\documentclass[9pt]{amsart}
\usepackage{amsmath}
\usepackage{amsxtra}
\usepackage{amscd}
\usepackage{amsthm}
\usepackage{amsfonts}
\usepackage{amssymb}
\usepackage{eucal}
\usepackage[all]{xy}
\usepackage{graphicx}
\usepackage[usenames]{color}

\newtheorem{cor}[subsubsection]{Corollary}
\newtheorem{maincor}[subsubsection]{Main Corollary}
\newtheorem{corconj}[subsubsection]{Corollary-of-Conjecture}
\newtheorem{lem}[subsubsection]{Lemma}
\newtheorem{prop}[subsubsection]{Proposition}

\newtheorem{conj}[subsubsection]{Conjecture}
\newtheorem{mainconj}[subsubsection]{Main Conjecture}

\newtheorem{thm}[subsubsection]{Theorem}
\newtheorem{mainthm}[subsubsection]{Main Theorem}

\newtheorem{defn}[subsubsection]{Definition}

\newtheorem{obs}[subsubsection]{Observation}
\newtheorem{hypoth}[subsubsection]{Hypothesis}

\theoremstyle{remark}
\newtheorem{rem}[subsubsection]{Remark}

\theoremstyle{definition}

\theoremstyle{remark}

\newcommand{\thmref}[1]{Theorem~\ref{#1}}

\newcommand{\secref}[1]{Sect.~\ref{#1}}
\newcommand{\lemref}[1]{Lemma~\ref{#1}}
\newcommand{\propref}[1]{Proposition~\ref{#1}}
\newcommand{\corref}[1]{Corollary~\ref{#1}}
\newcommand{\conjref}[1]{Conjecture~\ref{#1}}

\numberwithin{equation}{section}

\newcommand{\nc}{\newcommand}
\nc{\renc}{\renewcommand}
\nc{\ssec}{\subsection}
\nc{\sssec}{\subsubsection}
\nc{\on}{\operatorname}

\nc{\ips}{{\iota_P^{(S)}}}
\nc{\ipms}{{\iota_{P^-}^{(S)}}}
\nc{\sfpps}{{\sfp_P^{(S)}}}
\nc{\sfppms}{{\sfp_{P^-}^{(S)}}}

\nc\ol{\overline}
\nc\ul{\underline}
\nc\wt{\widetilde}
\nc\tboxtimes{\wt{\boxtimes}}
\nc\tstar{\wt{\star}}
\nc{\alp}{\alpha}

\nc{\ZZ}{{\mathbb Z}}
\nc{\NN}{{\mathbb N}}
\nc{\OO}{{\mathbb O}}
\renc{\SS}{{\mathbb S}}
\nc{\DD}{{\mathbb D}}
\nc{\GG}{{\mathbb G}}

\nc{\Fq}{{\mathbb F}_q}
\nc{\Fqb}{\ol{\mathbb F}_q}
\nc{\Ql}{{\mathbb Q}_\ell}
\nc{\Qlb}{{\ol{\mathbb Q}_\ell}}
\nc{\id}{\text{id}}
\nc\X{\mathcal X}

\nc{\red}{\on{red}}
\nc{\Ho}{\on{Ho}}
\nc{\Hom}{\on{Hom}}
\nc{\coHom}{\ul{\on{coHom}}}
\nc{\coMaps}{{\bf{coMaps}}}
\nc{\coef}{\on{coef}}
\nc{\Lie}{\on{Lie}}
\nc{\Loc}{\on{Loc}}
\nc{\Pic}{\on{Pic}}
\nc{\Bun}{\on{Bun}}
\nc{\IC}{\on{IC}}
\nc{\Aut}{\on{Aut}}
\nc{\rk}{\on{rk}}
\nc{\Sh}{\on{Sh}}
\nc{\Perv}{\on{Perv}}
\nc{\pos}{{\on{pos}}}
\nc{\Conv}{\on{Conv}}
\nc{\Sph}{\on{Sph}}
\nc{\Sym}{\on{Sym}}
\nc{\BunBb}{\overline{\Bun}_B}
\nc{\BunNb}{\overline{\Bun}_N}
\nc{\BunTb}{\overline{\Bun}_T}
\nc{\BunBbm}{\overline{\Bun}_{B^-}}
\nc{\BunBbel}{\overline{\Bun}_{B,el}}
\nc{\BunBbmel}{\overline{\Bun}_{B^-,el}}
\nc{\Buno}{\overset{o}{\Bun}}
\nc{\BunPb}{{\overline{\Bun}_P}}
\nc{\BunBM}{\Bun_{B(M)}}
\nc{\BunBMb}{\overline{\Bun}_{B(M)}}
\nc{\BunPbw}{{\widetilde{\Bun}_P}}
\nc{\BunBP}{\widetilde{\Bun}_{B,P}}
\nc{\GUb}{\overline{G/U}}
\nc{\GUPb}{\overline{G/U(P)}}

\nc{\Hhom}{\underline{\on{Hom}}}
\nc\syminfty{\on{Sym}^{\infty}}
\nc\lal{\ol{\lambda}}
\nc\xl{\ol{x}}
\nc\thl{\ol{\theta}}
\nc\nul{\ol{\nu}}
\nc\mul{\ol{\mu}}
\nc\Sum\Sigma
\nc{\oX}{\overset{o}{X}{}}
\nc{\hl}{\overset{\leftarrow}h{}}
\nc{\hr}{\overset{\rightarrow}h{}}
\nc{\M}{{\mathcal M}}
\nc{\N}{{\mathcal N}}
\nc{\F}{{\mathcal F}}
\nc{\D}{{\mathcal D}}
\nc{\Q}{{\mathcal Q}}
\nc{\Y}{{\mathcal Y}}
\nc{\G}{{\mathcal G}}
\nc{\E}{{\mathcal E}}
\nc{\CalC}{{\mathcal C}}
\nc\Dh{\widehat{\D}}

\nc{\C}{{\mathcal C}}
\nc{\K}{{\mathcal K}}
\renewcommand{\H}{{\mathcal H}}

\nc{\T}{{\mathcal T}}
\nc{\V}{{\mathcal V}}
\renc{\P}{{\mathcal P}}
\nc{\A}{{\mathcal A}}
\nc{\B}{{\mathcal B}}
\nc{\U}{{\mathcal U}}

\nc{\Gr}{{\on{Gr}}}

\nc{\frn}{{\check{\mathfrak u}(P)}}

\nc{\fC}{\mathfrak C}
\nc{\fT}{\mathfrak T}
\nc{\p}{\mathfrak p}
\nc{\q}{\mathfrak q}
\nc\f{{\mathfrak f}}

\nc{\qo}{{\mathfrak q}}
\nc{\po}{{\mathfrak p}}
\nc{\s}{{\mathfrak s}}
\nc\w{\text{w}}

\renewcommand{\mod}{{\on{-mod}}}

\nc\mathi\iota
\nc\Spec{\on{Spec}}
\nc\Proj{\on{Proj}}
\nc\Mod{\on{Mod}}
\nc{\tw}{\widetilde{\mathfrak t}}
\nc{\pw}{\widetilde{\mathfrak p}}
\nc{\qw}{\widetilde{\mathfrak q}}
\nc{\jw}{\widetilde j}

\nc{\grb}{\overline{\Gr}}
\nc{\I}{\mathcal I}

\nc{\lambdach}{{\check\lambda}}
\nc{\Lambdach}{{\check\Lambda}{}}
\nc{\much}{{\check\mu}}
\nc{\omegach}{{\check\omega}}
\nc{\nuch}{{\check\nu}}
\nc{\etach}{{\check\eta}}
\nc{\alphach}{{\check\alpha}}
\nc{\oblvtach}{{\check\oblvta}}
\nc{\rhoch}{{\check\rho}}
\nc{\ch}{{\check h}}

\nc{\Hb}{\overline{\H}}

\nc{\BA}{{\mathbb{A}}}
\nc{\BC}{{\mathbb{C}}}
\nc{\BE}{{\mathbb{E}}}
\nc{\BF}{{\mathbb{F}}}
\nc{\BG}{{\mathbb{G}}}
\nc{\BL}{{\mathbb{L}}}
\nc{\BM}{{\mathbb{M}}}
\nc{\BO}{{\mathbb{O}}}
\nc{\BD}{{\mathbb{D}}}
\nc{\BN}{{\mathbb{N}}}
\nc{\BP}{{\mathbb{P}}}
\nc{\BQ}{{\mathbb{Q}}}
\nc{\BR}{{\mathbb{R}}}
\nc{\BZ}{{\mathbb{Z}}}
\nc{\BS}{{\mathbb{S}}}
\nc{\Deep}{{\bf{deep}}}
\nc{\deep}{deep}

\nc{\CA}{{\mathcal{A}}}
\nc{\CB}{{\mathcal{B}}}

\nc{\CE}{{\mathcal{E}}}
\nc{\CF}{{\mathcal{F}}}
\nc{\CH}{{\mathcal{H}}}

\nc{\CL}{{\mathcal{L}}}
\nc{\CC}{{\mathcal{C}}}
\nc{\CG}{{\mathcal{G}}}
\nc{\CalD}{{\mathcal{D}}}
\nc{\CM}{{\mathcal{M}}}
\nc{\CN}{{\mathcal{N}}}
\nc{\CK}{{\mathcal{K}}}
\nc{\CO}{{\mathcal{O}}}
\nc{\CP}{{\mathcal{P}}}
\nc{\CQ}{{\mathcal{Q}}}
\nc{\CR}{{\mathcal{R}}}
\nc{\CS}{{\mathcal{S}}}
\nc{\CT}{{\mathcal{T}}}
\nc{\CU}{{\mathcal{U}}}
\nc{\CV}{{\mathcal{V}}}
\nc{\CW}{{\mathcal{W}}}
\nc{\CX}{{\mathcal{X}}}
\nc{\CY}{{\mathcal{Y}}}
\nc{\CZ}{{\mathcal{Z}}}
\nc{\CI}{{\mathcal{I}}}

\nc{\csM}{{\check{\mathcal A}}{}}
\nc{\oM}{{\overset{\circ}{\mathcal M}}{}}
\nc{\obM}{{\overset{\circ}{\mathbf M}}{}}
\nc{\oCA}{{\overset{\circ}{\mathcal A}}{}}
\nc{\obA}{{\overset{\circ}{\mathbf A}}{}}
\nc{\ooM}{{\overset{\circ}{M}}{}}
\nc{\osM}{{\overset{\circ}{\mathsf M}}{}}
\nc{\vM}{{\overset{\bullet}{\mathcal M}}{}}
\nc{\nM}{{\underset{\bullet}{\mathcal M}}{}}
\nc{\oD}{{\overset{\circ}{\mathcal D}}{}}
\nc{\obD}{{\overset{\circ}{\mathbf D}}{}}
\nc{\oA}{{\overset{\circ}{A}}{}}
\nc{\op}{{\overset{\bullet}{\mathbf p}}{}}
\nc{\cp}{{\overset{\circ}{\mathbf p}}{}}
\nc{\oU}{{\overset{\bullet}{\mathcal U}}{}}
\nc{\oZ}{{\overset{\circ}{\mathcal Z}}{}}
\nc{\ofZ}{{\overset{\circ}{\mathfrak Z}}{}}
\nc{\oF}{{\overset{\circ}{\fF}}}

\nc{\fa}{{\mathfrak{a}}}
\nc{\ofa}{\overset{\circ}{\mathfrak{a}}}
\nc{\fb}{{\mathfrak{b}}}
\nc{\fd}{{\mathfrak{d}}}
\nc{\ff}{{\mathfrak{f}}}
\nc{\fg}{{\mathfrak{g}}}
\nc{\fgl}{{\mathfrak{gl}}}
\nc{\fh}{{\mathfrak{h}}}
\nc{\fj}{{\mathfrak{j}}}
\nc{\fl}{{\mathfrak{l}}}
\nc{\fm}{{\mathfrak{m}}}
\nc{\ofm}{\overset{\circ}{\mathfrak{m}}}
\nc{\fn}{{\mathfrak{n}}}
\nc{\fu}{{\mathfrak{u}}}
\nc{\fp}{{\mathfrak{p}}}
\nc{\fr}{{\mathfrak{r}}}
\nc{\fs}{{\mathfrak{s}}}
\nc{\ft}{{\mathfrak{t}}}
\nc{\oft}{\overset{\circ}{\mathfrak{t}}}
\nc{\fz}{{\mathfrak{z}}}
\nc{\fsl}{{\mathfrak{sl}}}
\nc{\hsl}{{\widehat{\mathfrak{sl}}}}
\nc{\hgl}{{\widehat{\mathfrak{gl}}}}
\nc{\hg}{{\widehat{\mathfrak{g}}}}
\nc{\chg}{{\widehat{\mathfrak{g}}}{}^\vee}
\nc{\hn}{{\widehat{\mathfrak{n}}}}
\nc{\chn}{{\widehat{\mathfrak{n}}}{}^\vee}

\nc{\fA}{{\mathfrak{A}}}
\nc{\fB}{{\mathfrak{B}}}
\nc{\fD}{{\mathfrak{D}}}
\nc{\fE}{{\mathfrak{E}}}
\nc{\fF}{{\mathfrak{F}}}
\nc{\fG}{{\mathfrak{G}}}
\nc{\fK}{{\mathfrak{K}}}
\nc{\fL}{{\mathfrak{L}}}
\nc{\fM}{{\mathfrak{M}}}
\nc{\fN}{{\mathfrak{N}}}
\nc{\fP}{{\mathfrak{P}}}
\nc{\fU}{{\mathfrak{U}}}
\nc{\fV}{{\mathfrak{V}}}
\nc{\fZ}{{\mathfrak{Z}}}

\nc{\ba}{{\mathbf{a}}}
\nc{\bb}{{\mathbf{b}}}
\nc{\bc}{{\mathbf{c}}}
\nc{\bd}{{\mathbf{d}}}
\nc{\bbf}{{\mathbf{f}}}
\nc{\be}{{\mathbf{e}}}
\nc{\bi}{{\mathbf{i}}}
\nc{\bj}{{\mathbf{j}}}
\nc{\bh}{{\mathbf{h}}}
\nc{\bm}{{\mathbf{m}}}
\nc{\bn}{{\mathbf{n}}}
\nc{\bo}{{\mathbf{o}}}
\nc{\bp}{{\mathbf{p}}}
\nc{\bq}{{\mathbf{q}}}
\nc{\brr}{{\mathbf{r}}}
\nc{\bu}{{\mathbf{u}}}
\nc{\bv}{{\mathbf{v}}}
\nc{\bx}{{\mathbf{x}}}
\nc{\bs}{{\mathbf{s}}}
\nc{\by}{{\mathbf{y}}}
\nc{\bw}{{\mathbf{w}}}
\nc{\bA}{{\mathbf{A}}}
\nc{\bK}{{\mathbf{K}}}
\nc{\bB}{{\mathbf{B}}}
\nc{\bC}{{\mathbf{C}}}
\nc{\bG}{{\mathbf{G}}}
\nc{\bD}{{\mathbf{D}}}
\nc{\bE}{{\mathbf{E}}}
\nc{\bH}{{{\mathbf{H}}}}
\nc{\bM}{{\mathbf{M}}}
\nc{\bN}{{\mathbf{N}}}
\nc{\bO}{{\mathbf{O}}}
\nc{\bP}{{\mathbf{P}}}
\nc{\bQ}{{\mathbf{Q}}}
\nc{\bV}{{\mathbf{V}}}
\nc{\bW}{{\mathbf{W}}}
\nc{\bX}{{\mathbf{X}}}
\nc{\bZ}{{\mathbf{Z}}}
\nc{\bS}{{\mathbf{S}}}

\nc{\sA}{{\mathsf{A}}}
\nc{\sB}{{\mathsf{B}}}
\nc{\sC}{{\mathsf{C}}}
\nc{\sD}{{\mathsf{D}}}
\nc{\sF}{{\mathsf{F}}}
\nc{\sG}{{\mathsf{G}}}
\nc{\sH}{{\mathsf{H}}}
\nc{\sK}{{\mathsf{K}}}
\nc{\sM}{{\mathsf{M}}}
\nc{\sN}{{\mathsf{N}}}
\nc{\sO}{{\mathsf{O}}}
\nc{\sV}{{\mathsf{V}}}
\nc{\sW}{{\mathsf{W}}}
\nc{\sQ}{{\mathsf{Q}}}
\nc{\sP}{{\mathsf{P}}}
\nc{\sR}{{\mathsf{R}}}
\nc{\sT}{{\mathsf{T}}}
\nc{\sZ}{{\mathsf{Z}}}
\nc{\sfp}{{\mathsf{p}}}
\nc{\sfq}{{\mathsf{q}}}
\nc{\sft}{{\mathsf{t}}}
\nc{\sr}{{\mathsf{r}}}
\nc{\bk}{{\mathsf{k}}}
\nc{\sa}{{\mathsf{s}}}
\nc{\sg}{{\mathsf{g}}}
\nc{\sn}{{\mathsf{n}}}
\nc{\sh}{{\mathsf{h}}}
\nc{\sff}{{\mathsf{f}}}
\nc{\sfb}{{\mathsf{b}}}
\nc{\sfc}{{\mathsf{c}}}
\nc{\sfe}{{\mathsf{e}}}
\nc{\sd}{{\mathsf{d}}}

\nc{\BK}{{\bar{K}}}

\nc{\tA}{{\widetilde{\mathbf{A}}}}
\nc{\tB}{{\widetilde{\mathcal{B}}}}
\nc{\tg}{{\widetilde{\mathfrak{g}}}}
\nc{\tG}{{\widetilde{G}}}
\nc{\TM}{{\widetilde{\mathbb{M}}}{}}
\nc{\tO}{{\widetilde{\mathsf{O}}}{}}
\nc{\tU}{{\widetilde{\mathfrak{U}}}{}}
\nc{\TZ}{{\tilde{Z}}}
\nc{\tx}{{\tilde{x}}}
\nc{\tbv}{{\tilde{\bv}}}
\nc{\tfP}{{\widetilde{\mathfrak{P}}}{}}
\nc{\tz}{{\tilde{\zeta}}}
\nc{\tmu}{{\tilde{\mu}}}

\nc{\urho}{\underline{\rho}}
\nc{\uB}{\underline{B}}
\nc{\uC}{{\underline{\mathbb{C}}}}
\nc{\ui}{\underline{i}}
\nc{\uj}{\underline{j}}
\nc{\ofP}{{\overline{\mathfrak{P}}}}
\nc{\oB}{{\overline{\mathcal{B}}}}
\nc{\og}{{\overline{\mathfrak{g}}}}
\nc{\oI}{{\overline{I}}}

\nc{\eps}{\varepsilon}
\nc{\hrho}{{\hat{\rho}}}

\nc{\one}{{\mathbf{1}}}
\nc{\two}{{\mathbf{t}}}

\nc{\Rep}{{\mathop{\operatorname{\rm Rep}}}}
\nc{\Tot}{{\mathop{\operatorname{\rm Tot}}}}
\nc{\Ker}{{\mathop{\operatorname{\rm Ker}}}}
\nc{\im}{{\mathop{\operatorname{\rm Im}}}}
\nc{\Hilb}{{\mathop{\operatorname{\rm Hilb}}}}
\nc{\End}{{\mathop{\operatorname{\rm End}}}}
\nc{\Ext}{{\mathop{\operatorname{\rm Ext}}}}
\nc{\CHom}{{\mathop{\operatorname{{\mathcal{H}}\it om}}}}
\nc{\CEnd}{{\mathop{\operatorname{{\mathcal{E}}\it nd}}}}
\nc{\GL}{{\mathop{\operatorname{\rm GL}}}}
\nc{\gr}{{\mathop{\operatorname{\rm gr}}}}
\nc{\HN}{{\mathop{\operatorname{\rm HN}}}}
\nc{\Id}{{\mathop{\operatorname{\rm Id}}}}
\nc{\de}{{\mathop{\operatorname{\rm def}}}}
\nc{\length}{{\mathop{\operatorname{\rm length}}}}
\nc{\supp}{{\mathop{\operatorname{\rm supp}}}}

\nc{\Cliff}{{\mathsf{Cliff}}}
\nc{\Fl}{\on{Fl}}
\nc{\Fib}{{\mathsf{Fib}}}
\nc{\Coh}{{\on{Coh}}}
\nc{\QCoh}{{\on{QCoh}}}
\nc{\IndCoh}{{\on{IndCoh}}}
\nc{\FCoh}{{\mathsf{FCoh}}}

\nc{\reg}{{\text{\rm reg}}}

\nc{\cplus}{{\mathbf{C}_+}}
\nc{\cminus}{{\mathbf{C}_-}}
\nc{\cthree}{{\mathbf{C}_\bullet}}
\nc{\Qbar}{{\bar{Q}}}
\nc\Eis{\on{Eis}}
\nc\Eisb{\ol\Eis{}}
\nc\Eisr{\on{Eis}^{rat}{}}
\nc\wh{\widehat}
\nc{\Def}{\on{Def_{\check{\fb}}(E)}}
\nc{\barZ}{\overline{Z}{}}
\nc{\barbarZ}{\overline{\barZ}{}}
\nc{\barpi}{\overline\pi}
\nc{\barbarpi}{\overline\barpi}
\nc{\barpip}{\overline\pi{}^+}
\nc{\barpim}{\overline\pi{}^-}

\nc{\fq}{\mathfrak q}

\nc{\fqb}{\ol{\sfq}{}}
\nc{\fpb}{\ol{\sfp}{}}
\nc{\fpr}{{\sfp^{rat}}{}}
\nc{\fqr}{{\sfq^{rat}}{}}

\nc{\hattimes}{\wh\otimes}

\nc{\bOmega}{{\overline{\Omega(\check \fn)}}}

\nc{\seq}[1]{\stackrel{#1}{\sim}}

\nc{\cT}{{\check{T}}}
\nc{\cG}{{\check{G}}}
\nc{\cM}{{\check{M}}}
\nc{\cB}{{\check{B}}}

\nc{\ct}{{\check{\mathfrak t}}}
\nc{\cg}{{\check{\fg}}}
\nc{\cb}{{\check{\fb}}}
\nc{\cn}{{\check{\fn}}}

\nc{\cLambda}{{\check\Lambda}}

\nc{\cla}{{\check\lambda}}
\nc{\cmu}{{\check\mu}}
\nc{\cnu}{{\check\nu}}
\nc{\ceta}{{\check\eta}}

\nc{\DefbE}{{\on{Def}_{\cB}(E_\cT)}}

\nc{\imathb}{{\ol{\imath}}}
\nc{\rlr}{\overset{\longrightarrow}{\underset{\longrightarrow}\longleftarrow}}

\nc{\oBun}{\overset{\circ}\Bun}
\nc{\LocSys}{\on{LocSys}}
\nc{\BunBbb}{\ol{\ol{Bun}}_B}
\nc{\BunBr}{\Bun_B^{rat}}
\nc{\BunBrsg}{\Bun_B^{rat,\on{s.g.}}}
\nc{\BunBrp}{\Bun_B^{rat,polar}}
\nc{\BunBrpbg}{\Bun_B^{rat,polar,\on{b.g.}}}
\nc{\BunBrpsg}{\Bun_B^{rat,polar,\on{s.g.}}}
\nc{\BunTrp}{\Bun_T^{rat,polar}}
\nc{\BunTrpbg}{\Bun_T^{rat,polar,\on{b.g.}}}
\nc{\BunTrpsg}{\Bun_T^{rat,polar,\on{s.g.}}}
\nc{\BunNr}{\Bun_N^{rat}}
\nc{\BunNre}{\Bun_N^{enh,rat}}
\nc{\BunTr}{\Bun_T^{rat}}
\nc{\Vect}{\on{Vect}}
\nc{\Whit}{\on{Whit}}
\nc{\CTb}{\ol{\on{CT}}}
\nc{\Ran}{{\on{Ran}}}
\nc{\fSet}{{\on{fSet}}}
\nc{\CTr}{\on{CT}^{rat}{}}
\nc\jmathr{\jmath^{rat}{}}
\nc{\ux}{\underline{x}}
\nc{\clambda}{{\check\lambda}}
\nc{\calpha}{{\check\alpha}}
\nc{\ind}{{\mathbf{ind}}}
\nc{\coinv}{{\mathbf{coinv}}}
\nc{\inv}{{\mathbf{inv}}}
\nc{\oblv}{{\mathbf{oblv}}}
\nc{\free}{{\mathbf{free}}}
\nc{\ox}{{\overline{x}}}
\nc{\cLa}{\check{\Lambda}}
\nc{\StinftyCat}{\on{DGCat}}
\nc{\inftyCat}{\infty\on{-Cat}}
\nc{\inftygroup}{\infty\on{-Grpd}}
\nc{\Dmod}{\on{D-mod}}
\nc{\CMaps}{{\mathcal Maps}}
\nc{\Maps}{\on{Maps}}
\nc{\affSch}{\on{Sch}^{\on{aff}}}
\nc{\dr}{{\on{dR}}}
\nc{\oCF}{\overset{\circ}\CF}
\nc{\oCY}{\overset{\circ}\CY}
\nc{\opi}{\overset{\circ}\pi}
\nc{\leqG}{\underset{G}\leq}
\nc{\leqM}{\underset{M}\leq}
\nc{\leqGad}{\underset{G_{ad}}\leq}
\nc{\leqMad}{\underset{M_{ad}}\leq}
\nc{\Tr}{\on{Tr}}
\nc{\Frob}{{\on{Frob}}}
\nc{\DGCat}{\on{DGCat}}
\nc{\tDGCat}{2\on{-DGCat}_{\on{u.g.}}}
\nc{\ev}{\on{ev}}
\nc{\mmod}{\on{-}\mathbf{mod}}
\nc{\sotimes}{\overset{!}\otimes}
\nc{\Shv}{\on{Shv}}
\nc{\Spc}{\on{Spc}}
\nc{\LS}{\qLisse}
\nc{\Res}{\on{Res}}
\nc{\bDelta}{{\mathbf{\Delta}}}
\nc{\bMaps}{{\mathbf{Maps}}}
\nc{\cD}{\mathcal D}
\nc{\ocD}{\overset{\circ}\cD}
\nc{\ppart}{(\!(t)\!)}
\nc{\qqart}{[\![t]\!]}
\nc{\oCU}{\overset{\circ}{\CU}}
\nc{\Exc}{{\mathcal{E}xc}}
\nc{\Sht}{\on{Sht}}
\nc{\Nilp}{{\on{Nilp}}}
\nc{\Drinf}{\on{Drinf}}
\nc{\Sing}{\on{Sing}}
\nc{\IndLisse}{\qLisse}
\nc{\Shvl}{\on{Shv}_{\{0\}}} 
\nc{\Lisse}{\on{Lisse}}
\nc{\qLisse}{\on{QLisse}}
\nc{\iLisse}{\on{IndLisse}}
\nc{\hbarp}{{k[\hbar]}}
\nc{\hbarlp}{{k[\hbar,\hbart^{-1}]}}
\nc{\hbart}{{k[\![\hbar]\!]}} 
\nc{\hbarl}{{k(\!(\hbar)\!)}} 
\nc{\bPhi}{{\mathbf \Phi}}
\nc{\bPsi}{{\mathbf \Psi}}

\begin{document}


\title[Geometric Langlands with nilpotent singular support]{The stack of local systems with restricted variation and \\
geometric Langlands theory with nilpotent singular support}

\author{D.~Arinkin, D.~Gaitsgory, D.~Kazhdan, 
S.~Raskin, N.~Rozenblyum, Y.~Varshavsky}

\begin{abstract}
We define a new geometric object--the stack of local systems with restricted variation. We formulate a version of the 
categorical geometric Langlands conjecture that makes sense for any constructible sheaf theory (such as $\ell$-adic sheaves). 
We formulate a conjecture that makes precise the connection between the category of automorphic sheaves and the space of automorphic functions.
\end{abstract}

\date{\today}

\maketitle

\tableofcontents

\section*{Introduction}

\ssec{Starting point}


Classically, Langlands proposed a framework for understanding irreducible
automorphic representations for a reductive group $G$ via spectral data involving
the dual group $\cG$. 

\medskip

P.~Deligne (for $GL_1$), V.~Drinfeld (for $GL_2$) and G.~Laumon (for $GL_n$)
realized Langlands-style phenomena in algebraic geometry. 
In their setting, the fundamental
objects of interest are \emph{Hecke eigensheaves}.
This theory works over an arbitrary ground field $k$, and takes
as an additional input a \emph{sheaf theory} for varieties over that
field. 
Specializing to $k = \overline{\BF}_q$ and \'etale sheaves,
one recovers special cases of Langlands's conjectures by taking 
the trace of Frobenius. 

\medskip

Inspired by these works, Beilinson and Drinfeld proposed
the \emph{categorical Geometric Langlands Conjecture}
\begin{equation} \label{e:BD GLC}
\Dmod(\Bun_G) \simeq \IndCoh_\Nilp(\LocSys_\cG(X))
\end{equation} 
for $X$ a smooth projective curve over a field $k$ of 
characteristic zero. Here the left-hand side $\Dmod(\Bun_G)$ is a 
sheaf-theoretic analogue of the space of unramified automorphic functions,
and the right hand side is defined in \cite{AG}.

\medskip

There are (related) discrepancies between this categorical conjecture
and more classical conjectures.

\medskip 

\noindent(i) Hecke eigensheaves make sense in any sheaf theory, while 
the Beilinson-Drinfeld conjecture applies only in the setting of D-modules.

\medskip 

\noindent(ii) Hecke eigensheaves categorify the arithmetic 
Langlands correspondence through the trace of Frobenius construction,
while the Beilinson-Drinfeld conjecture bears no 
direct relation to automorphic functions.

\medskip

\noindent(iii) Langlands's conjecture considers \emph{irreducible} automorphic
representations, while the Beilinson-Drinfeld conjecture provides
a spectral decomposition of (a sheaf-theoretic
analogue of) the whole space of (unramified) automorphic functions.

\medskip 

These differences provoke natural questions:

\medskip

\noindent--Is there a categorical geometric Langlands conjecture that
applies in any sheaf-theoretic context, in particular, in the 
\'etale setting over finite fields? 

\medskip

\noindent--The trace construction attaches automorphic functions to
particular \'etale sheaves on $\Bun_G$; is there a direct relationship between the \emph{category} 
of \'etale sheaves on $\Bun_G$ and the \emph{space} of automorphic functions?

\medskip

\noindent--Is it possible to give a spectral description of the space 
of classical automorphic functions, not merely its irreducible constituents? 

\ssec{Summary}

\sssec{}\label{sss:yes}

In this paper, we 
provide positive answers to the three questions raised above.

\medskip 

\noindent Our \conjref{c:restr GLC} provides an analogue of 
the categorical Geometric Langlands Conjecture that is suited
to any ground field and any sheaf theory. 

\medskip 

\noindent Our \conjref{c:Trace conj} proposes a closer relationship between
sheaves on $\Bun_G$ and unramified automorphic functions than 
was previously considered. As such, it allows one to extract new,
concrete conjectures on automorphic functions from our categorical
Geometric Langlands Conjecture, see right below. 

\medskip 

\noindent Our \conjref{c:autom-omega} describes
the space of unramified automorphic functions over a function field
in spectral terms, refining Langlands's conjectures in this setting. 

\medskip 

In sum, the main purpose of this work is to propose
a variant of the categorical Beilinson-Drinfeld conjecture that makes
sense over finite fields, and in that setting, 
to connect it with the arithmetic Langlands program.

\sssec{}

This paper contains two main ideas. The first of them is the introduction of a space 
$$\LocSys^{\on{restr}}_\sG(X)$$ 
of \emph{$\sG$-local systems with restricted variation} 
on $X$. In our \conjref{c:restr GLC}, 
$\LocSys^{\on{restr}}_\cG(X)$ replaces 
$\LocSys_\cG(X)$ from the original conjecture of Beilinson and Drinfeld.

\medskip 

We discuss $\LocSys^{\on{restr}}_\sG(X)$ 
in detail later in the introduction. For now, let us
admit it into the discussion as a black box. 

\medskip 

Then our \conjref{c:restr GLC} asserts
\begin{equation} \label{e:Nilp GLC intro}
\Shv_\Nilp(\Bun_G) \simeq \IndCoh_\Nilp(\LocSys^{\on{restr}}_\cG(X)),
\end{equation} 
where the left-hand side is the category of
ind-constructible sheaves on $\Bun_G$ with nilpotent singular support;
we study this category in detail in \secref{s:Nilp}.

\sssec{}

In addition, we make some progress toward 
\conjref{c:restr GLC}.

\medskip 

Our \thmref{t:spectral decomp} 
provides an action of $\QCoh(\LocSys^{\on{restr}}_\cG(X))$
on $\Shv_\Nilp(\Bun_G)$ compatible with Hecke functors. 
We regard this result as a 
spectral decomposition of the category $\Shv_\Nilp(\Bun_G)$ over
$\LocSys^{\on{restr}}_\cG(X)$.
This theorem is a counterpart of \cite[Corollary 4.5.5]{Ga5}, which applied in the D-module setting
and whose proof used completely different methods.

\medskip

Using these methods, we settle long-standing conjectures
on the structure of Hecke eigensheaves. 
Our \corref{c:eigensheaves nilp} shows that Hecke eigensheaves
have nilpotent singular support, as predicted by
G.~Laumon in \cite[Conjecture 6.3.1]{laumon}.
In addition, our \corref{c:RS}
shows that in the D-module setting, 
any Hecke eigensheaf has regular singularities, as predicted
by Beilinson-Drinfeld in \cite[Sect. 5.2.7]{BD1}.

\sssec{}

The second main idea of this paper is that of categorical trace. It appears
in our \conjref{c:Trace conj}, which we title the \emph{Trace Conjecture}.
This conjecture predicts a stronger link between
geometric and arithmetic Langlands than was previously 
considered: 

\medskip

Suppose $k=\ol\BF_q$ and that $X$ and $G$ are defined over $\BF_q$, 
and therefore carry 
Frobenius endomorphisms. The Trace Conjecture asserts
that the categorical trace of the functor $(\Frob_{\Bun_G})_*$
on $\Shv_\Nilp(\Bun_G)$ maps isomorphically to the space of 
(compactly supported) unramified automorphic functions
$$\on{Autom}:=\on{Funct}_c(\Bun_G(\BF_q)).$$

\medskip

More evocatively: we conjecture that
a trace of Frobenius construction recovers the \emph{space} of automorphic
forms from the \emph{category} $\Shv_\Nilp(\Bun_G)$, 
much as one classically extracts a automorphic functions
from an automorphic sheaf by a trace of Frobenius construction.

\medskip 

Combined with our \thmref{t:spectral decomp}, the Trace Conjecture 
gives rise to the spectral decomposition of $\on{Autom}$ along the set of isomorphism
classes of semi-simple Langlands parameters, recovering the (unramified case of)
V.~Lafforgue's result. 

\medskip

Moreover, if we combine the Trace Conjecture with our version of the categorical
Geometric Langlands Conjecture (i.e., \conjref{c:restr GLC}), we obtain a full
description of the space of (unramified) automorphic functions in terms of Langlands
parameters (and not just the spectral decomposition): 
$$\on{Autom} \simeq 
\Gamma(\LocSys^{\on{arithm}} _\cG(X),\omega_{\LocSys^{\on{arithm}} _\cG(X)}),$$
where $\LocSys^{\on{arithm}} _\cG(X)$ is the algebraic stack of Frobenius-fixed points, i.e., 
$$\LocSys^{\on{arithm}} _\cG(X):=(\LocSys^{\on{restr}}_\cG(X))^\Frob,$$
where $\Frob$ is the automorphism of $\LocSys^{\on{restr}}_\cG(X)$
induced by the geometric Frobenius on $X$. 
This is our \conjref{c:autom-omega}, as referenced in \secref{sss:yes}.

\ssec{Some antecedents}

Before discussing the contents of this paper in more detail,
we highlight two points that are \emph{not} original to our work.

\sssec{Work of Ben-Zvi and Nadler}

Observe that in \conjref{c:restr GLC} we consider the subcategory
$\Shv_\Nilp(\Bun_G) \subset \Shv(\Bun_G)$ of 
(ind-constructible) sheaves with
nilpotent singular support, a hypothesis with no counterpart
in the Beilinson-Drinfeld setting of D-modules. 

\medskip 

The idea of considering this subcategory, which is so crucial to our work, is due to 
D.~Ben-Zvi and D.~Nadler, who did so in their setting 
of \emph{Betti} Geometric Langlands, see \cite{BN}. 

\sssec{}

Let us take a moment to clarify the relationship between our work and \cite{BN}.

\medskip 

For $k = \BC$, Ben-Zvi and Nadler consider the larger category 
$\Shv^{\on{all}}(\Bun_G)$ of \emph{all} (possibly 
not ind-constructible) sheaves on 
$\Bun_G(\BC)$, considered as a complex stack via its analytic topology.
Let $\Shv_\Nilp^{\on{all}}(\Bun_G)\subset \Shv^{\on{all}}(\Bun_G)$ be the full
subcategory consisting of objects with nilpotent singular support. Ben-Zvi and Nadler
conjectured an equivalence
\begin{equation} \label{e:Betti GLC}
\Shv^{\on{all}}_{\Nilp}(\Bun_G)\simeq \IndCoh_\Nilp(\LocSys_\cG(X)),
\end{equation} 
where in the right-hand side $\LocSys_\cG(X)$ is the Betti version of the stack of $\cG$-local
systems on $X$. 

\medskip

Let us compare this conjectural equivalence with the Beilinson-Drinfeld version \eqref{e:BD GLC}.
The latter is particular to D-modules, while \eqref{e:Betti GLC} is particular to topological sheaves. 
Our \eqref{e:Nilp GLC intro} sits in the middle between the two: when $k=\BC$ our $\Shv_\Nilp(\Bun_G)$
can be thought of as a full subcategory of both $\Dmod(\Bun_G)$ and $\Shv^{\on{all}}_{\Nilp}(\Bun_G)$.

\medskip

Similarly, our $\LocSys^{\on{restr}}_\cG(X)$ is an algebro-geometric object that is embedded into both
the de Rham and Betti versions of $\LocSys^{\on{restr}}_\cG(X)$. Now, the point is that $\LocSys^{\on{restr}}_\cG(X)$
can be defined abstractly, so that it makes sense in any sheaf-theoretic context, along with the conjectural 
equivalence \eqref{e:Nilp GLC intro}. 

\begin{rem} 
We should point out another source of initial evidence towards the relationship between $\Shv_\Nilp(\Bun_G)$
and $\LocSys^{\on{restr}}_\cG(X)$: 

\medskip

It was a discovered by D.~Nadler and Z.~Yun in \cite{NY1} that when we apply Hecke functors to objects from $\Shv_\Nilp(\Bun_G)$, 
we obtain objects in $\Shv(\Bun_G\times X)$ that \emph{behave like local systems along $X$}; see 
\thmref{t:NY} for a precise assertion.

\end{rem} 

\begin{rem}

We should also emphasize that what enabled us to even talk about $\Shv_\Nilp(\Bun_G)$ in the
context of $\ell$-adic sheaves was the work of A.~Beilinson \cite{Be2} and T.~Saito \cite{Sai}, where the singular support
of \'etale sheaves over any ground field was defined and studied. 

\end{rem} 

\sssec{Work of V.~Lafforgue}

Our Trace Conjecture is inspired by the
work \cite{VLaf1} of V.~Lafforgue on the arithmetic Langlands
correspondence for function fields.

\medskip

A distinctive feature of Geometric Langlands is that
Hecke functors are defined not merely at points $x \in X$
of a curve, but extend over all of $X$, and moreover,
extend over $X^I$ for any finite set $I$. These considerations
lead to the distinguished role played by the \emph{factorization algebras} of 
\cite{BD2} and \emph{Ran space} in geometric Langlands theory.

\medskip 

In his work, V.~Lafforgue showed that the existence of Hecke \emph{functors} over powers
of a curve has implications for automorphic \emph{functions}.
Specifically, he used the existence of these 
functors to construct \emph{excursion operators},
and used these excursion operators to define the spectral decomposition
of automorphic functions (over function fields) as predicted by the Langlands conjectures.

\sssec{}

In \cite{GKRV}, a subset of the authors of this paper attempted to reinterpret V.~Lafforgue's constructions using categorical 
traces. It provided a toy model for the spectral decomposition in \cite{VLaf1} in the following sense:

\medskip

In {\it loc.cit.} one starts with an abstract category $\CC$ equipped with an action of 
Hecke functors \emph{in the Betti setting} and an endofunctor $\Phi:\CC\to \CC$ (to be thought of
as a prototype of Frobenius), and obtains a spectral decomposition of the vector space
$\Tr(\Phi,\CC)$ along a certain space, which could be thought of as a Betti analog of the coarse
moduli space of arithmetic Langlands parameters. 

\medskip

Now, the present work allows to carry the construction of \cite{GKRV} in the actual setting of 
applicable to the study of automorphic functions: we take our $\CC$ to be the $\ell$-adic version of
$\Shv_\Nilp(\Bun_G)$ (for a curve $X$ over $\ol\BF_q$). 

\medskip 

In \secref{s:spectral}, we revisit V.~Lafforgue's work,
and show how our Trace Conjecture recovers and (following
ideas of V.~Drinfeld) refines
the main results of \cite{VLaf1} in the unramified case.

\ssec{Contents}

This paper consists of four parts. 

\medskip

In Part I we define and study the properties of the stack $\LocSys^{\on{restr}}_\cG(X)$.

\medskip

In Part II we establish a general spectral decomposition result that produces an action
of $\QCoh(\LocSys^{\on{restr}}_\cG(X))$ on a category $\bC$, equipped with what one
can call a \emph{lisse Hecke action}. 

\medskip

In Part III we study the properties of the category $\Shv_\Nilp(\Bun_G)$. We should say right away that
in this Part we prove two old-standing conjectures: that Hecke eigensheaves have nilpotent singular support,
and that (in the case of D-modules) all sheaves with nilpotent singular support have regular singularities. 
 
\medskip

In Part IV we study the applications of the theory developed hereto to the
Langlands theory.

\medskip

Below we will review the main results of each of the Parts.

%
%

\ssec{Overview: the stack $\LocSys^{\on{restr}}_\sG(X)$}

Let $\sG$ be an arbitrary affine algebraic group over a field of coefficients $\sfe$ of
characteristic $0$. 

\sssec{}

Let us start by recalling the definition of the (usual) algebraic stack $\LocSys^{\on{Betti}}_\sG(X)$ of
$\sG$-local systems on $X$ in the context of sheaves in the classical topology
(to be referred to as the \emph{Betti} context). 

\medskip

On the first pass, let us take $\sG=GL_n$. 

%

\medskip

Choose a base point $x\in X$. For an affine test scheme $S=\Spec(A)$ over $\sfe$, an $S$-point of 
$\LocSys_{GL_n}(X)$ is an $A$-module $E_S$, locally free of rank $n$, equipped with an action of 
$\pi_1(X,x)$.

\medskip

For an arbitrary $\sG$, the definition is obtained from the one for $GL_n$ via Tannakian formalism.

\sssec{}

We now give the definition of $\LocSys^{\on{restr}}_{GL_n}(X)$, still in the Betti context. Namely 
$\LocSys^{\on{restr}}_{GL_n}(X)$ is a subfunctor of $\LocSys_{GL_n}(X)$ that corresponds to the
following condition: 

\medskip

We require that the action of $\pi_1(X,x)$ on $E_S$ be $\sfe$-locally finite, i.e., each element of $E_S$
is contained in a finite-dimensional $\sfe$-vector subspace, preserves by the action of $\pi_1(X,x)$. 

\medskip

For an arbitrary $\sG$, one imposes this condition for each finite-dimensional representation
$\sG\to GL_n$ (or, equivalently, for one faithful representation). 

\medskip

When $A$ is Artinian, the above condition is automatic, so the formal completions of 
$\LocSys_\sG(X)$ and $\LocSys^{\on{restr}}_\sG(X)$ at any point are the same.
The difference appears for $A$ that have positive Krull dimension. 

\medskip

With that we should emphasize
that $\LocSys^{\on{restr}}_\sG(X)$ is \emph{not} entirely formal, i.e., it is \emph{not}
true that any $S$-point of $\LocSys^{\on{restr}}_\sG(X)$
factors though an $S'$-point with $S'$ Artinian. For example, for
$\sG=\BG_a$, the map 
$$\LocSys^{\on{restr}}_\sG(X)\to \LocSys_\sG(X)$$
is an isomorphism. 

\sssec{}

Let us explain the terminology ``restricted variation", again in the example of $\sG=GL_n$,

\medskip

The claim is that when we move along $S=\Spec(A)$, the 
corresponding representation of $\pi_1(X,x)$ does not change too much, in the sense that
the isomorphism class of its semi-simplification is constant (as long as $S$ is connected).

\medskip

Indeed, let us show that for every $\gamma\in \pi_1(X,x)$ and every $\lambda\in \sfe$,
the generalized $\lambda$-eigenspace of $\gamma$ on $E_s:=E_S\underset{A,s}\otimes \sfe$ 
has a constant dimension as $s$ moves along $S$. 

\medskip

Indeed, due to the locally finiteness condition, we can decompose $E_S$ into a direct
sum of generalized eigenspaces for $\gamma$
$$E_S=\underset{\lambda}\oplus\, E_S^{(\lambda)},$$
where each $E_S^{(\lambda)}$ is an $A$-submodule, and being a direct summand of 
a locally free $A$-module, it is itself locally free. 

\medskip

The same phenomenon will happen for any $\sG$: 
an $S$-point of $\LocSys_\sG(X)$ factors through $\LocSys^{\on{restr}}_\sG(X)$
if and only for all $\sfe$-points of $S$, the resulting $\sG$-local systems on $X$ all
have the same semi-simplification. 

\sssec{}

We are now ready to give the general definition of $\LocSys^{\on{restr}}_\sG(X)$. 

\medskip

Within the given sheaf theory, we consider the full subcategory 
$$\Lisse(X)\subset \Shv(X)^{\on{constr}}$$
of local systems (of finite rank).

\medskip

Consider its ind-completion, denoted $\iLisse(X)$. Finally, let $\qLisse(X)$ be
the left completion of $\iLisse(X)$ in the natural t-structure\footnote{The last step of
left completion is unnecessary if $X$ is a \emph{categorical $K(\pi,1)$}, see \secref{sss:Kpi1}, which is
the case of curves of genus $>0$. However, left completion is \emph{non-trivial} for $X=\BP^1$, i.e.,
$\iLisse(\BP^1)\neq \qLisse(\BP^1)$, see \secref{sss:analyze P1}}.
Now, for an
affine test scheme $S=\Spec(A)$, an $S$-point of $\LocSys^{\on{restr}}_\sG(X)$ is a symmetric 
monoidal functor
$$\Rep(\sG)\to A\mod\otimes \qLisse(X),$$
required to be right t-exact with respect to the natural t-structures. 

\medskip

By definition, $\sfe$-points of $\LocSys^{\on{restr}}_\sG(X)$ are just $\sG$-local systems on $X$. 

\medskip

Two remarks are in order:

\medskip

\noindent{(i)} In the definition of $\LocSys^{\on{restr}}_\sG(X)$ one can (and should!)
allow $S$ to be a \emph{derived} affine scheme over $\sfe$ (i.e., we allow $A$ to be a
connective commutative $\sfe$-algebra). Thus, $\LocSys^{\on{restr}}_\sG(X)$ is inherently
an object of derived algebraic geometry\footnote{In fact, our definition of the usual
$\LocSys^{\on{Betti}}_\sG(X)$ in the Betti context was a bit of a euphemism: for the correct definition 
in the context of derived algebraic geometry, one has to use the entire fundamental
groupoid of $X$, and not just $\pi_1$; the difference does not matter, however, when 
we evaluate on classical test affine schemes, while the distinction between  
$\LocSys^{\on{Betti}}_\sG(X)$ and $\LocSys^{\on{restr}}_\sG(X)$ happens at the classical level.}.

\medskip

\noindent{(ii)} The definition of $\LocSys^{\on{restr}}_\sG(X)$ uses the \emph{large}
category $A\mod\otimes \qLisse(X)$. When we evaluate our functor on \emph{truncated}
affine schemes, we can replace $\qLisse(X)$ by $\iLisse(X)=\on{Ind}(\Lisse(X))$
(see \propref{p:replace by Ind-Lisse}), and so we can express the definition in terms of small categories. But for an arbitrary $S$,
it is essential to work with the entire $\qLisse(X)$, to ensure \emph{convergence} (see \secref{ss:convergence}). 

\sssec{} 

As defined above, $\LocSys^{\on{restr}}_\sG(X)$ is just a functor on (derived) affine schemes,
so is just a prestack. But what kind of prestack is it? I.e., can we say something about its geometric
properties?  

\medskip

The majority of Part I is devoted to investigating this question. 

\medskip

While the geometric properties we find are exotic, 
this study plays a key role in Part III, where the geometry of $\LocSys^{\on{restr}}_\sG(X)$ (for $\sG=\cG$,
the Langlands dual of $G$) has concrete consequences for the category the category $\Shv_\Nilp(\Bun_G)$. 

\sssec{} \label{sss:Intro Betti}

First, let us illustrate the shape that $\LocSys^{\on{restr}}_\sG(X)$ has in the Betti context. Recall
that in this case we have the usual moduli stack $\LocSys^{\on{Betti}}_\sG(X)$, which is a quotient
of the affine scheme $\LocSys^{\on{Betti},\on{rigid}_x}_\sG(X)$ (that classifies local systems with
a trivialization at $x$) by $\sG$. 

\medskip

Assume that $\sG$ is reductive, and let 
$$\LocSys^{\on{Betti,coarse}}_\sG(X):=\LocSys^{\on{Betti},\on{rigid}_x}_\sG(X)/\!/\sG:=
\Spec(\Gamma(\LocSys^{\on{Betti}}_\sG(X),\CO_{\LocSys^{\on{Betti}}_\sG(X)}))$$
be the corresponding coarse moduli space. We have the tautological map
\begin{equation} \label{e:r coarse}
\brr:\LocSys^{\on{Betti}}_\sG(X)\to \LocSys^{\on{coarse}}_\sG(X),
\end{equation}
and recall that two $\sfe$-points of $\LocSys^{\on{Betti}}_\sG(X)$ lie in the same fiber of this map
if and only if they have isomorphic semi-simplifications.

\medskip

We can describe $\LocSys^{\on{restr}}_\sG(X)$ as the disjoint union
of formal completions of the fibers of $\brr$ over $\sfe$-points of $\LocSys^{\on{Betti},\on{coarse}}_\sG(X)$
(see \thmref{t:coarse fibers}). 

\medskip

In particular, we note one thing that $\LocSys^{\on{restr}}_\sG(X)$ is \emph{not}: it is \emph{not}
an algebraic stack (or union of such), because it has all these formal directions. 

\begin{rem} 

The above explicit description of $\LocSys^{\on{restr}}_\sG(X)$ in the Betti case 
may suggest that it is in general a ``silly" object. Indeed, why would we want
a moduli space in which all irreducible local systems belong to different connected components? 

\medskip

However,
as the results in Parts III and IV of this paper show, $\LocSys^{\on{restr}}_\sG(X)$ is actually a natural object
to consider, in that it is perfectly adapted to the study of $\Shv_\Nilp(\Bun_G)$, and thereby to applications
to the arithmetic theory. 

\medskip

For example, formula \eqref{e:dir sum Intro} below is the reflection on the
automorphic side of the above decomposition of $\LocSys^{\on{restr}}_\sG(X)$ as a disjoint union. 
See also \eqref{e:GLC Intro} for a version of the Geometric Langlands Conjecture with nilpotent 
singular support. Finally, see formula \eqref{e:Autom Intro} for an expression for the space of
automorphic functions in terms of Frobenius-fixed locus on $\LocSys^{\on{restr}}_\sG(X)$. 

\medskip

Looked at from a different angle, in the Betti and de Rham contexts, there are the ``honest" moduli 
spaces of local systems, denoted $\LocSys^{\on{Betti}}_\sG(X)$ and $\LocSys^\dr_\sG(X)$, respectively.
However, in the \'etale context, $\LocSys^{\on{restr}}_\sG(X)$ is the best algebro-geometric approximation
to the moduli of local systems that we can imagine. 

\end{rem}

\sssec{} \label{sss:formal affine Intro}

For a general sheaf theory, we prove the following theorem concerning the structure of $\LocSys^{\on{restr}}_\sG(X)$. 
Let $\LocSys^{\on{restr},\on{rigid}_x}_\sG(X)$ be the fiber product
$$\LocSys^{\on{restr}}_\sG(X)\underset{\on{pt}/\sG}\times \on{pt},$$
where 
$$\LocSys^{\on{restr}}_\sG(X)\to \on{pt}/\sG$$
is the map corresponding to taking the fiber at a chosen base point $x\in X$. So
$$\LocSys^{\on{restr}}_\sG(X)\simeq \LocSys^{\on{restr},\on{rigid}_x}_\sG(X)/\sG.$$

\medskip

We prove (in \thmref{t:main 1}) that $\LocSys^{\on{restr},\on{rigid}_x}_\sG(X)$ is a disjoint union of ind-affine ind-schemes $\CY$
(locally almost of finite type), each of which is a \emph{formal affine scheme}.

\medskip

We recall that a prestack $\CY$ is a formal affine scheme if it an be written as
a formal completion
$$\Spec(R)^\wedge_Y,$$
where $R$ is a connective $\sfe$-algebra (but not necessarily almost of finite type over $\sfe$) and
$Y\simeq \Spec(R')$ is a Zariski closed subset in $\Spec(R)$, where $R'$ is a (classical, reduced)
$\sfe$-algebra of finite type. 

\medskip

This all may sound technical, but the upshot is that the $\LocSys^{\on{restr}}_\sG(X)$ fails to be
an algebraic stack precisely to the same extent as in the Betti case, and the extent of this failure is such 
that we can control it very well. 

\medskip

To illustrate the latter point, in \secref{s:qcoh formal} we study the category
$\QCoh(\CY)$ on formal affine schemes (or quotients of these by groups) and show that its behavior
is very close to that of $\QCoh(-)$ on affine schemes (which is \emph{not at all} the case of $\QCoh(-)$ on
arbitrary ind-schemes). 

\sssec{}

As we have seen in \secref{sss:Intro Betti}, in the Betti context, the prestack $\LocSys^{\on{restr}}_\sG(X)$ 
splits into a disjoint union of prestacks $\CZ_\sigma$ parameterized by isomorphism classes of
semi-simple $\sG$-local systems\footnote{When $\sG$ is not reductive, the parameterization is by
the same set for the maximal reductive quotient of $\sG$.} $\sigma$ on $X$. Moreover, the underlying reduced prestack of each
$\CZ_\sigma$ is an algebraic stack.

\medskip

In \secref{s:uniformization} we prove that the same is true in any sheaf theory.  Furthermore, for each $\sigma$,
we construct a \emph{uniformization map} 
$$\underset{\sP}\sqcup\, \LocSys^{\on{restr}}_{\sP,\sigma_\sM}(X)\to \CZ_\sigma,$$
which is proper and surjective on geometric points, where: 

\begin{itemize}

\item The disjoint union runs over the set over parabolic subgroups $\sP$, such that $\sigma$
can be factored via an \emph{irreducible} local system $\sigma_\sM$ for some/any Levi splitting $\sP\hookleftarrow \sM$
(here $\sM$ is the Levi quotient of $\sP$);

\medskip

\item $\LocSys^{\on{restr}}_{\sP,\sigma_\sM}(X)$ is the \emph{algebraic stack}
$$\LocSys^{\on{restr}}_\sP(X)\underset{\LocSys^{\on{restr}}_\sM(X)}\times \on{pt}/\on{Aut}(\sigma_\sM).$$

\end{itemize} 

\sssec{} \label{sss:formal coarse Intro}

Let $\sG$ be again reductive. For a general sheaf theory, we do not have the picture involving \eqref{e:r coarse} that we had in the Betti case.
However, we do have a formal part of it.

\medskip

Namely, let $\CZ$ be a connected component of $\LocSys^{\on{restr}}_\sG(X)$. This is an ind-algebraic stack,
which can be written as
$$\underset{i}{\on{colim}}\, \CZ_i,$$
where each $\CZ_i$ is an algebraic stack isomorphic to the quotient of a (derived) affine scheme by $\sG$.

\medskip

We can consider the ind-affine ind-scheme
$$\CZ^{\on{coarse}}:=\underset{i}{\on{colim}}\, \Spec(\Gamma(\CZ_i,\CO_{\CZ_i})),$$
and the map
\begin{equation} \label{e:r Intro}
\brr:\CZ\to \CZ^{\on{coarse}}.
\end{equation} 

In \thmref{t:coarse restr} we prove that:

\smallskip

\noindent{(i)} $\CZ^{\on{coarse}}$ is a formal affine scheme (see \secref{sss:formal affine Intro} for what this means)
whose underlying reduced scheme is $\on{pt}$; 

\smallskip

\noindent{(ii)} The map \eqref{e:r Intro} makes $\CZ$ into a \emph{relative algebraic stack} over $\CZ^{\on{coarse}}$. 

%
%
%
%

\ssec{Overview: spectral decomposition}

Part II contains one of the two the main results of this paper, \thmref{t:action}. 

\sssec{}

We again start with a motivation in the Betti context. 

\medskip

Let $\CX$ be a connected space, and let $\bC$ be a DG category. 

\medskip

In this case, we have the notion of action of $\Rep(\sG)^{\otimes \CX}$ on $\bC$, see \cite[Sect. 1.7]{GKRV}. It consists of a compatible 
family of functors
$$\Rep(\sG)^{\otimes I}\to \End(\bC)\otimes (\Vect^\CX_\sfe)^{\otimes I}, \quad I\in \on{fSet},$$
where 
$\Vect^\CX_\sfe$ is the DG category $\on{Funct}(\CX,\Vect_\sfe)$ (it can be thought of as the category of
local systems of vector spaces on $\CX$, \emph{not necessarily} of finite rank), and $\on{fSet}$ is the category 
of finite sets. 

\medskip

Now, we have the stack of Betti local systems $\LocSys^{\on{Betti}}_\sG(\CX)$ and we can consider actions 
of the symmetric monoidal category $\QCoh(\LocSys^{\on{Betti}}_\sG(\CX))$ on $\bC$. 

\medskip

The tautological defined symmetric monoidal functor
$$\Rep(\sG)\otimes \QCoh(\LocSys^{\on{Betti}}_\sG(\CX))\to \Vect^\CX_\sfe$$ gives rise to a map (of $\infty$-groupoids)
\begin{equation} \label{e:action Betti Intro}
\{\text{Actions of $\QCoh(\LocSys^{\on{Betti}}_\sG(\CX))$ on $\bC$}\}\to \{\text{Actions of $\Rep(\sG)^{\otimes \CX}$ on $\bC$}\}.
\end{equation}

A relatively easy result (see \cite[Theorem 1.5.5]{GKRV}) says that the map \eqref{e:action Betti Intro} is an equivalence
(of $\infty$-groupoids). 

\sssec{}

We now transport ourselves to the context of algebraic geometry. Let $X$ be a connected scheme over $k$ and $\bC$ be a $\sfe$-linear
DG category. By an action of $\Rep(\sG)^{\otimes X\on{-lisse}}$ on $\bC$ we shall mean a compatible collection of functors
$$\Rep(\sG)^{\otimes I}\to \End(\bC)\otimes \qLisse(X)^{\otimes I}, \quad I\in \on{fSet}.$$

As before, we have the tautological symmetric monoidal functor
$$\Rep(\sG)\otimes \QCoh(\LocSys^{\on{restr}}_\sG(X))\to \qLisse(X),$$ 
and we obtain a map
\begin{equation} \label{e:action lisse Intro}
\{\text{Actions of $\QCoh(\LocSys^{\on{restr}}_\sG(\CX))$ on $\bC$}\}\to \{\text{Actions of $\Rep(\sG)^{\otimes X\on{-lisse}}$ on $\bC$}\}.
\end{equation}

One can can ask whether the map \eqref{e:action lisse Intro}
is an isomorphism as well, and our Spectral Decomposition theorem, namely,
\thmref{t:spectral} in the main body of the paper, says that it is. 

\sssec{}

Unfortunately, our proof of \thmref{t:spectral} is not aesthetically very satisfactory. In fact,
we conjecture that a more general statement along the same lines holds (see \conjref{c:Hom abs}),
when we replace the category $\qLisse(X)$ by what we call a \emph{gentle} Tannakian category $\bH$. 

\medskip

Our proof of \thmref{t:spectral} is very specific to $\bH$ being $\qLisse(X)$, where $X$ is a smooth
proper curve. 

\medskip

Namely, we use the fact that \eqref{e:action Betti Intro} is an equivalence to prove that the assertion of 
\thmref{t:spectral} holds in the Betti context (i.e., when $\qLisse(X)$ is the left
completion of the ind-completion of the category of \emph{finite-dimensional} Betti local systems on $X$). 

\medskip

Using Riemann-Hilbert, this formally implies the assertion of \thmref{t:spectral} holds in the de Rham
context (i.e., when $\qLisse(X)$ is the left
completion of the ind-completion of the category of de Rham local systems on $X$). 

\medskip 

Finally, we show that in the \'etale context, the assertion of \thmref{t:spectral} 
follows formally from its validity in the Betti context, essentially because the \'etale $\qLisse(X)$
(over any algebraically closed ground field) can be realized as a direct factor of the Betti version of $\qLisse(X')$ for some
complex curve $X'$. 

\begin{rem}
The particularly troublesome aspect of our proof of \thmref{t:spectral} is that it is not
applicable to the case when $X$ is a non-complete curve, while this case is of
interest if we have an eye on extending our theory to the ramified case.
\end{rem}

\ssec{Overview: the category $\Shv_\Nilp(\Bun_G)$}

In Part III of the paper, we take $G$ to be a reductive group and we will study the category 
$\Shv_\Nilp(\Bun_G)$ of sheaves on $\Bun_G$ (within any of our contexts) with singular
support in the nilpotent cone $\Nilp\subset T^*(\Bun_G)$.

\sssec{}

The stack $\Bun_G$ is non quasi-compact, and what allows us to work efficiently with the category $\Shv(\Bun_G)$
is the fact that we can simultaneously think of it as a \emph{limit}, taken over poset of quasi-compact open substacks $\CU\subset \Bun_G$,
$$\underset{\CU}{\on{lim}}\, \Shv(\CU),$$
with transition functors given by restriction, and \emph{also as a colimit}
$$\underset{\CU}{\on{colim}}\, \Shv(\CU),$$
with transition functors given by !-extension.

\medskip

We now take $\Shv_{\Nilp}(\Bun_G)$. More or less by definition, we still have 
$$\Shv_{\Nilp}(\Bun_G):=\underset{\CU}{\on{lim}}\, \Shv_\Nilp(\CU),$$
but we run into trouble with the colimit presentation:

\medskip

In order for such presentation to exist, we should be able to find a cofinal
family of quasi-compact opens, such that for every pair $\CU_1\overset{j}\hookrightarrow \CU_2$
from this family, the functor $j_!$ sends 
$$\Shv_\Nilp(\CU_1)\to \Shv_\Nilp(\CU_2).$$

\medskip

Fortunately, we can find such a family; its existence is guaranteed by \thmref{t:preserve Nilp Sing Supp prel}. 

\sssec{}

Thus, we can access the category $\Shv_{\Nilp}(\Bun_G)$ via the corresponding categories on 
quasi-compact open substacks. But our technical troubles are not over: 

\medskip

We do not know whether the categories $\Shv_{\Nilp}(\CU)$ are compactly generated. 
Such questions (for an arbitrary algebraic stack or even
scheme $\CY$, with a fixed $\CN\subset T^*(\CY)$) may be non-trivial. For example, it is \emph{not}
true in general that $\Shv_\CN(\CY)$ is generated by objects that are compact in $\Shv(\CY)$. 
We refer the reader to \secref{s:shvs on stacks} where some general facts pertaining to these
issues are summarized. 

\medskip

Although we conjecture that $\Shv_{\Nilp}(\Bun_G)$ is generated by objects
that are compact in the ambient category $\Shv(\Bun_G)$, we were not able to
prove this in full generality. We do, however, prove this in the de Rham and
Betti contexts. 

\medskip

That said, we were able to prove that $\Shv_{\Nilp}(\Bun_G)$ is compactly generated
as a DG category, and hence is dualizable. The latter is important for Part IV of the paper,
in order for the trace of the Frobenius endofunctor on $\Shv_{\Nilp}(\Bun_G)$ to be
well-defined. 

\sssec{}

We now proceed to formulating the other results in Part III. 

\medskip 

We consider the Hecke action on $\Shv(\Bun_G)$. Now, the subcategory 
$$\Shv_{\Nilp}(\Bun_G)\subset \Shv(\Bun_G)$$
has the following key feature with respect to this action:

\medskip

According to \cite{NY1}, combined with \cite[Theorem A.3.8]{GKRV}, the Hecke functors 
\begin{equation} \label{e:Hecke Intro}
\on{H}(-,-):\Rep(\cG)^{\otimes I}\otimes \Shv(\Bun_G)\to \Shv(\Bun_G\times X^I), \quad I\in\on{fSet},
\end{equation} 
send the subcategory
$$\Rep(\cG)^{\otimes I}\otimes \Shv_\Nilp(\Bun_G)\subset \Rep(\cG)^{\otimes I}\otimes \Shv(\Bun_G)$$
to
$$\Shv_\Nilp(\Bun_G)\otimes \qLisse(X)^{\otimes I}\subset \Shv(\Bun_G\times X^I).$$

\medskip

This means that $\Shv_\Nilp(\Bun_G)$ carries an action of $\Rep(\cG)^{\otimes X\on{-lisse}}$, i.e.,
we find ourselves in the setting of the Spectral Decomposition theorem. 

\medskip 

Thus, combined with \thmref{t:spectral} described above, we obtain the following assertion
(it appears as \thmref{t:spectral decomp} in the main body of the paper):

\begin{thm} \label{t:spectral decomp intro}
The category $\Shv_\Nilp(\Bun_G)$ has a natural structure of module category over
$\QCoh(\LocSys_\cG^{\on{restr}}(X))$. 
\end{thm} 

\thmref{t:spectral decomp intro} has an obvious ideological significance. For example, it immediately
implies that the category $\Shv_\Nilp(\Bun_G)$ splits as a direct sum
\begin{equation} \label{e:dir sum Intro}
\underset{\sigma}\oplus\, \Shv_\Nilp(\Bun_G)_\sigma,
\end{equation}
indexed by isomorphism classes of semi-simple $\cG$-local systems.

\medskip

However, in addition, we use \thmref{t:spectral decomp intro} extensively to prove a number of structural
results about $\Shv_\Nilp(\Bun_G)$. For example, we use it to prove: (i) the compact generation
of $\Shv_\Nilp(\Bun_G)$ (this is \thmref{t:Nilp comp gen abs});
(ii) the fact that in the de Rham context,  objects from $\Shv_\Nilp(\Bun_G)$ have regular singularities 
(this is \corref{c:RS all}); (iii) the tensor product property of $\Shv_\Nilp(\Bun_G)$ (\thmref{t:tensor product}, see below). 

\sssec{}

We now come to the second main result of this paper (it appears as \thmref{t:lisse} in the main body of the paper), 
which is in some sense a converse to the assertion of \cite{NY1} mentioned above: 

\medskip

\begin{thm} \label{t:lisse intro}
Let $\CF$ be an object of $\Shv(\Bun_G)$, such that the Hecke functors
\eqref{e:Hecke Intro} send it to
$$\Shv(\Bun_G)\otimes \qLisse(X),$$
then $\CF\in \Shv_\Nilp(\Bun_G)$. 
\end{thm} 

\medskip

A particular case of this assertion was conjectured by G.~Laumon. Namely, 
\cite[Conjecture 6.3.1]{laumon} says that Hecke eigensheaves have nilpotent singular support.

\sssec{}

The combination of Theorems \ref{t:spectral decomp intro} and \ref{t:lisse intro} allows us to establish a whole 
array of results about $\Shv_\Nilp(\Bun_G)$, in conjunction with another tool: Beilinson's spectral 
projector, whose definition we will now recall.

\medskip

Let us first start with a single $\cG$-local system $\sigma$. We can consider the category 
$$\on{Hecke}_\sigma(\Shv(\Bun_G))$$
of Hecke eigensheaves on $\Bun_G$ with respect to $\sigma$. 

\medskip 

We have a tautological forgetful functor
\begin{equation} \label{e:forget Hecke intro}
\oblv_{\on{Hecke}_\sigma}:\on{Hecke}_\sigma(\Shv(\Bun_G))\to \Shv(\Bun_G),
\end{equation} 
and Beilinson's spectral projector is a functor
$$\sP^{\on{enh}}_\sigma:\Shv(\Bun_G)\to \on{Hecke}_\sigma(\Shv(\Bun_G)),$$
left adjoint to \eqref{e:forget Hecke intro}. 

\medskip

A feature of the functor $\sP^{\on{enh}}_\sigma$ is that the composition
$$\Shv(\Bun_G)\overset{\sP^{\on{enh}}_\sigma}\longrightarrow \on{Hecke}_\sigma(\Shv(\Bun_G))
\overset{\oblv_{\on{Hecke}_\sigma}}\longrightarrow \Shv(\Bun_G)$$
is given by an explicit \emph{integral Hecke functor}\footnote{I.e., a colimit of functors \eqref{e:Hecke Intro}
for explicit objects of $\Rep(\cG)^{\otimes I}$, as $I$ ranges over the category of finite sets.}. 

\medskip

However, now that we have $\LocSys^{\on{restr}}_\cG(X)$, we can consider a version of the functor
$\sP^{\on{enh}}_\sigma$ is families: 

\sssec{}

Let $\CZ$ be a prestack over the field of coefficients $\sfe$, equipped with a map 
$$f:\CZ\to \LocSys^{\on{restr}}_\cG(X).$$

Then it again makes sense to consider the category of \emph{Hecke eigensheaves} parametrized by $S$: 
$$\on{Hecke}(\CZ,\Shv(\Bun_G)).$$

It is endowed with a forgetful functor
$$\oblv_{\on{Hecke},\CZ}:\on{Hecke}(\CZ,\Shv(\Bun_G))\to \QCoh(\CZ)\otimes \Shv(\Bun_G)$$
(i.e., forget the eigenproperty). 

\medskip

We have a version of Beilinson's spectral projector, which is now a functor, denoted in this paper by
$$\sP^{\on{enh}}_\CZ:\Shv(\Bun_G)\to \on{Hecke}(\CZ,\Shv(\Bun_G)),$$
left adjoint to the composition
$$\on{Hecke}(\CZ,\Shv(\Bun_G))\to \QCoh(\CZ)\otimes \Shv(\Bun_G)\to \Shv(\Bun_G).$$

\medskip

Let us note that the definition of functor $\sP^{\on{enh}}_\CZ$ only uses the existence of $\LocSys^{\on{restr}}_\cG(X)$. 
We do not need to use Theorems \ref{t:spectral decomp intro} and \ref{t:lisse intro} to prove its existence
or to establish its properties.

\sssec{}

However, let us now use the functor $\sP^{\on{enh}}_\CZ$ in conjunction with Theorems \ref{t:spectral decomp intro} and \ref{t:lisse intro}.
 
 \medskip

First, \thmref{t:lisse intro} implies that the inclusion
$$\on{Hecke}(\CZ,\Shv_\Nilp(\Bun_G)) \subset \on{Hecke}(\CZ,\Shv(\Bun_G))$$
is an equality.  

\medskip

And \thmref{t:spectral decomp intro} implies that the category 
$\on{Hecke}(\CZ,\Shv_\Nilp(\Bun_G))$ identifies with
$$\QCoh(\CZ)\underset{\QCoh(\LocSys^{\on{restr}}_\cG(X))}\otimes \Shv_\Nilp(\Bun_G).$$

\medskip

Thus, we obtain that the functor $\sP^{\on{enh}}_\CZ$ provides a left adjoint to the functor
\begin{multline*} 
\QCoh(\CZ)\underset{\QCoh(\LocSys^{\on{restr}}_\cG(X))}\otimes \Shv_\Nilp(\Bun_G)
\overset{f_*\otimes \on{Id}}\longrightarrow 
\QCoh(\CZ)\underset{\QCoh(\LocSys^{\on{restr}}_\cG(X))}\otimes \Shv_\Nilp(\Bun_G)\simeq \\
\simeq \Shv_\Nilp(\Bun_G)\hookrightarrow \Shv(\Bun_G),
\end{multline*} 
provided that $\CO_\CZ$ is compact. 

\medskip

This construction has a number of consequences:

\medskip

\noindent(i) It allows us to prove the compact generation of $\Shv_\Nilp(\Bun_G)$
(left adjoints can be used to construct compact generators); this is \thmref{t:Nilp comp gen abs}.

\medskip

\noindent(ii) We construct explicit generators of $\Shv_\Nilp(\Bun_G)$ by applying the
functor $\sP^{\on{enh}}_\CZ$ (for some particularly chosen $f:\CZ\to \LocSys^{\on{restr}}_\cG(X)$)
to $\delta$-function objects in $\Shv_\Nilp(\Bun_G)$. This leads to the theorem that all
objects in $\Shv_\Nilp(\Bun_G)$ have regular singularities (in the de Rham context); this is
Main Corollary \ref{c:RS all}. Combined with \corref{c:eigensheaves nilp}, we obtain that all Hecke eigensheaves have
regular singularities; this is Main Corollary \ref{c:RS}. 

\medskip

\noindent(iii) We use the above generators of $\Shv_\Nilp(\Bun_G)$ to prove the (unexpected,
but important for future applications) property that the tensor product functor
$$\Shv_\Nilp(\Bun_G)\otimes  \Shv_\Nilp(\Bun_G)\to \Shv_{\Nilp\times \Nilp}(\Bun_G\times \Bun_G)$$
is an equivalence; this is \thmref{t:tensor product}  
(see the discussion in \secref{sss:when tensor product} regarding why such an equivalence is not something
we should expect on general grounds). 

\begin{rem}

The properties of $\Shv_\Nilp(\Bun_G)$ mentioned above indicate that this category exhibits
behavior similar to that of $\Shv(\CY)$, where $\CY$ is an algebraic stack (equal to the 
union of open substacks) with a finite number of isomorphism classes of $k$-points 
(e.g., $N\backslash G/B$ or its affine counterparts), or to the category of character sheaves on $G$. 

\medskip

The analogy is in fact not too far-fetched, as for $X=\BP^1$, our $\Shv_\Nilp(\Bun_G)$ is all
of $\Shv(\Bun_G)$, and $\Bun_G$ is indeed an affine parabolic version of $N\backslash G/B$.

\end{rem} 

\ssec{Overview: Langlands theory}

Let $X$ be a curve over a ground field $k$, and we will work with any of the sheaf-theoretic contexts
from our list. 

\sssec{}

Having set up the theories of $\LocSys^{\on{restr}}_\cG(X)$ and $\Shv_\Nilp(\Bun_G)$, we are now in the position 
to state a version of the (categorical) Geometric Langlands Conjecture, with nilpotent singular support: this is
\conjref{c:restr GLC}. It says that we have an equivalence 
\begin{equation} \label{e:GLC Intro}
\Shv_\Nilp(\Bun_G)\simeq \IndCoh_\Nilp(\LocSys_\cG^{\on{restr}}(X)),
\end{equation}
as categories equipped with an action of $\QCoh(\LocSys_\cG^{\on{restr}}(X))$.

\medskip

Here in the right-hand side, $\IndCoh_{?}(-)$ stands for the category of ind-coherent sheaves with prescribed
\emph{coherent} singular support, a theory developed in \cite{AG}. (In {\it loc.cit.}, this theory was developed 
for quasi-smooth schemes/algebraic stacks, but in \secref{ss:coh sing supp} we show that it equally applicable
to objects such as our $\LocSys_\cG^{\on{restr}}$.) In our case $?=\Nilp$, the \emph{global nilpotent cone} 
in $\Sing(\LocSys_\cG^{\on{restr}}(X))$, see \secref{sss:Nilp spec}\footnote{It should not be confused with $\Nilp\subset T^*(\Bun_G)$:
the two uses of $\Nilp$ have different meanings, and occur on different sides of Langlands duality.}. 

\sssec{}

Note that \conjref{c:restr GLC} may be the first instance when a categorical statement is suggested
for automorphic sheaves in the context of $\ell$-adic sheaves.

\medskip

That said, both the de Rham and Betti contexts have their own forms of the (categorical) Geometric Langlands Conjecture.
In the de Rham context, this is an equivalence
\begin{equation} \label{e:GLC de Rham Intro}
\Dmod(\Bun_G)\simeq \IndCoh_\Nilp(\LocSys^\dr_\cG(X)),
\end{equation}
as categories equipped with an action of $\QCoh(\LocSys^\dr_\cG(X))$.

\medskip

In the Betti context, this is an equivalence 
\begin{equation} \label{e:GLC Betti Intro}
\Shv^{\on{all}}_\Nilp(\Bun_G)\simeq \IndCoh_\Nilp(\LocSys^{\on{Betti}}_\cG(X)),
\end{equation}
as categories equipped with an action of $\QCoh(\LocSys^{\on{Betti}}_\cG(X))$,
where $\Shv^{\on{all}}_{?}(-)$ stands for the category of all sheaves
(i.e., not necessarily ind-constructible ones) with a prescribed singular support. 

\medskip

We show that in each of these contexts, our \conjref{c:restr GLC} is a formal consequence 
of \eqref{e:GLC de Rham Intro} (resp., \eqref{e:GLC Betti Intro}), respectively. In fact, we show
that the two sides in \conjref{c:restr GLC} are obtained from the two sides in 
\eqref{e:GLC de Rham Intro} (resp., \eqref{e:GLC Betti Intro}) by
$$\QCoh(\LocSys^{\on{restr}}_\cG(X))\underset{\QCoh(\LocSys^?_\cG(X))}\otimes -$$
for ?= $\dr$ or Betti. 

\medskip

That said, we show (assuming Hypothesis \ref{h:Lde Rham GLC}) that the restricted version
of GLC (i.e., \eqref{e:GLC Intro}) actually implies the full de Rham version, i.e.,
\eqref{e:GLC de Rham Intro}. Probably, a similar argument can show that \eqref{e:GLC Intro}
implies the full Betti version (i.e., \eqref{e:GLC Betti Intro}) as well. 

\sssec{}

For the rest of this subsection we will work over the ground field $k=\ol\BF_q$, but assume that our geometric
objects (i.e., $X$ and $G$) are defined over $\BF_q$, so that they carry the geometric Frobenius 
endomorphism.

\medskip

We now come to the second main theme of this paper, the Trace Conjecture. 

\medskip

For any (quasi-compact) algebraic stack $\CY$ over $\ol\BF_q$, but defined over $\BF_q$. we can consider
the endomorphism (in fact, automorphism) of $\Shv(\CY)$ given by Frobenius pushforward,
$(\Frob_\CY)_*$. Since $\Shv(\CY)$ is a compactly generated (and, hence, dualizable)
category, we can consider the categorical trace of $(\Frob_\CY)_*$ on $\Shv(\CY)$:
$$\Tr((\Frob_\CY)_*,\Shv(\CY))\in \Vect_\sfe.$$

\medskip

The Grothendieck passage from Weil sheaves on $\CY$ to functions on $\CY(\BF_q)$ can
be upgraded to a map
$$\on{LT}:\Tr((\Frob_\CY)_*,\Shv(\CY))\to \on{Funct}(\CY(\BF_q)),$$
compatible with *-pullbacks and !-pushforwards, see \thmref{t:GLTF}.

\medskip

However, the map $\on{LT}$ is \emph{not at all} an isomorphism
(unless $\CY$ has finitely many isomorphism classes of $\ol\BF_q$-points). 

\sssec{} \label{sss:finally Trace}

We apply the above discussion to $\CY=\Bun_G$. Since $\Bun_G$ is not quasi-compact, the local term map is 
in this case a map
\begin{equation} \label{e:Trace map all Bun Intro}
\on{LT}:\Tr((\Frob_{\Bun_G})_*,\Shv(\Bun_G))\to \on{Funct}_c(\Bun_G(\BF_q)),
\end{equation}
where $\on{Funct}_c(-)$ stands for functions with finite support. 

\medskip

In what follows we will denote
$$\on{Autom}:=\on{Funct}_c(\Bun_G(\BF_q)).$$

This is the space of compactly supported unramified automorphic functions. 

\medskip

As we just mentioned, the map \eqref{e:Trace map all Bun Intro} is \emph{not} an isomorphism
(unless $X$ is of genus $0$).

\sssec{}

We now consider the full category 
$$\Shv_\Nilp(\Bun_G)\hookrightarrow \Shv(\Bun_G).$$

It is stable under the action of the Frobenius, and is dualizable as a DG category. Hence, it makes sense to consider
the object 
$$\Tr((\Frob_{\Bun_G})_*,\Shv_\Nilp(\Bun_G)) \in \Vect_\sfe.$$

Our Trace Conjecture (\conjref{c:Trace conj}) says that there exists a canonical isomorphism 
\begin{equation} \label{e:Trace Intro}
\Tr((\Frob_{\Bun_G})_*,\Shv_\Nilp(\Bun_G)) \simeq \on{Autom}.
\end{equation}

\begin{rem}

The sheaves-functions correspondence has been part of the geometric Langlands program since 
its inception by V.~Drinfeld: in his 1983 paper \cite{Dri}, he constructed a Hecke eigenfunction corresponding
to a 2-dimensional local system on $X$ by first constructing the corresponding sheaf and then
taking the associated functions.

\medskip

Constructions of this sort allow to produce particular elements in $\on{Autom}$ that satisfy some
desired properties. 

\medskip

Our Trace Conjecture is an improvement in that it, in principle, allows to deduce statements 
about the \emph{space} $\on{Autom}$ from statements of $\Shv_\Nilp(\Bun_G)$ as a \emph{category}. 

\end{rem}

\sssec{}

In fact, the Trace Conjecture is a particular case of a more general statement, \conjref{c:Trace conj legs},
which we refer to as the Shtuka Conjecture.

\medskip

Namely, for a finite set $I$ and $V\in \Rep(\cG)^{\otimes I}$ consider the Hecke functor
$$\on{H}(V,-): \Shv_\Nilp(\Bun_G)\to \Shv_\Nilp(\Bun_G)\otimes \IndLisse(X^I).$$

Generalizing the categorical trace construction, we can consider the trace of this functor,
precomposed with $(\Frob_{\Bun_G})_*$. The result will be an object that we denote
$$\wt{\on{Sht}}_{I,V}\in \IndLisse(X^I)\subset \Shv(X^I).$$

\medskip

Our Shtuka Conjecture says that we have a canonical isomorphism
\begin{equation} \label{e:shtuka conj Intro}
\wt{\on{Sht}}_{I,V} \simeq \on{Sht}_{I,V},
\end{equation}
where $\on{Sht}_{I,V}\in \Shv(X^I)$ is the shtuka cohomology, see \secref{sss:shtukas} where we recall
the definition.

\medskip

Note that the validity of \eqref{e:shtuka conj Intro} implies that the objects $\on{Sht}_{I,V}$ belong to $\qLisse(X^I)\subset \Shv(X^I)$. 

\medskip

The latter fact has been unconditionally established by C.~Xue in \cite{Xue2}, which provides a reality check for our 
Shtuka Conjecture. 

\sssec{}

We will now explain how the Trace Conjecture recovers V.~Lafforgue's spectral decomposition of $\on{Autom}$
along the arithmetic Langlands parameters.

\medskip

The ind-stack $\LocSys^{\on{restr}}_\cG(X)$ (which is an algebro-geometric object over $\sfe=\ol\BQ_\ell$) 
carries an action of Frobenius, by transport of structure; we denote it by $\Frob$. Denote
$$\LocSys^{\on{arithm}}_\cG:=\left(\LocSys^{\on{restr}}_\cG(X)\right)^\Frob.$$

A priori, $\LocSys^{\on{arithm}}_\cG(X)$ is also an \emph{ind}-algebraic stack, but we prove
(see \thmref{t:Frob-finite}) that $\LocSys^{\on{arithm}}_\cG(X)$ is an actual algebraic stack
(locally almost of finite type).  We also prove that it is quasi-compact
(i.e., even though $\LocSys^{\on{restr}}_\cG(X)$ had infinitely many connected components,
only finitely many of them are Frobenius-invariant). 

\medskip

The algebra
$$\Exc:=\Gamma(\LocSys^{\on{arithm}}_\cG(X),\CO_{\LocSys^{\on{arithm}}_\cG(X)})$$
receives a map from V.~Lafforgue's algebra of excursion operators; this map is surjective
at the level oh $H^0$, see \secref{sss:naive Exc}. 

\sssec{}

The categorical meaning of $\LocSys^{\on{arithm}}_\cG(X)$ is that the category
$\QCoh(\LocSys^{\on{arithm}}_\cG(X))$ identifies with the category of Hochschild chains
of $\Frob^*$ acting on $\QCoh(\LocSys^{\on{restr}}_\cG(X))$.

\medskip

We will now apply the relative version of the trace construction from \cite[Sect. 3.8]{GKRV}, and 
attach to the pair 
$$(\Shv_\Nilp(\Bun_G),(\Frob_{\Bun_G})_*),$$
viewed as acted on by the pair 
$$(\QCoh(\LocSys^{\on{restr}}_\cG(X)),\Frob^*),$$
its class 
$$\on{cl}(\Shv_\Nilp(\Bun_G),(\Frob_{\Bun_G})_*)\in 
\on{HH}_\bullet(\Frob^*,\QCoh(\LocSys^{\on{restr}}_\cG(X)))\simeq \QCoh(\LocSys^{\on{arithm}}_\cG(X)).$$

\medskip

We denote the resulting object of $\QCoh(\LocSys^{\on{arithm}}_\cG(X))$ by 
$$\Drinf\in \QCoh(\LocSys^{\on{arithm}}_\cG(X)).$$ 

Applying a version of \cite[Theorem 3.8.5]{GKRV}, we have
$$\Gamma(\LocSys^{\on{arithm}}_\cG(X),\Drinf)\simeq \Tr((\Frob_{\Bun_G})_*,\Shv_\Nilp(\Bun_G)).$$

Combining with the Trace Conjecture (see \eqref{e:Trace Intro}) we thus obtain an isomorphism
\begin{equation} \label{e:Drinf vs Autom}
\Gamma(\LocSys^{\on{arithm}}_\cG(X),\Drinf)\simeq \on{Autom}.
\end{equation} 

In particular, the tautological action of $\Exc$ on 
$\Gamma(\LocSys^{\on{arithm}}_\cG(X),\Drinf)$ gives rise to an action of $\Exc$ on $\on{Autom}$.
This recovers V.~Lafforgue's spectral decomposition.

\sssec{}

The ideological significance of the isomorphism \eqref{e:Drinf vs Autom} is that it 
provides a \emph{localization} picture for $\on{Autom}$.  

\medskip

Namely, it says that behind the vector space $\on{Autom}$ stands a finer object, namely, a quasi-coherent
sheaf (this is our $\Drinf$) on the moduli \emph{stack} of Langlands parameters (this is our $\LocSys^{\on{arithm}}_\cG(X)$),
such that $\on{Autom}$ is recovered as its global sections.

\medskip

In other words, $\on{Autom}$ is something that lives over the \emph{coarse moduli space}
$$\LocSys^{\on{arithm,coarse}}_\cG(X):=\Spec(\Exc),$$
and it is obtained as direct image along the tautological map 
$$\brr:\LocSys^{\on{arithm}}_\cG(X)\to \LocSys^{\on{arithm,coarse}}_\cG(X)$$
from a finer object, namely $\Drinf$, on the \emph{moduli stack}. 

\begin{rem}

The notation $\Drinf$ has the following origin: upon learning of V.~Lafforgue's work \cite{VLaf1}, 
V.~Drinfeld suggested that the objects
$$\wt{\on{Sht}}_{I,V}$$
mentioned above should organize themselves into an object of $\QCoh(\LocSys^{\on{arithm}}_\cG(X))$.
(However, at the time there was not yet a definition of $\LocSys^{\on{arithm}}_\cG(X)$.)

\medskip

Now, with our definition of $\LocSys^{\on{arithm}}_\cG(X)$, the Shtuka Conjecture, i.e., \eqref{e:shtuka conj Intro},
is precisely the statement
that the object $\Drinf$ constructed above realizes Drinfeld's vision. 

\end{rem}

\sssec{}

A particular incarnation of the localization phenomenon of $\on{Autom}$ is the following. 

\medskip

Fix an
$\sfe$-point of $\LocSys^{\on{arithm}}_\cG(X)$ corresponding to an \emph{irreducible}
Weil $\cG$-local system $\sigma$. In \thmref{t:irred Weil} we show that such a point
corresponds to a connected component of $\LocSys^{\on{arithm}}_\cG(X)$, isomorphic to
$\on{pt}/\on{Aut}(\sigma)$. 

\medskip

The restriction of $\Drinf$ to this connected component is then a representation of the
(finite) group $\on{Aut}(\sigma)$. The corresponding direct summand on $\on{Autom}$
is obtained by taking $\on{Aut}(\sigma)$-invariants in this representation. 

\sssec{}

Finally, let us juxtapose the Trace Conjecture with the Geometric Langlands Conjecture \eqref{e:GLC Intro}.
We obtain an isomorphism
$$\on{Autom}\simeq \Tr(\Frob^!,\IndCoh_\Nilp(\LocSys_\cG^{\on{restr}}(X))).$$

\medskip

Now, a (plausible, and much more elementary) \conjref{c:ignore N} says that the inclusion
$$\IndCoh_\Nilp(\LocSys_\cG^{\on{restr}}(X))\hookrightarrow \IndCoh(\LocSys_\cG^{\on{restr}}(X))$$
induces an isomorphism
\begin{equation} \label{e:ignore Nilp Intro}
\Tr(\Frob^!,\IndCoh_\Nilp(\LocSys_\cG^{\on{restr}}(X)))\simeq 
\Tr(\Frob^!,\IndCoh(\LocSys_\cG^{\on{restr}}(X))). 
\end{equation}

\medskip

Now, for any quasi-smooth stack $\CY$ with an endomorphism $\phi$, we have
$$\Tr(\phi^!,\IndCoh(\CY))\simeq \Gamma(\CY^\phi,\omega_{\CY^\phi}).$$

Hence, the right-hand side in \eqref{e:ignore Nilp Intro} identifies with 
$$\Gamma(\LocSys^{\on{arithm}}_\cG(X),\omega_{\LocSys^{\on{arithm}}_\cG(X)}).$$

Summarizing, we obtain that the combination of the above three conjectures yields an isomorphism
\begin{equation} \label{e:Autom Intro}
\on{Autom}\simeq \Gamma(\LocSys^{\on{arithm}}_\cG(X),\omega_{\LocSys^{\on{arithm}}_\cG(X)}).
\end{equation}

This gives a conjectural expression for the space of (unramified) automorphic functions purely in terms
of the stack of arithmetic Langlands parameters. 

\ssec{Notations and conventions}

The notations in this paper will largely follow those adopted in \cite{GKRV}. 

\sssec{Algebraic geometry}

There will be two algebraic geometries present in this paper. 

\medskip

On the one hand, we fix a ground field $k$ (assumed algebraically closed, but of arbitrary characteristic) and we 
will consider algebro-geometric objects over $k$. This algebraic geometry will occur on the \emph{geometric/automorphic} 
side of Langlands correspondence. 

\medskip

Thus, $X$ will be a scheme over $k$ (in Parts III and IV of the paper, $X$ will be a complete curve), $G$ will be a reductive 
group over $k$, $\Bun_G$ will be the stack of $G$-bundles on $X$, etc. 

\medskip

The algebro-geometric objects over $k$ will be \emph{classical}, i.e., \emph{non-derived}; this is because we will study
sheaf theories on them that are insensitive to the derived structure (such as $\ell$-adic sheaves, or D-modules). 

\medskip

All algebro-geometric objects over $k$ will be \emph{locally of finite type} (see \cite[Chapter 2, Sect. 1.6.1]{GR1} for what
this means). We will denote by $\on{Sch}_{\on{ft}/k}$ the category of schemes of finite type over $k$. 

\medskip

On the other hand, we will have a field of coefficients $\sfe$ (assumed algebraically closed \emph{and of characteristic zero}),
and we will consider \emph{derived} algebro-geometric objects over $\sfe$, see \secref{sss:DAG intro} below.

\medskip

The above two kinds of algebro-geometric objects do not generally mix unless we work with D-modules, 
in which case $k=\sfe$ is a field of characteristic zero.

\sssec{Higher categories}  \label{sss:higher categories}

This paper will substantially use the language of $\infty$-categories\footnote{We will often omit the adjective ``infinity" and refer to
$\infty$-categories simply as ``categories".}, as developed in \cite{Lu1}.

\medskip

We let $\Spc$ denote the $\infty$-category of spaces.

\medskip

Given an $\infty$-category $\bC$, and a pair of objects $\bc_1,\bc_2\in \bC$, we let $\Maps_\bC(\bc_1,\bc_2)\in \Spc$
the mapping space between them. 

\medskip

Recall that given an $\infty$-category $\bC$ that contains filtered colimits, an object $\bc\in \bC$ is said to be compact
if the Yoneda functor $\Maps_\bC(\bc,-):\bC\to \Spc$ preserves filtered colimits. We let $\bC^c\subset \bC$ denote
the full subcategory spanned by compact objects.

\medskip

Given a functor $F:\bC_1\to \bC_2$ between $\infty$-categories, we will denote by $F^R$ (resp., $F^L$) its right (resp., left) adjoint,
provided that it exists.

\sssec{Higher algebra}

Throughout this paper we will be concerned with \emph{higher algebra} over a field of coefficients, denoted $\sfe$
(as was mentioned above, throughout the paper, $\sfe$ will be assumed algebraically closed and of characteristic zero). 

\medskip

We will denote by $\Vect_\sfe$ the stable $\infty$-category of chain complexes of $\sfe$-modules,
see, e.g., \cite[Example 2.1.4.8]{GaLu}.

\medskip

We will regard $\Vect_\sfe$ as equipped with a symmetric monoidal structure (in the sense on $\infty$-categories), see, e.g., \cite[Sect. 3.1.4]{GaLu}.
Thus, we can talk about commutative/associative algebra objects in $\Vect_\sfe$, see, e.g., \cite[Sect. 3.1.3]{GaLu}.

%

\sssec{DG categories}

We will denote by $\DGCat$ the $\infty$-category of presentable stable $\infty$-categories, \emph{equipped with a module
structure over $\Vect_\sfe$ with respect to the symmetric monoidal structure on the $\infty$-category of presentable stable $\infty$-categories
given by the Lurie tensor product}, see \cite[Sect. 4.8.1]{Lu2}. We will refer to objects of $\DGCat$ as ``DG categories".
We emphasize that 1-morphisms in $\DGCat$ are in particular colimit-preserving.

\medskip

For a given DG category $\bC$, and a pair of objects $\bc_1,\bc_2\in \bC$, we have a well-defined ``inner Hom" object
$\CHom_\bC(\bc_1,\bc_2)\in \Vect_\sfe$, characterized by the requirement that
$$\Maps_{\Vect_\sfe}(V,\CHom_\bC(\bc_1,\bc_2))\simeq \Maps_\bC(V\otimes \bc_1,\bc_2), \quad V\in \Vect_\sfe.$$

\medskip

The category $\DGCat$ itself carries a symmetric monoidal structure, given by Lurie tensor product over $\Vect$:
$$\bC_1,\bC_2\rightsquigarrow \bC_1\otimes \bC_2.$$

\medskip

In particular, we can talk
about the $\infty$-category of associative/commutative algebras in $\DGCat$, which we denote by $\DGCat^{\on{Mon}}$ (resp., $\DGCat^{\on{SymMon}}$),
and refer to as monoidal (resp., symmetric monoidal) DG categories.

\medskip

Unless specified otherwise, all monoidal/symmetric monoidal DG categories will be assumed unital. Given a monoidal/symmetric monoidal DG category
$\CA$, we will denote by $\one_\CA$ its unit object.

\sssec{t-structures}

Given a DG category $\bC$, we can talk about a t-structures on it. 
For example, the category $\Vect_\sfe$ carries a natural t-structure. 

\medskip

Given a t-structure on $\bC$, we will denote by
$$\bC^{\leq n},\,\, \bC^{\geq n},\,\, \bC^\heartsuit$$
the corresponding subcategories (according to \emph{cohomological} indexing conventions), and also 
$$\bC^{<\infty}=\underset{n}\cup\,  \bC^{\leq n}, \,\, \bC^{>-\infty}=\underset{n}\cup\,  \bC^{\geq -n}.$$

We will refer to the objects of $\bC^{\leq 0}$ (resp.,  $\bC^{\geq 0}$) as \emph{connective}
(resp., \emph{coconnective}) with respect to the given t-structure.

%
%
%

\sssec{Derived algebraic geometry over $\sfe$}  \label{sss:DAG intro}

Most of the work in the present paper involves algebraic geometry on the spectral side of the Langlands correspondence. 
This is somewhat atypical to most work in geometric Langlands.

\medskip

As was mentioned above, algebraic geometry on the spectral side occurs over the field $\sfe$ and is derived. 
The starting point of derived algebraic geometry over $\sfe$ is the category $\affSch_{/\sfe}$ of \emph{derived affine schemes}
over $\sfe$, which is by definition the opposite category of the category of connective commutative algebras in $\Vect_\sfe$. 

\medskip

All other algebro-geometric objects over $\sfe$ will be \emph{prestacks}, i.e., accessible functors
$$(\affSch_{/\sfe})^{\on{op}}\to \Spc.$$

\medskip

Inside the category $\on{PreStk}_{/\sfe}$ of all prestacks, one singles out various subcategories. One such
subcategory is $\on{PreStk}_{\on{laft}/\sfe}$ that consists of prestacks \emph{locally almost of finite type}
(see \cite[Chapter 2, Sect. 1.7]{GR1}).  We set
$$\affSch_{\on{aft}/\sfe}:=\on{PreStk}_{\on{laft}/\sfe}\cap \affSch_{/\sfe}.$$

\medskip

We refer the reader to \cite[Chapter 3]{GR1} for the assighment
$$\CY\in \on{PreStk}_{/\sfe} \rightsquigarrow \QCoh(\CY)\in \DGCat.$$

\ssec{Acknowledgements}

We would like to express our gratitude to A.~Beilinson, D.~Ben-Zvi, D.~Beraldo, R.~Bezrukavnikov, J.~Campbell, 
V.~Drinfeld, M.~Kashiwara,
V.~Lafforgue, J.~Lurie, D.~Nadler, L.~Positselski, A.~Premet, T.~Saito, P.~Schapira and C.~Xue for illuminating discussions
and for helping us resolve numerous difficulties that we encountered.

\medskip

The entire project was supported by David Kazhdan's ERC grant No 669655
and BFS grant 2020189. The work of D.K. and Y.V. was supported BSF grant 2016363. 

\medskip

The work of D.A. was supported by NSF grant DMS-1903391. 
The work of D.G. was supported by NSF grant DMS-2005475. 
The work of D.K. was partially supported by ISF grant 1650/15. 
The work of S.R. was supported by NSF grant DMS-2101984.
The work of Y.V. was partially supported by ISF grants 822/17 and 2019/21.

\newpage

\centerline{\bf Part I: the (pre)stack $\LocSys^{\on{restr}}_{\sG}(X)$ and its properties}

\bigskip

Let us briefly describe the contents of this Part. 

\medskip

In \secref{s:restr} we define the prestack $\LocSys^{\on{restr}}_{\sG}(X)$ and state \thmref{t:main 1},
pertaining to its geometric properties. We study $\LocSys^{\on{restr}}_{\sG}(X)$ in the following 
general context: we consider prestacks of the form
$\bMaps(\sG),\bH)$, where $\bH$ is a \emph{gentle Tannakian category}. 
We recover $\LocSys^{\on{restr}}_{\sG}(X)$ by taking $\bH$ to be the category
$\qLisse(X)$ of lisse sheaves on $X$.

\medskip

In \secref{s:def} we establish the deformation theory properties of $\LocSys^{\on{restr}}_{\sG}(X)$
(and of its variant $\LocSys^{\on{restr},\on{rigid}_x}_{\sG}(X)$), 
leading to the conclusion that $\LocSys^{\on{restr},\on{rigid}_x}_{\sG}(X)$ is an ind-affine ind-scheme. 

\medskip

In \secref{s:uniformization} we finish the proof of \thmref{t:main 1} by combining the following
two results. One is \thmref{t:underlying reduced}, which says that the underlying \emph{reduced} prestack of 
$\LocSys^{\on{restr},\on{rigid}_x}_{\sG}(X)$ is a disjoint union of affine schemes. The other
is a general result due to J.~Lurie (we quote it as \thmref{t:formal}), which gives a deformation
theory criterion for an ind-scheme to be a formal scheme (for completeness, we supply
a proof in \secref{s:formal}). We prove \thmref{t:underlying reduced} by constructing a uniformization
of $\LocSys^{\on{restr}}_{\sG}(X)$ using parabolic subgroups of $\sG$
and \emph{irreducible} local systems for their Levi subgroups. In the process, we show that 
the set of connected components of $\LocSys^{\on{restr}}_{\sG}(X)$ is in bijection with the set
of isomorphism classes of \emph{semi-simple} $\sG$-local systems on $X$.

\medskip

In \secref{s:Betti and dR} we compare $\LocSys^{\on{restr}}_{\sG}(X)$ with the \emph{usual}
$\LocSys_{\sG}(X)$ in the two contexts in which the latter is defined, i.e., de Rham and Betti.
We show that in both cases, the resulting map $\LocSys^{\on{restr}}_{\sG}(X)\to \LocSys_{\sG}(X)$ 
identifies $\LocSys^{\on{restr}}_{\sG}(X)$ with the disjoint union of formal completions of
closed substacks, each corresponding to $\sG$-local systems with a fixed semi-simplification. 
Additionally, in the Betti context, we show that these closed substacks are exactly the fibers of
the map 
$$\brr:\LocSys_{\sG}(X)\to \LocSys^{\on{coarse}}_{\sG}(X),$$
where $\LocSys^{\on{coarse}}_{\sG}(X)$ is the corresponding coarse moduli space. 

\medskip

In \secref{s:geom properties}, we assume that $\sG$ is reductive. First,
we establish two more geometric properties of $\LocSys^{\on{restr}}_{\sG}(X)$: namely,
that it is \emph{mock-affine} and \emph{mock-proper}. We then state 
another structural result, \thmref{t:coarse restr}, which says that $\LocSys^{\on{restr}}_{\sG}(X)$
admits a coarse moduli space whose connected components are \emph{formal affine schemes}. 

\medskip

In \secref{s:formal coarse} we prove \thmref{t:coarse restr}. 

\medskip

In \secref{s:qcoh formal} we show that the category $\QCoh(-)$ of a formal affine scheme has
properties largely analogous to that of $\QCoh(-)$ of an affine scheme. The material from this
section will be applied when we will study the action of $\QCoh(\LocSys^{\on{restr}}_{\sG}(X))$
on $\Shv_{\Nilp}(\Bun_G)$.

\bigskip

\section{The restricted version of the stack of local systems}  \label{s:restr}

Let $X$ be a connected scheme over the ground field $k$. 

\medskip

We will begin this section by discussing several versions of the category of lisse
sheaves on $X$. 

\medskip

We will then proceed to the definition of our main object of study--the prestack
$\LocSys^{\on{restr}}_{\sG}(X)$, formulate a structural result about it, \thmref{t:main 1},
and consider a few examples.

\medskip

We will also introduce a rigidified version $\LocSys^{\on{restr},\on{rigid}_x}_{\sG}(X)$ of 
$\LocSys^{\on{restr}}_{\sG}(X)$, which gets rid of the ``stackiness".

\ssec{Sheaves}  

In this subsection we will define what we mean by the category of sheaves. 

\sssec{} \label{sss:Shv}

We will work in one of the following sheaf-theoretic contexts:
$$X\mapsto \Shv(X)^{\on{constr}}, \quad X\in \on{Sch}_{\on{ft}/k}.$$

\begin{itemize}

\item Constructible sheaves of $\sfe$-vector spaces on the topological space underlying $X$, when $k=\BC$
(to be referred to as the \emph{Betti context});

\smallskip

\item Holonomic D-modules on $X$, when $\on{char}(k)=0$
(to be referred to as the \emph{de Rham context});

\smallskip

\item Constructible $\ol\BQ_\ell$-adic \'etale sheaves on $X$, when $\on{char}(k)\neq \ell$
(to be referred to as the \emph{\'etale context}). 

\end{itemize}

In the Betti context, we will sometimes consider also the category of all sheaves of vector spaces on $X$;
we will denote it by $\Shv^{\on{all}}(X)$, see \secref{sss:all sheaves sch}.

\medskip

In the de Rham context we will sometimes consider also the category of all D-modules on $X$,
denoted $\Dmod(X)$. 

\sssec{}

We will denote by $\sfe$ the field of coefficients of our sheaves, which will always be algebraically closed and 
of characteristic $0$. In the three contexts above, this is
$\sfe$ (an arbitrary algebraically closed field of characteristic $0$), $k$ and $\ol\BQ_\ell$, respectively.

\sssec{}

The category $\Shv(X)^{\on{constr}}$ carries two symmetric monoidal structures. One is given by the ``usual" tensor product, denoted 
$\overset{*}\otimes$, for which the unit is the constant sheaf, denoted $\ul\sfe_X$. The other is given by the $\sotimes$ tensor product,
and its unit is the dualizing sheaf, denoted $\omega_X$. 

\medskip

By contrast, in the Betti context, the category $\Shv^{\on{all}}(X)$ only carries the $\overset{*}\otimes$ symmetric monoidal structure,
and in the de Rham context, the category $\Dmod(X)$ only carries the $\sotimes$ symmetric monoidal structure. 

\sssec{}

In any of the above three contexts, we set 
$$\Shv(X):=\on{Ind}(\Shv(X)^{\on{constr}}).$$

The two symmetric monoidal structures on $\Shv(X)^{\on{constr}}$ define by ind-extension the corresponding
two symmetric monoidal structures on $\Shv(X)$. 

\medskip

Unless explicitly stated otherwise, when talking about a symmetric monoidal structure on $\Shv(X)$, we will
be referring to the $\sotimes$ one. 

\sssec{}

The perverse t-structure on $\Shv(X)^{\on{constr}}$ uniquely extends to a t-structure on $\Shv(X)$
compatible with filtered colimits. Its heart is the ind-completion $\on{Ind}(\on{Perv}(X))$ of the category $\on{Perv}(X)$
of perverse sheaves on $X$.

\medskip

Since the above t-structure on $\Shv(X)$ is compactly generated (see \secref{sss:comp gen t} for what
this means), it is automatically right-complete.  

\medskip

We record the following result, proved in \secref{ss:shvs on schemes}:

\begin{thm} \label{t:Shv left-complete}
The category $\Shv(X)$ is left-complete in its t-structure.
\end{thm} 

\sssec{}  \label{sss:usual left-complete}

In addition to the perverse t-structure on $\Shv(X)^{\on{constr}}$, one can consider the \emph{usual} t-structure.

\medskip

It is a characterized by the requirement that the functors of *-fiber at closed points of $X$ are t-exact. 
By ind-extension, the usual t-structure on $\Shv(X)^{\on{constr}}$ defines a t-structure on 
$\Shv(X)$, which we will refer to as the ``usual" t-structure. 

\medskip

We note that the analog of \thmref{t:Shv left-complete} remains valid for the usual t-structure, due to the fact
that the two t-structures are a finite
distance apart (bounded by $\dim(X)$).

\medskip

That said, unless explicitly stated otherwise, we will work with the perverse t-structure. 

\ssec{Lisse sheaves}   \label{ss:lisse}

In this subsection we will introduce one of our main actors--the category of lisse sheaves on $X$. 

\sssec{} \label{sss:new lisse}

We define the full (abelian) subcategory 
$$\Lisse(X)^\heartsuit\subset \Shv(X)^{\on{constr}}$$
to consist of objects in the heart of the \emph{usual} t-structure that are dualizable
in the $\overset{*}\otimes$ symmetric monoidal structure. 

\medskip

We define the full DG subcategory
\begin{equation} \label{e:lisse subcat}
\Lisse(X)\subset \Shv(X)^{\on{constr}}
\end{equation}
to consist of objects whose cohomologies with respect to the \emph{usual} t-structure
belong to $\Lisse(X)^\heartsuit$.

\begin{rem}

Note that one can also characterize the subcategory 
$$\Lisse(X) \subset \Shv(X)$$
as the subcategory of objects dualizable with respect to the $\overset{*}\otimes$ symmetric monoidal structure on $\Shv(X)$.

\end{rem}

\sssec{Examples}

In the sheaf-theoretic contexts of \secref{sss:Shv}, the subcategory \eqref{e:lisse subcat}
can be characterized as follows:

\begin{itemize}

\item In the Betti context, $\Lisse(X)^\heartsuit$ is the abelian category of topological local systems on $X$ of finite rank;

\item In the \'etale context, $\Lisse(X)^\heartsuit$ consists of $\ell$-adic \'etale local systems on $X$; 

\item In the de Rham context (if $X$ is smooth of dimension $n$), $\Lisse(X)^\heartsuit[n]$ 
consists of $\CO$-coherent (right) D-modules;

\end{itemize}

\sssec{}

Suppose for a moment that $X$ is smooth of dimension $n$.  Set
$$\Lisse(X)^\heartsuit[n]=:\on{Perv}_{\on{lisse}}(X)\subset \on{Perv}(X).$$
be the full subcategory of \emph{lisse} objects. 

\medskip

In other words, in each of these sheaf-theoretic contexts of \secref{sss:Shv}, the condition that an object $\CF\in \on{Perv}(X)$ 
belong to $\on{Perv}_{\on{lisse}}(X)$ means that $\on{SingSupp}(\CF)=\{0\}\subset T^*(X)$, see \secref{sss:sing supp schemes perv} 
for the notations involving singular support. 

\medskip

It is easy to see that the subcategory \eqref{e:lisse subcat} can be alternatively 
characterized as consisting of objects whose cohomologies with respect to the \emph{perverse}
t-structure belong to $\on{Perv}_{\on{lisse}}(X)$. 

\medskip

Thus, in the notations of \secref{sss:sing supp schemes constr},
$$\Lisse(X)=\Shv_{\{0\}}(X)^{\on{constr}}.$$

\sssec{}

We define the full (abelian) subcategory 
$$\qLisse(X)^\heartsuit \subset \Shv(X)$$
to consist of objects in the heart of the \emph{usual} t-structure that could
be written as filtered colimits of objects from $\Lisse(X)^\heartsuit$.

\medskip

The following definition is central for this paper:

\begin{defn} 
Let 
\begin{equation} \label{e:q lisse subcat}
\qLisse(X)\subset \Shv(X)
\end{equation} 
be the full DG subcategory consisting of objects whose cohomologies with respect to the \emph{usual} t-structure
belong to $\qLisse(X)^\heartsuit$.
\end{defn} 

\sssec{}

Each of the above categories: 
$$\Lisse(X)^\heartsuit,\,\, \Lisse(X),\,\, \qLisse(X)^\heartsuit,\,\, \qLisse(X)$$
acquires a symmetric monoidal structures, induced by
the $\overset{*}\otimes$ symmetric monoidal structure on $\Shv(X)$. 

\medskip

The categories $\Lisse(X)$ and $\qLisse(X)$ carry t-structures, inherited from the \emph{usual} t-structure
on $\Shv(X)$, and their hearts identify with $\Lisse(X)^\heartsuit$ and $\qLisse(X)^\heartsuit$, respectively. 

\sssec{} \label{sss:fiber functor x}

Given a point $x\in X$, consider the (symmetric monoidal) functor
\begin{equation} \label{e:fiber functor}
\qLisse(X) \overset{\on{ev}_x}\longrightarrow \Vect_\sfe,
\end{equation}
given by taking the *-fiber at $x$.

\medskip

Note that the functor $\on{ev}_x$ is t-exact and conservative. We will regard \eqref{e:fiber functor} 
as a fiber functor on $\qLisse(X)$, viewed as a symmetric monoidal category. 

\sssec{}

Assume for a moment again that $X$ is smooth. Let $\on{Ind}(\on{Perv}_{\on{lisse}}(X))$ be the ind-completion 
of $\on{Perv}_{\on{lisse}}(X)$, viewed as a full abelian subcategory in $\on{Ind}(\on{Perv}(X))$.

\medskip

It is easy to see that the subcategory \eqref{e:q lisse subcat} can be alternatively 
characterized as consisting of objects whose cohomologies with respect to the \emph{perverse}
t-structure belong to $$\on{Ind}(\on{Perv}_{\on{lisse}}(X))\subset \on{Ind}(\on{Perv}(X)).$$

\medskip

In other words,
$$\qLisse(X)=\Shv_{\{0\}}(X)$$
in the notations of \secref{sss:sing supp schemes}. 

\sssec{}

We can define a \emph{different} embedding
\begin{equation} \label{e:!-emb}
\qLisse(X)\hookrightarrow \Shv(X), \quad E\mapsto E\overset{*}\otimes \omega_X.
\end{equation}

This embedding is a symmetric monoidal functor, when we regard $\Shv(X)$ as a symmetric 
monoidal category via the $\sotimes$ operation. 

\medskip

That said, when $X$ is smooth, the above two embeddings $\qLisse(X)\rightrightarrows \Shv(X)$ 
differ by a cohomological shift (by $\dim(X)$), and in particular, they have the same
essential image. 

\sssec{}

By \thmref{t:Shv left-complete}, the category $\qLisse(X)$ is left-complete in its t-structure. 

\medskip

Unfortunately, for a general $X$ we will be able to say very little about
the general categorical properties of $\qLisse(X)$. For example, we do not know whether
it is compactly generated or even dualizable.

\medskip

That said, our main application is when $X$ is a smooth algebraic curve, in which case we do
know that $\qLisse(X)$ is compactly generated, see Sects. \ref{sss:analyze P1}-\ref{sss:curves}. 

\ssec{Another version of lisse sheaves} \label{ss:iLisse}

In addition to $\qLisse(X)$ we can consider its variant, denoted $\iLisse(X)$, introduced below. 
The main advantage of $\iLisse(X)$ is that it is compactly generated, by definition.

\sssec{} \label{sss:ind-lisse}

Set
$$\iLisse(X):=\on{Ind}(\Lisse(X)).$$

In other words, $\iLisse(X)$ is the full subcategory in $\Shv(X):=\on{Ind}(\Shv(X)^{\on{constr}})$ generated by $\Lisse(X)$. 

\medskip

The t-structure on $\Lisse(X)$ uniquely extends to a t-structure on $\iLisse(X)$ compatible with filtered colimits. 

\sssec{} \label{sss:qLisse left compl iLisse}

We have a tautologically defined fully faithful functor
\begin{equation} \label{e:Ind lisse}
\iLisse(X)\to \qLisse(X).
\end{equation} 

\medskip

The functor \eqref{e:Ind lisse} sends compact generators of $\iLisse(X)$ to compact objects 
of $\qLisse(X)$. This implies that \eqref{e:Ind lisse} is fully faithful. 

\medskip

The functor \eqref{e:Ind lisse} is t-exact since the t-structure on $\qLisse(X)$ is also 
compatible with filtered colimits. Moreover, it is easy to see that the functor \eqref{e:Ind lisse} 
induces an equivalence on the hearts. Hence, it induces an equivalence $$(\iLisse(X))^{\geq -n}\to (\qLisse(X))^{\geq -n}$$
for any $n$. From here it follows that the functor \eqref{e:Ind lisse} identifies $\qLisse(X)$ with
the left completion of $ \iLisse(X)$.

\sssec{}

Note, however, the functor \eqref{e:Ind lisse} is \emph{not} always an equivalence. For example, it fails to be such for $X=\BP^1$,
see \secref{sss:analyze P1}. 

\medskip

Equivalently, the category $\iLisse(X)$ is \emph{not} 
necessarily left-complete in its t-structure.

\medskip

That said, as we will see in \secref{sss:curves}, the functor \eqref{e:Ind lisse} is an equivalence
for all smooth connected curves $X$ (projective or affine) different from $\BP^1$. 

\begin{rem}
The procedure by which we obtained $\iLisse(X)$ from $\qLisse(X)$ is similar to the procedure by which one 
produces $\IndCoh(S)$ from $\QCoh(S)$ (where $S$ is a scheme almost of finite type). 

\medskip

In that situation we also have a functor
\begin{equation} \label{e:IndCoh to QC}
\IndCoh(S)\to \QCoh(S). 
\end{equation} 

However, unlike \eqref{e:Ind lisse}, the functor \eqref{e:IndCoh to QC} is only fully faithful when $S$ is smooth, in which case
it is an equivalence. 

\medskip

If $S$ is not smooth but eventually coconnective, the functor \eqref{e:IndCoh to QC} is a co-localization
(i.e., admits a fully faithful left adjoint). So, in a sense, the functor \eqref{e:IndCoh to QC} exhibits a behavior
opposite to that of \eqref{e:Ind lisse}. 

\end{rem} 

\sssec{}

One should consider $ \iLisse(X)$ as a ``really nice" symmetric monoidal category,
in that it is compactly generated and rigid (see \cite[Chapter 1, Sect. 9.2]{GR1} for what this means). 

\medskip

Moreover, one can pick compact generators that belong to the heart of $ \iLisse(X)$,
and they will have cohomological dimension bounded by $\dim(X)$. 

\medskip

One thing that $ \iLisse(X)$ is \emph{not} is that it is \emph{not} the derived category
of its heart, see \secref{sss:Kpi1}. 

\ssec{Definition of $\LocSys^{\on{restr}}_{\sG}(X)$ as a functor}  \label{ss:restr}

For the duration of Part I, we let $\sG$ be a connected algebraic group. 

\sssec{} \label{sss:ten prod t}

Recall that if $\bC$ and $\bC'$ are DG categories, each equipped with a 
t-structure\footnote{Recall that according to our conventions, we require that t-structures be 
compatible with filtered colimits.}, then the tensor product $$\bC\otimes \bC'$$
carries a naturally defined t-structure, characterized by the property that its connective part 
$(\bC\otimes \bC')^{\leq 0}$ is generated by objects 
$$\bc\otimes \bc', \quad \bc\in \bC^{\leq 0},\,\, \bc'\in (\bC')^{\leq 0}.$$

Suppose for a moment that $\bC'$ is of the form $R\mod$, where $R\in \on{AssocAlg}(\Vect_\sfe^{\leq 0})$. 
Then the above t-structure is characterized by the property that the forgetful functor
$$\bC\otimes R\mod\simeq R\mod(\bC)\to \bC$$
is t-exact. 

\sssec{}

Let $X$ be a connected scheme of finite type over $k$. We define the prestack $$\LocSys^{\on{restr}}_\sG(X)$$ 
over the field $\sfe$ of coefficients, by sending 
$S\in \affSch_{/\sfe}$ to the space of \emph{right t-exact symmetric monoidal functors}
\begin{equation} \label{e:sym mon}
\Rep(\sG) \to \QCoh(S)\otimes  \qLisse(X).
\end{equation} 

Note that $\sfe$-points of $\LocSys^{\on{restr}}_\sG(X)$ are what we usually call $\sG$-local
systems on $X$ (within our sheaf theory). 

\sssec{Example}\label{sss:Tannaka BG}
Let $X = \on{pt}$.  Then $\qLisse(X) = \on{Vect}_\sfe$ and we have, by Tannaka duality (see e.g. \cite[Corollary 9.4.4.7]{Lu3}), 
that $\LocSys^{\on{restr}}(X) = \on{pt}/\sG$, the classifying stack of $\sG$.

\sssec{}

The main result of this Part is the following:

\begin{mainthm} \label{t:main 1} 
The prestack $\LocSys^{\on{restr}}_\sG(X)$ is an \'etale stack, equal to 
a disjoint union of \'etale stacks each of which can be written as an 
\'etale-sheafified quotient by $\sG$ of an \'etale stack $\CY$ with the following properties:

\medskip

\noindent{\em(a)} $\CY$ is locally almost of finite type\footnote{See \cite[Chapter 2, Sect. 1.7.2]{GR1} for what this means.};

\medskip

\noindent{\em(b)} $^{\on{red}}\CY$ is an affine (classical, reduced) scheme; 

\medskip

\noindent{\em(c)} $\CY$ is an ind-scheme;

\medskip

\noindent{\em(d)} $\CY$ can be written as
\begin{equation} \label{e:presentation A}
\underset{n\geq 0}{\on{colim}}\, \Spec(R_n)
\end{equation} 
where $R_n$ are connective commutative $\sfe$-algebras of the following form:
there exists a connective commutative $\sfe$-algebra $R$ and 
elements\footnote{By an element of $R$ we mean a point in the space corresponding to $R$ via Dold-Kan correspondence.} 
$f_1,...,f_m\in R$ so that
$$R_n=R\underset{\sfe[t_1,...,t_m]}\otimes \sfe[t_1,...,t_m]/(t_1^n,...,t_m^n), \quad t_i\mapsto f_i\in R,\,\, i=1,....,m.$$
\end{mainthm} 

%

\begin{rem}  \label{r:ind}

Points (a,b,c) of \thmref{t:main 1} can be combined to the following statement: $\CY$ can be written as \emph{filtered}
colimit 
\begin{equation} \label{e:ind prel}
\CY \simeq \underset{i}{\on{colim}}\, Y_n,
\end{equation} 
where all $Y_n$ are affine schemes almost of finite type\footnote{See \cite[Chapter 2, Sect. 1.7.1]{GR1} for what this means.}, 
and the maps $Y_{n_1}\to Y_{n_2}$ are closed nil-isomorphisms 
(i.e., closed embeddings that induce isomorphisms of the underlying reduced prestacks), see \cite[Chapter 2,  Corollary 1.8.6(a)]{GR2}. 

\end{rem}

\begin{rem} \label{r:formal affine}

Note, however, that points (a,b,c) of \thmref{t:main 1} do \emph{not} include the assertion contained in (d). For example,
if we take 
$$\CY:=\underset{n}{\on{colim}}\, (\BA^n)^\wedge_0,$$
then this $\CY$ admits a presentation as in \eqref{e:ind prel} (with the specified properties), but it does \emph{not} admit
a presentation \eqref{e:presentation A} (the reason is that prestacks of the latter form admit \emph{cotangent spaces},
while the former only \emph{pro-cotangent spaces}, see \secref{ss:deform}). 

\medskip

The property of admitting a presentation as in \eqref{e:presentation A} insures, among other things, that the category 
$\QCoh(\CY)$ is particularly well-behaved 
(has many properties similar to those of $\QCoh(-)$ of an affine scheme, see \secref{s:qcoh formal}). 

\medskip

Prestacks $\CY$ satisfying (d) are called \emph{formal affine schemes}. 

\medskip

Finally, we emphasize that the commutative algebra $R$ that appears in (d) is \emph{not} necessarily 
almost of finite type over $\sfe$. However, according to \corref{c:R is Noeth}, we can choose $R$ so that
it is Noetherian. 

\end{rem}

%
%
%
%
%
%

\begin{rem}
For the validity of \thmref{t:main 1} in the Betti context, we can work more generally: instead of starting with an algebraic variety $X$
over $\BC$, we can let $X$ be a topological space homotopy equivalent to a (retract of a) finite CW complex.
\end{rem} 

\ssec{Some examples}

\sssec{}  \label{sss:Gm}

Let $\sG=\BG_m$. As we shall see in \corref{c:irred}, in this case the \emph{underlying reduced prestack} of $\LocSys^{\on{restr}}_\sG(X)$
is a disjoint union, over the set of isomorphism classes of one-dimensional local systems on $X$, of copies of $\on{pt}/\BG_m$. 

\medskip

In \secref{ss:deform} we will see that for each 1-dimensional local system (i.e., an $\sfe$-point of $\LocSys^{\on{restr}}_{\BG_m}(X)$),
the tangent space to $\LocSys^{\on{restr}}_\sG(X)$ at this point identifies with 
\begin{equation} \label{e:cohomology}
\on{C}^\cdot(X,\sfe_X)[1]\in \Vect_\sfe,
\end{equation}
i.e., it looks like the tangent space of the ``usual would-be" $\LocSys_{\BG_m}(X)$. 
(Tangent spaces are defined for prestacks that admit deformation theory and are locally almost of finite type, see 
\cite[Chapter 1, Sect. 4.4]{GR2}.) 

\sssec{} \label{sss:Ga}

Let $\sG=\BG_a$. We claim that in this case $\LocSys^{\on{restr}}_\sG(X)$ is the algebraic stack associated with the 
object \eqref{e:cohomology}, i.e.,
\begin{equation} \label{e:cohomology 1}
\Maps(\Spec(R),\LocSys^{\on{restr}}_{\BG_a}(X))=\tau^{\leq 0}(R\otimes \on{C}^\cdot(X,\sfe_X)[1]),
\end{equation}
where we view an object of $(\Vect_\sfe)^{\leq 0}$ as a space by the Dold-Kan functor (see \cite[Chapter 1, Sect. 10.2.3]{GR1}). 

\medskip

Indeed, 
$$\Rep(\BG_a)\simeq \sfe[\xi]\mod, \quad \deg(\xi)=1,$$  
so the space of symmetric monoidal functors from $\Rep(\BG_a)$ to any symmetric monoidal category $\bA$ identifies with
$$\tau^{\leq 0}(\CEnd_\bA(\one_\bA)[1]).$$

In our case $\bA=R\mod\otimes \qLisse(X)$ so $\one_\bA=R\otimes \sfe_X$, whence \eqref{e:cohomology 1}.  

\sssec{Notation}

In what follows, for $V\in \Vect_\sfe$ we will use the notation $\on{Tot}(V)$ for the corresponding prestack, i.e.,
\begin{equation} \label{e:Tot}
\Hom(S,\on{Tot}(V))=\tau^{\leq 0}(V\otimes \Gamma(S,\CO_S)), \quad S\in \affSch_{/\sfe}. 
\end{equation} 

For example, when $V\in \Vect_\sfe^\heartsuit\cap \Vect_\sfe^c$, we have
$$\on{Tot}(V)=\Spec(\Sym(V^\vee)).$$

Thus,  \eqref{e:cohomology 1} is saying that 
$$\LocSys^{\on{restr}}_{\BG_a}(X)\simeq \on{Tot}(\on{C}^\cdot(X,\sfe_X)[1]).$$

%
%
%

\sssec{}

This is a preview of the material in \secref{s:Betti and dR}:

\medskip

Let our sheaf-theoretic context be either Betti or de Rham, so in both cases we have the usual algebraic stack
$\LocSys_\sG(X)$. In this case we will show that there exists a forgetful map
\begin{equation} \label{e:restr to all}
\LocSys^{\on{restr}}_\sG(X)\to \LocSys_\sG(X),
\end{equation} 
which identifies $\LocSys^{\on{restr}}_\sG(X)$ with the disjoint union of formal completions of a collection
of pairwise non-intersecting Zariski-closed reduced substacks of $\LocSys_\sG(X)$, such that every $\sfe$-point
of $\LocSys_\sG(X)$ belongs to (exactly) one of these substacks. 
Furthermore, we will be able to describe the corresponding reduced substacks explicitly.  

\ssec{Rigidification}

Let us choose a base point $x\in X$. We will introduce a cousin of 
$\LocSys^{\on{restr}}_\sG(X)$, denoted $\LocSys^{\on{restr},\on{rigid}_x}_\sG(X)$, that
has to do with choosing a trivialization of our local systems at $x$. 

%

\sssec{} \label{sss:fiber functor x again}

Given a base point $x\in X$, consider the functor 
\begin{equation} \label{e:fiber functor again}
\qLisse(X) \overset{\on{ev}_x}\longrightarrow \Vect_\sfe,
\end{equation}
of \eqref{e:fiber functor}. 

\medskip

Consider the prestack $\LocSys^{\on{restr},\on{rigid}_x}_\sG(X)$ that sends $S\in \affSch_{/\sfe}$ to the space of symmetric
monoidal functors \eqref{e:sym mon}, equipped with an isomorphism between the composition
$$\Rep(\sG) \to \QCoh(S)\otimes  \qLisse(X) \overset{\on{Id}\otimes \on{ev}_x}\longrightarrow \QCoh(S)$$
and
\begin{equation} \label{e:forget with S}
\Rep(\sG) \overset{\oblv_\sG}\longrightarrow \Vect_\sfe \overset{\on{unit}}\longrightarrow \QCoh(S)
\end{equation}
(as symmetric monoidal functors). Note that the latter identification implies that the functor \eqref{e:sym mon} is right t-exact:
this is due to the fact that $\on{ev}_x$ is t-exact and conservative. 

\medskip

In other words,
\begin{equation} \label{e:rigid as fiber prod}
\LocSys^{\on{restr},\on{rigid}_x}_\sG(X)\simeq \LocSys^{\on{restr}}_\sG(X)\underset{\on{pt}/\sG}\times \on{pt},
\end{equation}
where the map
$$\LocSys^{\on{restr}}_\sG(X)\to \on{pt}/\sG$$ 
is given by \eqref{e:fiber functor} (see Example \ref{sss:Tannaka BG}). 

\sssec{}

From the above description of $\LocSys^{\on{restr},\on{rigid}_x}_\sG(X)$ as a fiber product, we obtain that there is 
a natural action of $\sG$ on $\LocSys^{\on{restr},\on{rigid}_x}_\sG(X)$, and we will show in \corref{c:etale descent} that
$\LocSys^{\on{restr}}_\sG(X)$ identifies with the \'etale sheafification of the quotient of $\LocSys^{\on{restr},\on{rigid}_x}_\sG(X)$ by
this action.  

\medskip

Hence, in order to prove \thmref{t:main 1}, it will suffice to show the following:

\begin{thm} \label{t:main 1 bis}
The prestack $\LocSys^{\on{restr},\on{rigid}_x}_\sG(X)$ is an \'etale stack equal 
to a disjoint union of \'etale stacks $\CY$ with properties \em{(a)-(d)} listed in  \thmref{t:main 1}. 
\end{thm}

\ssec{Gentle Tannakian categories} \label{ss:gentle}

For the proof of \thmref{t:main 1} it will be convenient to replace $\qLisse(X)$ by an abstract
symmetric monoidal category, to be denoted $\bH$, that possesses certain properties.

\medskip

In this subsection we will introduce the notion of \emph{gentle Tannakian category}.
This will be the class of symmetric monoidal categories for which will state and prove
an appropriate generalization of \thmref{t:main 1}. 

\sssec{}  \label{sss:pre-Tannakian}

Let $\bH$ be a symmetric monoidal category, equipped with a t-structure and a 
conservative t-exact symmetric monoidal functor\footnote{Recall that all our functors are
by default assumed to commute with all colimits.} $\oblv_\bH$ to $\Vect_\sfe$.

\medskip

Note that the assumptions on $\oblv_\bH$ imply that the t-structure on $\bH$
is compatible with filtered colimits and right-complete. 

\medskip

We will additionally assume that $\bH$ is left-complete in its t-structure.  

\sssec{}  \label{sss:conditions}

Let $\bH$ be as above. We will make the following assumptions, which can be summarized in saying that $\bH$ is a 
particularly well-behaved Tannakian category: 

\begin{itemize}

\item The following classes of objects in $\bH^\heartsuit$ coincide: 

\smallskip

(i) Objects contained in $\bH^c\cap \bH^\heartsuit$;

\smallskip

(ii) Objects that are sent to compact objects in $\Vect_\sfe$ by $\oblv_\bH$;

\smallskip

(iii) Dualizable objects.

\medskip

\item The object $\one_\bH$ has the following properties:

\smallskip

(i) The functor $\CHom_\bH(\one_\bH,-)$ has a finite cohomological amplitude;

\smallskip

(ii) For any $\bh\in \bH^c\cap \bH^\heartsuit$, the cohomologies of $\CHom_\bH(\one_\bH,\bh)\in \Vect_\sfe$
are finite-dimensional.

\medskip

\item The category $\bH^\heartsuit$ is generated under colimits by 
$\bH^c\cap \bH^\heartsuit$.

\end{itemize}

\medskip

We will call a symmetric mononidal category with the above properties a \emph{gentle Tannakian category}. 

\sssec{}

Note that the above conditions imply that there exists a bound $n$, such that 
for any pair of objects $\bh_1,\bh_2\in \bH^c\cap \bH^\heartsuit$, the object 
$\CHom_\bH(\bh_1,\bh_2)$ lives in cohomological degrees $\leq n$, and each
cohomology is finite-dimensional. 

\sssec{Example}

The main is example of interest for us is of course
$$\bH=\qLisse(X).$$

In this case, we take the fiber functor $\oblv_\bH$ to be $\on{ev}_x$, see
\secref{sss:fiber functor x}.

\medskip

Another example to keep in mind is $\bH=\Rep(\sH)$, where $\sH$ is an algebraic group of finite type
over $\sfe$, and $\oblv_\bH$ is the tautological forgetful functor (to be henceforth denoted $\oblv_\sH$). 

\sssec{}

Let $\bH^{\on{access}}$ be the full subcategory of $\bH$ generated by $\bH^c\cap \bH^\heartsuit$.
In other words, $\bH^{\on{access}}$ is the ind-completion of the small DG subcategory $\bH^{\on{access},c}\subset \bH$
consisting of \emph{cohomologically bounded} objects all of whose cohomologies belong to 
$\bH^c\cap \bH^\heartsuit$.

\medskip

Since $\bH^c\cap \bH^\heartsuit$ is closed under the monoidal operation, the category $\bH^{\on{access}}$
inherits a symmetric monoidal structure. 

\medskip

By construction, $\bH^{\on{access}}$ is \emph{rigid} as a symmetric monoidal category
(see \cite[Chapter 1, Sect. 9.2]{GR1} for what this means)\footnote{Note that $\bH$ was not
necessarily rigid.}. 

\sssec{}

Consider the tautological embedding
\begin{equation} \label{e:un-ren H}
\bH^{\on{access}}\to \bH.
\end{equation} 

\medskip

The t-structure on $\bH$ restricts to a t-structure on $\bH^{\on{access},c}$, and the latter gives
rise to a t-structure $\bH^{\on{access}}$. The functor \eqref{e:un-ren H} is t-exact and induces an equivalence
\begin{equation} \label{e:ren coconn equiv}
(\bH^{\on{access}})^{\geq -n}\to \bH^{\geq -n}.
\end{equation} 

It follows that the functor \eqref{e:un-ren H} realizes $\bH$ as the left-completion of $\bH^{\on{access}}$.

\medskip

Since the functor \eqref{e:un-ren H} sends the compact generators of $\bH^{\on{access}}$ to compact
objects of $\bH$, we obtain that \eqref{e:un-ren H} is fully faithful. 

\begin{rem}
The mechanism by which \eqref{e:un-ren H} fails to be an equivalence is that the subcategory 
$\bH^c\cap \bH^\heartsuit$ does not necessarily generate $\bH$ under colimits.

\medskip

This may happen even if $\bH$ itself compactly generated: its compact generators may be
unbounded below. 

\end{rem} 

\sssec{Example}

When $\bH=\qLisse(X)$, the category $\bH^{\on{access}}$ is the category $\iLisse(X)$ introduced
in \secref{sss:ind-lisse}. 

\medskip

In the example $\bH=\Rep(\sH)$, where $\sH$ is an algebraic group of finite type, the functor
\eqref{e:un-ren H} is an equivalence. 

\sssec{}

The conditions on $\bH$ imply that the unit object $\one_\bH\in \bH$ is compact. Let $\inv_\bH$
denote the functor
$$\CHom_\bH(\one_\bH,-):\bH\to \Vect,$$
i.e., the right adjoint of the unit functor $\Vect_\sfe\overset{\sfe\mapsto \one_\bH}\longrightarrow \bH$. 

\medskip

For the sequel, we record the following:

\begin{lem} \label{l:coinv}
The unit functor $\Vect_\sfe\to \bH$ admits a left adjoint (to be denoted $\coinv_\bH$). 
\end{lem}

\begin{proof} 

First, the functor $\coinv_\bH$ is defined on objects from $\bH^c\cap \bH^\heartsuit$. Indeed,
$$\coinv_\bH(\bh)=\left(\CHom_\bH(\bh,\one_\bH)\right)^\vee.$$

Since $\bH^{\geq -n}$ is generated under filtered colimits by $\bH^c\cap \bH^\heartsuit$, we obtain 
that $\coinv_\bH$ is defined on $\bH^{\geq -n}$ for any $n$.

\medskip

We now claim that for an arbitrary $\bh\in \bH$, the value of $\coinv_\bH$ on it is given by
$$\underset{n}{\on{lim}}\, \coinv_\bH(\tau^{\geq -n}(\bh)).$$

Indeed, for every $m$, the $m$-th cohomology of the system
$$n\mapsto  \coinv_\bH(\tau^{\geq -n}(\bh))$$
stabilizes (since $\coinv_\bH$ is right t-exact), and the above object has the required adjunction 
property by the left-completeness of $\bH$.

\end{proof} 

\sssec{Example}

For $\bH=\qLisse(X)$, the functor $\coinv_\bH$ identifies with the functor 
of ``cochains with compact supports"
\begin{equation} \label{e:coinv as Cc}
E\mapsto \on{C}^\cdot_c(X,E\overset{*}\otimes \omega_X).
\end{equation}

\medskip

For $\bH=\Rep(\sH)$, the functors $\inv_\bH$ and $\coinv_\bH$ are the usual functors $\sH$-invariants
and coinvariants, respectively, to be henceforth denoted $\inv_\sH$ and $\coinv_\sH$. 

\ssec{The prestack $\Maps(\Rep(\sG),\bH)$ and the abstract version of \thmref{t:main 1}}  \label{ss:abs LocSys}

In this subsection we will introduce an abstract version of the prestacks $\LocSys^{\on{restr}}_{\sG}(X)$.

\sssec{}

Let $\bH$ be as in \secref{sss:pre-Tannakian}, and let $\sG$ be a connected algebraic group.

\medskip

We define the prestack $\bMaps(\Rep(\sG),\bH)$ by sending 
an affine scheme $S$ to the space of right t-exact symmetric monoidal functors
\begin{equation} \label{e:defn coMaps}
\Rep(\sG)\to \QCoh(S)\otimes \bH.
\end{equation} 

\sssec{}

We are now ready to state an abstract version of \thmref{t:main 1}:

\begin{thm} \label{t:main 1 abs}
Assume that $\bH$ is a gentle Tannakian category. Then the prestack $\bMaps(\Rep(\sG),\bH)$
has the properties listed in \thmref{t:main 1}.
\end{thm}

\sssec{A rigidified version} Along with $\bMaps(\Rep(\sG),\bH)$, we can consider its rigidified version. 
We define the prestack $\bMaps(\Rep(\sG),\bH)^{\on{rigid}}$ by sending an affine scheme $S$ to the space of
symmetric monoidal functors
\begin{equation} \label{e:to be rigidified}
\Rep(\sG)\to \QCoh(S)\otimes \bH,
\end{equation}
equipped with an identification of the composition
$$\Rep(\sG)\to \QCoh(S)\otimes \bH \overset{\on{Id}_{\QCoh(S)}\otimes \oblv_\bH}\longrightarrow \QCoh(S)$$
with the forgetful functor
$$\Rep(\sG) \overset{\oblv_\sG}\to \Vect_\sfe\overset{\CO_S}\longrightarrow \QCoh(S),$$
as symmetric monoidal functors. In other words,
\begin{equation} \label{e:rigid as fiber prod abs}
\bMaps(\Rep(\sG),\bH)^{\on{rigid}}\simeq \bMaps(\Rep(\sG),\bH)\underset{\on{pt}/\sG}\times \on{pt}.
\end{equation} 

\medskip

The remarks pertaining to the replacement 
$$\LocSys^{\on{restr}}_{\sG}(X) \rightsquigarrow \LocSys^{\on{restr},\on{rigid}_x}_{\sG}(X)$$
apply verbatim to the replacement
$$\bMaps(\Rep(\sG),\bH)\to \bMaps(\Rep(\sG),\bH)^{\on{rigid}}.$$

In particular, \thmref{t:main 1 abs} will follow once we prove its rigidified version:

\begin{thm} \label{t:main 1 rigid abs}
Let $\bH$ be as in \thmref{t:main 1 abs}. Then the prestack $\bMaps(\Rep(\sG),\bH)^{\on{rigid}}$
is an \'etale stack equal to a disjoint union of \'etale stacks $\CY$ with 
properties \em{(a)-(d)} listed in  \thmref{t:main 1}. 
\end{thm}

\section{Ind-representability and the beginning of proof of \thmref{t:main 1 abs}} \label{s:def}

In this section we will study deformation theory (=infinitesimal)
properties of $\LocSys^{\on{restr}}_{\sG}(X)$, or more generally, $\bMaps(\Rep(\sG),\bH)$. 
Most of these properties follow easily from the definitions, apart from some issues of convergence. 

\medskip

We will conclude that $\bMaps(\Rep(\sG),\bH)^{\on{rigid}}$ is an ind-scheme. We will show this
by reducing to the assertion that for a pair of algebraic groups $\sH$ and $\sG$, the prestack
of maps $\bMaps_{\on{Grp}}(\sH,\sG)$ is an ind-affine ind-scheme. 

\ssec{Convergence} \label{ss:convergence}

In this subsection we begin the investigation of infinitesimal properties of $\bMaps(\Rep(\sG),\bH)$
(resp., $\bMaps(\Rep(\sG),\bH)^{\on{rigid}}$). 

\medskip

We start with the most basic one--the property of being convergent\footnote{See \cite[Chapter 2, Sect. 1.4]{GR1} for what this means.},
which is one of the ingredients in the condition of being almost of finite type, stated in \thmref{t:main 1 abs}.

\sssec{} \label{sss:conv to check}

By the definition of convergence, we need to show that for a (derived) affine test-scheme $S$, the map
\begin{equation} \label{e:conv}
\Maps(S,\bMaps(\Rep(\sG),\bH)) \to \underset{n}{\on{lim}}\, \Maps({}^{\leq n}\!S,\bMaps(\Rep(\sG),\bH))
\end{equation} 
is an isomorphism, where $S\mapsto  {}^{\leq n}\!S$ denotes the $n$-th coconnective truncation, i.e., the operation
$$R\mapsto \tau^{\geq -n}(R)$$
at the level of rings. 

\sssec{}

In what follows we will use the following assertion:

\begin{lem} \label{l:left-compl convergent}
Let $\bC$ be a category equipped with a t-structure. 
Then for $S$ as above we have: 

\medskip

\noindent{\em(a)} If $\bC$ is left-complete, then $\bC\otimes \QCoh(S)$ is also left-complete. 

\medskip

\noindent{\em(a')} More generally, if $\bC^\wedge$ is the left completion of $\bC$, then the functor
$$\QCoh(S)\otimes \bC \to \QCoh(S)\otimes \bC^\wedge$$
identifies $\QCoh(S)\otimes \bC^\wedge$ with the left completion of $\QCoh(S)\otimes \bC$. 

\medskip

\noindent{\em(b)} If $\bC$ is left-complete, then the functor
$$(\bC\otimes \QCoh(S))^{\leq 0} \to \underset{n}{\on{lim}}\, (\bC\otimes \QCoh({}^{\leq n}\!S))^{\leq 0}$$
is an equivalence.
\end{lem}

\begin{proof}

Points (a) and (a') follow from the fact that the functor
$$\bC\otimes \QCoh(S) \overset{\on{Id}\otimes \Gamma(S,-)}\longrightarrow \bC$$
is t-exact, conservative and commutes with limits.

\medskip

Point (b) follows from point (a) and the fact that for any $n$, the functor
$$(\bC\otimes \QCoh(S))^{\leq 0,\geq -n} \to (\bC\otimes \QCoh({}^{\leq m}\!S))^{\leq 0,\geq -n}$$
is an equivalence, whenever $m\geq n$.

\end{proof}

\sssec{} \label{sss: proof of conv}

We are now ready to prove that \eqref{e:conv} is an equivalence. 

\begin{proof}

Since $\Rep(\sG)$ is the derived category of its heart and the monoidal operation is t-exact,
the space of right t-exact symmetric monoidal functors
$$\Rep(\sG) \to \QCoh(S)\otimes \bH$$
is isomorphic to the space of symmetric monoidal
functors
$$\Rep(\sG)^\heartsuit\to (\QCoh(S)\otimes \bH)^{\leq 0},$$
and similarly for every $^{\leq n}\!S$.

\medskip

The assertion now follows from the assumption that $\bH$ is left-complete in its t-structure and \lemref{l:left-compl convergent}(b). 

\end{proof}

\begin{rem}

Note that the above argument used the fact that $\bH$ is left-complete in its t-structure
(in order to be able to apply \lemref{l:left-compl convergent}). 

\medskip

This is the reason that we have to work with
$\qLisse(X)$ (resp., $\bH$) rather than with the more manageable category $\iLisse(X)$ (resp., $\bH^{\on{access}}$).

\end{rem}

\sssec{}

Despite the previous remark, we will now show that one can work $\bH^{\on{access}}$ instead of $\bH$ as long
as we evaluate our prestack on \emph{eventually coconnective} affine schemes. 

\medskip

Recall that an affine scheme $S$ is said to be eventually coconnective if it equals $^{\leq n}\!S$ for some $n$
(i.e., if its structure ring has cohomologies in finitely many degrees). 

\begin{prop} \label{p:replace by Ind-Lisse}
Suppose that $S$ is eventually coconnective. Then the functor \eqref{e:un-ren H} defines 
an isomorphism from the space of (right t-exact) symmetric monoidal
functors
$$\Rep(\sG) \to \QCoh(S) \otimes \bH^{\on{access}}$$
to the space of (right t-exact) symmetric monoidal functors
$$\Rep(\sG) \to \QCoh(S) \otimes \bH.$$
\end{prop} 

\begin{proof}

The space of (continuous) symmetric monoidal functors
$$\Rep(\sG) \to \QCoh(S) \otimes \bH^{\on{access}}$$
maps isomorphically to the space of symmetric monoidal functors
$$\Rep(\sG)^c\to \QCoh(S) \otimes \bH^{\on{access}},$$
and similarly for $\bH^{\on{access}}$ replaced by $\bH$.

\medskip

Since every object of $\Rep(\sG)^c$ is dualizable, it suffices to show that the embedding 
\begin{equation} \label{e:embed IndLisse S}
\QCoh(S) \otimes \bH^{\on{access}}\hookrightarrow \QCoh(S)\otimes \bH
\end{equation} 
induces an equivalence on the subcategories of dualizable objects. The functor \eqref{e:embed IndLisse S}
is a priori fully faithful because the functor \eqref{e:un-ren H} is such, while $\QCoh(S)$ is dualizable. 

\medskip

Note that by \lemref{l:left-compl convergent}(a'), the functor \eqref{e:embed IndLisse S}
identifies $\QCoh(S)\otimes \bH$ with the left completion of 
$\QCoh(S) \otimes \bH^{\on{access}}$. 

\medskip

Hence, it suffices to show that (for $S$ eventually coconnective), any dualizable object in the category 
$\QCoh(S) \otimes \bH$ is
bounded below (in the sense of the t-structure). 

\medskip

Since the functor $\oblv_\bH$ is conservative and t-exact, the functor 
$$\QCoh(S)\otimes \bH \overset{\on{Id}\otimes \oblv_\bH}\longrightarrow \QCoh(S)$$
is also t-exact and conservative. 

\medskip

Hence, it is enough to show that $\on{Id}\otimes \oblv_\bH$ 
sends dualizable objects to objects bounded below. However, $\on{Id}\otimes \oblv_\bH$ is symmetric monoidal,
the assertion follows from the fact that dualizable objects in $\QCoh(S)$ (for $S$ 
eventually coconnective) are bounded below.

\end{proof}

\ssec{Deformation theory: statements}  \label{ss:deform}

In this subsection we will formulate the deformation theory properties of 
$\bMaps(\Rep(\sG),\bH)$, along with 
its version $\bMaps(\Rep(\sG),\bH)^{\on{rigid}}$. 

\medskip

This will
be an ingredient in the proof of that fact that $\bMaps(\Rep(\sG),\bH)^{\on{rigid}}$ 
is an ind-scheme, stated in \thmref{t:main 1 rigid abs}. 

\sssec{}

We will prove:

\begin{prop}  \label{p:def} \hfill

\smallskip

\noindent{\em(a)}
The prestack $\bMaps(\Rep(\sG),\bH)$ admits a $(-1)$-connective corepresentable deformation theory.

\smallskip

\noindent{\em(b)}
For $S\in \affSch_{/\sfe}$ and an $S$-point 
$$\sF:\Rep(\sG)\to \QCoh(S)\otimes \bH$$
of $\bMaps(\Rep(\sG),\bH)$, the cotangent 
space $T^*_\sF(\bMaps(\Rep(\sG),\bH))\in \QCoh(S)^{\leq 1}$ identifies with 
$$(\on{Id}\otimes \coinv_\bH)(\sF(\sg^\vee))[-1],$$
where $\sg$ is the Lie algebra of $\sG$. 

\end{prop}

\begin{rem} 

We refer the reader to \cite[Chapter 1, Definition 7.1.5(a)]{GR2}, where it is explained what it means to admit
a $(-n)$-connective corepresentable deformation theory. In fact, there are three conditions:

\smallskip

\noindent(i) The first one is that the prestack admits deformation theory (i.e. admits pro-cotangent spaces
that are functorial in the test-scheme, and is infinitesimally cohesive);

\smallskip

\noindent(ii) The adjective ``corepresentable" refers to the fact that the pro-cotangent spaces are
actually objects (of $\QCoh(S)^{<\infty}$, where $S$ is the test-scheme), and not only pro-objects.

\smallskip

\noindent(iii) The adjective ``$(-n)$-connective" refers to the fact that cotangent spaces 
actually belong to $\QCoh(S)^{\leq n}$.

\end{rem}

\sssec{}
As a consequence, we deduce:

\begin{cor}\label{c:etale descent}
The prestack $\bMaps(\Rep(\sG),\bH)$ satisfies \'etale descent. In particular, it identifies with the \'etale quotient 
$\bMaps(\Rep(\sG),\bH)^{\on{rigid}}/\sG$.
\end{cor}

\begin{proof}
By \propref{p:def}(a) and \cite[Chapter 1, Proposition 8.2.2]{GR2}, it suffices to show that the underlying classical prestack satisfies \'etale descent.  
Thus, by \propref{p:replace by Ind-Lisse}, it suffices to show that the functor
$$ S \mapsto \{ \text{right $t$-exact symmetric monoidal functors } \Rep(\sG) \to \QCoh(S) \otimes  \bH^{\on{access}}\}$$
satisfies \'etale descent, for $S$ a classical affine scheme.  

\medskip

Since $\bH^{\on{access}}$ is compactly generated (and in particular dualizable), the functor
$ - \otimes \bH^{\on{access}}$ preserves limits.  The result now follows from the fact that $\QCoh(S)$ satisfies
\'etale descent.

\medskip

Now, the assertion that
$$\bMaps(\Rep(\sG),\bH) \simeq \bMaps(\Rep(\sG),\bH)^{\on{rigid}}/\sG$$
follows from the fact that $\bMaps(\Rep(\sG),\bH)$ satisfies \'etale descent and the identification
\eqref{e:rigid as fiber prod abs}. 

\end{proof}

Using the presentation of $\bMaps(\Rep(\sG),\bH)^{\on{rigid}}$ as \eqref{e:rigid as fiber prod abs}, from 
\propref{p:def}, we obtain:

\begin{cor}  \label{c:def}  \hfill

\noindent{\em(a)}
The prestack $\bMaps(\Rep(\sG),\bH)^{\on{rigid}}$ admits a connective corepresentable deformation theory.

\smallskip

\noindent{\em(b)}
For $S\in \affSch_{/\sfe}$ and an $S$-point $\sF$ of $\bMaps(\Rep(\sG),\bH)^{\on{rigid}}$, we have a canonical identification 
$$T^*_\sF(\bMaps(\Rep(\sG),\bH)^{\on{rigid}})\simeq \on{Fib}\left((\CO_S\otimes \sg^\vee) \to (\on{Id}\otimes \coinv_\bH)(\sF(\sg^\vee))\right).$$

\smallskip

\noindent{\em(b')} The object $T^*_\sF(\bMaps(\Rep(\sG),\bH)^{\on{rigid}})\in \QCoh(S)^{\leq 0}$ 
belongs to $\QCoh(S)^{\leq 0}\cap \QCoh(S)^c$. 

\end{cor}

\begin{proof}

The fact that $\bMaps(\Rep(\sG),\bH)^{\on{rigid}}$ admits corepresentable deformation theory
follows formally from \propref{p:def}(a) and \eqref{e:rigid as fiber prod abs}. 

\medskip

The functoriality of the identification in \propref{p:def}(b) with respect to $\bH$ implies that 
the codifferential of the map 
$$\bMaps(\Rep(\sG),\bH)\to \on{pt}/\sG$$
at a point $\sF\in \bMaps(\Rep(\sG),\bH)$ is given by
\begin{equation} \label{e:codiff F}
(\on{Id}\otimes \oblv_\bH)(\sF(\sg^\vee))[-1]\to (\on{Id}\otimes \coinv_\bH)(\sF(\sg^\vee))[-1].
\end{equation}

This implies the assertion of point (b). Since the map \eqref{e:codiff F} induces a \emph{surjection}
$$H^0\left((\on{Id}\otimes \oblv_\bH)(\sF(\sg^\vee))\right)\to 
H^0\left((\on{Id}\otimes \coinv_\bH)(\sF(\sg^\vee))\right),$$
this implies that the cotangent spaces of $\bMaps(\Rep(\sG),\bH)^{\on{rigid}}$ are connective. 

\medskip

For point (b') we will show that for any dualizable $V\in \Rep(\sG)$, the object
$$(\on{Id}\otimes \coinv_\bH)(\sF(V))\in \QCoh(S)^{\leq 0}$$ 
belongs to $\QCoh(S)^{\leq 0}\cap \QCoh(S)^c$. 

\medskip

It suffices to show that this happens after restriction to any truncation of $S$. 
Hence, we can assume that $S$ is eventually coconnective. In this case, 
by \propref{p:replace by Ind-Lisse}, we can regard $\sF$ as a functor  
$$\Rep(\sG)\to \QCoh(S)\otimes \bH^{\on{access}}.$$

\medskip

The object $\sF(V)\in \QCoh(S)\otimes \bH^{\on{access}}$ is dualizable. Since both $\QCoh(S)$ 
and $\bH^{\on{access}}$ are rigid, the category $\QCoh(S)\otimes \bH^{\on{access}}$ is also rigid.
Hence $\sF(V)$ is compact as an object of $\QCoh(S)\otimes \bH^{\on{access}}$. 

\medskip

The composite functor
$$\bH^{\on{access}} \overset{\text{\eqref{e:un-ren H}}}\to \bH\overset{\coinv_\bH}\to \Vect_\sfe$$
is the left adjoint of the unit functor. Hence, it preserves compactness. Hence, so does
the functor 
$$\on{Id}\otimes \coinv_\bH: \bH^{\on{access}}\otimes \QCoh(S)\to \QCoh(S).$$

\end{proof}


\ssec{Establishing deformation theory}

This subsection is devoted to the proof of the fact that $\bMaps(\Rep(\sG),\bH)$ admits deformation
theory, which is part of the assertion of \propref{p:def}(a). 

\sssec{}

By \cite[Chapter 1, Proposition 7.2.5]{GR2},  we need to show that for a push-out diagram of affine schemes
\begin{equation} \label{e:pushout}
\CD
S_1 @>>> S_2 \\
@VVV  @VVV  \\
S'_1 @>>> S'_2,
\endCD
\end{equation}
where $S_1\to S'_1$ is a nilpotent embedding, the diagram
$$
\CD
\Maps(S'_2,\bMaps(\Rep(\sG),\bH))  @>>> \Maps(S'_1,\bMaps(\Rep(\sG),\bH))  \\
@VVV  @VVV  \\
\Maps(S_2,\bMaps(\Rep(\sG),\bH))  @>>> \Maps(S_1,\bMaps(\Rep(\sG),\bH))
\endCD
$$
is a pullback square of spaces. 

\medskip

Since $\bMaps(\Rep(\sG),\bH)$ is convergent, 
we can assume that affine schemes in \eqref{e:pushout} are eventually coconnective. 

\sssec{}

Using \propref{p:replace by Ind-Lisse}, for $S$ eventually coconnective, 
we interpret $\Maps(S,\bMaps(\Rep(\sG),\bH))$ as the space of right t-exact symmetric monoidal functors
$$\Rep(\sG)^c\to (\QCoh(S)\otimes \bH^{\on{access}})^{\on{dualizable}}.$$

\medskip

Hence, it suffices to show that in the situation of \eqref{e:pushout}, the diagram
\begin{equation} \label{e:pullback sym mon}
\CD
(\QCoh(S'_2)\otimes \bH^{\on{access}})^{\on{dualizable}}  @>>> (\QCoh(S'_1)\otimes \bH^{\on{access}})^{\on{dualizable}}  \\
@VVV   @VVV   \\
(\QCoh(S_2)\otimes \bH^{\on{access}})^{\on{dualizable}}  @>>> (\QCoh(S_1)\otimes \bH^{\on{access}})^{\on{dualizable}}
\endCD
\end{equation}
is a pullback square of (small, symmetric monoidal) categories. 

\sssec{}

Note that by \cite[Chapter 1, Proposition 1.4.2]{GR2}, the functor
$$\QCoh(S'_2) \to \QCoh(S_2)\underset{\QCoh(S_1)}\times \QCoh(S'_1)$$
is fully faithful (but \emph{not} necessarily an equivalence, see \cite[Chapter 1, Remark 1.4.3]{GR2}). 

\medskip

Since $\bH^{\on{access}}$ is dualizable as a DG category, the functor
\begin{equation} \label{e:pushout iLisse}
\QCoh(S'_2)\otimes \bH^{\on{access}} \to 
(\QCoh(S_2)\otimes \bH^{\on{access}})
\underset{\QCoh(S_1)\otimes \bH^{\on{access}}}\times (\QCoh(S'_1)\otimes \bH^{\on{access}})
\end{equation} 
is also fully faithful. 

\medskip

Hence, the functor
\begin{multline} \label{e:pushout iLisse dual}
(\QCoh(S'_2)\otimes \bH^{\on{access}})^{\on{dualizable}} \to \\
\to (\QCoh(S_2)\otimes \bH^{\on{access}})^{\on{dualizable}}
\underset{(\QCoh(S_1)\otimes \bH^{\on{access}})^{\on{dualizable}}}\times (\QCoh(S'_1)\otimes \bH^{\on{access}})^{\on{dualizable}}
\end{multline} 
is fully faithful.

\medskip

It remains to prove that \eqref{e:pushout iLisse dual} is essentially surjective. The argument that follows is applicable
to $\bH^{\on{access}}$ replaced by any proper compactly generated rigid symmetric monoidal category $\bA$. 

\sssec{}

For an affine scheme $S$, the monoidal category $\QCoh(S)$ is rigid (see \cite[Chapter 1, Sect. 9]{GR1} for what this means).
Since $\bA$ was assumed rigid as well, we obtain that so are the categories of the form
$$\QCoh(S)\otimes \bA.$$

In particular,
$$(\QCoh(S)\otimes \bA)^{\on{dualizable}} = (\QCoh(S)\otimes \bA)^c$$
as subcategories of $\QCoh(S)\otimes \bA$.

\medskip
Since $\bA$ is rigid, it is in particular self-dual (see \cite[Chapter 1, Sect. 9.2]{GR1}); i.e. we have a canonical equivalence
$$ \bA \simeq \bA^{\vee} .$$

\medskip

Now, since $\bA$ is \emph{proper} (i.e., $\CHom$'s between compact objects lie in $\Vect_\sfe^c$), 
the equivalence
$$ \QCoh(S) \otimes \bA^{\vee} \simeq \on{Funct}(\bA, \QCoh(S)) $$
restricts to a fully faithful embedding
$$ (\QCoh(S) \otimes \bA^{\vee})^c \hookrightarrow \on{Funct}(\bA^c, \QCoh(S)^c) .$$
Thus,
we have a fully
faithful embedding
\begin{multline*}
(\QCoh(S)\otimes \bA)^c \simeq (\QCoh(S)\otimes \bA^\vee)^c
\hookrightarrow \on{Funct}(\bA^c,\QCoh(S)^c)=
\on{Funct}(\bA^c,\QCoh(S)^{\on{dualizable}}).
\end{multline*}

We now apply \cite[Chapter 8, Proposition 3.3.2]{GR2}, which implies that for the diagram \eqref{e:pushout}, the diagram
of categories
$$
\CD
\QCoh(S'_2)^{\on{dualizable}}  @>>> \QCoh(S'_1)^{\on{dualizable}}  \\
@VVV   @VVV   \\
\QCoh(S_2)^{\on{dualizable}}  @>>> \QCoh(S_1)^{\on{dualizable}}
\endCD
$$
is a pull-back square. Hence, so is the diagram
$$
\CD
\on{Funct}(\bA^c,\QCoh(S'_2)^{\on{dualizable}})  @>>> \on{Funct}(\bA^c,\QCoh(S'_1)^{\on{dualizable}})  \\
@VVV   @VVV   \\
\on{Funct}(\bA^c,\QCoh(S_2)^{\on{dualizable}})  @>>> \on{Funct}(\bA^c,\QCoh(S_1)^{\on{dualizable}}). 
\endCD
$$

\sssec{}

Hence, given an object $M$ in the right-hand side of 
\begin{multline} \label{e:pushout iLisse dual C}
(\QCoh(S'_2)\otimes \bA)^{\on{dualizable}} \to \\
\to (\QCoh(S_2)\otimes \bA)^{\on{dualizable}}
\underset{(\QCoh(S_1)\otimes \bA)^{\on{dualizable}}}\times (\QCoh(S'_1)\otimes \bA)^{\on{dualizable}},
\end{multline} 
we can create an object $M'$
in $\on{Funct}(\bA^c,\QCoh(S'_2)^{\on{dualizable}})$, so that $M$ and $M'$ have the same image in 
$$\on{Funct}(\bA^c,\QCoh(S_2)^{\on{dualizable}}) \underset{\on{Funct}(\bA^c,\QCoh(S_1)^{\on{dualizable}})}
\times \on{Funct}(\bA^c,\QCoh(S'_1)^{\on{dualizable}}).$$

Using the fully faithful embedding
$$\on{Funct}(\bA^c,\QCoh(S)^{\on{dualizable}})\hookrightarrow 
\QCoh(S)\otimes \bA,$$
we obtain that there exists an object $M''$ in the left-hand side of 
\begin{equation} \label{e:pushout iLisse C}
\QCoh(S'_2)\otimes \bA \to 
(\QCoh(S_2)\otimes \bA)
\underset{\QCoh(S_1)\otimes \bA}\times (\QCoh(S'_1)\otimes \bA),
\end{equation} 
whose image
in the right-hand side is isomorphic to that of $M$.

\medskip

Thus, it remains to show that $M''$ is compact. 

\sssec{}

Since the functor \eqref{e:pushout iLisse C} commutes with colimits and is fully faithful, an object 
in the left-hand side of \eqref{e:pushout iLisse C} is compact if its image 
in the right-hand side of \eqref{e:pushout iLisse C} is compact. 

\medskip

This implies that $M''$ is compact, since $M$, viewed as an object of the right-hand side of \eqref{e:pushout iLisse C},
is compact.

\qed

%
%
%
%
%
%
\ssec{Calculating the (co)tangent spaces} \label{ss:calc tang}

In this subsection we will prove the remaining assertions of \propref{p:def}. To do so, 
it suffices to perform the calculation of point (b). 

\sssec{} \label{sss:calc tang}

Let 
$\CM$ be an object of $\QCoh(S)^{\leq 0}$, and let $S_\CM\in \affSch_{/\sfe}$ be the corresponding split square-zero
extension of $S$. Unwinding the definitions, we obtain that we need to construct an isomorphism
\begin{multline*}
\Maps(S_\CM,\bMaps(\Rep(\sG),\bH))\underset{\Maps(S,\bMaps(\Rep(\sG),\bH))}\times \{*\}\simeq \\
\simeq \tau^{\leq 0}\left(\CHom_{\QCoh(S)\otimes \bH}\left(\CO_S\otimes \one_\bH, (\CM\otimes \one_\bH) \otimes \sF(\sg)\right)[1]\right).
\end{multline*} 

\sssec{}

Let $\bA$ be a symmetric monoidal DG category and let $\ba\in \bA$ be an object. We regard
$\one_\bA\oplus \ba$ as an object of $\on{ComAlg}(\bA)$, the square-zero extension of $\one_\bA$ 
by means of $\ba$. Consider the category
$$(\one_\bA\oplus \ba)\mod(\bA)$$
as a symmetric monoidal category, equipped with a symmetric monoidal functor back to $\bA$,
given by
$$-\underset{\one_\bA\oplus \ba}\otimes \one_\bA.$$

\medskip

We have the following general assertion:

\begin{lem} \label{l:funct from Rep}
Given a symmetric
monoidal functor
$$\sF:\Rep(\sG)\to \bA,$$
the space of its lifts to a functor 
$$\Rep(\sG)\to (\one_\bA\oplus \ba)\mod(\bA)$$
identifies canonically with
$$\Maps_\bA(\one_\bA,\ba\otimes \sF(\sg)[1]).$$
\end{lem}

\medskip

Applying this to 
$$\bA:=\QCoh(S)\otimes \bH, \quad \ba=\CM\otimes \one_\bH,$$
we obtain the result stated in \secref{sss:calc tang}.

\ssec{Proof of ind-representability}

In this subsection we will begin the proof that the prestack $\bMaps(\Rep(\sG),\bH)^{\on{rigid}}$
is an ind-affine ind-scheme locally almost of finite type. 

\medskip

We will do so by reducing to the case when instead of the category $\bH$ we
are dealing with the category $\Rep(\sH)$ of representations of an algebraic group $\sH$. 

\sssec{}

Recall (see \cite[Chapter 2, Definition 1.1.2]{GR2}) that prestack $\CY$ is an \emph{ind-scheme} if it can be written as a \emph{filtered} colimit
\begin{equation} \label{e:ind}
\CY \simeq \underset{i}{\on{colim}}\, Y_i,
\end{equation} 
where $Y_i$ are (quasi-compact) schemes and the transition maps $Y_i\to Y_j$ are closed embeddings. 

\medskip

An ind-scheme is \emph{ind-affine} if all $Y_i$ can be chosen to be affine. 

\medskip

An ind-scheme is locally almost of finite type as a prestack if all $Y_i$ can be chosen to be
almost of finite type see (\cite[Chapter 2, Corollary 1.7.5(a)]{GR2}).

%
%
%
%
%
%

\sssec{} \label{sss:proof of ind-rep}

By \corref{c:def}(a) combined with \cite[Chapter 2, Corollary 1.3.13]{GR2}, in order to show that $\bMaps(\Rep(\sG),\bH)^{\on{rigid}}$
is an ind-affine ind-scheme, it suffices to show that its classical truncation $^{\on{cl}}\bMaps(\Rep(\sG),\bH)^{\on{rigid}}$ is a classical
ind-affine ind-scheme.

\medskip

Similarly, by \corref{c:def}(b') combined with \cite[Chapter 1, Theorem 9.1.2]{GR2}, in order to show that the prestack 
$\bMaps(\Rep(\sG),\bH)^{\on{rigid}}$ is locally almost of finite type, 
it suffices to show that $^{\on{cl}}\bMaps(\Rep(\sG),\bH)^{\on{rigid}}$ is locally of finite type as a classical
prestack. 

\sssec{}  \label{sss:Galois pre bis}

Let $\sH$ be a pro-algebraic group. Consider the prestack 
$$\bMaps_{\on{Grp}}(\sH,\sG)$$
that sends an affine scheme $S$ to the space of homomorphisms of group-schemes over $S$
\begin{equation} \label{e:group hom}
\phi:S\times \sH\to S\times \sG.
\end{equation} 

We have a naturally defined map
$$\bMaps_{\on{Grp}}(\sH,\sG)\to \bMaps(\Rep(\sG),\Rep(\sH))^{\on{rigid}},$$
and it follows from Tannaka duality that it is actually an isomorphism. 

\medskip

We will prove:

\begin{prop} \label{p:maps of groups pro} 
The prestack 
$\bMaps_{\on{Grp}}(\sH,\sG)$
is an ind-affine ind-scheme locally almost of finite type. 
\end{prop}

The proof is given in \secref{ss:proof alg grp} below. 

\begin{rem}

Let us note that for a general pro-algebraic group $\sH$ (as opposed to an algebraic group
of finite type, the category $\Rep(\sH)$ is generally \emph{not} a gentle Tannakian category
(for example, $\Ext^1_{\Rep(\sH)}(\sfe,\sfe)$ can be infinite-dimensional). 

\medskip

As a result, the connected components of the ind-scheme $\bMaps_{\on{Grp}}(\sH,\sG)$ are \emph{not}, in general,
formal affine schemes.

\end{rem}

\sssec{} \label{sss:Galois}

Consider the \emph{abelian} symmetric monoidal category $\bH^\heartsuit$, equipped
with the fiber functor $\oblv_\bH$. 

\medskip

By Tannaka duality, there exists a pro-algebraic group,
to be denoted $\sH$, such that $\bH^\heartsuit$ identifies with the abelian category of algebraic 
representations of $\sH$ and
$\oblv_\bH$ corresponds to the tautological forgetful functor. 

\begin{rem} \label{r:Galois}
For $\bH=\qLisse(X)$, the resulting group $\sH$ is $\on{Gal}(X,x)_{\on{Pro-alg}}$,
the pro-algebraic completion of the (unramified) Galois group $\on{Gal}(X,x)$ of $X$
with base point $x$. 
\end{rem} 

\sssec{}

We claim:

\begin{prop} \label{p:Tannaka red cl}
There exists a canonical isomorphism of \emph{classical} prestacks
$$^{\on{cl}}\bMaps(\Rep(\sG),\bH)^{\on{rigid}} \simeq {}^{\on{cl}}\bMaps_{\on{Grp}}(\sH,\sG).$$
\end{prop}

Note that this proposition, combined with \propref{p:maps of groups pro}, implies that
$^{\on{cl}}\bMaps(\Rep(\sG),\bH)^{\on{rigid}}$ is representable by a classical ind-affine ind-scheme locally of finite type.  
By \secref{sss:proof of ind-rep}, this implies that $\bMaps(\Rep(\sG),\bH)^{\on{rigid}}$
is an ind-affine ind-scheme locally almost of finite type. 

\begin{proof}[Proof of \propref{p:Tannaka red cl}]

Let $S=\Spec(R)$ be a classical affine scheme. As in \secref{sss: proof of conv}, the value of 
$\bMaps(\Rep(\sG),\bH)^{\on{rigid}}$ on $S$ is the category of symmetric monoidal functors
\begin{equation} \label{e:functors F ab}
\sF:\Rep(\sG)^\heartsuit\to (R\mod\otimes \bH)^{\leq 0},
\end{equation} 
equipped with an identification of the composition 
$$\Rep(\sG)^\heartsuit\to (R\mod\otimes \bH)^{\leq 0}\overset{\on{Id}\otimes \oblv_\bH}\longrightarrow R\mod$$
with 
$$\Rep(\sG)^\heartsuit \overset{\oblv_\sG}\to \Vect_\sfe^\heartsuit \overset{\on{unit}}\to R\mod.$$

Such functors $\sF$ necessarily take values in the abelian monoidal category
$$(R\mod\otimes \bH)^\heartsuit\simeq R\mod(\bH^\heartsuit).$$

Similarly, $S$-values of $\bMaps_{\on{Grp}}(\sH,\sG)$ are symmetric monoidal functors 
$$\Rep(\sG)^\heartsuit\to R\mod(\Rep(\sH)^\heartsuit),$$
equipped with an identification of the composition
$$\Rep(\sG)^\heartsuit\to R\mod(\Rep(\sH)^\heartsuit)\overset{\oblv_\sH}\to R\mod$$
with 
$$\Rep(\sG)^\heartsuit \overset{\oblv_\sG}\to \Vect_\sfe^\heartsuit \overset{\on{unit}}\to R\mod.$$

Thus, the two sets of data are manifestly isomorphic.

\end{proof}

\ssec{Proof in the case of algebraic groups} \label{ss:proof alg grp}

In this subsection we will prove \propref{p:maps of groups pro}.

\sssec{}

We claim that in order to prove \propref{p:maps of groups pro}, it suffices to consider the case when $\sH$ is an
algebraic group of finite type. Indeed, write $\sH$ as an cofiltered limit
$$\sH:=\underset{\alpha}{\on{lim}}\, \sH^\alpha,$$
where $\sH^\alpha$ are algebraic groups of finite type and the transition maps are surjective.

\medskip

Then,
$$\bMaps_{\on{Grp}}(\sH,\sG)\simeq \underset{\alpha}{\on{colim}}\, 
\bMaps_{\on{Grp}}(\sH^\alpha,\sG).$$
and the transition maps are closed embeddings (see Remark \ref{r:closed emb maps} for the latter statement). 

\medskip

Hence, it suffices to show the following:

\begin{prop} \label{p:maps of groups} 
For a pair of algebraic groups $\sH$ and $\sG$ of finite type, the prestack 
$\bMaps_{\on{Grp}}(\sH,\sG)$
is an ind-affine ind-scheme locally almost of finite type. 
\end{prop}

The rest of this subsection is devoted to the proof of \propref{p:maps of groups}.

\sssec{}

We will show that the \emph{classical prestack underlying} $\bMaps_{\on{Grp}}(\sH,\sG)$ is an 
ind-affine ind-scheme locally of finite type. This will imply that $\bMaps_{\on{Grp}}(\sH,\sG)$ 
is an ind-affine ind-scheme locally almost of finite type by the argument in \secref{sss:proof of ind-rep}
applied to $\bH=\Rep(\sH)$. 

\begin{rem} \label{r:cotan Hom}

The description of the cotangent space in \corref{c:def}(b) can be translated as follows: 

\medskip

For an affine test scheme $S$,
and an $S$-point $\phi$ of $\bMaps_{\on{Grp}}(\sH,\sG)$, the cotangent space
$$T^*_\phi(\bMaps_{\on{Grp}}(\sH,\sG))\in \QCoh(S)$$
identifies with
$$\on{Fib}\left(\sg^\vee\otimes \CO_S \to \on{coinv}_{\sH}(\sg^\vee\otimes \CO_S)\right),$$
where $\on{coinv}_{\sH}$ stands for $\sH$-coinvariants, and $\sg^\vee\otimes \CO_S$ acquires
a structure of $\sH$-module via $\phi$. 

\end{rem} 

\sssec{}

From now on until the end of this section, we will consider the underlying classical prestacks and omit the superscript ``cl" from the notation. 

\medskip

For a pair of affine schemes of finite type $Y_1,Y_2$, consider the prestack 
$$\bMaps_{\on{Sch}}(Y_1,Y_2), \quad S\mapsto  \Hom(S\times Y_1,Y_2).$$

We claim that $\bMaps_{\on{Sch}}(Y_1,Y_2)$ is representable by an ind-affine ind-scheme locally of finite type. 

\medskip

Indeed, the formation of $\bMaps_{\on{Sch}}(Y_1,Y_2)$ commutes with limits in $Y_2$, and every affine scheme
of finite type
can be written as a (finite) limit of copies of $\BA^1$. This reduces the assertion to the case when $Y_2=\BA^1$. 

\medskip

However, for any prestack $\CY$
$$\bMaps_{\on{Sch}}(\CY,\BA^1)\simeq \on{Tot}(W),\quad W:=\Gamma(\CY,\CO_{\CY}),$$
while for any $W\in (\Vect_\sfe)^{\geq 0}$, the prestack $\on{Tot}(W)$ is indeed a classical ind-affine ind-scheme locally of finite type: write
$$W\simeq \underset{i}{\on{colim}}\, W_i,$$
with $W_i$ finite dimensional, and we have
$$\on{Tot}(W)\simeq \underset{i}{\on{colim}}\, \on{Tot}(W_i),$$
while
$$\on{Tot}(W_i)\simeq \Spec(\Sym(W_i^\vee)).$$

\sssec{}

Setting $Y_1=\sH$, $Y_2=\sG$, we obtain that 
$$\bMaps_{\on{Sch}}(\sH,\sG)$$ 
is an ind-affine ind-scheme locally of finite type. 

\medskip

Now, $\bMaps_{\on{Grp}}(\sH,\sG)$ can be expressed as a fiber product of copies of $\bMaps_{\on{Sch}}(\sH,\sG)$
and $\bMaps_{\on{Sch}}(\sH^2,\sG)$. This implies the assertion of \propref{p:maps of groups}. 

\qed[\propref{p:maps of groups}]

\begin{rem} \label{r:closed emb maps}
For future reference, we note that if $Y_2\hookrightarrow Y'_2$ is a closed embedding of affine schemes, then 
the corresponding map
$$\bMaps_{\on{Sch}}(Y_1,Y_2)\to \bMaps_{\on{Sch}}(Y_1,Y'_2)$$
is a closed embedding of functors. 

\medskip

In particular, if $\sG'\hookrightarrow \sG$ is a closed subgroup, then the map
$$\bMaps_{\on{Grp}}(\sH,\sG')\to \bMaps_{\on{Grp}}(\sH,\sG)$$
is a closed embedding.

\medskip

Similarly, for a surjection $\sH'\twoheadrightarrow \sH$, the map
$$\bMaps_{\on{Grp}}(\sH,\sG)\to \bMaps_{\on{Grp}}(\sH',\sG)$$
is a closed embedding, since
$$\bMaps_{\on{Grp}}(\sH,\sG) \simeq \bMaps_{\on{Grp}}(\sH',\sG)\underset{\bMaps_{\on{Grp}}(\sH'',\sG)}\times \on{pt},$$
where $\sH'':=\on{ker}(\sH'\to \sH)$. 

\end{rem} 

\section{Uniformization and the end of proof of \thmref{t:main 1 abs}} \label{s:uniformization}

In this section we will finish the proof of \thmref{t:main 1 abs}, while introducing a tool of independent interest:
a uniformization of $\LocSys^{\on{restr}}_\sG(X)$ by \emph{algebraic stacks} associated to parabolic
subgroups in $\sG$ and \emph{irreducible} local systems for their Levi quotients. 

\ssec{What is there left to prove?}

\sssec{}

We claim that in order to finish the proof of \thmref{t:main 1 rigid abs}, it remains to show the following:

\begin{thm}  \label{t:underlying reduced}
The underlying reduced prestack of $\bMaps(\Rep(\sG),\bH)^{\on{rigid}}$ is a disjoint union 
of affine schemes, sheafified in the Zariski/\'etale topology.
\end{thm}

Let us show how \thmref{t:underlying reduced} implies \thmref{t:main 1 rigid abs}. 

\sssec{}

Indeed, we have already shown that $\bMaps(\Rep(\sG),\bH)^{\on{rigid}}$ is an 
ind-affine ind-scheme locally almost of finite type. Combined with \thmref{t:underlying reduced},
this implies points (a,b,c) of \thmref{t:main 1}. 

\medskip

To prove point (d), we quote the following result, which is a particular case of \cite[Theorem 18.2.3.2]{Lu3}
(combined with \cite[Proposition. 6.7.4]{GR3}):

\begin{thm}  \label{t:formal}
Let $\CY$ be an ind-scheme locally almost of finite type 
with the following properties:

\medskip

\noindent{\em(i)} $^{\on{red}}\CY$ is an affine scheme;

\medskip

\noindent{\em(ii)} For any $(S,y)\in \affSch_{/\CY}$, the cotangent space $T^*_y(\CY)\in \on{Pro}(\QCoh(S)^{\leq 0})$
actually belongs to $\QCoh(S)^{\leq 0}$.

\medskip

Then $\CY$ can be written in the form \eqref{e:presentation A}. 

\end{thm}

For the sake of completeness, we will outline the proof of \thmref{t:formal} in \secref{s:formal}. 
%
%
%
%

\ssec{Uniformization} \label{ss:uniformization}

In this subsection we will begin the proof of \thmref{t:underlying reduced}. The method is based on 
constructing an algebraic stack that maps dominantly onto $^{\on{red}}\bMaps(\Rep(\sG),\bH)$. This
construction will also shed some light on ``what $\bMaps(\Rep(\sG),\bH)$ looks like". 

\sssec{} \label{sss:properties of uniformization}

Having proved that $\bMaps(\Rep(\sG),\bH)^{\on{rigid}}$ is an ind-affine ind-scheme locally almost of finite type, we know that
each connected component $^{\on{red}}\bMaps(\Rep(\sG),\bH)^{\on{rigid}}$ is a filtered colimit of reduced affine schemes
along closed embeddings. Hence, in order to prove \thmref{t:underlying reduced}, it suffices to show that these colimits stabilize. 

\medskip

We will achieve this by the following construction. We will find an \emph{algebraic stack} locally almost of finite type
$\wt\bMaps(\Rep(\sG),\bH)$, equipped with a map
$$\pi:\wt\bMaps(\Rep(\sG),\bH)\to \bMaps(\Rep(\sG),\bH)$$
with the following properties:

\begin{enumerate}

\item Each connected component of $\wt\bMaps(\Rep(\sG),\bH)$ is quasi-compact and irreducible; 

\item The map $\pi$ is schematic and proper on every connected component of $\wt\bMaps(\Rep(\sG),\bH)$; 

\smallskip

\item The map $\pi$ is surjective on geometric points; 

\item The set of connected components of $\wt\bMaps(\Rep(\sG),\bH)$ splits as a union of finite clusters,
and elements from different clusters have non-intersecting images in $\bMaps(\Rep(\sG),\bH)$. 

\end{enumerate}

\medskip

It is clear that an existence of such a pair $(\wt\bMaps(\Rep(\sG),\bH),\pi)$ would imply the required properties of
$^{\on{red}}\bMaps(\Rep(\sG),\bH)^{\on{rigid}}$. 

\medskip

Properties (1) and (2) will be established in \secref{ss:uniform is proper}; Property (3) in \secref{ss:uniform surj},
and Property (4) in \secref{ss:irred comps}.

\medskip

We will now proceed to the construction of $\wt\bMaps(\Rep(\sG),\bH)$. 

\sssec{}

Let $\on{Par}(\sG)$ be the (po)set of standard parabolics in $\sG$. 
For every $\sP\in \on{Par}(\sG)$, let $\sM$ denote its Levi quotient. Note that by convention, for $\sP=\sG$,
the corresponding Levi quotient is $\sG_{\on{red}}$, the quotient of $\sG$ by its unipotent radical. 

\medskip

The maps
$$\sG \leftarrow \sP\to \sM$$ 
induce the maps
$$\bMaps(\Rep(\sG),\bH)\overset{\sfp_\sP}\leftarrow \bMaps(\Rep(\sP),\bH)\overset{\sfq_\sP}\to \bMaps(\Rep(\sM),\bH).$$

\sssec{}

Let $\sfe'$ be an algebraically closed field containing $\sfe$. Let us call an $\sfe'$-point of $\bMaps(\Rep(\sG),\bH)$ \emph{irreducible} 
if it does not factor through the above map $\sfq_\sP$ for any \emph{proper} parabolic $\sP\subsetneq \sG$.

\sssec{}

Let $\sigma_\sM$ be an \emph{irreducible} $\sfe$-point of $\bMaps(\Rep(\sM),\bH)$. Choose its lift to an $\sfe$-point of 
the (ind)-scheme $\bMaps(\Rep(\sM),\bH)^{\on{rigid}}$. Let $\on{Stab}_\sM(\sigma_\sM)\subset \sM$ be its 
stabilizer with respect to the $\sM$-action on $\bMaps(\Rep(\sM),\bH)^{\on{rigid}}$.  (Note that the subgroup 
$\on{Stab}_\sM(\sigma_\sM)$ depends on the choice of a lift, and a change of the choice by $m\in \sM(\sfe)$
results in conjugating $\on{Stab}_\sM(\sigma_\sM)$ by $m$.) 

\medskip

We obtain a locally closed embedding
\begin{equation} \label{e:embed irred}
\on{pt}/\on{Stab}_\sM(\sigma_\sM)\hookrightarrow \bMaps(\Rep(\sM),\bH).
\end{equation} 

We claim, however: 

\begin{prop} \label{p:embed irred}
The map \eqref{e:embed irred} is a \emph{closed} embedding.
\end{prop}

The proof will be given in \secref{sss:proof embed irred}.

\begin{rem}
As we will see in \corref{c:conn irred}, the map \eqref{e:embed irred} is actually the embedding
of a connected component \emph{at the reduced level}.
\end{rem}

\sssec{}

Denote
$$\bMaps(\Rep(\sP),\bH)_{\sigma_\sM}:=\bMaps(\Rep(\sP),\bH)\underset{\bMaps(\Rep(\sM),\bH)}\times \on{pt}/\on{Stab}_\sM(\sigma_\sM).$$

\medskip

Finally, we set
$$\wt\bMaps(\Rep(\sG),\bH):=\underset{\sP\in \on{Par}(\sG)}\sqcup\, \underset{\sigma_\sM}\sqcup\, \bMaps(\Rep(\sP),\bH)_{\sigma_\sM}.$$

The maps $\sfp_\sP$ define the sought-for map 
$$\pi: \wt\bMaps(\Rep(\sG),\bH)\to \bMaps(\Rep(\sG),\bH).$$

\begin{rem}

The prestacks $\bMaps(\Rep(\sP),\bH)$ and $\bMaps(\Rep(\sP),\bH)_{\sigma_\sM}$ have a very transparent meaning
in the main example of $\bH=\qLisse(X)$.

\medskip

In this case, 
$$\bMaps(\Rep(\sP),\bH)=:\LocSys^{\on{restr}}_{\sP}(X)$$
is the prestack classifying local systems with a reduction to $\sP$, and
$$\bMaps(\Rep(\sP),\bH)_{\sigma_\sM}=:\LocSys^{\on{restr}}_{\sP,\sigma_\sM}(X)$$
is the prestack of local systems with a reduction to $\sP$, whose induced $\sM$-local
system is isomorphic to $\sigma_\sM$. 

\medskip

So the properties of the resulting map
$$\pi:\wt\LocSys{}^{\on{restr}}_\sG(X)\to \LocSys^{\on{restr}}_\sG(X)$$
say that $\LocSys^{\on{restr}}_\sG(X)$ is uniformized by the disjoint union of the prestacks 
$\LocSys^{\on{restr}}_{\sP,\sigma_\sM}(X)$. 

\end{rem}

\ssec{Properness of the uniformization morphism} \label{ss:uniform is proper}

In this subsection we will show that the prestack $\wt\bMaps(\Rep(\sG),\bH)$ is an algebraic stack, each of whose connected components is 
quasi-compact and irreducible, and that the map $\pi$ is proper on each connected component. This will establish Properties (1) and (2)
from \secref{sss:properties of uniformization}. 

\sssec{} \label{sss:algebraic stack 1}

We will first show that each $\bMaps(\Rep(\sP),\bH)_{\sigma_\sM}$ is an algebraic stack, which is 
quasi-compact and irreducible. 

\medskip

For this, it is sufficient to show that the map 
$$\sfq_\sP:\bMaps(\Rep(\sP),\bH)\to \bMaps(\Rep(\sM),\bH)$$
is a relative algebraic stack (i.e., its base change by a derived affine
scheme yields an algebraic stack) with fibers that are 
quasi-compact and irreducible. 

\medskip

The property of a map between prestacks to be a relative algebraic stack with fibers that are quasi-compact and irreducible is stable under compositions. 
Filtering the unipotent radical of $\sP$ we reduce the assertion to the following:

\begin{prop}  \label{p:radical} 
Let 
\begin{equation} \label{e:SES groups}
1\to \on{Tot}(V) \to \sG_1\to \sG_2\to 1
\end{equation}
be a short exact sequence of algebraic groups, where $\on{Tot}(V)$ is the vector group
associated with a finite-dimensional $\sG_2$-representation $V$. Then the resulting map
$$\bMaps(\Rep(\sG_1),\bH)\to \bMaps(\Rep(\sG_2),\bH)$$
is a relative algebraic stack whose fibers are quasi-compact and irreducible.
\end{prop}

\sssec{} \label{sss:obstruction}

Before we prove \propref{p:radical}, we make the following observation. 

\medskip

First, the datum of \eqref{e:SES groups} is equivalent to that of an object
$$\on{cl}_{\sG_1}\in \Maps_{\Rep(\sG_2)}(\on{triv},V[2]).$$

\medskip

Let $\bA$ be a symmetric monoidal category,
and let us be given a symmetric monoidal functor
$$\sF:\Rep(\sG_2)\to \bA.$$

Consider the object $\sF(V)\in \bA$ and
$$\sF(\on{cl}_{\sG_1})\in \Maps(\one_\bA,\sF(V)[2]).$$

\begin{lem} \label{l:obstruction}
Under the above circumstances, the space of lifts of $\sF$ to a functor
$$\Rep(\sG_1)\to \bA$$ 
identifies with the space of null-homotopies of $\sF(\on{cl}_{\sG_1})$. 
\end{lem} 

\begin{proof}

To simplify the notation, we will assume that $V$ is the trivial 1-dimensional representation
of $\sG_2$. Then the datum of $\sG_1$ is equivalent to that of a map of prestacks 
$$s:B(\sG_2)\to B^2(\BG_a),$$
so that
$$B(\sG_1) \simeq B(\sG_2)\underset{B^2(\BG_a)}\times \on{pt}.$$

It follows that
\begin{equation} \label{e:QCoh on gerbe}
\QCoh(B(\sG_1)) \simeq \QCoh(B(\sG_2))\underset{\QCoh(B^2(\BG_a))}\otimes \Vect_\sfe.
\end{equation} 

Note that $$\QCoh(B^2(\BG_a))\simeq \sfe[\eta]\mod, \quad \deg(\eta)=2.$$
The pullback of $\eta$, viewed as a point in
$$\Maps_{\QCoh(B^2(\BG_a))}(\CO_{B^2(\BG_a)},\CO_{B^2(\BG_a)}[2])$$
by means of $s$ is our 
$$\on{cl}_{\sG_1}\in \Maps_{\QCoh(B(\sG_2))}(\CO_{B(\sG_2)},\CO_{B(\sG_2)}[2]).$$

\medskip

Note that for a symmetric monoidal category $\bA'$, and $A\in \on{ComAlg}(\Vect_\sfe)$, 
the space of symmetric monoidal functors
$$A\mod\to \bA'$$
is isomorphic to the space of maps in $\on{ComAlg}(\Vect_\sfe)$ 
$$A \to \Maps_{\bA'}(\one_{\bA'},\one_{\bA'}).$$

For $A=\sfe[\eta]$, this space is further isomorphic to 
$$\Maps_{\bA'}(\one_{\bA'},\one_{\bA'}[2]).$$

\medskip

Combining this with \eqref{e:QCoh on gerbe}, we obtain that
$$\Rep(\sG_1)\simeq \Rep(\sG_2)\underset{\sfe[\eta]\mod}\otimes \Vect_\sfe,$$
where $\sfe[\eta]\mod\to \Rep(\sG_2)$ is given by
$$\eta \mapsto \on{cl}_{\sG_1}\in \Maps_{\Rep(\sG_2)}(\one_{\Rep(\sG_2)},\one_{\Rep(\sG_2)}[2]).$$

I.e., symmetric monoidal functors $\Rep(\sG_1)\to \bA$ are the same as symmetric monoidal functors
$\sF:\Rep(\sG_2)\to \bA$, equipped with the homotopy of the induced symmetric functor
$$\sfe[\eta]\mod \to \bA, \quad \eta\mapsto \sF(\on{cl}_{\sG_1})\in \Maps_\bA(\one_\bA,\one_\bA[2])$$
with
$$\sfe[\eta]\mod \to \Vect_\sfe \overset{\sfe\mapsto \one_\bA}\longrightarrow \bA,$$
as required. 

\end{proof}

\sssec{Proof of \propref{p:radical}} \label{sss:algebraic stack 2}

Let us be given an affine scheme $S$ and an $S$-point 
$$\sF:\Rep(\sG_2)\to \QCoh(S)\otimes \bH$$
of $\bMaps(\Rep(\sG_2),\bH)$. Consider the object 
$$\sF(V)\in \QCoh(S)\otimes \bH,$$
and the object
$$\CE:=(\on{Id}_{\QCoh(S)}\otimes \CHom_\bH(\one_\bH,-))(\sF(V))\in \QCoh(S).$$

According to \lemref{l:obstruction}, we have a point
$$\sF(\on{cl}_{\sG_1})\in \Gamma(S,\CE[2]),$$
and the fiber product 
$$S\underset{\bMaps(\Rep(\sG_2),\bH)}\times \bMaps(\Rep(\sG_1),\bH)$$
is the functor that sends $S'\to S$ to the space of null-homotopies
of $\sF(\on{cl}_{\sG_1})|_{S'}$. 

\medskip

Hence, it remains to show that the above functor of null-homotopies
is indeed an algebraic stack over $S$ with fibers that are quasi-compact and irreducible. For that
it suffices to show that, locally on $S$, the object $\CE$ is perfect of amplitude $[0,d]$
for some $d$, i.e., can be represented by a finite complex
$$\CE_0 \to \CE_1\to ...\to \CE_d,$$
where each $\CE_i$ is locally free of finite rank. 

\sssec{}

Note that
$$\QCoh(S)^{\on{perf}}\simeq \underset{n}{\on{lim}}\, \QCoh({}^{\leq n}\!S)^{\on{perf}}.$$

Hence, we can assume that $S$ is eventually coconnective. By \propref{p:replace by Ind-Lisse},
we can view $\sF$ as a functor with values in $\QCoh(S)\otimes \bH^{\on{access}}$. Moreover, 
$\sF(V)$ is compact in $\QCoh(S)\otimes \bH^{\on{access}}$ (see the proof of \corref{c:def}(b')).

\medskip

The functor $\CHom_{\bH^{\on{access}}}(\one_\bH,-)$
preserves compactness, hence, the object $\CE$ belongs to 
$$\QCoh(S)^c= \QCoh(S)^{\on{perf}},$$ 
and we only need to estimate its cohomological amplitude. 

\medskip

It is easy to see that if $\CE\in \QCoh(S)^{\on{perf}}$ is such that $\CE|_{^{\on{cl}}\!S}$ has
amplitude $[d_1,d_2]$, then so does $\CE$ itself. Hence, we can assume that $S$
is classical. Furthermore, since the prestacks
involved are locally (almost) of finite type, we can assume that $S$ is of finite type. 

\medskip

In this case, if $\CE\in \QCoh(S)^{\on{perf}}$ is such that its *-fiber at any $\sfe$-point
of $S$ lives in degrees $[d_1,d_2]$, then $\CE$ has amplitude $[d_1,d_2]$. Hence, we
have reduced the assertion to the case when $S=\on{pt}=\Spec(\sfe)$. 

\medskip

Now, the required property follows from the fact that for
$$\sF(V)=:\bh\in \bH^\heartsuit,$$
we have
$$\CHom_\bH(\one_\bH,\bh)\in (\Vect_\sfe)^{\geq 0,\leq d}$$
for some $d$ (i.e., $d$ is the cohomological amplitude of the functor $\CHom_\bH(\one_\bH,-)$, which 
is finite by the assumption on $\bH$).

\qed[\propref{p:radical}]

\sssec{}

We will now show that $\pi$ is schematic and proper when restricted to every connected component of $\wt\bMaps(\Rep(\sG),\bH)$.
Given \propref{p:embed irred} (which will be proved independently), 
it is sufficient to show that at the level of the underlying reduced prestacks, the map
$$\sfp_\sP:\bMaps(\Rep(\sP),\bH)\to \bMaps(\Rep(\sG),\bH)$$
is schematic, quasi-compact and proper.

\medskip

However, using the fact that $\sG/\sP$ is proper, the assertion follows from the next observation:

\begin{prop} \label{p:repr of fibers}
Let $\sG'$ be a subgroup of $\sG$. Then the map
$$\bMaps(\Rep(\sG'),\bH)\to \bMaps(\Rep(\sG),\bH)\underset{\on{pt}/\sG}\times \on{pt}/\sG',$$
given by $\oblv_\bH$, is a closed embedding. 
\end{prop} 

\begin{proof}

The statement is equivalent to the assertion that
$$\bMaps(\Rep(\sG'),\bH)^{\on{rigid}}\to \bMaps(\Rep(\sG),\bH)^{\on{rigid}}$$
is a closed embedding. 

\medskip

By \propref{p:Tannaka red cl}, it suffices to show that for an algebraic group $\sH$, the map
$$\bMaps_{\on{Grp}}(\sH,\sG')\to \bMaps_{\on{Grp}}(\sH,\sG)$$
is a closed embedding.

\medskip

However, this is the content of Remark \ref{r:closed emb maps}.

\end{proof}

\ssec{Surjectivity of the uniformization morphism} \label{ss:uniform surj}

In this subsection we will prove that $\pi$ is surjective, i.e., Property (3) from \secref{sss:properties of uniformization}. 

\sssec{}

We need to show that for an algebraically closed field $\sfe'$ containing $\sfe$, any $\sfe'$-point $\sigma_\sG$ of 
$\bMaps(\Rep(\sG),\bH)$ equals the image of an $\sfe'$-point of $\wt\bMaps(\Rep(\sG),\bH)$.

%

\medskip

We will argue by induction on the semi-simple rank of $\sG$. If $\sigma_\sG$ equals the image of a $\sfe'$-point
$\sigma_\sP$ of $\bMaps(\Rep(\sP),\bH)$ for a \emph{proper} parabolic $\sP\subsetneq \sG$, 
we are done by the induction hypothesis. 

\medskip

Hence, we can assume that $\sigma_\sG$ is irreducible, and we need to show the following:

\begin{prop} \label{p:dom irred}
Let $\sG$ be reductive. Then for an algebraically closed field extension $\sfe'\supseteq \sfe$, 
any irreducible $\sfe'$-point of $\bMaps(\Rep(\sG),\bH)$ factors through an $\sfe$-point. 
\end{prop}

Note that as a particular case, we obtain: 

\begin{cor}  \label{c:irred}
For $\sG=\sT$ being a torus, the prestack $^{\on{red}}\bMaps(\Rep(\sT),\bH)$
is the disjoint union of copies of $\on{pt}/\sT$ over the set of isomorphism classes of 
$\sfe$-points of $\bMaps(\Rep(\sT),\bH)$.
\end{cor} 

\sssec{} We proceed with the proof of \propref{p:dom irred}. By \propref{p:Tannaka red cl}, 
the statement of \propref{p:dom irred} is equivalent to the following: 

\begin{prop} \label{p:group irred}
Let $\sG,\sH$ be algebraic groups with $\sG$ reductive. Let $\sfe'\supseteq \sfe$ be an 
algebraically closed field extension, and let $\phi:\sH\to \sG$ be a homomorphism 
defined over $\sfe'$. Assume that the image of $\phi$ is not contained in any proper
parabolic of $\sG$ defined over $\sfe'$. Then $\phi$ is $\sG$-conjugate to a homomorphism defined over $\sfe$.
\end{prop} 

\ssec{Irreducible homomorphisms of reductive groups}

This subsection is devoted to the proof of \propref{p:group irred}. 

\sssec{}

Consider the Levi decomposition of $\sH$
$$1\to \sH_u\to \sH\to \sH_{\on{red}}\to 1.$$

We claim that $\phi$ factors via a homomorphism
$$\sH_{\on{red}}\to \sG.$$

Let $\sH' = \sH/\on{ker}(\phi)$.  Thus, we have an injective homomorphism
$$ \sH' \to \sG .$$
We need to show that $\sH'$ is reductive.
We now recall the following assertion from \cite[Proposition 4.2]{Se}: 

\begin{thm} \label{t:Serre}
For a connected reductive group $\sG$ and a subgroup $\sH\subset \sG$ the following conditions are equivalent:

\smallskip

\noindent{\em(i)} $\sH$ is reductive;

\smallskip

\noindent{\em(ii)} Whenever there exists a parabolic $\sP\subset \sG$ that contains $\sH$, there also exists a Levi splitting 
$\sP\leftrightarrows \sM$ such that $\sH\subset \sM$.

\end{thm}

By the irreducibility assumption, our subgroup $\sH'$ satisfies (ii) in \thmref{t:Serre}. Hence, it is reductive as claimed.

\sssec{}

Thus, in order to prove \propref{p:group irred}, it suffices to establish the following:

\begin{prop} \label{p:red group hom}
Let $\sH$ and $\sG$ be a pair of algebraic groups with $\sH$ reductive. Then the ind-scheme $\bMaps_{\on{Grp}}(\sH,\sG)$
is the disjoint union over isomorphism classes of homomorphisms
$$\phi:\sH\to \sG$$
of the (classical) schemes $\sG/\on{Stab}_\sG(\phi)$, where the stabilizer is taken with respect to the action of $\sG$ on
$\bMaps_{\on{Grp}}(\sH,\sG)$ by conjugation. 
\end{prop}

This proposition is well-known. We will supply a proof for completeness. 

\begin{proof}

It is enough to show that for any $\sfe$-point of $\bMaps_{\on{Grp}}(\sH,\sG)$, the resulting locally closed embedding
$$\sG/\on{Stab}_\sG(\phi)\to \bMaps_{\on{Grp}}(\sH,\sG)$$
induces an isomorphism at the level of tangent spaces. 

\medskip

Let our $\sfe$-point correspond to a homomorphism $\phi:\sH\to \sG$. Thus, we have to show that
the map
\begin{equation} \label{e:orbit}
\on{coFib}(H^0(\inv_{\sH}(\sg)) \to \sg) \to T_\phi(\bMaps_{\on{Grp}}(\sH,\sG))
\end{equation} 
is an isomorphism, where $\inv_{\sH}$ stands for $\sH$-invariants, and $\sg$ is 
viewed as a $\sH$-representation via $\phi$ and the adjoint action. 

\medskip

Comparing with the formula for $T^*_\phi(\bMaps_{\on{Grp}}(\sH,\sG))$ in Remark \ref{r:cotan Hom}, 
we obtain that we need to show that 
$$\inv_{\sH}(\sg)$$
is concentrated in cohomological degree $0$. 

\medskip

However, this follows from the assumption that $\sH$ is reductive (and hence the category $\Rep(\sH)$ is semi-simple). 

\medskip

\noindent NB: note that validity of \propref{p:red group hom} depends on the assumption 
that we work over a field of characteristic $0$ (in our case this is the field of coefficients $\sfe$).

\end{proof} 

\sssec{Proof of \propref{p:embed irred}} \label{sss:proof embed irred}

It suffices to show that if $\sG$ is reductive and $\phi:\sH\to \sG$ is a homomorphism of algebraic groups that does not
factor through a parabolic, then the $\on{Ad}(\sG)$-orbit of $\phi$ is closed in $\bMaps(\sH,\sG)$.

\medskip

Set $\sH':=\sH/\on{ker}(\phi)$, so that $\phi$ factors through a map $\phi':\sH'\to \sG$. It is enough to
check that the $\on{Ad}(\sG)$-orbit of $\phi'$ is closed in $\bMaps(\sH',\sG)$. We will show that it
is in fact a connected component.

\medskip

Indeed, by \thmref{t:Serre}, $\sH'$ is reductive, and the assertion follows from \propref{p:red group hom}. 

\qed[\propref{p:embed irred}]

\ssec{Associated pairs and semi-simple local systems}  \label{ss:assoc}

In this subsection we will make a digression and discuss the classification of \emph{semi-simple}  
points of $\bMaps(\Rep(\sG),\bH)$ in terms of irreducible ones of $\bMaps(\Rep(\sM),\bH)$,
where $\sM$ is a Levi subgroups of $\sG$.

\sssec{Terminology} \label{sss:term loc sys}

In what follows, for an algebraic group $\sG'$, 
it will be convenient to refer to $\sfe$-points of $\bMaps(\Rep(\sG'),\bH)$ as ``$\sG'$-local systems". 
They are literally such in the key example $\bH=\qLisse(X)$, in which case 
$$\bMaps(\Rep(\sG'),\bH)=\LocSys^{\on{restr}}_{\sG'}(X).$$

\medskip

For a homomorphism $\phi:\sG'_1\to \sG'_2$, and a $\sG'_1$-local system 
$$\sigma_{\sG'_1}\in \bMaps(\Rep(\sG'_1),\bH)(\sfe),$$
we will refer to its image 
$$\sigma_{\sG'_2}\in \bMaps(\Rep(\sG'_2),\bH)(\sfe)$$
as ``the $\sG'_2$-local system induced by $\sigma_{\sG'_1}$ by means of $\phi$".

\medskip

Vice versa, given $\sigma_{\sG'_2}$, we will refer to $\sigma_{\sG'_1}$ as a ``reduction of $\sigma_{\sG'_2}$
to $\sG'_1$". 

\sssec{}  \label{sss:assoc}

Let $\sP$ be a standard parabolic in $\sG$. We will view the partial flag variety $\sG/\sP$ as the space of 
parabolics conjugate to $\sP$.

\medskip

Let $\sP_1$ and $\sP_2$ be a pair of standard parabolics in $\sG$, each equipped with an irreducible local system, 
$\sigma_{\sM_i}$ with respect to the corresponding Levi quotient $\sM_i$. 

\medskip

We shall say that the pairs $(\sP_1,\sigma_{\sM_1})$ and $(\sP_2,\sigma_{\sM_2})$ are \emph{associated}
if there exists a $\sG$-orbit $\sO$ in $\sG/\sP_1\times \sG/\sP_2$, such that for some/any pair of points 
$(\sP'_1,\sP'_2)\in \sO$ the following holds: 

\begin{itemize}

\item The maps
$$\sM_1\leftarrow \sP'_1 \leftarrow \sP'_1\cap \sP'_2 \to \sP'_2\to \sM_2$$
identify $\sM_i$, $i=1,2$, with the Levi quotient of $\sP'_1\cap \sP'_2$;

\smallskip

\item Under the resulting isomorphism $\sM_1\simeq \sM_2$, the local systems $\sigma_{\sM_1}$
and $\sigma_{\sM_2}$ are isomorphic.

\end{itemize}

\sssec{}

We claim:

\begin{lem} \label{l:assoc}
Two pairs $(\sP_1,\sigma_{\sM_1})$ and $(\sP_2,\sigma_{\sM_2})$ are associated if and only if
there exists a $\sG$-local system $\sigma_\sG$, equipped with reductions to both $\sP_1$ and $\sP_2$,
such that the induced $\sM_i$-local systems are $\sigma_{\sM_i}$, respectively.
\end{lem}

\begin{proof}

One direction is clear: if $(\sP_1,\sigma_{\sM_1})$ and $(\sP_2,\sigma_{\sM_2})$ are associated, 
choose a pair $(\sP'_1,\sP'_2)$ on the corresponding orbit, and choose a Levi splitting 
of $\sP'_1\cap \sP'_2$. Then the resulting local system with respect to $\sP'_1\cap \sP'_2$
projects to $\sigma_{\sM_1}$ and $\sigma_{\sM_2}$, respectively.

\medskip

For the other implication, the two reductions of $\sigma_\sG$ correspond to a $\sG$-orbit $\sO$ on $\sG/\sP_1\times \sG/\sP_2$.
We will show that this orbit satisfies the two conditions of \secref{sss:assoc}. 

\medskip

By assumption, we can choose parabolics $\sP'_1$ and $\sP'_2$, conjugate to $\sP_1$ and $\sP_2$, respectively,
and lying on $\sO$, so that $\sigma_\sG$ admits a reduction to $\sP'_1\cap \sP'_2$; denote this reduction by $\sigma_{1,2}$. 
Furthermore, $\sigma_{\sM_i}$, $i=1,2$, is induced from $\sigma_{1,2}$ along the map 
\begin{equation} \label{e:P12}
\sP'_1\cap \sP'_2\hookrightarrow \sP'_i\twoheadrightarrow \sM_i.
\end{equation}

We note that (for any pair of parabolics) the image of \eqref{e:P12}
is a parabolic subgroup in $\sM_i$, $i=1,2$. 

\medskip

Hence, 
by the assumption on $\sigma_{\sM_i}$, the maps \eqref{e:P12} are surjective, and hence identify
$\sM_i$ as a Levi quotient of $\sP'_1\cap \sP'_2$.

\end{proof} 

\sssec{}  \label{sss:semi-simple}

We will say that a $\sG$-local system $\sigma_\sG$ is \emph{semi-simple} if whenever it admits a reduction to 
a local system $\sigma_\sP$ with $\sP\subset \sG$ a parabolic, then $\sigma_\sP$ admits a further reduction 
to a local system $\sigma_\sM$ for $\sM$ for \emph{some} Levi splitting 
$$\sP\leftrightarrows \sM.$$

\sssec{}

By \thmref{t:Serre}, when we think of $\sigma_\sG$ as a conjugacy class of homomorphisms 
$$\phi:\sH\to \sG$$
(for $\sH$ as in \secref{sss:Galois}), semi-simplicity is equivalent to the condition that the image of $\phi$ 
be reductive. 

\medskip

The latter interpretation has the following consequence:

\begin{lem} \label{l:semi-simple}
Let $\sG'\hookrightarrow \sG$ be an injection of algebraic groups. Then a 
$\sG'$-local system is semi-simple if and only if the induced 
$\sG$-local system is semi-simple.
\end{lem} 

\sssec{}

Let $\sP\subset \sG$ be a standard parabolic and choose a Levi splitting $\sP\leftrightarrows \sM$. 
Let $\sigma_\sM$ be an irreducible $\sM$-local system, and let $\sigma_\sG$ denote the induced $\sG$-local system 
via $\sM\to \sP\to \sG$. From \lemref{l:semi-simple} we obtain that $\sigma_\sG$ is semi-simple. 

\medskip

From \lemref{l:assoc}, we obtain: 

\begin{cor} \label{c:assoc vs semisimple} 
For two pairs $(\sP_1,\sigma_{\sM_1})$ and $(\sP_2,\sigma_{\sM_2})$, the $\sG$-local systems 
$\sigma_{\sG,1}$ and $\sigma_{\sG,2}$ are isomorphic if and only if $(\sP_1,\sigma_{\sM_1})$ and $(\sP_2,\sigma_{\sM_2})$ are associated.
\end{cor}  

And further: 

\begin{cor} \label{c:assoc vs semisimple bis} \hfill

\smallskip

\noindent{\em(a)} Association is an equivalence relation on the set of isomorphism classes of pairs 
$(\sP,\sigma_\sM)$, where $\sP$ is a parabolic, $\sM$ is its Levi quotient, and $\sigma_\sM$
is an irreducible local system with respect to $\sM$.

\smallskip

\noindent{\em(b)} The assignment $(\sP,\sigma_\sM)\mapsto \sigma_\sG$ 
establishes a bijection between classes of association
of pairs $(\sP,\sigma_\sM)$ and isomorphism classes of semi-simple $\sG$-local systems.

\end{cor} 

\sssec{}

For future use we notice:

\begin{lem} \label{l:assoc finite}
Each equivalence class of associated pairs $(\sP,\sigma_\sM)$ contains only
finitely many elements.
\end{lem}

\begin{proof}
Follows from the fact that for every pair of standard parabolics $\sP_1$ and $\sP_2$, 
there are finitely many $\sG$-orbits on $\sG/\sP_1\times \sG/\sP_2$.
\end{proof} 

\sssec{}

Given two $\sG$-local systems $\sigma_1$ and $\sigma_2$, we shall say that $\sigma_2$ is a \emph{semi-simplification} of $\sigma_1$ if 

\begin{itemize}

\item $\sigma_2$ is semi-simple;

\item there exists a parabolic $\sP$ and reductions $\sigma_{1,\sP}$ and $\sigma_{2,\sP}$
of $\sigma_1$ and $\sigma_2$, respectively, such that the induced local systems with respect to the Levi quotient of $\sP$
are isomorphic. 

\end{itemize} 

It is clear that every local system $\sigma_\sG$ admits a semi-simplification: take the minimal standard
parabolic $\sP$ to which $\sigma_\sG$ can be reduced, and let $\sigma'_\sG$ be the $\sG$-local system induced 
from the reduction of $\sigma_\sG$ to $\sP$ via the homomorphism
$$\sP\to \sM\to \sP\to \sG$$
for some Levi splitting of $\sM$.

\medskip 

From \lemref{l:assoc} and \corref{c:assoc vs semisimple bis}(b) we obtain:

\begin{cor} \label{c:semi-simplification}
For a given local system, its semi-simplification is well-defined up to isomorphism.
\end{cor} 

\ssec{Analysis of connected/irreducible components}  \label{ss:irred comps}

In this subsection, we will establish Property (4) of the map
$$\pi:\wt\bMaps(\Rep(\sG),\bH)\to \bMaps(\Rep(\sG),\bH).$$
from \secref{sss:properties of uniformization}, namely, 
that the set of connected components of $\wt\bMaps(\Rep(\sG),\bH)$ is a union of finite clusters,
and elements from different clusters have non-intersecting images in $\bMaps(\Rep(\sG),\bH)$ along $\pi$. 

\medskip

In addition, we will describe explicitly the set of connected components of $\bMaps(\Rep(\sG),\bH)$. 

\sssec{}

We will prove:

\begin{prop} \label{p:conn comps LocSys}
There exists a bijection between the set of connected components of the prestack $\bMaps(\Rep(\sG),\bH)$
and the set of isomorphisms classes of semi-simple $\sG$-local systems, characterized by either
of the following two properties: 

\medskip

\noindent{\em(a)} 
Two $\sfe$-points of $\bMaps(\Rep(\sG),\bH)$ belong to the same connected component 
if and only if they have isomorphic semi-simplifications. 

\medskip

\noindent{\em(b)} For a standard parabolic $\sP$ and an irreducible $\sM$-local system $\sigma_\sM$, the 
map $$\bMaps(\Rep(\sP),\bH)_{\sigma_\sM}\to \bMaps(\Rep(\sG),\bH)$$ lands
in the connected component corresponding via the bijection of \corref{c:assoc vs semisimple bis}(b) to the
association class of $(\sP,\sigma_\sM)$. 

\end{prop}

Note that the assertion of \propref{p:conn comps LocSys} combined with that of \lemref{l:assoc finite} implies
Property (4) from \secref{sss:properties of uniformization}, which was the last
one remaining to establish.

\sssec{}

We note that point (a) of \propref{p:conn comps LocSys} contains the following statement:

\begin{cor} \label{c:conn comps LocSys}
Each connected component of $\bMaps(\Rep(\sG),\bH)$ contains a unique isomorphism class of $\sfe$-points corresponding 
to a semi-simple $\sG$-local system.
\end{cor} 

Note also the following consequence of \propref{p:conn comps LocSys}:

\begin{cor}   \label{c:conn irred}
Let $\sG$ be reductive. Then a connected component of $\bMaps(\Rep(\sG),\bH)$ containing an irreducible local system
has a unique isomorphism class of $\sfe$-valued points.
\end{cor} 

\sssec{Proof of \propref{p:conn comps LocSys}}

For a given class of association $\fP$ of pairs $(\sP,\sigma_\sM)$, let
$$\bMaps(\Rep(\sG),\bH)_\fP$$
be the union of images of the maps
$$\pi:\bMaps(\Rep(\sP),\bH)_{\sigma_\sM}\to \bMaps(\Rep(\sG),\bH)$$
from $\fP$. This union is finite by \lemref{l:assoc finite}. Since the maps $\pi$
are proper, this is a closed substack of $\bMaps(\Rep(\sG),\bH)$.

\medskip

We will show that the substacks $\bMaps(\Rep(\sG),\bH)_\fP$ are:

\medskip

\noindent(i) Pairwise disjoint;

\medskip

\noindent(ii) Connected;

\medskip

\noindent(iii) The semi-simplification of every $\sfe$-point of $\bMaps(\Rep(\sG),\bH)_\fP$
corresponds under the bijection of \corref{c:assoc vs semisimple bis}(b) to $\fP$. 

\medskip

Point (i) follows readily from \lemref{l:assoc}. 

\medskip

Let $\sigma_\sG$ be the semi-simple $\sG$-local system corresponding to $\fP$, and consider the corresponding map
$$\pi:\bMaps(\Rep(\sP),\bH)_{\sigma_\sM}\to \bMaps(\Rep(\sG),\bH)$$
for $(\sP,\sigma_\sM)\in \fP$. 

\medskip

By definition, all $\sfe$-points in the image of this map have $\sigma_\sG$ as their semi-simplification. This proves 
point (iii).

\medskip

Further, $\sigma_\sG$ itself is contained in the image of the above map $\pi$. This proves point (ii),
since each $\bMaps(\Rep(\sP),\bH)_{\sigma_\sM}$ is irreducible (and hence connected), and their
images for $(\sP,\sigma_\sM)\in \fP$ all intersect at $\sigma_\sG$. 

\qed[\propref{p:conn comps LocSys}]

\begin{rem}
It is easy to see from the above argument that for given a local system $\sigma$, the map
$$\on{pt}/\on{Stab}_{\sG}(\sigma)\to \bMaps(\Rep(\sG),\bH)$$
is a closed embedding if $\sigma$ is semi-simple (cf. Propositions \ref{p:irred closed} and \ref{p:irred closed Betti}).
Furthermore, if $\sG$ is reductive, then the above assertion is ``if and only if". 

\medskip

Indeed, for every $(\sP,\sigma_\sM)$, it is clear that the $\sP$-local system $\sigma_P^0$ induced from $\sigma_\sM$
via a Levi splitting $\sM\to \sP$ is a closed point in 
$$\bMaps(\Rep(\sP),\bH)\underset{\bMaps(\Rep(\sM),\bH)}\times \on{pt},$$
and the assertion follows from the fact that $\pi$ is proper. 

\medskip

Further, if $\sG$ is reductive, the action of the center $\sM$ contracts any $\sfe$-point of the above fiber product 
to $\sigma_P^0$. 

\end{rem}

\begin{rem}
It is clear that the image of each $\bMaps(\Rep(\sP),\bH)_{\sigma_\sM}$ along $\pi$ is irreducible.
However, it is \emph{not} true that the images of different of $\bMaps(\Rep(\sP_1),\bH)_{\sigma_{\sM_1}}$
and $\bMaps(\Rep(\sP_2),\bH)_{\sigma_{\sM_2}}$ 
in $\bMaps(\Rep(\sG),\bH)$
will always produce different irreducible components:

\medskip

For example, take $G=GL_2$, $\sP_1=\sP_2=\sB$, so $\sM_1=\sM_2=\BG_m\times \BG_m$. Take
$\bH=\qLisse(X)$ and let $\sigma_{\sM_1}$ and $\sigma_{\sM_2}$ be given by
$$(E_1,E_2) \text{ and } (E_2,E_1),$$
where $E_1$ and $E_2$ are non-isomorphic one-dimensional local systems.

\medskip

Then if $X$ is a curve of genus $\geq 2$, the images of $\LocSys^{\on{restr}}_{\sP,\sigma_{\sM_1}}(X)$
and $\LocSys^{\on{restr}}_{\sP,\sigma_{\sM_2}}(X)$ are two distinct irreducible components of 
$\LocSys^{\on{restr}}_\sG(X)$.

\medskip

By contrast, if $X$ is a curve of genus $1$, both $\LocSys^{\on{restr}}_{\sP,\sigma_{\sM_1}}(X)$
and $\LocSys^{\on{restr}}_{\sP,\sigma_{\sM_2}}(X)$ are set-theoretically isomorphic to $\on{pt}/(\BG_m\times \BG_m)$,
and they map onto the same closed subset of $\LocSys^{\on{restr}}_\sG(X)$.

\end{rem}

\section{Comparison with the Betti and de Rham versions of $\LocSys^\dr_\sG(X)$}  \label{s:Betti and dR}

In this section we study the relationship between $\LocSys^{\on{restr}}_\sG(X)$ with 
$\LocSys^\dr_\sG(X)$ in the two contexts when the latter is defined: de Rham and Betti.

\medskip

We will show that in both cases, the map
$$\LocSys^{\on{restr}}_\sG(X)\to \LocSys^\dr_\sG(X)$$
is a \emph{formal isomorphism} with an explicit image at the reduced level. 

\ssec{Relation to the Rham version}  \label{ss:dR}

In this subsection we will take our ground field $k$ to be of characteristic $0$.
We will take $\sfe=k$ and let $\Shv(-)$ to be the sheaf theory of ind-holonomic D-modules. 

\medskip

We will study the relationship between $\LocSys^{\on{restr}}_\sG(X)$
and the ``usual" stack $\LocSys^\dr_\sG(X)$ classifying de Rham local systems. 

\sssec{}

Let $X$ be a scheme of finite type over $k$. Recall (see, e.g., \cite[Sects. 10.1-2]{AG}), that the 
prestack of de Rham local systems on $X$, denoted $\LocSys^\dr_\sG(X)$, is defined by sending 
$S\in \affSch_{/\sfe}$ to the space of right t-exact symmetric monoidal functors
$$\Rep(\sG)\to \QCoh(S)\otimes \Dmod(X).$$

\medskip

It is shown in \cite[Proposition 10.3]{AG} that $\LocSys^\dr_\sG(X)$ is laft (=locally almost of finite type)
and admits $(-1)$-connective corepresentable deformation theory.

\begin{rem}

Note that the prestack $\LocSys^\dr_\sG(X)$ is of the form $\bMaps(\Rep(\sG),\bH)$, where 
$\bH$ is the symmetric monoidal category $\Dmod(X)$.

\end{rem}

\sssec{}

We have a tautologically defined symmetric monoidal functor
\begin{equation} \label{e:lisse into Dmod}
\qLisse(X) \hookrightarrow \Shv(X)\to \Dmod(X),
\end{equation} 
which gives rise to a map of prestacks 
\begin{equation} \label{e:restr to all dR}
\LocSys^{\on{restr}}_\sG(X)\to \LocSys^\dr_\sG(X).
\end{equation} 

We observe: 

\begin{lem}
The map \eqref{e:restr to all dR} is a monomorphism 
(i.e., is a monomorphism of spaces when evaluated on any affine scheme).
\end{lem}

\begin{proof} 

Note that objects of $\Shv(X)^{\on{constr}}$ are compact as objects of $\Dmod(X)$. Hence,
the functor
$$\Shv(X)\to \Dmod(X),$$
obtained by ind-extending the tautological embedding is fully faithful. 

\medskip

Therefore, so is the composite functor \eqref{e:lisse into Dmod}. Since $\QCoh(S)$ is dualizable, the functor 
\begin{equation} \label{e:LS in Dmod}
\QCoh(S)\otimes \qLisse(X)\to \QCoh(S)\otimes \Dmod(X)
\end{equation} 
is also fully faithful. This implies the assertion of the lemma.

\end{proof}

\sssec{Example}  \label{sss:Gm bis}

Let us explain how the difference between $\LocSys^\dr_{\BG_m}(X)$ and $\LocSys^{\on{restr}}_{\BG_m}(X)$ plays out in the 
simplest cases when $\sG=\BG_m$ and $\sG=\BG_a$ . 

\medskip

Take $S$ to be classical. Then $S$-points of $\LocSys^\dr_{\BG_m}(X)$ are line bundles over $S\times X$,
equipped with a connection along $X$. Trivializing this line bundle locally, the connection corresponds to a section of
$$\CO_S\boxtimes \Omega^{1,\on{cl}}_X,$$
i.e., a function on $S$ with values in closed 1-forms on $X$.

\medskip

By contrast, if our $S$-point lands in $\LocSys^{\on{restr}}_{\BG_m}(X)$, and 
if we further assume that $S$ is integral, by Example \ref{sss:Gm}, our line bundle with connection
is necessarily pulled back from $X$. 

\medskip

Let us now take $\sG=\BG_a$. Then it follows from \secref{sss:Ga} that the map 
$$\LocSys^{\on{restr}}_{\BG_a}(X)\to \LocSys^\dr_{\BG_a}(X)$$
is an isomorphism. 

\sssec{}

Recall that a map of prestacks $\CY_1\to \CY_2$ is said to be a \emph{formal isomorphism} if
$\CF_1$ identifies with its own formal completion inside $\CY_2$, i.e., if the map
$$\CY_1\to (\CY_1)_{\on{dR}}\underset{(\CY_2)_{\on{dR}}}\times \CY_2$$
is an isomorphism.

\sssec{}

We claim:

\begin{prop} \label{p:formal compl dr}
The map \eqref{e:restr to all dR} is a formal isomorphism, i.e., identifies $\LocSys^{\on{restr}}_\sG(X)$ with
its formal completion inside $\LocSys^\dr_\sG(X)$.
\end{prop}

\begin{proof}

We need to show that 
%
%
for $S\in \affSch_{/\sfe}$ and a map
\begin{equation} \label{e:S to LocSys}
S\to \LocSys^\dr_\sG(X),
\end{equation} 
such that the composite map
$$^{\on{red}}\!S\to S\to \LocSys^\dr_\sG(X)$$
factors through $\LocSys^{\on{restr}}_\sG(X)$, the initial map \eqref{e:S to LocSys} factors though $\LocSys^{\on{restr}}_\sG(X)$
as well.

\medskip

Since both $\LocSys^\dr_\sG(X)$ and $\LocSys^{\on{restr}}_\sG(X)$
are prestacks locally almost of finite type,
we can assume that $S$ is eventually coconnective
and almost of finite type.

\medskip

Thus, we need to show that given a functor
\begin{equation} \label{e:S to LocSys bis}
\sF:\Rep(\sG)\to \QCoh(S)\otimes \Dmod(X),
\end{equation} 
such that the composite functor
$$\Rep(\sG)\to \QCoh(S)\otimes \Dmod(X)\to \QCoh({}^{\on{red}}\!S)\otimes \Dmod(X),$$ 
lands in 
\begin{equation} \label{e:LS in Dmod again red}
\QCoh({}^{\on{red}}S)\otimes \qLisse(X)\subset \QCoh({}^{\on{red}}\!S)\otimes \Dmod(X),
\end{equation} 
the functor \eqref{e:S to LocSys bis} also lands in 
\begin{equation} \label{e:LS in Dmod again}
\QCoh(S)\otimes \qLisse(X)\subset \QCoh(S)\otimes \Dmod(X).
\end{equation} 

\medskip 

Since $S$ was assumed eventually coconnective, 
by \propref{p:replace by Ind-Lisse}, in \eqref{e:LS in Dmod again red} and 
\eqref{e:LS in Dmod again}, we can replace $\qLisse(X)$ by $\iLisse(X)$.

\medskip

Let $\iota$ denote the embedding $\iLisse(X)\hookrightarrow \Dmod(X)$.
It sends compacts to compacts, hence admits a continuous right adjoint, to be denoted 
$\iota^R$. 

\medskip

We need to show that the natural transformation
$$(\on{Id}\otimes \iota)\circ (\on{Id}\otimes \iota^R)\circ \sF \to \sF$$
is an isomorphism. 

\medskip

Let $f$ denote the embedding $^{\on{red}}\!S\to S$. We know that 
$$(f^*\otimes \on{Id})\circ (\on{Id}\otimes \iota)\circ (\on{Id}\otimes \iota^R)\circ \sF \simeq 
(\on{Id}\otimes \iota)\circ (\on{Id}\otimes \iota^R)\circ (f^*\otimes \on{Id})\circ \sF \to (f^*\otimes \on{Id})\circ  \sF$$
is an isomorphism. 

\medskip

This implies the assertion since for $S\in {}^{<\infty}\!\affSch_{\on{aft}/\sfe}$, the functor 
$$f^*\otimes \on{Id}:\QCoh(S)\otimes \bC\to \QCoh({}^{\on{red}}\!S)\otimes \bC$$
is conservative for any DG category $\bC$ (indeed, $\QCoh(S)$ is generated under finite
limits by the essential image of $f_*$). 

\end{proof} 

\sssec{}

From now on, until the end of this subsection we will assume that $X$ is proper. In this case by 
\cite[Sects. 10.3.8 and 10.4.3]{AG}, we know that $\LocSys^\dr_\sG(X)$ is an algebraic stack locally almost of finite type.

\medskip

We claim:

\begin{thm} \label{t:compare dR}
The map 
$$\LocSys^{\on{restr}}_\sG(X)\to \LocSys^\dr_\sG(X)$$
is a closed embedding \emph{at the reduced level} for \emph{each connected component of} $\LocSys^{\on{restr}}_\sG(X)$.
\end{thm} 

This theorem will be proved in \secref{ss:proof of compare dR}. In the course of the proof we will also describe the closed
substacks of $^{\on{red}}\!\LocSys^\dr_\sG(X)$ that arise as images of connected components of $^{\on{red}}\!\LocSys^{\on{restr}}_\sG(X)$. 

\medskip

Combined with \propref{p:formal compl dr}, we obtain:

\begin{cor}  \label{c:compare dR}
The subfunctor $$\LocSys^{\on{restr}}_\sG(X)\subset \LocSys^\dr_\sG(X)$$
is the disjoint union\footnote{Sheafified in the Zariski/\'etale topology.} 
of formal completions of a collection of pairwise
non-intersecting closed substacks of ${}^{\on{red}}\!\LocSys^\dr_\sG(X)$.
\end{cor} 

\begin{rem}
The closed substacks of ${}^{\on{red}}\!\LocSys^\dr_\sG(X)$ appearing in \corref{c:compare dR}
will be explicitly described in Remark \ref{r:image}.
\end{rem}

\begin{rem} \label{r:G' to G deR}
Let $\sG'\to \sG$ be a closed embedding. It is not difficult to show that the diagram
$$
\CD
\LocSys^{\on{restr}}_{\sG'}(X) @>>> \LocSys^\dr_{\sG'}(X) \\
@VVV @VVV \\
\LocSys^{\on{restr}}_\sG(X) @>>> \LocSys^\dr_\sG(X) 
\endCD
$$
is a fiber square.
\end{rem} 

\ssec{A digression: ind-closed embeddings}

\sssec{}  \label{sss:ind-closed}

Let us recall the notion of \emph{ind-closed embedding} of prestacks (see \cite[Sect. 2.7.2]{GR3}). 

\medskip

First, if $S$ is an affine scheme and $\CY$ is a prestack mapping to it, we shall say that this map is an \emph{ind-closed} embedding
if $\CY$ is an ind-scheme and for some/any presentation of $\CY$ as \eqref{e:ind}, the resulting maps
$$Y_i\to S$$
are closed embeddings.

\medskip

We shall say that a map of prestacks $\CY_1\to \CY_2$ is an ind-closed embedding if its base change by
an affine scheme yields an ind-closed embedding.

\begin{rem}
Let us emphasize the difference between ``ind-closed embedding" and ``closed embedding". For example,
the inclusion of the disjoint union of infinitely many copies of $\on{pt}$ onto $\BA^1$ is an 
ind-closed embedding but not a closed embedding. Similarly, the map
$$\on{Spf}(\sfe\qqart)\to \BA^1$$
is an ind-closed embedding but not a closed embedding.
\end{rem}

\sssec{}

From \corref{c:compare dR} we obtain:

\begin{cor}  \label{c:compare dR ind}
The map $\LocSys^{\on{restr}}_\sG(X)\to \LocSys^\dr_\sG(X)$ is an ind-closed embedding.
\end{cor}

\begin{rem}  \label{r:union formal compl}

Let $f:\CY_1\to \CY_2$ be a map, where $\CY_2$ is an algebraic stack, locally almost of finite type. 
Assume that $f$ is a \emph{formal isomorphism}. It is not difficult to see that the following conditions on $f$ are equivalent: 

\medskip

\noindent{(i)} It is an ind-closed embedding;

\smallskip

\noindent{(ii)} $^{\on{red}}\CY_1$ is a union of 
closed subfunctors of $^{\on{red}}\CY_2$;

\smallskip

\noindent{(iii)} $\CY_1$ is obtained as the completion of $\CY_2$ along a 
subfunctor of $\CZ\subset {}^{\on{red}}\CY_2$ equal to a union of 
closed subfunctors.

\end{rem}

%


\ssec{Uniformization and the proof of \thmref{t:compare dR}}  \label{ss:proof of compare dR}

\sssec{}

For a standard parabolic $\sP$ consider the diagram
$$\LocSys^\dr_{\sG}(X)\overset{\sfp_\sP}\leftarrow \LocSys^\dr_{\sP}(X)\overset{\sfq_\sP}\to \LocSys^\dr_{\sM}(X).$$

Fix an irreducible local system $\sigma_\sM$ for $\sM$ and denote
$$\LocSys^\dr_{\sP,\sigma_\sM}(X):=\LocSys^\dr_{\sP}(X)\underset{\LocSys^\dr_{\sM}(X)}\times \on{pt}/\on{Stab}_\sM(\sigma_\sM).$$

\sssec{}

We have a commutative diagram
$$
\CD
\LocSys^{\on{restr}}_{\sP,\sigma_\sM}(X) @>>> \LocSys^\dr_{\sP,\sigma_\sM}(X) \\
@VVV   @VVV  \\
\LocSys^{\on{restr}}_{\sG}(X)  @>>>  \LocSys^\dr_{\sG}(X). 
\endCD
$$

Consider the composite morphism
\begin{equation} \label{e:uniformaze LocSys}
\CD
\LocSys^{\on{restr}}_{\sP,\sigma_\sM}(X) @>>> \LocSys^\dr_{\sP,\sigma_\sM}(X) \\
& & @VVV   \\
& &  \LocSys^\dr_{\sG}(X)
\endCD
\end{equation} 

Given that the map $\LocSys^{\on{restr}}_\sG(X)\to \LocSys^\dr_\sG(X)$ is a monomorphism, an easy diagram chase,
using properties (2) and (3) of the uniformization morphism in \secref{sss:properties of uniformization},  
shows that in order to prove \thmref{t:compare dR}, it suffices to show that the composite morphism \eqref{e:uniformaze LocSys}
is schematic and proper. 

\medskip

This follows from the combination of the next three assertions:

\begin{prop}  \label{p:comp de Rham P sigma}
The map 
$$\LocSys^{\on{restr}}_{\sP,\sigma_\sM}(X) \to \LocSys^\dr_{\sP,\sigma_\sM}(X)$$
is an isomorphism.
\end{prop}

\begin{prop} \label{p:p spec prop}
The map 
$$\sfp:\LocSys^\dr_{\sP}(X)\to \LocSys^\dr_{\sG}(X)$$
is schematic and proper.
\end{prop}

\begin{prop} \label{p:irred closed}
For a reductive group $\sG$ and an irreducible local system $\sigma$, the resulting map
$$\on{pt}/\on{Stab}_\sG(\sigma)\to \LocSys^\dr_{\sG}(X)$$
is a closed embedding. 
\end{prop}

\begin{rem} \label{r:image}
Note that the combination of the above three propositions describes the ind-closed substack
$$^{\on{red}}\!\LocSys^{\on{restr}}_{\sG}(X)\subset {}^{\on{red}}\!\LocSys^\dr_{\sG}(X).$$

Namely, it equals the disjoint union over classes of association of $(\sP,\sigma_\sM)$ of the unions
of the images of the maps
$$^{\on{red}}\!\LocSys^\dr_{\sP,\sigma_\sM}(X)\to {}^{\on{red}}\!\LocSys^\dr_{\sG}(X)$$
within a given class. 

\end{rem} 

\sssec{}

We now prove the above three propositions.

\medskip

The assertion of \propref{p:comp de Rham P sigma} follows by tracing the proof of \propref{p:radical}: namely, 
in the situation of {\it loc.cit.}, for an $S$-point of $\LocSys_{\sG_2}^{\on{restr}}$, the map
$$S\underset{\LocSys_{\sG_2}^{\on{restr}}}\times \LocSys_{\sG_1}^{\on{restr}}\to
S\underset{\LocSys^\dr_{\sG_2}(X)}\times \LocSys^\dr_{\sG_1}(X)$$
is an isomorphism. Indeed, in both cases, this fiber product classifies null-homotopies for the same class.

\medskip

\propref{p:p spec prop} is well-known: it follows from the fact that the map
$$\LocSys^\dr_{\sP}(X)\to \LocSys^\dr_{\sG}(X)\underset{\on{pt}/\sG}\times \on{pt}/\sP$$
is a closed embedding, where $\LocSys^\dr_{\sG}(X)\to \on{pt}/\sG$ is obtained by taking
the fiber at some point $x\in X$.

\medskip

It remains to prove \propref{p:irred closed}.

\begin{proof}[Analytic proof]

We can assume that $k=\BC$. Clearly, 
$$\on{pt}/\on{Stab}_\sG(\sigma)\to \LocSys^\dr_{\sG}(X)$$
is a locally closed embedding. To prove that it is a closed embedding, it is enough to show that its image is closed
\emph{in the analytic topology}. 

\medskip

Using Riemann-Hilbert, we identify the \emph{analytic stack} underlying ${}^{\on{cl}}\!\LocSys^\dr_{\sG}(X)$ with its Betti version 
${}^{\on{cl}}\!\LocSys^{\on{Betti}}_{\sG}(X)$
(see \secref{sss:loc const shv}). Hence, 
it is enough to show that the map
$$\on{pt}/\on{Stab}_\sG(\sigma)\to \LocSys^{\on{Betti}}_{\sG}(X)$$
is a closed embedding. 

\medskip

However, in this case the assertion follows from \propref{p:irred closed Betti} below. 

\end{proof}

\ssec{Algebraic proof of \propref{p:irred closed}}

\sssec{}

The map $\on{pt}/\on{Stab}_\sG(\sigma) \to \LocSys^\dr_{\sG}(X)$ is a priori
a locally closed embedding. Hence, in order to prove that it is actually a closed
embedding, it is enough to show that it is proper. We will do so by applying the
valuative criterion.

\medskip

Thus, is enough to show that for a smooth affine curve $C$ over $\sfe$ and a point $c\in C$, 
given a map
\begin{equation} \label{e:DVR map}
C\to \LocSys^\dr_{\sG}(X),
\end{equation} 
such that the composite map
$$(C-c)\to C\to \LocSys^\dr_{\sG}(X)$$
factors as
$$(C-c) \to \on{pt}/\on{Stab}_\sG(\sigma) \to \LocSys^\dr_{\sG}(X),$$
then the initial map \eqref{e:DVR map} also factors as
\begin{equation} \label{e:DVR map bis}
C\to \on{pt}/\on{Stab}_\sG(\sigma) \to \LocSys^\dr_{\sG}(X).
\end{equation} 

\medskip

Furthermore, it is enough to show that there exists a covering
$$\wt{C}\to C,$$
allowed to be branched at $c$, such that the composition
\begin{equation} \label{e:DVR map branched}
\wt{C}\to C\to \LocSys^\dr_{\sG}(X)
\end{equation} 
factors as
$$\wt{C} \to \on{pt}/\on{Stab}_\sG(\sigma) \to \LocSys^\dr_{\sG}(X).$$

\sssec{}

The assertion is easy if $\sG$ is a torus.  Hence, we obtain that for the induced bundles 
with respect to $\sG/[\sG,\sG]$, the given isomorphism indeed extends over over all $C\times X$. 
Modifying by means of a local system with respect to $Z_\sG$, we can thus assume that 
the induced local systems for $\sG/[\sG,\sG]$ are trivial. Hence, we can replace $\sG$
by $[\sG,\sG]$, i.e., we can assume that $\sG$ is semi-simple. 

\sssec{}

Since $\sigma$ was assumed irreducible and $\sG$ semi-simple, the group $\on{Stab}_\sG(\sigma)$ is finite. 
The given map $(C-c) \to \on{pt}/\on{Stab}_\sG(\sigma)$ corresponds to an \'etale covering
of $C-c$. Let $\wt{C}$ denote its normalization over $C$; let $\wt{c}$ be the preimage
of $c$ in $\wt{C}$.

\medskip

By construction, the map
$$(\wt{C}-\wt{c})\to (C-c) \to \on{pt}/\on{Stab}_\sG(\sigma) \to \LocSys^\dr_{\sG}(X)$$
factors as 
$$(\wt{C}-\wt{c}) \to \on{pt} \overset{\sigma}\to \LocSys^\dr_{\sG}(X).$$

We will show that the map \eqref{e:DVR map branched} also factors 
\begin{equation} \label{e:DVR map bis branched}
\wt{C} \to \on{pt} \overset{\sigma}\to \LocSys^\dr_{\sG}(X),
\end{equation} 
in a way compatible with the restriction to $\wt{C}-\wt{c}$. 

\sssec{}

The maps \eqref{e:DVR map branched} and \eqref{e:DVR map bis branched} correspond to $\sG$-bundles
$\CF^1_\sG$ and $\CF^2_\sG$ on $\wt{C}\times X$, each equipped with a connection,
and we are given an isomorphism of these data over $(\wt{C}-\wt{c})\times X$. We wish to show that
this isomorphism extends over all $\wt{C}\times X$. 

\medskip

Let $\eta_X$ denote the generic point of $X$. Then the relative position of $\CF^1_\sG$ and $\CF^2_\sG$ at 
$\wt{c}\times \eta_X$ is a cell of the affine Grassmannian of $\sG$, and hence corresponds to a coweight $\lambda$ of
$\sG$, which is $0$ if and only if the isomorphism between $\CF^1_\sG$ and $\CF^2_\sG$ extends 
over all $\wt{C}\times X$. 

\medskip

Furthermore, the restrictions of both $\CF^1_\sG$ and $\CF^2_\sG$ to $\wt{c}\times \eta_X$
acquire a reduction to the corresponding standard parabolic $\sP$ (it corresponds to those vertices $i$ of the
Dynkin diagram, for which $\langle \check\alpha_i,\lambda\rangle=0$).  
These reductions to $\sP$ are horizontal with respect to the connection along $\eta_X$. 

\sssec{}

Note that $\CF^2_\sG$ is isomorphic to the constant family corresponding to $\sigma$, so $\CF^2_\sG|_{\wt{c}\times X}$ is also given by 
$\sigma$. By the valuative criterion for $\sG/\sP$, its reduction to $\sP$ over $\wt{c}\times \eta_X$ extends to all of $\wt{c}\times X$.
However, since $\sigma$ was assumed irreducible, we have $\sP=\sG$. Hence, $\lambda=0$, as required.

\qed[\propref{p:irred closed}]

\ssec{The Betti version of $\LocSys_\sG(X)$}  \label{ss:Betti}

From this point until the end of this section we let $\sfe$ be an arbitrary algebraically closed field
of characteristic $0$. 

\sssec{}

Let $\CX$ be a connected object of $\Spc$.  We define the prestack $\LocSys^{\on{Betti}}_\sG(\CX)$ to be
$$(\on{pt}/\sG)^\CX=\bMaps(\CX,\on{pt}/\sG).$$

I.e., for $S\in \affSch_{/\sfe}$,
$$\Maps_{\on{PreStk}}(S,\LocSys^{\on{Betti}}_\sG(\CX))=\Maps_{\Spc}(\CX,\Maps_{\on{PreStk}}(S,\on{pt}/\sG)).$$

The fact that $\on{pt}/\sG$ admits $(-1)$-connective corepresentable 
deformation theory formally implies that the same is true for $\LocSys^{\on{Betti}}_\sG(\CX)$.

\sssec{} \label{sss:X compact}

Assume for a moment that $\CX$ is compact, i.e., is a retract of a space that 
can be obtained by a finite operation of taking push-outs from $\{*\}\in \Spc$. 

\medskip

In this case, it is clear from the definitions that $\LocSys^{\on{Betti}}_\sG(\CX)$
is locally almost of finite type.

\sssec{}

We claim: 

\begin{prop}  \label{p:LocSys Betti}
The prestack $\LocSys^{\on{Betti}}_\sG(\CX)$ is a derived algebraic stack. It can be realized as a quotient of an affine scheme
(to be denoted $\LocSys_\sG^{\on{Betti},\on{rigid}_x}(\CX)$) by an action of $\sG$.
\end{prop}

\begin{proof}

Choose a base point $x\in \CX$. Denote
$$\LocSys_\sG^{\on{Betti},\on{rigid}_x}(\CX):=\LocSys^{\on{Betti}}_\sG(\CX)\underset{\on{pt}/\sG}\times \on{pt},$$
where the map $\LocSys^{\on{Betti}}_\sG(\CX)\to \on{pt}/\sG$ is given by restriction to $x$. 

\medskip

We have a natural action of $\sG$ on $\LocSys_\sG^{\on{Betti},\on{rigid}_x}(\CX)$ so that
$$\LocSys^{\on{Betti}}_\sG(\CX)\simeq \LocSys_\sG^{\on{Betti},\on{rigid}_x}(\CX)/\sG.$$

We will show that $\LocSys_\sG^{\on{Betti},\on{rigid}_x}(\CX)$ is an affine scheme. The fact that 
$\LocSys^{\on{Betti}}_\sG(\CX)$ admits $(-1)$-connective corepresentable deformation theory implies that 
$\LocSys_\sG^{\on{Betti},\on{rigid}_x}(\CX)$ admits connective corepresentable deformation theory (we argue
as in \corref{c:def}(b) and use the assumption that $\CX$ is connected).

\medskip

Hence, by \cite[Theorem 18.1.0.1]{Lu3}, in order to show that $\LocSys_\sG^{\on{Betti},\on{rigid}_x}(\CX)$ is an affine scheme,
it suffices to show that
$$^{\on{cl}}\!\LocSys_\sG^{\on{Betti},\on{rigid}_x}(\CX)$$
is a classical affine scheme.

\medskip

Denote $$\Gamma:=\pi_1(\CX,x).$$ It follows from the definitions that for $S\in {}^{\on{cl}}\!\affSch_{/\sfe}$, the space 
$\Maps(S,\LocSys_\sG^{\on{Betti},\on{rigid}_x}(\CX))$ is a set of homomorphisms $\Gamma\to \sG$,
parameterized by $S$. 

\medskip

I.e., $\LocSys_\sG^{\on{Betti},\on{rigid}_x}(\CX)$ is a subfunctor of
$$S\mapsto \Maps(S,\sG)^\Gamma\simeq \Maps(S,\sG^{\Gamma}),$$
consisting of elements that obey the group law, i.e.,
$$\sG^{\Gamma} \underset{\sG^{\Gamma\times \Gamma}}\times \on{pt}.$$

\medskip

Since $\sG^{\Gamma}$ and $\sG^{\Gamma\times \Gamma}$ are affine schemes, we obtain that so is
$\LocSys_\sG^{\on{Betti},\on{rigid}_x}(\CX)$. 

\end{proof} 

\begin{rem}

It follows from \secref{sss:X compact} that if $\CX$ compact, then $\LocSys_\sG^{\on{Betti},\on{rigid}_x}(\CX)$
is almost of finite type.

\end{rem} 

\sssec{} \label{sss:Vect X}

Let us now rewrite the definition of $\LocSys^{\on{Betti}}_\sG(\CX)$ slightly differently. Consider the DG category
$$\Vect_\sfe^\CX\simeq \on{Funct}(\CX,\Vect_\sfe),$$
see \cite[Sects. 1.4.1-2]{GKRV}.

\medskip

For any DG category $\bC$, we have a tautological functor
\begin{equation} \label{e:tensor by Funct}
\bC\otimes \Vect_\sfe^\CX \to \bC^\CX,
\end{equation} 
which is an equivalence if $\bC$ is dualizable (or if $\CX$ is compact). 

\medskip

Furthermore $\Vect_\sfe^\CX$ has a natural symmetric monoidal structure, and if $\bC$ is also
symmetric monoidal, the functor \eqref{e:tensor by Funct} is symmetric monoidal. 

\medskip

Assume for a moment that $\bC$ has a t-structure. Then $\bC^\CX$ also acquires a t-structure
(an object is connective/coconnetive if its value for any $x\in X$ is connective/coconnective). In particular, $\Vect_\sfe^\CX$
has a t-structure.  With respect to these t-structures, the functor \eqref{e:tensor by Funct} is t-exact. 

\sssec{} \label{sss:Betti coMaps}

By the definitions of $\LocSys^{\on{Betti}}_\sG(\CX)$ and of $\on{pt}/\sG$, the value of $\LocSys^{\on{Betti}}_\sG(\CX)$ on 
an affine scheme $S$ is the space of functors
$$\CX \times \Rep(\sG)\to \QCoh(S)$$
that are symmetric monoidal and right t-exact in the second variable. By the above, this is the same 
as the space of right t-exact symmetric monoidal functors
$$\Rep(\sG)\to \QCoh(S)\otimes \Vect_\sfe^\CX.$$

\medskip

Thus, the prestack $\LocSys^{\on{Betti}}_\sG(\CX)$ is also of the form $\bMaps(\Rep(\sG),\bH)$, for 
$\bH:=\Vect_\sfe^\CX$.

\sssec{}  \label{sss:loc const shv}

Let now $X$ be CW complex. Let $\Shv^{\on{all}}_{\on{loc.const.}}(X)$ be the category of sheaves of 
$\sfe$-vector spaces \emph{with locally constant cohomologies}.

\medskip

We define the prestack $\LocSys^{\on{Betti}}_\sG(X)$ as follows. It sends an affine scheme $S$ to the space 
of right t-exact symmetric monoidal functors
$$\Rep(\sG)\to \QCoh(S)\otimes \Shv^{\on{all}}_{\on{loc.const.}}(X).$$

In other words,
$$\LocSys^{\on{Betti}}_\sG(X):=\bMaps(\Rep(\sG),\bH)$$
for $\bH:=\Shv^{\on{all}}_{\on{loc.const.}}(X)$. 

\sssec{}  \label{sss:geom real}

Let us write $X$ as a geometric realization of an object $\CX$ of $\Spc$. 

\medskip

In this case we have a canonical t-exact equivalence of symmetric monoidal categories
$$\Shv^{\on{all}}_{\on{loc.const.}}(X)\simeq  \Vect_\sfe^\CX.$$
Hence, we obtain that in this case we have a canonical isomorphism of prestacks
$$\LocSys_\sG^{\on{Betti}}(X)\simeq \LocSys^{\on{Betti}}_\sG(\CX).$$

Thus, the results pertaining to $\LocSys^{\on{Betti}}_\sG(\CX)$ that we have reviewed above carry over to 
$\LocSys_\sG^{\on{Betti}}(X)$ as well. 

\ssec{The coarse moduli space of Betti local systems} \label{ss:coarse}

In this subsection we will make a digression and discuss the coarse 
moduli space of Betti local systems (a.k.a. character variety).

\medskip

In this subsection we will assume that $\sG$ is reductive. 

\sssec{}

Let $\CX$ be as in \secref{ss:Betti}. Consider the object of $\on{ComAlg}(\Vect_\sfe)$ given by
\begin{equation} \label{e:functions on coarse}
\Gamma(\LocSys^{\on{Betti}}_\sG(\CX),\CO_{\LocSys^{\on{Betti}}_\sG(\CX)}).
\end{equation} 

Note that it is connective: this follows from the presentation of $\LocSys^{\on{Betti}}_\sG(\CX)$ as 
$$\LocSys_\sG^{\on{Betti},\on{rigid}_x}(\CX)/\sG,$$
so that
$$\Gamma(\LocSys^{\on{Betti}}_\sG(\CX),\CO_{\LocSys^{\on{Betti}}_\sG(\CX)})=\inv_\sG\left(\Gamma(\LocSys^{\on{rigid}_x}_\sG(\CX),\CO_{\LocSys^{\on{rigid}_x}_\sG(\CX)})\right),$$
and 
using the fact that $\LocSys_\sG^{\on{Betti},\on{rigid}_x}(\CX)$ is an affine scheme and $\sG$ is reductive, so the functor $\inv_\sG$
is t-exact. 

\medskip

Note that if $\CX$ is compact, so that $\LocSys^{\on{Betti}}_\sG(\CX)$ is almost of finite type, the algebra \eqref{e:functions on coarse} 
is also almost of finite type.

\sssec{} \label{sss:coarse Betti}

Set
$$\LocSys_\sG^{\on{Betti,coarse}}(\CX):=\Spec\left(\Gamma(\LocSys^{\on{Betti}}_\sG(\CX),\CO_{\LocSys^{\on{Betti}}_\sG(\CX)})\right).$$

We have a tautologically defined map
\begin{equation} \label{e:projection to coarse}
\brr:\LocSys^{\on{Betti}}_\sG(\CX)\to \LocSys_\sG^{\on{Betti,coarse}}(\CX).
\end{equation}

\sssec{}

Let us describe the classical affine scheme underlying $\LocSys_\sG^{\on{Betti,coarse}}(\CX)$. Recall the notation
$$\Gamma:=\pi_1(X,x),$$ 
and consider the affine scheme
$$\bMaps_{\on{Grp}}(\Gamma,\sG),$$
which is acted on by $\sG$ by conjugation.

\medskip

Set
$$\bMaps_{\on{Grp}}(\Gamma,\sG)/\!/\!\on{Ad}(\sG)=\Spec\left(\inv_\sG\left(\Gamma(\bMaps_{\on{Grp}}(\Gamma,\sG),\CO_{\bMaps_{\on{Grp}}(\Gamma,\sG)})\right)\right).$$

As we have seen in the course of the proof of \propref{p:LocSys Betti}, we have the isomorphisms
$$^{\on{cl}}\!\LocSys^{\on{rigid}_x}_\sG(\CX)\simeq {}^{\on{cl}}\bMaps_{\on{Grp}}(\Gamma,\sG)$$
and
$$^{\on{cl}}\!\LocSys^{\on{Betti}}_\sG(\CX)\simeq {}^{\on{cl}}\bMaps_{\on{Grp}}(\Gamma,\sG)/\on{Ad}(\sG).$$

Hence, since $\sG$ reductive, we have
$$^{\on{cl}}\!\LocSys_\sG^{\on{Betti,coarse}}(\CX)\simeq {}^{\on{cl}}\bMaps_{\on{Grp}}(\Gamma,\sG)/\!/\!\on{Ad}(\sG).$$
%
%

\sssec{}

For future use we now quote the following fundamental result of \cite{Ri}:

\begin{thm} \label{t:Rich}  Let $\Gamma$ be an abstract group. Then 
%
%
two $\sfe$-points of the stack $\bMaps_{\on{Grp}}(\Gamma,\sG)/\on{Ad}(\sG)$ 
map to the same point in the affine scheme $\bMaps_{\on{Grp}}(\Gamma,\sG)/\!/\!\on{Ad}(\sG)$ if and only if they have isomorphic semi-simplifications.
%
%
%
\end{thm} 

By the above, we immediately obtain:

\begin{cor} \label{c:Rich}  
%
%
Two $\sfe$-points of the stack $\LocSys^{\on{Betti}}_\sG(\CX)$ 
get sent by $\brr$ to the same point in the affine scheme $\LocSys_\sG^{\on{Betti,coarse}}(\CX)$ if and only if they have isomorphic semi-simplifications.
%
%
\end{cor} 

\ssec{Relationship of the restricted and Betti versions}  \label{ss:restr vs Betti}

In this subsection we let $X$ be a smooth connected algebraic variety\footnote{The material of this and the next
 subsection is equally applicable, when instead of $X$ we take a 
connected finite CW complex. In this case we let $\qLisse(X)$ be the full subcategory of $\Shv^{\on{all}}_{\on{loc.const.}}(X)$
consisting of objects such that each of their cohomologies (with respect to the usual t-structure)
is locally finite as a representation of $\pi_1(X,x)$.}  
over $\BC$.

\sssec{}

Consider the functor 

\begin{equation} \label{e:LS to loc const}
\qLisse(X)\to \Shv^{\on{all}}_{\on{loc.const.}}(X)
\end{equation}

We claim:

\begin{prop} \label{p:lisse to loc const}
The functor \eqref{e:LS to loc const} is fully faithful.
\end{prop} 

%

\begin{rem}
Note that, unlike the de Rham case, in the Betti setting, the fully faithulness of \eqref{e:LS to loc const}
is not a priori evident (because objects from $\Shv(X)^{\on{constr}}$ are \emph{not} compact as objects
in the category of \emph{all} sheaves of $\sfe$-vector spaces on $X$). 
\end{rem} 

\begin{proof}

Since both categories are left-complete and \eqref{e:LS to loc const} is t-exact, 
it is sufficient to show that it induces fully faithful functors
\begin{equation} \label{e:LS to loc const n}
(\qLisse(X))^{\geq -n}\to (\Shv^{\on{all}}_{\on{loc.const.}}(X))^{\geq -n}.
\end{equation}

Now, 
$$(\iLisse(X))^{\geq -n}\to (\qLisse(X))^{\geq -n}$$
is an equivalence, and hence the functor \eqref{e:LS to loc const n} sends compacts to compacts.

\medskip

Since $(\iLisse(X))^{\geq -n}$ is compactly generated (by $(\Lisse(X))^{\geq -n}$) and 
$$\Lisse(X)\to \Shv^{\on{all}}_{\on{loc.const.}}(X)$$
is fully faithful, we obtain that \eqref{e:LS to loc const n} is fully faithful.

\end{proof} 

\sssec{}

The functor \eqref{e:LS to loc const} defines a map
\begin{equation} \label{e:restr to Betti}
\LocSys^{\on{restr}}_\sG(X)\to \LocSys_\sG^{\on{Betti}}(X).
\end{equation} 

As in the de Rham case, from \propref{p:lisse to loc const} we obtain that the map \eqref{e:restr to Betti}
is a monomorphism. 

\begin{rem}
This remark is parallel to Remark \ref{sss:Gm bis}. Let us explain how the difference between $\LocSys_\sG^{\on{Betti}}(X)$ and $\LocSys^{\on{restr}}_\sG(X)$ 
plays out in the simplest cases when $\sG=\BG_m$ and $\sG=\BG_a$. Take $S=\Spec(R)$ to be classical.  

\medskip

In this case, an $S$-point of $\LocSys_{\BG_m}(X)$ is a homomorphism
$$\pi_1(X)\to R^\times.$$

By contrast, if we further assume $S$ to be reduced, then an $S$-point of $\LocSys^{\on{restr}}_{\BG_m}(X)$ is a homomorphism
$$\pi_1(X)\to \sfe^\times.$$

\medskip

Take now $\sG=\BG_a$. In this case, by Example \secref{sss:Ga}, the map
$$\LocSys^{\on{restr}}_{\BG_a}(X)\to \LocSys_{\BG_a}(X)$$
is an isomorphism.
\end{rem}

\begin{rem}
A remark parallel to Remark \ref{r:G' to G deR} holds in the Betti context as well, i.e., for a closed
embedding $\sG'\to \sG$, the diagram
$$
\CD
\LocSys^{\on{restr}}_{\sG'}(X) @>>> \LocSys_{\sG'}(X) \\
@VVV @VVV \\
\LocSys^{\on{restr}}_\sG(X) @>>> \LocSys_\sG^{\on{Betti}}(X) 
\endCD
$$
is a fiber square.
\end{rem}

\sssec{}

We have also the following statements that are completely parallel with the de Rham situation
(with the same proofs):

\begin{prop} \label{p:formal compl Betti}
The map \eqref{e:restr to Betti} is a formal isomorphism, i.e., identifies $\LocSys^{\on{restr}}_\sG(X)$ with
its formal completion inside $\LocSys_\sG^{\on{Betti}}(X)$.
\end{prop}

\begin{thm} \label{t:compare Betti}
The map 
$$\LocSys^{\on{restr}}_\sG(X)\to \LocSys_\sG^{\on{Betti}}(X)$$
is a closed embedding \emph{at the reduced level} for \emph{each connected component of} $\LocSys^{\on{restr}}_\sG(X)$.
\end{thm} 

\begin{cor}  \label{c:compare Betti}
The subfunctor $$\LocSys^{\on{restr}}_\sG(X)\subset \LocSys_\sG^{\on{Betti}}(X)$$
is the disjoint union of formal completions of a collection of pairwise
non-intersecting closed substacks of ${}^{\on{red}}\!\LocSys_\sG^{\on{Betti}}(X)$. 
\end{cor} 

\begin{cor}   \label{c:compare Betti ind}
The map $\LocSys^{\on{restr}}_\sG(X)\to \LocSys_\sG^{\on{Betti}}(X)$ is an ind-closed embedding.
\end{cor}

Note, however, that we still have to supply a proof of the Betti version of \propref{p:irred closed}:

\begin{prop} \label{p:irred closed Betti}
For a reductive group $\sG$ and an irreducible local system $\sigma$, the resulting map
$$\on{pt}/\on{Stab}_\sG(\sigma)\to \LocSys^{\on{Betti}}_{\sG}(X)$$
is a closed embedding. 
\end{prop}

The proof is given in \secref{sss:proof of irred Betti closed} below. 

\begin{rem} \label{r:image Betti}

Note that as in Remark \ref{r:image}, we obtain that the 
image of 
\begin{equation} \label{e:restr to Betti red}
^{\on{red}}\!\LocSys^{\on{restr}}_\sG(X)\to  {}^{\on{red}}\!\LocSys_\sG^{\on{Betti}}(X)
\end{equation}
is the ind-closed substack
that equals the disjoint union over classes of association of $(\sP,\sigma_\sM)$ of the unions
of the images of the maps
$$^{\on{red}}\!\LocSys^{\on{Betti}}_{\sP,\sigma_\sM}(X)\to {}^{\on{red}}\!\LocSys^{\on{Betti}}_{\sG}(X)$$
within a given class. 

\medskip

In \secref{ss:char} below we will give an alternative description of the image of
\eqref{e:restr to Betti red}, which is specific to the Betti situation. 

\end{rem}

\ssec{Comparison of $\LocSys^{\on{restr}}_\sG(X)$ vs $\LocSys_\sG^{\on{Betti}}(X)$ via the coarse moduli space} \label{ss:char}

Let $X$ be as in \secref{ss:restr vs Betti}. 
We will give a more explicit description of $\LocSys^{\on{restr}}_\sG(X)$ as a subfunctor 
of $\LocSys_\sG^{\on{Betti}}(X)$. 

\sssec{}

Let $\sG_{\on{red}}$ denote the reductive quotient of $\sG$. We have a fiber square
\begin{equation} \label{e:red to red}
\CD
\LocSys^{\on{restr}}_\sG(X)  @>>> \LocSys_\sG^{\on{Betti}}(X)  \\
@VVV  @VVV  \\
\LocSys^{\on{restr}}_{\sG_{\on{red}}}(X)  @>>> \LocSys_{\sG_{\on{red}}}^{\on{Betti}}(X). 
\endCD
\end{equation} 

Hence, in order to describe $\LocSys^{\on{restr}}_\sG(X)$ as a subfunctor of $\LocSys_\sG^{\on{Betti}}(X)$,
it is enough to do so for $\sG$ replaced by $\sG_{\on{red}}$. So, from now until the end of this
subsection we will assume that $\sG$ is reductive.

\sssec{} \label{sss:proof of irred Betti closed}

First, we are  going to deduce \propref{p:irred closed Betti} from \corref{c:Rich}: 

\begin{proof}

Let $\sigma$ be irreducible, and consider the closed substack
$$\on{pt}\underset{\LocSys_\sG^{\on{Betti,coarse}}(X)}\times \LocSys_\sG^{\on{Betti}}(X) \subset \LocSys_\sG^{\on{Betti}}(X),$$
where 
$$\on{pt}\to \LocSys_\sG^{\on{Betti,coarse}}(X)$$
is given by $\brr(\sigma)$. 

\medskip

By \corref{c:Rich} and the irreducibility assumption on $\sigma$, the above stack contains
a unique isomorphism class of $\sfe$-points. Hence, the map
$$\on{pt}/\on{Stab}_\sG(\sigma)\to \on{pt}\underset{\LocSys_\sG^{\on{Betti,coarse}}(X)}\times \LocSys_\sG^{\on{Betti}}(X)$$
is an isomorphism of the underlying reduced substacks. In particular, it is a closed embedding.

\end{proof} 

\sssec{}

We now claim:

\begin{thm}  \label{t:coarse fibers}
The subfunctor 
$$^{\on{red}}\!\LocSys^{\on{restr}}_\sG(X)\subset {}^{\on{red}}\!\LocSys_\sG^{\on{Betti}}(X)$$
is the disjoint union of the fibers of the map $\brr$ of \eqref{e:projection to coarse}.
\end{thm} 

Combining with \corref{c:compare Betti}, we obtain: 

\begin{cor}  \label{c:coarse fibers}
The subfunctor $\LocSys^{\on{restr}}_\sG(X)\subset \LocSys_\sG^{\on{Betti}}(X)$ is
the disjoint union of formal completions of the fibers of the map
$$\LocSys_\sG^{\on{Betti}}(X)\to \LocSys_\sG^{\on{Betti,coarse}}(X).$$
\end{cor} 

The rest of this subsection is devoted to the proof of \thmref{t:coarse fibers}.

\sssec{}  \label{sss:coarse fiber}

We will prove the following slightly more precise statement (which would imply \thmref{t:coarse fibers}
in view of Remark \ref{r:image Betti}): 

\medskip

Fix a class of association of pairs $(\sP,\sigma_\sM)$. For each element in this class pick a Levi splitting
$$\sP\leftrightarrows \sM,$$
and consider the induced $\sG$-local system. Note, however, that by \corref{c:assoc vs semisimple bis}, 
these $\sG$-local systems are all isomorphic (for different elements $(\sP,\sigma_\sM)$ in the given
class); denote the resulting local system by $\sigma_\sG$. 

\medskip

We will show that the reduced substack underlying
\begin{equation}  \label{e:coarse fiber}
\on{pt}\underset{\LocSys_\sG^{\on{Betti,coarse}}(X)}\times \LocSys_\sG^{\on{Betti}}(X) 
\end{equation} 
(where $\on{pt}\to \LocSys_\sG^{\on{Betti,coarse}}(X)$ is given by $\brr(\sigma_\sG)$), equals 
the union of the images of the maps
\begin{equation} \label{e:uniform again} 
\LocSys^{\on{Betti}}_{\sP,\sigma_\sM}(X) \to \LocSys_\sG^{\on{Betti}}(X),
\end{equation}
where the union is taken over the pairs $(\sP,\sigma_\sM)$ in our chosen class of association. 

\sssec{}

We claim: 

\begin{prop} \label{p:assoc char}
Let $\sP\leftrightarrows \sM$ be a parabolic with a Levi splitting. Let $\sigma_\sM$ be an
irreducible $\sM$-local system, 
and let $\sigma_\sG$ be the induced $\sG$-local system. Then the composite 
$$^{\on{red}}\!\LocSys^{\on{Betti}}_{\sP,\sigma_\sM}(X) \to {}^{\on{red}}\!\LocSys_\sG^{\on{Betti}}(X)\overset{\brr}\to
{}^{\on{red}}\!\LocSys_\sG^{\on{Betti,coarse}}(X)$$
factors as 
$$^{\on{red}}\!\LocSys^{\on{Betti}}_{\sP,\sigma_\sM}(X) \to \on{pt} \overset{\brr(\sigma_\sG)}\longrightarrow {}^{\on{red}}\!\LocSys_\sG^{\on{Betti,coarse}}(X).$$
\end{prop} 

\begin{proof}

Note that all $\sfe$-points of $\LocSys_\sG^{\on{Betti}}(X)$
obtained from $\sfe$-points of $\LocSys^{\on{Betti}}_{\sP,\sigma_\sM}(X)$ have $\sigma_\sG$ as their semi-simplification.

\medskip

Hence, the assertion of the proposition follows from \thmref{t:Rich}.

\end{proof}

We will now deduce from \propref{p:assoc char} the description of \eqref{e:coarse fiber} as the union of the images
of the maps \eqref{e:uniform again}. 

\sssec{}

Indeed, on the one hand, \propref{p:assoc char} implies that the images of the maps \eqref{e:uniform again} 
(at the reduced level) indeed lie in the fiber \eqref{e:coarse fiber}. 

\medskip

On the other hand, take an $\sfe$-point $\sigma'_\sG$ in the fiber \eqref{e:coarse fiber}, and let $(\sP',\sigma_{\sM'})$ be 
a pair such $\sigma'_\sG$ lies in the image of 
$$\LocSys^{\on{Betti}}_{\sP',\sigma_{\sM'}}(X) \to \LocSys_\sG^{\on{Betti}}(X).$$

We need to show that $(\sP',\sigma_{\sM'})$ lies in our class of association. However, by \propref{p:assoc char},
the $\sG$-local system, induced from $\sigma_{\sM'}$, is isomorphic to $\sigma_\sG$. This implies the result by
\corref{c:assoc vs semisimple bis}. 

\qed[\thmref{t:coarse fibers}]

\section{Geometric properties of $\LocSys^{\on{restr}}_\sG(X)$} \label{s:geom properties}

In this section we will assume that $\sG$ is reductive. The goal of this section is to establish a version, adapted to
$\LocSys^{\on{restr}}_\sG(X)$, of the picture
$$\brr:\LocSys_\sG(X)\to \LocSys^{\on{coarse}}_\sG(X)$$
that we have in the Betti case (see \secref{sss:coarse Betti}). This will be stated as \thmref{t:coarse restr}, which constructs the desired picture
$$\brr:\CZ \to \CZ^{\on{restr}}$$
for each connected component $\CZ$ of $\LocSys^{\on{restr}}_\sG(X)$. 

\medskip

Prior to doing so, we show that $\LocSys^{\on{restr}}_\sG(X)$ has the following two geometric properties:
it is \emph{mock-affine} and \emph{mock-proper}.

\ssec{``Mock-properness" of $^{\on{red}}\!\LocSys_\sG^{\on{restr}}(X)$}

\sssec{} 

Let $\CZ$ be a quasi-compact algebraic stack locally almost of finite type over $\sfe$. Let
$$\Coh(\CZ)\subset \QCoh(\CZ)$$
be the full subcategory consisting of objects whose pullback under a smooth cover (equivalently, any map)
$$S\to \CZ, \quad S\in \affSch_{\on{aft}/\sfe}$$
belongs to $\Coh(S)\subset \QCoh(S)$.

\medskip

We shall say that $\CZ$ is \emph{mock-proper} if the functor
$$\Gamma(\CZ,-):\QCoh(\CZ)\to \Vect_\sfe$$
sends $\Coh(\CZ)$ to $\Vect_\sfe^c$.

\begin{rem}
This definition is equivalent to one in \cite[Sect. 6.5]{Ga3}. Indeed, the subcategory
$$\Dmod(\CZ)^c\subset \Dmod(\CZ)$$
is generated under finite colimits by the image of $\Coh(\CZ)$ along induction functor
$$\ind_{\Dmod}:\QCoh(\CZ)\to \Dmod(\CZ).$$
\end{rem} 

\sssec{Examples} \label{sss:ex mock} \hfill

\medskip

\noindent(i) If $\CZ$ is a scheme, then it is mock-proper as a stack if and only it is proper as a scheme.

\medskip

\noindent(ii) The stack $\on{pt}/\sH$ is mock-proper for any algebraic group $\sH$.

\medskip

\noindent(iii) For a (finite-dimensional) vector space $V$, the stack $\on{Tot}(V)/\BG_m$ is mock-proper.
(This is just the fact that for a finitely generated graded $\Sym(V^\vee)$-module, its degree $0$ component
is finite-dimensional as a vector space.) 

\sssec{}

Let $\bH$ be a gentle Tannakian category, and let $\CZ$ be a connected component of 
$\bMaps(\Rep(\sG),\bH)$. Recall that according to
\thmref{t:main 1 abs}, its underlying reduced prestack $^{\on{red}}\CZ$ is actually a quasi-compact algebraic stack.

\medskip

We will prove:

\begin{thm} \label{t:mock proper}
The algebraic stack $^{\on{red}}\CZ$ is mock-proper.
\end{thm}

Of course, our main application is when $\bH=\qLisse(X)$, so that $\CZ$ is a connected component of 
$\LocSys_\sG^{\on{restr}}(X)$.

\medskip

The rest of the subsection is devoted to the proof of \thmref{t:mock proper}.

\sssec{}

Recall (see \secref{ss:assoc}) that to $\CZ$ there corresponds a class of association of pairs
$(\sP,\sigma_\sM)$, where $\sP$ is a parabolic in $\sG$ and $\sigma_\sM$ is an irreducible local system
with respect to the Levi quotient $\sM$ of $\sP$.

\medskip

Moreover, the resulting morphism 
$$\pi:\underset{(\sP,\sigma_\sM)}\sqcup\, {}^{\on{red}}\bMaps(\Rep(\sP),\bH)_{\sigma_\sM}\to {}^{\on{red}}\CZ$$
(the union is taken over the given class of association)
is proper and surjective at the level of geometric points.

\medskip

We claim that the category $\Coh({}^{\on{red}}\CZ)$ is generated under finite colimits and retracts
by the essential image of 
$$\Coh\left(\underset{(\sP,\sigma_\sM)}\sqcup\, {}^{\on{red}}\bMaps(\Rep(\sP),\bH)_{\sigma_\sM}\right)$$
along $\pi_*$.

\medskip

Indeed, this follows from the next general assertion:

\begin{lem}
Let $\pi:\CZ'\to \CZ$ be a proper map between algebraic stacks, surjective at the level of geometric points. 
Then $\Coh(\CZ)$ is generated 
under finite colimits and retracts by the essential image of $\Coh(\CZ')$ along $\pi_*$.
\end{lem}

\begin{proof}
First, since $\pi$ is proper, the functor $\pi_*$ does indeed send $\Coh(\CZ')$ to $\Coh(\CZ)$.
Since $\IndCoh(\CZ')$ is generated by $\Coh(\CZ')$ (see \cite[Proposition 3.5.1]{DrGa1}), the assertion 
of the lemma is equivalent to the fact that the essential image of $\IndCoh(\CZ')$ along
$$\pi^{\IndCoh}_*:\IndCoh(\CZ')\to \IndCoh(\CZ)$$
generates $\IndCoh(\CZ)$. This is equivalent to the fact that the right adjoint
$$\pi^!:\IndCoh(\CZ)\to \IndCoh(\CZ')$$
is conservative. However, the latter is \cite[Proposition 8.1.2]{Ga4}. 
\end{proof}

\sssec{} \label{sss:rigid P}

Thus, we obtain that it suffices to show that for a parabolic $\sP$
with Levi quotient $\sM$ and a $\sM$-local system $\sigma_\sM$, the algebraic stack
$$\bMaps(\Rep(\sP),\bH)_{\sigma_\sM}$$
is mock-proper. 

\medskip

We will consider separately two cases: when $\sP=\sG$ and when $\sP$ is a proper parabolic.
If $\sP=\sG$, 
$$\bMaps(\Rep(\sP),\bH)_{\sigma_\sM}=\on{pt}/\on{Aut}(\sigma_\sM)$$
and the assertion obvious. Hence, from now on we will assume that $\sP$ is a proper parabolic. 

\sssec{}

Let $\bMaps(\Rep(\sP),\bH)^{\on{rigid}}_{\sigma_\sM}$ be the following (algebraic) stack: it classifies the data of
$$(\sigma_\sP,\alpha,\epsilon),$$
where:

\begin{itemize} 

\item $\sigma_\sP$ is a point of $\bMaps(\Rep(\sP),\bH)$, 

\item $\alpha$ is an identification $\sM\overset{\sP}\times \sigma_\sP\simeq \sigma_\sM$, 
so that the pair $(\sigma_\sP,\alpha)$ is a point of
$$\on{pt}\underset{\bMaps(\Rep(\sM),\bH)}\times \bMaps(\Rep(\sP),\bH);$$

\item $\epsilon$ is an identification  
$$\oblv_\bH(\sigma_\sP) \simeq \sP\overset{\sM}\times \oblv_\bH(\sigma_{\sM}),$$
as points of $\on{pt}/\sP$, compatible with the datum of $\alpha$.

\end{itemize}

\medskip

The stack $\bMaps(\Rep(\sP),\bH)^{\on{rigid}}_{\sigma_\sM}$ carries an action of $\on{Aut}(\sigma_\sM)$
(by changing the datum of $\alpha$); in particular, it is acted on by $Z(\sM)^0$, the connected component of
the center of $\sM$. In addition, it carries an action of the (unipotent) group 
$$(\sN_\sP)_{\oblv_\bH(\sigma_\sM)}$$
(by changing the datum of $\epsilon$), where:

\begin{itemize} 

\item $\sN_\sP$ is the unipotent radical of $\sP$;

\item $(\sN_\sP)_{\oblv_\bH(\sigma_\sM)}$ is the twist of $\sN_\sP$ by the $\sM$-torsor $\oblv_\bH(\sigma_\sM)$,
using the adjoint action of $\sM$ on $\sN_\sP$.

\end{itemize}

Combining, we obtain an action on $\bMaps(\Rep(\sP),\bH)^{\on{rigid}}_{\sigma_\sM}$ of the semi-direct product 
$$\on{Aut}(\sigma_\sM)\ltimes (\sN_\sP)_{\oblv_\bH(\sigma_\sM)}.$$

We have:
$$\bMaps(\Rep(\sP),\bH)^{\on{rigid}}_{\sigma_\sM}/\on{Aut}(\sigma_\sM)\ltimes (\sN_\sP)_{\oblv_\bH(\sigma_\sM)}
\simeq \bMaps(\Rep(\sP),\bH)_{\sigma_\sM}.$$

\sssec{}

Choose a coweight $\BG_m\to Z(\sM)^0$, dominant and regular with respect to $\sP$ (i.e., one such that the adjoint action
of $\BG_m$ on $\sn_\sP$ has positive eigenvalues). Such a coweight exists by the assumption that $\sP$ is a proper
parabolic. 
 
\medskip

We claim that it suffices to show that the algebraic stack
\begin{equation} \label{e:semidir}
\bMaps(\Rep(\sP),\bH)^{\on{rigid}}_{\sigma_\sM}/\BG_m\ltimes (\sN_\sP)_{\oblv_\bH(\sigma_\sM)}
\end{equation}
is mock-proper. Indeed, the space global sections of an object in $\CF\in \QCoh(\bMaps(\Rep(\sP),\bH)_{\sigma_\sM})$ can be expressed
as invariants with respect to the quotient group $\on{Aut}(\sigma_\sM)/\BG_m$ on the space of global
sections of the pullback of $\CF$ to \eqref{e:semidir}.

\medskip

Furthermore, we claim that it suffices to show that the algebraic stack
$$\bMaps(\Rep(\sP),\bH)^{\on{rigid}}_{\sigma_\sM}/\BG_m$$
is mock proper. Indeed, let $\CF'$ be a quasi-coherent sheaf on \eqref{e:semidir}, and let $\CF''$ denote 
its pullback to $\bMaps(\Rep(\sP),\bH)^{\on{rigid}}_{\sigma_\sM}/\bG_m$.
Since the group $(\sN_\sP)_{\oblv_\bH(\sigma_\sM)}$ is unipotent, using the Chevalley complex that
computes Lie algebra cohomology, we obtain that the space of global sections of $\CF'$ admits a finite
filtration with subquotients of the form
$$\Gamma\left(\bMaps(\Rep(\sP),\bH)^{\on{rigid}}_{\sigma_\sM}/\BG_m,\CF''\otimes 
\Lambda^\cdot((\sn_\sP)_{\oblv_\bH(\sigma_\sM)})\right),$$
where $(\sn_\sP)_{\oblv_\bH(\sigma_\sM)}$ is the Lie algebra of $(\sN_\sP)_{\oblv_\bH(\sigma_\sM)}$.

\sssec{} \label{sss:contr N}

Note that the proof in Sects. \ref{sss:algebraic stack 1}-\ref{sss:algebraic stack 2} of the fact that the morphism
$$\bMaps(\Rep(\sP),\bH)\to \bMaps(\Rep(\sM),\bH)$$
is a relative algebraic stack implies that $\bMaps(\Rep(\sP),\bH)^{\on{rigid}}_{\sigma_\sM}$
is actually an affine (derived) scheme.

\medskip

Furthermore, the fact that $\BG_m$ acts on $\sn_\sP$ with positive eigenvalues implies that the action
of $\BG_m$ on $\bMaps(\Rep(\sP),\bH)^{\on{rigid}}_{\sigma_\sM}$ is \emph{contracting}: 

\medskip

Recall (see \cite[Sect. 1.4.4]{DrGa3}) that an action of $\BG_m$ on an affine scheme $Z$ is said to be contracting
if it can be extended to an action of the monoid $\BA^1$, so that the action of $0\in \BA^1$ factors as
$$Z\to \on{pt}\to Z.$$

\medskip

The required result follows now from the next general assertion, which generalizes Example (iii) in
\secref{sss:ex mock}:

\begin{lem}
Let $Z$ be an affine scheme almost of finite type, equipped with a contracting action of $\BG_m$. 
Then the algebraic stack $Z/\BG_m$ is mock-proper.
\end{lem}

\begin{proof}

Write $Z=\Spec(A)$. The $\BG_m$-action on $Z$ equips $A$ with a grading. The fact that the $\BG_m$-action
is contracting is equivalent to the fact that the grading on $A$ is non-negative and that the map $\sfe\to A^0$
is an isomorphism.

\medskip

The category $\QCoh(Z/\BG_m)$ consists of complexes $M$ of graded $A$-modules. 
The subcategory $\Coh(Z/\BG_m)\subset \QCoh(Z/\BG_m)$ corresponds to the condition
that $M$ is cohomologically bounded and all $H^i(M)$ are finitely generated over $H^0(A)$.

\medskip

The functor 
$$\Gamma(Z/\BG_m,-):\Coh(Z/\BG_m)\to \Vect_\sfe$$
takes $M$ to its degree $0$ component $M^0$, which 
is finite-dimensional. 

\medskip

This implies the assertion of the lemma. 

\end{proof}

\ssec{A digression: ind-algebraic stacks} \label{ss:ind alg}

\sssec{}

Let $\CZ$ be a prestack. 

\medskip

We shall say that $\CZ$ is an \emph{ind-algebraic stack} if it is \emph{convergent} and 
for every $n$, the $n$th coconnective truncation $^{\leq n}\CZ$, can be written as 
\begin{equation} \label{e:presentation stack}
^{\leq n}\CZ\simeq \underset{i\in I}{\on{colim}}\, \CZ_{i,n},
\end{equation} 
where:

\begin{itemize}

\item Each $\CZ_{i,n}$ is a quasi-compact $n$-coconnective algebraic stack locally of finite type;

\item The category $I$ of indices is filtered;

\item The transition maps $\CZ_{i,n}\to \CZ_{j,n}$ are closed embedding.

\end{itemize} 

We claim:

\begin{lem}  \label{l:closed in ind}
Let $\CZ$ be an $n$-coconnective ind-algebraic stack. Then:

\medskip

\noindent{\em(a)} The maps $\CZ_{i,n}\to \CZ$ are closed embeddings. 

\medskip

\noindent{\em(b)} The family 
$$i\mapsto  (\CZ_{i,n}\to \CZ)$$
is cofinal in the category of $n$-coconnective algebraic quasi-compact stacks equipped with a closed
embedding into $\CZ$.
\end{lem}  

The proof is parallel to \cite[Lemma 1.3.6]{GR3}\footnote{The $n$-coconnectivity condition is important here:
we use it when we say that an $n$-coconnective quasi-compact algebraic stack can be written as a \emph{finite}  
colimit of affine schemes, sheafified in the \'etale/fppf topology.}.

\sssec{}

We now claim:

\begin{lem} \label{l:ind-alg as a quotient}
Let a prestack $\CZ$ be equal to the quotient $\CY/\sG$, where $\CY$ is an ind-scheme locally almost of finite type, 
and $\sG$ is an algebraic group. Then $\CZ$ is an ind-algebraic stack.
\end{lem} 

\begin{proof}

The convergence condition easily follows from the fact that both $\CY$ and $\on{pt}/\sG$ are convergent.
Thus, we may assume that $\CY$ is $n$-coconnective. We need to show that we can write $\CY$ as
a filtered colimit
$$\CY\simeq \underset{i\in I}{\on{colim}}\,  Y_i,$$
where:

\begin{itemize}

\item Each $Y_i$ is a quasi-compact scheme almost of finite type, stable under the $\sG$-action;

\item The transition maps $Y_i\to Y_j$ are closed embeddings, compatible with the action of $\sG$. 

\end{itemize} 

We will first show that such a presentation exists but without the condition that $Y_i$ be 
almost of finite type. 

\medskip

Recall (see \cite[Sect. 3.1.6]{GR3}) that for a map of schemes $f:Y\to Z$, it makes sense to consider the
closure of the image of $Y$ inside $Z$, to be denoted $\ol{\on{Im}(f)}$. This is the universal
closed subscheme of $Z$ for which there exists a factorization of $f$ as
$$Y\to \ol{\on{Im}(f)} \to Z.$$

Furthermore, if $i:Z\to Z'$
is a closed embedding and $f':=i\circ f$, then the natural map $$\ol{\on{Im}(f)}\to \ol{\on{Im}(f')}$$
is an isomorphism. In particular, it makes sense to talk about the closure of the image in the target 
that is an ind-scheme.

\medskip

Write $\CY$ as a filtered colimit of closed (but not necessarily $\sG$-invariant) subschemes 
$$\CY\simeq \underset{i\in I}{\on{colim}}\,  Y'_i.$$

Now, let $Y_i$ be the closure of the image of the map
$$\sG\times Y'_i\to \sG \times \CY\to \CY,$$
where the last arrow is the action map. By construction, the closed subschemes $Y_i$ are $\sG$-invariant,
and the resulting map
$$\underset{i\in I}{\on{colim}}\,  Y_i\to \CY$$
is an isomorphism. 

\medskip

Now, starting from the family of subschemes constructed above, we apply a $\sG$-equivariant version 
of \cite[Proposition 1.7.7]{GR3} (proved in {\it loc.cit.} Sect. 3.5.2) to produce a family that consists of 
$\sG$-equivariant schemes almost of finite type. 

\end{proof}

\sssec{}

As corollary, combining with \thmref{t:main 1 abs}, we obtain:

\begin{cor} \label{c:restr ind-alg}
Every connected component of $\bMaps(\Rep(\sG),\bH)$ is an ind-algebraic stack.
\end{cor}

\ssec{Mock-affineness and coarse moduli spaces} \label{ss:mock-affine}

\sssec{} \label{sss:mock-affine}

Let $\CZ$ be an algebraic stack. We shall say that $\CZ$ is \emph{mock-affine} if the functor of global sections
$$\Gamma(\CZ,-):\QCoh(\CZ)\to \Vect_\sfe$$
is t-exact.

\medskip

Clearly, $\CZ$ is mock-affine if and only if its underlying classical stack $^{\on{cl}}\CZ$ is mock-affine. 

\sssec{Example}  \label{sss:mock aff as quot}
Let $\CZ$ be of the form $Y/\sG$, where $Y$ is affine scheme and $\sG$ is a \emph{reductive} algebraic group. 
Then (assuming that $\sfe$ has characteristic zero) the stack $\CZ$ is mock-affine.

\sssec{}

Let $\CZ$ be an ind-algebraic stack. We shall say that $\CZ$ is mock-affine if $^{\on{cl}}\CZ$ admits a presentation
\eqref{e:presentation stack} whose terms are mock-affine. 

\medskip

By \lemref{l:closed in ind}, this is equivalent to requiring that for every algebraic stack $\CZ'$ equipped with a closed embedding $\CZ'\to \CZ$,
the stack $\CZ'$ is mock-affine.

\sssec{}

From \thmref{t:main 1 abs}, combined with \lemref{l:ind-alg as a quotient} and Example \ref{sss:mock aff as quot},
we obtain that each connected component of $\bMaps(\Rep(\sG),\bH)$ is mock-affine.

\sssec{}

Let $\CZ$ be a mock-affine algebraic stack. In particular, the $\sfe$-algebra
$$\Gamma(\CZ,\CO_\CZ)$$
is connective. 

\medskip

Further, for every $n$, 
$$\tau^{\geq -n}(\Gamma(\CZ,\CO_\CZ))\simeq \Gamma({}^{\leq n}\CZ,\CO_{^{\leq n}\CZ}).$$

\medskip

We define the coarse moduli space $\CZ^{\on{coarse}}$ of $\CZ$ to be the affine scheme 
$$\Spec(\Gamma(\CZ,\CO_\CZ)).$$

By construction, we have a canonical projection
$$\brr:\CZ\to \CZ^{\on{coarse}}.$$

\sssec{} \label{sss:map brr}

Let $\CZ$ be a mock-affine ind-algebraic stack. For every $n$ consider the $n$-coconnective ind-affine ind-scheme
$$^{\leq n}\CZ^{\on{coarse}}:=\underset{i}{\on{colim}}\, \Spec(\Gamma(\CZ_{i,n},\CO_{\CZ_{i,n}}))$$
for $^{\leq n}\CZ$ written as in \secref{e:presentation stack} (by \lemref{l:closed in ind}, this definition is independent of
the presentation). 

\medskip

We define the ind-affine ind-scheme $\CZ^{\on{coarse}}$ to be the 
convergent completion\footnote{See \cite[Chapter 2, Sect. 1.4.8]{GR1} for what this means.} 
of
\begin{equation} \label{e:course moduli conv}
\underset{n}{\on{colim}}\, {}^{\leq n}\CZ^{\on{coarse}}.
\end{equation} 

I.e., this is a convergent prestack whose value on eventually coconnective affine schemes is given by the colimit
\eqref{e:course moduli conv}. 

\medskip

We have a canonical projection 
$$\brr:\CZ\to \CZ^{\on{coarse}}.$$

\sssec{}

We claim:

\begin{lem} \label{l:coarse pt}
Let $\CZ$ be a mock-affine ind-algebraic stack satisfying:

\begin{itemize}

\item $\CZ$ is locally almost of finite type;

\item $^{\on{red}}\CZ$ is a mock-proper algebraic stack;

\item $^{\on{red}}\CZ$ is connected.

\end{itemize}

Then $\CZ^{\on{coarse}}$ has the following properties:

\begin{itemize}

\item It is locally almost of finite type;

\item $^{\on{red}}(\CZ^{\on{coarse}})\simeq \on{pt}$.

\end{itemize}

\end{lem} 

\begin{proof}

To prove that $\CZ^{\on{coarse}}$ is locally almost of finite type, it suffices to show that for every $n$, and a presentation of $^{\leq n}\CZ$ as in 
\secref{e:presentation stack}, the rings $\Gamma(\CZ_{i,n},\CO_{\CZ_{i,n}})$ are finite-dimensional over $\sfe$. However, this follows from the
mock-properness assumption. 

\medskip

This also implies that $^{\on{red}}(\CZ^{\on{coarse}})$ is Artinian, i.e., is the union of finite many copies of $\on{pt}$.
The connectedness assumption on $^{\on{red}}\CZ$ implies that there is only one copy.

\end{proof}

\ssec{Coarse moduli spaces for connected components of $\LocSys_\sG^{\on{restr}}(X)$} \label{ss:coarse for comps}

\sssec{}

Let $\bH$ be again a gentle Tannakian category. We will apply the discussion from \secref{ss:mock-affine} to $\CZ$ being a connected component of
$\bMaps(\Rep(\sG),\bH)$. 

\medskip

Note that $\CZ$ satisfies the conditions of \lemref{l:coarse pt} by the combination of Theorems \ref{t:main 1 abs} and \ref{t:mock proper}. 
In particular, we obtain that
$\CZ^{\on{coarse}}$ is an ind-affine ind-scheme locally almost of finite type, and $^{\on{red}}(\CZ^{\on{coarse}})\simeq \on{pt}$.

\medskip

We are now ready to state the main result of this subsection: 

\begin{mainthm} \label{t:coarse restr}
Let $\CZ$ being a connected component of $\bMaps(\Rep(\sG),\bH)$, and consider the corresponding
map
$$\brr:\CZ \to \CZ^{\on{coarse}}.$$ 
We have:

\smallskip

\noindent{\em(a)}
The map $\brr$ makes $\CZ$ into a relative algebraic stack over $\CZ^{\on{coarse}}$, i.e., the base change of $\brr$
by an affine scheme yields an algebraic stack. 

\smallskip

\noindent{\em(b)} The ind-scheme $\CZ^{\on{coarse}}$ is a formal affine scheme
(see Remark \ref{r:formal affine} for what this means).

\end{mainthm} 

The proof of \thmref{t:coarse restr} (for a general gentle Tannakian category $\bH$) will given in \secref{s:formal coarse}.
In the particular case when $\bH=\qLisse(X)$ when $X$ is a smooth and complete algebraic curve, a simpler 
argument will be given in \secref{ss:simple coarse}. 

\sssec{}

Our main application is when $\bH=\qLisse(X)$, and so $\bMaps(\Rep(\sG),\bH)=\LocSys_\sG^{\on{restr}}(X)$.

\medskip

Denote by $\LocSys_\sG^{\on{restr,coarse}}(X)$ the disjoint union of the formal affine schemes $\CZ^{\on{coarse}}$
over the connected components $\CZ$ of $\LocSys_\sG^{\on{restr}}(X)$, and consider the corresponding map
$$\brr:\LocSys_\sG^{\on{restr}}(X)\to \LocSys_\sG^{\on{restr,coarse}}(X).$$

In \secref{ss:coarse Betti again} we will show that in the Betti context, this map can be obtained as a formal completion of the map
$$\brr:\LocSys_\sG^{\on{Betti}}(X)\to \LocSys_\sG^{\on{Betti,coarse}}(X)$$
of \eqref{e:projection to coarse} at the disjoint union of $\sfe$-points of $\LocSys_\sG^{\on{Betti,coarse}}(X)$.

\sssec{}

For a connected component $\CZ$ of $\bMaps(\Rep(\sG),\bH)$, set
$$\CZ^{\on{rigid}}:=\CZ\underset{\bMaps(\Rep(\sG),\bH)}\times \bMaps(\Rep(\sG),\bH)^{\on{rigid}}.$$

A consequence of \thmref{t:coarse restr} of particular importance for the sequel is:

\begin{cor} \label{c:coarse restr bis}
The fiber product
$$\on{pt} \underset{\CZ^{\on{coarse}}}\times \CZ$$
is an algebraic stack\footnote{It follows automatically that it is quasi-compact and locally almost of finite type}.
\end{cor} 

From here we obtain:

\begin{cor} \label{c:coarse restr}
The fiber product
$$\on{pt} \underset{\CZ^{\on{coarse}}}\times \CZ^{\on{rigid}}$$
is an affine scheme.
\end{cor} 

(Indeed, it is easy to see that a prestack that is simultaneously an ind-affine
ind-scheme and an algebraic stack is actually an affine scheme.)

\begin{rem}
We emphasize that the assertion of \corref{c:coarse restr} (resp., \corref{c:coarse restr bis}) is that the corresponding
fiber products do \emph{not} have ind-directions. 

\medskip 

They may be non-reduced, but the point is that they are (locally) schemes, as opposed to 
formal schemes.

\end{rem}

\section{The formal coarse moduli space} \label{s:formal coarse}

This section is devoted to the proof of \thmref{t:coarse restr}, and we continue to assume that $\sG$ is reductive. 

\medskip

In the course of the proof we will encounter another fundamental feature of $\bMaps(\Rep(\sG),\bH)$
(\thmref{t:top fin gen}):

\medskip

Recall that at the classical level, when we can think of $\bMaps(\Rep(\sG),\bH)$ as the prestack of homomorphisms
$\sH\to \sG$, where $\sH$ is the pro-algebraic Tannakian group attached to $(\bH^\heartsuit,\oblv_\bH)$,
see \propref{p:Tannaka red cl}. The claim is that on each component of $\bMaps(\Rep(\sG),\bH)$, these homomorphisms
factor via a particular quotient of $\sH$ which is \emph{topologically finitely generated}. 

\ssec{The coarse moduli space in the Betti setting} \label{ss:coarse Betti again}

In this subsection we return to the context of \secref{ss:Betti}. 
We will illustrate what \thmref{t:coarse restr} says in this case. 

\sssec{} \label{sss:coarse Betti again}

Let $X$ be a compact connected CW complex. 

\medskip

Recall the setting of \secref{sss:coarse Betti}: we have the affine scheme $\LocSys^{\on{Betti,coarse}}_\sG(X)$ and a map
\begin{equation} \label{e:to coarse again}
\brr: \LocSys^{\on{Betti}}_\sG(X)\to \LocSys^{\on{Betti,coarse}}_\sG(X).
\end{equation} 

%

\medskip


Let
$$\LocSys^{\on{restr},\on{coarse}}_\sG(X)$$ be the disjoint union of formal completions of $\LocSys^{\on{Betti,coarse}}_\sG(X)$ at its
$\sfe$-points.

\medskip

Note that \corref{c:coarse fibers} can be reformulated as saying that we have a Cartesian diagram
\begin{equation} \label{e:cart coarse}
\CD
\LocSys^{\on{restr}}_\sG(X) @>>> \LocSys^{\on{Betti}}_\sG(X) \\
@VVV @VVV \\
\LocSys^{\on{restr},\on{coarse}}_\sG(X) @>>> \LocSys^{\on{Betti,coarse}}_\sG(X). 
\endCD
\end{equation} 

\sssec{} \label{sss:betti coarse clear}

For a fixed $\sigma\in \LocSys^{\on{Betti,coarse}}_\sG(X)$, let $\CZ_\sigma\subset \LocSys^{\on{restr}}_\sG(X)$
be the corresponding connected component of $\LocSys^{\on{restr}}_\sG(X)$. 

\medskip

It is clear from \eqref{e:cart coarse} that
\begin{equation} \label{e:coarse restr Betti}
(\CZ_\sigma)^{\on{coarse}}\simeq (\LocSys^{\on{Betti,coarse}}_\sG(X))^\wedge_\sigma,
\end{equation} 
where the right-hand side is the formal completion of $\LocSys^{\on{Betti,coarse}}_\sG(X)$ at $\sigma$.

\medskip

The isomorphism \eqref{e:coarse restr Betti} makes both assertions of \thmref{t:coarse restr} manifest.
Indeed, point (a) follows from the fact that the projection
$$\CZ_\sigma \overset{r}\to (\CZ_\sigma)^{\on{coarse}}$$
is a base change of the map \eqref{e:to coarse again}, while $\LocSys^{\on{Betti}}_\sG(X)$ is an algebraic stack. 

\ssec{Property W} 

The rest of this section is devoted to the proof of \thmref{t:coarse restr}.

\sssec{} \label{sss:formal coarse situation}

Let $\CZ$ be an \'etale stack of the form $\CZ^{\on{rigid}}/\sG$, where $\CZ^{\on{rigid}}$ is an ind-affine ind-scheme
and $\sG$ is a reductive group.

\medskip

Assume that $^{\on{red}}\CZ$ is connected and mock-proper, so that \lemref{l:coarse pt} applies. In particular, 
$^{\on{red}}\CZ^{\on{rigid}}$ has a unique closed $\sG$-orbit, which corresponds to a unique closed point of $\CZ$,
\begin{equation} \label{e:dist point}
\on{pt}\to \CZ.
\end{equation}

\medskip

Consider the corresponding map
$$\brr:\CZ\to \CZ^{\on{coarse}},$$
and the unique point
$$\on{pt}\to \CZ^{\on{coarse}}.$$

\sssec{}

We shall say that $\CZ$ has Property W if the prestack 
$$\CW:=\on{pt}\underset{\CZ^{\on{coarse}}}\times \CZ$$ is an algebraic stack
(as opposed to an ind-algebraic stack).

\medskip

This is equivalent to requiring that 
\begin{equation} \label{e:W rigid}
\CW^{\on{rigid}}:=\on{pt}\underset{\CZ^{\on{coarse}}}\times \CZ^{\on{rigid}}
\end{equation}
is an affine scheme (as opposed to an ind-affine ind-scheme).

\sssec{}

We claim:

\begin{lem} \label{l:property W}
The following conditions are equivalent:

\smallskip

\noindent{\em(i)} The map $\brr$ makes $\CZ$ into a relative algebraic stack over $\CZ^{\on{coarse}}$;

\smallskip

\noindent{\em(ii)} $\CZ$ has Property W.

\end{lem} 

\begin{proof}

Clearly, we have (i) $\Rightarrow$ (ii). For the opposite implication, it suffices to show that if
$\CZ^{\on{coarse}}$ is written as
$$\underset{i}{\on{colim}}\, \Spec(A_i)$$ 
with $A_i$ Artinian, then each 
$$\Spec(A_i)\underset{\CZ^{\on{coarse}}}\times \CZ^{\on{rigid}}$$
is an affine scheme. 

\medskip

However, since $A_i$ is Artinian, this would follow once we know that the further base change
\begin{equation} \label{e:A pt fiber}
\on{pt} \underset{\Spec(A_i)}\times \left(\Spec(A_i)\underset{\CZ^{\on{coarse}}}\times \CZ^{\on{rigid}}\right)
\end{equation} 
is an affine scheme (indeed, given an ind-scheme $\CY$ over $\Spec(A)$ with $A$ an Artinian ring, if 
$\CY\underset{\Spec(A)}\times \on{pt}$ is a scheme, then $\CY$ is a scheme). 
However the fiber product \eqref{e:A pt fiber} is the same as $\CW^{\on{rigid}}$. 

\end{proof}

\sssec{}

Here is how Property W will be used. Let $\CZ$ be as in \secref{sss:formal coarse situation}. 

\begin{prop}  \label{p:property W}
Assume that the ind-affine ind-scheme $\CZ^{\on{rigid}}$ is a formal affine scheme, and assume that $\CZ$ has property $W$. 
Then $\CZ^{\on{coarse}}$ is also a formal affine scheme.
\end{prop}

\begin{proof}

By \thmref{t:formal}, it suffices to show that the tangent space to $\CZ^{\on{coarse}}$ at its unique
closed point, viewed as an object of $\Vect_\sfe^{\geq 0}$, is finite-dimensional in each degree. For that it suffices to check
that the !-pullback of $T_{\on{pt}}(\CZ^{\on{coarse}})$ to $\CW^{\on{rigid}}$, viewed as an object of 
$\IndCoh(\CW^{\on{rigid}})$, is such that all its cohomologies are in 
$\Coh(\CW^{\on{rigid}})^\heartsuit$.

\medskip

We have a fiber sequence
$$T(\CW^{\on{rigid}})\to T(\CZ^{\on{rigid}})|_{\CW^{\on{rigid}}}\to T_{\on{pt}}(\CZ^{\on{coarse}})|_{\CW^{\on{rigid}}}.$$

The cohomologies of $T(\CW^{\on{rigid}})\in \IndCoh(\CW^{\on{rigid}})$ lie in $\Coh(\CW^{\on{rigid}})^\heartsuit$ 
because $\CW^{\on{rigid}}$ is an affine scheme (locally almost of finite type). Now, $T(\CZ^{\on{rigid}})|_{\CW^{\on{rigid}}}$ also has cohomologies
lying in $\Coh(\CW^{\on{rigid}})^\heartsuit$ since $\CZ^{\on{rigid}}$ is formal affine scheme.

\end{proof}

\sssec{} \label{sss:return to coarse restr 1}

Let us return to the setting of \thmref{t:coarse restr}. As in \secref{sss:term loc sys}, we will refer to 
points of $\bMaps(\Rep(\sG),\bH)$ as ``local systems".

\medskip

Let $\sigma$ be a semi-simple local system, and let $\CZ_\sigma$ be the sconnected 
component of $\bMaps(\Rep(\sG),\bH)$ corresponding to $\sigma$. 

\medskip

Combining \lemref{l:property W} and \propref{p:property W}, we obtain that in order to prove 
\thmref{t:coarse restr}, it suffices to show that $\CZ_\sigma$ has Property W. 

\ssec{Property A}

Let $\CZ$ be as in \secref{sss:formal coarse situation}. 

\sssec{}

We shall say that $\CZ$ has Property A if there exists a classical affine scheme $\Spec(A)$ (\emph{not} necessarily almost of finite type) 
and a map $$\brr_A:{}^{\on{cl}}\CZ\to \Spec(A)$$
such that the following holds:

\medskip

\noindent The classical prestack underlying the fiber product 
\begin{equation} \label{e:A fiber product}
\on{pt}\underset{\Spec(A)}\times {}^{\on{cl}}\CZ^{\on{rigid}},
\end{equation} 
where $\on{pt}\to \Spec(A)$ is the map
\begin{equation} \label{e:pt to A}
\on{pt}\overset{\text{\eqref{e:dist point}}}\longrightarrow \CZ^{\on{rigid}}\to \Spec(A),
\end{equation}
is a (classical) affine scheme (as opposed to an ind-affine ind-scheme). 

\sssec{}

We claim:

\begin{prop} \label{p:property A}
If $\CZ$ has property $A$, then it has property $W$.
\end{prop}

\sssec{} \label{sss:return to coarse restr 2}

Assuming \propref{p:property A} for a moment, we obtain that in order to prove 
\thmref{t:coarse restr}, it suffices to show that the prestack
$\CZ_\sigma$ as in \secref{sss:return to coarse restr 1} has Property A. 

\sssec{}

The rest of this subsection is devoted to the proof of \propref{p:property A}. 

\medskip

Let $\Spec(A)$ and $\brr_A$ be as above. First, we claim that we can extend $\brr_A$ to a map at the derived level,
$$\CZ\to \Spec(A)$$
which we will denote by the same symbol $\brr_A$.

\medskip

Indeed, with no restriction of generality, we can assume that $A$ is a classical polynomial algebra,
so the datum of $\brr_A$ amounts to a collection $\sG$-invariant elements in $\Gamma({}^{\on{cl}}\CZ^{\on{rigid}},\CO_{^{\on{cl}}\CZ^{\on{rigid}}})$
or $\Gamma(\CZ^{\on{rigid}},\CO_{\CZ^{\on{rigid}}})$ for the classical and derived versions of $\brr_A$, respectively.

\medskip

Now, since $\CZ^{\on{rigid}}$ is a formal affine scheme, the map 
$$\Gamma(\CZ^{\on{rigid}},\CO_{\CZ^{\on{rigid}}})\to \Gamma({}^{\on{cl}}\CZ^{\on{rigid}},\CO_{^{\on{cl}}\CZ^{\on{rigid}}})$$
is an isomorphism on $H^0$. Hence, so is the map
$$\Gamma(\CZ^{\on{rigid}}_\phi,\CO_{\CZ^{\on{rigid}}})^{\sG}\to 
\Gamma({}^{\on{cl}}\CZ^{\on{rigid}}_\phi,\CO_{^{\on{cl}}\CZ^{\on{rigid}}})^{\sG},$$
since $\sG$ is reductive. Hence every element can be lifted.

\sssec{} \label{sss:fiber derived bis}

We now claim that the fiber product 
$$\on{pt}\underset{\Spec(A)}\times \CZ^{\on{rigid}}$$
itself is an affine scheme.

\medskip

Indeed, its underlying classical prestack is a classical affine scheme, by assumption. 
Further, it has a connective co-representable deformation theory, because $\CZ^{\on{rigid}}$ has this
property. Hence, it is indeed an affine scheme by \cite[Theorem 18.1.0.1]{Lu3}.  

\sssec{} 

We are now ready to prove that $\CW^{\on{rigid}}$ is an affine scheme. 

\medskip

Note that for $\Spec(A)$ as above, 
the map 
$$\brr_A:\CZ\to \Spec(A)$$
canonically factors as
$$\CZ \overset{\brr}\to \CZ^{\on{coarse}}\to \Spec(A).$$

Let us base change these maps by \eqref{e:pt to A}. Thus, from $\brr$, we obtain a map
\begin{equation} \label{e:base changed coarse}
\on{pt}\underset{\Spec(A)}\times \CZ \to \on{pt}\underset{\Spec(A)}\times \CZ^{\on{coarse}}.
\end{equation}

The map \eqref{e:base changed coarse} realizes 
$\on{pt}\underset{\Spec(A)}\times \CZ^{\on{coarse}}$ 
as 
$$\left(\on{pt}\underset{\Spec(A)}\times \CZ\right){}^{\on{coarse}}.$$

\medskip

The left-hand side in \eqref{e:base changed coarse} is 
$$(\on{pt}\underset{\Spec(A)}\times \CZ^{\on{rigid}})/\sG,$$
and hence, by \secref{sss:fiber derived bis}, is a mock-affine \emph{algebraic stack}
(as opposed to ind-algebraic stack). 

\medskip

From here, we obtain that the right-hand side in 
\eqref{e:base changed coarse} is an affine \emph{scheme} (as opposed to ind-scheme). 
Therefore, the map
$$\on{pt}\underset{\Spec(A)}\times \CZ^{\on{rigid}} \to \on{pt}\underset{\Spec(A)}\times \CZ^{\on{coarse}}.$$
is a map between affine schemes. Hence, its further pullback with respect to
$$\on{pt}\to \on{pt}\underset{\Spec(A)}\times \CZ^{\on{coarse}}$$
is still an affine scheme. But the latter pullback is the prestack $\CW^{\on{rigid}}$ of \eqref{e:W rigid}.

\medskip

Thus, $\CW^{\on{rigid}}$ is an affine scheme, as required.  

\qed[\propref{p:property A}]

\ssec{A digression: the case of algebraic groups}  \label{ss:maps of alg groups again}

In this subsection we will establish a particular case of \thmref{t:coarse restr}. Namely, we will 
show that it holds for $\bH=\Rep(\sH)$, where $\sH$ be a (finite-dimensional) algebraic group. 

\sssec{}

Let $\sH$ be an affine algebraic group of finite type, and consider the prestacks
$$\bMaps_{\on{Grp}}(\sH,\sG) \text{ and } 
\bMaps_{\on{Grp}}(\sH,\sG)/\on{Ad}(\sG).$$

\medskip

In \propref{p:maps of groups} we have already established that $\bMaps_{\on{Grp}}(\sH,\sG)$
is an ind-affine ind-scheme. Furthermore, if $\sH$ is reductive, we know by \propref{p:red group hom}
that $\bMaps_{\on{Grp}}(\sH,\sG)$ is a disjoint union of (classical smooth) affine schemes.

%

\sssec{} \label{sss:Maps phi}

Choose a Levi splitting of $\sH$, i.e.,
$$\sH:=\sH_{\on{red}}\ltimes \sH_u.$$

We have a natural projection
$$\bMaps_{\on{Grp}}(\sH,\sG)\to \bMaps_{\on{Grp}}(\sH_{\on{red}},\sG).$$

Fix a point $\phi\in \bMaps_{\on{Grp}}(\sH_{\on{red}},\sG)$, and set
$$\bMaps_{\on{Grp}}(\sH,\sG)_\phi:=\bMaps_{\on{Grp}}(\sH,\sG)\underset{\bMaps_{\on{Grp}}(\sH_{\on{red}},\sG)}\times \{\phi\}.$$

This is an ind-scheme, equipped with an action of $\on{Stab}_{\sG}(\phi)$. 
We are going to exhibit $\bMaps_{\on{Grp}}(\sH,\sG)_\phi$ as the completion of an affine scheme along a Zariski closed subset,
such that the entire situation carries an action of $\on{Stab}_{\sG}(\phi)$. 

\sssec{} \label{sss:comps for maps of grps}

Note that $\bMaps_{\on{Grp}}(\sH,\sG)_\phi$ identifies with
$$\bMaps_{\on{Grp}}(\sH_u,\sG)^{\sH_{\on{red}}},$$
where $\sH_{\on{red}}$ acts on $\sG$ via $\phi$ and on $\sH_u$ by conjugation. 

\medskip

Consider the affine scheme
$$\bMaps_{\Lie}(\sh_u,\sg)$$
(see \secref{sss:maps of Lie algs} below), and its closed subscheme
$$\bMaps_{\Lie}(\sh_u,\sg)^{\sH_{\on{red}}}.$$

We have a naturally defined map
\begin{equation} \label{e:from Lie to grp}
\bMaps_{\on{Grp}}(\sH_u,\sG)^{\sH_{\on{red}}}\to \bMaps_{\Lie}(\sh_u,\sg)^{\sH_{\on{red}}}.
\end{equation} 

We claim:

\begin{prop} \label{p:from Lie to grp}
The map \eqref{e:from Lie to grp} realizes $\bMaps_{\on{Grp}}(\sH_u,\sG)^{\sH_{\on{red}}}$
as the formal completion of the affine scheme $\bMaps_{\Lie}(\sh_u,\sg)^{\sH_{\on{red}}}$ along
the closed subset $\bMaps_{\Lie}(\sh_u,\sg)^{\sH_{\on{red}}}_{\on{nilp.im.}}$  consisting of those maps
$$\sh_u\to \sg$$ 
whose image is contained in the nilpotent cone of $\sg$.
\end{prop}

\begin{proof}

We interpret 
$$\bMaps_{\on{Grp}}(\sH_u,\sG)$$
as $\bMaps(\Rep(\sG),\Rep(\sH_u))^{\on{rigid}}$, 
and $\bMaps_{\Lie}(\sh_u,\sg)$ as 
$$\bMaps(\Rep(\sG),\sh_u\mod)^{\on{rigid}},$$
see \propref{p:Lie maps}.

\medskip

The fact that \eqref{e:from Lie to grp} is an ind-closed embedding and a formal isomorphism follows now from
the fact that the restriction functor
$$\Rep(\sH_u)\to \sh_u\mod$$
is fully faithful, whose essential image consists of objects, all of whose
cohomologies are such that the action of $\sh_u$ on them is locally nilpotent. 

\medskip

This description also implies the stated description of the essential image at the reduced level.

\end{proof}

\begin{cor} \label{c:from Lie to grp conn}
The formal affine scheme $\bMaps_{\on{Grp}}(\sH,\sG)_\phi$ is connected.
\end{cor}

\begin{proof} 
The action of $\BG_m$ by dilations contracts $\bMaps_{\Lie}(\sh_u,\sg)^{\sH_{\on{red}}}$ to a single point,
and this action preserves the closed subset $\bMaps_{\Lie}(\sh_u,\sg)^{\sH_{\on{red}}}_{\on{nilp.im.}}$. 
\end{proof} 

\sssec{} \label{sss:coarse for grps}

Let $\CZ$ be a connected component of $\bMaps_{\on{Grp}}(\sH,\sG)/\on{Ad}(\sG)$. 
From \corref{c:from Lie to grp conn} we obtain that $\CZ$ has the form 
$$\bMaps_{\on{Grp}}(\sH,\sG)_\phi/\on{Ad}(\on{Stab}_{\sG}(\phi))$$
for some $\phi:\sH_{\on{red}}\to \sG$. Denote such $\CZ$ by $\CZ_\phi$. 

\medskip

Note that the unique closed point of $\CZ_\phi$ identifies with
$$\on{pt}\to \on{pt}/\on{Stab}_{\sG}(\phi) \hookrightarrow \bMaps_{\on{Grp}}(\sH_u,\sG)^{\sH_{\on{red}}}/\on{Ad}(\on{Stab}_{\sG}(\phi))
\simeq \bMaps_{\on{Grp}}(\sH,\sG)_\phi/\on{Ad}(\on{Stab}_{\sG}(\phi)),$$
where the middle arrow corresponds to the trivial homomorphism $\sH_u\to \sG$.

\medskip

In other other words, it corresponds to the locus of homomorphisms $\sH\to \sG$ that factor as 
$$\sH\twoheadrightarrow \sH_{\on{red}}\overset{\phi}\to \sG.$$

\sssec{} \label{sss:Chev}

Let us show that the stack $\CZ_\phi$ has Property A, thereby establishing that \thmref{t:coarse restr} holds
for $\bH=\Rep(\sH)$, see \secref{sss:return to coarse restr 2}. 

\medskip 

Let
$$\fa:=\sg/\!/\on{Ad}(\sG)$$
be the Chevalley space of $\sg$. This is an affine scheme equipped with an action of $\BG_m$. 

\medskip

We let $\Spec(A)$ be the affine scheme 
$$\bMaps_{\on{Sch}}(\sh_u,\fa)^{\BG_m},$$
where $\BG_m$ acts on $\sh_u$ by dilations, and on $\fa$ via its action on $\sg$ (also by dilations). 
It is easy to see that this is indeed an affine scheme. 

\medskip

We define map $\brr_A$ as the composition
\begin{multline*} 
\CZ_\phi\to \bMaps_{\on{Grp}}(\sH,\sG)/\on{Ad}(\sG) \to
\bMaps_{\on{Grp}}(\sH_u,\sG)/\on{Ad}(\sG)\to
\bMaps_{\on{Lie}}(\sh_u,\sg)/\on{Ad}(\sG) \to \\
\to \bMaps_{\on{Sch}}(\sh_u,\sg)^{\BG_m}/\on{Ad}(\sG) \to 
\bMaps_{\on{Sch}}(\sh_u,\fa)^{\BG_m}
\end{multline*} 

\medskip

Let us show that the fiber product 
$$\on{pt}\underset{\bMaps_{\on{Sch}}(\sh_u,\fa)^{\BG_m}}\times \CZ_\phi$$
is an algebraic stack (as opposed to an ind-algebroac stack). We will do so
right away at the derived level. 

\medskip

We will establish an equivalent fact, namely, that
\begin{equation} \label{e:fib prod phi}
\on{pt} \underset{\bMaps_{\on{Sch}}(\sh_u,\fa)^{\BG_m}}\times \bMaps_{\on{Grp}}(\sH,\sG)_\phi
\end{equation}
is an affine scheme.

\sssec{}

We rewrite \eqref{e:fib prod phi} as 
$$\on{pt} \underset{\bMaps_{\on{Sch}}(\sh_u,\fa)^{\BG_m}}\times \bMaps_{\on{Grp}}(\sH_u,\sG)^{\sH_{\on{red}}},$$
and consider the fiber product
$$\on{pt} \underset{\bMaps_{\on{Sch}}(\sh_u,\fa)^{\BG_m}}\times \bMaps_{\on{Lie}}(\sh_u,\sg)^{\sH_{\on{red}}},$$
which is an affine scheme, because $\bMaps_{\on{Lie}}(\sh_u,\sg)^{\sH_{\on{red}}}$ is such.

\medskip

Hence, it suffices to show that the map
\begin{equation} \label{e:Chev map u}
\on{pt} \underset{\bMaps_{\on{Sch}}(\sh_u,\fa)^{\BG_m}}\times \bMaps_{\on{Grp}}(\sH_u,\sG)^{\sH_{\on{red}}}\to
\on{pt} \underset{\bMaps_{\on{Sch}}(\sh_u,\fa)^{\BG_m}}\times \bMaps_{\on{Lie}}(\sh_u,\sg)^{\sH_{\on{red}}}
\end{equation}
is schematic. We claim that \eqref{e:Chev map u} is in fact an isomorphism.

\sssec{}

By \propref{p:from Lie to grp}, a priori, the map \eqref{e:Chev map u} realizes the left-hand side as the formal
completion of the right-hand side along the closed subset
\begin{equation} \label{e:Chev map u fiber}
\on{pt} \underset{\bMaps_{\on{Sch}}(\sh_u,\fa)^{\BG_m}}\times\bMaps_{\Lie}(\sh_u,\sg)^{\sH_{\on{red}}}_{\on{nilp.im.}}
\end{equation}
where 
$$\bMaps_{\Lie}(\sh_u,\sg)^{\sH_{\on{red}}}_{\on{nilp.im.}}\subset \bMaps_{\Lie}(\sh_u,\sg)^{\sH_{\on{red}}}$$
is the locus of maps 
\begin{equation} \label{e:unip map}
\sh_u\to \sg
\end{equation} 
whose image is contained in the nilpotent cone. 

\medskip

However, if a map \eqref{e:unip map} is such that the composition
$$\sh_u\to \sg\to \fa$$
is zero, then this map automatically lands in the nilpotent cone. 

\medskip

This implies that the closed subset \eqref{e:Chev map u fiber} is all of 
$$\on{pt} \underset{\bMaps_{\on{Sch}}(\sh_u,\fa)^{\BG_m}}\times \bMaps_{\on{Lie}}(\sh_u,\sg)^{\sH_{\on{red}}},$$
and hence \eqref{e:Chev map u} is an isomorphism. 

\ssec{The case of pro-algebraic groups} \label{ss:pro-alg again}

In this subsection we will study connected components of the ind-algebraic stack
$\bMaps_{\on{Grp}}(\sH,\sG)/\on{Ad}(\sG)$, where $\sH$ is a pro-algebraic group. 

\sssec{} \label{sss:pro-alg again}

Choose a Levi splitting
$$\sH\simeq \sH_{\on{red}} \ltimes \sH_u,$$
see \cite[Theorem 3.2]{HM}.  

\medskip

The description of connected components of $\bMaps_{\on{Grp}}(\sH,\sG)/\on{Ad}(\sG)$
in the case when $\sH$ is of finite type given in \secref{ss:maps of alg groups again} applies
verbatim to the present situation:

\medskip

The connected components are in bijection with conjugacy classes of homomorphisms
$\phi:\sH_{\on{red}}\to \sG$, and for a given $\phi$, the corresponding connected component
$\CZ_\phi$ identifies with 
$$\bMaps_{\on{Grp}}(\sH,\sG)_\phi/\on{Ad}(\on{Stab}_{\sG}(\phi))$$
and
$$\bMaps_{\on{Grp}}(\sH,\sG)_\phi\simeq \bMaps_{\on{Grp}}(\sH_u,\sG)^{\sH_{\on{red}}}.$$

However, we do not know, in general, whether such $\CZ_\phi$ satisfies Property A.

\begin{rem}
Note that in the above discussion, $\sH$ is an arbitrary pro-algebraic group, so it is \emph{not}
true, in general, that its category of representation $\Rep(\sH)$ is a gentle Tannakian category. 
Hence, it is \emph{not} true that the ind-affine ind-scheme 
$$\CZ_\phi^{\on{rigid}}:=\CZ_\phi\underset{\on{pt}/\sG}\times \on{pt}$$
is a formal affine scheme.
\end{rem}

\sssec{} 

Let $\on{Free}_n$ be the free group on $n$ letters, and let $\on{Free}^{\on{Pro-alg}}_n$ be its pro-algebraic envelope
over $\sfe$, i.e., 
\begin{equation} \label{e:free group}
\Hom_{\on{Grp}}(\on{Free}^{\on{Pro-alg}}_n,\sH')\simeq (\sH'(\sfe))^{\times n}, \quad \sH'\in \on{Alg. Groups}.
\end{equation}

\sssec{}

Let $\sH$ be a pro-algebraic group, written as $\underset{\alpha}{\on{lim}}\, \sH_\alpha$ with surjective transition maps. 
A map $\on{Free}^{\on{Pro-alg}}_n\to \sH$ is then the same as an $n$-tuple $\ul{g}$ of elements in $\sH(\sfe)$. 

\medskip

We shall say that an $n$-tuple $\ul{g}$ \emph{topologically generates} $\sH$ if the corresponding map
$\on{Free}^{\on{Pro-alg}}_n\to \sH$ is such that all the composite maps
$$\on{Free}^{\on{Pro-alg}}_n\to \sH\to \sH_\alpha$$
are surjective.

\medskip

This is equivalent to the condition that the Zariski closure of the abstract group generated by the images
of the elements of $\ul{g}$ in $\sH_\alpha$ is all of $\sH_\alpha$. 

\sssec{}  \label{sss:top fin gen}

We will say that $\sH$ is \emph{topologically finitely generated} if it admits a finite set of topological generators.

\sssec{}

We will prove:

\begin{thm} \label{t:coarse groups pro}
Assume that $\sH$ is topologically finitely generated. Then every connected component of 
$\bMaps_{\on{Grp}}(\sH,\sG)/\on{Ad}(\sG)$ has Property A.
\end{thm}

\begin{rem}
Note that since any algebraic group of finite type (over a field of characteristic $0$) is topologically finitely generated,
\thmref{t:coarse groups pro} provides an alternative proof of \thmref{t:coarse restr} for $\bH=\Rep(\sH)$, where 
$\sH$ is an algebraic group of finite type.
\end{rem}

\ssec{Proof of \thmref{t:coarse groups pro}}

\sssec{} \label{sss:surj coarse}

Write
$$\sH\simeq \underset{\alpha}{\on{lim}}\, \sH_\alpha.$$

Let $\sH'\to \sH$ be a homomorphism of pro-algebraic groups, such that for every $\alpha$
the composite map
$$\sH'\to \sH\to \sH_\alpha$$
is surjective. 

\medskip

We claim that if every connected component of $\bMaps_{\on{Grp}}(\sH',\sG)/\on{Ad}(\sG)$
has Property A, then so does every connected component of $\bMaps_{\on{Grp}}(\sH,\sG)/\on{Ad}(\sG)$
(for a given $\sG$). 

\sssec{}

Let $\CZ_\phi$ be a connected component of $\bMaps_{\on{Grp}}(\sH,\sG)/\on{Ad}(\sG)$ containing 
a given map $$\phi:\sH_{\on{red}}\to \sG.$$ 

The map $\sH'\to \sH$ induces a map 
\begin{equation} \label{e:maps H H'}
\bMaps_{\on{Grp}}(\sH,\sG)/\on{Ad}(\sG)\to \bMaps_{\on{Grp}}(\sH',\sG)/\on{Ad}(\sG).
\end{equation} 
The surjectivity property of the map of groups implies that \eqref{e:maps H H'} 
is a \emph{closed embedding}.

\medskip

Let $\CZ'_\phi$ be the connected component of $\bMaps_{\on{Grp}}(\sH',\sG)/\on{Ad}(\sG)$
containing the image of $\CZ_\phi$.
Since $\CZ'_\phi$ has Property A, we can find a map 
$$\brr'_A:\CZ'_\phi\to \Spec(A),$$
such that (the classical prestack underlying)
$\on{pt}\underset{\Spec(A)}\times \CZ'_\phi$
is an algebraic stack. 

\medskip

Define a map $\brr_A:\CZ_\phi\to \Spec(A)$
to be the composition
$$\CZ_\phi\to \CZ'_\phi\overset{\brr'_A}\longrightarrow\Spec(A).$$

Since the map
$$\on{pt}\underset{\Spec(A)}\times \CZ_\phi\to \on{pt}\underset{\Spec(A)}\times \CZ'_\phi$$
is a closed embedding, and we obtain that (the classical prestack underlying) 
$\on{pt}\underset{\Spec(A)}\times \CZ_\phi$ is also an algebraic stack, as required. 

\sssec{}

Thus, by the assumption on $\sH$ and \secref{sss:surj coarse}, we can replace the original $\sH$ by $\on{Free}^{\on{Pro-alg}}_n$. 

\medskip

Note that $$\Rep(\on{Free}^{\on{Pro-alg}}_n) \simeq \qLisse(X),$$
where $X$ is the bouquet of $n$ copies of $S^1$.

\medskip

Hence, the prestack $$\bMaps_{\on{Grp}}(\on{Free}^{\on{Pro-alg}}_n,\sG)/\on{Ad}(\sG)$$ is the same
as (the Betti version of) $\LocSys_\sG^{\on{restr}}(X)$. 

\medskip

Now, the fact that connected components of (the Betti version of) $\LocSys_\sG^{\on{restr}}(X)$
have Property A follows from \secref{ss:coarse Betti again}: we can take
$$\Spec(A):=\LocSys_\sG^{\on{Betti,coarse}}(X).$$

\qed[\thmref{t:coarse groups pro}]

\begin{rem}

Let us emphasize that the pro-algebraic group $\on{Free}^{\on{Pro-alg}}_n$ satisfies its universal
property \eqref{e:free group} for individual target groups $\sH'$, but \emph{not} in families. So, 
the ind-scheme prestack $\bMaps_{\on{Grp}}(\on{Free}^{\on{Pro-alg}}_n,\sG)$ is $\LocSys_\sG^{\on{restr},\on{rigid}_x}(X)$
(for $X$ the bouquet of $n$ copies of $S^1$), which is different from
$$\bMaps_{\on{Grp}}(\on{Free}_n,\sG)\simeq \LocSys_\sG^{\on{Betti,rigid}_x}(X)\simeq \sG^{\times n}.$$

\end{rem} 

\ssec{Proof of \thmref{t:coarse restr}}

\sssec{}

Let $\CZ_\sigma$ be a connected component of $\bMaps(\Rep(\sG),\bH)$. According to 
\secref{sss:return to coarse restr 2}, it suffices to show that $\CZ_\sigma$ has Property A. 

\medskip

Recall that, according to \propref{p:Tannaka red cl}, the prestack $^{\on{cl}}\bMaps(\Rep(\sG),\bH)$ identifies with the classical prestack
underlying
$$\bMaps_{\on{Grp}}(\sH,\sG)/\on{Ad}(\sG),$$
where $\sH$ is as in \secref{sss:Galois}. The local system $\sigma$ corresponds to the conjugacy class
of a homomorphism $\phi:\sH_{\on{red}}\to \sG$, so that
$$^{\on{cl}}\CZ_\sigma \simeq {}^{\on{cl}}\CZ_\phi.$$

Hence, it suffices to show that $\CZ_\phi$ has Property A, 

%

\begin{rem}
If we knew that $\sH$ is topologically finitely generated, then the fact that $\CZ_\phi$ has Property A
would follow from \thmref{t:coarse groups pro}. 

\medskip 

However, we do not know whether $\sH$ is topologically finitely generated. Instead, we will show that 
for every $\sigma$, there exists a particular quotient of (the unipotent part of) $\sH$ that is topologically
finitely generated, such that the passage to this quotient does not change $^{\on{cl}}\CZ_\phi$. This will effectively
reduce us to the situation of \thmref{t:coarse groups pro}. 

\end{rem}

\sssec{}

Let 
$$\bMaps_{\on{Grp}}(\sH,\sG)_\phi\simeq
\bMaps_{\on{Grp}}(\sH_u,\sG)^{\sH_{\on{red}}}$$
be as in \secref{sss:pro-alg again}. 

\medskip

Being a pro-unipotent group, we can write $\sH_u$ as
$$\underset{\alpha}{\on{lim}}\, \sH_\alpha,$$
where $\alpha$ runs over a filtered family of indices, the groups $\sH_\alpha$ are
finite-dimensional and unipotent and the transition maps 
$$\sH_{\alpha_2}\to \sH_{\alpha_1}$$
are surjective. 

\medskip

With no restriction of generality, we can assume that the $\sH_{\on{red}}$-action on
the pro-algebraic group $\sH_u$ comes from a compatible family of actions on the $\sH_\alpha$'s. 

\medskip

We have:
$$\bMaps_{\on{Grp}}(\sH_u,\sG)^{\sH_{\on{red}}}\simeq
\underset{\alpha}{\on{colim}}\, \bMaps_{\on{Grp}}(\sH_\alpha,\sG)^{\sH_{\on{red}}}.$$

\sssec{}

For each index $\alpha$, let $\sh_{\alpha,\phi\on{-isotyp}}$ be the maximal Lie algebra quotient of
$$\sh_{u,\alpha}:=\on{Lie}(\sH_\alpha)$$
on which the action of $\sH_{\on{red}}$ has only the same isotypic components as those that
appear in $\sg:=\on{Lie}(\sG)$, where the latter is acted on by $\sH_{\on{red}}$ via $\phi$.

\medskip

Let $\sH_{\alpha,\phi\on{-isotyp}}$ denote the corresponding quotient of $\sH_\alpha$. 

\begin{lem} \label{l:isotyp matters}
The map
\begin{equation} \label{e:map isotyp}
\bMaps_{\on{Grp}}(\sH_{\alpha,\phi\on{-isotyp}},\sG)^{\sH_{\on{red}}}\to 
\bMaps_{\on{Grp}}(\sH_\alpha,\sG)^{\sH_{\on{red}}}
\end{equation}
induces an isomorphism of the underlying classical prestacks. 
\end{lem}

\begin{proof}

Follows from \propref{p:from Lie to grp}.

\end{proof}

\sssec{}

Set
$$\sH_{\phi\on{-isotyp}}:=\underset{\alpha}{\on{lim}}\, \sH_{\alpha,\phi\on{-isotyp}}.$$

From \lemref{l:isotyp matters} we obtain that the map
$$\bMaps_{\on{Grp}}(\sH_{\phi\on{-isotyp}},\sG)^{\sH_{\on{red}}}\to 
\bMaps_{\on{Grp}}(\sH_u,\sG)^{\sH_{\on{red}}}$$
induces an isomorphism of the underlying classical prestacks. 

\medskip

Hence, it is sufficient to show that 
$$\bMaps_{\on{Grp}}(\sH_{\phi\on{-isotyp}},\sG)^{\sH_{\on{red}}}/\on{Ad}(\on{Stab}_{\sG}(\phi))$$
has Property A. 

\sssec{}

Consider the maps
$$\bMaps_{\on{Grp}}(\sH_{\phi\on{-isotyp}},\sG)^{\sH_{\on{red}}}/\on{Ad}(\on{Stab}_{\sG}(\phi)) \to
\bMaps_{\on{Grp}}(\sH_{\phi\on{-isotyp}},\sG)/\on{Ad}(\on{Stab}_{\sG}(\phi)) \to$$
$$\to \bMaps_{\on{Grp}}(\sH_{\phi\on{-isotyp}},\sG)/\on{Ad}(\sG).$$

In the above composition, the first map is a closed embedding, and the second map is schematic.

\medskip

Hence, it is sufficient to show that every connected component of 
$$\bMaps_{\on{Grp}}(\sH_{\phi\on{-isotyp}},\sG)/\on{Ad}(\sG)$$
has Property A.

\medskip

This follows by combining \thmref{t:coarse groups pro} and the following result: 

\begin{thm} \label{t:top fin gen}
The pro-algebraic group $\sH_{\phi\on{-isotyp}}$ is topologically finitely 
generated.
\end{thm}

\ssec{Proof of \thmref{t:top fin gen}}

\sssec{}

Let $\sH'$ be a pro-unipotent group
$$\sH'\simeq \underset{\beta}{\on{lim}}\,\sH'_\beta,$$
where $\sH'_\beta$ are unipotent algebraic groups of finite type. 

\medskip

Consider $\sh':=\on{Lie}(\sH')$ as a pro-finite dimensional vector space.
The following is elementary:

\begin{lem} \label{l:Lie top gener}
Let $\sV$ be a finite-dimensional subspace of $\sh'$ such that for every $\beta$, the image of $\sV$ in
$\sh'_\beta:=\on{Lie}(\sH'_\beta)$ generates it as a Lie algebra. Then $\sH'$ is topologically finitely
generated.
\end{lem} 

\sssec{}

Let $\sH'$ be a pro-unipotent group as above. We claim:

\begin{prop} \label{p:commut top gen}
Assume that $\sh'/[\sh',\sh']$ is finite-dimensional. Then $\sH'$ is topologically finitely
generated.
\end{prop} 

\begin{proof}

Let $\sV\subset \sh'$ be a finite-dimensional vector space that projects surjectively onto
$\sh'/[\sh',\sh']$. By \lemref{l:Lie top gener}, it suffices to see that for any $\beta$, the image
of $\sV$ in $\sh'_\beta$ generates it as a Lie algebra.

\medskip

But this follows from the next property of nilpotent Lie algebras: if a subspace $\wt\sV$ in a 
nilpotent finite-dimensional Lie algebra $\sh''$ projects surjectively onto $\sh''/[\sh'',\sh'']$,
then $\wt\sV$ generates $\sh''$ as a Lie algebra. 

\end{proof} 

\sssec{} \label{sss:estim generators}

We will prove \thmref{t:top fin gen} by applying \propref{p:commut top gen} to $\sH_{\phi\on{-isotyp}}$. 

\medskip

Note that the quotient
$$\sh_{\phi\on{-isotyp}}/[\sh_{\phi\on{-isotyp}},\sh_{\phi\on{-isotyp}}]$$ is the maximal pro-abelian quotient of 
$\on{Lie}(\sH_u)$ on which
$\sH_{\on{red}}$ acts via isotypic components that appear in its action on $\sg$ via $\sigma$. 

\medskip

Hence, it is enough to show that the vector space
$$\Hom\left(\on{Lie}(\sH_u)/[\on{Lie}(\sH_u),
\on{Lie}(\sH_u)],\sg\right)^{\sH_{\on{red}}}$$
is finite-dimensional.

\medskip

Note, however, that the above vector space is the same as
$$H^1\left(\on{Lie}(\sH_u),\sg\right)^{\sH_{\on{red}}},$$
which is the same as
$$H^1\left(\inv_{\sH_u}(\sg)\right)^{\sH_{\on{red}}}\simeq H^1\left(\inv_\sH(\sg)\right).$$

\sssec{}

We have 
$$H^1\left(\inv_\sH(\sg)\right) \simeq H^0\left(T_\phi(\bMaps_{\on{Grp}}(\sH,\sG)/\on{Ad}(\sG))\right),$$
and we have an exact triangle
$$\sg \to T_\phi(\bMaps_{\on{Grp}}(\sH,\sG))\to T_\phi(\bMaps_{\on{Grp}}(\sH,\sG)/\on{Ad}(\sG)),$$
hence it is enough to show that $H^0\left(T_\phi(\bMaps_{\on{Grp}}(\sH,\sG))\right)$ is finite-dimensional.

\sssec{}

We claim that we have a canonical isomorphism
\begin{equation} \label{e:tang space class}
H^0\left(T_\phi(\bMaps_{\on{Grp}}(\sH,\sG))\right) \simeq H^0\left(T_\sigma(\bMaps(\Rep(\sG),\bH)^{\on{rigid}})\right).
\end{equation} 

This would imply the finite-dimensionality claim by \corref{c:def}(b').

\sssec{}

To prove \eqref{e:tang space class}, we note that by \corref{c:def}(a), both
$$T^*_\phi(\bMaps_{\on{Grp}}(\sH,\sG)) \text{ and } T^*_\sigma(\bMaps(\Rep(\sG),\bH)^{\on{rigid}})$$
belong to $\Vect^{\leq 0}_\sfe$. 

\medskip

Hence, for $V\in \Vect^\heartsuit_\sfe$,
$$\Maps_{\Vect_\sfe}(T^*_\phi(\bMaps_{\on{Grp}}(\sH,\sG)),V)\simeq H^0\left(T_\phi(\bMaps_{\on{Grp}}(\sH,\sG))\right) \otimes V$$
and
$$\Maps_{\Vect_\sfe}(T^*_\sigma(\bMaps(\Rep(\sG),\bH)^{\on{rigid}}),V) \simeq H^0\left(T_\sigma(\bMaps(\Rep(\sG),\bH)^{\on{rigid}})\right) \otimes V.$$

\medskip

Now, by deformation theory
\begin{multline} \label{e:dual numbers 1}
\Maps_{\Vect_\sfe}(T^*_\phi(\bMaps_{\on{Grp}}(\sH,\sG)),V)\simeq \\
\simeq \Maps(\Spec(\sfe\oplus V),\bMaps_{\on{Grp}}(\sH,\sG))\underset{\Maps(\on{pt},\bMaps_{\on{Grp}}(\sH,\sG))}\times \{\phi\}
\end{multline} 
and
\begin{multline}  \label{e:dual numbers 2}
\Maps_{\Vect_\sfe}(T^*_\sigma(\bMaps(\Rep(\sG),\bH)^{\on{rigid}}),V) \simeq \\
\simeq \Maps(\Spec(\sfe\oplus V),\bMaps(\Rep(\sG),\bH)^{\on{rigid}}) \underset{\Maps(\on{pt},\bMaps(\Rep(\sG),\bH)^{\on{rigid}})}\times \{\sigma\},
\end{multline} 
where $\sfe\oplus V$ is a square-zero extension of $\sfe$ by means of $V$.

\medskip

However, since $V$ is classical, in the right-hand sides in \eqref{e:dual numbers 1} and \eqref{e:dual numbers 2} we can replace 
$$\bMaps_{\on{Grp}}(\sH,\sG) \text{ and } \bMaps(\Rep(\sG),\bH)^{\on{rigid}}$$
by
$$^{\on{cl}}\bMaps_{\on{Grp}}(\sH,\sG) \text{ and } ^{\on{cl}}\bMaps(\Rep(\sG),\bH)^{\on{rigid}},$$
respectively. Now, the assertion follows from the fact that
$$^{\on{cl}}\bMaps_{\on{Grp}}(\sH,\sG) \simeq {}^{\on{cl}}\bMaps(\Rep(\sG),\bH)^{\on{rigid}},$$
by \propref{p:Tannaka red cl}. 

\qed[\thmref{t:top fin gen}]

\section{Quasi-coherent sheaves on a formal affine scheme} \label{s:qcoh formal}

In this section we will study properties of the category of quasi-coherent sheaves on a formal
affine scheme, and then apply the results to 
$\QCoh(\bMaps(\Rep(\sG),\bH))$, where $\bH$ is a gentle Tannakian category. 

\medskip

The special feature of formal schemes among general ind-schemes is the following: for
an ind-scheme $\CY$, the category $\QCoh(\CY)$ is by definition the inverse limit of the categories 
$\QCoh(Y_i)$ for closed subschemes $Y_i\hookrightarrow \CY$. The functors in this inverse systems
are given by *-pullback and they do not generally admit left adjoints. So we do not in general know
whether $\QCoh(\CY)$ is compactly generated. 

\medskip

However, in the case of formal affine schemes, the situation is much better. 

\ssec{Formal affine schemes: basic properties} \label{ss:qcoh formal}

Let $\CY$ be a formal affine scheme. I.e., $\CY$ is a prestack that can be written as 
\begin{equation} \label{e:presentation A qcoh}
\underset{n\geq 1}{\on{colim}}\, \Spec(R_n)
\end{equation} 
as in \thmref{t:main 1}(d). I.e., $R_n$ are connective commutative $\sfe$-algebras of the form
$$R_n=R\underset{\sfe[t_1,...,t_m]}\otimes \sfe[t_1,...,t_m]/(t_1^n,...,t_m^n), \quad t_i\mapsto f_i\in R,\,\, i=1,....,m,$$
where $R$ is a connective commutative $\sfe$-algebra and $f_1,...,f_m$ is a collection of elements in $R$.

\medskip

Equivalently, we can write 
$$R_n=R\underset{\sfe[t_1,...,t_m]}\otimes \sfe, \quad t_i\mapsto f_i^n,$$

\medskip

In this subsection we will describe some favorable properties enjoyed by $\QCoh(\CY)$ for such $\CY$.
In general, $\QCoh$ of an ind-scheme is unwieldy, but \propref{p:ind vs pro} below allows one to get one's hand on 
$\QCoh(\CY)$ for $\CY$ a formal affine scheme. 

\sssec{} \label{sss:i n Y}

Fix a presentation of $\CY$ as in \eqref{e:presentation A qcoh}; denote by $i_\infty$ the resulting map $\CY\to \Spec(R)$. 
Set 
$$Y_n:=\Spec(R_n)\overset{i_n}\hookrightarrow \Spec(R).$$ 
For $n_1\leq n_2$, let $i_{n_1,n_2}$ denote the corresponding map $Y_{n_1}\to Y_{n_2}$. 

\medskip

Let $U\overset{j}\hookrightarrow \Spec(R)$ be the (open) complement of $\Spec(R_1)$. 

\begin{rem} \label{r:R is can}

Note that by the proof of \thmref{t:formal}, the $\sfe$-algebra $R$ and the map $i_\infty:\CY\to \Spec(R)$
can be constructed canonically starting from $\CY$, namely
$$R=\Gamma(\CY,\CO_\CY).$$ 

\medskip

However, there is a choice involved in choosing
the elements $f_1,...,f_n\in R$, and hence of the subschemes $Y_n$.

\end{rem} 

\sssec{}

Let 
$$\QCoh(\Spec(R))_\CY\overset{(i_\infty)_!}\hookrightarrow \QCoh(\Spec(R))$$ be the inclusion of the full subcategory consisting of objects 
with \emph{set-theoretic} support on $Y_1$ (i.e., these are objects whose restriction to $U$ vanishes). 
This inclusion 
admits a right adjoint, denoted $(i_\infty)^!$; explicitly, for every $\CF\in \QCoh(\Spec(R))$ we have the Cousin exact triangle
$$(i_\infty)_!\circ (i_\infty)^!(\CF)\to \CF\to j_*\circ j^*(\CF).$$

Furthermore, we can explicitly write the functor $(i_\infty)_!\circ (i_\infty)^!$ as
\begin{equation} \label{e:w supports as colimit} 
\underset{n}{\on{colim}}\, (i_n)_*\circ (i_n)^!,
\end{equation}
where we note that each $i_n^!$ is continuous because $i_n$ is a regular embedding.
(Note, however, that for fixed $n_1,n_2$, 
the functor $(i_{n_1,n_2})^!$, right adjoint to $(i_{n_1,n_2})_*$, is discontinuous.) 

\sssec{}

Consider the composite functor 
\begin{equation}   \label{e:ind vs pro}
\QCoh(\Spec(R))_\CY  \overset{(i_\infty)_!}\hookrightarrow \QCoh(\Spec(R)) \overset{(i_\infty)^*}\to \QCoh(\CY).
\end{equation} 

The following is established in \cite[Proposition 7.1.3]{GR3}:

\begin{prop}  \label{p:ind vs pro}
The functor \eqref{e:ind vs pro} is an equivalence.
\end{prop}

From here we formally obtain:

\begin{cor} \label{c:ind vs pro}  \hfill

\smallskip

\noindent{\em(a)} 
There exists a (unique) equivalence $\QCoh(\Spec(R))_\CY  \simeq  \QCoh(\CY)$, under which
the functor 
$$(i_\infty)^!: \QCoh(\Spec(R))\to \QCoh(\Spec(R))_\CY$$
goes over to the functor 
$$(i_\infty)^*:\QCoh(\Spec(R))\to \QCoh(\CY).$$

\smallskip

\noindent{\em(b)} 
The functor $(i_\infty)^*$ realizes $\QCoh(\CY)$ both as a co-localization and a localization of $\QCoh(\Spec(R))$ with respect to
the essential image of $\QCoh(U)$ along $j_*$.

\end{cor} 

\sssec{}

We observe:

\begin{lem} \label{l:gen of categ supp 1}  Let $\CY$ and $Y_n$ be as above.

\smallskip

\noindent{\em(a)} For $\CF\in \QCoh(\Spec(R))_\CY$, the map
$$\underset{n}{\on{colim}}\, (i_n)_*\circ (i_n)^!(\CF)\to \CF$$
is an isomorphism. 

\smallskip

\noindent{\em(b)} 
The category $\QCoh(\Spec(R))_\CY$ is compactly generated by the objects $(i_n)_*(\CO_{Y_n})$. 

\smallskip

\noindent{\em(c)} The subcategory of compact
objects in $\QCoh(\Spec(R))_\CY$ is closed under the monoidal operation. 

\end{lem} 

\begin{proof}

Point (a) follows from \eqref{e:w supports as colimit}. 

\medskip

The fact that the objects $(i_n)_*(\CO_{Y_n})$ generate $\QCoh(\Spec(R))_\CY$ follows from point (a). 
The fact that they are compact follows from the fact that they are compact as objects of $\QCoh(\Spec(R))$. 
This proves point (b).

\medskip

The fact that the 
subcategory of compact objects is closed under the monoidal operation follows from the corresponding fact for $\QCoh(\Spec(R))$. 
This proves point (c).

\end{proof}

\sssec{} \label{sss:i n infty}

Let $i_{n,\infty}$ denote the map $Y_n\to \CY$. Note that by \corref{c:ind vs pro}, the 
functor 
$$(i_{n,\infty})_*:\QCoh(Y_n)\to \QCoh(\CY),$$
right adjoint to 
$$(i_{n,\infty})^*:\QCoh(\CY)\to \QCoh(Y_n),$$ identifies
with $(i_\infty)^*\circ (i_n)_*$; in particular, it is continuous. 

\medskip

Furthermore, the above functor $(i_{n,\infty})_*$ admits a right adjoint,
to be denoted $(i_{n,\infty})^!$, which under the equivalence of 
\eqref{e:ind vs pro} corresponds to
$$(i_n)^!:\QCoh(\Spec(R))_\CY\to \QCoh(Y_n).$$

\medskip

Hence, from \lemref{l:gen of categ supp 1}, we obtain:

\begin{cor}  \label{c:colim An 1}  \hfill

\smallskip

\noindent{\em(a)} For $\CF\in \QCoh(\CY)$, the map
$$\underset{n}{\on{colim}}\, (i_{n,\infty})_*\circ (i_{n,\infty})^!(\CF)\to \CF$$
is an isomorphism. 

\smallskip

\noindent{\em(b)} 
The category $\QCoh(\CY)$ is compactly generated by the objects $(i_{n,\infty})_*(\CO_{Y_n})$. 

\smallskip

\noindent{\em(c)} The subcategory of compact
objects in $\QCoh(\CY)$ is closed under the monoidal operation. 

\end{cor}

\sssec{}

Finally, we claim:

\begin{prop} \label{p:gen of categ supp 2} 
The functor
\begin{equation} \label{e:colim An}
\underset{n}{\on{colim}}\, \QCoh(Y_n) \to \QCoh(\Spec(R))_\CY,
\end{equation} 
given by $\{(i_n)_*\}$, is an equivalence.
\end{prop}

Combining with \propref{p:ind vs pro}, we obtain: 

\begin{cor}  \label{c:colim An 2}  
The functor
$$\underset{n}{\on{colim}}\, \QCoh(Y_n) \to \QCoh(\CY),$$
given by $\{(i_{n,\infty})_*\}$, is an equivalence.
\end{cor}

\ssec{Proof of \propref{p:gen of categ supp 2}}

\sssec{}

For an index $n_0$, let 
$$\on{ins}_{n_0}:\QCoh(Y_{n_0})\to \underset{n}{\on{colim}}\, \QCoh(Y_n)$$
denote the corresponding tautological functor.

\medskip

For any object 
$$\CF\in \underset{n}{\on{colim}}\, \QCoh(Y_n),$$
we have a tautological isomorphism
\begin{equation} \label{e:object as evals}
\CF\simeq \underset{n}{\on{colim}}\, \on{ins}_n\circ (\on{ins}_n)^R(\CF).
\end{equation} 

\sssec{}

Denote the functor \eqref{e:colim An} by $\Psi$ and its right adjoint by $\Phi$
(note that we do not yet know that $\Phi$ is continuous).

\medskip

Let us rewrite $\underset{n}{\on{colim}}\, \QCoh(Y_n)$
as
\begin{equation} \label{e:QCoh as lim Y}
\underset{n}{\on{lim}}\,  \QCoh(Y_n),
\end{equation} 
where the limit is formed using the \emph{discontinuous} functors 
$$(i_{n_1,n_2})^!:\QCoh(Y_{n_2})\to \QCoh(Y_{n_1}),$$
see \cite[Chapter 1, Proposition 2.5.7]{GR1}.

\medskip

In terms of \eqref{e:QCoh as lim Y}, the functor $\Phi$ is given by the compatible collection of functors
$$\{(i_n)^!\}:\QCoh(\Spec(R)) \to \underset{n}{\on{lim}}\,  \QCoh(Y_n),$$
precomposed with $\QCoh(\Spec(R))_\CY\hookrightarrow \QCoh(\Spec(R))$. 

\medskip

In other words,
$$(\on{ins}_n)^R\circ \Phi\simeq (i_n)^!.$$

\sssec{}

Using \eqref{e:object as evals}, we obtain that the composition $\Psi\circ \Phi$ identifies
with the functor \eqref{e:w supports as colimit}. Hence, the counit of the adjunction
$$\Psi\circ \Phi\to \on{Id}$$
is an isomorphism, by \lemref{l:gen of categ supp 1}(a). 

\medskip

Hence, $\Phi$ is fully faithful. 

\sssec{}

We now show that the essential image of $\Phi$ generates the colimit category. 
It suffices to show that for every fixed $n_0$, and $\CF_0\in \QCoh(Y_{n_0})$
the object $\on{ins}_{n_0}(\CF_0)$ lies in the essential 
image of $\Phi$. We will show that
\begin{equation} \label{e:ident colim n0}
\on{ins}_{n_0}(\CF_0) \simeq \Phi\circ (i_{n_0})_*(\CF_0).
\end{equation}

Using \eqref{e:object as evals}, the desired isomorphism \eqref{e:ident colim n0} translates as 
\begin{equation} \label{e:ident colim n1}
\on{ins}_{n_0}(\CF_0) \simeq 
\underset{n\geq n_0}{\on{colim}}\, \on{ins}_n\circ i_n^! \circ (i_{n_0})_*(\CF_0).
\end{equation}

Consider yet another object: 
\begin{equation} \label{e:ident colim n2}
\underset{n\geq n_0}{\on{colim}}\, \underset{N\geq n}{\on{colim}}\, \on{ins}_n\circ i_{n,N}^! \circ (i_{n_0,N})_*(\CF_0).
\end{equation}

We will show that \eqref{e:ident colim n2} is isomorphic both to the left-hand side and the right-hand side 
of \eqref{e:ident colim n1}.

\sssec{}

The isomorphism with the left-hand side follows by replacing the index category in \eqref{e:ident colim n2} 
by the cofinal category with $N=n$. 

\sssec{}

For the isomorphism with the right-hand side, we will show that for every fixed $n\geq n_0$, the natural map 
\begin{equation} \label{e:ident colim n3}
\underset{N\geq n}{\on{colim}}\, i_{n,N}^! \circ (i_{n_0,N})_*(\CF_0) \to
i_n^! \circ (i_{n_0})_*(\CF_0),
\end{equation}
is an isomorphism (taking place in $\QCoh(Y_n)$). 

\medskip

Since we are dealing with affine schemes, it suffices to show that the isomorphism takes place 
at the level of global sections. By base change, the latter is equivalent to the fact that the map
$$\underset{N\geq n}{\on{colim}}\, 
\CHom_{\QCoh(Y_{n_0})}(\CO_{Y_{n_0}\underset{Y_N}\times Y_n},\CF_0)
\to \CHom_{\QCoh(Y_{n_0})}(\CO_{Y_{n_0}\underset{\Spec(R)}\times Y_n},\CF_0)$$
is an isomorphism in $\Vect_\sfe$. This follows from the next assertion:

\begin{lem}
The map from $\CO_{Y_{n_0}\underset{\Spec(R)}\times Y_n}$ to 
$$M\mapsto \CO_{Y_{n_0}\underset{Y_N}\times Y_n},$$
as a pro-object of $\QCoh(Y_{n_0})$, is an isomorphism.
\end{lem} 

\begin{proof}

The assertion immediately reduces to the case when $R=\sfe[t_1,...,t_m]$, and further to the
case when $m=1$. In this case, it becomes a calculation similar to \cite[Lemma 7.1.5]{GR3}.

\end{proof}

\qed[\propref{p:gen of categ supp 2}]

%
%

\ssec{Mapping affine schemes into a formal affine scheme}

Let $\CY$ be a formal affine scheme. 

\sssec{}

First, we notice:

\begin{lem} \label{l:aff diag}
The diagonal map $\Delta_\CY:\CY\to \CY\times \CY$ is affine.
\end{lem}

\begin{proof}
Fix a presentation of $\CY$ as in \eqref{e:presentation A qcoh}.  Then the map $\Delta_\CY$ can be obtained
as the base change of the diagonal map $\Delta_{\Spec(R)}:\Spec(R)\to \Spec(R)\times \Spec(R)$, i.e.,
the square
\begin{equation} \label{e:base change diag}
\CD
\CY  @>{\Delta_\CY}>>   \CY\times \CY  \\
@VVV   @VVV  \\
\Spec(R)  @>{\Delta_{\Spec(R)}}>>  \Spec(R)\times \Spec(R)
\endCD
\end{equation} 
is Cartesian.
\end{proof} 

\sssec{}  \label{sss:map from affine}

Let $S$ be an affine scheme, equipped with a map $f$ to $\CY$. Note that $f$ is \emph{affine} as a map of prestacks
(by \lemref{l:aff diag}). Hence, the functor $f_*$, right adjoint to $f^*$ is continuous. 



\begin{cor} \label{c:on pro as colim}
The functor
\begin{equation} \label{e:colim S}
\underset{(S,f)}{\on{colim}}\, \QCoh(S) \to \QCoh(\CY),
\end{equation} 
is an equivalence, where:

\begin{itemize}

\item The index category is either of the following:
$$\affSch_{/\CY},\,\, \affSch_{/\CY,\on{closed}},$$ where the subscript ``closed" indicates that we consider only closed 
embeddings\footnote{When $\CY$ locally almost of finite
type as a prestack, we can further allow $(\affSch_{\on{aft}/\sfe})_{/\CY}$ and $(\affSch_{\on{aft}/\sfe})_{/\CY,\on{closed}}$ 
as index categories in the above colimit.}
$S\to \CY$;

\item The colimit is formed using the pushforward functors $(f_{1,2})_*:\QCoh(S_1)\to \QCoh(S_2)$ for 
$$f_{1,2}:S_1\to S_2, \quad f_2\circ f_{1,2}=f_1.$$

\item The map in \eqref{e:colim S} is given by $\{\QCoh(S)\overset{f_*}\to \QCoh(\CY)\}$.

\end{itemize} 

\end{cor} 

\begin{proof}
Fix a presentation of $\CY$ as in \eqref{e:presentation A}. The assertion follows from \corref{c:colim An 2}
and the fact that the family $Y_n\overset{i_{n,\infty}}\longrightarrow \CY$ is cofinal in any of the above categories.
\end{proof} 

\sssec{} \label{sss:regular into formal}

Let $S\overset{f}\to \CY$ be as above. We shall say that $f$
is a \emph{regular closed embedding} if there exists a map $\CY\to \BA^m$, and a Cartesian diagram
$$
\CD
S @>{f}>> \CY \\
@VVV @VVV \\
\on{pt} @>{\bi_0}>> \BA^{m}. 
\endCD
$$

In this case, we have
$$\QCoh(S)\simeq \QCoh(\on{pt})\underset{\QCoh(\BA^m)}\otimes \QCoh(\CY).$$

Therefore, the adjunction $((\bi_0)_*,(\bi_0)^!)$ implies that the right adjoint $f^!$ of $f_*$
is continuous and is strictly compatible with the $\QCoh(\CY)$-actions.  In particular, the 
functor $f_*$ preserves compactness.

\medskip

Furthermore, the isomorphism
$$(\bi_0)^!\simeq  (\bi_0)^*[-m]$$ 
implies that we have an isomorphism 
\begin{equation} \label{e:upper ! regular}
f^!(\CF)\simeq f^*(\CF)[-m]. 
\end{equation} 

\sssec{} \label{sss:regular into formal bis}

Let $\CY$ be realized as
$$\Spec(R)^\wedge_{\Spec({}^{\on{cl}}R/I)},$$
where $I\subset {}^{\on{cl}}R$ is a finitely generated ideal. Let us be given a map $\Spec(R)\to  \BA^m$,
such that 
$$^{\on{red}}(\on{pt}\underset{\BA^m}\times \Spec(R))={}^{\on{red}}\CY$$
as subsets of $^{\on{red}}\Spec(R)$.

\medskip

Then 
\begin{equation} \label{e:R and Y BC}
\on{pt}\underset{\BA^m}\times \CY\to \on{pt}\underset{\BA^m}\times \Spec(R)
\end{equation}
is an isomorphism. Indeed, the left-hand side in \eqref{e:R and Y BC} is a priori the completion of the right-hand side along 
a closed subset, which is actually the whole thing. 

\medskip

In particular, we obtain that in the above situation, we have a regular closed embedding
$$\on{pt}\underset{\BA^m}\times \Spec(R)\to \CY.$$

\sssec{} \label{sss:regular into formal bis bis}

The situation of \secref{sss:regular into formal bis} is realized for $\CY$ written as in \eqref{e:presentation A qcoh},
with the map $\Spec(R)\to  \BA^m$ given by the $m$-tuple $(t_1,....,t_m)\in R$.

\medskip

Hence, we obtain that the maps $i_{n,\infty}:Y_n\to \CY$ of \secref{sss:i n Y}
are regular closed embeddings.

%
%
%
%

\ssec{Semi-rigidity and semi-passable prestacks} \label{ss:formalq-rigid}

\sssec{}

In \secref{s:semi-rigid} we introduce the notion of \emph{semi-rigid} symmetric monoidal category. We observe:

\begin{lem} \label{l:QCoh s-rigid}
Let $\CY$ be a prestack such that:

\smallskip

\noindent{\em(i)} The diagonal morphism $\Delta_{\CY}:\CY\to \CY\times \CY$ is schematic\footnote{In this paper, all schemes are
assumed quasi-separated and quasi-compact.};

\smallskip

\noindent{\em(ii)} The category $\QCoh(\CY)$ is dualizable.

\medskip

Then $\QCoh(\CY)$ is semi-rigid.

\end{lem}

\begin{proof}

If $\QCoh(\CY)$ is dualizable, for any prestack $\CZ$, the functor
$$\QCoh(\CY)\otimes \QCoh(\CZ)\to \QCoh(\CY\times \CZ)$$
is an equivalence, \cite[Chapter 3, Proposition 3.1.7]{GR1}. 

\medskip

In particular, the functor 
$$\QCoh(\CY)\otimes \QCoh(\CY)\to \QCoh(\CY\times \CY)$$
is an equivalence. Hence, we can identify the functor 
$$\on{mult}_{\QCoh(\CY)}:\QCoh(\CY)\otimes \QCoh(\CY)\to \QCoh(\CY)$$
with $\Delta_\CY^*$. 

\medskip

Hence, from the fact that the diagonal morphism of $\CY$ is schematic, we obtain that the functor
$$\on{mult}_{\QCoh(\CY)}:\QCoh(\CY)\otimes \QCoh(\CY)\to \QCoh(\CY),$$
admits a continuous right adjoint, namely, $(\Delta_\CY)_*$, see \cite[Chapter 3, Proposition 2.2.2]{GR1}.
Moreover, by \cite[Chapter 3, Lemma 3.2.4]{GR1} the functor $(\Delta_\CY)_*$ satisfies the projection
formula; hence, the structure of right-lax compatibility on $(\on{mult}_{\QCoh(\CY)})^R$ 
with the $\QCoh(\CY)$-bimodule structure is strict.

\end{proof}

\sssec{}

Let us say call a prestack $\CY$ \emph{semi-passable}\footnote{This condition is slightly weaker 
than ``quasi-passable" from \cite[Sect. 1.5.7]{GKRV}} if it satisfies the assumptions of \lemref{l:QCoh s-rigid}
(cf. \cite[Chapter 3, Sect. 3.5.1]{GR1} for the choice of the terminology). 

\medskip

We obtain:

\begin{cor}
Let $\CY$ be a formal affine scheme. Then $\CY$ is semi-passable.
\end{cor}

\sssec{}

Semi-rigid categories enjoy some very favorable 2-categorical properties. For example,
a module category over a semi-rigid category is dualizable if and only if it is dualizable
as a plain DG category, see \lemref{l:dual module over s-rigid}. 

\ssec{Duality for semi-passable prestacks} \label{ss:self-duality passable}

%
%

\sssec{}

Recall that if $Y$ is an affine scheme, the functors
\begin{equation} \label{e:naive pairing}
\QCoh(Y)\otimes \QCoh(Y) \overset{\otimes}\to \QCoh(Y) \overset{\Gamma(Y,-)}\to \Vect_\sfe
\end{equation} 
and
\begin{equation} \label{e:naive unit}
\Vect_\sfe \overset{\sfe\mapsto \CO_Y}\longrightarrow \QCoh(Y) \overset{(\Delta_*)_Y}\longrightarrow 
\QCoh(Y\times Y)\simeq \QCoh(Y)\otimes \QCoh(Y)
\end{equation} 
define an identification
$$\QCoh(Y)\simeq \QCoh(Y)^\vee.$$

\sssec{}

Assume now that $\CY$ is a semi-passable prestack. In this case, $\Gamma(\CY,-)$ may be discontinuous
(this happens when $\CY$ is a formal affine scheme). So, \eqref{e:naive pairing} cannot serve as a counit
of a self-duality. Yet, we will see that \eqref{e:naive unit} does form the unit of a self-duality. 

\sssec{}

We claim:

\begin{prop} \label{p:self-duality for semi-pass}
Let $\CY$ be a semi-passable prestack. Then the object
$$(\Delta_\CY)_*(\CO_\CY)\in \QCoh(\CY\times \CY)\simeq \QCoh(\CY)\otimes \QCoh(\CY)$$
is the unit of a duality.
\end{prop}

\begin{proof}

This is a special case of \lemref{l:semi-rigid self dual}.

\end{proof} 

For future needs, we observe:

\begin{lem}   \label{l:dual of f*}
Let $\CY$ be a semi-passable prestack. Then 
for an affine scheme $S$ and a map $S\overset{f}\to \CY$, with respect to the above self-duality 
on $\QCoh(\CY)$ and the canonical self-duality on $\QCoh(S)$, the functor $f^*$ is the dual of the functor $f_*$.
\end{lem}

\begin{proof}

We need to establish an isomorphism
\begin{equation} \label{e:dual of f*}
(f\times \on{id}_\CY)^*\circ (\Delta_\CY)_*(\CO_\CY)\simeq (\on{id}_S\times f)_*\circ (\Delta_S)_*(\CO_S).
\end{equation} 

\medskip

Consider the Cartesian diagram
$$
\CD 
S @>{f}>> \CY \\
@V{\on{Graph}_f}VV @VV{\Delta_\CY}V \\
S\times \CY @>{f\times \on{id}_\CY}>> \CY\times \CY. 
\endCD
$$

Since the vertical arrows are schematic, we obtain a commutative diagram
\begin{equation} \label{e:BC}
\CD
\QCoh(S)   @<{f^*}<< \QCoh(\CY) \\
@V{(\on{Graph}_f)_*}VV  @VV{(\Delta_\CY)_*}V  \\
\QCoh(S\times \CY)  @<{(f\times \on{id}_\CY)^*}<<  \QCoh(\CY\times \CY). 
\endCD
\end{equation} 

\medskip

Evaluating the two circuits of \eqref{e:BC} on $\CO_\CY\in \QCoh(\CY)$, we obtain the desired isomorphism in \eqref{e:dual of f*}. 

\end{proof} 

\ssec{The functor of !-global sections} \label{ss:!-sect semi-pass}

In this subsection, we will let $\CY$ be a semi-passable prestack.

\medskip

As was mentioned above, for a formal affine scheme, the functor of global sections
$$\Gamma(\CY,-)=\CHom_{\QCoh(\CY)}(\CO_\CY,-), \quad \QCoh(\CY)\to \Vect_\sfe$$
is discontinuous.

\medskip

In this subsection we will introduce its substitute, denoted $\Gamma_!(\CY,-)$. 

\sssec{} \label{sss:!-sect semi-pass}

Let $\Gamma_!(\CY,-)$ denote the functor
$$\QCoh(\CY)\to \Vect_\sfe,$$
dual to the functor
$$\Vect_\sfe \overset{\CO_\CY}\longrightarrow \QCoh(\CY)$$
with respect to the self-duality
\begin{equation} \label{e:self-duality Y}
\QCoh(\CY)^\vee\simeq \QCoh(\CY)
\end{equation}
of \propref{p:self-duality for semi-pass}. 

\medskip

Note that according to \secref{sss:Gamma ! is lax mon} possesses a natural (non-unital) right-lax
symmetric monoidal structure. 

\begin{rem}

According to Remark \ref{r:Gamma ! as ren}, if $\QCoh(\CY)$ compactly generated, 
the functor $\Gamma_!(\CY,-)$ can  be characterized as follows:
it is the unique continuous functor
$$\QCoh(\CY)\to  \Vect_\sfe$$
that restricts to $\Gamma(\CY,-)$ on the subcategory of compact objects.

\end{rem}

\sssec{}

We claim:

\begin{prop} \label{p:counit self-duality for semi-pass}  The counit for the self-duality \eqref{e:self-duality Y} is given by 
\begin{equation} \label{e:counit self-duality for semi-pass}
\QCoh(\CY)\otimes \QCoh(\CY) \overset{\otimes}\to \QCoh(\CY) \overset{\Gamma_!(\CY,-)}\longrightarrow \Vect_\sfe.
\end{equation} 
\end{prop}

\begin{proof}

This is a special case of \lemref{l:semi-rigid self dual counit}. 

\end{proof} 

\sssec{}

Here is one more property of the functor $\Gamma_!(\CY,-)$:

\begin{prop}  \label{p:Gamma c via S}
For an affine scheme $S$ and a map $S\overset{f}\to \CY$, there is a canonical isomorphism
$$\Gamma_!(\CY,-)\circ f_*\simeq \Gamma(S,-): \QCoh(S)\to \Vect_\sfe.$$
\end{prop}

\begin{proof} 
By \lemref{l:dual of f*}, the functors dual to both sides identify with 
$$\Vect_\sfe\overset{\CO_S}\to \QCoh(S).$$
\end{proof} 

\begin{rem}
One can show that the isomorphism of \propref{p:Gamma c via S} is compatible
with the right-lax symmetric monoidal structures on both sides. 
\end{rem}

\ssec{The functor of !-global sections on a formal affine scheme}

In this subsection we specialize again to the case when $\CY$ is a formal affine scheme.

\sssec{}

Note that \propref{p:Gamma c via S} allows us to describe the functor $\Gamma_!(\CY,-)$ as follows: 
in terms of the presentation \eqref{e:colim S}, it corresponds to the compatible collection of functors
$$\Gamma(S,-):\QCoh(S)\to  \Vect_\sfe.$$

This functor should \emph{not} be confused with the \emph{discontinuous} functor
$$\Gamma(\CY,-):\QCoh(\CY)\to \Vect_\sfe,$$
corepresented by $\CO_\CY$.

\sssec{}  \label{sss:Gamma !}

For a choice of the presentation of $\CY$ as in \eqref{e:presentation A qcoh}, 
in terms of the identification $$\QCoh(\CY)\simeq \QCoh(\Spec(R))_\CY,$$ the functor 
$\Gamma_!(\CY,-)$ corresponds to the composition
\begin{equation} \label{e:Gamma ! via ambient bis}
\Gamma(\Spec(R),-)\circ (i_\infty)_!.
\end{equation} 

In other words,
\begin{equation} \label{e:Gamma ! via ambient}
\Gamma_!(\CY,-) \simeq \Gamma(\Spec(R),-)\circ ((i_\infty)^*)^L.
\end{equation} 

\begin{rem}
Note that the functor $(i_\infty)_!$ is (non-unital) symmetric monoidal. It is easy to see that
the isomorphism \eqref{e:Gamma ! via ambient bis} is compatible with the right-lax 
symmetric monoidal structures on both sides. 
\end{rem}

\sssec{}

Here is an explicit expression for the functor $\Gamma_!(\CY,-)$ in terms of \secref{sss:i n infty}:

\medskip

For $\CF\in \QCoh(\CY)$, we have
\begin{equation} \label{e:Gamma ! expl}
\Gamma_!(\CY,\CF) \simeq \underset{n}{\on{colim}}\, \Gamma(Y_n,i_{n,\infty}^!(\CF)),
\end{equation}
where we also note that
$$i_{n,\infty}^!(\CF)\simeq i_{n,\infty}^*(\CF)\otimes i_{n,\infty}^!(\CO_\CY),$$
by \eqref{e:upper ! regular}. 

\sssec{}

Finally, we claim:

\begin{prop} \label{p:t on Y}
The category $\QCoh(\CY)$ carries a t-structure, uniquely characterized by the requirement that the
functor $\Gamma_!(\CY,-)$ is t-exact. Furthermore, $\QCoh(\CY)$ is left-complete in this t-structure. 
\end{prop}

\begin{proof}
Choose a presentation as in \eqref{e:presentation A qcoh}. Then the assertion of the proposition follows 
from \propref{p:ind vs pro}:

\medskip

The corresponding t-structure
on $\QCoh(\Spec(R))_\CY$ is the unique one for which the functor $(i_\infty)_!$ is t-exact. 
\end{proof}

\ssec{Applications of 1-affineness}

\sssec{} \label{sss:1-affine formal}

Recall what it means for a prestack to be 1-affine, see \cite[Definition 1.3.7]{Ga2}.

\medskip

From \cite[Theorem 2.3.1]{Ga2}, we obtain:

\begin{prop} \label{p:formal affine 1-affine}
A formal affine scheme is 1-affine.
\end{prop}

\sssec{}

We also have (see \cite[Theorems 1.5.7 and 2.2.2]{Ga2}):

\begin{thm} \label{t:pt/G is 1-affine}
Let $\sG$ be an algebraic group. Then the (pre)stack $\on{pt}/\sG$ is 1-affine.
\end{thm}

\sssec{}

Here is the concrete meaning of \thmref{t:pt/G is 1-affine}. It says that the operations
$$\bC \mapsto \bC\underset{\Rep(\sG)}\otimes \Vect_\sfe, \quad \Rep(\sG)\mod \to \QCoh(\sG)\mod$$
and 
$$\bC'\mapsto (\bC')^\sG:=\on{Funct}_{\QCoh(\sG)\mod}(\Vect_\sfe,\bC'), \quad \QCoh(\sG)\mod \to \Rep(\sG)\mod$$
define mutually inverse equivalences of categories. 

\medskip

In the above formulas, we regard $\QCoh(\sG)$ as a monoidal DG category with respect to 
the operation \emph{convolution}, i.e., by means of taking pushfoward along the group law
$\sG\times \sG\to \sG$. 

\sssec{}

In particular, an object $\bC\in \Rep(\sG)\mod$ is dualizable (this is equivalent to being dualizable as a plain
DG category, since $\Rep(\sG)$ is rigid, see \cite[Chapter 1, Proposition 9.4.4]{GR1}) if and only if 
$$\bC':=\bC\underset{\Rep(\sG)}\otimes \Vect_\sfe$$ is dualizable as an object of $\QCoh(\sG)\mod$ (this is equivalent
to being dualizable as a plain DG category, see \cite[Proposition 1.4.5]{Ga2}). 

\medskip

As another consequence of \thmref{t:pt/G is 1-affine}, we obtain that the functor
\begin{equation} \label{e:de-eq}
\bC \mapsto \bC\underset{\Rep(\sG)}\otimes \Vect_\sfe, \quad \Rep(\sG)\mod \to \DGCat
\end{equation} 
is conservative. 

\sssec{}

Let $\CY'$ be a prestack acted on by $\sG$, and set $\CY=\CY'/\sG$. We can regard 
$\QCoh(\CY')$ as a $\QCoh(\sG)$-module category and $\QCoh(\CY)$ as a $\Rep(\sG)$-module category so that we have
$$\QCoh(\CY)\simeq \QCoh(\CY')^\sG$$
and 
$$\QCoh(\CY')\simeq \Vect_\sfe\underset{\Rep(\sG)}\otimes \QCoh(\CY).$$

From \thmref{t:pt/G is 1-affine}, we obtain:

\begin{cor} \label{c:descent do quot} 
For $\CY'$ and $\CY$ as above, we have:

\smallskip

\noindent{\em(a)} If $\QCoh(\CY')$ is dualizable, then so is $\QCoh(\CY)$. 

\smallskip

\noindent{\em(b)} If $\CY'$ is semi-passable, then so is $\CY$. 

\smallskip

\noindent{\em(c)} If $\CY'$ is 1-affine, then so is $\CY$. 

\end{cor}

As a particular case, we obtain:

\begin{cor} \label{c:descent do quot form aff} 
If $\CY$ is a prestack of the form $\CY'/\sG$, where $\CY'$ is a formal affine scheme, then:

\smallskip

\noindent{\em(a)} $\CY$ is semi-passable;

\smallskip

\noindent{\em(b)} $\CY$ is 1-affine.

\end{cor}

\sssec{}

Here is an application of 1-affinenness that we will need:

\begin{lem}  \label{l:base change QCoh Z}
Let $\CY$ be a 1-affine prestack, and let 
$$\CZ\to \CY\leftarrow \CZ'$$
be a diagram of prestacks. Assume that $\QCoh(\CZ')$ is dualizable as a $\QCoh(\CY)$-module. 
Then the functor 
$$\QCoh(\CZ)\underset{\QCoh(\CY)}\otimes \QCoh(\CZ')\to
\QCoh(\CZ\underset{\CY}\times \CZ')$$
is an equivalence.
\end{lem}

\begin{proof}

Write
$$\QCoh(\CZ)\simeq \underset{f:S\to \CZ}{\on{lim}}\, \QCoh(S),$$
where $S$ are affine schemes.  Since $\QCoh(\CZ')$ was assumed dualizable as a $\QCoh(\CY)$-module, the functor
$$\QCoh(\CZ)\underset{\QCoh(\CY)}\otimes \QCoh(\CZ') \to
\underset{f:S\to \CZ}{\on{lim}}\, \left(\QCoh(S)\underset{\QCoh(\CY)}\otimes \QCoh(\CZ')\right)$$
is an equivalence. 

\medskip

The functor
$$\QCoh(\CZ\underset{\CY}\times \CZ')\to \underset{f:S\to \CZ}{\on{lim}}\, \QCoh(S\underset{\CY}\times \CZ')$$
is tautologically an equivalence. 

\medskip

This reduces the assertion of the lemma to the case when $\CZ=S$ is an affine scheme. In this case, it follows
from \cite[Proposition 3.1.9]{Ga2}.

\end{proof}

\begin{cor} \label{c:base change QCoh Z}
Let $\CY$ be of the form $\CY'/\sG$, where $\CY'$ is a formal affine scheme. Then for a diagram of prestacks
$$\CZ\to \CY\leftarrow \CZ',$$
if either $\CZ$ or $\CZ'$ is dualizable as a plain DG category, then the functor
\begin{equation} \label{e:base change QCoh Z}
\QCoh(\CZ)\underset{\QCoh(\CY)}\otimes \QCoh(\CZ')\to
\QCoh(\CZ\underset{\CY}\times \CZ')
\end{equation}
is an equivalence.
\end{cor}

\begin{proof}

Follows by combining \lemref{l:base change QCoh Z}, \corref{c:descent do quot form aff}
and \lemref{l:dual module over s-rigid}.

\end{proof}

\ssec{Compact generation of $\QCoh(\bMaps(\Rep(\sG),\bH))$}  \label{ss:comp gen restr}

\sssec{} \label{sss:formal union quot}

Recall that the prestack $\CY:=\bMaps(\Rep(\sG),\bH)$ can be written as $\CY'/\sG$,
where $\CY'$ is a disjoint union of formal affine schemes, equipped with an action of an algebraic
group $\sG$. 

\medskip

The results of the preceding subsections apply to prestacks of this form as well. In particular,
such $\CY$ is semi-passable, 1-affine, and an analog of \corref{c:base change QCoh Z}
holds. 

\medskip

In particular, from \lemref{l:QCoh s-rigid}, we obtain:

\begin{cor} \label{c:LocSys semi-rigid}
The category $\QCoh(\bMaps(\Rep(\sG),\bH))$ is semi-rigid.
\end{cor}

\sssec{}

However, there is one property of $\QCoh(\bMaps(\Rep(\sG),\bH))$ that does not follow from the
preceding results, namely, that $\QCoh(\bMaps(\Rep(\sG),\bH))$ is compactly generated. The goal
of this subsection is to establish this.

\begin{rem} \label{r:non-lift action}

Let $\CY'$ be a formal affine scheme acted on by $\sG$, and set 
$\CY\simeq \CY'/\sG$. 

\medskip

Recall that according to Remark \ref{r:R is can} we have a canonical choice for an affine 
scheme $\Spec(R)$ such that $\CY'$ can be obtained as its formal completion.
By canonicity, $\sfe$-points of $\sG$ act on $R$. However, we are not
guaranteed to have an action on $R$ of $\sG$ as an algebraic group; this
is because the construction of $R$ involves the procedure of passing to the inverse limit. 

\medskip

Hence, it is not clear that we can find a $\sG$-equivariant model for a presentation of 
$\CY'$ as in \eqref{e:presentation A qcoh}.

\medskip

If we had such a presentation, we could give an easy proof of the fact that $\QCoh(\CY)$ is
compactly generated. 

\medskip

In the case when $\CY$ is a connected component of $\bMaps(\Rep(\sG),\bH)$,
we will take a different route, namely, one supplied by \thmref{t:coarse restr}.

\end{rem}

\sssec{} 

Let $\CZ$ be a connected component of $\bMaps(\Rep(\sG),\bH)$. Our current goal is to prove the following:

\begin{thm} \label{t:comp gen restr rough}
The category $\QCoh(\CZ)$ is compactly generated.
\end{thm} 

\begin{rem}
We will actually prove a slightly more precise version of \thmref{t:comp gen restr rough},
see \thmref{t:comp gen restr fine}, in which we will explicitly describe compact generators
of $\QCoh(\CZ)$.
\end{rem} 

\sssec{}  \label{sss:almost lift action a}

Consider the coarse moduli space $\CZ^{\on{coarse}}=:\CS$ corresponding to $\CZ$ and the map 
$$\brr:\CZ\to \CS,$$
see \secref{sss:map brr}. 

\medskip

Recall that according to \thmref{t:coarse restr}(b), the ind-scheme $\CS$ is actually a formal affine scheme. 
Write 
$$\CS\simeq \underset{n}{\on{colim}}\, \Spec(R_n)$$
as in \eqref{e:presentation A qcoh}. Denote by $i_{n,\infty}$ the corresponding maps 
$$\Spec(R_n)=:S_n\to \CS.$$

\medskip

By Sects. \ref{sss:regular into formal}-\ref{sss:regular into formal bis bis}, for every $n$, we have an adjunction
\begin{equation} \label{e:adj coarse}
(i_{n,\infty})_*:\QCoh(S_n)\rightleftarrows \QCoh(\CS):(i_{n,\infty})^!
\end{equation} 
as $\QCoh(\CS)$-module categories. Moreover, by \corref{c:colim An 1}(a), the map 
\begin{equation} \label{e:adj isom coarse}
\underset{n}{\on{colim}}\, (i_{n,\infty})_*\circ (i_{n,\infty})^! \to \on{Id}_{\QCoh(\CS)}
\end{equation} 
is an isomorphism.

\sssec{}  \label{sss:almost lift action b}

Set
$$Z_n:=S_n\underset{\CS}\times \CZ,$$
and let $\wt{i}_{n,\infty}$ denote the resulting maps
$$Z_n\to \CZ.$$

By \corref{c:base change QCoh Z}, we have
$$\QCoh(Z_n)\simeq \QCoh(S_n)\underset{\QCoh(\CS)}\otimes \QCoh(\CZ).$$

Hence, from \eqref{e:adj coarse} we obtain an adjunction 
\begin{equation} \label{e:adj upstairs}
(\wt{i}_{n,\infty})_*:\QCoh(Z_n)\rightleftarrows \QCoh(\CZ):(\wt{i}_{n,\infty})^!
\end{equation} 
as $\QCoh(\CZ)$-module categories. 

\medskip

In particular, the functors $(\wt{i}_{n,\infty})_*$ preserve compactness. Moreover, from 
\eqref{e:adj isom coarse} we obtain that the map
\begin{equation} \label{e:adj isom upstairs}
\underset{n}{\on{colim}}\, (\wt{i}_{n,\infty})_*\circ (\wt{i}_{n,\infty})^! \to \on{Id}_{\QCoh(\CZ)}
\end{equation} 
is an isomorphism.

\sssec{}  \label{sss:almost lift action c}

Let $$\CZ^{\on{rigid}}:=\CZ\underset{\on{pt}/\sG}\times \on{pt}$$
be the preimage of $\CZ$ in $\bMaps(\Rep(\sG),\bH)^{\on{rigid}}$. 

\medskip

Set
$$Z_n^{\on{rigid}}:=Z_n\underset{\CZ}\times \CZ^{\on{rigid}}\simeq S_n\underset{\CS}\times \CZ^{\on{rigid}},$$
so that
$$Z_n\simeq Z_n^{\on{rigid}}/\sG.$$

Note now that by \thmref{t:coarse restr}, $Z_n^{\on{rigid}}$ is an affine scheme. Hence, $Z_n$ is an algebraic stack. 

\sssec{Proof of \thmref{t:comp gen restr rough}} \label{sss:proof of rough}

Since the functors $(\wt{i}_{n,\infty})_*$ preserve compactness, and by \eqref{e:adj isom upstairs},
it suffices to show that each of the categories $\QCoh(Z_n)$ is compactly generated.

\medskip

However, for any algebraic stack $Z$ equal to the quotient of an affine scheme $Z^{\on{rigid}}$ by an action 
of the algebraic group $\sG$, the category $\QCoh(Z)$ is compactly generated by objects of the form
$$\CO_Z\otimes p^*(V), \quad V\in \Rep(\sG)^c,$$
where $p$ denote the map
$$Z\to \on{pt}/\sG,$$
corresponding to the $\sG$-torsor $Z^{\on{rigid}}\to Z$. 

\qed[\thmref{t:comp gen restr rough}]

\sssec{}

We will now give a slightly more precise form of the generation assertion. Let $p$
denote the map
$$\CZ\to \on{pt}/\sG,$$
corresponding to the $\sG$-torsor $\CZ^{\on{rigid}}\to \CZ$. 

\begin{thm} \label{t:comp gen restr fine}
The category $\QCoh(\CZ)$ is compactly generated by a family of objects of the form $\CF\otimes p^*(V)$,
where:

\begin{itemize}

\item $V\in \Rep(\sG)^c$;

\item $\CF$ can be expressed as a finite colimit in terms of $\CO_\CZ$.

\end{itemize} 

\end{thm}

\begin{rem}
We emphasize that the object $\CO_\CZ\in \QCoh(\CZ)$ itself is \emph{not} compact.
\end{rem}

\begin{proof}

The proof of \thmref{t:comp gen restr rough} in \secref{sss:proof of rough} shows
that in order to prove \thmref{t:comp gen restr fine}, we only need to prove that
$$(\wt{i}_{n,\infty})_*(\CO_{Z_n})\in \QCoh(\CZ)$$
can be expressed as a finite colimit in terms of $\CO_\CZ$.

\medskip

For that, it is sufficient to show that each
$$(i_{n,\infty})_*(\CO_{S_n})\in \QCoh(\CS)$$
can be expressed as a finite colimit in terms of $\CO_\CS$.

\medskip

However, this follows from the expression for the ring $R_n$ as
$$R_n\simeq R\underset{\sfe[t_1,...,t_m]}\otimes \sfe.$$

\end{proof} 

\ssec{Enhanced categorical trace}

In this subsection, we will prove an assertion that will be used in \secref{ss:enh Tr}.

\sssec{} \label{sss:enhanced trace setup}

Recall the set-up of \cite[Sects. 3.6-3.8]{GKRV}. 
We start with a symmetric monoidal category $\bA$ 
(assumed dualizable as a DG category), equipped with a symmetric monoidal endofunctor $F_\bA$. 
Let $\bM$ be an $\bA$-module category (assumed dualizable as such), equipped with
an endofunctor $F_\bM$, compatible with $F_\bA$. 

\medskip

Consider the category
$$\on{HH}_\bullet(F_\bA,\bA),$$
i.e., the category of Hochschild chains on $\bA$ twisted by $F_\bA$, see \cite[Sect. 3.7.2]{GKRV}. 
The fact that the monoidal structure on $(\bA,F_\bA)$ is symmetric allows us to define 
a symmetric monoidal structure on $\on{HH}_\bullet(F_\bA,\bA)$.

\medskip

Further, to $(\bM,F_\bM)$ we can attach an object
$$\Tr^{\on{enh}}_\bA(F_\bM,\bM)\in \on{HH}_\bullet(F_\bA,\bA),$$
see \cite[Sect. 3.8.2]{GKRV}. 

\sssec{}

Under the assumption that $\bA$ is rigid, we have the following assertion
(\cite[Theorem 3.8.5]{GKRV}): 

\medskip

There exists a canonical isomorphism in $\Vect_\sfe$:
%
\begin{equation} \label{e:enhanced trace abs}
\Tr(F_\bM,\bM)\simeq \CHom_{\on{HH}_\bullet(F_\bA,\bA)}\left(\one_{\on{HH}_\bullet(F_\bA,\bA)},\Tr^{\on{enh}}_\bA(F_\bM,\bM)\right),
\end{equation} 
where $\one_{\on{HH}_\bullet(F_\bA,\bA)}$ is the monoidal unit in $\on{HH}_\bullet(F_\bA,\bA)$.

\medskip

For example, if $\bA=\QCoh(\CY)$, where $\CY$ is an algebraic stack, and $F_\bA$ is given by $\phi^*$,
where $\phi$ is an endomorphism of $\CY$, we have
$$\on{HH}_\bullet(F_\bA,\bA) \simeq \QCoh(\CY^\phi),$$
where 
$$\CY^\phi:=\CY\underset{\Delta_\CY,\CY\times \CY,(\on{id}\times \phi)\circ \Delta_\CY}\times \CY,$$
and the right-hand side in \eqref{e:enhanced trace abs} is
\begin{equation} \label{e:Gamma stack}
\Gamma(\CY^\phi,\Tr^{\on{enh}}_\bA(F_\bM,\bM)).
\end{equation} 

\sssec{}

Our current goal is to generalize \eqref{e:enhanced trace abs} when instead of requiring that $\bA$
be rigid, we only require that $\bA$ be semi-rigid.

\medskip

The appropriate generalization is stated in \thmref{t:Tr enh abstract}, and proved in \secref{ss:proof Tr enh abstract}. 

\medskip

Here we will formulate its particular case, pertaining to the geometric situation, when 
$\bA:=\QCoh(\CY)$, for $\CY=\CY'/\sG$ where $\CY'$ and $\sG$ are as in \secref{sss:formal union quot}.

\sssec{} \label{sss:enhanced trace setup Y}

By \corref{c:descent do quot form aff}(a), $\bA$ is semi-rigid.
In particular, by \lemref{l:dual module over s-rigid}, in this case, 
an $\bA$-module category $\bM$ is dualizable if and only if the underlying DG category is dualizable.

\medskip

Assume that $\CY'$ is equipped with an endomorphism that commutes with the $\sG$-action, 
so that $\phi$ induces an endomorphism of $\CY$. We will denote both these endomorphisms by $\phi$.

\medskip

Note that 
$$\CY^\phi\simeq \left((\sG\times \CY')\underset{\on{act},\CY'\times \CY',(\on{id}\times \phi)\circ \Delta_\CY'}\times \CY'\right)/\sG,$$
where $(\sG\times \CY')\underset{\on{act},\CY'\times \CY',\Delta_\CY'}\times \CY'$ is also a formal affine scheme. 

\medskip

Let $F_\bA:=\phi^*$. By definition
$$\on{HH}_\bullet(F_\bA,\bA) \simeq 
\QCoh(\CY)\underset{\on{mult}_{\QCoh(\CY)},\QCoh(\CY)\otimes \QCoh(\CY),(\on{id}\otimes \phi^*)\circ \on{mult}_{\QCoh(\CY)}}\otimes \QCoh(\CY).$$

By \corref{c:base change QCoh Z}, the latter category maps isomorphically to
$$\QCoh(\CY \underset{\Delta_{\CY},\CY\times \CY,(\on{id}\times \phi)\circ \Delta_\CY}\times \CY)=\QCoh(\CY^\phi).$$

\medskip

Hence, in the setting of \secref{sss:enhanced trace setup}, we can think of $\Tr^{\on{enh}}_\bA(F_\bM,\bM)$ as an object of $\QCoh(\CY^\phi)$. 

\sssec{}

We claim:

\begin{thm} \label{t:enhanced Tr Y}
In the setting of \secref{sss:enhanced trace setup Y}, there 
is a canonical isomorphism 
\begin{equation} \label{e:enhanced trace abs Y}
\Tr(F_\bM,\bM)\simeq \Gamma_!\left(\CY^\phi,\Tr^{\on{enh}}_{\QCoh(\CY)}(F_\bM,\bM)\right).
\end{equation} 
\end{thm}

The proof is a generalization of the argument in \cite[Theorem 3.10.6]{GKRV}, and will be given
in \secref{s:semi-rigid} (see \secref{sss:proof enhanced Tr Y}). 

\begin{rem}
Note the difference between the assertion of \thmref{t:enhanced Tr Y} and a similar assertion 
when $\CY$ is an algebraic stack: in the latter case, instead of the right-hand side in \eqref{e:enhanced trace abs Y} we have
\eqref{e:Gamma stack}. 

\medskip

By contrast, in the case of formal schemes/stacks, the (discontinuous) functor $\Gamma(\CY^\phi,-)$
gets replaced by the functor $\Gamma_!(\CY^\phi,-)$.
\end{rem}

\begin{rem} \label{r:enhanced Tr Y}

Let us take $(\bM,F_\bM)$ to be $(\bA,F_\bA)$ itself. Then $\Tr^{\on{enh}}_\bA(F_\bA,\bA)=\one_{\on{HH}_\bullet(F_\bA,\bA)}$. 
So, in the setting of \thmref{t:enhanced Tr Y}, $\Tr^{\on{enh}}_\bA(F_\bA,\bA)\simeq \CO_{\CY}$. 

\medskip

Recall that for any semi-passable prestack $\wt{\CY}$,
the functor $\Gamma_!(\wt\CY,-)$ is \emph{non-unital} right-lax symmetric monoidal
(see \secref{sss:!-sect semi-pass}). So, $\Gamma_!(\wt\CY,\CO_{\wt\CY})$
acquires a structure of (not necessarily unital) commutative algebra, and for any $\CF\in \QCoh(\wt\CY)$, the object 
$\Gamma_!(\wt\CY,\CF)$ is naturally a module for $\Gamma_!(\wt\CY,\CO_{\wt\CY})$. 

\medskip

Applying this to $\wt\CY=\CY^\phi$, we obtain that, on the one hand,
$\Gamma_!\left(\CY^\phi,\Tr^{\on{enh}}_\bA(F_\bA,\bA)\right)$ acquires a structure of (not necessarily unital) 
commutative algebra,
and $\Gamma_!\left(\CY^\phi,\Tr^{\on{enh}}_\bA(F_\bM,\bM)\right)$ acquires a structure of module
over this commutative algebra. 

\medskip

On the other hand, \cite[Sects. 3.3.2 and 3.3.3]{GKRV}
implies that $\Tr(F_\bA,\bA)$ is naturally also a (not necessarily unital) commutative algebra and $\Tr(F_\bM,\bM)$ is a module 
over $\Tr(F_\bA,\bA)$. 

\medskip

As in \cite[Theorem 3.8.5]{GKRV}, the statement of \thmref{t:enhanced Tr Y} should be complemented
as follows:

\begin{itemize}

\item The isomorphism 
$$\Tr(F_\bA,\bA)\simeq \Gamma_!(\CY^\phi,\CO_{\CY^\phi})$$
of \eqref{e:enhanced trace abs Y} is compatible with commutative algebra structures
on both sides. 

\smallskip

\item The isomorphism 
$$\Tr(F_\bM,\bM)\simeq \Gamma_!\left(\CY^\phi,\Tr^{\on{enh}}_{\QCoh(\CY)}(F_\bM,\bM)\right)$$
respects the module structures for these algebras.

\end{itemize}

The proof of these compatibility assertions follows formally from \thmref{t:enhanced Tr Y}
as in \cite[Sects. 3.12.7-3.12.8]{GKRV}. 

\end{rem}

\newpage

\centerline{\bf Part II: Lisse actions and the spectral decomposition over $\LocSys^{\on{restr}}_{\sG}(X)$}

\bigskip

Let us make a brief overview of the contents of this Part. 

\medskip

In \secref{s:spectral decomp} we describe the set-up for the following question: what does it take to have an action 
of the monoidal category $\QCoh(\LocSys^{\on{restr}}_\sG(X))$ on a DG category $\bC$? It turns
that the appropriate input datum is what one can call \emph{an action of $\Rep(\sG)^{\otimes X\on{-lisse}}$ on $\bC$}.
In \thmref{t:action} we state that these two pieces of data are indeed
in bijection. We introduce an abstract framework for this result, where instead of $\qLisse(X)$
we are dealing with a general symmetric monoidal category $\bH$, equipped with a t-structure. The object of study
becomes the functor between symmetric monoidal categories
$$\coHom(\Rep(\sG),\bH) \to \QCoh(\Maps(\Rep(\sG),\bH))$$
(see \eqref{e:Hom abs}). We call a symmetric monoidal category \emph{adapted for spectral decomposition}
if the above functor is an equivalence. We state \conjref{c:Hom abs} to the effect that any \emph{gentle Tannakian category}
(see \secref{ss:gentle}) is adapted for spectral decomposition. A reformulation of \thmref{t:action}, stated 
as \thmref{t:spectral}, says that this conjecture holds for $\bH:=\qLisse(X)$.

\medskip

In \secref{s:part cases act} we prove \thmref{t:spectral}. We first show that the category 
$\bH:=\Shv_{\on{loc.const.}}(X)$, where $X$ is a connected CW complex is 
adapted for spectral decomposition. Next, we show that if $\bH$ is a adapted for spectral decomposition,
and $\bH'\subset \bH$ is a full subcategory, then under certain conditions, $\bH'$ is also
adapted for spectral decomposition. We then use a series of reductions showing that the category
$\qLisse(X_1)$ (where $X_2$ is a smooth and complete curve over a ground field $k$ of any
characteristic) can be realized as a full subcategory in one of the form $\Shv_{\on{loc.const.}}(X_2)$,
where $X_2$ is a curve over $\BC$, thereby deducing \thmref{t:spectral} from the Betti case. 

\medskip
%
%

\secref{s:Lie} is not needed for the rest of the paper. Here we consider another class
of symmetric monoidal categories adapted for spectral decomposition, namely, 
categories of the form $\sh\mod$, where $\sh$ is a connective Lie algebra.
We prove it by a method that we hope can be useful for the proof of \conjref{c:Hom abs}.  


\medskip

In Sects. \ref{s:progenitor projector}, \ref{s:Loc on LocSys} and \ref{s:projector abstract} 
we introduce a tool that will be extensively used in Part III of the 
paper--Beilinson's spectral projector.

\bigskip

\section{The spectral decomposition theorem} \label{s:spectral decomp}

In this section we the state the main theorem of Part II, \thmref{t:action} that
describes what it takes to have an action of $\QCoh(\LocSys^{\on{restr}}_\sG(X))$ 
on a DG category $\bC$.

\medskip

We will introduce an abstract framework in which \thmref{t:action} will be proved, and discuss
several reformulations. 


\ssec{Actions of $\Rep(\sG)^{\otimes X\on{-lisse}}$}

Let $X$ be a smooth, connected and complete curve. 

\medskip

In this subsection we define what it means to have an action $\Rep(\sG)^{\otimes X\on{-lisse}}$ on $\bC$, 
and state the main theorem of this part, \thmref{t:action}, which says that the datum of such an action on a category $\bC$ is
equivalent to the datum of an action on $\bC$ of the category $\QCoh(\LocSys^{\on{restr}}_\sG(X))$. 

%

\sssec{}  \label{sss:lisse action}

Let $\bC$ be a DG category. We define the notion of \emph{action of $\Rep(\sG)^{\otimes X\on{-lisse}}$ on $\bC$} 
by imitating \cite[Sects. C.1.2 and C.2.2]{GKRV}. Namely, this is a natural transformation between the following 
two functors $\on{fSet}\to \DGCat^{\on{Mon}}$: 

\medskip

From the functor
$$I\mapsto \Rep(\sG)^{\otimes I}$$
to the functor
$$I\mapsto \End(\bC)\otimes \qLisse(X)^{\otimes I}.$$


\medskip

In other words, informally, for every finite set $I$ we need to specify a monoidal functor 
$$\Rep(\sG)^{\otimes I}\to \End(\bC)\otimes \qLisse(X)^{\otimes I},$$
and for every map of finite sets $I\to J$, we need to supply a data of commutativity for
$$
\CD
\Rep(\sG)^{\otimes I} @>>>  \End(\bC)\otimes \qLisse(X)^{\otimes I} \\
@VVV @VVV \\
\Rep(\sG)^{\otimes J} @>>>  \End(\bC)\otimes \qLisse(X)^{\otimes J},
\endCD
$$
along with a homotopy-coherent system of compatibilities for compositions. 

\sssec{} \label{sss:action universal}

Consider the symmetric monoidal category $\QCoh(\LocSys^{\on{restr}}_\sG(X))$. We claim that there is
a canonically defined natural transformation between the following 
two functors $\on{fSet}\to \DGCat^{\on{SymMon}}$:

\medskip

From the functor
\begin{equation} \label{e:Rep to I}
I\mapsto \Rep(\sG)^{\otimes I}
\end{equation}
to the functor
\begin{equation} \label{e:QCoh to I}
I\mapsto \QCoh(\LocSys^{\on{restr}}_\sG(X))\otimes \qLisse(X)^{\otimes I}.
\end{equation}

Indeed, since $ \qLisse(X)$ is dualizable (and hence, tensoring by it commutes with limits), a datum of such a natural transformation is equivalent to a compatible system of 
natural transformations from \eqref{e:Rep to I} to
$$I\mapsto \QCoh(S)\otimes \qLisse(X)^{\otimes I} \text{ for } S\in \affSch_{/\LocSys^{\on{restr}}_\sG(X)}.$$
 
\medskip

By definition, the datum of a map $S\to \LocSys^{\on{restr}}_\sG(X)$ is a (right t-exact) symmetric monoidal functor
$$\sF:\Rep(\sG)\to \QCoh(S) \otimes \qLisse(X).$$

The required functor
$$\Rep(\sG)^{\otimes I}\to \QCoh(S)\otimes \qLisse(X)^{\otimes I}$$ 
is then the composition
\begin{equation} \label{e:E S}
\Rep(\sG)^{\otimes I}\overset{\sF^{\otimes I}}\longrightarrow \QCoh(S)^{\otimes I} \otimes \qLisse(X)^{\otimes I}\to
\QCoh(S)\otimes \qLisse(X)^{\otimes I},
\end{equation} 
where
$$\QCoh(S)^{\otimes I}\to \QCoh(S)$$
is the tensor product map.

\sssec{}  \label{sss:construction from universal}

From \secref{sss:action universal} we obtain that for any DG category $\bC$, equipped with an action of $\QCoh(\LocSys^{\on{restr}}_\sG(X))$,
we obtain an action of $\Rep(\sG)^{\otimes X\on{-lisse}}$ on $\bC$. I.e., we obtain a map of spaces
\begin{equation} \label{e:action}
\{\text{Actions of $\QCoh(\LocSys^{\on{restr}}_\sG(X))$ on $\bC$} \} \to 
\{\text{Actions of $\Rep(\sG)^{\otimes X\on{-lisse}}$ on $\bC$} \}. 
\end{equation}

\medskip

The main result of Part II of this paper is the following:

\begin{mainthm} \label{t:action}
The map \eqref{e:action} is an isomorphism. 
\end{mainthm}

We can regard this theorem as saying that a category $\bC$ equipped with an action of $\Rep(\sG)^{\otimes X\on{-lisse}}$, 
admits a spectral decomposition with respect to $\LocSys^{\on{restr}}_\sG(X)$. 

\sssec{}

The proof of \thmref{t:action} will be given in \secref{s:part cases act}. 

\medskip

In the next few subsections we will set up an abstract framework for \thmref{t:action}.

%
%

\ssec{The coHom symmetric monoidal category}

In this subsection we will make preparations for an abstract framework in which \thmref{t:action} can be stated. 

\sssec{} \label{sss:coHom first}

Let $\bH$ be a symmetric monoidal category. Assume that it is dualizable as a DG category. In this case,
there exists a monoidal category, to be denoted $\coHom(\Rep(\sG),\bH)$, defined by the universal property
that for a target symmetric monoidal category $\bA$, we have
$$\Maps_{\DGCat^{\on{SymMon}}}(\coHom(\Rep(\sG),\bH),\bA)\simeq 
\Maps_{\DGCat^{\on{SymMon}}}(\Rep(\sG),\bA\otimes \bH).$$

The construction of $\coHom(\Rep(\sG),\bH)$ fits into the following general paradigm. 

\sssec{}

Let $\bO$ be a symmetric monoidal category. We will assume that $\bO$ admits all colimits and that 
the monoidal operation commutes with colimits in each variable. 

\medskip

Let $A$ (resp., $C$) be a unital commutative algebra (resp., cocommutative coalgebra) object in $\bO$. 
In this case one can form a unital commutative algebra object
$$\on{coEnd}(A,C)\in \bO,$$
with the following universal property: for a unital commutative algebra object $A'\in \bO$, the space of maps of
(unital) commutative algebras $\on{coEnd}(A,C)\to A'$ is the space of maps in $\bO$
$$\phi:A\otimes C\to A',$$
equipped with a datum of commutativity for the diagrams
$$
\CD
A\otimes A\otimes C  @>{\on{mult}_A\otimes \on{id}_C}>>  A \otimes C  @>{\phi}>>  A' \\
@V{\on{id}_{A\otimes A}\otimes \on{comult}_C}VV & &   @AA{\on{mult}_{A'}}A   \\
A\otimes A\otimes C \otimes C  & @>{\phi\otimes \phi}>>  & A'\otimes A',
\endCD
$$
and
$$
\CD
C @>{\on{unit}_A}>> A \otimes C \\
@V{\on{counit}_C}VV @VV{\phi}V \\
\sfe @>{\on{unit}_{A'}}>> A',
\endCD
$$
along with a homotopy-coherent system of higher compatibilities.

\sssec{} \label{sss:coEnd}

The formal definition of $\on{coEnd}(A,C)$ is as follows: 

\medskip

Let $\on{fSet}$ be the category of finite sets, and let $\on{TwArr}(\on{fSet})$ be the 
corresponding twisted arrows category, see \cite[Sect. 1.2.2]{GKRV}.

\medskip

Consider the functor $\on{TwArr}(\on{fSet})\to O$ that at the level of objects sends 
$$(I\to J)\in \on{TwArr}(\on{fSet}) \to A^{\otimes I}\otimes C^{\otimes J}.$$
At the level of 1-morphisms, it sends the morphism 
$$
\CD
I_0 @>>> J_0 \\
@VVV @AAA \\
I_1 @>>> J_1
\endCD
$$
in $\on{TwArr}(\on{fSet})$, to the corresponding map
$$A^{\otimes I_0}\otimes C^{\otimes J_0}\to A^{\otimes I_1}\otimes C^{\otimes J_1}$$
given by the maps $A^{\otimes I_0}\to A^{\otimes I_1}$ (resp., $C^{\otimes J_0}\to C^{\otimes J_1}$),
given by the commutative algebra structure on $A$ (resp., cocommutative coalgebra structure on $C$).  

\medskip

Consider the colimit
\begin{equation} \label{e:coEnd}
\underset{(I\to J)\in \on{TwArr}(\on{fSet})}{\on{colim}}\, A^{\otimes I}\otimes C^{\otimes J}.
\end{equation} 

We endow \eqref{e:coEnd} with a structure of commutative algebra via the operation of disjoint union on $\on{fSet}$.



\sssec{} \label{sss:coEnd dualizable}

Let $A$ be as above, and let $B$ be another unital commutative algebra. Let $\coHom(A,B)$
be the commutative algebra in $\bO$ (if it exists) that has the following universal property
$$\Maps_{\on{ComAlg}(\bO)}(\coHom(A,B),A')\simeq \Maps_{\on{ComAlg}(\bO)}(A,B\otimes A'), \quad A'\in \on{ComAlg}(\bO).$$

We claim: 

\begin{lem} \label{l:coEnd abs}
Suppose $B$ is dualizable as a plain object of $\bO$. 
Then the object $\coHom(A,B)\in \on{ComAlg}(\bO)$ exists and is canonically isomorphic to
$\on{coEnd}(A,C)$, where $C:=B^\vee$, viewed as a cocommutative coalgebra in $\bO$. 
\end{lem}

The proof will be given in \secref{sss:proof coEnd}.

\sssec{}

Thus, the symmetric monoidal category $\coHom(\Rep(\sG),\bH)$ introduced in \secref{sss:coHom first} fits into the above paradigm with
$\bO:=\DGCat$ and $A:=\Rep(\sG)$. 

\sssec{}

We return to the setting of \secref{sss:coEnd dualizable}. As in \cite[Theorem 1.2.4]{GKRV}, one shows:

\begin{lem}   \label{l:maps from coend}  Assume that $B$ be a unital commutative algebra that is
dualizable as an object of $\bO$. 

\smallskip

\noindent{\em(a)} For an associative/commutative algebra object $D\in \bO$, the space of maps of associative/commutative algebras
$$\coHom(A,B)\to D$$ identifies with the space of compatible collections of maps of associative/commutative algebras
$$A^{\otimes I}\to D\otimes B^{\otimes I}, \quad I\in \on{fSet}.$$

\smallskip

\noindent{\em(b)} For a \emph{plain} object $D\in \bO$, the space $\Maps_\bO(\coHom(A,B),D)$
identifies with the space of compatible collections of maps in $\bO$
$$A^{\otimes I}\to D\otimes B^{\otimes I}, \quad I\in \on{fSet}.$$

\end{lem} 

\sssec{Sketch of proof of \lemref{l:maps from coend}} \label{sss:proof tw Arr}

Here we will sketch a proof. A full argument will be given in \secref{ss:proof coEnd assoc}.

\medskip

Consider the colimit \eqref{e:coEnd}. 
We first consider it as a plain object of $O$. Then for $D\in O$,
the space of maps from \eqref{e:coEnd} to $D$ is, by definition,
$$\underset{(I\to J)\in \on{TwArr}(\on{fSet})}{\on{lim}}\, \Maps_O(A^{\otimes I},D\otimes B^{\otimes J}).$$

However (see, e.g., \cite[Lemma 1.3.12]{GKRV}), the latter expression identifies with the space of natural transformations
between the functors
$$\on{fSet}\to \bO, \quad (I\mapsto A^{\otimes I})\, \Rightarrow\, (I\mapsto D\otimes B^{\otimes I}).$$

\medskip

Suppose now that $D$ is an associative/commutative algebra in $\bO$. We claim that the space of maps from \eqref{e:coEnd} to $D$ 
that are upgraded to maps of algebras correspond to compatible systems of maps of algebras
\begin{equation} \label{e:compat collection}
A^{\otimes I}\to D\otimes B^{\otimes I}, \quad I\in \on{fSet}.
\end{equation} 

\medskip

Indeed, for a map from \eqref{e:coEnd} to $D$, the data of compatibility with an associative/commutative algebra structure
translates into the data of commutativity of the diagrams 
\begin{equation} \label{e:compat via disj}
\CD
A^{\otimes I_1}\otimes...\otimes A^{\otimes I_k} @>>> 
(D\otimes B^{\otimes I_1})\otimes...\otimes (D\otimes B^{\otimes I_k}) \\
@V{\sim}VV @VVV \\
A^{\otimes (I_1 \sqcup...\sqcup I_k)} @>>> D\otimes B^{\otimes (I_1\sqcup...\sqcup I_k)}
\endCD
\end{equation} 
for unordered/ordered collections of finite sets $I_1,...,I_k$. 

\medskip

If the maps in \eqref{e:compat collection} are maps of algebras, the commutativity for the diagrams
\eqref{e:compat via disj} arises from 
$$
\CD
A^{\otimes I_1}\otimes...\otimes A^{\otimes I_k} @>>> (D\otimes B^{\otimes I_1})\otimes...\otimes (D\otimes B^{\otimes I_k}) \\ 
@VVV @VVV  \\
(A^{\otimes I})^{\otimes k} @>>> (D\otimes B^{\otimes I})^{\otimes k}  \\
@VVV @VVV \\
A^{\otimes I} @>>> D\otimes B^{\otimes I},
\endCD
$$
where $I=I_1\sqcup...\sqcup I_k$, 
and where the upper vertical maps are given by inclusions $I_j\to I$. 

\medskip

Vice versa, given the commutative diagrams \eqref{e:compat via disj}, we construct the data of
compatibility for \eqref{e:compat collection} by
$$
\CD
(A^{\otimes I})^{\otimes k} @>>> (D\otimes B^{\otimes I})^{\otimes k} \\
@V{\sim}VV @VVV  \\
A^{\otimes I'} @>>> D \otimes B^{\otimes I'} \\
@VVV @VVV \\
A^{\otimes I} @>>> D \otimes B^{\otimes I},
\endCD
$$
where $I'$ is the disjoint union of $k$ copies of $I$, and the lower vertical maps are given by the natural projection $I'\to I$.

\qed

\ssec{Maps vs coHom}

In this subsection we study the relationship of the category $\coHom(\Rep(\sG),\bH)$ introduced
above and its algebro-geometric counterpart, the prestack $\bMaps(\Rep(\sG),\bH)$.

\sssec{}  \label{sss:map abs}

Let $\bH$ be a symmetric monoidal category, equipped with a t-structure and a 
fiber functor satisfying the assumptions of \secref{sss:pre-Tannakian}.  

\medskip

From now on we will add the assumption that $\bH$ is dualizable as a DG category.

\sssec{}  \label{sss:map abs bis}

On the one hand, we can consider the symmetric monoidal category $\coHom(\Rep(\sG),\bH)$, introduced above.
On the other hand, we can consider the prestack $\bMaps(\Rep(\sG),\bH)$, see \secref{ss:abs LocSys}. 

\medskip

We claim that we have a canonically defined 
symmetric monoidal functor
\begin{equation} \label{e:Hom abs}
\coHom(\Rep(\sG),\bH) \to \QCoh(\bMaps(\Rep(\sG),\bH)).
\end{equation} 

Indeed, the datum of such a functor is by definition equivalent to the datum of a symmetric monoidal functor
\begin{equation} \label{e:Hom abs pre}
\Rep(\sG) \to \QCoh(\bMaps(\Rep(\sG),\bH)) \otimes \bH.
\end{equation} 

The functor in \eqref{e:Hom abs pre} is obtained by passing to the limit from the tautological functors
$$\Rep(\sG) \to \QCoh(S) \otimes \bH, \quad S\in \affSch_{/\bMaps(\Rep(\sG),\bH)},$$
using the fact that
$$\QCoh(\bMaps(\Rep(\sG),\bH)) \otimes \bH=\left(\underset{S\in \affSch_{/\bMaps(\Rep(\sG),\bH)}}{\on{lim}}\, \QCoh(S)\right)\otimes \bH\to$$
$$\to \underset{S\in \affSch_{/\bMaps(\Rep(\sG),\bH)}}{\on{lim}}\, \left(\QCoh(S)\otimes \bH\right)$$
is an equivalence, the latter since $\bH$ is dualizable as a DG category. 

\sssec{}

We shall say that $\bH$ is \emph{adapted for spectral decomposition} (for a given $\sG$) if the functor 
\eqref{e:Hom abs} is an equivalence.

\medskip

\begin{rem}
In the course of the next two sections we will see examples of symmetric monoidal categories $\bH$
that are adapted for spectral decomposition. These examples include 
$\bH:=\Vect_\sfe^\CX$, where $\CX$ is a connected object of $\Spc$ (see \secref{sss:Vect X}), and $\bH:=\sh\mod$, where
$\sh$ is a connective Lie algebra. 
%
%
\medskip

Note that in the above two examples, $\bH$ is not gentle (see \secref{sss:conditions} for what this means). 

\medskip

One can also show it holds for $\bH=\Rep(\sH)$, where $\sH$ is an affine algebraic group of finite type
(but we will not prove this in the present paper). 

\end{rem} 

\sssec{}

We propose: 

\begin{conj} \label{c:Hom abs}
If $\bH$ is a gentle Tannakian category, then it is adapted for spectral decomposition.
\end{conj}

We will prove:

\begin{mainthm} \label{t:spectral}
\conjref{c:Hom abs} holds when $\bH=\qLisse(X)$, where $X$ is a smooth and complete algebraic curve.
\end{mainthm} 

We will see shortly that \thmref{t:spectral} is equivalent to \thmref{t:action}. 

\ssec{Spectral decomposition vs actions}

In this subsection we will reformulate the property of being adapted for spectral decomposition
in terms of actions on a module category $\bC$.

\sssec{}

Let $\bH$ be a (dualizable) symmetric monoidal category and let $\bC$ be a DG category. By an 
$\bH$-family of actions of $\Rep(\sG)$ on $\bC$ we will mean a natural transformation between the following 
two functors $\on{fSet}\to \DGCat^{\on{Mon}}$: 

\medskip

From the functor
$$I\mapsto \Rep(\sG)^{\otimes I}$$
to the functor
$$I\mapsto \End(\bC)\otimes \bH^{\otimes I}.$$

\medskip

\lemref{l:maps from coend}(a) implies that this data is equivalent to that of an action of 
$\coHom(\Rep(\sG),\bH)$, viewed as a monoidal category, on $\bC$.

\sssec{} \label{sss:X to the Lisse}

Taking $\bH=\qLisse(X)$, we obtain that the notion of a family of $\qLisse(X)$-actions on $\bC$
just defined is the same as a the notion of action of $\Rep(\sG)^{\otimes X\on{-lisse}}$ on $\bC$ from
\secref{sss:lisse action}.

\medskip

In particular, the symbol $\Rep(\sG)^{\otimes X\on{-lisse}}$ stands for an actual (symmetric) monoidal category, namely, 
$$\Rep(\sG)^{\otimes X\on{-lisse}}\simeq \coHom(\Rep(\sG),\qLisse(X)).$$

\medskip

Thus, we can regard \thmref{t:spectral} as saying that the functor
\begin{equation}  \label{e:Hom concr}
\Rep(\sG)^{\otimes X\on{-lisse}}\to \QCoh(\LocSys^{\on{restr}}_\sG(X)),
\end{equation}
described in \secref{sss:action universal}, is an equivalence. 

\sssec{}

Let now $\bH$ be endowed with a t-structure and a fiber functor satisfying the assumptions of \secref{sss:pre-Tannakian}.  

\medskip

The map \eqref{e:Hom abs} gives rise to a map of spaces 
\begin{equation} \label{e:action abs}
\{\text{Actions of $\QCoh(\bMaps(\Rep(\sG),\bH))$ on $\bC$} \} \to 
\{\text{$\bH$-families of actions of $\Rep(\sG)$ on $\bC$} \}. 
\end{equation}

From here, we obtain:

\begin{lem} \label{l:action vs spectral}
The category $\bH$ is adapted for spectral decomposition if and only if the map \eqref{e:action abs}
is an equivalence for any $\bC$.
\end{lem}

\sssec{}

Unwinding the constructions, it is easy to see that for $\bH=\qLisse(X)$, 
the map  \eqref{e:action abs} is the same as the map \eqref{e:action}. 

\medskip

Hence, \lemref{l:action vs spectral} implies that Theorems \ref{t:action} and \ref{t:spectral}
are logically equivalent.

\sssec{}

Let $\bH$ be again a (dualizable) symmetric monoidal category and let $\bC$ be a DG category. 
By an 
$\bH$-family of functors $\Rep(\sG)\to \bC$ we will mean a natural transformation between the following 
two functors $\on{fSet}\to \DGCat$: 

\medskip

From the functor
$$I\mapsto \Rep(\sG)^{\otimes I}$$
to the functor
$$I\mapsto \bC\otimes \bH^{\otimes I}.$$

\medskip

From \lemref{l:maps from coend}(b) we obtain that this data is equivalent to the data of a functor
$$\coHom(\Rep(\sG),\bH)\to \bC.$$

\sssec{}

Let $\bH$ be endowed with t-structure and a fiber functor satisfying the assumptions of \secref{sss:pre-Tannakian}.  

\medskip

The map \eqref{e:Hom abs} gives rise to a map of spaces 
\begin{equation} \label{e:functors abs}
\{\text{Functors $\QCoh(\bMaps(\Rep(\sG),\bH))\to \bC$} \} \to 
\{\text{$\bH$-families of functors $\Rep(\sG)\to \bC$}\}. 
\end{equation}

As in \lemref{l:action vs spectral}, we have:

\begin{lem} \label{l:functors vs spectral}
The category $\bH$ is adapted for spectral decomposition if and only if the map \eqref{e:functors abs}
is an equivalence for any $\bC$.
\end{lem}

\sssec{} \label{sss:pre-shtuka}

Let us write out the map \eqref{e:functors abs} explicitly. 

\medskip

We start with the functors \eqref{e:Hom abs pre}. Then for $I\in \on{fSet}$ we obtain a functor
\begin{equation} \label{e:E univ}
\Rep(\sG)^{\otimes I}\to \QCoh(\bMaps(\Rep(\sG),\bH))^{\otimes I} \otimes \bH^{\otimes I} \to
\QCoh(\bMaps(\Rep(\sG),\bH)) \otimes \bH^{\otimes I},
\end{equation} 
where the last arrow uses the tensor product functor $$\QCoh(\bMaps(\Rep(\sG),\bH))^{\otimes I}\to \QCoh(\bMaps(\Rep(\sG),\bH)).$$ 

\medskip

We will denote the functor \eqref{e:E univ} by $\CE^I$. For a fixed $V\in \Rep(\sG)^{\otimes I}$, we denote the resulting
object of $\QCoh(\bMaps(\Rep(\sG),\bH)) \otimes \bH^{\otimes I}$ by $\CE^I_V$. 

\medskip

Now, given a functor 
$$\CS:\QCoh(\bMaps(\Rep(\sG),\bH)) \to \bC,$$
the resulting system of functors
$$\CS^I:\Rep(\sG)^{\otimes I}\to  \bC\otimes \bH^{\otimes I}$$
sends 
$$V\in  \Rep(\sG)^{\otimes I} \, \mapsto\, (\CS\otimes \on{Id})(\CE^I_V)\in \bC\otimes \bH^{\otimes I}.$$

\begin{rem}

Note that for $\bH=\qLisse(X)$ and $\bC=\Vect_\sfe$, a system of functors
$$\Rep(\sG)^{\otimes I}\to \qLisse(X)^{\otimes I}, \quad I\in \on{fSet}$$
is exactly the structure that arises from the \emph{Shtuka} construction. 

\medskip

We will explore this in \secref{ss:shtukas} to relate Shtukas to objects in $\QCoh(\LocSys_\sG^{\on{restr}}(X))$. 

\end{rem}


%

%

%
%
%
%
%
%
%
%

\ssec{A rigidified version} \label{ss:rigidified coMaps}

Let $\bH$ be as in \secref{sss:pre-Tannakian}. Recall that along with the prestack $\bMaps(\Rep(\sG),\bH)$ we considered
its rigidified version $\bMaps(\Rep(\sG),\bH)^{\on{rigid}}$. In this subsection we will introduce a counterpart of this
rigidification for $\coHom(\Rep(\sG),\bH)$. 

\sssec{}

Let $\bH$ be a (dualizable) symmetric monoidal category, equipped with a symmetric monoidal functor 
$\oblv_\bH:\bH\to \Vect_\sfe$.

\medskip

Composition with $\oblv_\bH$ defines a symmetric monoidal functor
$$\Rep(\sG)\simeq \coHom(\Rep(\sG),\Vect_\sfe) \to \coHom(\Rep(\sG),\bH).$$

Denote 
$$\coHom(\Rep(\sG),\bH)^{\on{rigid}}:=\coHom(\Rep(\sG),\bH)\underset{\Rep(\sG)}\otimes \Vect_\sfe.$$

\sssec{}

By construction, for a symmetric monoidal category $\bA$, the datum of a symmetric monoidal functor
$$\coHom(\Rep(\sG),\bH)^{\on{rigid}}\to \bA$$
is equivalent to that of a symmetric monoidal functor
$$\Rep(\sG)\to \bA\otimes \bH,$$
equipped with an identification of the composition 
$$\Rep(\sG)\to \bA\otimes \bH \overset{\on{Id}_\bA\otimes \oblv_\bH}\longrightarrow \bA,$$
with the forgetful functor
$$\Rep(\sG) \overset{\oblv_\sG}\to \Vect_\sfe\overset{\one_\bA}\longrightarrow \bA.$$

\sssec{}

As in \secref{sss:map abs bis}, we have a symmetric monoidal functor
\begin{equation} \label{e:Hom abs rigid}
\coHom(\Rep(\sG),\bH)^{\on{rigid}} \to \QCoh(\bMaps(\Rep(\sG),\bH)^{\on{rigid}}).
\end{equation} 

Since
$$\QCoh(\bMaps(\Rep(\sG),\bH)^{\on{rigid}})
\simeq \QCoh(\bMaps(\Rep(\sG),\bH))\underset{\Rep(\sG)}\otimes \Vect_\sfe,$$
we obtain that $\bH$ is adapted for spectral decomposition if and only if the functor \eqref{e:Hom abs rigid}
is an equivalence. Indeed, this follows from the fact that the functor \eqref{e:de-eq}
is conservative. 

\section{Categories adapted for spectral decomposition} \label{s:part cases act}

The goal of this section is to prove \thmref{t:spectral}. Our strategy will be as follows:

\medskip

We will first show that the category $\bH:=\Vect^\CX_\sfe$ is adapted for spectral decomposition
(where $\CX$ is a connected object of $\Spc$). From this we will then formally deduce that the
category $\bH:=\qLisse(X)$ is also adapted for spectral decomposition, in the particular case when
$X$ is a smooth and compete curve. 

\ssec{The Betti case}  \label{ss:full Betti}

In this subsection we let $\CX$ be a connected object of $\Spc$. 

\sssec{}

Consider the symmetric monoidal category $\Vect_\sfe^\CX$, equipped with its natural t-structure 
and the fiber functor (the latter is given by a choice of a base point $x\in X$). 

\medskip

We claim:

\begin{thm} \label{t:spectral Betti}
The symmetric monoidal category $\Vect_\sfe^\CX$ is adapted for spectral decomposition.
\end{thm}

This result is stated and proved in \cite[Theorem 1.5.5]{GKRV}. In fact, the category that we denote
$\coHom(\Rep(\sG),\Vect_\sfe^\CX)$ is exactly the category denoted $\Rep(\sG)^{\otimes \CX}$
in {\it loc.cit.}, and 
$$\bMaps( \Rep(\sG),\Vect_\sfe^\CX)=\LocSys^{\on{Betti}}_\sG(\CX),$$
see \secref{sss:Betti coMaps}.

\begin{rem}

In \secref{ss:another Betti} we will give another proof of \thmref{t:spectral Betti}, which has a potential 
for generalization for other symmetric monoidal categories $\bH$.

\end{rem}

\ssec{The heriditary property of being adapted}

In this subsection we will perform a crucial step towards the proof of \thmref{t:spectral}: we will
show that the property of being adapted for spectral decomposition is, under certain conditions,
inherited by full subcategories.

\sssec{}

Let $\bH$ be a (dualizable) symmetric monoidal category as in \secref{sss:pre-Tannakian}, and let
$\bH'\subset \bH$ be a full symmetric monoidal subcategory. Assume that $\bH'$ inherits a t-structure
(i.e., it is preserved by the truncation functors). 

\medskip

We will prove:

\begin{thm} \label{t:adapted inherited} Suppose that:

\begin{itemize}
  
\item The embedding $\iota:\bH'\hookrightarrow \bH$, considered as a functor between plain DG categories,
admits a continuous right adjoint;

\smallskip

\item The prestack $\bMaps(\Rep(\sG),\bH)^{\on{rigid}}$ is an eventually coconnective affine scheme almost of finite type;

\smallskip

\item The map $\bMaps(\Rep(\sG),\bH')\to \bMaps(\Rep(\sG),\bH)$ is a formal isomorphism and an 
ind-closed embedding\footnote{See Remark \ref{r:union formal compl} where it is explained what the combination
of these two conditions amounts to.}. 

\end{itemize}

Then if $\bH$ is adapted for spectral decomposition, then so is $\bH'$.
\end{thm} 

The rest of this subsection is devoted to the proof of this theorem.

\sssec{}

Let $\bA$ be a target symmetric monoidal category. We wish to show that the map
$$\Maps_{\DGCat^{\on{SymMon}}}(\QCoh(\bMaps(\Rep(\sG),\bH')),\bA) \to 
\Maps_{\DGCat^{\on{SymMon}}}(\Rep(\sG),\bA\otimes \bH')$$
is an isomorphism of spaces. 

\medskip

We have a commutative diagram
\begin{equation}  \label{e:diagram of funct}
\CD
\Maps_{\DGCat^{\on{SymMon}}}(\QCoh(\bMaps(\Rep(\sG),\bH)),\bA) @>>> 
\Maps_{\DGCat^{\on{SymMon}}}(\Rep(\sG),\bA\otimes \bH) \\
@AAA @AAA \\
 \Maps_{\DGCat^{\on{SymMon}}}(\QCoh(\bMaps(\Rep(\sG),\bH')),\bA) @>>> 
\Maps_{\DGCat^{\on{SymMon}}}(\Rep(\sG),\bA\otimes \bH'),
 \endCD
\end{equation}
where the top horizontal arrow is an isomorphism by assumption. We will show that both vertical arrows are fully faithful, and that their
essential images match up under the equivalence given by the top horizontal arrow.  

\sssec{}

The right vertical arrow in \eqref{e:diagram of funct} is a fully faithful because the functor
$$\bA\otimes \bH'\to \bA\otimes \bH$$
is fully faithful. The latter is true because the inclusion functor $\bH'\hookrightarrow \bH$ is fully faithful and
admits a continuous right adjoint. 

\sssec{} \label{sss:union formal compl again}

%

Let 
$$\CY_1=(\CY_2)^\wedge_\CZ\to \CY_2$$
be as in Remark \ref{r:union formal compl}. 

\medskip
 
A simple colimit argument, combined with \corref{c:ind vs pro}(a), shows that the restriction functor
\begin{equation} \label{e:restr to formal compl}
\QCoh(\CY_2)\to \QCoh(\CY_1)
\end{equation}
admits a fully faithful left adjoint, whose essential image is the full subcategory
$$\QCoh(\CY_2)_\CZ\subset \QCoh(\CY_2)$$
consisting of objects with set-theoretic support on $\CZ$, i.e., those objects 
$\CF\in \QCoh(\CY_2)$ such that for every affine scheme $S$ mapping to $\CY_2$,
the pullback $\CF_S$ of $\CF$ to $S$ vanishes on the localization of $S$ at every scheme-theoretic
point not contained in $S\underset{\CY_2}\times \CZ$ (equivalently, $\CF_S$ is such that its 
cohomologies are unions of subsheaves supported on closed subsets of $S$ that comprise
$S\underset{\CY_2}\times \CZ$).

\sssec{}  \label{sss:when action factors}

This implies that in the situation of \secref{sss:union formal compl again}, 
for any DG category $\bD$,
restriction along \eqref{e:restr to formal compl} defines a fully faithful embedding
\begin{equation} \label{e:restr to formal compl funct}
\Maps_{\DGCat}(\QCoh(\CY_1),\bD)\to 
\Maps_{\DGCat}(\QCoh(\CY_2),\bD),
\end{equation}
with essential image consisting of those functors that vanish on the full subcategory 
\begin{equation} \label{e:vanish Z alpha}
\{\CF\in \QCoh(\CY_2),\,\, \CF|_{\CY_1}=0\}
\end{equation} 
of $\QCoh(\CY_2)$.

\medskip

This formally implies that for any symmetric monoidal category $\bA$ the map
$$\Maps_{\DGCat^{\on{SymMon}}}(\QCoh(\CY_1),\bA)\to 
\Maps_{\DGCat^{\on{SymMon}}}(\QCoh(\CY_2),\bA)$$
is fully faithful, whose essential image consists of those functors that vanish on the subcategory 
\eqref{e:vanish Z alpha} of $\QCoh(\CY_2)$. 

\sssec{}

We will need the following assertion: 

\begin{lem} \label{l:descr ker}
Suppose that $\CY_2$
is an eventually coconnective affine scheme almost of finite type. 
Then the subcategory \eqref{e:vanish Z alpha}
is generated by objects of the form $f_*(\wt\sfe)$, where $\wt\sfe$ is a field
extension of $\sfe$ and $f$ is a map $\Spec(\wt\sfe)\to \CY_2$ that \emph{does not} factor
through $\CZ$. 
\end{lem} 

\begin{proof}
%
%

Since $\CY_2$ is eventually coconnective, every object is a (finite) colimit of objects
obtained as direct images along $^{\on{cl}}\CY_2\to \CY_2$. 

\medskip

Hence, we can assume that $\CY_2$ is classical. In this case, the statement follows
by a standard Cousin complex argument.

\end{proof}

\sssec{} 

We apply the discussion in \secref{sss:when action factors} to the embedding
\begin{equation} \label{e:compare Maps}
\bMaps(\Rep(\sG),\bH')\to \bMaps(\Rep(\sG),\bH). 
\end{equation}

We obtain that the left vertical arrow in \eqref{e:diagram of funct} is fully faithful.

\begin{rem} \label{r:surj on QCoh}

Note that the above argument shows that in the situation of \thmref{t:adapted inherited}, 
the functor
$$\coHom(\Rep(\sG),\bH')\to \QCoh(\bMaps(\Rep(\sG),\bH'))$$
is a localization (admits a fully faithful (but not necessarily continuous) right adjoint, 
even without the assumption that $\bMaps(\Rep(\sG),\bH)^{\on{rigid}}$ is eventually coconnective. 

\end{rem} 

\sssec{}

To prove \thmref{t:adapted inherited}, it remains to show that the essential images of the vertical arrows in \eqref{e:diagram of funct} 
match under the equivalence given by the top horizontal arrow. 

\ssec{Proof of \thmref{t:adapted inherited}: identifying the essential image}

\sssec{}

Applying base change
$$\Vect_\sfe\underset{\Rep(\sG)}\otimes -,$$
we can assume that we are given a functor
$$\Phi:\QCoh(\bMaps(\Rep(\sG),\bH)^{\on{rigid}})\to \bA,$$
such that the corresponding functor
$$\sF:\Rep(\sG)\to \bA\otimes \bH$$
factors as
$$\Rep(\sG)\overset{\sF'}\to \bA\otimes \bH'\to \bA\otimes \bH.$$

We wish to show that $\Phi$ factors as 
$$\QCoh(\bMaps(\Rep(\sG),\bH)^{\on{rigid}})\to \QCoh(\bMaps(\Rep(\sG),\bH')^{\on{rigid}})\to \bA.$$

I.e., we wish to show that $\Phi$ vanishes on
$$\on{ker}\left(\QCoh(\bMaps(\Rep(\sG),\bH)^{\on{rigid}})\to \QCoh(\bMaps(\Rep(\sG),\bH')^{\on{rigid}})\right).$$ 

\sssec{}

By \lemref{l:descr ker}, it suffices to show the following;

\medskip

Let $\wt\sfe$ be a field extension of $\sfe$, and let us be given a map
$$f:\Spec(\wt\sfe)\to \bMaps(\Rep(\sG),\bH)^{\on{rigid}}.$$

We wish to show that if $\Phi(f_*(\wt\sfe))\neq 0$, then $f$ factors through 
$\bMaps(\Rep(\sG),\bH')^{\on{rigid}}$.

\sssec{}

Consider the tensor product category
$$\wt\bA:=\Vect_{\wt\sfe}\underset{\QCoh(\bMaps(\Rep(\sG),\bH)^{\on{rigid}})}\otimes \bA.$$

Since the morphism $f$ is affine, we have 
$$\Vect_{\wt\sfe}\simeq f_*(\wt\sfe)\mod(\QCoh(\bMaps(\Rep(\sG),\bH)^{\on{rigid}})).$$

Hence, 
$$\wt\bA\simeq \Phi(f_*(\wt\sfe))\mod(\bA).$$
Therefore, if $\Phi(f_*(\wt\sfe))\neq 0$, then $\wt\bA\neq 0$. 

\sssec{}

Denote by $\wt\Phi$ the composite functor 
$$\QCoh(\bMaps(\Rep(\sG),\bH)^{\on{rigid}})\overset{\Phi}\to \bA\to \wt\bA$$
and by $\wt\sF$ the corresponding functor
$$\Rep(\sG) \overset{\sF}\to \bA\otimes \bH \to \wt\bA\otimes \bH.$$

By assumption, the functor $\wt\sF$ takes values in the full subcategory 
$$\wt\bA\otimes \bH'\subset \wt\bA\otimes \bH.$$

Note, however, that $\wt\Phi$ factors as 
$$\QCoh(\bMaps(\Rep(\sG),\bH)^{\on{rigid}})\overset{f^*}\to \Vect_{\wt\sfe}\to \wt\bA,$$
and $\wt\sF$ factors as 
$$\Rep(\sG) \overset{\sF_f}\to \Vect_{\wt\sfe}\otimes \bH \to \wt\bA\otimes \bH,$$
where $\sF_f$ is the functor corresponding to $f$ in the definition of 
$\bMaps(\Rep(\sG),\bH)^{\on{rigid}}$.

\sssec{}

Thus, we wish to show that $\sF_f$ takes values in
$$\Vect_{\wt\sfe}\otimes \bH'\subset  \Vect_{\wt\sfe}\otimes \bH.$$

\sssec{}

Recall that $\iota$ denotes the embedding $\bH'\hookrightarrow \bH$. 
For an object $V\in \Rep(\sG)$ consider the counit of the adjunction
$$(\on{Id}_{\Vect_{\wt\sfe}}\otimes (\iota\circ\iota^R))(\sF_f(V))\to \sF_f(V).$$

We wish to show that it is an isomorphism. We know that this map becomes an isomorphism 
after applying the functor
$$\Vect_{\wt\sfe}\otimes \bH\to \wt\bA\otimes \bH.$$

Hence, it is enough to show that the latter functor is conservative. Since $\bH$ is dualizable,
it suffices to show that the functor
$$\Vect_{\wt\sfe}\to \wt\bA$$
is conservative. 

\medskip

However, the latter is evident: up to a cohomological shift, a non-zero object of $\Vect_{\wt\sfe}$ has a copy of $\wt\sfe$
as a retract, and $\wt\sfe \mapsto \one_{\wt\bA}$, which is non-zero, since $\wt\bA$ was assumed non-zero. 

\qed[\thmref{t:adapted inherited}]

\ssec{Proof of \thmref{t:spectral}, Betti and de Rham contexts} \label{ss:spectral Betti and de Rham}

In this subsection we will prove \thmref{t:spectral} in the Betti and de Rham contexts.

\medskip

Our method will consist of combining 
Theorems \ref{t:spectral Betti} and \ref{t:adapted inherited}. Throughout this section,
$X$ will be a smooth, complete and connected curve over a ground field $k$.

\sssec{} \label{sss:spectral proof Betti}

We will first consider the case when $k=\BC$ and our sheaf-theoretic context
is Betti (see \secref{sss:Shv} for what this means). In this case, our curve $X$ does 
not need to be complete. 

\medskip

We take $\bH:=\Shv^{\on{all}}_{\on{loc.const.}}(X)$ and $\bH':=\qLisse(X)$. We know that the category
$\Shv^{\on{all}}_{\on{loc.const.}}(X)$ is adapted for spectral decomposition by \thmref{t:spectral Betti},
which would imply \thmref{t:spectral} in this case. 

\medskip

The corresponding functor 
\begin{equation} \label{e:qLisse to Betti again}
\iota:\qLisse(X)\hookrightarrow \Shv^{\on{all}}_{\on{loc.const.}}(X)
\end{equation} 
is fully faithful by \propref{p:lisse to loc const}. We will show that 
it satisfies the requirements of \thmref{t:adapted inherited}. 

\sssec{}

We first show that $\iota$ admits a continuous right adjoint. We will distinguish two cases:

\medskip

\noindent Case 1: $X$ has genus\footnote{If $X$ not complete, then we stipulate that we are in Case 2 below.}
$0$. In this case $\iota$ is an equivalence. 

\medskip

\noindent Case 2: $X$ has genus $\geq 1$. In this case, by \thmref{t:curves}(a) and \corref{c:Kpi1}, the functor
$$\iLisse(X)\to \qLisse(X)$$ is an equivalence. Now, for any $X$, the composite functor
$$\iLisse(X)\to \qLisse(X)\to \Shv^{\on{all}}_{\on{loc.const.}}(X)$$
sends compacts to compacts (indeed, for a finite CW complex, objects from $\Lisse(X)$ are
compact in $\Shv^{\on{all}}_{\on{loc.const.}}(X)$), and since $\iLisse(X)$ is compactly generated,
it admits a continuous right adjoint.

\sssec{}

We now show that 
$$\bMaps(\Rep(\sG),\bH)^{\on{rigid}} \simeq \LocSys_\sG^{\on{Betti,rigid}_x}(X)$$
is eventually coconnective and almost of finite type. The aft condition holds for any $X$
that is homotopy-equivalent to a finite CW complex. To show that it is eventually coconnective,
we will show that it is quasi-smooth. 

\medskip

Let $\overset{\circ}{X}$ be obtained by removing from $X$ one point (different from $x$). Then 
$\overset{\circ}{X}$ is homotopy-equivalent to a bouquet of $n$ circles, and 
$$\LocSys_\sG^{\on{Betti,rigid}_x}(\overset{\circ}{X})\simeq \sG^{\times n}$$
We have an isomorphism of homotopy types
$$X\simeq \overset{\circ}{X}\underset{S^1}\sqcup\, \on{pt},$$
hence
$$\LocSys_\sG^{\on{Betti,rigid}_x}(X)\simeq
\LocSys_\sG^{\on{Betti,rigid}_x}(\overset{\circ}{X})\underset{\LocSys_\sG^{\on{Betti,rigid}_x}(S^1)}\times \on{pt}
\simeq \sG^{\times n}\underset{\sG}\times \on{pt},$$
and hence is manifestly quasi-smooth. 

\sssec{}

Finally, the fact that the map
$$\LocSys_\sG^{\on{restr}}(X)\to  \LocSys_\sG^{\on{Betti}}(X)$$
is a formal isomorphism and an ind-closed embedding is the content of
\propref{p:formal compl Betti} and \thmref{t:compare Betti}. 

\sssec{}

Next we consider the de Rham context. We wish to show that the category 
$\qLisse(X)$,
which is a symmetric monoidal category over $k$ is adapted for spectral decomposition.

\medskip 

By Lefschetz principle, we can assume that $k=\BC$. In this case, by Riemann-Hilbert,
the de Rham version of $\qLisse(X)$ is equivalent as a symmetric monoidal category to 
its Betti counterpart with $\sfe=\BC$. 

\medskip 

Hence, the assertion follows from \secref{sss:spectral proof Betti}. 

\ssec{Proof of \thmref{t:spectral}, \'etale context over a field of characteristic $0$} \label{ss:spectral etale char 0}

In this subsection we will prove \thmref{t:spectral} in the \'etale context, but for $k$ being an algebraically closed field 
of characteristic $0$. 

\sssec{}

We wish to show that the symmetric monoidal category $\qLisse(X)$ is adapted for spectral decomposition. 
We will consider separately the cases when $X$ has genus $0$, and when the genus of $X$ is $\geq 1$. 

\medskip

When $X$ has genus $0$, the statement follows from the description of the category $\qLisse(X)$ in 
\secref{sss:analyze P1}. 

\medskip

Hence, from now on we will assume that the genus of $X$ is $\geq 1$. In this case, by \thmref{t:curves}(a) and \corref{c:Kpi1}, 
\begin{equation} \label{e:qLisse via Ab}
\qLisse(X)\simeq \iLisse(X) \text{ and } \Lisse(X)\simeq D^b(\Lisse(X)^\heartsuit).
\end{equation}

\sssec{}

Recall (see \cite[Expos\'e X, Corollary 1.8]{SGA1}) that if $k\hookrightarrow k'$
is an extension of algebraically closed fields, for a proper\footnote{When $k$ has characteristic $0$ (which is the case here), the assertion 
formulated below holds without the properness assumption.} scheme $Y$ over $k$ and
its base change $Y'$ to $k'$, the assignment
$$(\wt{Y}\to Y) \mapsto (\wt{Y}'\to Y'), \quad \wt{Y}':=\wt{Y}\underset{\Spec(k)}\times \Spec(k')$$
is an equivalence of categories between finite \'etale covers of $Y$ and those of $Y'$. Hence, 
pullback defines an equivalence of categories
$$\Lisse(Y)^\heartsuit\simeq \Lisse(Y')^\heartsuit.$$

Taking $Y=X$, by \eqref{e:qLisse via Ab} we obtain 
$$\qLisse(X)\simeq \qLisse(X').$$

Thus, embedding $k$ into a larger
algebraically closed field that also contains $\BC$, we can assume that $k=\BC$. 

\begin{rem} \label{r:indep of alg closed field}
The fact that pullback defines an isomorphism
$$\on{C}^\cdot_{\on{et}}(Y,\CF)\to \on{C}^\cdot_{\on{et}}(Y',\CF|_{Y'}), \quad \CF\in \Lisse(Y)$$
implies that the pullback functor 
$$\qLisse(Y)\to \qLisse(Y')$$
is fully faithful, and being an equivalence at the abelian level, is actually an equivalence for any $Y$
(not just a curve of genus $\geq 1$). 

\end{rem} 

\sssec{}

For $k=\BC$, on the one hand, we can consider the symmetric monoidal category 
$$\qLisse_{\on{et}}(X),$$
and on the other hand 
$$\qLisse_{\on{Betti}}(X),$$
for $\sfe=\ol\BQ_\ell$. 

\medskip

We will construct a (symmetric monoidal) functor
\begin{equation} \label{e:l-adic to Betti}
\jmath:\qLisse_{\on{et}}(X)\to \qLisse_{\on{Betti}}(X),
\end{equation} 
which is fully faithful, admits a continuous right adjoint. Moreover, the induced map
\begin{equation} \label{e:l-adic to Betti LocSys}
\LocSys_\sG^{\on{restr}_{\on{et}}}(X)\to \LocSys_\sG^{\on{restr}_{\on{Betti}}}(X)
\end{equation} 
will be an isomorphism on each connected component of the source (i.e., this map
identifies the source with the union of some of the connected components of the target).

\medskip

Once we show this, composing with the functor \eqref{e:qLisse to Betti again}, we will
establish that $\qLisse_{\on{et}}(X)$ is adapted for spectral decomposition in view
of what we have shown already in \secref{ss:spectral Betti and de Rham}, combined
with Theorems \ref{t:spectral Betti} and \ref{t:adapted inherited}.

\sssec{} 

Since the genus of $X$ is $\geq 1$, in addition to \eqref{e:qLisse via Ab}, we also have the equivalences
\begin{equation} \label{e:qLisse via Ab Betti}
\qLisse_{\on{Betti}}(X)\simeq \iLisse_{\on{Betti}}(X) \text{ and } \Lisse_{\on{Betti}}(X)\simeq D^b(\Lisse_{\on{Betti}}(X)^\heartsuit).
\end{equation}

Hence, in order to construct \eqref{e:l-adic to Betti}, it suffices to construct the corresponding functor  
\begin{equation} \label{e:l-adic to Betti lisse Ab}
\Lisse_{\on{et}}(X)^\heartsuit\to \Lisse_{\on{Betti}}(X)^\heartsuit.
\end{equation} 

The sought-for functor \eqref{e:l-adic to Betti lisse Ab} is obtained from the equivalence
$$\{\text{Finite \'etale covers of $X$}\} \leftrightarrow \{\text{Finite covers of the topological space underlying $X$}\}.$$

The functor \eqref{e:l-adic to Betti lisse Ab} is fully faithful. In fact, its essential image
consists of those those representations of the fundamental group of (the topological space underlying) $X$
on finite-dimensional $\ol\BQ_\ell$-vector spaces that admit
a $\ol\BZ_\ell$-lattice.

\sssec{} \label{sss:etale Betti ff}

The resulting functor $\jmath$ in \eqref{e:l-adic to Betti} maps
\begin{equation} \label{e:l-adic to Betti lisse}
\Lisse_{\on{et}}(X)\to \Lisse_{\on{Betti}}(X),
\end{equation} 
by construction. Hence, it preserves
compactness, and hence admits a continuous right adjoint.

\medskip

We claim that $\jmath$ is fully faithful. This is standard, but we will give an elementary proof for completeness:

\medskip

It is enough to show that the functor \eqref{e:l-adic to Betti lisse} is fully faithful.  
I.e., we have to show that for $\CF_1,\CF_2\in \Lisse_{\on{et}}(X)^\heartsuit$,
the maps
\begin{equation} \label{e:Hom 0}
\Hom_{\Lisse_{\on{et}}(X)^\heartsuit}(\CF_1,\CF_2)\to
\Hom_{\Lisse_{\on{Betti}}(X)^\heartsuit}(\jmath(\CF_1),\jmath(\CF_2))
\end{equation} 
\begin{equation} \label{e:Ext 1}
\Ext^1_{\Lisse_{\on{et}}(X)^\heartsuit}(\CF_1,\CF_2)\to
\Ext^1_{\Lisse_{\on{Betti}}(X)^\heartsuit}(\jmath(\CF_1),\jmath(\CF_2))
\end{equation} 
and
\begin{equation} \label{e:Ext 2}
\Ext^2_{\Lisse_{\on{et}}(X)^\heartsuit}(\CF_1,\CF_2)\to
\Ext^2_{\Lisse_{\on{Betti}}(X)^\heartsuit}(\jmath(\CF_1),\jmath(\CF_2))
\end{equation} 
are isomorphisms.

\medskip

The map \eqref{e:Hom 0} is an isomorphism since \eqref{e:l-adic to Betti lisse Ab} is fully faithful.

\medskip

To prove that \eqref{e:Ext 2} is an isomorphism, replacing $\CF_1$ by $\CF_1\otimes \CF_2^\vee$, we can assume 
that $\CF_2=\ul\sfe_X$. For $\CF\in \qLisse(X)$ (in either context) let $\CF_0$ be its maximal trivial quotient. 
Verdier duality implies that the map $$\CF\to \CF_0\simeq V\otimes \ul\sfe_X$$ defines an isomorphism 
$$V^* \otimes H^2(X,\sfe) \simeq \Hom_{\Lisse(X)}(\CF_0,\ul\sfe_X[2]) \simeq \Hom_{\Lisse(X)}(\CF,\ul\sfe_X[2])\simeq \Ext^2_{\Lisse(X)^\heartsuit}(\CF,\ul\sfe_X).$$

Since the functor \eqref{e:l-adic to Betti lisse Ab} is fully faithful, we have $\jmath(\CF)_0\simeq \jmath(\CF_0)$. This implies that
\eqref{e:Ext 2} is an isomorphism, since the functor \eqref{e:l-adic to Betti} induces an isomorphism
$$H^2_{\on{et}}(X,\sfe)\to H^2_{\on{Betti}}(X,\sfe).$$
(Indeed, both sides are 1-dimensional vector spaces, and the above map is easily seen to be non-zero.) 

\medskip

The map \eqref{e:Ext 1} is injective again because \eqref{e:l-adic to Betti lisse Ab} is fully faithful. Hence, in order to show
that it is surjective, it suffices to show that both sides have the same dimension. However, the latter follows from the 
Grothendieck-Ogg-Shafarevich formula. 

\begin{rem}
One can show that for an algebraic variety $Y$ over $\BC$, we have a well-defined
fully faithful functor 
$$\jmath:\Shv_{\on{et}}(Y)^\heartsuit\to \Shv_{\on{Betti}}(Y)^\heartsuit.$$
\end{rem}

\sssec{}

Finally, let us show that the map \eqref{e:l-adic to Betti LocSys} is an isomorphism on every connected component of
the source. 

\medskip

Given that we already know that \eqref{e:l-adic to Betti} is fully faithful, the arguments proving \propref{p:formal compl dr}
and \thmref{t:compare dR} imply that \eqref{e:l-adic to Betti LocSys} is a formal isomorphism and an ind-closed embedding. 
Hence, it suffices to show that \eqref{e:l-adic to Betti LocSys} defines a surjection at the level of $\sfe$-points from 
a connected component of $\LocSys_\sG^{\on{restr}_{\on{et}}}(X)$ to the corresponding connected component of
$\LocSys_\sG^{\on{restr}_{\on{Betti}}}(X)$.

\medskip

By \secref{ss:uniform surj}, it suffices to show that if $\sM$ is a Levi quotient of a parabolic and 
$\sigma_\sM$ is an irreducible \'etale $\sM$-local system on $X$, then the corresponding map
\begin{equation} \label{e:l-adic to Betti LocSys M}
\LocSys^{\on{restr}_{\on{et}}}_{\sP,\sigma_\sM}(X)\to \LocSys^{\on{restr}_{\on{Betti}}}_{\sP,\sigma_\sM}(X)
\end{equation}
is surjective at the level of $\sfe$-points. 

\medskip

However, as in \propref{p:comp de Rham P sigma}, the map \eqref{e:l-adic to Betti LocSys M} is in fact
an isomorphism.

\begin{rem}
The last argument shows the mechanism by which the map \eqref{e:l-adic to Betti LocSys} fails to be
an isomorphism: it misses those connected components that correspond to semi-simple Betti local
systems that do not come from the \'etale ones.
\end{rem} 

\ssec{Proof of \thmref{t:spectral}, \'etale context over a field of positive characteristic} \label{ss:spectral etale char p}

In this subsection we will show that \thmref{t:spectral} holds in the \'etale context, when 
$k$ is an algebraically closed field of positive characteristic.

\sssec{}

Let $R$ denote the ring of Witt vectors on $k$, and let $k'$ be an algebraic closure of the field of fractions of $R$.

\medskip

Since $X$ is proper and smooth of dimension $1$, it can be lifted to a smooth proper relative curve $X_{R^\wedge}$ over $\on{Spf}(R)$.

\medskip

By GAGA, $X_{R^\wedge}$ comes from a uniquely defined curve $X_R$ over $\Spec(R)$. Let $X'$ denote the base change of $X_R$
to $\Spec(k')$. 

\sssec{}

We claim that there exists a symmetric monoidal functor
$$\qLisse(X)\to \qLisse(X'),$$
with the same properties as the functor \eqref{e:l-adic to Betti}. 

\medskip

Once we show this, composing with the functors \eqref{e:l-adic to Betti} and \eqref{e:qLisse to Betti again}, we will
know that $\qLisse(X)$ is adapted for spectral decomposition by \secref{sss:spectral proof Betti}. 

\medskip

If $X$ has genus $0$, the assertion follows from the explicit description of the category $\qLisse(X)$ in 
this case, see \secref{sss:analyze P1}. Hence, we can assume that the genus of $X$ is $\geq 1$. In this case, 
it suffices to construct an exact symmetric monoidal functor at the abelian level
\begin{equation} \label{e:fund grp surj}
\Lisse(X)^\heartsuit \to \Lisse(X')^\heartsuit,
\end{equation}
and show that it is fully faithful (the fully faithfulness at the derived level will then follow by the argument in 
\secref{sss:etale Betti ff}). 

\sssec{}

Now, the existence of the functor \eqref{e:fund grp surj} with the required properties 
follows from the fact that we have a canonically defined
surjection at the level of \'etale fundamental groups
$$\pi_{1,\on{et}}(X')\to  \pi_{1,\on{et}}(X).$$

This follows from \cite[Expos\'e X, Corollary 2.3]{SGA1}.

\qed[\thmref{t:spectral}]

\ssec{A simple proof of \thmref{t:coarse restr}} \label{ss:simple coarse}

In this subsection we let $X$ be a smooth and complete algebraic curve. We will revisit
\thmref{t:coarse restr} in the case $\bH=\qLisse(X)$. 

\sssec{}

The assertion of \thmref{t:coarse restr} in the Betti context was established in
\secref{ss:coarse Betti again}.

\medskip

The assertion of \thmref{t:coarse restr} in the de Rham context follows from the Betti case
by Riemann-Hilbert.

\sssec{} \label{sss:coarse change}

Thus, it remains to treat the \'etale context. 

\medskip

However, as we have seen in Sects. \ref{ss:spectral etale char 0} and \ref{ss:spectral etale char p},
there exists a curve $X'$ over $\BC$, such that every connected component of the \'etale $\LocSys^{\on{restr}}_\sG(X)$ is 
isomorphic to a connected component of the Betti $\LocSys^{\on{restr}}_\sG(X')$ for $\sfe=\ol\BQ_\ell$. 

\bigskip

\qed[\thmref{t:coarse restr}]

\ssec{Complements: de Rham and Betti spectral actions}

In this subsection we will make a brief digression, and consider the de Rham or Betti contexts,
in which the ``usual" (i.e., not restricted) $\LocSys_\sG(X)$ is defined.  Let us be given a category $\bC$ equipped with 
an action of $\QCoh(\LocSys_\sG(X))$.

\medskip

Let $X$ be a smooth and complete curve. 
We will explicitly describe the full subcategory 
$$\QCoh(\LocSys_\sG^{\on{restr}}(X))\underset{\QCoh(\LocSys_\sG(X))}\otimes \bC \subset
\QCoh(\LocSys_\sG(X))\underset{\QCoh(\LocSys_\sG(X))}\otimes \bC=\bC,$$
where we view
\begin{equation} \label{e:qcoh restr via amb}
\QCoh(\LocSys_\sG^{\on{restr}}(X))\simeq \QCoh(\LocSys_\sG(X))_{\LocSys_\sG^{\on{restr}}(X)}
\end{equation}
as a co-localization of $\QCoh(\LocSys_\sG(X))$, see \secref{sss:union formal compl again}. 

\sssec{} \label{sss:fin mon}

Let is first specialize to the Betti context. Consider the algebraic stack $\QCoh(\LocSys^{\on{Betti}}_\sG(X))$.

\medskip

For a given $V\in \Rep(\sG)$, let 
$$\CE_V\in \QCoh(\LocSys_\sG(X))\otimes \Shv^{\on{all}}_{\on{loc.const.}}(X)$$
be as in \secref{sss:pre-shtuka}. 

\medskip

Let $\bC$ be a DG category, equipped with an action of $\QCoh(\LocSys_\sG(X))$. 
In particular, for
$V\in \Rep(\sG)$, we have the functor
$$\on{H}(V,-):\bC \to \bC\otimes \Shv^{\on{all}}_{\on{loc.const.}}(X),$$ 
corresponding to the action of the object $\CE_V$ above.

\medskip

Let 
$$\bC^{\on{fin.mon.}}\subset \bC$$
be the full subcategory consisting of objects $\bc\in \bC$, for
which
$$\on{H}(V,\bc)\in \bC\otimes \qLisse(X) \subset \bC\otimes \Shv^{\on{all}}_{\on{loc.const.}}(X).$$

\medskip

As in \cite[Proposition C.2.5]{GKRV}, one shows that the category $\bC^{\on{fin.mon.}}$ is stable under the action of 
$\QCoh(\LocSys_\sG(X))$ and it carries an action of $\Rep(\sG)^{\otimes X\on{-lisse}}$.

\sssec{}

We claim:

\begin{prop} \label{p:mon fin subcategory}
The full subcategory $\bC^{\on{fin.mon.}}\subset \bC$ equals
$$\QCoh(\LocSys_\sG^{\on{restr}}(X))\underset{\QCoh(\LocSys^{\on{Betti}}_\sG(X))}\otimes \bC \subset
\QCoh(\LocSys^{\on{Betti}}_\sG(X))\underset{\QCoh(\LocSys^{\on{Betti}}_\sG(X))}\otimes \bC=\bC.$$
\end{prop} 

\begin{proof}

The inclusion 
$$\QCoh(\LocSys_\sG^{\on{restr}}(X))\underset{\QCoh(\LocSys^{\on{Betti}}_\sG(X))}\otimes \bC \subset \bC^{\on{fin.mon.}}$$
is clear.

\medskip

For the opposite inclusion, we can assume that $\bC^{\on{fin.mon.}}=\bC$, and we will need to show that 
$$\QCoh(\LocSys_\sG^{\on{restr}}(X))\underset{\QCoh(\LocSys^{\on{Betti}}_\sG(X))}\otimes \bC \to \bC$$
is an equivalence. 

\medskip

The assumption on $\bC$ implies that the action of $\Rep(\sG)^{\otimes X}$ factors through an action 
of $\Rep(\sG)^{\otimes X\on{-lisse}}$. Hence, by \thmref{t:spectral}, the action of 
$\QCoh(\LocSys^{\on{Betti}}_\sG(X))$ on $\bC$ factors through $\QCoh(\LocSys^{\on{restr}}_\sG(X))$.

\medskip

Hence,
$$\QCoh(\LocSys_\sG^{\on{restr}}(X))\underset{\QCoh(\LocSys^{\on{Betti}}_\sG(X))}\otimes \bC \simeq$$
$$\simeq \QCoh(\LocSys_\sG^{\on{restr}}(X))\underset{\QCoh(\LocSys^{\on{Betti}}_\sG(X))}\otimes 
\QCoh(\LocSys_\sG^{\on{restr}}(X)) \underset{\QCoh(\LocSys^{\on{restr}}_\sG(X))}\otimes \bC,$$
while 
$$\QCoh(\LocSys_\sG^{\on{restr}}(X))\underset{\QCoh(\LocSys^{\on{Betti}}_\sG(X))}\otimes 
\QCoh(\LocSys_\sG^{\on{restr}}(X)) \simeq \QCoh(\LocSys_\sG^{\on{restr}}(X)),$$
since $\QCoh(\LocSys_\sG^{\on{restr}}(X))$ is a monoidal co-localization of 
$\QCoh(\LocSys_\sG^{\on{Betti}}(X))$.

\end{proof} 

\sssec{} \label{sss:action de Rham}

Let us now specialize to the de Rham context. 

\medskip

For a given $V\in \Rep(\sG)$ consider the corresponding object 
$$\CE_V\in \QCoh(\LocSys^{\dr}_\sG(X))\otimes \Dmod(X),$$
see \secref{sss:pre-shtuka}. 

\medskip

Let $\bC$ be a DG category, equipped with an action of $\QCoh(\LocSys^{\dr}_\sG(X))$. 
In particular, for
$V\in \Rep(\sG)$, we have the functor
$$\on{H}(V,-):\bC \to \bC\otimes \Dmod(X),$$ 
corresponding to the action of the object $\CE_V$ above.

\medskip

Let 
$$\bC^{\on{Lisse}}\subset \bC$$
be the full subcategory consisting of objects $\bc\in \bC$, for
which
$$\on{H}(V,\bc)\in \bC\otimes \qLisse(X) \subset \bC\otimes \Dmod(X).$$

\medskip

As in \cite[Proposition C.2.5]{GKRV}, one shows that the category $\bC^{\on{Lisse}}$ is stable under the action of 
$\QCoh(\LocSys^{\dr}_\sG(X))$ and it carries an action of $\Rep(\sG)^{\otimes X\on{-lisse}}$.

\sssec{}

We claim: 

\begin{prop} \label{p:lisse subcategory}
The full subcategory $\bC^{\on{Lisse}}\subset \bC$ equals
$$\QCoh(\LocSys_\sG^{\on{restr}}(X))\underset{\QCoh(\LocSys_\sG(X))}\otimes \bC \subset
\QCoh(\LocSys^{\dr}_\sG(X))\underset{\QCoh(\LocSys^{\dr}_\sG(X))}\otimes \bC=\bC.$$
\end{prop} 

The proof repeats that of \propref{p:mon fin subcategory}. 

\begin{rem}
A statement somewhat weaker than \propref{p:lisse subcategory} appeared in \cite{GKRV} as 
Conjecture C.5.5 of {\it loc. cit.} 
\end{rem}

\begin{rem} \label{r:rep G Ran dr prel} 

Consider the category 
$$\coHom(\Rep(\sG),\Dmod(X)).$$

This is the category that appears, e.g., in \cite[Sect. 4.2.7]{Ga7}; in this paper
we denote it\footnote{Our version of $\Rep(\sG)^\dr_{\Ran}$ is a slightly different from the one
in \cite[Sect. 4.2.7]{Ga7} in that it is the unital version of the category 
considered in {\it loc.cit.}}
$$\Rep(\sG)^\dr_{\Ran},$$
see Remark \ref{r:rep G Ran dr}.

\medskip

The corresponding functor
$$\Rep(\sG)^\dr_{\Ran}\simeq \coHom(\Rep(\sG),\Dmod(X))\to \QCoh(\LocSys^{\dr}_\sG(X))$$
is a localization (i.e., admits fully faithful right adjoint, which is, moreover, continuous),
see \cite[Proposition 4.3.4]{Ga7} . 

\medskip

However, it is \emph{not} an equivalence. Hence, the arrow
$$\{\text{Actions of $\QCoh(\LocSys^{\dr}_\sG(X))$ on $\bC$} \} \to 
\{\text{Actions of $\Rep(\sG)_{\Ran}$ on $\bC$} \}$$
is fully faithful, but not an equivalence.

\medskip

A key result of \cite[Corollary 4.5.5]{Ga7} says that for $\bC=\Dmod(\Bun_G)$, the
action of $\Rep(\cG)_{\Ran}$ given by Hecke functors lies in the essential image of the above map. 

\end{rem} 

\section{Other examples of categories adapted for spectral decomposition} \label{s:Lie}

The contents of this section are \emph{not} needed for the rest of the paper. 

\medskip

We will find another class of symmetric monoidal categories adapted for spectral decomposition,
namely, categories of modules over connective Lie algebras. 

\medskip

The method of proof will allow us to give an alternative argument also for the proof \thmref{t:spectral Betti}, 
and potentially sheds some light on the nature of the ``adapted for spectral decomposition" condition. 

\ssec{The case of Lie algebras}  \label{ss:Lie}

We are going to establish a variant of \thmref{t:spectral Betti}, 
where instead of an object $\CX\in \Spc$ we have a Lie algebra $\sh\in \on{LieAlg}(\Vect^{\leq 0}_\sfe)$. 

\sssec{}

Consider the category
$$\bH:=\sh\mod.$$

It carries a symmetric monoidal structure (given by tensor product of modules over $\sh$), and a fiber functor
$$\oblv_\sh:\sh\mod\to \Vect_\sfe,$$
given by forgetting the action of the Lie algebra. Since $\sh$ was assumed connective, the category 
$\sh\mod$ carries a t-structure, for which $\oblv_\sh$ is t-exact. This t-structure is left-complete: indeed
$$\sh\mod \simeq U(\sh)\mod,$$
and it is is known that the category of modules over a connective associative algebra is left-complete
in its t-structure. 

\medskip

Hence, $\sh\mod$ is a category that satisfies the requirements of \secref{sss:pre-Tannakian}.
Furthermore, $\sh\mod$ is dualizable (in fact, compactly generated). 

\sssec{}

We will prove:

\begin{thm} \label{t:Lie adapted}
The category $\sh\mod$ is adapted for spectral decomposition.
\end{thm} 

The proof will occupy the next two subsections. As a first step, we will reinterpret the prestack $\bMaps(\Rep(\sG),\sh\mod)$.

\begin{rem} \label{r:Lie non-Tann}
Note that when $H^0(\sh)$ is nilpotent, \thmref{t:Lie adapted} is a particular case of \thmref{t:spectral Betti}:
indeed, rational homotopy type theory implies that there exists a pointed space $\CX$ such that the 
pair $(\Vect_\sfe^\CX,\on{ev}_x)$ is equivalent to $(\sh\mod,\oblv_\sh)$. 
\end{rem}

\ssec{The space of maps of Lie algebras}

As a first step towards the proof of \thmref{t:Lie adapted}, we will reinterpret the prestack
$\bMaps(\Rep(\sG),\sh\mod)$, or rather its version $\bMaps(\Rep(\sG),\sh\mod)^{\on{rigid}}$,
as the space of maps of Lie algebras. 

\sssec{} \label{sss:maps of Lie algs}

Let $\sg$ denote the Lie algebra of $\sG$. Consider the prestack, denoted $\bMaps_{\Lie}(\sh,\sg)$,
that sends an affine scheme $S=\Spec(A)$ to the space of maps of Lie algebras in $A\mod$. 
$$\sh\otimes A\to \sg\otimes A.$$

%
%
%
%
%
%

\sssec{} \label{sss:maps of Lie algs old}

Let $S$ be an affine scheme, and let us be given an $S$-point of $\bMaps_{\Lie}(\sh,\sg)$.
It gives rise to a symmetric monoidal functor
$$\sg\mod\to \QCoh(S)\otimes \sh\mod,$$
such that the composition
$$\sg\mod\to \QCoh(S)\otimes \sh\mod \overset{\on{Id}\otimes \oblv_\sh}\longrightarrow \QCoh(S)$$
identifies with
$$\sg\mod \overset{\oblv_\sg}\longrightarrow \Vect_\sfe \overset{\on{unit}}\longrightarrow \QCoh(S).$$

\medskip

We have the (symmetric monoidal) restriction functor
\begin{equation} \label{e:retsr from G to g}
\Rep(\sG)\to \sg\mod,
\end{equation} 
which commutes with the forgetful functors to $\Vect_\sfe$.

\medskip

Composing, we obtain a map of prestacks
\begin{equation} \label{e:Lie}
\bMaps_{\Lie}(\sh,\sg)\to \bMaps(\Rep(\sG),\sh\mod)^{\on{rigid}}.
\end{equation}

We claim:

\begin{prop} \label{p:Lie maps}
The map \eqref{e:Lie} is an isomorphism.
\end{prop} 

\sssec{}

For the proof of \propref{p:Lie maps} (as well as that of \propref{p:Lie} below) we recall that any object in $\on{LieAlg}(\Vect^{\leq 0}_\sfe)$ can 
be written as a \emph{sifted} colimit of objects of the form
\begin{equation} \label{e:free Lie}
\free_{\on{Lie}}(V),\quad V\in \Vect_\sfe^{\leq 0}.
\end{equation} 

Note also that for $\sh$ as in \eqref{e:free Lie}, we have:
\begin{equation} \label{e:map out of free}
\bMaps_{\Lie}(\free_{\on{Lie}}(V),\sg)\simeq \Spec(\Sym(V\otimes \sg^\vee))
\end{equation}
is an affine scheme. 

\medskip

We also note the following lemma:

\begin{lem} \label{l:colim Lie} 
The assignment 
$$\sh\mapsto \sh\mod, \quad \on{LieAlg}(\Vect_\sfe)\to \DGCat$$
sends sifted colimits to limits.
\end{lem}

\begin{proof}

Consider the functor of universal enveloping algebra
$$\on{LieAlg}(\Vect_\sfe)\to \on{AssocAlg}(\Vect_\sfe).$$

Being a left adjoint, this functor sends colimits to colimits. We have
$$\sh\mod\simeq U(\sh)\mod.$$

Now, the assertion follows from \cite[Lemma 2.5.5]{GKRV}.

\end{proof} 

\sssec{Proof \propref{p:Lie maps}}

We need to show that for an affine scheme $S$, the map
\begin{equation} \label{e:Lie S}
\Maps(S,\bMaps_{\Lie}(\sh,\sg))\to \Maps(S,\bMaps(\Rep(\sG),\sh\mod)^{\on{rigid}})
\end{equation}
is an isomorphism. 

\medskip

The left-hand side in \eqref{e:Lie S} sends colimits in $\sh$ to limits in $\Spc$. The right-hand side
in \eqref{e:Lie S}, sends sifted colimits in $\sh$ to limits in $\Spc$, by \lemref{l:colim Lie}. Hence, we 
can assume that $\sh$ is of the form \eqref{e:free Lie}. Moreover, we can assume that 
$$V\in \Vect_\sfe^{\leq 0}\cap \Vect^c_\sfe.$$

Note that, by \eqref{e:map out of free}, for $S=\Spec(A)$, 
$$\Maps(S,\bMaps_{\Lie}(\free_{\on{Lie}}(V),\sg))\simeq \Maps_{\Vect_\sfe}(V,A\otimes \oblv_{\Lie}(\sg)).$$

\medskip

Let $A'$ be the split square-zero extension of $A$ equal to
$$A'=A\otimes (\sfe \oplus \epsilon\cdot V^*), \quad \epsilon^2=0.$$

Note that the category 
$$\QCoh(S)\otimes \free_{\on{Lie}}(V)\mod$$
identifies as a symmetric monoidal category with the category of triples $(\CM,t,\alpha)$, where:

\begin{itemize}

\item $\CM\in A\mod$;

\smallskip

\item $t$ is an automorphism of $\CM':=A'\underset{A}\otimes \CM$;

\item $\alpha$ is a trivialization of the induced automorphism on $\CM\simeq A\underset{A'}\otimes \CM'$
induced by $t$.

\end{itemize}

Hence, the space
$$\Maps\left(S,\bMaps(\Rep(\sG),\free_{\on{Lie}}(V)\mod)^{\on{rigid}}\right)$$
identifies with the space of automorphisms of the symmetric monoidal functor
$$\Rep(\sG) \overset{\oblv_\sG}\longrightarrow \Vect_\sfe \overset{A'}\longrightarrow A'\mod,$$ 
equipped with the trivialization of the automorphism of the composite functor 
$$\Rep(\sG) \overset{\oblv_\sG}\longrightarrow \Vect_\sfe \overset{A'}\longrightarrow A'\mod\to A\mod.$$ 

\medskip

By Tannaka duality, the latter is the same as the space of maps 
$$\Spec(A')\to \sG,$$
equipped with the trivialization of the composite
$$S=\Spec(A)\to \Spec(A')\to \sG.$$

By deformation theory, we rewrite the latter as 
$$\Maps_{\Vect_\sfe}(\sg^\vee,A\otimes V^*)\simeq \Maps_{\Vect_\sfe}(V,A\otimes \oblv_{\Lie}(\sg)).$$

Unwinding the above identifications, it is easy to see that the map \eqref{e:Lie S} corresponds to the
identity map on $\Maps_{\Vect_\sfe}(V,A\otimes \oblv_{\Lie}(\sg))$; in particular, it is an isomorphism.

\qed[\propref{p:Lie maps}]

\sssec{}

We now record the following: 

\begin{prop} \hfill \label{p:Lie} 

\smallskip

\noindent{\em(a)} The prestack $\bMaps_{\Lie}(\sh,\sg)$ is an affine scheme.

\smallskip

\noindent{\em(b)} The affine scheme $\bMaps_{\Lie}(\sh,\sg)$ is almost of finite type
if $\sh$ is finite-dimensional in each degree.

\end{prop}

\begin{proof}

Since a limit of affine schemes is an affine scheme, for the proof of point (a) 
we can assume that $\sh$ is of the form \eqref{e:free Lie}. Then the statement 
follows from \eqref{e:map out of free}.

\medskip

For the proof of point (b) we argue as follows. Assume that $\sh$ is finite-dimensional in each degree. Then
$^{\on{cl}}\bMaps_{\Lie}(\sh,\sg)$ is the classical scheme that classifies maps of (classical) Lie algebras $H^0(\sh)\to \sg$;
in particular it is a closed subscheme in the affine space of the vector space Hom 
from $H^0(\sh)$ to $\fg$, i.e., 
$$\on{Tot}(\Hom_{\Vect_\sfe}(\oblv_{\Lie}(H^0(\sh)),\oblv_{\Lie}(\sg))),$$
and so is of finite type. By \cite[Chapter 1, Theorem 9.1.2]{GR2}, it remains to show that for a classical scheme $S$ of finite type 
and an $S$-point $\sF$ of $\bMaps_{\Lie}(\sh,\sg)$, the cotangent space 
$$T^*_\sF(\bMaps_{\Lie}(\sh,\sg))\in \QCoh(S)^{\leq 0}$$
has coherent cohomologies. 

\medskip

Using \propref{p:Lie maps} above and \corref{c:def}(b), we obtain that for $S=\Spec(A)$ and an $S$-point 
$\phi\in \bMaps_{\Lie}(\sh,\sg)$, 
$$T^*_\phi(\bMaps_{\Lie}(\sh,\sg))\simeq \on{Fib}\left(\sg^\vee\otimes A\to \on{C}_\cdot(\sh,\sg^\vee\otimes A)\right),$$
where $\sg^\vee\otimes A$ is viewed as a $\sh$-module via $\phi$. This implies the required assertion as $\on{C}_\cdot(\sh,-)$
can be computed by the standard Chevalley complex. 

\end{proof}

Combining with \propref{p:Lie maps}, we obtain:

\begin{cor} \hfill \label{c:Lie} 

\smallskip

\noindent{\em(a)} The prestack $\bMaps(\Rep(\sG),\sh\mod)^{\on{rigid}})$ is an affine scheme.

\smallskip

\noindent{\em(b)} The affine scheme $\bMaps(\Rep(\sG),\sh\mod)^{\on{rigid}})$ is almost of finite type
if $\sh$ is finite-dimensional in each degree.

\end{cor}

\ssec{Proof of \thmref{t:Lie adapted}} \label{ss:proof Lie}

We are now ready to prove \thmref{t:Lie adapted}. 

\sssec{}

We are going to prove an equivalent statement, namely that the functor
$$\coHom(\Rep(\sG),\sh\mod)^{\on{rigid}}\to \QCoh(\bMaps(\Rep(\sG),\sh\mod)^{\on{rigid}})$$
is an equivalence.

\medskip

I.e., we have to show that for a target symmetric monoidal category $\bA$, restriction along 
\eqref{e:Hom abs rigid} defines an equivalence from the space of symmetric monoidal functors
\begin{equation} \label{e:LHS Lie}
\QCoh(\bMaps(\Rep(\sG),\sh\mod)^{\on{rigid}})\to \bA
\end{equation} 
to the space of symmetric monoidal functors
\begin{equation} \label{e:RHS Lie}
\Rep(\sG)\to \bA\otimes \sh\mod,
\end{equation}
equipped with an identification of the composition
$$\Rep(\sG)\to \bA\otimes \sh\mod \overset{\on{Id}_\bA\otimes \oblv_\sh}\longrightarrow \bA$$
with
$$\Rep(\sG)\overset{\oblv_\sG}\longrightarrow \Vect_\sfe \overset{\one_\bA}\to \bA.$$

\sssec{Step 1}  \label{sss:step 1 Lie}

We will first show that we can assume that $\bA$ is of the form $A\mod$ for some $A\in \on{ComAlg}(\Vect_\sfe)$. 
Namely, we will show that both \eqref{e:LHS Lie} and \eqref{e:RHS Lie} factor canonically via $A\mod$,
where
$$A:=\CEnd_\bA(\one_\bA).$$

For \eqref{e:LHS Lie} this follows from \corref{c:Lie}(a): for any affine scheme $Y=\Spec(R)$, symmetric monoidal functors
$$\QCoh(Y)=R\mod\to \bA$$
are in bijection with maps of commutative algebras $R\to \CEnd_\bA(\one_\bA)=:A$, and the latter are the
same as symmetric monoidal functors
$$R\mod\to A\mod.$$

\medskip

For \eqref{e:RHS Lie} we argue as follows. We have a tautological (symmetric monoidal) functor
\begin{equation} \label{e:C' to C}
A\mod =:\bA' \to \bA,
\end{equation} 
(it is not necessarily fully faithful because $\one_\bA\in \bA$ is not necessarily compact). We claim that \eqref{e:C' to C}
induces an isomorphism between the data of \eqref{e:RHS Lie} for $\bA'$ and $\bA$, respectively. 

\medskip

Indeed, the datum of \eqref{e:RHS Lie} for $\bA$ (or $\bA'$) 
amounts to defining an action of $\sh$ on 
\begin{equation} \label{e:actions Lie}
\ul{V}\otimes \one_\bA
\end{equation} 
(viewed as an object of either $\bA$ or $\bA'$) for every $V\in \Rep(\sG)^c$, in a way compatible with the tensor structure
(here $V\to \ul{V}$ denotes the forgetful functor $\Rep(\sG)\to \Vect_\sfe$). 

\medskip

The datum of such action consists of a compatible family of diagrams
\begin{equation} \label{e:actions Lie IJ}
U(\sh)^{\otimes I}\otimes \ul{V}\otimes \one_\bA\to \ul{V}\otimes \one_\bA, \quad V\in \Rep(\sG^J)^c,\,\, I,J\in\on{fSet}.
\end{equation} 
The assertion that \eqref{e:C' to C} induces an isomorphism on the data \eqref{e:RHS Lie}
follows now from the fact that the functor \eqref{e:C' to C} \emph{does} induce an isomorphisms on the mapping
space from objects of the form
$$W_1\otimes \one_\bA,\quad W_1\in \Vect_\sfe$$
to objects of the form
$$W_2\otimes \one_\bA,  \quad W_2\in \Vect^c_\sfe.$$

\sssec{Step 2} \label{sss:step 2 Lie}

Thus, we can assume that $\bA=A\mod$ for $A\in \on{ComAlg}(\Vect_\sfe)$. Next we claim that we 
can assume that $A$ is connective. More precisely, we claim that \eqref{e:LHS Lie} and \eqref{e:RHS Lie} factor 
canonically via $A'\mod$, where $A':=\tau^{\leq 0}(A)$. 

\medskip

This is again obvious for \eqref{e:LHS Lie}: for $R\in \on{ComAlg}(\Vect^{\leq 0}_\sfe)$, 
a map $R\to A$ factors canonically through a map $R\to A'$. 

\medskip

For \eqref{e:RHS Lie} we argue as follows: since $\Rep(\sG)$ is the derived category of its heart and the 
tensor product operation is t-exact, in \eqref{e:actions Lie} we can assume that $V\in \Rep(\sG)^\heartsuit\cap \Rep(\sG)^c$.
Hence, in \eqref{e:actions Lie IJ} we can also assume that
$$V\in  \Rep(\sG^J)^\heartsuit\cap  \Rep(\sG^J)^c.$$ 

Now, in this case, maps in  \eqref{e:actions Lie IJ}, which correspond to points in
$$\Maps_{\Vect_\sfe}(U(\sh)^{\otimes I}\otimes \ul{V},A\otimes \ul{V})$$
factor canonically via 
$$\Maps_{\Vect_\sfe}(U(\sh)^{\otimes I}\otimes \ul{V},A'\otimes \ul{V}).$$

\sssec{Step 3} 

Thus, we can assume that $\bA=\QCoh(S)$, where $S$ is an affine scheme. However, in this case, 
the spaces \eqref{e:LHS Lie} and  \eqref{e:RHS Lie} are just the same. 

\qed[\thmref{t:Lie adapted}]

\ssec{Back to the Betti case} \label{ss:another Betti}

We will now show how to prove \thmref{t:spectral Betti}, along the lines of the proof of \thmref{t:Lie adapted}. 

\sssec{}

Let $\CX$ be a connected object of $\Spc$, and let $x\in \CX$ be a base point.

\medskip

It suffices to show that functor
$$\Rep(\sG)^{\otimes \CX}_\sfe \underset{\Rep(\sG)}\otimes \Vect_\sfe\to \QCoh(\LocSys^{\on{Betti}}_\sG(X)^{\on{rigid}_x})$$
is an equivalence.

\medskip

I.e., we have to show that for a target symmetric monoidal category $\bA$, the space of symmetric monoidal functors
\begin{equation} \label{e:LHS Betti}
\QCoh(\LocSys^{\on{Betti}}_\sG(X)^{\on{rigid}_x})\to \bA
\end{equation}
maps isomorphically to the space of symmetric monoidal functors
\begin{equation} \label{e:RHS Betti}
\Rep(\sG) \to \bA\otimes \Vect^\CX_\sfe
\end{equation}
equipped with an identification of the composition
$$\Rep(\sG) \to \bA\otimes \Vect^\CX_\sfe \overset{\on{Id}\otimes \on{ev}_x}\longrightarrow \bA$$
with
$$\Rep(\sG) \overset{\oblv_\sG}\to \Vect_\sfe\overset{\one_\bA}\to\bA.$$

\sssec{}

As in the proof of \thmref{t:Lie adapted}, it suffices to show that we can replace the category $\bA$
by a category $A\mod$, where $A$ is a connective commutative algebra.

\medskip

For \eqref{e:LHS Betti}, this follows by the same argument as in \secref{ss:proof Lie}, since 
$\LocSys^{\on{Betti}}_\sG(X)^{\on{rigid}_x}$ is an affine scheme (by \propref{p:LocSys Betti}). 

\medskip

For \eqref{e:RHS Betti} we argue as follows.

\sssec{}

Let 
$$\Omega(\CX,x)\in \on{Grp}(\Spc)$$ be the loop space of $\CX$ based at $x$. 
We interpret $\Vect_\sfe^\CX$ as 
$$\Omega(\CX,x)\mod(\Vect_\sfe).$$

\medskip

Then the datum in \eqref{e:RHS Betti} amounts to the datum of tensor-compatible collection
of actions of $\Omega(\CX,x)$ on the objects 
$$\ul{V}\otimes \one_\bA, \quad V\in \Rep(\sG).$$

Then the arguments in Sects. \ref{sss:step 1 Lie}-\ref{sss:step 2 Lie} apply verbatim, allowing to replace
$$\bA\rightsquigarrow \CEnd_\bA(\one_\bA)$$ and
$$A \rightsquigarrow \tau^{\leq 0}(A).$$

\qed

\begin{rem}
Note that an analog of \propref{p:Lie maps} holds also in the present context: we can interpret (the affine scheme)
$\LocSys^{\on{Betti}}_\sG(X)^{\on{rigid}_x}$ as the affine scheme 
$$\bMaps_{\on{Grp}}(\Omega(\CX,x),\sG),$$ 
where
$$\Maps(S,\bMaps_{\on{Grp}}(\Omega(\CX,x),\sG))):=
\Maps_{\on{Grp}}(\Omega(\CX,x),\Maps_{\affSch}(S,\sG)).$$

This follows from the fact that for $S=\Spec(A)$ and $\CY\in \Spc$, the datum of a system tensor-compatible maps
$$\CY\mapsto \on{Aut}_{A\mod}(\ul{V}\otimes A), \quad V\in \Rep(\sG),$$
is equivalent to the datum of a map 
$$\CY\to \Maps_{\affSch}(S,\sG),$$
by Tannaka duality. 

\end{rem}

\section{Ran version of $\Rep(\sG)$ and Beilinson's spectral projector} \label{s:progenitor projector}

In this section and the next sections we develop a tool that we will use in the sequel in order to
produce Hecke eigensheaves. 

\medskip

This tool is Beilinson's spectral projector, which is an object of the \emph{Ran version} of the
category $\Rep(\sG)$. 


\ssec{The category $\Rep(\sG)_{\Ran}$} \label{ss:C Ran}

In this subsection we introduce the Ran version of the category $\Rep(\sG)$,
to be denoted $\Rep(\sG)_{\Ran}$. 

\sssec{}

Let $X$ be an arbitrary scheme, and let $\CC$ be a symmetric monoidal category. 

\medskip

Recall that $\on{fSet}$ denotes the category of finite sets, and $\on{TwArr}(\on{fSet})$ its twisted
arrows category.

\medskip

Consider the functor 
\begin{equation}\label{e:Ran diagram}
\on{TwArr}(\on{fSet})\to \DGCat
\end{equation}
that at the level of objects sends
$$(I\to J) \mapsto \CC^{\otimes I}\otimes \Shv(X^J).$$

At the level of morphisms,  for a map
\begin{equation} \label{e:mor Tw arr}
\CD
I_1 @>>> J_1 \\
@V{\phi_I}VV @AA{\phi_J}A \\
I_2 @>>> J_2,
\endCD
\end{equation} 
in $\on{TwArr}(\on{fSet})$,
the corresponding functor
\begin{equation} \label{e:transition functors}
\CC^{\otimes I_1}\otimes \Shv(X^{J_1})\to \CC^{\otimes I_2}\otimes \Shv(X^{J_2})
\end{equation} 
is given by the tensor product functor along the fibers of $\phi_I$
\begin{equation} \label{e:mult phi}
\on{mult}_\CC^{\phi_I}:\CC^{\otimes I_1}\to \CC^{\otimes I_2}
\end{equation} 
and the functor
\begin{equation} \label{e:Delta phi}
(\Delta_{\phi_J})_*:\Shv(X^{J_1})\to \Shv(X^{J_2}),
\end{equation} 
where $\Delta_{\phi_J}:X^{J_1}\to X^{J_2}$ is the diagonal map induced by $\phi_J$. 

\sssec{}

We define the key actor in section, the category $\CC_\Ran$, as the colimit of the functor 
\eqref{e:Ran diagram}.

\medskip

Our main example of interest is when $\CC=\Rep(\sG)$, where $\sG$ is an algebraic group.

\begin{rem}
Note that the definition of $\CC_\Ran$ makes sense without the assumption that $X$ be proper.
\end{rem}

\sssec{}

Let $(I \to J) \in \on{TwArr}(\fSet)$ be given. We will denote by 
$$
\on{ins}_{I\to J}:\CC^{\otimes I} \otimes \Shv(X^J) \to 
\CC_\Ran
$$ 
the corresponding functor. 

\sssec{}

The functor \eqref{e:Ran diagram} is naturally right-lax symmetric monoidal.  Therefore, the colimit $\CC_\Ran$
carries a natural symmetric monoidal structure.  Explicitly, this symmetric monoidal structure can be described
as follows.
For
$$V_1\otimes \CF_1\in \CC^{\otimes I_1}\otimes \Shv(X^{J_1}) \text{ and }
V_2\otimes \CF_2\in \CC^{\otimes I_2}\otimes \Shv(X^{J_2}),$$
the tensor product of their images in $\CC_\Ran$ is the image of the object
$$(V_1\otimes V_2)\otimes (\CF_1\boxtimes \CF_2)\in \CC^{\otimes (I_1\sqcup I_2)}\otimes \Shv(X^{J_1\sqcup J_2}).$$

\medskip

We will denote the resulting monoidal operation on $\CC_\Ran$ by
$$\CV_1,\CV_2\mapsto \CV_1\star \CV_2.$$

\medskip

We denote the unit object by $\one_{\CC_\Ran}$. It is given by
$$\on{ins}_{\emptyset \to \emptyset}(\sfe), \quad \sfe\in \Vect_\sfe\simeq \CC^{\emptyset}\otimes \Shv(X^\emptyset).$$

\sssec{Example}

Let $\CC=\Vect_\sfe$. Then $\CC_\Ran\simeq \Vect$. For example, an object of the form
$$\on{ins}_\psi(\sfe\otimes \CF), \quad \CF\in \Shv(X^J)$$
is canonically isomorphic to 
$$\on{C}^\cdot(X^J,\CF)\otimes \one_{\CC_\Ran}$$
via the following morphisms in $\on{TwArr}(\fSet)$: 
$$
\CD
I @>{\psi}>> J \\
@AAA @VV{\on{id}}V \\
\emptyset @>>>  J \\
@VVV @AAA \\
\emptyset @>>> \emptyset.
\endCD
$$

\sssec{} \label{sss:maps out of Ran new}

As in Sects. \ref{sss:coEnd} and \ref{sss:proof tw Arr}, given a target symmetric monoidal/monoidal/plain DG category $\bA$
there is a naturally defined map

\begin{itemize}

\item From the space of natural transformation of functors from $\fSet$ to $$\DGCat^{\on{SymMon}}/\DGCat^{\on{Mon}}/\DGCat,$$
$$(I\mapsto \CC^{\otimes I})\,\to \,(I\mapsto \bA\otimes \Shv(X^I)),$$

\item To the space of symmetric monoidal/monoidal/plain continuous functors $\CC_\Ran\to \bA$. 

\end{itemize}

\noindent (Here we use the fact that the DG category $\Shv(X^I)$ (or in fact $\Shv(Y)$ on any scheme $Y$)
is canonically self-dual by means of Verdier duality, see \secref{sss:duality on stack}.)

\medskip 

As in \lemref{l:maps from coend}(b), the above map is an isomorphism for plain DG categories,
However, unlike \lemref{l:maps from coend}(a), this map fails to be an isomorphism in the associative and 
commutative cases. 

\begin{rem}

Along with $\CC_\Ran$, one can also consider the category
$$\coHom(\CC,\Shv(X)),$$
where $\Shv(X)$ is viewed as a symmetric monoidal category via the $\sotimes$ operation.

\medskip

The duals (which are also the right adjoints) of the functors
$$\Shv(X)^{\otimes J}\overset{\boxtimes}\to \Shv(X^J)$$
define a symmetric monoidal functor
\begin{equation} \label{e:Ran to coHom}
\CC_\Ran\to \coHom(\CC,\Shv(X)).
\end{equation}

\end{rem} 

\begin{rem}  \label{r:rep G Ran dr} 

We can apply the construction of $\CC_\Ran$ verbatim, when instead of $\Shv(X)$ we use the
category $\Dmod(X)$ (when the ground field $k$ has characteristic $0$); denote the resulting 
symmetric monoidal category by $\CC^\dr_\Ran$. 

\medskip

However, in this case, the functors
$$\Dmod(X)^{\otimes J}\overset{\boxtimes}\to \Dmod(X^J)$$
are equivalences. Hence, the counterpart of the functor \eqref{e:Ran to coHom}
$$\CC^\dr_\Ran\to \coHom(\CC,\Dmod(X))$$
is an equivalence. 

\medskip

In a constructible de Rham context, the (fully faithful) functors $\Shv(-)\to \Dmod(-)$ 
induce a fully faithful symmetric monoidal functor
$$\CC_\Ran\to \CC^\dr_\Ran.$$

%
%
%
%

\end{rem}

\begin{rem}  \label{r:rep G Ran Betti} 

Similarly to Remark \ref{r:rep G Ran dr}, we can also consider the category $\CC^{\on{Betti}}_\Ran$,
where now instead of $\Shv(X)$ we use $\Shv^{\on{all}}(X)$ and instead of the functors $(\Delta_\phi)_*$
(which fail to be continuous unless $X$ is proper), we use the functors $(\Delta_\phi)_!$. 

\medskip

In this case, the functors
$$\Shv^{\on{all}}(X)^{\otimes J}\overset{\boxtimes}\to \Shv^{\on{all}}(X^J)$$
are also equivalences (see \secref{sss:prod Betti}).  

\medskip 

The categories $\Shv^{\on{all}}(X^J)$ (or more generally, $\Shv^{\on{all}}(Y)$ on any finite CW complex $Y$)
are also canonically self-dual via the pairing
$$\CF_1,\CF_2\mapsto \on{C}_c^\cdot(Y,\CF_1\overset{*}\otimes \CF_2).$$

With respect to this duality, the functors $(\Delta_\phi)_!$ are dual to the functors 
$(\Delta_\phi)^*$. So, the (symmetric monoidal) category $\CC^{\on{Betti}}_\Ran$
identifies with
$$\coHom(\CC,\Shv^{\on{all}}(X)),$$
where $\Shv^{\on{all}}(X)$ is viewed as a symmetric monoidal category via the 
$\overset{*}\otimes$ operation.

%

\end{rem}

\ssec{Relation to the lisse version} \label{ss:Ran vs Lisse}

In this subsection, we will relate the category $\CC_\Ran$ to its lisse counterpart.

\sssec{} 

Recall the category
$$\Rep(\sG)^{\otimes X\on{-lisse}}\simeq \coHom(\Rep(\sG),\qLisse(X)),$$
see \secref{sss:X to the Lisse}.  By the same token, we can consider the category
$$\CC^{\otimes X\on{-lisse}}:= \coHom(\CC,\qLisse(X))$$
for an arbitrary symmetric monoidal $\CC$.

\sssec{}

Let us note the difference between $\CC^{\otimes X\on{-lisse}}$ and $\CC_\Ran$.
In the former the terms of the colimit are 
$$\CC^{\otimes I}\otimes (\qLisse(X)^\vee)^{\otimes J}$$
and in the latter
$$\CC^{\otimes I}\otimes \Shv(X^J).$$

\sssec{} \label{sss:Ran vs Lisse new}

We claim that there is a naturally defined symmetric monoidal functor
\begin{equation} \label{e:Ran to lisse}
\CC_\Ran\to \CC^{\otimes X\on{-lisse}}.
\end{equation} 

In fact, there is a natural transformation between the corresponding right-lax symmetric 
monoidal functors $\on{TwArr}(\on{fSet})\to \DGCat$.

\medskip

Namely, for every $I\to J$, the corresponding functor is induced by the functor
$$\Shv(X^J) \to (\qLisse(X)^\vee)^{\otimes J},$$
\emph{dual} with respect to the Verdier self-duality on $\Shv(X^J)$ to the functor 
$$\qLisse(X)^{\otimes J}\to \Shv(X)^{\otimes J}\overset{\boxtimes}\to \Shv(X^J),$$
where 
$$\qLisse(X)\to \Shv(X)$$
is the embedding \eqref{e:!-emb}. 

\begin{rem}

Note that the functor \eqref{e:Ran to lisse} factors as 
$$\CC_\Ran\overset{\text{\eqref{e:Ran to coHom}}}\longrightarrow \coHom(\CC,\Shv(X))\to \coHom(\CC,\qLisse(X)),$$
where the second arrow comes from the symmetric monoidal embedding $\qLisse(X)\to \Shv(X)$ of \eqref{e:!-emb}.

\end{rem}

\sssec{}

We claim: 

\begin{prop} \label{p:Ran to lisse} 
Let $\bA$ be a dualizable DG category. Then the functor 
\begin{equation} \label{e:Ran to lisse bis} 
\on{Funct}_{\on{cont}}(\CC^{\otimes X\on{-lisse}},\bA)\to 
\on{Funct}_{\on{cont}}(\CC_\Ran,\bA),
\end{equation}
given by precomposition with \eqref{e:Ran to lisse} is fully faithful.
\end{prop} 

\begin{proof}

%

By \secref{sss:maps out of Ran new}, it suffices to show that for every $I\in \fSet$, the functor 
$$\bA\otimes \qLisse(X)^{\otimes I}\to \bA\otimes \Shv(X^I)$$
is fully faithful. However, this follows from the fact that $\qLisse(X)^{\otimes I}\to \Shv(X^I)$
is fully faithful, since $\bA$ is dualizable. 

\end{proof}

\begin{cor} \label{c:Ran to lisse prel}
For any DG category $\bA'$, a natural number $n$ and a dualizable $\bA$, the functor 
$$\on{Funct}_{\on{cont}}((\CC^{\otimes X\on{-lisse}})^{\otimes n}\otimes \bA',\bA)\to 
\on{Funct}_{\on{cont}}((\CC_\Ran)^{\otimes n}\otimes \bA',\bA)$$
is fully faithful.
\end{cor}

\begin{proof}

The assertion for $n=1$ follows from \propref{p:Ran to lisse} by taking $\on{Funct}_{\on{cont}}(\bA',-)$
into both sides of \eqref{e:Ran to lisse bis}. 

\medskip

The assertion for $n>1$ follows by iteration. 

\end{proof}


%

\sssec{}

From \corref{c:Ran to lisse prel} we obtain:

\begin{cor} \label{c:Ran to lisse} \hfill

\smallskip

\noindent{\em(a)} For any monoidal category $\bA$, precomposition 
with \eqref{e:Ran to lisse} defines a \emph{monomorphism} 
$$\Maps_{\DGCat^{\on{Mon}}}(\CC^{\otimes X\on{-lisse}},\bA) \to 
\Maps_{\DGCat^{\on{Mon}}}(\CC_\Ran,\bA),$$
provided $\bA$ is \emph{dualizable} as a DG category. 

\smallskip

\noindent{\em(b)}
For any symmetric monoidal category $\bA$, precomposition 
with \eqref{e:Ran to lisse} defines a \emph{monomorphism} 
$$\Maps_{\DGCat^{\on{SymMon}}}(\CC^{\otimes X\on{-lisse}},\bA) \to 
\Maps_{\DGCat^{\on{SymMon}}}(\CC_\Ran,\bA),$$
provided $\bA$ is \emph{dualizable} as a DG category. 

\smallskip

\noindent{\em(c)} For a pair of $\CC^{\otimes X\on{-lisse}}$-module categories $\bM_1,\bM_2$, the map
$$\on{Funct}_{\CC^{\otimes X\on{-lisse}}\mod}(\bM_1,\bM_2)\to 
\on{Funct}_{\CC_\Ran\mod}(\bM_1,\bM_2)$$
is an isomorphism, provided $\bM_2$ is dualizable as a DG category. 
\end{cor} 

\begin{rem} \label{r:Ran to lisse Betti}

An analog of the situation described in this subsection takes place for $\CC_\Ran^{\on{Betti}}$. In this
case, we have a tautological embedding
\begin{equation}  \label{e:loc const into all Ran }
\Shv^{\on{all}}_{\on{loc.const.}}(X) \hookrightarrow \Shv^{\on{all}}(X),
\end{equation}
which gives rise to a symmetric monoidal functor
$$\CC_\Ran^{\on{Betti}}\simeq \coHom(\CC,\Shv^{\on{all}}(X))\to
\coHom(\CC,\Shv^{\on{all}}_{\on{loc.const.}}(X))\simeq \CC^{\otimes X}.$$

The assertions and proofs of \propref{p:Ran to lisse} and \corref{c:Ran to lisse}
remain valid in this context as well. 

\end{rem}

\ssec{Rigidity} \label{ss:Ran rigid}

In this subsection we will show that $\CC_\Ran$ is \emph{rigid} as a monoidal category,
and as a result, is canonically self-dual. We will also describe the resulting datum
of self-duality explicitly. 

\sssec{}

We now reimpose the condition that $X$ be proper, for the duration of this section. 

\medskip

We will also assume that $\CC$ is compactly 
generated and \emph{rigid}. Given compact generation, the latter condition means
that compact generators of $\CC$ are dualizable (in the sense of the 
symmetric monoidal structure).

\sssec{} 

We claim that the above conditions imply that $\CC_\Ran$ is also
compactly generated and rigid. 

\medskip

First, since the transition functors 
\eqref{e:transition functors} preserve compactness, a set of compact generators
of $\CC_\Ran$ is provided by objects of the form 
\begin{equation} \label{e:comp gen Ran}
\on{ins}_{I\to J}(V\otimes \CF), \quad V\in (\CC^{\otimes I})^c,\,\, \CF\in\Shv(X^J)^c,
\end{equation} 

\sssec{}

We now show that $\CC_\Ran$ is rigid. To do so, it is enough to show that its
compact generators are dualizable (in the sense of the symmetric monoidal structure
on $\CC_\Ran$). We will exhibit the duality data for compact generators explicitly.

\medskip

Namely, for an object \eqref{e:comp gen Ran}, 
its monoidal dual is given by 
$$\on{ins}_{I\to J}(V^\vee\otimes \BD(\CF)),$$
where $\BD$ denotes Verdier duality on $\Shv(X^J)^c$.

\medskip

The unit and counit maps are defined as follows. 

\sssec{}

The counit is:

$$\on{ins}_{I\to J}(V\otimes \CF)\otimes \on{ins}_{I\to J}(V^\vee\otimes \BD(\CF))
\simeq 
\on{ins}_{I\sqcup I\to J\sqcup J}((V\boxtimes V^\vee)\otimes (\CF\boxtimes \BD(\CF)))\to $$
$$\to \on{ins}_{I\sqcup I\to J\sqcup J}((V\boxtimes V^\vee)\otimes (\Delta_{X^J})_*(\omega_{X^J}))
\simeq 
\on{ins}_{I\sqcup I\to J}((V\boxtimes V^\vee)\otimes \omega_{X^J})\simeq$$
$$\simeq \on{ins}_{I\to J}((V\otimes V^\vee)\otimes \omega_{X^J}) \to
\on{ins}_{I\to J}(\one_{\CC^{\otimes I}}\otimes \omega_{X^J})\simeq
\on{ins}_{\emptyset\to J}(\sfe \otimes \omega_{X^J})\simeq $$
$$\simeq \on{ins}_{\emptyset\to \emptyset}(\sfe \otimes \on{C}^\cdot(X^J,\omega_{X^J}))
\simeq \one_{\CC_\Ran}\otimes \on{C}^\cdot(X^J,\omega_{X^J})\to \one_{\CC_\Ran},$$
where the last arrow is the trace map, well-defined due to the fact that $X$ is proper.

\begin{rem}

In the above formula, we have used the notation 
$$V\boxtimes V^\vee\in \CC^{\otimes I}\otimes \CC^{\otimes I},$$
to be distinguished from
$$V\otimes V^\vee\in \CC^{\otimes I}.$$
I.e., the latter object is obtained from the former by applying the monoidal operation
$$\CC^{\otimes I}\otimes \CC^{\otimes I}\to \CC^{\otimes I}.$$

\end{rem}

\sssec{}

The unit is given by 
$$\one_{\CC_\Ran}\to \one_{\CC_\Ran}\otimes \on{C}^\cdot(X^J,\ul\sfe_{X^J})
 \simeq \on{ins}_{\emptyset \to \emptyset}(\sfe\otimes \on{C}^\cdot(X^J,\ul\sfe_{X^J}))\simeq \on{ins}_{\emptyset \to J}(\sfe\otimes \ul\sfe_{X^J})\simeq $$
$$\simeq \on{ins}_{I\to J}(\one_{\CC^{\otimes I}}\otimes \ul\sfe_{X^J})\to \on{ins}_{I\to J}((V\otimes V^\vee)\otimes \ul\sfe_{X^J})\simeq 
\on{ins}_{I\sqcup I\to J}((V\boxtimes V^\vee)\otimes \ul\sfe_{X^J})\simeq$$
$$\simeq \on{ins}_{I\sqcup I\to J\sqcup J}((V\boxtimes V^\vee)\otimes (\Delta_{X^J})_*(\ul\sfe_{X^J}))\simeq
\on{ins}_{I\sqcup I\to J\sqcup J}((V\boxtimes V^\vee)\otimes (\Delta_{X^J})_!(\ul\sfe_{X^J}))\to$$
$$\to \on{ins}_{I\sqcup I\to J\sqcup J}((V\boxtimes V^\vee)\otimes (\CF\boxtimes \BD(\CF)))
\simeq \on{ins}_{I\to J}(V\otimes \CF)\otimes \on{ins}_{I\to J}(V^\vee\otimes \BD(\CF)).$$

\sssec{} \label{sss:self-dual naive} 

Recall that if $\bA$ is a compactly generated rigid symmetric monoidal category, then it is canonically self-dual
as a DG category. Namely, the corresponding anti-equivalence
$$(\bA^c)^{\on{op}}\to \bA^c$$
is given by monoidal dualization. (For another description of this self-duality see 
\secref{sss:self-duality rigid} below.)  

\medskip

In the next subsection we will describe explicitly the resulting self-duality on $\CC_\Ran$. 

\begin{rem} \label{r:rigid dr non-rigid Betti}

The material in this subsection can be applied ``as-is" to $\CC_\Ran$ replaced by $\CC^\dr_\Ran$. 

\medskip

However, the situation is different for $\Shv^{\on{all}}(-)$ in that 
$\CC^{\on{Betti}}_\Ran$ is \emph{not} rigid (the unit object is no longer compact). Yet, it retains some
features, which will make the key construction work, see \secref{sss:proj Betti}.

\end{rem}

\ssec{Self-duality} 

\sssec{} \label{sss:limits and colimits new}

Let $I$ be an index category, and let 
$$i\mapsto \CC_i, \quad (i_1\to i_2) \rightsquigarrow \CC_{i_1}\overset{\phi_{i_1,i_2}}\longrightarrow \CC_{i_2}. $$
is a functor $I\to \DGCat$. Denote 
$$\cD:=\underset{i\in I}{\on{colim}}\, \CC_i.$$

For $i\in I$, let $\on{ins}_i$ denote the tautological functor $\CC_i\to \cD$. 

\sssec{}

Assume that for every 1-morphism $i_1\to i_2$ in $I$, the
transition functor $\phi_{i_1,i_2}:\CC_{i_1}\to \CC_{i_2}$ admits a continuous
right adjoint. 

\medskip

In this case we can form a functor 
$$I^{\on{op}}\to \DGCat, \quad i\mapsto \CC_i, \quad (i_1\to i_2) \rightsquigarrow \CC_{i_2}\overset{\phi^R_{i_1,i_2}}\longrightarrow \CC_{i_1}.$$

\medskip

According to \cite[Chapter 1, Proposion 2.5.7]{GR1}, the functors $\on{ins}_i$ also admit continuous
right adjoints. Furthermore, the resulting functor
\begin{equation} \label{e:colimit as limit gen}
\cD\to \underset{i\in I^{\on{op}}}{\on{lim}}\, \CC_i,
\end{equation} 
whose components are the right adjoints $(\on{ins}_i)^R$, is an equivalence.

\sssec{}

For future reference note that the equivalence \eqref{e:colimit as limit gen} implies that 
for $d\in \cD$, the canonical map
\begin{equation} \label{e:object as colimit}
\underset{i\in I}{\on{colim}}\, \on{ins}_i\circ \on{ins}_i^R(d)\to d
\end{equation} 
is an isomorphism.

\sssec{}

Assume now that each $\CC_i$ is dualizable. We can form a new functor
$$I\to \DGCat, \quad \quad i\mapsto \CC_i^\vee, \quad (i_1\to i_2) 
\rightsquigarrow \CC_{i_1}^\vee \overset{(\phi^R_{i_1,i_2})^\vee}\longrightarrow  \CC_{i_1}^\vee.$$

Note that if the $\CC_i$ are compactly generated, the functor $(\phi^R_{i_1,i_2})^\vee$, 
when restricted to compact objects, viewed as the functor $(\CC^c_{i_1})^{\on{op}}\to (\CC^c_{i_2})^{\on{op}}$,
is the opposite of $\phi_{i_1,i_2}:\CC^c_{i_1}\to \CC^c_{i_2}$, see \cite[Chapter 1, Proposition 7.3.5]{GR1}.

\medskip

According to \cite[Proposition 1.8.3]{DrGa2}, the category $\cD$ is also dualizable, and the functor 
$$\underset{i\in I}{\on{colim}}\, \CC^\vee_i\to \cD^\vee$$
comprised of the functors
$$(\on{ins}_i^R)^\vee:\CC^\vee_i\to \cD^\vee$$
is an equivalence. 

\sssec{}

To summarize, we obtain that there is a canonical duality between
$$\cD:=\underset{i\in I}{\on{colim}}\, \CC_i \text{ and } \cD':=\underset{i\in I}{\on{colim}}\, \CC^\vee_i,$$
under which, the functor 
$$\on{ins}'_i: \CC^\vee_i \to \cD'$$
identifies with the dual of 
$$\on{ins}_i^R:\cD\to \CC_i.$$

\medskip

In particular, if $\on{u}_\cD\in \cD\otimes \cD'$ denote the unit of the duality, formula \eqref{e:object as colimit} implies
that we have a canonical isomorphism
\begin{equation} \label{e:unit for colimit}
\on{u}_\cD \simeq \underset{i}{\on{colim}}\, (\on{ins}_i\otimes \on{ins}'_i)(\on{u}_{\CC_i}),
\end{equation}
where $\on{u}_{\CC_i}\in \CC_i\otimes \CC_i^\vee$ is the unit of the $(\CC_i,\CC_i^\vee)$ duality. 

\sssec{} \label{sss:colim self-dual}

Suppose now that for every $i\in I$ we are given a data of self-duality 
$$\CC_i^\vee \simeq \CC_i,$$
so that the functor
$$I^{\on{op}}\to \DGCat, \quad i\mapsto \CC_i^\vee, \quad (i_1\to i_2) \rightsquigarrow \CC_{i_2}^\vee \overset{\phi^\vee_{i_1,i_2}}\longrightarrow  \CC_{i_1}^\vee$$
is identifies with the functor 
$$I^{\on{op}}\to \DGCat, \quad i\mapsto \CC_i, \quad (i_1\to i_2) \rightsquigarrow \CC_{i_2}\overset{\phi^R_{i_1,i_2}}\longrightarrow \CC_{i_1}.$$

We obtain that in this case there is a canonical self-duality
$$\cD^\vee \simeq \cD,$$
with respect to which we have
$$\on{ins}_i^R\simeq \on{ins}_i^\vee.$$

\sssec{}

Applying this to $I:= \on{TwArr}(\fSet)$ and the functor \eqref{e:Ran diagram}, we obtain a self-duality
\begin{equation} \label{e:Ran self-dual}
(\CC_\Ran)^\vee \simeq \CC_\Ran.
\end{equation}

Indeed, for an individual object $(I\overset{\psi}\to J)\in \on{TwArr}(\fSet)$, the category 
$$\CC^{\otimes I}\otimes \Shv(X^J)$$
is canonically self-dual due to:

\begin{itemize}

\item The canonical self-duality on $\CC$ arising from the fact that $\CC$ is rigid
(i.e., it acts on compact objects as monoidal dualization);

\item Verdier self-duality on $\Shv(X^J)$.

\end{itemize}

For a 1-morphism \eqref{e:mor Tw arr}, the functor \eqref{e:transition functors} identifies
with the dual of its right adjoint because:

\begin{itemize}

\item The functor \eqref{e:mult phi} is monoidal and hence commutes with monoidal dualization on
dualizable (hence, compact) objects;

\item The functor $(\Delta_{\phi_J})_*$ commutes with Verdier duality, due to the
assumption that $X$ is proper.

\end{itemize} 

\sssec{}

According to \secref{sss:colim self-dual}, with respect to the identification
$$(\CC_\Ran)^\vee \simeq \CC_\Ran$$
of \secref{e:Ran self-dual}, the dual of the functor
$$\on{ins}_{I\to J}:\CC^{\otimes I}\otimes \Shv(X^J)\to \CC_\Ran$$
identifies with
\begin{equation} \label{e:ins dual}
(\CC_\Ran)^\vee \simeq  \CC_\Ran \overset{(\on{ins}_{I\to J})^R}\longrightarrow 
\CC^{\otimes I}\otimes \Shv(X^J) \simeq (\CC^{\otimes I}\otimes \Shv(X^J))^\vee.
\end{equation}

\sssec{}

Unwinding the definitions, one can see that the self-duality \eqref{e:Ran self-dual}
coincides with the one in \secref{sss:self-dual naive}.

\ssec{The progenitor of the projector}  \label{ss:progenitor}

In this subsection we introduce an object $$\sR_{\CC,\Ran}\in \CC_\Ran\otimes \CC_\Ran,$$
which will ultimately give rise to Beilinson's spectral projector. 

\sssec{} \label{sss:self-duality rigid}

Recall again that if $\bA$ is a rigid symmetric monoidal category, it is canonically self-dual, see
\cite[Chapter 1, Sect. 9.2]{GR1}.

\medskip

Namely, the counit is given by
$$\bA\otimes \bA\overset{\on{mult}_\bA}\to \bA \overset{\CHom_\bA(\one_\bA,-)}\longrightarrow \Vect_\sfe,$$
and the unit is given by
$$\Vect_\sfe \overset{\one_\bA}\to \bA \overset{\on{comult}_\bA}\to \bA\otimes \bA,$$
where the functor 
\begin{equation} \label{e:comult rigid}
\on{comult}_\bA:\bA \to \bA\otimes \bA
\end{equation}
is the right adjoint to $\on{mult}_\bA$. 

%

\medskip

Let $\sR_\bA\in \bA\otimes \bA$ denote the unit of the above self-duality on $\bA$. I.e.,
$$\sR_\bA:=\on{comult}_\bA(\one_\bA).$$

\medskip

Being the right adjoint of a symmetric monoidal functor, the functor $\on{comult}_\bA$
carries a natural right-lax symmetric monoidal structure. Hence, $\sR_\bA$ is naturally
a commutative algebra object in $\bA\otimes \bA$. 

\sssec{} \label{sss:R Ran}

Let us apply the above discussion to $\CC_\Ran$. Let 
$$\sR_{\CC,\Ran}\in \CC_\Ran\otimes \CC_\Ran$$
denote the unit of the self-duality. 

\medskip

This object will play a key role in the sequel. By the above, $\sR_{\CC,\Ran}$ carries a natural structure
of commutative algebra in $\CC_\Ran\otimes \CC_\Ran$.

\medskip

Our next goal is to describe $\sR_{\CC,\Ran}$ explicitly as a colimit. 

\sssec{} 

Recall that if $\bA$ is a rigid symmetric monoidal category, then the functor 
$\on{comult}_\bA$ of \eqref{e:comult rigid} is \emph{strictly} compatible with the $\bA$-bimodule 
structure (being a right adjoint of a map of bimodules, the functor $\on{comult}_\bA$ is a 
priori right-lax compatible with the bimodule structure). 

\medskip

This means that the object $\sR_\bA\in \bA\otimes \bA$ naturally lifts to an object of
$$\on{HC}^\bullet(\bA,\bA\otimes \bA):=\on{Funct}_{(\bA\otimes \bA)\mod}(\bA,\bA\otimes \bA)$$
(here $\on{HC}^\bullet(\bA,-)$ stands for the ``Hochschild cohomology"  category with coefficients
in a given $\bA$-bimodule category). 

\medskip

In other words, we have a canonical system of isomorphisms
$$(\ba\otimes \one_\bA)\otimes \sR_\bA \simeq \sR_\bA \otimes (\one_\bA\otimes \ba),\quad \ba\in \bA,$$
compatible with the monoidal structure on $\bA$. 

\sssec{}

Applying this for $\bA:=\CC_\Ran$, we obtain a system of isomorphisms 
\begin{equation} \label{e:Hecke univ}
(\CV \otimes \one_{\CC_\Ran})\star \sR_{\CC,\Ran} \simeq 
\sR_{\CC,\Ran} \star (\one_{\CC_\Ran} \otimes \CV), \quad \CV\in \CC_\Ran.
\end{equation}

As we shall see, the system of isomorphisms \eqref{e:Hecke univ} is the source 
of Hecke eigen-property of various objects that we will establish in the sequel. 

\begin{rem}
The system of isomorphisms \eqref{e:Hecke univ} is equally valid when we work with $\CC^\dr_\Ran$.
\end{rem}

\ssec{The progenitor as a colimit}

\sssec{}

Applying \eqref{e:colimit as limit gen} to $\CC_\Ran$, we obtain that it can also be written as a \emph{limit}
\begin{equation} \label{e:Ran as limit}
\CC_\Ran\simeq \underset{(I\to J)\in (\on{TwArr}(\on{fSet}))^{\on{op}}}{\on{lim}}\, \CC^{\otimes I}\otimes \Shv(X^J),
\end{equation}
where the transition functor corresponding to \eqref{e:mor Tw arr} is the tensor product of 
$$(\on{mult}_\CC^{\phi_I})^R: \CC^{\otimes I_2}\to \CC^{\otimes I_1}$$
and 
$$(\Delta_{\phi_J})^!:\Shv(X^{J_2})\to \Shv(X^{J_1}).$$

\sssec{}

Let us apply \eqref{e:unit for colimit} to the object 
$$\sR_{\CC,\Ran}\in \CC_\Ran\otimes \CC_\Ran.$$

We claim:
\begin{equation}  \label{e:progenitor}
\sR_{\CC,\Ran}\simeq
\underset{(I\to J)\in \on{TwArr}(\fSet)}{\on{colim}}\,
(\on{ins}_{I\to J}\otimes \on{ins}_{I\to J})(\sR^{\boxtimes I}_\CC\otimes \on{u}_{\Shv(X^J)}),
\end{equation} 
where:

\begin{itemize}

\item $\sR_\CC\in \CC\otimes \CC$ denotes the unit of the self-duality on $\CC$,
arising from the fact that $\CC$ is a rigid symmetric monoidal category;

\smallskip

\item $\sR^{\boxtimes I}_\CC$ denotes the $I$-tensor power of $\sR_\CC$, viewed
as an object of $\CC^{\otimes I}\otimes \CC^{\otimes I}$;

\smallskip

\item For a scheme $Y$, we denote by 
$\on{u}_{\Shv(Y)}\in \Shv(Y)\otimes \Shv(Y)$ is the unit of the Verdier 
self-duality on $\Shv(Y)$. 

\end{itemize}

\medskip

Indeed, this follows from \eqref{e:unit for colimit} using the identification
$$(\on{ins}_{I\to J})^R\overset{\text{\eqref{e:ins dual}}}\simeq (\on{ins}_{I\to J})^\vee$$
and the fact that $\sR^{\boxtimes I}_\CC$ is the unit of the self-duality on $\CC^{\otimes I}$,
induced by the self-duality of $\CC$. 

\begin{rem} \label{r:quasi-diag}
For future use, let us observe that the object $\on{u}_{\Shv(Y)}\in \Shv(Y)\otimes \Shv(Y)$
introduced above can also be interpreted as the value on $(\Delta_Y)_*(\omega_Y)$ of the
right adjoint $\boxtimes^R$ to the external tensor product functor 
\begin{equation} \label{e:ext Y}
\Shv(Y)\otimes \Shv(Y)\overset{\boxtimes}\to \Shv(Y\times Y)
\end{equation}
on $(\Delta_Y)_*(\omega_Y)$, see \secref{sss:quasi-diag}. 

\medskip

By a slight abuse of notation, we will sometimes denote by the same symbol $\on{u}_{\Shv(Y)}$
the image of this object along the (fully faithful) functor \eqref{e:ext Y}. This is done in order to
avoid the somewhat awkward notation $\boxtimes(\on{u}_{\Shv(Y)})$.

\medskip

The counit of the adjunction 
$$\on{u}_{\Shv(Y)}\to (\Delta_Y)_*(\omega_Y)$$
has the following basic property: for $\CF\in \Shv(Y)$, the induced map 
\begin{equation} \label{e:proj formula u}
(p_2)_*(\on{u}_{\Shv(Y)}\sotimes (p_1)^!(\CF))\to (p_2)_*((\Delta_Y)_*(\omega_Y)\sotimes (p_1)^!(\CF))\simeq \CF
\end{equation} 
is an isomorphism. 

\end{rem} 

\begin{rem}
Formula \eqref{e:progenitor} holds also for the unit of the self-duality of $\CC^\dr_\Ran$, 
with the difference that now instead of the object $\on{u}_{\Shv(X^J)}$ we use
$$(\Delta_{X^J})_*(\omega_{X^J})\in \Dmod(X^J\times X^J)
\simeq \Dmod(X^J)\otimes \Dmod(X^J).$$
\end{rem} 

\sssec{} \label{sss:counit R}

We will now describe explicitly particular values of the unit and counit of the adjunction
$$\on{mult}_{\CC_\Ran}:\CC_\Ran\otimes \CC_\Ran\rightleftarrows \CC_\Ran:\on{comult}_{\CC_\Ran},$$
in terms of formula \eqref{e:progenitor}.

\medskip

The unit of the adjunction, when evaluated on $\one_{\CC_\Ran}\otimes \one_{\CC_\Ran}\in \CC_\Ran\otimes \CC_\Ran$, is a map
$$\one_{\CC_\Ran}\otimes \one_{\CC_\Ran} \to \sR_{\CC,\Ran}.$$

It corresponds to the term $(I\to J)=(\emptyset\to \emptyset)$ in the colimit \eqref{e:progenitor}. 

\medskip

The counit of the adjunction, when evaluated on $\one_{\CC_\Ran}$, is a map 
\begin{equation} \label{e:counit via prog}
\on{mult}_{\CC_\Ran}(\sR_{\CC,\Ran})\to \one_{\CC_\Ran}.
\end{equation} 

Here is an explicit description of this map in terms of \eqref{e:progenitor}.

\sssec{} \label{sss:counit via prog}

In order to describe \eqref{e:counit via prog}, we need to specify a compatible system of maps
\begin{equation} \label{e:counit progenitor}
\on{ins}_{I\sqcup I\to J\sqcup J}(\sR_\CC^{\boxtimes I}\otimes  \on{u}_{\Shv(X^J)})\to \one_{\CC_\Ran}, \quad
(I\to J)\in \on{TwArr}(\fSet),
\end{equation} 

The map in \eqref{e:counit progenitor} is the following composition:

\begin{itemize}

\item Using the counit of the adjunction $\on{u}_{\Shv(X^J)}\to (\Delta_{X^J})_*(\omega_{X^J})$ we map
the left-hand side in \eqref{e:counit progenitor} to
\begin{equation} \label{e:counit progenitor 1}
\on{ins}_{I\sqcup I\to J\sqcup J}(\sR_\CC^{\boxtimes I}\otimes  (\Delta_{X^J})_*(\omega_{X^J})).
\end{equation}

\item The expression in \eqref{e:counit progenitor 1} is isomorphic to
\begin{equation} \label{e:counit progenitor 2}
\on{ins}_{I\sqcup I\to J}(\sR_\CC^{\boxtimes I}\otimes  \omega_{X^J}).
\end{equation}

\item The expression in \eqref{e:counit progenitor 2} is isomorphic to
\begin{equation} \label{e:counit progenitor 3}
\on{ins}_{I\to J}((\on{mult}_\CC(\sR_\CC))^{\boxtimes I}\otimes  \omega_{X^J}),
\end{equation}
where $\on{mult}_\CC(\sR_\CC)\in \CC$ and $(\on{mult}_\CC(\sR_\CC))^{\boxtimes I}\in \CC^{\otimes I}$.

\smallskip

\item Using the counit of the adjunction $\on{mult}_\CC(\sR_\CC)\to \one_\CC$, we map \eqref{e:counit progenitor 3} to
\begin{equation} \label{e:counit progenitor 4}
\on{ins}_{I\to J}((\one_\CC)^{\boxtimes I}\otimes  \omega_{X^J}).
\end{equation}

\item The expression in \eqref{e:counit progenitor 4} is isomorphic to
\begin{equation} \label{e:counit progenitor 5}
\on{ins}_{\emptyset\to J}(\sfe\otimes  \omega_{X^J}).
\end{equation}

\item The expression in \eqref{e:counit progenitor 5} is isomorphic to
\begin{equation} \label{e:counit progenitor 6}
\on{ins}_{\emptyset\to \emptyset}(\sfe\otimes  \on{C}^\cdot(X^J,\omega_{X^J})).
\end{equation}

\item Using the trace map $\on{C}^\cdot(X^J,\omega_{X^J})\to \sfe$, we map 
\eqref{e:counit progenitor 6} to 
$$\on{ins}_{\emptyset\to \emptyset}(\sfe\otimes \sfe)=\one_{\CC_\Ran}.$$

\end{itemize}

\ssec{Explicit construction of the Hecke isomorphisms} \label{ss:Hecke abs}

We have deduced the system of isomorphisms \eqref{e:Hecke univ} from the rigidity property
of $\CC_\Ran$. However, one can prove it by a direct computation if we take formula \eqref{e:progenitor}
as the \emph{definition} of $\sR_{\CC,\Ran}$.

\sssec{} \label{sss:Hecke abs 1}

Let $\CV\in \CC_\Ran$ be of the form
$$\on{ins}_{I_0\to J_0}(V\otimes \CF), \quad V\in \CC^{\otimes I_0},\, \CF\in \Shv(X^{J_0}).$$

Let us construct the corresponding isomorphism
\begin{equation} \label{e:Hecke univ again}
(\CV \otimes \one_{\CC_\Ran})\star \sR_{\CC,\Ran} \simeq 
\sR_{\CC,\Ran} \star (\one_{\CC_\Ran} \otimes \CV), \quad \CV\in \CC_\Ran.
\end{equation}

\medskip

Namely, we claim that each side in \eqref{e:Hecke univ again} can be identified with the corresponding
side in 
\begin{multline} \label{e:progen tensored}
\underset{(I\to J)\in \on{TwArr}(\fSet)}{\on{colim}}\,
(\on{ins}_{I_0\sqcup I\to J_0\sqcup J}\otimes \on{ins}_{I_0\sqcup I\to J_0\sqcup J})
\left((V\otimes \sR^{\boxtimes I_0\sqcup I}_\CC)\otimes (\CF\sotimes  \on{u}_{\Shv(X^{J_0\sqcup J})})\right)\simeq \\
\simeq 
\underset{(I\to J)\in \on{TwArr}(\fSet)}{\on{colim}}\,
(\on{ins}_{I_0\sqcup I\to J_0\sqcup J}\otimes \on{ins}_{I_0\sqcup I\to J_0\sqcup J})
\left((\sR^{\boxtimes I_0\sqcup I}_\CC\otimes V)\otimes (\on{u}_{\Shv(X^{J_0\sqcup J})}\sotimes \CF)\right),
\end{multline} 
where:

\begin{itemize}

\item $V\otimes \sR^{\boxtimes I_0\sqcup I}_\CC$ and $\sR^{\boxtimes I_0\sqcup I}_\CC\otimes V$ are 
(isomorphic) objects
of $\CC^{\otimes I_0\sqcup I}\otimes \CC^{\otimes I_0\sqcup I}$ obtained by tensoring the object 
$\sR^{\boxtimes I_0\sqcup I}_\CC$ by
$$V\otimes \sfe\in \CC^{\otimes I_0}\otimes \CC^{\otimes I_0}\to \CC^{\otimes I_0 \sqcup I}\otimes \CC^{\otimes I_0 \sqcup I}
\text{ and }
\sfe\otimes V \in \CC^{\otimes I_0}\otimes \CC^{\otimes I_0} \to \CC^{\otimes I_0 \sqcup I}\otimes \CC^{\otimes I_0 \sqcup I},$$
respectively.

\item $\CF\sotimes  \on{u}_{\Shv(X^{J_0\sqcup J})}$ and $\on{u}_{\Shv(X^{J_0\sqcup J})}\sotimes \CF$ are (isomorphic) objects
of $\Shv(X^{J_0\sqcup J})\otimes \Shv(X^{J_0\sqcup J})$ obtained by tensoring the object $\on{u}_{\Shv(X^{J_0\sqcup J})}$ by the
!-pullback of $\CF$ along 
$$X^{J_0\sqcup J}\to X^{J_0}$$
on the left and right factor, respectively. 

\end{itemize}

\medskip

Let us show how to identify the left-hand side of \eqref{e:Hecke univ again} with the left-hand side of \eqref{e:progen tensored}
(the right-hand sides are handled by symmetry). 

\sssec{} \label{sss:Hecke abs 1.5}

Let us construct a map from the left-hand side of \eqref{e:Hecke univ again} to the left-hand side of \eqref{e:progen tensored}. 
Fix $(I\to J)\in \on{TwArr}(\fSet)$. 

\medskip

\noindent{{\bf Step 0.}}
We start with \begin{equation} \label{e:ins LLLHS}
\on{ins}_{I_0\to J_0}(V\otimes \CF) \star 
\left((\on{ins}_{I\to J}\otimes \on{ins}_{I\to J})(\sR^{\boxtimes I}_\CC\otimes \on{u}_{\Shv(X^J)})\right),
\end{equation} 
which is a term corresponding to $(I\to J)$ in the colimit expression in the left-hand side of \eqref{e:Hecke univ again} .

\medskip

\noindent{{\bf Step 1.}}
The object \eqref{e:ins LLLHS} is canonically isomorphic to the object 
\begin{equation} \label{e:ins LHS}
(\on{ins}_{I_0\sqcup I\to J_0\sqcup J}\otimes \on{ins}_{I\to J_0\sqcup J})
\left((V\boxtimes \sR_\CC^{\boxtimes I})\otimes (\CF\sotimes \on{u}_{\Shv(X^{J_0\sqcup J})})\right),
\end{equation} 
where we regard $V\boxtimes\sR_\CC^{\boxtimes I}$ as an object of 
$\CC^{\otimes I_0\sqcup I}\otimes \CC^{\otimes I}$. 

\medskip

Indeed, this isomorphism is induced by the 1-morphism in $\on{TwArr}(\fSet)\times \on{TwArr}(\fSet)$, which is identity along the first factor, and 
$$
\CD
I @>>> J_0 \sqcup J \\
@V{\on{id}}VV  @AAA \\
I @>>> J
\endCD
$$
along the second factor. 

\medskip

\noindent{{\bf Step 2.}} Next, we consider the object 
\begin{equation} \label{e:ins LLHS}
(\on{ins}_{I_0\sqcup I_0\sqcup I\to J_0\sqcup J}\otimes \on{ins}_{I_0\sqcup I\to J_0\sqcup J})
\left(((V\otimes \sfe \otimes \sfe) \boxtimes\sR_\CC^{\boxtimes I})\otimes (\CF\sotimes \on{u}_{\Shv(X^{J_0\sqcup J})})\right),
\end{equation} 
where we view $V\otimes \sfe \otimes \sfe$ as an object of $\CC^{\otimes I_0\sqcup I_0} \otimes \CC^{\otimes I_0}$, and 
$(V\otimes \sfe \otimes \sfe) \boxtimes\sR_\CC^{\boxtimes I}$ as an object of
$\CC^{\otimes I_0\sqcup I_0\sqcup I}\otimes \CC^{\otimes I_0\sqcup I}$. 

\medskip

We have a canonical isomorphism from \eqref{e:ins LHS} to \eqref{e:ins LLHS},
induced by the inclusions
$$I_0\sqcup I \hookrightarrow I_0\sqcup I_0\sqcup I \text{ and } I \hookrightarrow I_0\sqcup I.$$

\medskip

\noindent{{\bf Step 3.}} Consider the object 
\begin{equation} \label{e:ins RHS}
(\on{ins}_{I_0\sqcup I_0\sqcup I\to J_0\sqcup J}\otimes \on{ins}_{I_0\sqcup I\to J_0\sqcup J})
\left((V\boxtimes\sR_\CC^{\boxtimes I_0\sqcup I})\otimes (\CF\sotimes \on{u}_{\Shv(X^{J_0\sqcup J})})\right),
\end{equation} 
where we regard $V\boxtimes\sR_\CC^{\boxtimes I_0\sqcup I}$ as an object of 
$\CC^{\otimes I_0\sqcup I_0\sqcup I}\otimes \CC^{\otimes I_0\sqcup I}$. 

\medskip

We have a canonically defined map from \eqref{e:ins LLHS} to \eqref{e:ins RHS}, given by
$$\sfe \otimes \sfe\to\sR_\CC^{\boxtimes I_0}.$$

\medskip

\noindent{{\bf Step 4.}} The object \eqref{e:ins RHS} admits a canonical isomorphism to
\begin{equation} \label{e:ins RHS bis}
(\on{ins}_{I_0\sqcup I\to J_0\sqcup J}\otimes \on{ins}_{I_0\sqcup I\to J_0\sqcup J})
\left((V\otimes \sR^{\boxtimes I_0\sqcup I}_\CC)\otimes (\CF\sotimes  \on{u}_{\Shv(X^{J_0\sqcup J})})\right).
\end{equation}

This isomorphism is induced by the 1-morphism in $\on{TwArr}(\fSet)\times \on{TwArr}(\fSet)$, which is identity 
along the second factor and 
$$
\CD
I_0\sqcup I_0\sqcup I @>>>  J_0\sqcup J \\
@VVV @AA{\on{id}}A \\
I_0\sqcup I @>>>  J_0\sqcup J
\endCD
$$
along the first factor. 

\medskip

\noindent{{\bf Final step.}}  Finally, the object \eqref{e:ins RHS bis} is the term 
corresponding to $(I\to J)$ in the colimit expression in the left-hand side of \eqref{e:progen tensored}.

\sssec{} \label{sss:Hecke abs 2}

Let us now construct a map from the left-hand side of \eqref{e:progen tensored} to the left-hand side of \eqref{e:Hecke univ again}. 

\medskip

\noindent{{\bf Step 0.}} 
By Step 4 in \secref{sss:Hecke abs 1.5}, the term corresponding to $(I\to J)$ in the colimit expression in the left-hand side of 
\eqref{e:progen tensored} is isomorphic to \eqref{e:ins RHS}.

\medskip

\noindent{{\bf Step 1.}} We have a canonical isomorphism between \eqref{e:ins RHS} and 
\begin{multline} \label{e:ins RHS another}
(\on{ins}_{I_0\sqcup I_0\sqcup I\to J_0\sqcup J_0\sqcup J}\otimes \on{ins}_{I_0\sqcup I\to J_0\sqcup J}) \\ 
\left((V\boxtimes\sR_\CC^{\boxtimes I_0\sqcup I})\otimes ((\Delta_{X^{J_0}}\times \on{id}_{X^J})_*\otimes \on{Id})
(\CF\sotimes \on{u}_{\Shv(X^{J_0\sqcup J})})\right),
\end{multline}
where $(\Delta_{X^{J_0}}\times \on{id}_{X^J})_*\otimes \on{Id}$ is the functor
$$\Shv(X^{J_0\sqcup J})\otimes \Shv(X^{J_0\sqcup J})\to \Shv(X^{J_0\sqcup J_0\sqcup J})\otimes \Shv(X^{J_0\sqcup J}).$$

This isomorphism is defined using the 1-morphism in $\on{TwArr}(\fSet)\times \on{TwArr}(\fSet)$, which is is identity 
along the second factor and 
$$
\CD
I_0\sqcup I_0\sqcup I @>>> J_0\sqcup J \\
@V{\on{id}}VV @AAA \\
I_0\sqcup I_0\sqcup I @>>> J_0\sqcup J_0\sqcup J
\endCD
$$ 
along the first factor.  

\medskip

\noindent{{\bf Step 2.}} We have a canonically defined map \eqref{e:ins RHS another} to 
\begin{equation} \label{e:ins RRHS}
(\on{ins}_{I_0\sqcup I_0\sqcup I\to J_0\sqcup J_0\sqcup J}\otimes \on{ins}_{I_0\sqcup I\to J_0\sqcup J})
\left((V\boxtimes\sR_\CC^{\boxtimes I_0\sqcup I})\otimes (\CF\boxtimes \on{u}_{\Shv(X^{J_0\sqcup J})})\right),
\end{equation} 
where we regard $\CF\boxtimes \on{u}_{\Shv(X^{J_0\sqcup J})}$ as an object of 
$\Shv(X^{J_0\sqcup J_0\sqcup J})\otimes \Shv(X^{J_0\sqcup J})$. 

\medskip

This map is induced by the morphism 
$$((\Delta_{X^{J_0}}\times \on{id}_{X^J})_*\otimes \on{Id})(\CF\sotimes \on{u}_{\Shv(X^{J_0\sqcup J})})\to
(\CF\boxtimes \on{u}_{\Shv(X^{J_0\sqcup J})}),$$
arising by adjunction from the isomorphism
$$\CF\sotimes \on{u}_{\Shv(X^{J_0\sqcup J})} \simeq ((\Delta_{X^{J_0}}\times \on{id}_{X^J})^!\otimes \on{Id})
(\CF\boxtimes \on{u}_{\Shv(X^{J_0\sqcup J})}).$$

\medskip

\noindent{{\bf Final step.}} The object \eqref{e:ins RRHS} is the term 
corresponding to $(I_0\sqcup I\to J_0\sqcup J)$ in the colimit expression in the right-hand side of \eqref{e:Hecke univ again}.

\sssec{}

One shows that by a routine diagram chase that the two maps between the left-hand side of 
\eqref{e:Hecke univ again} and the left-hand side of \eqref{e:progen tensored} are mutually inverse.  

%
%
%
%

\section{The spectral projector and localization} \label{s:Loc on LocSys}

In this section we will relate the object $\sR_{\Rep(\sG),\Ran}$ to the localization functor
$$\Loc:\Rep(\sG)_{\Ran}\to \QCoh(\LocSys^{\on{restr}}_\sG(X)).$$

\ssec{The progenitor for coHom} \label{ss:proj coHom}

\sssec{} \label{sss:progen coHom}

Let $\CC$ be a rigid symmetric monoidal
category, and let $\bH$ be another symmetric monoidal category, assumed dualizable as a
plain DG category. 

\medskip

Consider the category $\coHom(\CC,\bH)$ (see \secref{sss:coHom first}). In \secref{sss:coHom proof}
we will prove:

\begin{prop} \label{p:coHom almost semi-rigid}
The functor 
$$\on{mult}_{\coHom(\CC,\bH)}:\coHom(\CC,\bH)\otimes \coHom(\CC,\bH) \to \coHom(\CC,\bH)$$
admits a continuous right adjoint, to be denoted $\on{comult}_{\coHom(\CC,\bH)}$. Moreover, the 
structure on $\on{comult}_{\coHom(\CC,\bH)}$ of right-lax compatibility with the 
$\coHom(\CC,\bH)$-bimodule structure structure is strict. 
\end{prop}

\begin{rem}
Note that if we knew that $\coHom(\CC,\bH)$ was dualizable as a plan DG category, \propref{p:coHom almost semi-rigid}
would mean that the category $\coHom(\CC,\bH)$ is \emph{semi-rigid} (see \secref{sss:semi-rigid defn} for what this means). 
\end{rem}

\sssec{}

Denote 
\begin{equation} \label{e:progenitor coHom}
\sR_{\coHom(\CC,\bH)}:=\on{comult}_{\coHom(\CC,\bH)}(\one_{\coHom(\CC,\bH)})\in \coHom(\CC,\bH)\otimes \coHom(\CC,\bH).
\end{equation}

Being the right adjoint of a symmetric monoidal functor, the functor $\on{comult}_{\coHom(\CC,\bH)}$ acquires a natural
right-lax symmetric monoidal structure. Hence, the object $\sR_{\coHom(\CC,\bH)}$ carries a natural structure of 
commutative algebra in $\coHom(\CC,\bH) \otimes \coHom(\CC,\bH)$.

\sssec{}

Recall that the category $\coHom(\CC,\bH)$ identifies with
$$\underset{(I\to J)\in \on{TwArr}(\fSet)}{\on{colim}}\, \CC^{\otimes I} \otimes (\bH^\vee)^{\otimes J},$$
see \lemref{l:coEnd abs}. Denote by
$$\on{ins}_{I\to J}:\CC^{\otimes I} \otimes (\bH^\vee)^{\otimes J}\to \coHom(\CC,\bH)$$
the resulting tautological functors. 

\medskip

We will now describe $\sR_{\coHom(\CC,\bH)}$ in terms of the above colimit presentation.

\sssec{} \label{sss:again assumptions on H new}

Assume that the functor $\one_\bH:\Vect_\sfe\to \bH$ admits a left adjoint, to be denoted $\coinv_\bH$. 
Let us view $\coinv_\bH$ as an object of $\bH^\vee$. Let $\sR_{\bH^\vee}\in \bH^\vee\otimes \bH^\vee$ denote the image
of $\coinv_\bH$ under the dual of the monoidal operation $\on{mult}_\bH:\bH\otimes \bH\to \bH$.

\medskip

In \secref{sss:R coHom proof} we will prove:

\begin{prop} \label{p:R coHom almost semi-rigid}
There exists a canonical isomorphism: 
\begin{multline} \label{e:progenitor coHom formula}
\sR_{\coHom(\CC,\bH)} \simeq 
\underset{(I\to J)\in \on{TwArr}(\fSet)}{\on{colim}}\, 
(\on{ins}_{I\to J}\otimes \on{ins}_{I\to J})\left((\sR_\CC)^{\boxtimes I} \otimes (\sR_{\bH^\vee})^{\boxtimes J}\right)\in \\
\in \coHom(\CC,\bH)\otimes \coHom(\CC,\bH),
\end{multline} 
where we view $(\sR_\CC)^{\boxtimes I} \otimes (\sR_{\bH^\vee})^{\boxtimes J}$ as an object of
$$(\CC\otimes \CC)^{\otimes I}\otimes (\bH^\vee\otimes \bH^\vee)^{\otimes J}\simeq 
(\CC^{\otimes I} \otimes (\bH^\vee)^{\otimes J})\otimes (\CC^{\otimes I} \otimes (\bH^\vee)^{\otimes J}).$$
\end{prop}

\ssec{Abstract version of factorization homology}

In this subsection we introduce an important tool: an abstraction version of the
procedure known as \emph{factorization homology}.

\medskip

It will be used for the proof of Propositions \ref{p:coHom almost semi-rigid} and \ref{p:R coHom almost semi-rigid} as well as other results. 

\sssec{}  \label{sss:fact homology}

Consider the following paradigm. Let $\bA$ and $\bA'$ be a pair of symmetric monoidal categories, and let $\Phi:\bA'\to \bA$ be a symmetric monoidal
functor that admits a left adjoint, denoted $\Phi^L$, as a functor of plain DG categories. 

\medskip

Then the induced functor
$$\Phi:\on{ComAlg}(\bA')\to \on{ComAlg}(\bA)$$ admits a left adjoint, to be denoted $\Phi^{L,\on{ComAlg}}$,
which is described as follows.

\medskip

Define the functor
$$\wt\Phi^{L,\on{ComAlg}}: \on{ComAlg}(\bA)\to \on{ComAlg}(\bA')$$
as follows:

\medskip

Its value on $R\in \on{ComAlg}(\bA)$, viewed as a plain object of $\bA'$, 
is given by the colimit over $\on{TwArr}(\on{fSet})$ of the functor that sends
\begin{equation} \label{e:tw arrows object}
(I \overset{\psi}\to J)\in \on{TwArr}(\on{fSet})
\end{equation} 
to 
$$\on{mult}^J_{\bA'}\circ (\Phi^L)^{\otimes J}\circ \on{mult}_\bA^\psi(R^{\otimes I}),$$
where $\on{mult}^J_{\bA'}$ is the $J$-fold tensor product functor
$$(\bA')^{\otimes J}\to \bA'.$$

\medskip

The structure on $\wt\Phi^{L,\on{ComAlg}}$ of commutative algebra is induced by the operation of disjoint union on $\on{fSet}$.

\begin{lem} \label{l:fact homology}
The functor $\wt\Phi^{L,\on{ComAlg}}$ is canonically isomorphic to the left adjoint, denoted $\Phi^{L,\on{ComAlg}}$,
of
$$\Phi:\on{ComAlg}(\bA')\to \on{ComAlg}(\bA).$$ 
\end{lem}

The proof will be given in \secref{ss:proof fact homology}.

%
%
%

\begin{rem}
The above description of the left adjoint to $\Phi:\on{ComAlg}(\bA')\to \on{ComAlg}(\bA)$ is most familiar in the context
of \emph{factorization homology}. Namely, take 
$$\bA=(\Shv(X),\sotimes), \,\, \bA'=\Vect_\sfe,\,\, \Phi(\sfe)=\omega_X.$$
Then the functor $\wt\Phi^{L,\on{ComAlg}}$ evaluated on $R\in \on{ComAlg}^!(\Shv(X))$ sends to
$$\underset{(I\overset{\psi}\to J)\in \on{TwArr}(\on{fSet})}{\on{colim}}\,
\on{C}_c^\cdot\left(X^J,\underset{j\in J}\boxtimes\, R^{\otimes \psi^{-1}(j)}\right),$$
which is the formula for the facorization homology of $R$ along $X$. 
\end{rem} 

\sssec{} \label{sss:left adj fact homology}

We now make the following observation: let 
$$
\CD
\bA @>{\Psi}>> \bA_1 \\
@A{\Phi}AA @AA{\Phi_1}A \\
\bA' @>{\Psi'}>> \bA'_1
\endCD
$$
be a commutative diagram of symmetric monoidal categories. Note that we have natural transformations
\begin{equation} \label{e:left adj before fact homology}
\Phi_1^{L}\circ \Psi\to 
\Psi'\circ \Phi^{L}, \quad 
\end{equation}
and 
\begin{equation} \label{e:left adj fact homology}
\Phi_1^{L,\on{ComAlg}}\circ \Psi \to 
\Psi'\circ \Phi^{L,\on{ComAlg}}.
\end{equation}

\begin{lem} \label{l:left adj fact homology}
If the natural transformation \eqref{e:left adj before fact homology} is an isomorphism, then so is 
\eqref{e:left adj fact homology}. 
\end{lem} 

\begin{proof}

The category $\on{ComAlg}(\bA)$ is generated under sifted colimits by free objects, i.e., 
objects of the form $\Sym(\ba)$, for $\ba\in \bA$. Since all functors in \eqref{e:left adj fact homology}
preserve colimits, it suffices to show that the map \eqref{e:left adj fact homology} is an isomorphism
when evaluated on objects of the above form.

\medskip

We have, tautologically:
$$\Phi^{L,\on{ComAlg}}(\Sym(\ba))\simeq \Sym(\Phi^L(\ba)).$$

And similarly, $\Phi_1^{L,\on{ComAlg}}(\Sym(\ba_1))\simeq \Sym(\Phi_1^L(\ba_1))$. Hence, the map
\eqref{e:left adj fact homology}, evaluated on $\Sym(\ba)$ identifies with
$$\Sym(\Phi_1^L(\Psi(\ba_1))) \to \Psi'(\Sym(\Phi^L(\ba)))\simeq \Sym(\Psi(\Phi^L(\ba))),$$
which is an isomorphism by assumption.

\end{proof}

\ssec{Proofs of Propositions \ref{p:coHom almost semi-rigid} and \ref{p:R coHom almost semi-rigid}}

\sssec{} \label{sss:producing R}

Let $\bA$ be a symmetric monoidal DG category, and let $R_\bA\in \bA$ be a commutative algebra object. 
Let $\bH$ be a symmetric monoidal DG category, assumed dualizable as a plain DG category. 

\medskip

Consider the symmetric monoidal DG categories $\coHom(\bA,\bH)$ and $\coHom(R_\bA\mod(\bA),\bH)$.
Define the commutative algebra object $R_{\coHom(\bA,\bH)}\in \coHom(\bA,\bH)$
as follows:

\medskip

It is the value of the functor
$$(\on{Id}\otimes \one_\bH)^{L,\on{ComAlg}}:\on{ComAlg}(\coHom(\bA,\bH)\otimes \bH)\to \on{ComAlg}(\coHom(\bA,\bH))$$
(see \lemref{l:fact homology} for the notations) 
on the image of $R_\bA$ along the tautological symmetric monoidal functor
$$\bA\to \coHom(\bA,\bH)\otimes \bH.$$

\sssec{}

Unwinding the definitions and using \lemref{l:left adj fact homology}, we obtain:

\begin{lem} \label{l:affine functor}
There is a canonical equivalence
$$\coHom(R_\bA\mod(\bA),\bH) \simeq R_{\coHom(\bA,\bH)}\mod(\coHom(\bA,\bH)),$$
so that the symmetric monoidal functor
$$\coHom(\bA,\bH)\to \coHom(R_\bA\mod(\bA),\bH),$$
attached by the functoriality of $\coHom(-,\bH)$ to the symmetric monoidal functor $\bA\to R_\bA\mod(\bA)$
corresponds to the symmetric monoidal functor
$$\coHom(\bA,\bH)\to R_{\coHom(\bA,\bH)}\mod(\coHom(\bA,\bH)).$$
\end{lem} 

\sssec{} \label{sss:R H colim}

Let $\bH$ be as in \secref{sss:again assumptions on H new}. 
By \lemref{l:fact homology} we obtain the following explicit description for the object $R_{\coHom(\bA,\bH)}$: 

\medskip

Let
us identify $\coHom(\bA,\bH)$ with
$$\underset{(I\to J)\in \on{TwArr}(\fSet)}{\on{colim}}\, \bA^{\otimes I} \otimes (\bH^\vee)^{\otimes J},$$
see \lemref{l:coEnd abs}. Denote by
$$\on{ins}_{I\to J}:\bA^{\otimes I} \otimes (\bH^\vee)^{\otimes J}\to \coHom(\bA,\bH)$$
the resulting tautological functors. 

\medskip

Unwinding the definitions, we obtain: 
\begin{equation} \label{e:A  C H}
R_{\coHom(\bA,\bH)} \simeq 
\underset{(I\to J)\in \on{TwArr}(\fSet)}{\on{colim}}\, \on{ins}_{I\to J}(R_\bA^{\otimes I}\otimes (\coinv_\bH)^{\otimes J}),
\end{equation}
where we regard $\coinv_H$ as an object of $\bH^\vee$. 

\sssec{}

We return to the setting of \secref{ss:proj coHom}. 
Set $\bA:=\CC\otimes \CC$ and $R_\bA=\sR_\CC$. 
By Barr-Beck-Lurie, the functor $\on{comult}_\CC$ identifies 
$$\CC\simeq \sR_\CC\mod(\CC\otimes \CC).$$

Denote the resulting commutative algebra object of $\coHom(\CC\otimes \CC,\bH)$, constructed in \secref{sss:producing R},
by $\wt\sR_{\coHom(\CC,\bH)}$. By \lemref{l:affine functor}, we obtain an equivalence 
$$\coHom(\CC,\bH) \simeq \wt\sR_{\coHom(\CC,\bH)}\mod(\coHom(\CC\otimes \CC,\bH)),$$
so that the symmetric monoidal functor
\begin{equation} \label{e:coHom mult}
\coHom(\CC\otimes \CC,\bH) \overset{\on{mult}_\CC}\longrightarrow \coHom(\CC,\bH)
\end{equation} 
corresponds to the symmetric monoidal functor
$$\coHom(\CC\otimes \CC,\bH)\to \wt\sR_{\coHom(\CC,\bH)}\mod(\coHom(\CC\otimes \CC,\bH)).$$

In particular, we obtain that the right adjoint of the functor \eqref{e:coHom mult} identifies with the forgetful functor
$$\wt\sR_{\coHom(\CC,\bH)}\mod(\coHom(\CC\otimes \CC,\bH))\to \coHom(\CC\otimes \CC,\bH),$$
and hence is continuous and 
respects the $\coHom(\CC\otimes \CC,\bH)$-module structure.

\medskip

Furthermore, the value of the right adjoint to \eqref{e:coHom mult} on $\one_{\coHom(\CC,\bH)}$
is $\wt\sR_{\coHom(\CC,\bH)}$. Hence, if $\bH$ satisfies the assumption of \secref{sss:again assumptions on H new}, by 
\secref{sss:R H colim}, it is given by 
$$\underset{(I\to J)\in \on{TwArr}(\fSet)}{\on{colim}}\, \on{ins}_{I\to J}((\sR_\CC)^{\otimes I}\otimes \coinv_\bH^{\otimes J})\in
\coHom(\CC\otimes \CC,\bH).$$

\sssec{} \label{sss:coHom proof}

Note that for a pair of symmetric monoidal categories $\CC_1$ and $\CC_2$
we have a natural identification
$$\coHom(\CC_1,\bH)\otimes \coHom(\CC_2,\bH)\simeq \coHom(\CC_1\otimes \CC_2,\bH).$$
Let us take $\CC_1=\CC_2=\CC$. We obtain an equivalence
\begin{equation} \label{e:coHom double}
\coHom(\CC,\bH)\otimes \coHom(\CC,\bH)\simeq \coHom(\CC\otimes \CC,\bH),
\end{equation}
so that the functor \eqref{e:coHom mult} identifies with
$$\on{mult}_{\coHom(\CC,\bH)}:\coHom(\CC,\bH)\otimes \coHom(\CC,\bH) \to \coHom(\CC,\bH).$$

Hence, we obtain that its right adjoint has the properties specified in \propref{p:coHom almost semi-rigid}.

\sssec{} \label{sss:R coHom proof} 

Note that under the equivalence \eqref{e:coHom double}, the object
$$\wt\sR_{\coHom(\CC,\bH)}\in \on{ComAlg}(\coHom(\CC\otimes \CC,\bH))$$
corresponds to the object 
$$\sR_{\coHom(\CC,\bH)}\in \on{ComAlg}(\coHom(\CC,\bH)\otimes \coHom(\CC,\bH)).$$

Finally, it is easy to see that under the equivalence \eqref{e:coHom double}, the colimit expression
$$\underset{(I\to J)\in \on{TwArr}(\fSet)}{\on{colim}}\, \on{ins}_{I\to J}((\sR_\CC)^{\otimes I}\otimes \coinv_\bH^{\otimes J})$$
coincides term-wise with
$$\underset{(I\to J)\in \on{TwArr}(\fSet)}{\on{colim}}\, 
(\on{ins}_{I\to J}\otimes \on{ins}_{I\to J})\left((\sR_\CC)^{\boxtimes I} \otimes (\sR_{\bH^\vee})^{\boxtimes J}\right).$$

This proves \propref{p:R coHom almost semi-rigid}. 

\ssec{Applications to $\CC_\Ran^{\on{Betti}}$}

\sssec{} \label{sss:proj Betti}

Let us apply the results of \secref{ss:proj coHom} to $\bH=\Shv^{\on{all}}(X)$. From \propref{p:coHom almost semi-rigid} we obtain that
(although the category $\CC_\Ran^{\on{Betti}}$ is not rigid) the functor
$$\on{mult}_{\CC_\Ran^{\on{Betti}}}:\CC_\Ran^{\on{Betti}}\otimes \CC_\Ran^{\on{Betti}} \to \CC_\Ran^{\on{Betti}}$$
admits a continuous right adjoint, to be denoted $\on{comult}_{\CC_\Ran^{\on{Betti}}}$.

\medskip

Moreover, we obtain that the structure on $\on{comult}_{\CC_\Ran^{\on{Betti}}}$ of right-lax compatibility with the 
$\CC_\Ran^{\on{Betti}}$-bimodule structure structure is strict. Denote 
$$\sR_{\CC,\Ran}^{\on{Betti}}:=\on{comult}_{\CC_\Ran^{\on{Betti}}}(\one_{\CC_\Ran^{\on{Betti}}})\in \CC_\Ran^{\on{Betti}}\otimes \CC_\Ran^{\on{Betti}}.$$

We obtain that $\sR_{\CC,\Ran}^{\on{Betti}}$ naturally lifts to an object of
$$\on{HC}^\bullet(\CC_\Ran^{\on{Betti}},\CC_\Ran^{\on{Betti}}\otimes \CC_\Ran^{\on{Betti}}),$$
i.e., it is equipped with Hecke isomorphisms \eqref{e:Hecke univ}.

\sssec{}

Note that with respect to the canonical self-duality of $\Shv^{\on{all}}(X)$ (see \secref{sss:Shv Betti self-dual}), the object
$$\sR_{\Shv^{\on{all}}(X)^\vee}\in \Shv^{\on{all}}(X)^\vee\otimes \Shv^{\on{all}}(X)^\vee \simeq 
\Shv^{\on{all}}(X)\otimes \Shv^{\on{all}}(X)\simeq \Shv^{\on{all}}(X\times X)$$
identifies with 
$$(\Delta_X)_!(\omega_X).$$

\medskip

Hence, we obtain that $\sR_{\CC,\Ran}^{\on{Betti}}$ is given by the formula similar to \eqref{e:progenitor}, namely
\begin{equation} \label{e:progenitor Betti}
\sR_{\CC,\Ran}^{\on{Betti}}\simeq
\underset{(I\to J)\in \on{TwArr}(\fSet)}{\on{colim}}\,
(\on{ins}_{I\to J}\otimes \on{ins}_{I\to J})(\sR^{\boxtimes I}_\CC\otimes (\Delta_{X^J})_!(\omega_{X^J})).
\end{equation} 

\begin{rem}
Note the difference between formulas \eqref{e:progenitor} and \eqref{e:progenitor Betti}:
in the latter we have the objects
$$(\Delta_{X^J})_!(\omega_{X^J})\simeq (\Delta_{X^J})_*(\omega_{X^J})\in \Shv(X^J)\subset \Shv^{\on{all}}(X^J),$$
while in the former we have
$$\on{u}_{\Shv(X^J)}\in \Shv(X^J)\otimes \Shv(X^J).$$

Note also that, unlike the constructible contexts and that of $\Dmod(-)$, the 
object $$(\Delta_X)_!(\omega_X)\in \Shv^{\on{all}}(X)\otimes \Shv^{\on{all}}(X)$$ is \emph{not}
the unit of the canonical self-duality on $\Shv^{\on{all}}(X)$. The unit is given by
$(\Delta_X)_!(\ul\sfe_X)$.

\end{rem} 

\ssec{Applications to $\CC_\Ran$}

\sssec{} \label{sss:R Ran vs R lisse}

Let us apply the results in \secref{ss:proj coHom} to $\bH=\qLisse(X)$. Recall that 
$$\coHom(\CC,\qLisse(X))\simeq \CC^{\otimes X\on{-lisse}},$$
and consider 
the corresponding object
$$\sR_{\coHom(\CC,\qLisse(X))}\in \coHom(\CC,\qLisse(X))\otimes \coHom(\CC,\qLisse(X))\simeq 
\CC^{\otimes X\on{-lisse}}\otimes \CC^{\otimes X\on{-lisse}}.$$

\sssec{}

Recall the functor 
\begin{equation} \label{e:Ran to lisse again} 
\CC_\Ran\to \CC^{\otimes X\on{-lisse}}
\end{equation} 
of \eqref{e:Ran to lisse}, and let us denote by
$$\wt\sR_{\coHom(\CC,\qLisse(X))}\in \CC^{\otimes X\on{-lisse}}\otimes \CC^{\otimes X\on{-lisse}}$$
the image of
$$\sR_{\CC,\Ran}\in \CC_\Ran\otimes \CC_\Ran$$
along the tensor square of \eqref{e:Ran to lisse again}
$$\CC_\Ran\otimes \CC_\Ran\to \CC^{\otimes X\on{-lisse}}\otimes \CC^{\otimes X\on{-lisse}}.$$

By adjunction, we obtain a map of commutative algebras in $\CC^{\otimes X\on{-lisse}}\otimes \CC^{\otimes X\on{-lisse}}$
\begin{equation} \label{e:R ran to R lisse}
\wt\sR_{\coHom(\CC,\qLisse(X))}\to \sR_{\coHom(\CC,\qLisse(X))}.
\end{equation}

\sssec{}

We will prove:

\begin{prop} \label{p:R ran to R lisse}
The map \eqref{e:R ran to R lisse} is an isomorphism.
\end{prop}

\begin{proof}

We will show that in terms of presentations of $\wt\sR_{\coHom(\CC,\qLisse(X))}$ and $\sR_{\coHom(\CC,\qLisse(X))}$
as colimits, given by formulas \eqref{e:progenitor} and \eqref{e:progenitor coHom formula}, respectively, the map 
\eqref{e:R ran to R lisse} is a term-wise isomorphism.

\medskip

For the latter, we need to show that for a given finite set $J$, the functor dual to
$$\qLisse(X)^{\otimes J}\otimes \qLisse(X)^{\otimes J} \to \Shv(X^J)\otimes \Shv(X^J)$$
sends 
$$\on{u}_{\Shv(X^J)}\in \Shv(X^J)\otimes \Shv(X^J)\simeq \Shv(X^J)^\vee\otimes \Shv(X^J)^\vee$$ to the object 
$$(\sR_{\qLisse(X)^\vee})^{\boxtimes J}\in (\qLisse(X)^\vee)^{\otimes J}\otimes (\qLisse(X)^\vee)^{\otimes J}.$$

Note that the object $\on{u}_{\Shv(X^J)}\in \Shv(X^J)^\vee\otimes \Shv(X^J)^\vee$ equals the value of the functor
dual to
$$\Shv(X^J)\otimes \Shv(X^J) \overset{\sotimes}\to \Shv(X^J)$$
on $\on{C}^\cdot(X^J,-)$, viewed as an object of $\Shv(X^J)^\vee$. 

\medskip

The required assertion follows now by passing to dual functors in the commutative diagram
$$
\CD
\Shv(X^J)\otimes \Shv(X^J) @>>{\sotimes}>   \Shv(X^J) \\
@A{\text{\eqref{e:!-emb}}\otimes \text{\eqref{e:!-emb}}}AA @AA{\text{\eqref{e:!-emb}}}A \\
\qLisse(X)^{\otimes J} \otimes \qLisse(X)^{\otimes J}  @>>> \qLisse(X)^{\otimes J},
\endCD
$$
using the fact that the composition
$$\qLisse(X)^{\otimes J} \overset{\text{\eqref{e:!-emb}}}\longrightarrow \Shv(X^J) \overset{\on{C}^\cdot(X^J,-)}\longrightarrow \Vect_\sfe$$
identifies with $(\coinv_{\qLisse(X)})^{\otimes J}$.

\end{proof}

\begin{rem} \label{r:R Ran vs R lisse Betti}
An analog of the construction in \secref{sss:R Ran vs R lisse} applies when instead of the pair $(\CC_\Ran,\coHom(\CC,\qLisse(X))$
we take $(\CC_\Ran^{\on{Betti}},\coHom(\CC,\Shv^{\on{all}}_{\on{loc.const}}(X))$. 

\medskip

An assertion parallel to \propref{p:R ran to R lisse} continues to hold in this context, with the same proof.

\end{rem} 

\ssec{Identification of the diagonal} \label{ss:diagonal}

\sssec{} \label{sss:diagonal}

In the setting of \secref{ss:proj coHom}, let us take $\CC=\Rep(\sG)$, and let $\bH$ be as
in \secref{sss:map abs}. Consider the functor
\begin{equation} \label{e:Hom abs again}
\coHom(\Rep(\sG),\bH)\to \QCoh(\bMaps(\Rep(\sG),\bH))
\end{equation}
of \eqref{e:Hom abs}.

\medskip

Consider the object
$$\sR_{\coHom(\Rep(\sG),\bH)} \in \coHom(\Rep(\sG),\bH)\otimes \coHom(\Rep(\sG),\bH),$$
see \eqref{e:progenitor coHom}. 

\sssec{}

Let us denote by
$$\sR_{\bMaps(\Rep(\sG),\bH)}\in \QCoh(\bMaps(\Rep(\sG),\bH))\otimes \QCoh(\bMaps(\Rep(\sG),\bH))$$
the image of $\sR_{\coHom(\Rep(\sG),\bH)}$ along the tensor square of the functor \eqref{e:Hom abs again} 
$$\coHom(\Rep(\sG),\bH)\otimes \coHom(\Rep(\sG),\bH)\to \QCoh(\bMaps(\Rep(\sG),\bH))\otimes \QCoh(\bMaps(\Rep(\sG),\bH)).$$

\medskip

By adjunction, we obtain a map
\begin{equation} \label{e:diag abs}
\sR_{\bMaps(\Rep(\sG),\bH)} \to (\Delta_{\bMaps(\Rep(\sG),\bH)})_*(\CO_{\bMaps(\Rep(\sG),\bH)})
\end{equation}
of commutative algebras in
\begin{multline*} 
\QCoh(\bMaps(\Rep(\sG),\bH))\otimes \QCoh(\bMaps(\Rep(\sG),\bH))\simeq \\
\simeq \QCoh(\bMaps(\Rep(\sG),\bH)\times \bMaps(\Rep(\sG),\bH)).
\end{multline*}

The goal of this subsection is to prove the following result:

\begin{thm} \label{t:diag abs}
The map \eqref{e:diag abs} is an isomorphism.
\end{thm}

\begin{rem}
Note that if $\bH$ is adapted for spectral decomposition, the statement of \thmref{t:diag abs} is tautological.

\medskip

In general, we hope that \thmref{t:diag abs} goes some way in the direction of the proof of \conjref{c:Hom abs}. 
\end{rem} 

\begin{rem} 

Our main interest is the case when $\bH=\qLisse(X)$ (provided that $\qLisse(X)$ is dualizable). 
Note, however, that by the previous remark, if $X$ is a complete algebraic curve, 
the assertion of \thmref{t:diag abs} in this case is already known, due to \thmref{t:spectral}. 

\end{rem}

\begin{rem}
Note that \thmref{t:diag abs} is applicable also to $\bH$ being  $\Dmod(X)$ or $\Shv_{\on{loc.const}}(X)$, and hence 
it gives rise to a description of
$$(\Delta_{\LocSys^?_\sG(X)})_*(\CO_{\LocSys^?_\sG(X)})\in 
\QCoh(\LocSys^?_\sG(X))\otimes \QCoh(\LocSys^?_\sG(X)$$
for $?$ being $\dr$ or $\on{Betti}$ in terms of $\sR_{\Rep(\sG),\Ran}$, viewed as an object in 
$\Rep(\sG)^\dr_\Ran$ or $\Rep(\sG)^{\on{Betti}}_\Ran$. 
respectively. 
\end{rem} 

\sssec{}

The rest of this subsection is devoted to the proof of \thmref{t:diag abs}. 

\medskip

Let $S$ be an affine scheme, equipped with two maps
$$\sigma_i:S\to \bMaps(\Rep(\sG),\bH),$$
corresponding to symmetric monoidal functors
$$\sF_i:\Rep(\sG)\to \QCoh(S)\otimes \bH.$$

Let us denote by
$$\sR_{\sigma_1,\sigma_2}\in \QCoh(S)$$ 
the commutative object equal to the pullback by means of
$$S \overset{\sigma_1,\sigma_2}\longrightarrow \bMaps(\Rep(\sG),\bH)) \times \bMaps(\Rep(\sG),\bH)$$
of the object $\sR_{\bMaps(\Rep(\sG),\bH)}$. 

\medskip

We need to show that the space of homomorphisms of commutative algebras
$$\sR_{\sigma_1,\sigma_2}\to \CO_S$$
identifies canonically with the space of isomorphisms of symmetric monoidal functors $\sF_1\simeq \sF_2$. 

\sssec{} \label{sss:to prove fact homology}

Consider the following general situation. Let $\bA$ be a symmetric monoidal category, and let us be given a pair of symmetric monoidal
functors
$$\sF_1,\sF_2:\Rep(\sG)\to \bA.$$

Consider the commutative algebra object $\sR_{\sF_1,\sF_2}\in \bA$, obtained by applying 
the  symmetric monoidal functor
$$\Rep(\sG)\otimes \Rep(\sG) \overset{\sF_1\otimes \sF_2}\longrightarrow \bA\otimes \bA \overset{\on{mult}_\bA}\longrightarrow \bA$$
to the regular representation $\sR_{\Rep(\sG)}\in \Rep(\sG)\otimes \Rep(\sG)$. 

\medskip

Then it is easy to see that the space of isomorphisms between $\sF_1$ and $\sF_2$ identifies canonically with the space of 
maps of commutative algebras
$$\sR_{\sF_1,\sF_2}\to \one_\bA.$$

\sssec{} \label{sss:reduce to fact homology}

Hence, we need to show that the space of homomorphisms of commutative algebras in $\QCoh(S)$
$$\sR_{\sigma_1,\sigma_2}\to \CO_S$$
is canonically isomorphic to the space of maps of commutative algebras in 
$\bH\otimes \QCoh(S)$ 
$$\on{mult}_{\QCoh(S)\otimes \bH}\circ (\sF_1\otimes \sF_2)(\sR_{\Rep(\sG)}) \to \CO_S\otimes \one_\bH.$$

\sssec{}

Recall the setting of \secref{sss:left adj fact homology}. Set
$$\bA':=\coHom(\Rep(\sG)\otimes \Rep(\sG),\bH), \quad \bA:=\coHom(\Rep(\sG)\otimes \Rep(\sG),\bH) \otimes \bH, \quad 
\Phi=-\otimes \one_\bH,$$
$$\bA'_1=\QCoh(S), \quad \bA_1:=\QCoh(S)\otimes \bH, \quad \Phi_1=-\otimes \one_\bH,$$
with $\Psi'$ being the functor
$$\coHom(\Rep(\sG)\otimes \Rep(\sG),\bH) \simeq \coHom(\Rep(\sG),\bH)\otimes \coHom(\Rep(\sG),\bH)\to \\$$
$$\to \QCoh(\bMaps(\Rep(\sG),\bH)))\otimes \QCoh(\bMaps(\Rep(\sG),\bH))) \overset{(\sigma_1,\sigma_2)^*}\longrightarrow \QCoh(S),$$
and $\Psi=\Psi'\otimes \on{Id}_\bH$. 

\medskip

Note that the functor 
$$\Rep(\sG)\otimes \Rep(\sG)\to \QCoh(S)\otimes \bH,$$
corresponding to $\Psi'$ by adjunction, is the functor $\on{mult}_{\QCoh(S)\otimes \bH}\circ (\sF_1\otimes \sF_2)$. 

\medskip

Take 
$$R\in \on{ComAlg}(\coHom(\Rep(\sG)\otimes \Rep(\sG),\bH) \otimes \bH)$$
to be equal to the image of $\sR_{\Rep(\sG)}$ under the tautological functor
$$\Rep(\sG)\otimes \Rep(\sG)\to \coHom(\Rep(\sG)\otimes \Rep(\sG),\bH)\otimes \bH.$$

Then
$$\Psi(R)\simeq \on{mult}_{\QCoh(S)\otimes \bH}\circ (\sF_1\otimes \sF_2)(\sR_{\Rep(\sG)}) .$$

\medskip

The required assertion follows now evaluating both sides of \lemref{l:left adj fact homology} on the above 
object $R$. 

\qed[\thmref{t:diag abs}]

\ssec{Localization on $\LocSys_\sG(X)$} \label{ss:Loc}

\sssec{}

Consider the symmetric monoidal functors
$$\Rep(\sG)_\Ran \to \coHom(\Rep(\sG),\qLisse(X))$$
of \eqref{e:Ran to lisse again}
and
$$\coHom(\Rep(\sG),\qLisse(X))\to \QCoh(\LocSys_\sG^{\on{restr}}(X))$$
of \eqref{e:Hom abs again}.

\medskip

Let us denote their composition by
$$\on{Loc}:\Rep(\sG)_\Ran\to \QCoh(\LocSys_\sG^{\on{restr}}(X)).$$

\sssec{}

Explicitly, the functor $\on{Loc}$ sends an object of $\Rep(\sG)_\Ran$ of the form 
$$\on{ins}_{I\overset{\psi}\to J}(V \otimes \CF), \quad V\in \Rep(\sG)^{\otimes I},\, \CF\in \Shv(X^J)$$
to
$$(\on{Id}_{\QCoh(\LocSys_\sG^{\on{restr}}(X))}\otimes \on{C}^\cdot(X^I,-))\left(\CE^I_V\otimes (\Delta_\psi)_*(\CF)\right),$$
where
$$\CE^I_V\in \QCoh(\LocSys_\sG^{\on{restr}}(X))\otimes \qLisse(X)^{\otimes I}$$
is as in \secref{sss:pre-shtuka}. 

\sssec{}

As in \secref{sss:diagonal}, by adjunction, we obtain a map of commutative algebras:
\begin{equation} \label{e:Loc diag}
(\on{Loc}\otimes \on{Loc})(\sR_{\Rep(\sG),\Ran}) \to (\Delta_{\LocSys^{\on{restr}}_\sG(X)})_*(\CO_{\LocSys^{\on{restr}}_\sG(X)}).
\end{equation}

Combining \propref{p:R ran to R lisse} and \thmref{t:diag abs}, we obtain: 

\begin{thm} \label{t:diagonal}
The above map \eqref{e:Loc diag} is an isomorphism.
\end{thm}

\ssec{Tensor products over $\Rep(\sG)$ vs. $\QCoh(\LocSys^{\on{restr}}_\sG(X))$}

In this subsection we will use \thmref{t:diagonal} to deduce some results on the relationship
between module categories $\Rep(\sG)$ vs. those over $\QCoh(\LocSys^{\on{restr}}_\sG(X))$.

\sssec{}

We claim:

\begin{prop} \label{p:LocSys over ran}
The tensor product functor 
\begin{equation} \label{e:LocSys over ran}
\QCoh(\LocSys^{\on{restr}}_\sG(X))\underset{\Rep(\sG)_\Ran}\otimes 
\QCoh(\LocSys^{\on{restr}}_\sG(X)) \to \QCoh(\LocSys^{\on{restr}}_\sG(X))
\end{equation}
is an equivalence.
\end{prop}

\begin{proof}

This is a formal consequence of \thmref{t:diagonal}. Namely, let $\bA,\bA'$ be a pair of symmetric monoidal 
categories such that the functors
$$\on{mult}_\bA:\bA\otimes \bA\to \bA \text{ and } 
\on{mult}_{\bA'}:\bA'\otimes \bA'\to \bA'$$
both admit continuous right adjoints, which also respect the bimodule structure.

\medskip

Let $\Phi:\bA\to \bA'$ be a symmetric monoidal functor, such that the resulting map 
$$(\Phi\otimes \Phi)(\on{mult}^R_\bA(\one_\bA))\to \on{mult}^R_{\bA'}(\one_{\bA'}),$$
obtained by adjunction, is an isomorphism.

\medskip

Then we claim that the functor
$$(\bA'\otimes \bA')\underset{\bA\otimes \bA}\otimes \bA \simeq \bA'\underset{\bA}\otimes \bA' \to \bA'$$
is an equivalence.

\medskip

Indeed, the pairs
$$\bA'\otimes \bA' \simeq (\bA'\otimes \bA')\underset{\bA\otimes \bA}\otimes (\bA\otimes \bA) \rightleftarrows
(\bA'\otimes \bA')\underset{\bA\otimes \bA}\otimes \bA$$
and 
$$\bA'\otimes \bA'  \rightleftarrows \bA'$$
are monadic, and the corresponding monads are given by tensoring with
$$(\Phi\otimes \Phi)(\on{mult}^R_\bA(\one_\bA)) \text{ and } \on{mult}^R_{\bA'}(\one_{\bA'}),$$
respectively.

\end{proof}  

\sssec{}

As a formal consequence of \propref{p:LocSys over ran}, we obtain:

\begin{cor} \label{c:LocSys over ran} 
Let $\bM_1$ and $\bM_2$ be module categories over $\QCoh(\LocSys^{\on{restr}}_\sG(X))$. 

\smallskip

\noindent{\em(a)} The functor
$$\bM_1\underset{\Rep(\sG)_\Ran}\otimes \bM_2\to
\bM_1\underset{\QCoh(\LocSys^{\on{restr}}_\sG(X))}\otimes \bM_2$$
is an equivalence.

\smallskip

\noindent{\em(b)} The functor
$$\on{Funct}_{\QCoh(\LocSys^{\on{restr}}_\sG(X))\mod}(\bM_1,\bM_2)\to
\on{Funct}_{\Rep(\sG)_\Ran\mod}(\bM_1,\bM_2)$$
is an equivalence.
\end{cor}

\begin{proof}

Both assertions hold in the general context in which we proved \propref{p:LocSys over ran}:

\medskip

For a pair of $\bA'$-module categories $\bM_1,\bM_2$, we have
$$\bM_1\underset{\bA}\otimes \bM_2\simeq
(\bM_1\otimes \bM_2)\underset{\bA\otimes \bA}\otimes \bA\simeq
(\bM_1\otimes \bM_2)\underset{\bA'\otimes \bA'}\otimes (\bA'\otimes \bA')\underset{\bA\otimes \bA}\otimes \bA
\overset{\propref{p:LocSys over ran}}\simeq $$
$$\simeq (\bM_1\otimes \bM_2)\underset{\bA'\otimes \bA'}\otimes \bA'\simeq \bM_1\underset{\bA'}\otimes \bM_2.$$

$$\on{Funct}_{\bA\mod}(\bM_1,\bM_2)\simeq
\on{Funct}_{(\bA\otimes \bA)\mod}(\bA,\on{Funct}(\bM_1,\bM_2))\simeq$$
$$\simeq \on{Funct}_{(\bA'\otimes \bA')\mod}((\bA'\otimes \bA')\underset{\bA\otimes \bA}\otimes \bA,\on{Funct}(\bM_1,\bM_2))\overset{\propref{p:LocSys over ran}}\simeq $$
$$\simeq \on{Funct}_{(\bA'\otimes \bA')\mod}(\bA',\on{Funct}(\bM_1,\bM_2)) \simeq \on{Funct}_{\bA'\mod}(\bM_1,\bM_2).$$

\end{proof} 

\begin{rem}

An analog of the functor $\Loc$ exists also in the context of $\Dmod(-)$, in which case this is the functor
$$\Rep(\sG)^\dr_\Ran \simeq \coHom(\Rep(\sG),\Dmod(X))\overset{\text{\eqref{e:Hom abs again}}}\longrightarrow \QCoh(\LocSys^\dr_\sG(X)).$$

\medskip

Note that in this context, the functor $\Loc$ admits a continuous and fully faithful right adjoint. 

\medskip

A counterpart of \thmref{t:diagonal} in this case is \thmref{t:diag abs}. Hence, analogs of \propref{p:LocSys over ran} and \corref{c:LocSys over ran}
continue to hold in this context as well. 

\end{rem}

\begin{rem}

An analog of the functor $\Loc$ exists also in the context of $\Shv^{\on{all}}(-)$, in which case this is the functor
$$\Rep(\sG)^{\on{Betti}}_\Ran \simeq \coHom(\Rep(\sG),\Shv^{\on{all}}(X))\to 
\coHom(\Rep(\sG),\Shv_{\on{loc.const}}^{\on{all}}(X)) \overset{\text{\eqref{e:Hom abs again}}}\longrightarrow $$
$$\to \QCoh(\LocSys^{\on{Betti}}_\sG(X)).$$
where the last arrow is an equivalence by \thmref{t:spectral Betti}. 

\medskip

An analog of \thmref{t:diagonal} continues to hold in this context, see Remark \ref{r:R Ran vs R lisse Betti}. 
Hence, analogs of \propref{p:LocSys over ran} and \corref{c:LocSys over ran}
continue to hold in this context as well. 

\end{rem}

\section{Spectral projector and Hecke eigen-objects} \label{s:projector abstract}

Having introduced the progenitor of the projector in \secref{ss:progenitor}, we now proceed to
the definition of the projector itself. 

\ssec{Beilinson's spectral projector--abstract form} \label{ss:proj abs}

In this subsection we will finally define what we mean by the category of Hecke
eigen-objects, and introduce Beilinson's spectral projector. 


\sssec{} 

Let $\CZ$ be a prestack (over the field of coefficients $\sfe$). Let us be given a symmetric monoidal functor
\begin{equation} \label{e:eigencond abs}
\sF:\CC\to \QCoh(\CZ)\otimes \qLisse(X).
\end{equation} 

\sssec{Example}
Note that if $\CC=\Rep(\sG)$, and if $\sF$ is right t-exact\footnote{For any prestack $\CZ$, the category $\QCoh(\CZ)$
carries a canonically defined t-structure in which an object is connective if and only if its pullback to any affine scheme 
is connected.}, the above datum is equivalent to that of a map
$$\CZ\to \LocSys^{\on{restr}}_\sG(X).$$

\sssec{} \label{sss:tilde F}

We can interpret $\sF$ as a symmetric monoidal functor
$$\CC^{\otimes X\on{-lisse}}:=\coHom(\CC,\qLisse(X))\to \QCoh(\CZ).$$

Precomposing with \eqref{e:Ran to lisse}, we obtain a functor
\begin{equation} \label{e:functor F tilde}
\wt\sF:\CC_\Ran\to \QCoh(\CZ).
\end{equation} 

Let us describe the functor $\wt\sF$ explicitly. Its value on an object of $\CC_\Ran$ of the form
$$\on{ins}_{I\overset{\psi}\to J}(V \otimes \CF), \quad V\in \CC^{\otimes I},\, \CF\in \Shv(X^J)$$
is
$$\left(\on{Id}_{\QCoh(\CZ)}\otimes \on{C}^\cdot(X^J,-)\right)
\left(\sF^J(\on{mult}^\psi_\CC(V))\otimes \CF\right),$$
where:

\begin{itemize}

\item $\on{mult}^\psi_\CC$ is the tensor product functor $\CC^{\otimes I}\to \CC^{\otimes J}$
along the fibers of $\psi$;

\item $\sF^J$ is the functor $\CC^{\otimes J}\to \QCoh(\CZ)\otimes \qLisse(X)^{\otimes J}$ obtained
from $\sF$. 

\item $-\otimes \CF$ refers to the action of $\qLisse(-)$ on $\Shv(-)$ by tensor products\footnote{This is
either the $\overset{*}\otimes$ tensor product or, equivalently, the $\sotimes$ tensor product precomposed
with \eqref{e:!-emb}.}.

\end{itemize} 

\sssec{} \label{sss:Hecke abs}

Let $\bM$ be a module category over $\CC_\Ran$. We will denote the action functor
$$\CC_\Ran\otimes \bM\to \bM$$
by 
$$\CV,\bm\mapsto \CV\star \bm.$$

\medskip

We define the category of Hecke eigen-objects 
in $\bM$ with respect to \eqref{e:eigencond abs}, to be denoted  
$$\on{Hecke}(\CZ,\bM)_\sF$$
(or simply $\on{Hecke}(\CZ,\bM)$ if no ambiguity is likely to occur), to be
$$\on{HC}^\bullet(\CC_\Ran,\bM\otimes \QCoh(\CZ)):=
\on{Funct}_{(\CC_\Ran\otimes \CC_\Ran)\mod}(\CC_\Ran,\bM\otimes \QCoh(\CZ)),$$
where we regard $\bM\otimes \QCoh(\CZ)$ as a bimodule over $\CC_\Ran$. 

\medskip

Note that we can also rewrite
\begin{equation} \label{e:Hecke category alt}
\on{Hecke}(\CZ,\bM) \simeq
\on{Funct}_{(\CC_\Ran\otimes \QCoh(\CZ))\mod}(\QCoh(\CZ),\bM\otimes \QCoh(\CZ)).
\end{equation}

\medskip

By definition, we can also interpret $\on{Hecke}(\CZ,\bM)$ as the category of objects
$\bm\in \bM\otimes \QCoh(\CZ)$ equipped with a tensor-compatible system of isomorphisms
\begin{equation} \label{e:Hecke M}
\CV\star \bm \simeq \bm\otimes \wt\sF(\CV), \quad \CV\in \CC_\Ran.
\end{equation} 

We can regard the system of isomorphisms \eqref{e:Hecke M} as an abstract form of the
Hecke eigen-property. 

\sssec{} \label{sss:Hecke left adj}

The adjunction
\begin{equation} \label{e:mult adj Ran}
\on{mult}_{\CC_\Ran}:\CC_\Ran\otimes \CC_\Ran \rightleftarrows \CC_\Ran:\on{comult}_{\CC_\Ran}
\end{equation} 
as $\CC_\Ran$-bimodule categories induces an adjunction
\begin{equation} \label{e:Hecke adj}
\ind_{\on{Hecke},\CZ}:\bM\otimes \QCoh(\CZ) \rightleftarrows \on{Hecke}(\CZ,\bM)_\sF:\oblv_{\on{Hecke},\CZ},
\end{equation} 
where $\oblv_{\on{Hecke},\CZ}$ is the tautological forgetful functor. 

\medskip

The functor $\oblv_{\on{Hecke},\CZ}$ is conservative, and hence monadic. 

\sssec{} \label{sss:R object new}

Let $\sR_{\CZ,\sF}$ (or simply $\sR_\CZ$ if no ambiguity is likely to occur) denote the object 
$$(\on{Id}\otimes \wt\sF)(\sR_{\CC,\Ran})\in \CC_\Ran\otimes \QCoh(\CZ),$$
where
$$\sR_{\CC,\Ran} \in \CC_\Ran\otimes \sR_\Ran$$
is as in \secref{sss:R Ran}. 

\medskip

Since $\sR_{\CC,\Ran}$ has a structure of commutative algebra in $\CC_\Ran\otimes \sR_\Ran$
(see \secref{sss:R Ran}), the object $\sR_\CZ$ is naturally a 
commutative algebra in $\sR_{\CC,\Ran}\otimes \QCoh(\CZ)$.

\medskip

We obtain that the monad on $\bM\otimes \QCoh(\CZ)$ corresponding to the adjunction
\eqref{e:Hecke adj} is given by the action of $\sR_{\CZ,\sF}$.

\sssec{} \label{sss:change Z abs}

Let $g:\CZ'\to \CZ$ be a map of prestacks, and let $\sF'$ denote the composite functor
$$\CC\overset{\sF}\to \QCoh(\CZ)\otimes \qLisse(X) \overset{g^*\otimes \on{Id}}\to
\QCoh(\CZ')\otimes \qLisse(X).$$

It follows from the definitions that in this case we have a naturally defined functor
$$g^*:\on{Hecke}(\CZ,\bM)_\sF\to \on{Hecke}(\CZ',\bM)_{\sF'}$$
that makes both diagrams 
$$
\CD
\on{Hecke}(\CZ',\bM)_{\sF'} @>{\oblv_{\on{Hecke},\CZ'}}>> \bM \otimes \QCoh(\CZ')  \\
@A{g^*}AA @AA{\on{Id}\otimes g^*}A \\
\on{Hecke}(\CZ,\bM)_{\sF} @>{\oblv_{\on{Hecke},\CZ}}>> \bM \otimes \QCoh(\CZ) 
\endCD
$$
and 
$$
\CD
\bM \otimes \QCoh(\CZ') @>{\ind_{\on{Hecke},\CZ'}}>> \on{Hecke}(\CZ',\bM)_{\sF'} \\
@A{\on{Id}\otimes g^*}AA @AA{g^*}A \\
\bM \otimes \QCoh(\CZ) @>{\ind_{\on{Hecke},\CZ}}>> \on{Hecke}(\CZ,\bM)_{\sF}
\endCD
$$
commute. 

\sssec{}

Let us denote by $\oblv_{\on{Hecke}}$ the (not necessarily continuous) functor
$$\on{Hecke}(\CZ,\bM)\overset{\oblv_{\on{Hecke},\CZ}}\longrightarrow  \bM\otimes \QCoh(\CZ) \to \bM,$$
where the second arrow is the (not necessarily continuous) right adjoint to
$$\bM \overset{\on{Id}\otimes \CO_\CZ}\longrightarrow \bM\otimes \QCoh(\CZ).$$

\begin{rem}
Suppose for a moment that $\CZ$ is such that $\CO_\CZ \in \QCoh(\CZ)$ is compact
(e.g., $\CZ$ is an algebraic stack), so that the functor 
$$\Gamma(\CZ,-):\QCoh(\CZ)\to \Vect_\sfe$$
is continuous. Then the functor $\oblv_{\on{Hecke}}$ is continuous. Indeed, in the case, 
the corresponding functor $\bM\otimes \QCoh(\CZ) \to \bM$ 
is given by 
$$\bM\otimes \QCoh(\CZ) \overset{\on{Id}_\bM\otimes \Gamma(\CZ,-)}\longrightarrow \bM.$$
\end{rem} 

\sssec{} \label{sss:P enh first}

Consider the functor, to be denoted $\sP_{\CZ,\sF}$ (or simply $\sP_\CZ$ if no ambiguity is likely to occur), 
$$\bM \overset{\on{Id}\otimes \CO_\CZ}\longrightarrow \bM\otimes \QCoh(\CZ)
\overset{\sR_{\CZ,\sF}\star -}\longrightarrow \bM\otimes \QCoh(\CZ).$$

We obtain that the functor $\sP_\CZ$ naturally upgrades to a functor
\begin{equation} \label{e:projector abs}
\sP^{\on{enh}}_\CZ:\bM \to \on{Hecke}(\CZ,\bM),
\end{equation} 
where
$$\sP^{\on{enh}}_\CZ=\ind_{\on{Hecke},\CZ}\circ (\on{Id}\otimes \CO_\CZ).$$


\medskip

By \secref{sss:R object new}, the functor $\sP^{\on{enh}}_\CZ$ provides a left adjoint to $\oblv_{\on{Hecke}}$.

\medskip

The functor \eqref{e:projector abs} is the abstract form of Beilinson's spectral projector: it produces
Hecke eigen-objects from plain objects of $\bM$. 

%
%
%
%
%
%
%
%
%
%
%

\sssec{} \label{sss:R as colim}

Let us write down the object 
$$\sR_\CZ\in \CC_\Ran\otimes \QCoh(\CZ)$$ explicitly as a colimit. By 
\eqref{e:progenitor} and the description of the functor $\wt\sF$ in \secref{sss:tilde F}, it identifies with
\begin{equation} \label{e:formula for R}
\underset{(I\overset{\psi}\to J)\in \on{TwArr}(\fSet)}{\on{colim}}\,
\on{ins}_{I\to J}\left(\left(\on{Id} \otimes (\sF^J\circ \on{mult}^\psi_{\CC})\right)(\sR_\CC^{\boxtimes I})\right),
\end{equation}
where:

\begin{itemize}

\item $\sR_\CC^{\boxtimes I} \in (\CC\otimes \CC)^{\otimes I}\simeq \CC^{\otimes I}\otimes \CC^{\otimes I}$;

\smallskip

\item 
$\sF^J\circ \on{mult}^\psi_{\CC}:\CC^{\otimes I} \to \QCoh(\CZ)\otimes \qLisse(X)^{\otimes J} \simeq 
\qLisse(X)^{\otimes J} \otimes  \QCoh(\CZ)$
is as in \secref{sss:tilde F};

\smallskip

\item We view $\qLisse(X)^{\otimes J}$ as a full subcategory of $\Shv(X^J)$ via the embedding \eqref{e:!-emb}, 

\end{itemize}

\medskip

\noindent so that 
$\left(\on{Id} \otimes (\sF^J\circ \on{mult}^\psi_{\CC})\right)(\sR_\CC^{\boxtimes I})$ is an object of
$$\CC^{\otimes I} \otimes \Shv(X^J)\otimes \QCoh(\CZ).$$

\sssec{} \label{sss:R dr and Betti}

The material in this subsection applies ``as-is'' to $(\CC_\Ran,\qLisse(X))$ replaced by either
$(\CC^\dr_\Ran,\Dmod(X))$ or $(\CC^{\on{Betti}}_\Ran,\Shv^{\on{all}}_{\on{loc.const.}}(X))$.

\ssec{A multiplicativity property of the projector}

In this subsection we will establish a certain multiplicativity property of Beilinson's projector,
namely, \propref{p:prod projectors}, that will be needed in \secref{ss:ten product}.

\sssec{}

Let us be given a pair symmetric monoidal categories $\CC_i$, $i=1,2$ and 
symmetric monoidal functors
$$\sF_i:\CC_i\to \QCoh(\CZ_i)\otimes \qLisse(X),$$
consider
$$\CC:=\CC_1\otimes \CC_2$$
and the corresponding functor
$$\sF:\CC\to \QCoh(\CZ)\otimes \qLisse(X), \quad \CZ=\CZ_1\times \CZ_2.$$

Consider the corresponding objects
$$\sR_{\CZ_i}\in \CC_{i,\Ran}\otimes  \QCoh(\CZ_i), \quad \sR_{\CZ}\in \CC_{\Ran}\otimes  \QCoh(\CZ).$$

Note that we have a naturally defined functor
\begin{equation} \label{e:mult C Ran}
\CC_{1,\Ran}\otimes  \CC_{2,\Ran}\to \CC_{\Ran}\otimes  \CC_{\Ran} \overset{\on{mult}_{\CC_{\Ran}}}\longrightarrow \CC_{\Ran},
\end{equation}

We claim:

\begin{prop} \label{p:prod projectors}
Under the above circumstances, the functor 
\begin{equation} \label{e:mult C Ran bis}
\CC_{1,\Ran}\otimes  \QCoh(\CZ_1) \otimes  \CC_{2,\Ran} \otimes \QCoh(\CZ_2)\to \CC_{\Ran}\otimes  \QCoh(\CZ),
\end{equation}
induced by \eqref{e:mult C Ran} sends 
$$\sR_{\CZ_1}\otimes \sR_{\CZ_2} \mapsto \sR_{\CZ}.$$
\end{prop}

The rest of this subsection is devoted to the proof of \propref{p:prod projectors}.

\sssec{}

First, given a pair of rigid symmetric monoidal categories $\bA_1$ and $\bA_2$ and
$$\bA:=\bA_1\otimes \bA$$
note that we have a canonical isomorphism
\begin{equation} \label{e:R ten prod}
\sR_\bA\simeq \sR_{\bA_1}\otimes \sR_{\bA_2}
\end{equation} 
as objects in 
$$\bA^{\otimes 2}\simeq \bA_1^{\otimes 2}\otimes \bA_2^{\otimes 2}.$$

\medskip

Further, for a symmetric monoidal functor
$$\phi:\bA\to \bA'$$
we have a canonical map 
$$(\phi\otimes \phi)(\sR_\bA)\to \sR_{\bA'}.$$

\sssec{} \label{sss:R R}

Hence, we obtain a map from the image of 
$$\sR_{\CC_1,\Ran}\otimes \sR_{\CC_2,\Ran} \in (\CC_{1,\Ran}\otimes  \CC_{2,\Ran})^{\otimes 2}$$
to
$$\sR_{\CC,\Ran}\in (\CC_\Ran)^{\otimes 2}$$
along the tensor square of the map \eqref{e:mult C Ran}. (Note, however, that this map itself is \emph{not}
an isomorphism; cf. Remark \ref{r:prod dR} below.)

\medskip

The latter map induces a map from the image of $\sR_{\CZ_1}\otimes \sR_{\CZ_2}$ along \eqref{e:mult C Ran bis} to
$\sR_{\CZ}$. We will show that this map is an isomorphism. 

\begin{rem} \label{r:prod dR}

Note that in the case of 
$$\CC^\dr_{\Ran}\simeq \coHom(\CC,\Dmod(X)),$$ the corresponding functor
$$\CC^\dr_{1,\Ran}\otimes  \CC^\dr_{2,\Ran}\to  \CC^\dr_\Ran$$
is already an equivalence. 

\medskip

By \eqref{e:R ten prod}, this implies that the image of 
$$\sR_{\CC_1,\Ran}\otimes \sR_{\CC_2,\Ran} \in (\CC^\dr_{1,\Ran}\otimes  \CC^\dr_{2,\Ran})^{\otimes 2}$$
in $(\CC^\dr_\Ran)^{\otimes 2}$ is canonically isomorphic to $\sR_{\CC,\Ran}$.

\medskip

This immediately implies the assertion of \propref{p:prod projectors} in this context. A similar
observation holds also for $\CC^{\on{Betti}}_{\Ran}$.

\end{rem} 

\sssec{}

By \secref{sss:R as colim}, the object $\sR_\CZ\in \CC_{\Ran}\otimes  \QCoh(\CZ)$ is the colimit 
\begin{equation}  \label{e:R again}
\underset{(I\overset{\psi}\to J)}{\on{colim}}\, 
\on{ins}_{I\to J}\left(\left(\on{Id} \otimes (\sF^J\circ \on{mult}^\psi_{\CC})\right)(\sR_\CC^{\boxtimes I})\right),
\end{equation} 

\medskip

Similarly, the image of $\sR_{\CZ_1}\otimes \sR_{\CZ_2}$ in $\CC_{\Ran}\otimes  \QCoh(\CZ)$ 
is the colimit
\begin{equation}  \label{e:R square}
\underset{(I_1\overset{\psi_1}\to J_1),(I_2\overset{\psi_2}\to J_2)}{\on{colim}}\,
\on{ins}_{I_1\sqcup I_2\to J_1\sqcup J_2} \left(\left(\on{Id} \otimes (\sF^{J_1\sqcup J_2}\circ \on{mult}^{\psi_1\sqcup \psi_2}_{\CC})\right)
(\sR_{\CC_1}^{\boxtimes I_1}\otimes \sR_{\CC_2}^{\boxtimes I_2})\right),
\end{equation} 
where we regard $\sR_{\CC_1}^{\boxtimes I_1}\otimes \sR_{\CC_2}^{\boxtimes I_2}$ as an object of $\CC^{\otimes (I_1\sqcup I_2)}$. 

\medskip

The map from \eqref{e:R square} to \eqref{e:R again} constructed in \secref{sss:R R} is given by the functor
$$\on{TwArr}(\fSet)\times \on{TwArr}(\fSet)\to \on{TwArr}(\fSet), \quad (I_1\to J_1) \times (I_2\to J_2) \mapsto
(I_1\sqcup I_2\to J_1\sqcup J_2)$$
and the maps
$$\sR_{\CC_1}^{\boxtimes I_1}\otimes \sR_{\CC_2}^{\boxtimes I_2}\to 
\sR_{\CC}^{\boxtimes I_1}\otimes \sR_{\CC}^{\boxtimes I_2} \simeq \sR_{\CC}^{\boxtimes (I_1\sqcup I_2)}.$$

We will now construct an inverse map.

\sssec{}

Consider the object 
\begin{equation}  \label{e:R doubled}
\underset{(I\overset{\psi}\to J)\in \on{TwArr}(\fSet)}{\on{colim}}\,
\on{ins}_{I\sqcup I\overset{\psi,\psi}\to J}\left(\left(\on{Id}  \otimes (\sF^J\circ \on{mult}^{\psi,\psi}_{\CC})\right)
(\sR_{\CC_1}^{\boxtimes I}\otimes \sR_{\CC_2}^{\boxtimes I})\right),
\end{equation} 
where we regard $\sR_{\CC_1}^{\boxtimes I}\otimes \sR_{\CC_2}^{\boxtimes I}$ as an object of $\CC^{\otimes (I\sqcup I)}$. 

\medskip

The maps in $\on{TwArr}(\fSet)$ given by the diagrams 
$$
\CD
I \sqcup I @>{\psi,\psi}>> J \\
@V{\on{id},\on{id}}VV @AA{\on{id}}A \\
I @>{\psi}>> J
\endCD
$$
define an isomorphism from \eqref{e:R doubled} to \eqref{e:R again}.

\medskip

We now define a map from \eqref{e:R doubled} to \eqref{e:R square} to be given by mapping
\begin{multline*}
\on{ins}_{I\sqcup I\overset{\psi,\psi}\to J}\left(\left(\on{Id} \otimes (\sF^J\circ \on{mult}^{\psi,\psi}_{\CC})\right)
(\sR_{\CC_1}^{\boxtimes I}\otimes \sR_{\CC_2}^{\boxtimes I}) \right) \to \\
\to \on{ins}_{I\sqcup I\overset{\psi\sqcup \psi}\to J\sqcup J}\left(\left(\on{Id}\otimes (\sF^{J\sqcup J}\circ \on{mult}^{\psi\sqcup \psi}_{\CC})\right)
(\sR_{\CC_1}^{\boxtimes I}\otimes \sR_{\CC_2}^{\boxtimes I})\right)
\end{multline*} 
using the diagram
$$
\CD
I \sqcup I @>{\psi,\psi}>> J  \\
@V{\on{id} \sqcup \on{id}}VV @AA{\on{id},\on{id}}A \\
I \sqcup I @>{\psi\sqcup \psi}>> J \sqcup J 
\endCD
$$
and the natural transformation
$$(\on{Id}\otimes (\Delta_{X^J})_*)\circ \sF^J\circ \on{mult}_{\CC^{\otimes J}}\to \sF^{J\sqcup J}.$$

\sssec{}

It is a straightforward verification that the two maps 
$$\text{\eqref{e:R square}} \leftrightarrow \text{\eqref{e:R again}},$$
constructed above, are mutually inverse.

\qed[\propref{p:prod projectors}]

\ssec{The spectral (sub)category} \label{ss:spectral category}

In this subsection we specialize to the case $\CC=\Rep(\sG)$. 

\sssec{} \label{sss:spectral subcateg}

Let $\bM$ be a module category over $\Rep(\sG)_{\Ran}$. Denote 
$$\bM^{\on{spec}}:=\on{Funct}_{\Rep(\sG)_\Ran}(\QCoh(\LocSys^{\on{restr}}_\sG(X)),\bM).$$

Let $\iota_\bM$ denote the forgetful functor 
$$\bM^{\on{spec}}=\on{Funct}_{\Rep(\sG)_\Ran}(\QCoh(\LocSys^{\on{restr}}_\sG(X)),\bM)\overset{-\circ \Loc}\longrightarrow
\on{Funct}_{\Rep(\sG)_\Ran}(\Rep(\sG)_\Ran,\bM)\simeq \bM.$$

\sssec{} \label{sss:spectral subcateg LocSys}

Note that if $\bM$ is such that the action of $\Rep(\sG)_{\Ran}$ on it factors through an action of 
$\QCoh(\LocSys^{\on{restr}}_\sG(X))$, then the functor $\iota_\bM$ is an equivalence. 

\medskip

Indeed, this follows from \corref{c:LocSys over ran}(b). 

\sssec{}

We now claim:

\begin{prop} \label{p:iota abs is ff}
Assume that $\bM$ is dualizable as a DG category. Then the functor
$$\iota_\bM:\bM^{\on{spec}}\to \bM$$
is fully faithful.
\end{prop}

\begin{proof}

The functor $\Loc$ as the composition
$$\Rep(\sG)_\Ran \to \Rep(\sG)^{\otimes X\on{-lisse}}= \coHom(\Rep(\sG),\qLisse(X))\to \QCoh(\LocSys^{\on{restr}}_\sG(X)),$$
where the last arrow is the map \eqref{e:Hom abs}, which is an equivalence by \thmref{t:spectral}.

\medskip

Hence, the functor
$$\on{Funct}_{\Rep(\sG)_\Ran}(\QCoh(\LocSys^{\on{restr}}_\sG(X)),\bM)\to
\on{Funct}_{\Rep(\sG)_\Ran}(\coHom(\Rep(\sG),\qLisse(X)),\bM)$$
is an equivalence. 

\medskip

Now, by \corref{c:Ran to lisse prel}, the functor
$$\on{Funct}_{\Rep(\sG)_\Ran}(\Rep(\sG)^{\otimes X\on{-lisse}},\bM)\to
\on{Funct}_{\Rep(\sG)_\Ran}(\Rep(\sG)_\Ran,\bM)\simeq \bM$$
is fully faithful, provided that $\bM$ is dualizable. 

\end{proof}

\begin{rem}
We can view \propref{p:iota abs is ff} as saying that, if $\bM$ is dualizable, $\bM^{\on{spec}}$
is the maximal full subcategory of $\bM$ on which the action of $\Rep(\sG)_\Ran$ factors via the
functor
$$\Loc:\Rep(\sG)_\Ran\to \QCoh(\LocSys^{\on{restr}}_\sG(X)).$$

The superscript ``spec" stands for ``spectral decomposition with respect to 
$\LocSys^{\on{restr}}_\sG(X)$".
\end{rem}

\sssec{} \label{sss:Hecke as base change new}

Let $\CZ$ be a prestack equipped with a map $f:\CZ\to \LocSys^{\on{restr}}_\sG(X)$. 
The map $f$ gives rise to a functor $\sF$ as in \eqref{e:eigencond abs}. 

\medskip

Following \secref{sss:Hecke abs}, we can consider the category 
$$\on{Hecke}(\CZ,\bM).$$

\sssec{}

Take $\CZ=\LocSys^{\on{restr}}_\sG(X)$ with $f$ being the identity map. Note that the 
resulting functor $\wt\sF$ identifies with the functor $\Loc$. 

\medskip

Consider the functor
\begin{equation} \label{e:iota abs Hecke}
\on{Hecke}(\LocSys^{\on{restr}}_\sG(X),\bM) \overset{\oblv_{\on{Hecke}},\LocSys^{\on{restr}}_\sG(X)}
\longrightarrow 
\bM \otimes \QCoh(\LocSys^{\on{restr}}_\sG(X)) \overset{\on{Id}_\bM\otimes \Gamma_!}\to \bM,
\end{equation} 
where 
$$\Gamma_!:\QCoh(\LocSys^{\on{restr}}_\sG(X))\to \Vect_\sfe$$
is as in \secref{ss:!-sect semi-pass}.  

\sssec{}

We claim:

\begin{prop} \label{p:M spec}
There exists a canonical equivalence 
$$\on{Hecke}(\LocSys^{\on{restr}}_\sG(X),\bM)  \simeq \bM^{\on{spec}},$$
under which the functor \eqref{e:iota abs Hecke} identifies with the functor
$\iota_\bM$.
\end{prop}

\begin{proof}

We interpret $\on{Hecke}(\LocSys^{\on{restr}}_\sG(X),\bM)$ as
$$\on{Funct}_{\Rep(\sG)_\Ran\otimes \QCoh(\LocSys^{\on{restr}}_\sG(X))}
(\QCoh(\LocSys^{\on{restr}}_\sG(X)),\bM\otimes \QCoh(\LocSys^{\on{restr}}_\sG(X))),$$
see \eqref{e:Hecke category alt}.

\medskip

Now the assertion follows from the fact that we have a canonical identification 
$$\QCoh(\LocSys^{\on{restr}}_\sG(X))^\vee \simeq \QCoh(\LocSys^{\on{restr}}_\sG(X))$$
as $\QCoh(\LocSys^{\on{restr}}_\sG(X))$-modules, with the counit given by
$$\QCoh(\LocSys^{\on{restr}}_\sG(X))\otimes \QCoh(\LocSys^{\on{restr}}_\sG(X))
\overset{\on{mult}}\longrightarrow \QCoh(\LocSys^{\on{restr}}_\sG(X)) \overset{\Gamma_!}\to\Vect_\sfe,$$
see \corref{c:LocSys semi-rigid} and \lemref{l:semi-rigid self dual counit}. 

\end{proof}

\ssec{Beilinson's spectral projector--the universal case}

We retain the setting of \secref{ss:spectral category}. 

\sssec{} \label{sss:object R}

Consider the object
$$\sR_{\LocSys^{\on{restr}}_\sG(X)}\in \Rep(\sG)_\Ran \otimes \QCoh(\LocSys^{\on{restr}}_\sG(X)),$$
see \secref{sss:R object new}. Denote
$$\sR:=(\on{Id}_{\Rep(\sG)_\Ran} \otimes \Gamma_!)(\sR_{\LocSys^{\on{restr}}_\sG(X)})\in \Rep(\sG)_\Ran.$$

Recall that $\sR_{\LocSys^{\on{restr}}_\sG(X)}$ carries a natural structure of commutative algebra
(see \secref{sss:R object new}). Recall also that the functor $\Gamma_!(\LocSys^{\on{restr}}_\sG(X),-)$
has a natural right-lax symmetric monoidal structure (see \secref{sss:!-sect semi-pass}). 

\medskip

Hence, the
object $\sR$ has a natural structure of commutative algebra in $\Rep(\sG)_\Ran$.

\sssec{}

From \thmref{t:diagonal}, we obtain:

\begin{cor} \label{c:R LocSys}
There exists a canonical isomorphism of commutative algebras
$$\Loc(\sR)\simeq \CO_{\LocSys^{\on{restr}}_\sG(X)}.$$
\end{cor}

\begin{proof}

By definition,
$$\Loc(\sR) \simeq
(\on{Id}_{\QCoh(\LocSys^{\on{restr}}_\sG(X))} \otimes \Gamma_!)\circ (\Loc\otimes \Loc)(\sR_{\Rep(\sG),\Ran}).$$

However, by \thmref{t:diagonal}, we have
$$(\Loc\otimes \Loc)(\sR_{\Rep(\sG),\Ran})\simeq 
(\Delta_{\LocSys^{\on{restr}}_\sG(X)})_*(\CO_{\LocSys^{\on{restr}}_\sG(X)}).$$

Finally, we have
$$(\on{Id}_{\QCoh(\LocSys^{\on{restr}}_\sG(X))} \otimes \Gamma_!) \circ 
(\Delta_{\LocSys^{\on{restr}}_\sG(X)})_*\simeq \on{Id}_{\QCoh(\LocSys^{\on{restr}}_\sG(X))},$$
as commutative algebras 
(indeed, this is a feature of any semi-rigid category, see Lemmas \ref{l:semi-rigid self dual} and \ref{l:semi-rigid self dual counit}). 

\end{proof}

\sssec{} \label{p:functor P}

Let $\bM$ be a module category over $\Rep(\sG)_\Ran$. 

\medskip

Let us denote by $\sP$ the functor $\sP^{\on{enh}}_{\LocSys^{\on{restr}}_\sG(X)}$
(see \secref{sss:P enh first}), which we now view as a functor 
$$\sP:\bM \to  \bM^{\on{spec}},$$
thanks to \propref{p:M spec}. 

\sssec{}

Furthermore, from \propref{p:M spec} (combined with the observation in \secref{sss:R object new}), we obtain: 

\begin{cor} \label{c:the projector univ}
The endofunctor $\iota_\bM\circ \sP$ of $\bM$ is given by the action of the object $\sR\in \Rep(\sG)_\Ran$.
\end{cor}

Finally, we claim:

\begin{prop} \label{p:the projector univ}
The endofunctor $\sP\circ \iota_\bM$ is canonically isomorphic to the identity.
\end{prop}

\begin{proof}

It is enough to consider the universal case, i.e., $\bM=\Rep(\sG)_\Ran$. So, we can assume
that $\bM$ is dualizable. In this case, by \propref{p:iota abs is ff}, we can view $\bM^{\on{spec}}$
as a full subcategory of $\bM$, and by \corref{c:the projector univ}, the endofunctor 
$\sP\circ \iota_\bM$ is induced by the action of $\sR\in \Rep(\sG)_\Ran$ on $\bM$.

\medskip

However, for $\CV\in \Rep(\sG)_\Ran$, its action on 
$$\bM^{\on{spec}}:=\on{Funct}_{\Rep(\sG)_\Ran}(\QCoh(\LocSys^{\on{restr}}_\sG(X)),\bM)$$
is given by the action of $\Loc(\sR)$ on $\QCoh(\LocSys^{\on{restr}}_\sG(X))$. Now, the required
assertion follows from \corref{c:R LocSys}.

\end{proof}

\begin{rem}

We can view the combination of \corref{c:the projector univ} and \propref{p:the projector univ} as saying that, if $\bM$ is dualizable,
the action of $\sR$ on $\bM$ acts as projector onto the full subcategory 
$$\bM^{\on{spec}}\subset \bM.$$

Thus, we can think of $\sR$ as a ``universal spectral projector". 

\end{rem}

\sssec{} \label{sss:R expl}

Let us write the object $\sR$ explicitly as a colimit. 

\medskip

Namely, let $\CZ$ be a connected component of $\LocSys^{\on{restr}}_\sG(X)$, and let 
$f_n:Z_n\to \CZ$ be as in Sects. \ref{sss:almost lift action a}-\ref{sss:almost lift action b}. 

\medskip

Using \eqref{e:Gamma ! expl}, we obtain: 
\begin{equation} \label{e:formula for R univ}
\sR \simeq \underset{\CZ}\oplus\, \underset{n}{\on{colim}}\, (\on{Id}_{\Rep(\sG)_\Ran}\otimes \Gamma(Z_n,\ell_{Z_n}\otimes -))(\sR_{Z_n}),
\end{equation}
where:

\begin{itemize}

\item $\sR_{Z_n}\in \Rep(\sG)_\Ran\otimes \QCoh(Z_n)$ is given by formula \eqref{e:formula for R};

\item $\ell_{Z_n}$ is the line bundle $f_n^!(\CO_{\LocSys^{\on{restr}}_\sG(X)})$ on $Z_n$ (so, non-canonically, 
$\ell_{Z_n}\simeq \CO_{Z_n}[-m_n]$, where the integer $m_n$ only depends on $\CZ$, see \eqref{e:upper ! regular}).

\end{itemize}

\ssec{Beilinson's spectral projector--the general case}

\sssec{} \label{sss:Hecke as base change old}

Let $f:\CZ\to \LocSys^{\on{restr}}_\sG(X)$ be as in \secref{sss:Hecke as base change new}, and let
$\bM$ be a module category over $\Rep(\sG)_\Ran$. 

\medskip

Note that since $\Rep(\sG)_\Ran$ is rigid, we can rewrite 
$$\on{Hecke}(\CZ,\bM)\simeq \bM\underset{\Rep(\sG)_\Ran}\otimes \QCoh(\CZ)$$
(see \cite[Chapter 1, Proposition 9.4.4]{GR1} or \propref{p:HH semi-rigid}).

\medskip

In particular, 
$$\bM^{\on{spec}}\simeq \bM \underset{\Rep(\sG)_\Ran}\otimes \QCoh(\LocSys^{\on{restr}}_\sG(X)).$$

Combining, we obtain an equivalence
\begin{equation} \label{e:Hecke as ten}
\on{Hecke}(\CZ,\bM) \simeq \bM^{\on{spec}}\underset{\QCoh(\LocSys^{\on{restr}}_\sG(X))}\otimes \QCoh(\CZ)
\end{equation}

\sssec{} \label{sss:oblv Hecke via spec}

Unwinding the definitions, we obtain that, in terms of the equivalence \eqref{e:Hecke as ten}, 
the forgetful functor 
$$\oblv_{\on{Hecke},\CZ}:\on{Hecke}(\CZ,\bM)\to \bM\otimes \QCoh(\CZ)$$ 
identifies canonically with
\begin{multline} \label{e:oblv Hecke via spec}
\bM^{\on{spec}} \underset{\QCoh(\LocSys^{\on{restr}}_\sG(X))}\otimes \QCoh(\CZ)
\simeq \\
\simeq (\bM^{\on{spec}} \otimes \QCoh(\CZ)) \underset{\QCoh(\LocSys^{\on{restr}}_\sG(X))\otimes \QCoh(\LocSys^{\on{restr}}_\sG(X))}\otimes
\QCoh(\LocSys^{\on{restr}}_\sG(X)) \overset{\on{Id}\otimes (\Delta_{\LocSys^{\on{restr}}_\sG(X)})_*}\longrightarrow \\
(\bM^{\on{spec}} \otimes \QCoh(\CZ)) \underset{\QCoh(\LocSys^{\on{restr}}_\sG(X))\otimes \QCoh(\LocSys^{\on{restr}}_\sG(X))}\otimes
\left(\QCoh(\LocSys^{\on{restr}}_\sG(X)\otimes \LocSys^{\on{restr}}_\sG(X))\right)\simeq \\
\simeq \bM^{\on{spec}} \otimes \QCoh(\CZ)\overset{\iota_\bM\otimes \on{Id}}\longrightarrow \bM\otimes \QCoh(\CZ).
\end{multline}

Similarly, by \secref{sss:change Z abs}, the functor 
$$\sP^{\on{enh}}_\CZ:\bM\to \on{Hecke}(\CZ,\bM)$$
identifies with
\begin{multline} \label{e:P enh via spec}
\bM \overset{\sP}\to \bM^{\on{spec}} \simeq \bM^{\on{spec}}  \underset{\QCoh(\LocSys^{\on{restr}}_\sG(X))}\otimes \QCoh(\LocSys^{\on{restr}}_\sG(X))
\overset{\on{Id}\otimes f^*}\longrightarrow \\
\to \bM^{\on{spec}}  \underset{\QCoh(\LocSys^{\on{restr}}_\sG(X))}\otimes \QCoh(\CZ).
\end{multline}

%

\sssec{} \label{sss:Hecke as base change cont}

Assume now that $\CO_\CZ$ is compact as an object of $\QCoh(\CZ)$. Factoring the morphism $f$ as
$$\CZ \simeq \LocSys^{\on{restr}}_\sG(X)\underset{\LocSys^{\on{restr}}_\sG(X)\times \LocSys^{\on{restr}}_\sG(X)}\times
(\LocSys^{\on{restr}}_\sG(X)\times \CZ)\to $$
$$\to \LocSys^{\on{restr}}_\sG(X)\times \CZ \to \LocSys^{\on{restr}}_\sG(X)$$
and using the fact that $\Delta_{\LocSys^{\on{restr}}_\sG(X)}$ is an affine morphism, 
we obtain that in this case the functor $f_*$ is continuous and compatible with $\QCoh(\LocSys^{\on{restr}}_\sG(X))$-module structures.

\medskip

Hence, in this case, using the expression for $\oblv_{\on{Hecke},\CZ}$ in \secref{sss:oblv Hecke via spec}, we obtain that  the functor 
$$\oblv_{\on{Hecke}}\simeq (\on{Id}_\bM\otimes \Gamma(\CZ,-))\circ \oblv_{\on{Hecke},\CZ}$$
identifies with 
\begin{multline} \label{e:oblv Hecke abs}
\bM^{\on{spec}} \underset{\QCoh(\LocSys^{\on{restr}}_\sG(X))}\otimes \QCoh(\CZ)
\overset{\on{Id}\otimes f_*}\longrightarrow \\
\to \bM^{\on{spec}} \underset{\QCoh(\LocSys^{\on{restr}}_\sG(X))}\otimes \QCoh(\LocSys^{\on{restr}}_\sG(X))
\simeq \bM^{\on{spec}} \overset{\iota_\bM}\to \bM.
\end{multline}

In particular, we obtain:

\begin{cor}
If $\CO_\CZ$ is compact, then the functor \eqref{e:oblv Hecke abs} admits a left adjoint, explicitly 
given by \eqref{e:P enh via spec}.
\end{cor}

\begin{rem}

The material of this subsection applies equally well when instead of $\CC_\Ran$ and $\LocSys^{\on{restr}}_\sG(X)$ we
consider $\CC^\dr_\Ran$ and $\LocSys^{\dr}_\sG(X)$ or
$\CC^{\on{Betti}}_\Ran$ and $\LocSys^{\on{Betti}}_\sG(X)$.

\end{rem}

\ssec{A version with parameters} \label{ss:Ran param}

%

When we work in a constructible sheaf-theoretic context (as opposed to $\Dmod(-)$
or $\Shv^{\on{all}}(-)$), the formalism of the (symmetric) 
monoidal category $\Rep(\sG)_\Ran$ is not sufficient to encode the pattern of the Hecke
action, to be studied in Part III of the paper (there $\sG$ will be the Langlands dual $\cG$
of ``our" group $G$). 

\medskip

The reason for this is that the Hecke functors map $\Shv(\Bun_G)$ to
$\Shv(\Bun_G\times X^I)$, which contains, but is not equivalent to $\Shv(\Bun_G)\otimes \Shv(X^I)$.

\medskip

In order to account for this, we will need to introduce a version of $\Rep(\cG)_\Ran$, where we 
allow an additional scheme as a parameter. We will continue to work in an abstract setting,
when instead of $\Rep(\cG)_\Ran$ we have an arbitrary rigid symmetric monoidal category $\CC$. 

\sssec{}

Let $X$ and $\CC$ be as in \secref{ss:Ran rigid}.

\medskip

Let $Y$ be a scheme over $k$. We introduce the category $\CC_{\Ran\times Y}$
by the same colimit procedure as in the case of $\CC_\Ran$, with the difference that
instead of $\Shv(X^J)$ we now use $$\Shv(X^J\times Y).$$

\medskip

We endow $\CC_{\Ran\times Y}$ with a symmetric monoidal structure using the operation
of disjoint of finite sets, where we now use the functors
$$\Shv(X^{J_1}\times Y)\otimes \Shv(X^{J_2}\times Y)\overset{\boxtimes}\to
\Shv((X^{J_1}\times Y)\times (X^{J_2}\times Y))\overset{!\text{-pullback}}\longrightarrow
\Shv(X^{J_1\sqcup J_2}\times Y).$$

\medskip

Tensoring by !-pullbacks of objects of $\Shv(Y)$, we obtain a (unital) symmetric monoidal functor
\begin{equation} \label{e:Y to Ran Y}
\Shv(Y)\to \CC_{\Ran\times Y}.
\end{equation}

\sssec{}  \label{sss:change Y}

Let $f:Y_1\to Y_2$ be a map of schemes. Then !-pullback along $f$ defines a (unital) symmetric monoidal
functor
$$\CC_{\Ran\times Y_2}\to \CC_{\Ran\times Y_1}.$$

\sssec{}

In particular for any $Y$, we have a (unital) symmetric monoidal functor
\begin{equation} \label{e:Ran to Ran Y}
\CC_\Ran\simeq \CC_{\Ran\times \on{pt}}\to \CC_{\Ran\times Y}.
\end{equation} 

Combining with \eqref{e:Y to Ran Y}, we obtain a symmetric monoidal functor
\begin{equation} \label{e:ten by Y}
\CC_\Ran \otimes \Shv(Y)\to \CC_{\Ran\times Y}.
\end{equation} 

Since the individual functors
$$\Shv(X^J)\otimes \Shv(Y)\to \Shv(X^J\times Y)$$
are fully faithful and the category $\Shv(Y)$ is dualizable, it follows from \eqref{e:Ran as limit} (and a similar presentation for $\CC_{\Ran\times Y}$)
that the functor \eqref{e:ten by Y} is fully faithful. 

\begin{rem}
Similar definitions apply when instead of $\CC_\Ran$ we use $\CC_\Ran^\dr$. However, in this case, the 
corresponding functor
$$\CC^\dr_\Ran \otimes \Shv(Y)\to \CC^\dr_{\Ran\times Y}$$
is an equivalence. So in this case, there is no point of introducing $\CC^\dr_{\Ran\times Y}$ as a separate 
entity.

\medskip

Similar definitions also apply to $\CC_\Ran^{\on{Betti}}$ (but the !-pullbacks replaced by *-pullbacks).
Here again, the corresponding functor
$$\CC^{\on{Betti}}_\Ran \otimes \Shv(Y)\to \CC^{\on{Betti}}_{\Ran\times Y}$$
is an equivalence. 

\end{rem}

\sssec{} \label{sss:F Y}

Let us be given a symmetric monoidal functor $\sF$ as in \eqref{e:eigencond abs}. From $\sF$ we produce a 
symmetric monoidal functor
$$\wt\sF_Y:\CC_{\Ran\times Y}\to \QCoh(\CZ)\otimes \Shv(Y).$$

Namely, $\wt\sF_Y$ sends an object
$$\on{ins}_{I\overset{\psi}\to J}(V \otimes \CF_Y), \quad V\in \CC^{\otimes I},\, \CF_Y\in \Shv(X^J\times Y)$$
to 
$$\left(\on{Id}_{\QCoh(\CZ)}\otimes (p_Y)_*\right)(\sF^J(\on{mult}^\psi_\CC(V))\otimes \CF_Y),$$
where $p_Y$ denotes the projection
$$X^J\times Y\to Y.$$

\sssec{} 

Let $\bM$ be a module category over $\CC_{\Ran\times Y}$. Let us regard
$$\bM\otimes \QCoh(\CZ)\simeq
\bM \underset{\Shv(Y)}\otimes (\QCoh(\CZ)\otimes \Shv(Y))$$
as a module category over 
$$\CC_{\Ran\times Y}\otimes \QCoh(\CZ)\simeq 
\CC_{\Ran\times Y}\underset{\Shv(Y)}\otimes (\QCoh(\CZ)\otimes \Shv(Y)).$$

\medskip

We can also view it as a module over $\CC_{\Ran\times Y}\underset{\Shv(Y)}\otimes  \CC_{\Ran\times Y}$
via $\wt\sF_Y$.

\sssec{} \label{sss:Hecke param}

We define
$$\on{Hecke}_Y(\CZ,\bM)$$ as the category of functors 
$$\CC_{\Ran\times Y}\to \bM\otimes \QCoh(\CZ).$$
of modules categories over $\CC_{\Ran\times Y}\underset{\Shv(Y)}\otimes  \CC_{\Ran\times Y}$.

\medskip

Equivalently, we can view $\on{Hecke}_Y(\CZ,\bM)$ as the category of functors 
\begin{equation} \label{e:Q Y M}
\QCoh(\CZ)\otimes \Shv(Y) \to 
\bM\otimes \QCoh(\CZ)
\end{equation} 
of modules categories over 
$\CC_{\Ran\times Y}\otimes \QCoh(\CZ)$, 
where $\QCoh(\CZ)\otimes \Shv(Y)$ in the left-hand side of \eqref{e:Q Y M}
is regarded as a module over $\CC_{\Ran\times Y}$ via $\wt\sF_Y$. 

\medskip

We have a naturally defined forgetful functor
\begin{equation} \label{e:pre-oblv Y}
\oblv_{\on{Hecke}_Y,\CZ}:\on{Hecke}_Y(\CZ,\bM)\to \bM\otimes \QCoh(\CZ).
\end{equation} 

\sssec{}

Note that given $\bM$ as above, we can regard it as a module category over $\CC_{\Ran}$ via the symmetric
monoidal functor 
\eqref{e:Ran to Ran Y}. We have a naturally defined forgetful functor
\begin{equation} \label{e:get rid of Y}
\on{Hecke}_Y(\CZ,\bM)\to\on{Hecke}(\CZ,\bM)
\end{equation} 
that makes the diagram 
$$
\CD
\on{Hecke}_Y(\CZ,\bM) @>{\oblv_{\on{Hecke}_Y,\CZ}}>> \bM\otimes \QCoh(\CZ) \\
@V{\text{\eqref{e:get rid of Y}}}VV @VV{\on{Id}}V \\
\on{Hecke}(\CZ,\bM) @>{\oblv_{\on{Hecke},\CZ}}>> \bM\otimes \QCoh(\CZ) 
\endCD
$$
commute.

\medskip

We will prove:

\begin{thm} \label{t:get rid of Y}
The functor \eqref{e:get rid of Y} is an equivalence.
\end{thm}

\sssec{}

Let us denote by 
$$\sR_{\CZ,Y}\in  \CC_{\Ran\times Y}\otimes \QCoh(\CZ)$$
the image of $\sR_\CZ\in \CC_\Ran\otimes \QCoh(\CZ)$ along 
$$\CC_\Ran\otimes \QCoh(\CZ) \to \CC_{\Ran\times Y}\otimes \QCoh(\CZ).$$

\medskip

From \thmref{t:get rid of Y} we obtain: 

\begin{cor} \label{c:get rid of Y}
The functor $\oblv_{\on{Hecke}_Y,\CZ}$ is monadic, and the resulting monad
on $\bM\otimes \QCoh(\CZ)$ is given by the action of $\sR_{\CZ,Y}$.
\end{cor}

\ssec{Proof of \thmref{t:get rid of Y}}

\sssec{}

Consider the functor 
\begin{multline} \label{e:left adjoint Y}
\CC_{\Ran\times Y}\underset{\Shv(Y)}\otimes (\QCoh(\CZ)\otimes \Shv(Y))
\overset{\wt\sF_Y\otimes \on{Id}}\longrightarrow \\
\to (\QCoh(\CZ)\otimes \Shv(Y)) \underset{\Shv(Y)}\otimes (\QCoh(\CZ)\otimes \Shv(Y)) \overset{\on{mult}}\to
\QCoh(\CZ)\otimes \Shv(Y).
\end{multline}

It is enough to show that 
$$\one_{\QCoh(\CZ)\otimes \Shv(Y)}\mapsto \sR_{\CZ,Y}$$
extends to a map of 
$\CC_{\Ran\times Y}\underset{\Shv(Y)}\otimes (\QCoh(\CZ)\otimes \Shv(Y))$-module categories 
\begin{equation} \label{e:right adjoint Y}
\QCoh(\CZ)\otimes \Shv(Y)\to \CC_{\Ran\times Y}\otimes \QCoh(\CZ)\simeq 
\CC_{\Ran\times Y}\underset{\Shv(Y)}\otimes (\QCoh(\CZ)\otimes \Shv(Y)),
\end{equation}
which is the right adjoint of \eqref{e:left adjoint Y}. 

\sssec{}

Let us denote the functor \eqref{e:ten by Y} by $\Phi$. 

\medskip

Note that the individual functors
$$\boxtimes:\Shv(X^J)\otimes \Shv(Y)\to \Shv(X^J\times Y)$$
preserve compactness, and hence admit continuous right adjoints,
to be denoted $\boxtimes^R$. Since $\boxtimes$ commutes with
Verdier duality, the functor $\boxtimes^R$ also identifies with the
dual of $\boxtimes$. The latter observation implies that for
$\phi:J_2\to J_1$, the diagram
$$
\CD
\Shv(X^{J_1}\times Y) @>{(\Delta_\phi\times \on{id})_*}>> \Shv(X^{J_2}\times Y) \\
@V{\boxtimes^R}VV @VV{\boxtimes^R}V \\
\Shv(X^{J_1})\otimes \Shv(Y) @>{(\Delta_\phi)_*\otimes \on{Id}}>> \Shv(X^{J_2})\otimes \Shv(Y),
\endCD
$$
which a priori commutes up to a natural transformation, commutes strictly.

\medskip

This implies that the functors  $\boxtimes^R$ assemble to a functor
$$\Psi:\CC_{\Ran\times Y}\to \CC_\Ran\otimes \Shv(Y),$$
right adjoint to $\Phi$. The unit map $\on{Id}\to \Psi\circ \Phi$ is an isomorphism
since $\Phi$ is fully faithful. 

\medskip

Being a right adjoint to a symmetric monoidal functor, the functor $\Psi$ carries
a canonically defined right-lax symmetric monoidal structure. It is easy to see, however, 
that $\Psi$ is strictly linear with respect to $\CC_\Ran\otimes \Shv(Y)$.

\sssec{}

Let us observe that the functor $\wt\sF_Y$ identifies 
with
$$\CC_{\Ran\times Y} \overset{\Psi}\to \CC_\Ran\otimes \Shv(Y)
\overset{\wt\sF\otimes \on{Id}}\longrightarrow \QCoh(\CZ)\otimes \Shv(Y),$$
as a right-lax symmetric monoidal functor (however, the functor $\wt\sF_Y$ itself 
is strict). 

\sssec{}

By \eqref{e:Hecke univ}, we have a tensor-compatible system of isomorphisms
$$\CV\star \sR_\CZ \simeq  \sR_\CZ\otimes \wt\sF(\CV), \quad \CV\in \CC_\Ran.$$

From here, the $(\Phi,\Psi)$-adjunction gives rise to a tensor-compatible system of \emph{morphisms}
\begin{equation} \label{e:Hecke univ Y}
\CV_Y \star \sR_{\CZ,Y} \leftarrow \sR_{\CZ,Y} \otimes \wt\sF_Y(\CV_Y), \quad \CV_Y\in \CC_{\Ran\times Y}.
\end{equation} 

The key observation is the following: 

\begin{prop} \label{p:tightness}
The maps \eqref{e:Hecke univ Y} are isomorphisms.
\end{prop}

The proof will be given in Sects. \ref{ss:tightness1} and \ref{ss:tightness2}
(we will give two proofs, each in the corresponding section). 

\medskip

Let us accept this proposition temporarily and
finish the proof of \thmref{t:get rid of Y}. 

\sssec{} \label{sss:proj Y}

By \propref{p:tightness}, the object $\sR_{\CZ,Y}$, defines a map 
\begin{equation} \label{e:R Z Y}
\QCoh(\CZ)\otimes \Shv(Y)\to \CC_{\Ran\times Y}\underset{\Shv(Y)}\otimes (\QCoh(\CZ)\otimes \Shv(Y))\simeq \CC_{\Ran\times Y}\otimes \QCoh(\CZ)
\end{equation} 
of module categories over $\CC_{\Ran\times Y}\otimes \QCoh(\CZ)$. 

\medskip

We now construct the adjunction datum between \eqref{e:left adjoint Y} and \eqref{e:R Z Y}. 

\medskip

\noindent Unit: Since \eqref{e:left adjoint Y} and \eqref{e:R Z Y} are functors of module categories over 
$\CC_{\Ran\times Y}\otimes \QCoh(\CZ)$, it suffices to specify the value of the unit on the object
$\one_{\CC_{\Ran\times Y}\otimes \QCoh(\CZ)}\in \CC_{\Ran\times Y}\otimes \QCoh(\CZ)$.

\medskip

It is obtained by applying the functor $\Phi$ to the unit of the adjunction
\begin{equation} \label{e:R Z adj}
\on{mult}_{\QCoh(\CZ)}\circ (\wt\sF\otimes \on{Id}):\CC_\Ran\otimes \QCoh(\CZ) 
\rightleftarrows \QCoh(\CZ):\sR_\CZ
\end{equation}
on $\one_{\CC_{\Ran}\otimes \QCoh(\CZ)}\in \CC_{\Ran}\otimes \QCoh(\CZ)$.

\medskip

\noindent Counit: it suffices to specify the value of the counit on $\one_{\QCoh(\CZ)\otimes \Shv(Y)}\in \QCoh(\CZ)\otimes \Shv(Y)$. 
I.e., we have to specify a map
\begin{equation} \label{e:counit Y Z}
\on{mult}_{\QCoh(\CZ)\otimes \Shv(Y)}\circ (\wt\sF_Y\otimes \on{Id}_{\QCoh(\CZ)\otimes \Shv(Y)})(\sR_{\CZ,Y})\to 
\one_{\QCoh(\CZ)\otimes \Shv(Y)}.
\end{equation}

Note that 
$$\wt\sF_Y\circ \Phi \simeq (\wt\sF\otimes \on{Id}_{\Shv(Y)}) \circ \Psi \circ \Phi\simeq \wt\sF\otimes \on{Id}_{\Shv(Y)}.$$

Hence, the left-hand side in \eqref{e:counit Y Z} identifies with 
$$(\wt\sF\otimes \on{Id}_{\QCoh(\CZ)})(\sR_\CZ)\otimes \one_{\Shv(Y)}.$$

The required map in \eqref{e:counit Y Z} is obtained from the counit of the adjunction \eqref{e:R Z adj}
by tensoring with $\one_{\Shv(Y)}=\omega_Y$.

\qed[\thmref{t:get rid of Y}]

\ssec{Direct proof of \propref{p:tightness}} \label{ss:tightness1}

\sssec{}

It is sufficient to prove that \eqref{e:Hecke univ Y} is an isomorphism for $\CV_Y$ of the form 
$$\on{ins}_{I_0\to J_0}(V\otimes \CF), \quad V\in \CC^{\otimes I_0},\quad \CF\in \Shv(X^{J_0}\times Y).$$

\sssec{}

Let us denote by 
$$\CC_{\Ran\times \Ran} \text{ and } \CC_{\Ran\times \Ran\times Y}$$ the categories defined as
$$\underset{(I_1\to J_1),(I_2\to J_2)}{\on{colim}}\,
\CC^{\otimes I_1}\otimes \CC^{\otimes I_2}\otimes \Shv(X^{J_1}\times X^{J_2})$$
and 
$$\underset{(I_1\to J_1),(I_2\to J_2)}{\on{colim}}\,
\CC^{\otimes I_1}\otimes \CC^{\otimes I_2}\otimes \Shv(X^{J_1}\times X^{J_2}\times Y),$$
respectively, where 
$$(I_1\to J_1),(I_2\to J_2)\in  \on{TwArr}(\fSet)\times  \on{TwArr}(\fSet).$$

Denote by $\on{ins}_{(I_1\to J_1),(I_2\to J_2)}$ the corresponding tautological functors 
$$\CC^{\otimes I_1}\otimes \CC^{\otimes I_2}\otimes \Shv(X^{J_1}\times X^{J_2})\to\CC_{\Ran\times \Ran}$$
and 
$$\CC^{\otimes I_1}\otimes \CC^{\otimes I_2}\otimes \Shv(X^{J_1}\times X^{J_2}\times Y)\to\CC_{\Ran\times \Ran\times Y},$$
respectively. 

\sssec{}

The operation of disjoint union on finite sets makes $\CC_{\Ran\times \Ran}$ and $\CC_{\Ran\times \Ran\times Y}$
into symmetric monoidal categories.

\medskip

We have naturally defined symmetric monoidal functors
$$\Upsilon:\CC_{\Ran}\otimes \CC_{\Ran}\to \CC_{\Ran\times \Ran}  \text{ and }
\Upsilon_Y:\CC_{\Ran\times Y}\underset{\Shv(Y)}\otimes \CC_{\Ran\times Y}\to \CC_{\Ran\times \Ran\times Y}.$$

We also have a symmetric monoidal functor 
\begin{equation} \label{e:double Ran pullback Y}
\CC_{\Ran\times \Ran} \to \CC_{\Ran\times \Ran\times Y},
\end{equation} 
given by pullback along $Y\to \on{pt}$.

\begin{rem}

Note the difference between $\CC_{\Ran\times \Ran}$ and $\CC_{\Ran}\otimes \CC_{\Ran}$. In the former the 
terms of the colimit have factors $\Shv(X^{J_1}\times X^{J_2})$, and in the latter 
$\Shv(X^{J_1})\otimes \Shv(X^{J_2})$.

\end{rem}

\sssec{}

Let $\sR_{\on{geom},\CC,\Ran}$ be the object of $\CC_{\Ran\times \Ran}$ defined as 
$$\underset{(I\to J)\in \on{TwArr}(\fSet)}{\on{colim}}\, \on{ins}_{(I\to J),(I\to J)}\,
\left((\sR_\CC)^{\boxtimes I}\otimes (\Delta_{X^J})_*(\omega_{X^J})\right).$$

The maps
$$\on{u}_{\Shv(X^J)}\to (\Delta_{X^J})_*(\omega_{X^J})$$ gives rise to a map
\begin{equation} \label{e:R to R geom}
\Upsilon(\sR_{\CC,\Ran})\to \sR_{\on{geom},\CC,\Ran}.
\end{equation} 

Let $\sR_{\on{geom},\CC,\Ran,Y}$ be the object of $\CC_{\Ran\times \Ran\times Y}$ equal to the image 
of $\sR_{\on{geom},\CC,\Ran}$ along \eqref{e:double Ran pullback Y}. 

\sssec{}

Let $\CV$ be an object of $\CC^{\otimes I_0} \otimes \Shv(X^{J_0})$ for some
$(I_0\to J_0)\in \on{TwArr}(\fSet)$. Let $\CV^l$ and $\CV^r$
denote its images in $\CC_{\Ran\times \Ran}$
along
$$\CC^{\otimes I_0} \otimes \Shv(X^{J_0}) 
\overset{\on{ins}_{I_0\to J_0}}\longrightarrow \CC_{\Ran} 
\overset{\on{Id}\otimes \one_{\CC_{\Ran}}} \longrightarrow \CC_{\Ran}\otimes \CC_{\Ran}\to \CC_{\Ran\times \Ran}$$
and
$$\CC^{\otimes I_0} \otimes \Shv(X^{J_0}) 
\overset{\on{ins}_{I_0\to J_0}}\longrightarrow \CC_{\Ran} 
\overset{\one_{\CC_{\Ran}}\otimes \on{Id}} \longrightarrow \CC_{\Ran}\otimes \CC_{\Ran}\to \CC_{\Ran\times \Ran},$$
respectively.

\medskip

The calculation performed in \secref{ss:Hecke abs} shows that we have a canonical isomorphism
$$\CV^l\star \sR_{\on{geom},\CC,\Ran} \simeq \sR_{\on{geom},\CC,\Ran}\star \CV^r$$
in $\CC_{\Ran\times \Ran}$. 

\medskip

Let $\CV_Y$ be an object of $\CC^{\otimes I_0} \otimes \Shv(X^{J_0}\times Y)$, and let
$$\CV_Y^l,\CV_Y^r\in \CC_{\Ran\times \Ran\times Y}$$
be defined in a way similar to the above.

\medskip

Then the same calculation shows that 
we have a canonical isomorphism
\begin{equation} \label{e:Hecke Y geom}
\CV_Y^l\star \sR_{\on{geom},\CC,\Ran,Y} \simeq \sR_{\on{geom},\CC,\Ran,Y}\star \CV_Y^r
\end{equation} 
in $\CC_{\Ran\times \Ran\times Y}$.

\sssec{}

Let $(\CZ,\sF)$ be as in \secref{sss:F Y}. We claim that we have a naturally defined symmetric monoidal functor
$$\wt\sF_{\Ran}:\CC_{\Ran\times \Ran}\to \CC_{\Ran}\otimes \QCoh(\CZ).$$

Namely, the functor $\wt\sF_{\Ran}$ sends an object 
$$\on{ins}_{(I_1\overset{\psi_1}\to J_1),(I_2\overset{\psi_2}\to J_2)}(V_1\otimes V_2\otimes \CF), \quad
V_1\in \CC^{\otimes I_1},\quad V_2\in \CC^{\otimes I_2},\quad \CF\in \Shv(X^{J_1}\times X^{J_2})$$
to 
$$\on{ins}_{I_1\to J_1}(V_1\otimes \CF_1),$$
where $\CF_1$ is the object of $\QCoh(\CZ)\otimes \Shv(X^{J_1})$ equal to
$$(\on{Id}_{\QCoh(\CZ)}\otimes (p_{J_1})_*)(\sF^{J_2}(\on{mult}^{\psi_2}_\CC(V_2))\otimes \CF),$$
where $p_{J_1}$ is the projection $X^{J_1}\times X^{J_2}\to X^{J_1}$. 

\sssec{}

Note that the composition
$$\CC_{\Ran} \otimes \CC_{\Ran} \overset{\Upsilon}\to \CC_{\Ran\times \Ran} \overset{\wt\sF_{\Ran}}\to \CC_{\Ran}\otimes \QCoh(\CZ)$$
identifies with the functor $\on{Id}_{\CC_\Ran}\otimes \wt\sF$, where $\wt\sF$ is an in \eqref{e:functor F tilde}. 

\medskip

We claim:

\begin{lem} \label{l:F on R geom}
The map 
$$\sR_\CZ = (\on{Id}_{\CC_\Ran}\otimes \wt\sF)(\sR_{\CC,\Ran})\simeq 
\wt\sF_{\Ran} \circ \Upsilon(\sR_{\CC,\Ran}) \overset{\text{\eqref{e:R to R geom}}}\longrightarrow \wt\sF_{\Ran}(\sR_{\on{geom},\CC,\Ran})$$
is an isomorphism.
\end{lem}

\begin{proof}

Follows from the isomorphism \eqref{e:proj formula u}. 

\end{proof} 

\sssec{}

By a similar token we construct a map
$$\wt\sF_{\Ran,Y}:\CC_{\Ran\times \Ran\times Y}\to \CC_{\Ran\times Y}\otimes \QCoh(\CZ)$$
such that the composition
$$\CC_{\Ran\times Y} \underset{\Shv(Y)}\otimes \CC_{\Ran\times Y} 
\overset{\Upsilon_Y}\to \CC_{\Ran\times \Ran\times Y} \overset{\wt\sF_{\Ran,Y}}\to \CC_{\Ran\times Y}\otimes \QCoh(\CZ)$$
identifies with the map 
$$\CC_{\Ran\times Y} \underset{\Shv(Y)}\otimes \CC_{\Ran\times Y}  
\overset{\on{Id}\otimes \wt\sF_Y}\longrightarrow \CC_{\Ran\times Y} \underset{\Shv(Y)}\otimes  (\QCoh(\CZ)\otimes \Shv(Y))\simeq
\CC_{\Ran\times Y}  \otimes  \QCoh(\CZ).$$

As in \lemref{l:F on R geom}, we obtain that the resulting map
$$\sR_{\CZ,Y}\to \wt\sF_{\Ran,Y}(\sR_{\on{geom},\CC,\Ran,Y})$$
is an isomorphism.

\sssec{}

Combining the latter isomorphism with the isomorphisms \eqref{e:Hecke Y geom}, we obtain isomorphisms
$$\CV_Y\star \sR_{\CZ,Y} \simeq \sR_{\CZ,Y}  \otimes \wt\sF_Y(\CV_Y).$$

Unwinding the definitions, we obtain that the resulting morphisms 
$$\CV_Y\star \sR_{\CZ,Y} \leftarrow \sR_{\CZ,Y}  \otimes \wt\sF_Y(\CV_Y)$$
are equal to those in \eqref{e:Hecke univ Y}.

\qed[\propref{p:tightness}]

\ssec{An indirect proof of \propref{p:tightness}} \label{ss:tightness2}

We will give a proof that works in the \'etale and constructible de Rham contexts;
the constructible Betti case will follow from the de Rham case by Riemann-Hilbert.

\sssec{}

We need to show that the maps \eqref{e:Hecke univ Y} are isomorphisms in
$$\CC_{\Ran\times Y}\otimes \QCoh(\CZ).$$

We can rewrite this category as a \emph{limit} with terms
$$\CC^{\otimes I}\otimes \Shv(X^J\times Y)\otimes \QCoh(\CZ).$$

So, we need to show that the map in question becomes an isomorphism in each 
of the above terms.

%
%
%
%
%

\sssec{}

Note that for a scheme $W$, an object $\CF\in \Shv(W\times Y)$ is zero if and only if
for every geometric point $\Spec(k')\to Y$, the pullback $\CF$ to 
$$W':=\Spec(k')\underset{Y}\times (W\times Y)$$
is zero. 

\medskip

The same remains true for $\Shv(W\times Y)\otimes \bC$ for any DG category $\bC$. 

\medskip

Thus, it is sufficient to show that the map \eqref{e:Hecke univ Y} becomes an isomorphism
after the base change $k\rightsquigarrow k'$ for $Y'=\Spec(k')$. 

\sssec{}

However, the base change of the map \eqref{e:Hecke univ Y} is a similar map over the ground field $k'$ for
$\sF'$ being the following functor:

\medskip

In the \'etale context, $\sF'$ is
$$\CC \overset{\sF}\to \QCoh(\CZ) \underset{\sfe}\otimes \qLisse(X)\to \QCoh(\CZ) \underset{\sfe}\otimes \qLisse(X'),$$
and in the de Rham context, $\sF'$ is 
$$\CC \overset{\sF}\to \QCoh(\CZ) \underset{k}\otimes \qLisse(X)\to  \QCoh(\CZ') \underset{k'}\otimes \qLisse(X'),$$
where in both cases 
$$X':=\Spec(k')\underset{\Spec(k)}\times X,$$ 
and in the de Rham context 
$$\CZ':=\Spec(k')\underset{\Spec(k)}\times \CZ.$$

\qed[\propref{p:tightness}]

\sssec{}

Thus, we have reduced the verification of fact that \eqref{e:Hecke univ Y} is an isomorphism
to the case when $Y=\on{pt}$. However, in this case, the assertion is already known by 
\eqref{e:Hecke univ}.

%
%
%
%

\newpage

\centerline{\bf Part III: The category of automorphic sheaves with nilpotent singular support} 

\bigskip

Let us make a brief overview of the contents of this Part.

\medskip

In \secref{s:Nilp} we introduce and study the category $\Shv_\Nilp(\Bun_G)$. First, we state a key technical result, 
\thmref{t:preserve Nilp Sing Supp prel}, which says that
$\Bun_G$ can be covered by quasi-compact open substacks, such that the functor of !-extension from each
of them preserves the nilpotence of singular support. Next we observe that the action of Hecke functors  
on the entire category $\Shv(\Bun_G)$ gives rise to an action of $\Rep(\cG)^{\otimes X\on{-lisse}}$ on 
$\Shv_\Nilp(\Bun_G)$. Applying our Spectral Decomposition theorem, we obtain $\Shv_\Nilp(\Bun_G)$ carries a monoidal
action of $\QCoh(\LocSys_\cG^{\on{restr}}(X))$. We also state the second main result of this paper,
\thmref{t:lisse}, which says that if an object $\CF\in \Shv(\Bun_G)$ is such that the Hecke action on it is
lisse, then $\CF$ belongs to $\Shv_\Nilp(\Bun_G)$. This implies, in particular, that Hecke eigensheaves
have nilpotent singular support.

\medskip

In \secref{s:projector and eigen} we introduce yet another tool in the study of $\Shv(\Bun_G)$--Beilinson's spectral
projector, denoted $\sP^{\on{enh}}_\CZ$, which is defined for a prestack $\CZ$ equipped with a map 
$f:\CZ \to \LocSys_\cG^{\on{restr}}(X)$. This is a functor, given by an explicit Hecke operator, and it provides a left 
adjoint to the forgetful functor
$$\on{Hecke}(\CZ,\Shv(\Bun_G))\to \QCoh(\CZ)\otimes \Shv(\Bun_G)\to \Shv(\Bun_G).$$
Using our Spectral Decomposition theorem and \thmref{t:lisse}, we interpret 
$\sP^{\on{enh}}_\CZ$ as the left adjoint to the functor\footnote{Provided $\QCoh(\CZ)$ is dualizable 
and $\CO_\CZ\in \QCoh(\CZ)$ is compact.}
\begin{multline*}
\QCoh(\CZ)\underset{\QCoh(\LocSys_\cG^{\on{restr}}(X))}\otimes \Shv_{\on{Nilp}}(\Bun_G)\overset{f_*\otimes \on{Id}}\longrightarrow \\
\to \QCoh(\LocSys_\cG^{\on{restr}}(X)) \underset{\QCoh(\LocSys_\cG^{\on{restr}}(X))}\otimes \Shv_{\on{Nilp}}(\Bun_G)
\simeq \Shv_{\on{Nilp}}(\Bun_G) \hookrightarrow \Shv(\Bun_G).
\end{multline*}

\medskip

In \secref{s:projector} we use Beilinson's spectral projector to prove an array of structural results about $\Shv_{\on{Nilp}}(\Bun_G)$: we will show that
the category $\Shv_{\on{Nilp}}(\Bun_G)$ is compactly generated, that the external tensor product functor
$$\Shv_{\on{Nilp}}(\Bun_{G_1})\otimes \Shv_{\on{Nilp}}(\Bun_{G_2})\to \Shv_{\on{Nilp}}(\Bun_{G_1\times G_2})$$
is an equivalence and that, in the de Rham context, 
all objects in $\Shv_\Nilp(\Bun_G)$ have regular singularities.

\medskip

In \secref{s:more} we will make several observations regarding a conjecture, initially formulated in \secref{s:Nilp},
which can be stated as saying that the right adjoint of the embedding $\Shv_{\on{Nilp}}(\Bun_G) \hookrightarrow \Shv(\Bun_G)$
is continuous. 

\medskip

In \secref{s:spectral Betti}, we establish analogs of the results of the preceding sections in Part III, when work
over the ground field $k=\BC$ and instead of $\Shv(-)$ we consider the category $\Shv^{\on{all}}(-)$ of all sheaves in the
classical topology. 

\medskip

In \secref{s:preserve sing} we prove \thmref{t:preserve Nilp Sing Supp prel} about the preservation
of nilpotence of singular support under the functor of direct image for certain open embeddings $\CU\overset{j}\hookrightarrow \Bun_G$.
The proof follows closely the strategy of \cite{DrGa2}: by the same method as in {\it loc.cit.}, it turns out that we can
control the singular support of the extension in a \emph{contractive} situation. 

\medskip 

In \secref{s:proof lisse} we prove \thmref{t:lisse}. We first consider the case of $G=GL_2$, which explains the
main idea of the argument. We then implement this idea in a slightly more involved case of $G=GL_n$
(where it is sufficient consider the minuscule Hecke functors).  Finally, we treat the case of an arbitrary $G$;
the proof reduces to the analysis of the local Hitchin map and affine Springer fibers. 

\bigskip

\section{Automorphic sheaves with nilpotent singular support and spectral decomposition} \label{s:Nilp}

In this section we introduce and study the category $\Shv_\Nilp(\Bun_G)$. 

\medskip

The central results of this section are: 

\medskip

\noindent--\thmref{t:preserve Nilp Sing Supp prel}, which expresses a locality property of the 
nilpotence of singular support condition; 

\medskip

\noindent--\thmref{t:spectral decomp}, which says that $\Shv_\Nilp(\Bun_G)$ carries a monoidal action 
of $\QCoh(\LocSys_\cG^{\on{restr}}(X))$;

\medskip

\noindent--\thmref{t:lisse} that any object of $\Shv(\Bun_G)$ on which the Hecke action
is lisse, belongs to $\Shv_\Nilp(\Bun_G)$. 

\ssec{Definition and basic properties}

In this subsection we define the category $\Shv_{\on{Nilp}}(\Bun_G)$ and 
formulate a key result (\thmref{t:preserve Nilp Sing Supp prel}) that ensures that
it is, in a certain sense, local with respect to $\Bun_G$. 

\sssec{}

From now on we let $X$ be a smooth, connected and complete curve and $G$ a reductive group, over a ground field $k$
(assumed algebraically closed). 

\medskip

Consider $\Bun_G$, the moduli space of principal $G$-bundles on $X$. Our object of study is the category 
$$\Shv(\Bun_G)$$
of sheaves on $\Bun_G$. (The basics of the theory of sheaves on algebraic stacks are 
reviewed in \secref{s:shvs on stacks}.)  

\sssec{}

Recall that $T^*(\Bun_G)$ can be identified with the moduli space of pairs
$(\CP_G,A)$, where $\CP_G$ is a $G$-bundle on $X$, and $A$ is a global section of $\fg^\vee_{\CP_G}\otimes \omega_X$. 

\medskip

Let $\on{Nilp}\subset T^*(\Bun_G)$ be the nilpotent cone, i.e., the closed subset consisting of those $(\CP_G,A)$, for
which $A$ is nilpotent (at the generic point of $X$). 

\medskip

When $\on{char}(k)=0$, it is well-known that $\Nilp$ is half-dimensional (and even Lagrangian). 

\medskip

When $\on{char}(k)$ is positive, but a ``very good" prime for $G$, the corresponding assertion
seems to have been known in the folklore, but we were not able to find a proof in the literature.  
For completeness, we will supply a proof in \secref{s:glob Nilp}. 

\medskip

From now on, we will assume that the above restrictions on $\on{char}(k)$ are satisfied, so that 
$\Nilp$ is half-dimensional.

\sssec{}

The main object of study in this part is the subcategory
$$\Shv_{\on{Nilp}}(\Bun_G)\subset \Shv(\Bun_G),$$
see \secref{sss:non qc N}. 

\medskip

We denote the tautological embedding $\Shv_{\on{Nilp}}(\Bun_G)\hookrightarrow \Shv(\Bun_G)$ by $\iota$.

%
%

\sssec{}

For an open substack $\CU\subset \Shv(\Bun_G)$, we can consider the full subcategory
$$\Shv_{\on{Nilp}}(\CU)\subset \Shv(\CU).$$

We have the following result, which insures that the category $\Shv_{\on{Nilp}}(\Bun_G)$ can be obtained as a colimit of the
corresponding categories on quasi-compact open substacks of $\Bun_G$, see \secref{sss:N preserved abs a}: 

\begin{mainthm} \label{t:preserve Nilp Sing Supp prel}
The stack $\Shv(\Bun_G)$ can be written as a filtered union of quasi-compact open substacks 
$$\CU_i\overset{j_i}\hookrightarrow \Bun_G$$
such that the extension functors
$$(j_i)_!,(j_i)_*:\Shv(\CU_i)\to \Shv(\Bun_G)$$
send $\Shv_{\on{Nilp}}(\CU_i)^{\on{constr}}\to \Shv_{\on{Nilp}}(\Bun_G)^{\on{constr}}$.
\end{mainthm} 

The proof will be given in \secref{s:preserve sing}.  

\begin{rem} \label{r:preserve Nilp Sing Supp}

In the terminology of \secref{sss:N preserved abs a}, \thmref{t:preserve Nilp Sing Supp prel} says that
the pair $\Bun_G$ is $\Nilp$-\emph{truncatable}.

\medskip

By \secref{sss:N preserved abs b}, the statement of \thmref{t:preserve Nilp Sing Supp prel}
can be reformulated as the assertion that that for $(\CU_i,j_i)$ as above, 
the functors $(j_i)_!$ and $(j_i)_*$ send 
$$\Shv_{\on{Nilp}}(\CU_i)\to \Shv_{\on{Nilp}}(\Bun_G).$$

\end{rem}


\sssec{}

Here is one property of the category $\Shv(\Bun_G)$ that we expect to hold, but at the 
moment are unable to prove in general (but we can prove it in the de Rham and Betti contexts, 
see Theorems \ref{t:Nilp comp gen dR} and \ref{t:Nilp comp gen Betti}):

\begin{conj} \label{c:Nilp comp gen}
The category $\Shv_{\on{Nilp}}(\Bun_G)$ is generated by objects that are compact in the ambient category
$\Shv(\Bun_G)$.
\end{conj} 

In the sequel, we will give several (equivalent) reformulations of \conjref{c:Nilp comp gen}.

%
%
%

\begin{rem}
In the terminology of \secref{sss:non qc N constraccess}, the above conjecture says that the pair $(\Bun_G,\on{Nilp})$
is \emph{renormalization-adapted} and \emph{constraccessible}. 

\medskip

Given \thmref{t:preserve Nilp Sing Supp prel}, the property of being \emph{renormalization-adapted} is known to hold since 
$\Bun_G$ is locally a quotient, see \corref{c:N preserved abs} (and is in fact expected to hold for any pair $(\CY,\CN)$ of an algebraic 
stack (with an affine diagonal) and a subset of its cotangent bundle, see \conjref{c:ren adapt}). 

\medskip

The property of being \emph{constraccessible} is much more mysterious: it reflects a particular feature of the pair $(\Bun_G,\on{Nilp})$
(for example, it fails for $(\BP^1,\{0\})$, see Remark \ref{r:bad P1}).

\end{rem}
%

\ssec{Hecke action on the category with nilpotent singular support}

In this subsection we recall the pattern of Hecke action on $\Shv(\Bun_G)$, and the particular feature that the subcategory
$$\Shv_{\on{Nilp}}(\Bun_G) \subset \Shv(\Bun_G)$$
has with respect to this action. 

\sssec{}

Let $\cG$ denote the Langlands dual group of $G$. 
The following result encodes the Hecke action of $\Rep(\cG)$ on $\Shv(\Bun_G)$
(see \cite[Proposition B.2.3]{GKRV}):

\begin{thm} \label{t:Hecke}
The Hecke functors combine to a compatible family of actions of 
\begin{equation} \label{e:Hecke action initial}
\Rep(\cG)^{\otimes I}\otimes \Shv(X^I) \text{ on } \Shv(\Bun_G\times X^I), \quad I\in \on{fSet},
\end{equation}
extending the tautological action of 
$$\Shv(X^I) \text{ on } \Shv(\Bun_G\times X^I), \quad I\in \on{fSet}.$$
\end{thm}

\sssec{}

We are going to combine \thmref{t:Hecke} with the following result, 
established in \cite[Theorem 5.2.1]{NY1} (see also \cite[Theorem B.5.2]{GKRV}):

\begin{thm} \label{t:NY}
The Hecke functor
\begin{equation} \label{e:1 step Hecke}
\on{H}(-,-):\Rep(\cG)\otimes \Shv(\Bun_G) \to \Shv(\Bun_G\times X)
\end{equation} 
sends 
$$\Shv_{\on{Nilp}}(\Bun_G) \subset \Shv(\Bun_G)$$
to the full subcategory
$$\Shv_{\on{Nilp}\times \{0\}}(\Bun_G \times X) \subset \Shv(\Bun_G\times X).$$
\end{thm} 

\sssec{}

Note that by \thmref{t:product thm stack 2} and \corref{c:curve Verdier compat},
the external tensor product functor
$$\Shv_{\Nilp}(\Bun_G)\otimes \qLisse(X) \to \Shv_{\on{Nilp}\times \{0\}}(\Bun_G\times X)$$
is an equivalence.

\medskip

Hence, the Hecke functor restricted to the category of sheaves with nilpotent singular support can be viewed as a functor
\begin{equation} \label{e:1 step Hecke bis}
\Rep(\cG)\otimes \Shv_\Nilp(\Bun_G)\to \Shv_{\Nilp}(\Bun_G)\otimes \qLisse(X).
\end{equation} 

\begin{rem} \label{r:why access prel}

In Remark \ref{r:why access} we will see that \eqref{e:1 step Hecke bis} can be somewhat refined: the essential
image of the functor \eqref{e:1 step Hecke bis} actually belongs to
$$\Shv_{\on{Nilp}}(\Bun_G)\otimes \iLisse(X) \subset \Shv_{\on{Nilp}}(\Bun_G)\otimes \qLisse(X).$$

\end{rem} 

\sssec{} \label{sss:Hecke action on Nilp}

Iterating, from \thmref{t:NY}, we obtain that for any $I\in \on{fSet}$, the Hecke functors 
\begin{equation} \label{e:full Hecke}
\on{H}(-,-):\Rep(\cG)^{\otimes I}\otimes \Shv(\Bun_G) \to \Shv(\Bun_G\times X^I), \quad I\in \on{fSet}
\end{equation} 
define a system of functors
\begin{equation} \label{e:full Hecke bis}
\Rep(\cG)^{\otimes I}\otimes \Shv_{\Nilp}(\Bun_G) \to \Shv_{\Nilp}(\Bun_G)\otimes \qLisse(X)^{\otimes I}.
\end{equation}

\medskip

Combining with \thmref{t:Hecke}, we obtain:

\begin{cor} \label{c:Hecke nilp}
The Hecke action gives rise to a compatible family of monoidal functors
$$\Rep(\cG)^{\otimes I}\to \End(\Shv_{\on{Nilp}}(\Bun_G))\otimes \qLisse(X)^{\otimes I}, \quad I\in \on{fSet}.$$
\end{cor}

\sssec{}

Thus, in the terminology of \secref{sss:lisse action}, we obtain that the Hecke action gives rise to an action of
$\Rep(\cG)^{\otimes X\on{-lisse}}$ on $\Shv_{\on{Nilp}}(\Bun_G)$. 

\ssec{Spectral decomposition of the category with nilpotent singular support}

We now come to the first main point of this paper: the spectral decomposition of $\Shv_\Nilp(\Bun_G)$
over $\LocSys_\cG^{\on{restr}}(X)$.

\sssec{}

Combining \corref{c:Hecke nilp} with \thmref{t:action}, we obtain: 

\begin{mainthm} \label{t:spectral decomp}
The action of $\Rep(\cG)^{\otimes X\on{-lisse}}$ on $\Shv_{\on{Nilp}}(\Bun_G)$ (arising from the Hecke action) 
factors via a (uniquely defined) action of $\QCoh(\LocSys_\cG^{\on{restr}}(X))$.
\end{mainthm}

\sssec{} \label{sss:EV again}

Let us emphasize the main feature of the action of $\QCoh(\LocSys_\cG^{\on{restr}}(X))$ on 
$\Shv_{\on{Nilp}}(\Bun_G)$.

\medskip

For a finite set $I$ and $V\in \Rep(\cG)^{\otimes I}$, let 
\begin{equation} \label{e:EV again}
\CE^I_V\in \QCoh(\LocSys_\cG^{\on{restr}}(X))\otimes \qLisse(X)^{\otimes I}
\end{equation}
be as in \secref{sss:pre-shtuka}. 

\medskip

Then the action of $\CE^I_V$ on $\Shv_{\on{Nilp}}(\Bun_G)$, viewed as a functor
$$\Shv_{\on{Nilp}}(\Bun_G) \to \Shv_{\on{Nilp}}(\Bun_G) \otimes \qLisse(X)^{\otimes I}$$
equals the Hecke functor 
$$\on{H}(V,-):\Shv_{\Nilp}(\Bun_G) \to \Shv_{\Nilp}(\Bun_G)\otimes \qLisse(X)^{\otimes I}$$
of \eqref{e:full Hecke bis}.



\sssec{}

As a first corollary of \thmref{t:spectral decomp} (combined with \propref{p:conn comps LocSys}), we obtain: 

\begin{cor} \label{c:spectral coarse}
The category $\Shv_{\on{Nilp}}(\Bun_G)$ splits canonically as a direct sum
$$\Shv_{\on{Nilp}}(\Bun_G) \simeq \underset{\sigma}\oplus\, \Shv_{\on{Nilp}}(\Bun_G)_\sigma,$$
where $\sigma$ runs over the set of isomorphism classes of semi-simple $\cG$-local systems on $X$.
\end{cor}

\sssec{Example} \label{sss:decomp Gm}
Let us explain what \corref{c:spectral coarse} says in concrete terms for $G=\BG_m$. 

\medskip

Recall that the geometric class field theory attaches to a 1-dimensional local system $\sigma$ on $X$
a local system $E_\sigma$ on $\on{Pic}$. Then \corref{c:spectral coarse} is the assertion that 
$\qLisse(\on{Pic})$ splits as a direct sum
$$\qLisse(\on{Pic})\simeq \underset{\sigma}\oplus\, \qLisse(\on{Pic})_\sigma,$$
where each $\qLisse(\on{Pic})_\sigma$ is generated by $E_\sigma$. 

\medskip

In the particular case of $\on{Pic}$, such a decomposition is not difficult to establish directly: it follows
from the fact that every lisse irreducible object in $\Shv(\on{Pic})$ is isomorphic to one of the $E_\sigma$ 
(this is the assertion that the \'etale fundamental group of $\on{Pic}$ is the abelianization of the \'etale
fundamental group of $X$) and the different $E_\sigma$ are mutually orthogonal. 

\sssec{}

From \corref{c:spectral coarse} (combined with \corref{c:eigensheaves nilp} below) we obtain the following result:

\begin{cor} \label{c:eigen orth}
Let $\CF_1$ and $\CF_2$ be Hecke eigensheaves corresponding to $\sG$-local systems $\sigma_1$
and $\sigma_2$ with non-isomorphic semi-simplifications. Then $\CF_1$ and $\CF_2$ are
mutually orthogonal, i.e.,
$$\CMaps(\CF_1,\CF_2)=0.$$
\end{cor} 

\begin{rem} \label{r:why access}

Let us show that the functor \eqref{e:1 step Hecke bis} has essential image in 
$$\Shv_{\on{Nilp}}(\Bun_G)\otimes \iLisse(X) \subset \Shv_{\on{Nilp}}(\Bun_G)\otimes \qLisse(X).$$

\medskip

This follows from the fact that the objects $\CE_V$ above in fact belong to
$$\QCoh(\LocSys_\cG^{\on{restr}}(X))\otimes \iLisse(X)\subset \QCoh(\LocSys_\cG^{\on{restr}}(X))\otimes \qLisse(X).$$

To prove this, it is enough to show that for a \emph{cofinal} family of maps
$S\to \LocSys_\cG^{\on{restr}}(X)$, the objects
$$\CE_{V}|_S\in \QCoh(S)\otimes \qLisse(X)$$
belong to
$$\QCoh(S)\otimes \iLisse(X)\subset \QCoh(S)\otimes \qLisse(X).$$

We now use the fact that for $X$ a curve, the prestack $\LocSys_\cG^{\on{restr}}(X)$ is 
\emph{eventually coconnective} (see \cite[Chapter 2, Sect. 1.3.5]{GR1}). 
This follows from the fact that the connected components of $\LocSys_\cG^{\on{restr},\on{rigid}_x}(X)$
are \emph{quasi-smooth formal affine schemes}, see Sects. \ref{sss:formal q-smooth} and \ref{sss:LocSys q-sm}. 

\medskip

Hence, inside the category $\affSch_{/\LocSys_\cG^{\on{restr}}(X)}$, a cofinal family is formed by those $S$ that
are eventually coconnective. Now, for $S$ eventually coconnective, the fact that
$\CE_{V}|_S$ belongs to the subcategory $\QCoh(S)\otimes \iLisse(X)$ is a reformulation of \propref{p:replace by Ind-Lisse}.

\end{rem}

\ssec{A converse to \thmref{t:NY}} \label{ss:lisse thm}

In this subsection we state the second main result of this paper, \thmref{t:lisse}, which says
that the statement of \thmref{t:NY} is ``if and only if". This theorem 
implies, among the rest, that Hecke eigensheaves have nilpotent singular support. 

\sssec{} \label{sss:assump char}

For the validity of \thmref{t:lisse} we will have to make the following assumptions on $\on{char}(k)$:

\smallskip

\begin{itemize}

\item
There exists a non-degenerate $G$-equivariant symmetric bilinear form on $\fg$, 
whose restriction to the center of any Levi subalgebra remains non-degenerate;

\smallskip

\item The projection $\ft/\!/W\to \fg/\!/\on{Ad}(G)$ is an isomorphism, and similarly 
for any Levi subgroup of $G$. 

\smallskip

\item
The centralizer in $G$ of a semi-simple element in $\fg$ is a Levi subgroup.

\end{itemize}

\medskip

From now on, we will assume
that the above assumptions on $\on{char}(k)$ are satisfied. 

\medskip

For example, for $G=GL_n$, these assumptions are equivalent to $\on{char}(k)>n$. 

\medskip

In general, it is known the above conditions are satisfied away from very small characteristics for
a given type of $G$. 

\sssec{} \label{sss:defn Hecke-lisse}

Let
$$\Shv(\Bun_G)^{\on{Hecke-lisse}}\subset \Shv(\Bun_G)$$
be the full subcategory consisting of objects $\CF$ such that for all $V\in \Rep(\cG)$, we have
$$\on{H}(V,\CF)\in \Shv(\Bun_G)\otimes \qLisse(X)\subset \Shv(\Bun_G\times X),$$
where $\on{H}(V,\CF)$ is the Hecke functor of \eqref{e:1 step Hecke}. 

\medskip

We can phrase \thmref{t:NY} as saying that 
\begin{equation} \label{e:Nilp and lisse}
\Shv_{\on{Nilp}}(\Bun_G)\subset \Shv(\Bun_G)^{\on{Hecke-lisse}}.
\end{equation}

\medskip

The following was proposed as a conjecture in \cite{GKRV} (it appears as Conjecture C.2.8 in {\it loc.cit.}):  

\begin{mainthm} \label{t:lisse prel}
The inclusion \eqref{e:Nilp and lisse} is an equality. 
\end{mainthm}

In fact, we will prove a stronger result: 

\begin{mainthm} \label{t:lisse}
Let $\CF\in \Shv(\Bun_G)$ be such that for all $V\in \Rep(\cG)$, the singular support of 
the object 
$$\on{H}(V,\CF)\in \Shv(\Bun_G\times X)$$
is contained in $T^*(\Bun_G)\times \{0\}\subset T^*(\Bun_G\times X)$. Then $\CF\in \Shv_{\on{Nilp}}(\Bun_G)$.
\end{mainthm}

\begin{rem} \label{r:lisse hol}
When $\on{char}(k)=0$, the assertion of \thmref{t:lisse} is actually equivalent to that of \thmref{t:lisse prel},
by \corref{c:1.75}.
\end{rem} 

\begin{rem} \label{r:lisse other dR}
Recall that in our notation $\Shv(\Bun_G)$ refers to a constructible sheaf theory. However, the statement
of \thmref{t:lisse prel} remains valid, when instead of $\Shv(-)$ we consider $\Dmod(-)$ 
(when $\on{char}(k)=0$ and $\sfe=k$). We will prove this in \secref{ss:lisse for D-mod}. \

\medskip

Note that in the case of D-modules, the inclusion
\begin{equation} \label{e:incl Nilp Dmod}
\Shv_{\on{Nilp}}(\Bun_G)\subset \Dmod_{\on{Nilp}}(\Bun_G)
\end{equation} 
is an equality: since $\Nilp$ is Lagrangian, every object of $\Dmod_{\on{Nilp}}(\Bun_G)$
is necessarily ind-holonomic. 
\end{rem} 

\begin{rem} \label{r:lisse other dR big}
We conjecture that the statement of (the stronger) \thmref{t:lisse} also remains valid for $\Dmod(-)$.
In fact, it would follow if we knew that \thmref{t:sing supp dir im} holds for $\Dmod(-)$, see
Remark \ref{r:lisse for Dmod gen}. 
\end{rem}


\sssec{}

From \thmref{t:lisse} we obtain\footnote{The conclusion of \corref{c:loose} appears in \cite{GKRV} as Conjecture C.2.10.}:

\begin{maincor} \label{c:loose}
Let $\CF\in \Shv(\Bun_G)$ be a \emph{loose} Hecke eigensheaf, i.e., for every $V\in \Rep(\cG)^\heartsuit$, the object
$$\on{H}(V,\CF)\in \Shv(\Bun_G\times X)$$
is of the form $\CF\boxtimes E_V$ for some $E_V\in \qLisse(X)$. Then $\CF\in \Shv_\Nilp(\Bun_G)$.
\end{maincor}

As a particular case, we obtain the following statement, which was conjectured by G.~Laumon
(\cite[Conjecture 6.3.1]{laumon}):

\begin{maincor} \label{c:eigensheaves nilp}
Hecke eigensheaves in $\Shv(\Bun_G)$ have nilpotent singular support.
\end{maincor}

\sssec{Example}  \label{sss:Hecke lisse Gm}

The assertion of \thmref{t:lisse prel} is easy for $G=\BG_m$. Note that in this case $\Bun_G=\on{Pic}$, and 
$$\Shv_{\on{Nilp}}(\Bun_G)=\qLisse(\on{Pic}).$$

\begin{proof}

The Hecke functor for the standard character of $\cG=\BG_m$ is the pullback functor with respect to the addition map
$$\on{add}:\on{Pic}\times X \to \on{Pic}.$$

\medskip

Let us be given an object $\CF\in \Shv(\on{Pic})$ such that 
$$\on{add}^!(\CF)\in \Shv(\on{Pic}\times X)$$
belongs to 
$$\Shv(\on{Pic})\otimes \qLisse(X)\subset \Shv(\on{Pic}\times X)$$. 

\medskip

We wish to show that $\CF$ belongs to $\qLisse(\on{Pic})$. It is easy to see that
it is enough to prove that $\CF|_{\on{Pic}^d}$ belongs to $\qLisse(\on{Pic}^d)$
for some/any $d$. 

\medskip

By \cite[Proposition C.2.5]{GKRV} quoted above, for any integer $d$ we have
$$\on{add}_d^!(\CF)\in \Shv(\on{Pic})\otimes \qLisse(X^d),$$
where $\on{add}_d$ is the $d$-fold addition map
$$\on{add}_d:\on{Pic}\times X^d \to \on{Pic}.$$

In particular, the !-pullback of $\CF$ along
\begin{equation} \label{e:add d}
X^d\simeq \one_{\on{Pic}}\times X^d \to \on{Pic}\times X^d \to \on{Pic}^d
\end{equation}
belongs to $\qLisse(X^d)$.  

\medskip

Note that the map \eqref{e:add d} factors as 
$$X^d\overset{\on{sym}^d}\longrightarrow X^{(d)} \overset{\on{AJ}_d}\to \on{Pic}^d,$$
where $\on{AJ}_d$ is the Abel-Jacobi map. For $d>2g-2$, the map $\on{AJ}_d$ is smooth and surjective. 
Let $\overset{\circ}X{}^d\subset X^d$ be the complement of the diagonal divisor. For $d\gg 0$, the composite map 
$$\overset{\circ}X{}^d\hookrightarrow X^d \overset{\text{\eqref{e:add d}}}\longrightarrow \on{Pic}^d,$$
is also surjective. It is smooth because the map $\on{sym}^d$ is \'etale when restricted to $\overset{\circ}X{}^d$. 

\medskip

Hence, for such $d$, if $\CF|_{\oX^d}$ is lisse, then so is $\CF$. 

\end{proof}

\sssec{}

Let us now prove (the stronger) \thmref{t:lisse} for $\BG_m$. In fact, the proof is even simpler
(and will be the prototype of the proof of \thmref{t:lisse} for any $G$): 

\medskip

Consider again the map
$$\on{add}:\on{Pic}\times X\to \on{Pic}.$$

We identify the cotangent space to $\on{Pic}$ at any $\CL\in \on{Pic}$ with $\Gamma(X,\omega_X)$.
Then the codifferential of $\on{add}$ at any $(\CL,x)\in \on{Pic}\times X$ 
is the map
$$\Gamma(X,\omega_X)\to \Gamma(X,\omega_X)\oplus T^*_x(X),$$
whose first component is the identity map and the second component is the evaluation map at $x$.

\medskip

Now, if for $\CF\in \Shv(\on{Pic})$ its singular support does not lie in $\Nilp=\{0\}$, we can find 
$\CL \in \on{Pic}$ and non-zero $\xi\in H^0(T^*_\CL(\on{Pic}))$ such that
$$(\xi,\CL)\in \on{SingSupp}(\CF).$$

Let $x\in X$ be such that the value $\xi|_x\in T^*_x(X)$ of $\xi$ at $x$ is non-zero. 
Since the map $\on{add}$ is smooth, the element
$$((\xi,\xi|_x),(\CL,x))\in T^*(\on{Pic}\times X)$$
belongs to $\on{SingSupp}(\on{add}^!(\CF))$. 

\medskip

However, $\xi|_x\neq 0$ by assumption, and we have obtained a contradiction with the fact that
$$\on{SingSupp}(\on{add}^!(\CF))\in T^*(\Bun_G)\times \{0\}.$$

\ssec{Spectral decomposition in the de Rham context}

In this subsection we will assume that our ground field $k$ has characteristic $0$,
and we will work with the entire category of D-modules, i.e., $\Dmod(-)$ instead of $\Shv(-)$. 

\sssec{}

Hecke action in the context of D-modules is a compatible family of functors
$$\Rep(\cG)^{\otimes I}\otimes \Dmod(\Bun_G) \to \Dmod(\Bun_G\times X^I)\simeq  \Dmod(\Bun_G)\otimes \Dmod(X)^{\otimes I}, \quad I\in \on{fSet}.$$

\medskip

It was shown in \cite[Corollary 4.5.5]{Ga5} that the above family of functors 
comes from a (uniquely defined) action of the category $\QCoh(\LocSys^{\dr}_\cG(X))$ on $\Dmod(\Bun_G)$. 

\medskip

Here again, for a fixed $V\in \Rep(\cG)$, the corresponding functor
$$\on{H}(V,-):\Dmod(\Bun_G) \to \Dmod(\Bun_G\times X)\simeq  \Dmod(\Bun_G)\otimes \Dmod(X)$$
is given by the action of the object
$$\CE_V\in \QCoh(\LocSys^{\dr}_\cG(X))\otimes \Dmod(X),$$
see \secref{sss:action de Rham}. 

\sssec{}

We now claim:

\begin{prop} \label{p:nilp via spectral}
The full subcategory
$$\Shv_{\on{Nilp}}(\Bun_G)\subset \Dmod(\Bun_G)$$
equals 
\begin{multline*} 
\Dmod(\Bun_G) \underset{\QCoh(\LocSys^{\dr}_\cG(X))}\otimes \QCoh(\LocSys_\cG^{\on{restr}}(X))
\subset  \\
\subset 
\Dmod(\Bun_G)\underset{\QCoh(\LocSys^{\dr}_\cG(X))}\otimes \QCoh(\LocSys^{\dr}_\cG(X))=\Dmod(\Bun_G),
\end{multline*} 
where we view
$$\QCoh(\LocSys_\cG^{\on{restr}}(X))\simeq \QCoh(\LocSys^{\dr}_\cG(X))_{\LocSys_\cG^{\on{restr}}(X)}$$
as a co-localization of $\QCoh(\LocSys^{\dr}_\cG(X))$. 
\end{prop}

\begin{proof}

This is obtained by combining \propref{p:lisse subcategory} with \thmref{t:lisse prel}, applied for all D-modules
(see Remark \ref{r:lisse other dR}). 

\end{proof} 

\sssec{}

From \propref{p:nilp via spectral} we deduce:

\begin{cor} \label{c:nilp via spectral}
In the de Rham context, the embedding $\iota:\Shv_{\on{Nilp}}(\Bun_G)\hookrightarrow \Dmod(\Bun_G)$
admits a \emph{continuous} right adjoint.
\end{cor}

\begin{proof}

The right adjoint is obtained by tensoring $-\underset{\QCoh(\LocSys^{\dr}_\cG(X))}\otimes \Dmod(\Bun_G)$ 
from the right adjoint to
$$\QCoh(\LocSys_\cG^{\on{restr}}(X))\simeq \QCoh(\LocSys^{\dr}_\cG(X))_{\LocSys_\cG^{\on{restr}}(X)}\hookrightarrow \QCoh(\LocSys^{\dr}_\cG(X)).$$

\end{proof} 

\section{Beilinson's spectral projector and Hecke eigensheaves} \label{s:projector and eigen}

In this section we recall the notion of Hecke eigensheaf and introduce 
Beilinson's spectral projector as a tool of constructing them. 


\ssec{Ran version of the Hecke action and Hecke eigensheaves} \label{ss:Hecke Ran}

In this subsection, we will show how the category $\Rep(\cG)_\Ran$ of \secref{ss:C Ran}
(and its version with a scheme of parameters, see \secref{ss:Ran param}) encode the formalism
of Hecke action. 

\sssec{} \label{sss:Ran version of the Hecke action new}

Recall the setting of \secref{ss:proj abs}. We claim that the category 
$\Shv(\Bun_G)$ carries an action of the (symmetric) monoidal category $\Rep(\cG)_{\Ran}$. 

\medskip

Namely, for 
\begin{equation} \label{e:Ran version of the Hecke action Bun}
(I\overset{\psi}\to J)\in \on{TwArr}(\on{fSet}), \quad V\in \Rep(\cG)^{\otimes I}, \quad \CS\in \Shv(X^J),
\end{equation} 
and 
\begin{equation} \label{e:object of Ran Bun}
\CV:=\on{ins}_{I\to J}(V\otimes \CS)\in \Rep(\cG)_{\Ran},
\end{equation} 
we let the action of $\CV$ on $\Shv(\Bun_G)$ be given by 
$$\Shv(\Bun_G) \overset{\on{H}(V,-)}\longrightarrow 
\Shv(\Bun_G\times X^I) \overset{(\on{id}\times \Delta_\psi)^!}\longrightarrow
\Shv(\Bun_G\times X^J)  \overset{\sotimes \CS}\to $$
$$\to \Shv(\Bun_G\times X^J) \overset{(p_{\Bun_G})_*}\longrightarrow 
\Shv(\Bun_G),$$ 
where:

\begin{itemize}

\item The symbol $\sotimes \CS$ means !-tensor product by the !-pullback of $\CS$;

\item $p_{\Bun_G}$ denotes the projection $\Bun_G\times X^J\to X^J$. 

\end{itemize}

\medskip

We will denote the resulting action of $\Rep(\cG)_{\Ran}$ on $\Shv(\Bun_G)$ by
$$\CV \in \Rep(\cG)_{\Ran}, \,\, \CF\in \Shv(\Bun_G) \mapsto \CV\star \CF,$$
following the notation of \secref{sss:Hecke abs}.



\sssec{} \label{sss:Hecke main}

Let now $\CZ$ be a prestack over $\sfe$, equipped with a map 
$$f:\CZ\to \LocSys^{\on{restr}}_\cG(X).$$

Let $\sF$ denote the resulting functor
$$\Rep(\cG)\to \QCoh(\CZ)\otimes \qLisse(X).$$

We will study the resulting category of Hecke eigen-objects
$$\on{Hecke}(\CZ,\Shv(\Bun_G)),$$
see \secref{sss:Hecke abs}. 

\medskip

We can think of $\on{Hecke}(\CZ,\Shv(\Bun_G))$ as objects $\CF\in \Shv(\Bun_G)\otimes \QCoh(\CZ)$
equipped with a system of isomorphisms
\begin{equation} \label{e:Hecke property Ran Bun}
\CV \star \CF\simeq \CF\otimes \wt\sF(\CV), \quad \CV\in \Rep(\cG)_{\Ran}.
\end{equation}

\sssec{} \label{sss:Ran version of the Hecke action param}

We will now rerun the above story, in the presence of a scheme of parameters $Y$. 

\medskip

Recall the setting of \secref{ss:Ran param}. We claim that for any test $k$-scheme $Y$, the category 
$\Shv(\Bun_G\times Y)$ carries an action of the (symmetric) monoidal category $\Rep(\cG)_{\Ran\times Y}$. 

\medskip

Namely, for 
\begin{equation} \label{e:Ran version of the Hecke action Y}
(I\overset{\psi}\to J)\in \on{TwArr}(\on{fSet}), \quad V\in \Rep(\cG)^{\otimes I}, \quad \CS\in \Shv(X^J\times Y),
\end{equation} 
and 
\begin{equation} \label{e:object of Ran Y}
\CV_Y:=\on{ins}_{I\to J}(V\otimes \CS)\in \Rep(\cG)_{\Ran\times Y},
\end{equation} 
we let the action of $\CV_Y$ on $\Shv(\Bun_G\times Y)$ be given by 
$$\Shv(\Bun_G\times Y) \overset{\on{H}(V,-)}\longrightarrow 
\Shv(\Bun_G\times X^I\times Y) \overset{(\on{id}\times \Delta_\psi\times \on{id})^!}\longrightarrow
\Shv(\Bun_G\times X^J\times Y)  \overset{\sotimes \CS}\to $$
$$\to \Shv(\Bun_G\times X^J\times Y) \overset{(p_{\Bun_G\times Y})_*}\longrightarrow 
\Shv(\Bun_G\times Y),$$ 
where:

\begin{itemize}

\item We denote by $\on{H}(V,-)$ the corresponding version of the Hecke functor, where
we have a scheme $Y$ of parameters;

\item The symbol $\sotimes \CS$ means !-tensor product by the !-pullback of $\CS$;

\smallskip

\item $p_{\Bun_G\times Y}$ denotes the projection $\Bun_G\times X^J\times Y\to \Bun_G\times Y$.

\end{itemize}

\sssec{}

For a pair $(\CZ,f:\CZ\to \LocSys^{\on{restr}}_\cG(X))$ as above, we can consider the resulting categories 
$$\on{Hecke}_Y(\CZ,\Shv(\Bun_G\times Y)) \text{ and } \on{Hecke}(\CZ,\Shv(\Bun_G\times Y)).$$
see \secref{sss:Hecke param}. 

\medskip

Note, however, that thanks to \thmref{t:get rid of Y}, the forgetful functor 
$$\on{Hecke}_Y(\CZ,\Shv(\Bun_G\times Y)) \to \on{Hecke}(\CZ,\Shv(\Bun_G\times Y))$$
is an equivalence. 

\sssec{}

For a map $Y_1\to Y_2$, the functor of !-pullback 
$$\Shv(\Bun_G\times Y_2)\to \Shv(\Bun_G\times Y_1)$$
is compatible with the action of $\Rep(\cG)_\Ran$. Hence, it induces a functor  
$$\on{Hecke}(\CZ,\Shv(\Bun_G\times Y_2))\to \on{Hecke}(\CZ,\Shv(\Bun_G\times Y_1)).$$

\medskip

In particular, for every $Y$, we have a canonically defined functor
$$\on{Hecke}(\CZ,\Shv(\Bun_G)) \to \on{Hecke}_Y(\CZ,\Shv(\Bun_G\times Y))$$
that fits into a commutative diagram
$$
\CD
\on{Hecke}(\CZ,\Shv(\Bun_G)) @>>>  \on{Hecke}_Y(\CZ,\Shv(\Bun_G\times Y)) \\
@VVV @VVV \\
\Shv(\Bun_G)\otimes \QCoh(\CZ) @>>> \Shv(\Bun_G\times Y)\otimes \QCoh(\CZ),
\endCD
$$
where the vertical arrows are forgetful functors, and the bottom horizontal arrow is the
!-pullback functor. 

\medskip

Hence, for $\CF\in \on{Hecke}(\CZ,\Shv(\Bun_G))$ and an object $\CV_Y\in \Rep(\cG)_{\Ran\times Y}$, 
we have a Hecke isomorphism in $\Shv(\Bun_G\times Y)\otimes \QCoh(\CZ)$:
\begin{equation} \label{e:Hecke property Ran Y}
\CV_Y \star (\CF\boxtimes \omega_Y)\simeq (\CF\boxtimes \omega_Y)\sotimes \wt\sF_Y(\CV_Y), \quad \CV_Y\in \Rep(\cG)_{\Ran\times Y},
\end{equation} 
for $\wt\sF_Y$ as in \secref{sss:F Y}, where $\sotimes$ denotes the natural action of $\QCoh(\CZ)\otimes \Shv(Y)$ on
$$\Shv(\Bun_G\times Y)\otimes \QCoh(\CZ).$$

\sssec{}

All of the above discussion is equally applicable when instead of $\Shv(-)$ we consider 
$\Dmod(-)$ (and $f$ is a map $\CZ\to \LocSys^\dr_\cG(X)$) or $\Shv^{\on{all}}(-)$ 
(and $f$ is a map $\CZ\to \LocSys^{\on{Betti}}_\cG(X)$). 

\ssec{Another notion of Hecke eigensheaf} \label{ss:usual Hecke}

In this subsection we relate the category $$\on{Hecke}(\CZ,\Shv(\Bun_G))$$ to another, probably more 
familiar, notion of Hecke eigensheaf.

\sssec{} \label{sss:usual Hecke}

We will show that the category $\on{Hecke}(\CZ,\Shv(\Bun_G))$ can be identified with the category of objects
$$\CF\in \Shv(\Bun_G)\otimes \QCoh(\CZ),$$ equipped with a system of isomorphisms
\begin{equation} \label{e:Hecke property usual}
\on{H}(V,\CF) \overset{\alpha_V}\simeq \CF \underset{\CO_\CZ}\boxtimes \sF^I(V), \quad V\in \Rep(\cG)^{\otimes I}, \quad I\in \fSet,
\end{equation} 
where:

\begin{itemize}

\item $\sF^I$ is the functor $\Rep(\cG)^{\otimes I}\to \QCoh(\CZ)\otimes \qLisse(X)^{\otimes I}$,
corresponding to $f$;

\item Both sides in \eqref{e:Hecke property usual} are viewed as objects of $\Shv(\Bun_G\times X^I)\otimes \QCoh(\CZ)$;

\item $\underset{\CO_\CZ}\boxtimes$ denotes the external tensor product functor
$$(\Shv(\Bun_G) \otimes  \QCoh(\CZ))\otimes (\QCoh(\CZ)\otimes \qLisse(X)^{\otimes I})\to \Shv(\Bun_G\times X^I)\otimes \QCoh(\CZ).$$

\end{itemize}

\bigskip

The isomorphisms \eqref{e:Hecke property usual} are required to be equipped with a homotopy-coherent system of 
compatibilities for maps between finite sets: 

\medskip

For $I_1\overset{\phi}\to I_2$ and $V_1\in \Rep(\cG)^{\otimes I_1}$, we must be given a data of commutativity for the diagram
$$
\CD
((\on{id}_{\Bun_G}\times \Delta_\phi)^!\otimes \on{Id}_{\QCoh(\CZ)})
(\on{H}(V_1,\CF)) @>{\alpha_{V_1}}>> ((\on{id}_{\Bun_G}\times \Delta_\phi)^!\otimes \on{Id}_{\QCoh(\CZ)})(\CF \underset{\CO_\CZ}\boxtimes \sF^{I_1}(V_1)) \\
@V{\sim}VV @VV{\sim}V \\
\on{H}(V_2,\CF) @>{\alpha_{V_2}}>> \CF \underset{\CO_\CZ}\boxtimes \sF^{I_2}(V_2),
\endCD
$$
where 
$$V_2:=\on{mult}^\phi_{\Rep(\cG)}(V_1)\in \Rep(\cG)^{\otimes I_2}.$$

\sssec{} \label{sss:usual Hecke one dir}

Indeed, let $\CF$ be an object of $\Shv(\Bun_G)\otimes \QCoh(\CZ)$ equipped with a lift to an object of 
$\on{Hecke}(\CZ,\Shv(\Bun_G))$ in the definition of \secref{sss:Hecke main}, 
and let us construct the isomorphisms \eqref{e:Hecke property usual}. 

\medskip

For a finite set $I$ and $V\in \Rep(\cG)^{\otimes I}$, take $Y=X^I$, and take in \eqref{e:Ran version of the Hecke action Y} 
$$\CS:=(\Delta_{X^I})_*(\omega_{X^I})\in \Shv(X^I\times X^I), \quad \CV_Y:=\on{ins}_{I\to I}(V\otimes \CS).$$

Then the isomorphism \eqref{e:Hecke property Ran Y} for the above $\CV_Y$ amounts to
the isomorphism \eqref{e:Hecke property usual}. 

\sssec{}

Vice versa, let $\CF$ be an object of $\Shv(\Bun_G)\otimes \QCoh(\CZ)$ equipped with a system of Hecke isomorphisms \eqref{e:Hecke property usual}. 
We recover the isomorphism \eqref{e:Hecke property Ran Bun} for
$$\CV=\on{ins}_{I \overset{\psi}\to J}(V \otimes \CS), \quad V \in \Rep(\cG)^{\otimes I},\,\, \CS\in \Shv(X^J),$$
by applying the functor
$$\left((p_{\Bun_G})_*\circ \left((\on{id}_X\times \Delta_\psi)^!(-)\sotimes \CS\right)\right) \otimes \on{Id}_{\QCoh(\CZ)}$$
to the two sides of \eqref{e:Hecke property usual}. 

\sssec{}

All of the above discussion is equally applicable when instead of $\Shv(-)$ we consider 
$\Dmod(-)$ (and $f$ is a map $\CZ\to \LocSys^\dr_\cG(X)$ or $\Shv^{\on{all}}(-)$ and 
(and $f$ is a map $\CZ\to \LocSys^{\on{Betti}}_\cG(X)$). 

\ssec{Creating Hecke eigensheaves}

Beilinson's projector is a functor that manifactures Hecke eigen-objects from arbitrary
objects of $\Shv(\Bun_G)$. In this subsection we will define it as the left adjoint of the
forgetful functor. But its crucial feature is that it can also be constructed as an explicit
integral Hecke functor. 

\sssec{} \label{sss:data for projector}

Let $f:\CZ\to \LocSys^{\on{restr}}_\cG(X)$ be as above. Recall that we denote by 
$\oblv_{\on{Hecke},\CZ}$ the forgetful functor 
\begin{equation} \label{e:forget Hecke Bun} 
\on{Hecke}(\CZ,\Shv(\Bun_G))\to
\Shv(\Bun_G) \otimes \QCoh(\CZ),
\end{equation} 
and by $\ind_{\on{Hecke},\CZ}$ its left adjoint
$$\Shv(\Bun_G) \otimes \QCoh(\CZ)\to \on{Hecke}(\CZ,\Shv(\Bun_G)).$$

\medskip

Recall also that the resulting monad on $\Shv(\Bun_G) \otimes \QCoh(\CZ)$
is given by the action of the (commutative) algebra object
$$\sR_\CZ\in \Rep(\cG)_\Ran\otimes \QCoh(\CZ),$$
see \secref{sss:R object new}. 

\sssec{} \label{sss:PZ}

Let $\oblv_{\on{Hecke}}$ denote the (not necessarily continuous) functor equal to the composition
$$\on{Hecke}(\CZ,\Shv(\Bun_G))\overset{\oblv_{\on{Hecke},\CZ}}\longrightarrow 
\Shv(\Bun_G)\otimes \QCoh(\CZ)\to \Shv(\Bun_G),$$
where the second arrow is the right adjoint to
\begin{equation} \label{e:ten OZ}
\Shv(\Bun_G) \overset{-\otimes \CO_\CZ}\to \Shv(\Bun_G) \otimes \QCoh(\CZ).
\end{equation} 

Let $\sP_\CZ$ denote the functor 
$$\Shv(\Bun_G)\to \Shv(\Bun_G) \otimes \QCoh(\CZ)$$
equal to the composition of \eqref{e:ten OZ} and the functor given by the action of $\sR_\CZ$.

\medskip

We obtain that the functor $\sP_\CZ$ naturally upgrades to a functor
$$\sP^{\on{enh}}_\CZ:\Shv(\Bun_G)\to \on{Hecke}(\CZ,\Shv(\Bun_G)),$$
which provides a left adjoint to $\oblv_{\on{Hecke}}$, i.e., it identifies with
$$\Shv(\Bun_G)\overset{-\otimes \CO_\CZ}\longrightarrow \Shv(\Bun_G) \otimes \QCoh(\CZ)
\overset{\ind_{\on{Hecke},\CZ}}\longrightarrow \on{Hecke}(\CZ,\Shv(\Bun_G)).$$

\medskip

We will refer to $\sP^{\on{enh}}_\CZ$ as \emph{Beilinson's spectral projector}. 

\sssec{} \label{sss:P Z as colim}

Note that the functor $\sP_\CZ$ can written down explicitly as a colimit, via the presentation
of $\sR_\CZ$ as a colimit, see \secref{sss:R as colim}.

\medskip

Namely, 
$$\sP_\CZ \simeq 
\underset{(I\overset{\psi}\to J)\in \on{TwArr}(\fSet)}{\on{colim}}\, \sP_\CZ^{I \overset{\psi}\to J},$$
where $\sP_\CZ^{I \overset{\psi}\to J}$ equals the composition 
\begin{multline*} 
\Shv(\Bun_G) \overset{\sR_{\Rep(\cG)}^{\otimes I}}\longrightarrow 
(\Rep(\cG)\otimes \Rep(\cG))^{\otimes I}\otimes \Shv(\Bun_G) \overset{\on{mult}_{\Rep(\cG)\otimes \Rep(\cG)}^\psi\otimes \on{Id}}\longrightarrow \\
\to (\Rep(\cG)\otimes \Rep(\cG))^{\otimes J}\otimes \Shv(\Bun_G)\simeq 
\Rep(\cG)^{\otimes J}\otimes \Shv(\Bun_G) \otimes \Rep(\cG)^{\otimes J}
\overset{\on{H}(-,-)\otimes \on{Id}}\longrightarrow \\
\to \Shv(\Bun_G\times X^J) \otimes \Rep(\cG)^{\otimes J} \overset{\on{Id}\otimes \sF^J}\longrightarrow 
\Shv(\Bun_G\times X^J) \otimes \qLisse(X)^{\otimes J} \otimes \QCoh(\CZ) \to \\
\to \Shv(\Bun_G\times X^J\times X^J) \otimes \QCoh(\CZ) \overset{(\on{id}\times \Delta_{X^J})^!\otimes \on{Id}}\longrightarrow \\
\to \Shv(\Bun_G\times X^J) \otimes \QCoh(\CZ) \overset{(p_{\Bun_G})_*\otimes \on{Id}}\longrightarrow \Shv(\Bun_G) \otimes \QCoh(\CZ), 
  \end{multline*} 
where:

\begin{itemize}

\item $\sR_{\Rep(\cG)}$ denotes the regular representation of $\cG$, regarded as an object of 
$\Rep(\cG)\otimes \Rep(\cG)$; 

\medskip

\item $\sF^J$ denotes the functor
$$\Rep(\cG)^{\otimes J}\to \QCoh(\CZ)\otimes \qLisse(X)^{\otimes J}$$ 
corresponding to the given map $f:\CZ\to \LocSys_\cG^{\on{restr}}(X)$;

\end{itemize}


%

\begin{rem} \label{r:projector dr and Betti new}
The above observations are equally applicable when instead of $\Shv(-)$ we consider 
$\Dmod(-)$ (and $f$ is a map $\CZ\to \LocSys^\dr_\cG(X)$ or $\Shv^{\on{all}}(-)$ and 
(and $f$ is a map $\CZ\to \LocSys^{\on{Betti}}_\cG(X)$). 

\end{rem}

\ssec{Beilinson's projector and nilpotence of singular support} \label{ss:Beil and Nilp}\label{ss:the projector new}

\sssec{} \label{sss:Ran and Nilp}

Consider the category $\Shv_\Nilp(\Bun_G)\subset \Shv(\Bun_G)$. It follows from \secref{sss:Hecke action on Nilp}
that the action of $\Rep(\cG)_\Ran$ on $\Shv(\Bun_G)$ preserves the subcategory $\Shv_\Nilp(\Bun_G)$.

\medskip

Thus, we can consider $\Shv_\Nilp(\Bun_G)$ as a $\Rep(\cG)_\Ran$-module category. 

\sssec{}

Let $\Shv(\Bun_G)^{\on{spec}}$ be as in \secref{sss:spectral subcateg}. Since $\Shv(\Bun_G)$
is dualizable as a DG category, the corresponding functor
$$\iota_{\Bun_G}:\Shv(\Bun_G)^{\on{spec}}\to \Shv(\Bun_G)$$
is fully faithful, by \propref{p:iota abs is ff}. 

\sssec{}

We claim:

\begin{prop} \label{p:Nilp as spec}
The subcategories
$$\Shv(\Bun_G)^{\on{spec}} \text{ and } \Shv_\Nilp(\Bun_G)$$
of $\Shv(\Bun_G)$ coincide.
\end{prop}

\begin{proof}

We have a commutative diagram
$$
\CD
\Shv_\Nilp(\Bun_G)^{\on{spec}}  @>>> \Shv(\Bun_G)^{\on{spec}}  \\
@V{\iota_{\Shv_\Nilp(\Bun_G)}}VV @VV{\iota_{\Bun_G}}V \\
\Shv_\Nilp(\Bun_G) @>>>  \Shv(\Bun_G),
\endCD
$$
where the left vertical arrow is an equivalence by \secref{sss:spectral subcateg LocSys}.
This implies the inclusion
$$\Shv_\Nilp(\Bun_G)\subset \Shv(\Bun_G)^{\on{spec}}.$$

\medskip

For the opposite inclusion, by \propref{p:M spec}, it suffices to show that the inclusion
$$\on{Hecke}(\LocSys^{\on{restr}}_\cG(X),\Shv_\Nilp(\Bun_G))\hookrightarrow
\on{Hecke}(\LocSys^{\on{restr}}_\cG(X),\Shv(\Bun_G))$$
is an equality.

\medskip

I.e., we have to show that the forgetful functor
$$\on{Hecke}(\LocSys^{\on{restr}}_\cG(X),\Shv(\Bun_G)) \to \Shv(\Bun_G)\otimes \QCoh(\LocSys^{\on{restr}}_\cG(X))$$
has essential image contained in
$$\Shv_\Nilp(\Bun_G)\otimes \QCoh(\LocSys^{\on{restr}}_\cG(X)) \subset \Shv(\Bun_G)\otimes \QCoh(\LocSys^{\on{restr}}_\cG(X)).$$

We will show that for any prestack $\CZ$ such that $\QCoh(\CZ)$ is dualizable, the forgetful functor 
$$\oblv_{\on{Hecke},\CZ}:\on{Hecke}(\CZ,\Shv(\Bun_G))\to \Shv(\Bun_G)\otimes \QCoh(\CZ)$$
is contained in $\Shv_\Nilp(\Bun_G)\otimes \QCoh(\CZ)$. 

\medskip

Since $\QCoh(\CZ)$ is dualizable, we have to show that
if $\sT$ is a continuous functor $\QCoh(\CZ)\to \Vect_\sfe$, the essential image of the composition
$$\on{Hecke}(\CZ,\Shv(\Bun_G))\to
\Shv(\Bun_G) \otimes \QCoh(\CZ)\overset{\on{Id}_{\Shv(\Bun_G)}\otimes \sT}\longrightarrow \Shv(\Bun_G)$$
is contained in $\Shv_\Nilp(\Bun_G)$.

\medskip

Let $\CF$ be an object of $\Shv_\Nilp(\Bun_G)\otimes \QCoh(\CZ)$ that can be upgraded to an object of 
the category $\on{Hecke}(\CZ,\Shv(\Bun_G))$. We have to show that
$$(\on{Id}_{\Shv(\Bun_G)}\otimes \sT)(\CF)\in \Shv_\Nilp(\Bun_G).$$

We will show that $(\on{Id}_{\Shv(\Bun_G)}\otimes \sT)(\CF) \in \Shv(\Bun_G)^{\on{Hecke-lisse}}$.
This would imply the required assertion by \thmref{t:lisse prel}.

\medskip

By \secref{sss:usual Hecke one dir}, for $V\in \Rep(\cG)$, we have: 
$$\on{H}(V,\CF) \simeq \CF\underset{\CO_\CZ}\boxtimes \sF(V)$$
as objects of
$$\Shv(\Bun_G\times X)\otimes \QCoh(\CZ).$$

Hence, 
$$\on{H}\left(V,(\on{Id}_{\Shv(\Bun_G)}\otimes \sT)(\CF)\right) \simeq
(\on{Id}_{\Shv(\Bun_G\times X)}\otimes \sT)(\CF\underset{\CO_\CZ}\boxtimes \sF(V))$$
as objects of $\Shv(\Bun_G\times X)$. 

\medskip

Now, the functor $\on{Id}_{\Shv(\Bun_G\times X)}\otimes \sT$ maps objects in the essential
image of the functor
$$\left(\Shv_\Nilp(\Bun_G)\otimes \QCoh(\CZ)\right)\otimes
\left(\QCoh(\CZ)\otimes \qLisse(X)\right) \overset{\underset{\CO_\CZ}\boxtimes}\longrightarrow \Shv(\Bun_G\times X)\otimes \QCoh(\CZ)$$
to $\Shv(\Bun_G)\otimes \qLisse(X)$. Hence, we conclude that
$$\on{H}\left(V,(\on{Id}_{\Shv(\Bun_G)}\otimes \sT)(\CF)\right) \in \Shv(\Bun_G)\otimes \qLisse(X),$$
as required. 

\end{proof}

\sssec{} \label{sss:the projector new}

Let $\sP$ be the functor from \secref{p:functor P}, i.e.,
$$\sP=\sP^{\on{enh}}_{\LocSys^{\on{restr}}_\cG(X)}.$$

Taking into account \propref{p:Nilp as spec}, we can view
it as a functor
\begin{equation} \label{e:P new}
\Shv(\Bun_G)\to \Shv_\Nilp(\Bun_G).
\end{equation}

By \propref{p:the projector univ}, the functor $\sP$ provides a left inverse to the embedding
$$\iota:\Shv_\Nilp(\Bun_G)\hookrightarrow \Shv(\Bun_G).$$

\begin{rem}

Unless a confusion is likely to occur, we will sometimes view $\sP$ as an endofunctor of $\Shv(\Bun_G)$,
by composing \eqref{e:P new} with the embedding $\iota$.

\medskip

When viewed as such, $\sP$ is a projector onto the full subcategory 
$\Shv_\Nilp(\Bun_G)\subset \Shv(\Bun_G)$

\end{rem}

\sssec{}

By \corref{c:the projector univ}, when we view $\sP$ as an endofunctor of $\Shv(\Bun_G)$,
it is given by the action of the object 
$$\sR\in \Rep(\cG)_\Ran,$$
see \secref{sss:object R}. 

\sssec{} \label{sss:proj int Hecke}

The expression for $\sR$ given by \eqref{e:formula for R univ} implies that the endofunctor $\sP$ 
of $\Shv(\Bun_G)$ 
is an \emph{integral Hecke functor}, i.e., a colimit
of functors, each of which is the composition of a Hecke functor
$$\on{H}(V,-):\Shv(\Bun_G)\to \Shv(\Bun_G\times X^I), \quad V\in \Rep(\cG)^{\otimes I}$$
and a functor $\Shv(\Bun_G\times X^I)\to \Shv(\Bun_G)$ given by 
\begin{equation} \label{e:int Hecke} 
\CF\mapsto (p_1)_*(\CF\sotimes p_2^!(\CS)), \quad \CS\in \Shv(X^I).
\end{equation} 

Furthermore, in the colimit expression for $\sP$, the above
objects $\CS$ belong to
$$\qLisse(X)^{\otimes I}\subset \Shv(X^I).$$

\ssec{Implications for Hecke eigensheaves} \label{ss:impl for eigen}

\sssec{}

Let $\CZ$ be a prestack over $\sfe$, and 
fix a map $f:\CZ\to  \LocSys^{\on{restr}}_\cG(X)$. Consider the resulting category 
$\on{Hecke}(\CZ,\Shv(\Bun_G))$, equipped with the pair of adjoint functors
$$\ind_{\on{Hecke},\CZ}:\Shv(\Bun_G)\otimes \QCoh(\CZ) \rightleftarrows \on{Hecke}(\CZ,\Shv(\Bun_G)):\oblv_{\on{Hecke},\CZ}.$$

\sssec{}

Combining \propref{p:Nilp as spec} with 
Sects. \ref{sss:Hecke as base change old}-\ref{sss:Hecke as base change cont}, we obtain:

\begin{cor} \label{c:right adj to projector} \hfill

\smallskip

\noindent{\em(a)} There exists a unique identification 
$$\on{Hecke}(\CZ,\Shv(\Bun_G)) \simeq \Shv_\Nilp(\Bun_G)\underset{\QCoh(\LocSys^{\on{restr}}_\cG(X))}\otimes  \QCoh(\CZ)$$
so that the forgetful functor 
$$\oblv_{\on{Hecke},\CZ}: \on{Hecke}(\CZ,\Shv(\Bun_G))\to \Shv(\Bun_G)\otimes \QCoh(\CZ)$$
identifies with 
\begin{multline} \label{e:right adj to projector}
\Shv_\Nilp(\Bun_G)\underset{\QCoh(\LocSys^{\on{restr}}_\cG(X))}\otimes  \QCoh(\CZ) \simeq \\
\simeq 
(\Shv_\Nilp(\Bun_G) \otimes \QCoh(\CZ)) \underset{\QCoh(\LocSys^{\on{restr}}_\sG(X))\otimes \QCoh(\LocSys^{\on{restr}}_\cG(X))}\otimes
\QCoh(\LocSys^{\on{restr}}_\sG(X)) \to \\
\overset{\on{Id}\otimes (\Delta_{\LocSys^{\on{restr}}_\sG(X)})_*}\longrightarrow 
\Shv_\Nilp(\Bun_G) \otimes \QCoh(\CZ)\overset{\iota\otimes \on{Id}}\longrightarrow  \Shv(\Bun_G) \otimes \QCoh(\CZ).
\end{multline} 

\smallskip

\noindent{\em(b)} Assume moreover that $\CO_\CZ\in \QCoh(\CZ)$ is compact. Then with respect to the identification 
of point (a), the forgetful functor 
$$\oblv_{\on{Hecke}}: \on{Hecke}(\CZ,\Shv(\Bun_G))\to \Shv(\Bun_G)$$
identifies with 
\begin{multline} \label{e:right adj to projector cont}
\Shv_\Nilp(\Bun_G)\underset{\QCoh(\LocSys^{\on{restr}}_\cG(X))}\otimes  \QCoh(\CZ)
\overset{\on{Id}\otimes f_*}\longrightarrow \\
\to \Shv_\Nilp(\Bun_G)\underset{\QCoh(\LocSys^{\on{restr}}_\cG(X))}\otimes \QCoh(\LocSys^{\on{restr}}_\cG(X))
\simeq \Shv_\Nilp(\Bun_G)\overset{\iota}\hookrightarrow \Shv(\Bun_G).
\end{multline}

\smallskip

\noindent{\em(c)} The functor $\sP^{\on{enh}}_\CZ$ identifies with
\begin{multline} \label{e:P enh expl}
\Shv(\Bun_G) \overset{\sP}\to \Shv_\Nilp(\Bun_G) \simeq 
 \Shv_\Nilp(\Bun_G)\underset{\QCoh(\LocSys^{\on{restr}}_\cG(X))}\otimes \QCoh(\LocSys^{\on{restr}}_\cG(X)) 
\overset{\on{Id}\otimes f^*}\longrightarrow \\
\to \Shv_\Nilp(\Bun_G)\underset{\QCoh(\LocSys^{\on{restr}}_\cG(X))}\otimes \QCoh(\CZ).
\end{multline}

\end{cor}

In particular, we obtain:

\begin{cor} \label{c:P Z}
The functor $\sP_\CZ$ takes values in
$$\Shv_\Nilp(\Bun_G)\otimes \QCoh(\CZ)
\subset \Shv(\Bun_G)\otimes \QCoh(\CZ).$$
\end{cor}

\begin{cor} \label{c:right adj to projector bis}
Under the assumption of \corref{c:right adj to projector}(b), the functor \eqref{e:right adj to projector cont}
admits a left adjoint, given by \eqref{e:P enh expl}.  
\end{cor}

\begin{rem}
Note that the functor $\iota: \Shv_\Nilp(\Bun_G)\hookrightarrow \Shv(\Bun_G)$ itself does does 
\emph{not} admit a left adjoint (see \secref{ss:iota left} for more details). This does not violate 
\corref{c:right adj to projector bis} since $\CO_{\LocSys^{\on{restr}}_\cG(X)}$ is not compact as an object
of $\QCoh(\LocSys^{\on{restr}}_\cG(X))$.
\end{rem} 

%

%
%

\sssec{} \label{sss:P Z O comp}

Furthermore, under the assumption of \corref{c:right adj to projector}(b), we have
\begin{equation} \label{e:P and P enh}
(\on{Id}\otimes \Gamma(\CZ,-))\circ \sP_\CZ \simeq \iota \circ (\on{Id}\otimes f_*)\circ \sP^{\on{enh}}_\CZ,
\end{equation}
where $\on{Id}\otimes f_*$ is the functor
$$\Shv_\Nilp(\Bun_G)\underset{\QCoh(\LocSys^{\on{restr}}_\cG(X))}\otimes  \QCoh(\CZ)\to$$
$$\to \Shv_\Nilp(\Bun_G)\underset{\QCoh(\LocSys^{\on{restr}}_\cG(X))}\otimes \QCoh(\LocSys^{\on{restr}}_\cG(X))
\simeq \Shv_\Nilp(\Bun_G).$$

Note that the left-hand side in \eqref{e:P and P enh} can also be identified with the functor
\begin{equation} \label{e:P and P enh bis}
\CV\star -, \quad \CV:=\left(\on{Id}_{\Rep(\cG)_\Ran}\otimes \Gamma(\CZ,-)\right)(\sR_{\CZ})\in \Rep(\cG)_\Ran.
\end{equation}

%

\begin{rem} \label{r:proj Nilp dr}

The material in this subsection applies when instead of $\Shv(-)$ we consider $\Dmod(-)$,
but $f$ is a map $f:\CZ\to \LocSys^{\on{restr}}_\cG(X)\subset \LocSys^{\dr}_\cG(X)$. Here we use the fact that the inclusion
\eqref{e:incl Nilp Dmod} is an equality. 

\end{rem}

\begin{rem} \label{r:proj Nilp Betti}
This is a preview of the material in  \secref{s:spectral Betti}:

\medskip

The discussion in this subsection applies when instead of $\Shv(-)$ we consider $\Shv^{\on{all}}(-)$.
In fact, here we have two variants:

\medskip

We can consider maps $f:\CZ\to \LocSys^{\on{restr}}_\cG(X)$, in which case the statements 
of the assertions in this subsection hold with no modifications.

\medskip

However, we can also consider maps $f:\CZ\to \LocSys^{\on{Betti}}_\cG(X)$. In this 
case, one has to replace $\Shv_\Nilp(\Bun_G)$ by $\Shv^{\on{all}}_\Nilp(\Bun_G)$.

\end{rem}

\section{Applications of Beilinson's spectral projector} \label{s:projector}

In this section we will combine our Theorems \ref{t:spectral decomp} and \ref{t:lisse} 
with Beilinson's spectral projector to prove some key theorems
about $\Shv_\Nilp(\Bun_G)$.

\ssec{Compact generation of $\Shv_{\on{Nilp}}(\Bun_G)$}

In this subsection we will use Beilinson's projector to prove the following assertion:

\begin{thm} \label{t:Nilp comp gen abs}
The category $\Shv_{\on{Nilp}}(\Bun_G)$ is compactly generated.
\end{thm}

\sssec{} \label{sss:proof of Nilp comp gen abs}

By Sects. \ref{sss:almost lift action a}-\ref{sss:almost lift action c}, for every connected component $\CZ$ of 
$\LocSys_\cG^{\on{restr}}(X)$, we can find a family of algebraic stacks mapping to $\CZ$
$$f_n:Z_n\to \CZ$$
with the following properties:

\begin{itemize}

\item Each $Z_n$ is of the form $S/\sH$ with $S$ an affine scheme and $\sH$ an algebraic group;

\smallskip

\item Each $f_n$ is affine (so that $(f_n)_*$ is continuous and $\QCoh(\LocSys_\cG^{\on{restr}}(X))$-linear),
and $(f_n)_*$ admits a continuous right adjoint $(f_n)^!$, which is also 
$\QCoh(\LocSys_\cG^{\on{restr}}(X))$-linear. 

\smallskip

\item The essential images of the functors $(f_n)_*$ generate $\QCoh(\LocSys_\cG^{\on{restr}}(X))$. 

\end{itemize} 

\begin{rem}
In what follows, in order to unburden the notations, we will group the connected components
of $\LocSys_\cG^{\on{restr}}(X)$ together, and consider the stacks $Z_n$ as mapping directly to  
$\LocSys_\cG^{\on{restr}}(X)$. So in the sequel, in this context, the index $n$ no longer refers to 
a natural number but rather to an element $\underset{\CZ}\cup\, \BN$.

\end{rem}

\sssec{} \label{sss:comp gen 1}

The adjunction 
$$(f_n)_*:\QCoh(Z_n) \rightleftarrows \QCoh(\LocSys_\cG^{\on{restr}}(X))):(f_n)^!$$
as $\QCoh(\LocSys_\cG^{\on{restr}}(X))$-module categories induces an adjunction 
$$\on{Id} \otimes (f_n)_*:
\Shv_{\on{Nilp}}(\Bun_G)\underset{\QCoh(\LocSys_\cG^{\on{restr}}(X))}\otimes\QCoh(Z_n)
\rightleftarrows \Shv_{\on{Nilp}}(\Bun_G):\on{Id}\otimes f_n^!.$$

\medskip

In particular, the functor 
\begin{equation} \label{e:from Zn prel}
\Shv_{\on{Nilp}}(\Bun_G)\underset{\QCoh(\LocSys_\cG^{\on{restr}}(X))}\otimes\QCoh(Z_n) 
\overset{\on{Id} \otimes (f_n)_*}\longrightarrow \Shv_{\on{Nilp}}(\Bun_G)
\end{equation} 
preserves compactness. Furthermore,
since the essential images of the functors $(f_n)_*$ generate $\QCoh(\LocSys_\cG^{\on{restr}}(X))$,
the essential images of the functors \eqref{e:from Zn prel} generate $\Shv_{\on{Nilp}}(\Bun_G)$.

%
%
%
%

\medskip

Hence, it is enough to show that each of the categories
$$\Shv_{\on{Nilp}}(\Bun_G)\underset{\QCoh(\LocSys_\cG^{\on{restr}}(X))}\otimes\QCoh(Z_n)$$
is compactly generated.

\sssec{} \label{sss:comp gen 2}

Consider the functor 
\begin{equation} \label{e:from Zn}
\Shv_{\on{Nilp}}(\Bun_G)\underset{\QCoh(\LocSys_\cG^{\on{restr}}(X))}\otimes\QCoh(Z_n) 
\overset{\text{\eqref{e:from Zn prel}}}\longrightarrow \Shv_{\on{Nilp}}(\Bun_G)\overset{\iota}\to \Shv(\Bun_G).
\end{equation} 

Since $\Shv(\Bun_G)$ is compactly generated, it suffices to show that this functor is conservative and admits a left adjoint.

\medskip

The existence of the left adjoint follows from \corref{c:right adj to projector bis} (the conditions of the
corollary are guaranteed by the assumption that $Z_n$ is of the form $S/\sH$). 

\sssec{} \label{sss:comp gen 3}

To prove that \eqref{e:from Zn} is conservative we argue as follows. It suffices to show that the functor \eqref{e:from Zn prel} is conservative.
The latter is equivalent to the fact that the essential image of the functor
\begin{multline*}
\Shv_{\on{Nilp}}(\Bun_G)\simeq 
\Shv_{\on{Nilp}}(\Bun_G)  \underset{\QCoh(\LocSys_\cG^{\on{restr}}(X))}\otimes \QCoh(\LocSys_\cG^{\on{restr}}(X)) 
\overset{\on{Id}\otimes (f_n)^*}\longrightarrow \\
\to \Shv_{\on{Nilp}}(\Bun_G)\underset{\QCoh(\LocSys_\cG^{\on{restr}}(X))}\otimes\QCoh(Z_n)
\end{multline*}
generates $\Shv_{\on{Nilp}}(\Bun_G)\underset{\QCoh(\LocSys_\cG^{\on{restr}}(X))}\otimes\QCoh(Z_n)$.

\medskip

To prove this, it is sufficient to show that the essential image of the functor
$$f_n^*:\QCoh(\LocSys_\cG^{\on{restr}}(X))\to \QCoh(Z_n)$$
generates $\QCoh(Z_n)$. This is equivalent to the functor $(f_n)_*$ being conservative. But this
is indeed the case since $f_n$ is affine. 

\qed[\thmref{t:Nilp comp gen abs}]

\sssec{}

Note that \thmref{t:Nilp comp gen abs} admits the following corollary:

\begin{cor} \label{c:Nilp comp gen U}
Let $\CU\overset{j}\hookrightarrow \Bun_G$ be an open substack such that
the functor $j_!$ (equivalently, $j_*$) sends $\Shv_{\on{Nilp}}(\CU)^{\on{constr}}$ to
$\Shv_{\on{Nilp}}(\Bun_G)^{\on{constr}}$. Then $\Shv_{\on{Nilp}}(\CU)$ is
compactly generated.
\end{cor}

\begin{proof}

Follows from the fact that the functor 
$$j^*:\Shv_{\on{Nilp}}(\Bun_G)\to \Shv_{\on{Nilp}}(\CU)$$
admits a conservative right adjoint, given by $j_*$ (see Remark \ref{r:preserve Nilp Sing Supp}). 

\end{proof} 

\ssec{A set of generators for $\Shv_{\on{Nilp}}(\Bun_G)$} \label{ss:gen Nilp}

In the previous subsection we showed that the category $\Shv_{\on{Nilp}}(\Bun_G)$
is compactly generated. In this subsection, we will sharpen this by writing down an explicit
set of generators for $\Shv_{\on{Nilp}}(\Bun_G)$. 

\sssec{}  \label{sss:indicators for N}

Let $\CY$ be an algebraic stack and $\CN\subset T^*(\CY)$ a half-dimensional conical subset. We claim that
we can find a collection of $k$-points $y_i$ on $\CY$ (finitely many in every quasi-compact open of $\CY$)
such that for every $0\neq \CF\in \Shv_\CN(\CY)$, the !-fiber of $\CF$ at least one $y_i$ will be non-zero.

\medskip

With no restriction of generality, we can assume that $\CY$ is a (quasi-compact) smooth scheme. 
We can partition $\CY$ into smooth, connected locally closed subschemes $\CY_i$ such that the 
dimensions of the fibers of the map 
\begin{equation} \label{e:cotangent fibers proj}
\CN \hookrightarrow T^*(\CY)\to \CY
\end{equation} 
are $\leq \on{codim}(\CY_i,\CY)$ over $\CY_i$. By refining the partition, we can assume that for 
every $n$, the union 
$$\underset{i,\on{codim}(\CY_i,\CY)\geq n}\cup\, \CY_i$$
is closed. 

\medskip

Choose a point $y_i\in \CY_i(k)$. We claim that these points will have the required property. 

\begin{proof} 

For a given $0\neq \CF\in  \Shv_\CN(\CY)$, let $i$ be an index with a minimal 
$\on{codim}(\CY_i,\CY)$, such that $\CF|_{\CY_i}\neq 0$.  It suffices to show that 
$\CF|_{\CY_i}$ is lisse. I.e., we want to show that
$$\on{SingSupp}(\CF|_{\CY_i})\subset T^*(\CY_i)$$
is the zero section.

\medskip

Removing the (closed) subscheme 
$$\underset{j,\on{codim}(\CY_j,\CY)> \on{codim}(\CY_i,\CY)}\cup\, \CY_j,$$
we can assume that $\CY_i$ is a closed subscheme of $\CY$, and $\CF$ is
supported on $\CY_i$. Denote by $\pi_i$ the projection
$$T^*(\CY)|_{\CY_i}\to T^*(\CY_i).$$
This is a smooth surjective map, with fibers of dimension $\on{codim}(\CY_i,\CY)$.

\medskip

We have
$$\on{SingSupp}(\CF)=\pi_i^{-1}(\on{SingSupp}(\CF|_{\CY_i})).$$

Hence, the dimension of the fibers of the map 
$$\on{SingSupp}(\CF) \hookrightarrow T^*(\CY)\to \CY$$
over $\CY_i$ equals $\on{codim}(\CY_i,\CY)$ plus the dimension of the fibers of the map 
\begin{equation} \label{e:cotangent fibers proj i}
\on{SingSupp}(\CF|_{\CY_i}) \hookrightarrow T^*(\CY_i)\to \CY_i.
\end{equation} 
Hence, since $\on{SingSupp}(\CF)\subset \CN$ and by the assumption on \eqref{e:cotangent fibers proj},
we obtain that the fibers of the map \eqref{e:cotangent fibers proj i} are zero-dimensional.

\medskip

Since $\on{SingSupp}(\CF|_{\CY_i})$ is conical, we obtain that it is necessarily the zero section.

\end{proof} 

\sssec{} \label{sss:generators Nilp}

We are now going to exhibit a particular set of generators for the category $\Shv_{\on{Nilp}}(\Bun_G)$.

\medskip

Recall that $\on{Nilp}\subset T^*(\Bun_G)$ is half-dimensional, see \secref{s:glob Nilp}.
For every $i$, let $\delta_{y_i}\in \Shv(\Bun_G)$ be the corresponding !-delta function object, i.e.,
$$\delta_{y_i}=(\bi_{y_i})_!(\sfe),$$
where
$$\Spec(k) \overset{\bi_{y_i}}\to \Bun_G$$
is the morphism corresponding to $y_i$. 

\medskip

Let 
$$f_n:Z_n\to \LocSys^{\on{restr}}_\cG(X)$$
be one of the substacks as in \secref{sss:proof of Nilp comp gen abs}. 

\medskip

Consider the objects
$$\sP_{Z_n}^{\on{enh}}(\delta_{y_i})\in \Shv_{\on{Nilp}}(\Bun_G)\underset{\QCoh(\LocSys^{\on{restr}}_\cG(X))}\otimes \QCoh(Z_n),$$
and
\begin{equation} \label{e:generators}
(\on{Id} \otimes (f_n)_*)(\sP_{Z_n}^{\on{enh}}(\delta_{y_i}))\in \Shv_{\on{Nilp}}(\Bun_G).
\end{equation} 


We claim that the objects \eqref{e:generators} provide a set of compact generators for $\Shv_{\on{Nilp}}(\Bun_G)$. 

\sssec{}

Since the functors \eqref{e:from Zn prel} preserve compactness and the union of their essential images generates
$\Shv_{\on{Nilp}}(\Bun_G)$, it suffices to show that for a fixed $n$, the objects $\sP_{Z_n}^{\on{enh}}(\delta_{y_i})$ generate 
$$\Shv_{\on{Nilp}}(\Bun_G)\underset{\QCoh(\LocSys_\cG^{\on{restr}}(X))}\otimes\QCoh(Z_n).$$

\medskip

By adjunction, for $\CF\in \Shv_{\on{Nilp}}(\Bun_G)\underset{\QCoh(\LocSys_\cG^{\on{restr}}(X))}\otimes\QCoh(Z_n)$,
we have
$$\CHom_{\Shv_{\on{Nilp}}(\Bun_G)\underset{\QCoh(\LocSys_\cG^{\on{restr}}(X))}\otimes\QCoh(Z_n)}
(\sP_{Z_n}^{\on{enh}}(\delta_{y_i}),\CF)\simeq
\CHom_{\Shv(\Bun_G)}\left(\delta_{y_i},\iota\circ (\on{Id} \otimes (f_n)_*)(\CF)\right).$$

Since the functor
$$\on{Id} \otimes (f_n)_*:\Shv_{\on{Nilp}}(\Bun_G)\underset{\QCoh(\LocSys_\cG^{\on{restr}}(X))}\otimes\QCoh(Z_n)\to
\Shv_\Nilp(\Bun_G)$$
is conservative, the latter assertion is equivalent to the statement that if $\CF'\in \Shv_{\on{Nilp}}(\Bun_G)$ is \emph{non-zero}, 
then \emph{not all}  
$$\CHom_{\Shv(\Bun_G)}(\delta_{y_i},\CF')\simeq (\bi_{y_i})^!(\CF')$$
are zero. 

\medskip

However, the latter was proved in \secref{sss:indicators for N}. 

\ssec{The tensor product property} \label{ss:ten product}

\sssec{} \label{sss:when tensor product}

Let $Y_1$ and $Y_2$ be a pair of quasi-compact schemes (or algebraic stacks) over the ground field $k$. Let $\Shv(-)$
be a constructible sheaf theory (i.e., in the case of D-modules, we will consider the subcategory of
holonomic D-modules or regular holonomic D-modules).

\medskip

Consider the external tensor product functor
\begin{equation} \label{e:external ten prod}
\Shv(Y_1)\otimes \Shv(Y_2) \to \Shv(Y_1\times Y_2).
\end{equation} 

It is fully faithful (see \cite[Lemma A.2.6]{GKRV}), but very rarely an equivalence. However, it is an equivalence, for example, if 
either $Y_1$ or $Y_2$ is an algebraic stack with a \emph{finite number of isomorphism classes of $k$-points}\footnote{Such as, for example,
$B\backslash G/N$, or the stack $\Bun_G$ for $X$ of genus $0$.}. 

\medskip

In fact, there is a clear obstruction for an object of $\Shv(Y_1\times Y_2)$ to belong to the essential image 
of \eqref{e:external ten prod}. Namely all objects in the essential image have their singular support contained
in a subset of $T^*(Y_1\times Y_2)\simeq T^*(Y_1)\times T^*(Y_2)$ of the form
$$\CN_1\times \CN_2, \quad \CN_i\subset T^*(Y_i).$$

\medskip

Thus, one can wonder whether, given $\CN_1$ and $\CN_2$ as above, the functor
\begin{equation} \label{e:external ten prod sing supp}
\Shv_{\CN_1}(Y_1)\otimes \Shv_{\CN_2}(Y_2) \to \Shv_{\CN_1\times \CN_2}(Y_1\times Y_2)
\end{equation} 
is an equivalence.

\medskip

Now, this happens to always be the case for constructible sheaves in the Betti context, at least after
passing to the left completion of the left-hand side (assuming $\CN_i$ are Lagrangian). However, this is \emph{not} the case of 
$\ell$-adic sheaves over a field of positive characteristic, and not for holonomic (but irregular) D-modules. 

\medskip

For example,
taking $Y_1=Y_2=\BA^1$ and $\CN_1=\CN_2=\{0\}$, the map
$$\qLisse(\BA^1)\otimes \qLisse(\BA^1)\to \qLisse(\BA^1\times \BA^1)$$
is \emph{not} an equivalence. Indeed, the object
$$\on{mult}^*(\on{A-Sch})\in \qLisse(\BA^1\times \BA^1)$$
does not lie in the essential image, where 
$$\on{mult}:\BA^1\times \BA^1\to \BA^1$$
is the product map and $\on{A-Sch}\in \Shv(\BA^1)$ is the Artin-Schreier local system. 

\medskip

That said, Theorems 
\ref{t:product thm stack 1b} and \thmref{t:product thm stack 2}
say that the functor \eqref{e:external ten prod sing supp}
is an equivalence in the case when $Y_1$ is a proper scheme and $\CN_1=\{0\}$, either 
up to left completions or under an additional assumption on $Y_1$. However, an assertion of this sort 
would still fail even when $Y_1$ is proper for a more general $\CN_1$.

\sssec{}

The main result of the present subsection is the following:

\begin{thm} \label{t:tensor product}
Let $G_1$ and $G_2$ be a pair of reductive groups. 
Then the functor
\begin{equation} \label{e:ten prod BunG}
\Shv_{\on{Nilp}}(\Bun_{G_1})\otimes \Shv_{\on{Nilp}}(\Bun_{G_2})\to \Shv_\Nilp(\Bun_{G_1\times G_2})
\end{equation} 
is an equivalence. 
\end{thm}

The rest of this subsection is devoted to the proof of this theorem.

\sssec{}

We first show that \eqref{e:ten prod BunG} is fully faithful. 

\medskip

Indeed, a standard colimit argument 
shows that the fully faithfulness of \eqref{e:external ten prod} implies the fully faithfulness of 
\eqref{e:external ten prod sing supp}, whenever one of the categories $\Shv_{\CN_i}(Y_i)$
is dualizable. 

\medskip

The dualizability assumption holds for $\Shv_{\on{Nilp}}(\Bun_{G_i})$ by \thmref{t:Nilp comp gen abs}. 

\sssec{}

We now show that the essential image of \eqref{e:ten prod BunG} generates the target category.

\medskip

Let $\{y_{i_1}\}\in \Bun_{G_1}(k)$ and 
$\{y_{i_2}\}\in \Bun_{G_2}(k)$ be collections of points chosen as in \secref{sss:indicators for N} with respect
to the subsets $\Nilp\subset T^*(\Bun_{G_1})$ and $\Nilp\subset T^*(\Bun_{G_2})$,
respectively. Then, by construction, the collection of points 
$$\{y_{i_1}\times y_{i_2}\}\in \Bun_{G_1}(k)\times \Bun_{G_2}(k)\simeq \Bun_{G_1\times G_2}(k)$$
will have the corresponding property with respect to
$$\Nilp\subset T^*(\Bun_{G_1\times G_2}).$$

\medskip

Consider the corresponding stacks
$$Z_{n_1}\overset{f_{n_1}}\to \LocSys_{\cG_1}^{\on{restr}}(X) \text{ and }
Z_{n_2}\overset{f_{n_2}}\to \LocSys_{\cG_2}^{\on{restr}}(X)$$
and
$$Z_{n_1}\times Z_{n_2} \overset{f_{n_1}\times f_{n_2}}\to \LocSys_{\cG_1\times \cG_2}^{\on{restr}}(X).$$

By \secref{sss:generators Nilp}, it suffices to show that for all quadruples
$i_1,i_2$, $n_1,n_2$, the object
\begin{equation} \label{e:proj delta}
(\on{Id}\otimes (f_{n_1}\times f_{n_2})_*)(\sP_{Z_{n_1}\times Z_{n_2}}^{\on{enh}}(\delta_{y_{i_1}\times y_{i_2}}))\in
\Shv_{\on{Nilp}\times \on{Nilp}}(\Bun_{G_1\times G_2})
\end{equation}
lies in the essential image of the functor \eqref{e:ten prod BunG}.

\sssec{}

We claim the object \eqref{e:proj delta} equals the external tensor product
$$(\on{Id}\otimes (f_{n_1})_*)(\sP_{Z_{n_1}}^{\on{enh}}(\delta_{y_{i_1}}))\boxtimes 
(\on{Id}\otimes (f_{n_2})_*)(\sP_{Z_{n_2}}^{\on{enh}}(\delta_{y_{i_2}})).$$

Using \eqref{e:P and P enh}, in order to prove the latter assertion, it suffices to show that
the diagram
$$
\CD
\Shv(\Bun_{G_1})\otimes \Shv(\Bun_{G_2}) @>{\sP_{Z_1}\boxtimes \sP_{Z_2}}>>  
\Shv(\Bun_{G_1})\otimes \Shv(\Bun_{G_2}) \otimes \QCoh(Z_1)\otimes \QCoh(Z_2) \\
@VVV @VVV \\
\Shv(\Bun_{G_1\times G_2}) @>{\sP_{Z_1\times Z_2}}>> \Shv(\Bun_{G_1\times G_2})\otimes \QCoh(Z_1\times Z_2)
\endCD
$$
commutes.

\medskip

However, this follows from \propref{p:prod projectors}. 

\qed[\thmref{t:tensor product}]

\ssec{Some consequences pertaining to \conjref{c:Nilp comp gen}}

\sssec{}

Consider the tautological embedding
\begin{equation} \label{e:embed nilp}
\iota:\Shv_{\on{Nilp}}(\Bun_G) \hookrightarrow \Shv(\Bun_G).
\end{equation} 
 
\medskip

Now that we know that the category $\Shv_{\on{Nilp}}(\Bun_G)$ is compactly generated, we can equivalently reformulate
\conjref{c:Nilp comp gen} as follows:

\medskip

\noindent(i) The functor $\iota$ preserves compactness;

\medskip

\noindent(ii) The right adjoint of $\iota$ is continuous. 

\sssec{}

We claim: 

\begin{thm} \label{t:Nilp comp gen dR}
\conjref{c:Nilp comp gen} holds in the de Rham context.
\end{thm} 

\begin{proof}

By \corref{c:nilp via spectral}, the composite functor
$$\Shv_\Nilp(\Bun_G) \overset{\iota}\hookrightarrow \Shv(\Bun_G)\hookrightarrow \Dmod(\Bun_G)$$
admits a continuous right adjoint, where $\Shv(-)$ is the category of holonomic D-modules.

\medskip

In particular, the above composite functor sends compact objects in $\Shv_\Nilp(\Bun_G)$ to objects
that are compact in $\Dmod(\Bun_G)$. However, since the embedding 
$$\Shv(\Bun_G)\hookrightarrow \Dmod(\Bun_G)$$
is fully faithful, we obtain that $\iota$ preserves compactness as well.

\end{proof} 

\begin{rem}
Recall the functor 
$$\sP:\Shv(\Bun_G)\to \Shv_{\on{Nilp}}(\Bun_G),$$
see \secref{sss:the projector new}. In \propref{p:projector as adjoint} below 
we will show that $\sP$
\emph{would be} the right adjoint of $\iota$,
if we knew that that right adjoint was continuous. 
\end{rem} 

\sssec{} \label{sss:access BunG}

Recall (see \secref{sss:non qc N access}) that 
\begin{equation} \label{e:access on Bun_G}
\Shv_{\on{Nilp}}(\Bun_G)^{\on{access}}\subset \Shv_{\on{Nilp}}(\Bun_G)
\end{equation} 
denotes the full subcategory generated by the essential image of
$$\Shv_{\on{Nilp}}(\Bun_G)\cap \Shv(\Bun_G)^c \subset \Shv_{\on{Nilp}}(\Bun_G).$$

\medskip

The above definition of $\Shv_{\on{Nilp}}(\Bun_G)^{\on{access}}$ coincides with the one given in
\secref{sss:non qc N}. This is due to the combination of \thmref{t:preserve Nilp Sing Supp prel} and \corref{c:N preserved abs}. 

\medskip

Note that \conjref{c:Nilp comp gen} can be reformulated as the assertion that the inclusion
$$\Shv_{\on{Nilp}}(\Bun_G)^{\on{access}} \subset \Shv_{\on{Nilp}}(\Bun_G)$$
is an equality.

\sssec{}

We now claim:

\begin{prop} \label{p:cohom estim 1}
The following statements are equivalent:

\smallskip

\noindent{\em(a)} \conjref{c:Nilp comp gen} holds;

\smallskip

\noindent{\em(b)} The endofunctor 
$$\iota \circ (\on{Id}_{\Shv_\Nilp(\Bun_G)} \otimes \Gamma(Z_n,-))\circ \sP_{Z_n}\simeq \oblv_{\on{Hecke}}\circ \sP^{\on{enh}}_{Z_n}$$
of $\Shv(\Bun_G)$ for $Z_n$ being each of the stacks from \secref{sss:proof of Nilp comp gen abs},
has the following properties:

\smallskip

\noindent{\em(bi)} It preserves compactness;

\smallskip

\noindent{\em(bii)} It sends compact objects to objects that are comologically bounded;

\smallskip

\noindent{\em(biii)} It sends compact objects to objects that are comologically bounded on the left.

\end{prop}

\begin{proof}

By \corref{c:right adj to projector}(b), we can rewrite the functor $\oblv_{\on{Hecke}}\circ \sP^{\on{enh}}_{Z_n}$
as
\begin{multline} \label{e:endofunctor point b}
\Shv(\Bun_G) \overset{\sP^{\on{enh}}_{Z_n}}\to 
\Shv_\Nilp(\Bun_G) \underset{\QCoh(\LocSys_\cG^{\on{restr}}(X))}\otimes \QCoh(Z_n)\overset{\on{Id}\otimes (f_n)_*}\longrightarrow 
\Shv_\Nilp(\Bun_G)\overset{\iota}\to \\
\to \Shv(\Bun_G).
\end{multline}

\medskip

Point (a) implies (bi) because the endofunctor in (b) is the composition of $\iota$ (which preserves compactness by the
assumption in (a)), the functor $(f_n)_*\otimes \on{Id}$ (which preserves compactness by the assumption on $(Z_n,f_n)$),
and the functor $\sP^{\on{enh}}_{Z_n}$ (which preserves compactness, being a left adjoint). 

\medskip

Vice versa, point (bi) implies point (a), because the images of the compacts under
$$\Shv(\Bun_G) \overset{\sP^{\on{enh}}_{Z_n}}\to 
\Shv_\Nilp(\Bun_G) \underset{\QCoh(\LocSys_\cG^{\on{restr}}(X))}\otimes \QCoh(Z_n)\overset{\on{Id}\otimes (f_n)_*}\longrightarrow 
\Shv_\Nilp(\Bun_G)$$
generate $\Shv_\Nilp(\Bun_G)$, see Sects. \ref{sss:comp gen 2}-\ref{sss:comp gen 3}. 

\medskip

We have the tautological implications (bi) $\Rightarrow$ (bii) $\Rightarrow$ (biii). Suppose that (biii) holds, and let us deduce (a). 

\medskip

The embedding
\begin{equation} \label{e:access in all}
\Shv_\Nilp(\Bun_G)^{\on{access}}\hookrightarrow \Shv_\Nilp(\Bun_G)
\end{equation} 
induces an equivalence on bounded below (eventually coconnective) subcategories (see \secref{sss:non qc N access}). 
Point (a) is equivalent to this embedding being an equivalence, see \secref{sss:access BunG}. 

\medskip

The assumption in (biii) implies that the (compact) generators of $\Shv_\Nilp(\Bun_G)$ are bounded on the left
(=eventually coconnective), when considered as objects of $\Shv(\Bun_G)$; hence they belong to the essential
image of \eqref{e:access in all}. This implies that \eqref{e:access in all} is an equivalence. 

\end{proof} 

\begin{rem}

Note that although we cannot prove \conjref{c:Nilp comp gen} in general, and as result we do not know that 
the endofunctors \eqref{e:endofunctor point b} preserve compactness, we do know that these functors
send compact objects to objects that are bounded above and such that all of their individual cohomologies 
are constructible. This is because the functor $\iota$ sends objects that are compact in $\Shv_\Nilp(\Bun_G)$
to objects of $\Shv(\Bun_G)$ with these properties. 

\medskip

This finiteness property of individual cohomologies is non-obvious from the presentation of the functors $\sP_{Z_n}$ as colimits,
see \secref{sss:P Z as colim}.

\end{rem} 

\sssec{}

We will now prove:

\begin{thm} \label{t:Nilp comp gen Betti}
\conjref{c:Nilp comp gen} holds in the Betti context.
\end{thm} 

\begin{proof} 

Consider the endofunctor \eqref{e:endofunctor point b}. By \propref{p:cohom estim 1}, it suffices to show 
that it sends objects that are compact in $\Shv_\Nilp(\Bun_G)$ (for $\Shv(-)$ being the constructible 
category in the classical topology) to objects that are cohomologically bounded.

\medskip

We can assume that our field of coefficients $\sfe$ is $\BC$, and let us apply the Riemann-Hilbert equivalence.
I.e., we can replace the initial $\Shv(-)$ by the \emph{regular} holonomic category $\Shv^{\on{reg.hol}}(-)$.

\medskip

We now embed $\Shv^{\on{reg.hol}}(-)$ into the entire holonomic category $\Shv^{\on{hol}}(-)$, i.e.,
without the regularity assumption. By \thmref{t:Nilp comp gen dR}, we know that the functor 
\eqref{e:endofunctor point b} sends compact objects in $\Shv^{\on{hol}}(\Bun_G)$ to objects that are
cohomologically bounded. 

\medskip

Hence, it suffices to show that the embedding
$$\Shv^{\on{reg.hol}}(\Bun_G)\hookrightarrow \Shv^{\on{hol}}(\Bun_G)$$
preserves compactness. However, this is true for any algebraic stack (e.g., this follows from 
the description of compact generators on a stack in \secref{sss:comp gen stack}). 

\end{proof}

\ssec{The de Rham context}

In the de Rham context, the category $\Dmod(\Bun_G)$ carries an action of $\QCoh(\LocSys_\cG(X))$.

\medskip

In this subsection we will recast some of the results of the preceding subsections in terms of this action.  

\sssec{}

We will use the functors $\sP_\CZ$ and $\sP^{\on{enh}}_\CZ$ on $\Dmod(\Bun_G)$, for a map
$$f:\CZ\to \LocSys^\dr_\cG(X),$$
see Remark \ref{r:projector dr and Betti new}.

\medskip

Unwinding the definitions, as in \secref{sss:Hecke as base change old}, we obtain that there is a canonical identification
$$\on{Hecke}(\CZ,\Dmod(\Bun_G))\simeq \Dmod(\Bun_G)\underset{\QCoh(\LocSys^\dr_\cG(X))}\otimes \QCoh(\CZ),$$
so that the functor $\sP^{\on{enh}}_\CZ$ corresponds to the pullback functor
\begin{multline*}
\Dmod(\Bun_G)\simeq 
\Dmod(\Bun_G) \underset{\QCoh(\LocSys^\dr_\cG(X))}\otimes \QCoh(\LocSys^\dr_\cG(X)) \overset{\on{Id}\otimes f^*}\longrightarrow \\
\to \Dmod(\Bun_G)\underset{\QCoh(\LocSys^\dr_\cG(X))}\otimes \QCoh(\CZ).
\end{multline*}

%
%
%

\sssec{} \label{sss:P Z dr}

Let us denote the action functor of $\QCoh(\LocSys^\dr_\cG(X))$ on $\Dmod(\Bun_G)$ by 
$$\CE\in \QCoh(\LocSys^\dr_\cG(X)),\,\, \CM\in \Dmod(\Bun_G) \mapsto \CE\star \CM.$$

\medskip

In particular, we obtain that if $f$ is such that the functor 
$$f_*:\QCoh(\CZ)\to \QCoh(\LocSys^\dr_\cG(X))$$ is continuous
(in which case it is automatically a map of $\QCoh(\LocSys^\dr_\cG(X))$-module categories since 
$\QCoh(\LocSys^\dr_\cG(X))$ is rigid), then the monad
$$\oblv_{\on{Hecke}}\circ \sP^{\on{enh}}_\CZ$$ of $\Dmod(\Bun_G)$ identifies with 
$$\CF\mapsto f_*(\CO_\CZ)\star \CF.$$

%


\sssec{}

Recall that we have an identification
\begin{equation} \label{e:restr as base change again}
\Shv_\Nilp(\Bun_G)\simeq 
\Dmod(\Bun_G) \underset{\QCoh(\LocSys^\dr_\cG(X))}\otimes \QCoh(\LocSys^\dr_\cG(X))_{\LocSys^{\on{restr}}_\cG(X)} ,
\end{equation}
see \propref{p:nilp via spectral}.

\medskip

In particular, taking $\CZ=\LocSys^{\on{restr}}_\cG(X)$, we obtain that the endofunctor 
$$\oblv_{\on{Hecke}}\circ \sP$$ 
of $\Dmod(\Bun_G)$ identifies with $\iota^R\circ \iota$, 
and is given by the action of the object
$$\iota\circ \iota^R(\CO_{\LocSys^\dr_\cG(X)})\in \QCoh(\LocSys^\dr_\cG(X)),$$
where by a slight abuse of notation we denote by $\iota$ the embedding
$$\QCoh(\LocSys^\dr_\cG(X))_{\LocSys^{\on{restr}}_\cG(X)}\hookrightarrow \QCoh(\LocSys^\dr_\cG(X)),$$
and $\iota^R$ by its right adjoint (i.e., the functor of ``local sections with set-theoretic support on $\LocSys^{\on{restr}}_\cG(X)$").

\sssec{}

Finally, we claim that the assertion of \thmref{t:tensor product} can also be easily obtained from this perspective.
Indeed, since $\Dmod(\Bun_G)$ is compactly generated, and hence dualizable, the functor
$$\Dmod(\Bun_G)\otimes \Dmod(\Bun_G) \to \Dmod(\Bun_G\times \Bun_G)$$
is an equivalence.

\medskip

Now, the equivalence in \thmref{t:tensor product} can be obtained by tensoring both sides over
$$\QCoh(\LocSys^\dr_\cG(X))\otimes \QCoh(\LocSys^\dr_\cG(X))\simeq \QCoh(\LocSys^\dr_\cG(X)\times \LocSys^\dr_\cG(X))$$
with 
$$\QCoh(\LocSys^{\on{restr}}_\cG(X))\otimes \QCoh(\LocSys^{\on{restr}}_\cG(X))\simeq \QCoh(\LocSys^{\on{restr}}_\cG(X)\times \LocSys^{\on{restr}}_\cG(X)).$$

\sssec{The regular singularity property}

There is, however, one new property of the category $\Shv_{\on{Nilp}}(\Bun_G)$ that one obtains by the 
methods of the spectral projector:

\begin{maincor} \label{c:RS all}
All compact objects of $\Shv_{\on{Nilp}}(\Bun_G)$ have regular singularities.
\end{maincor}

\begin{proof}

It suffices to show that all objects of $\Shv_{\on{Nilp}}(\Bun_G)$ lie in the ind-completion of the \emph{regular}
holonomic subcategory. For that it suffices to show that the compact generators of $\Shv_{\on{Nilp}}(\Bun_G)$
have this property.

\medskip

However, this follows from the description of the generators of $\Shv_{\on{Nilp}}(\Bun_G)$ given in 
\secref{sss:generators Nilp}:

\medskip

We have to show that the objects \eqref{e:generators} have regular singularities. The objects $\delta_{y_i}$ 
have regular singularities, so it suffices to show that the endofunctors 
$$\iota \circ (\on{Id} \otimes (f_n)_*)\circ \sP_{Z_n}^{\on{enh}} \simeq
(\on{Id} \otimes \Gamma(Z_n,-)) \circ \sP_{Z_n}$$ of $\Shv(\Bun_G)$
preserve the ind-completion of the regular holonomic subcategory.

\medskip

However, this follows from the description of the functors $\sP_{\CZ}$ in \secref{sss:P Z as colim}. 

\end{proof}

Combining with \corref{c:eigensheaves nilp}, we obtain: 

\begin{maincor} \label{c:RS}
All Hecke eigensheaves have regular singularities.
\end{maincor}

The above corollary was suggested as a conjecture in \cite[Sect. 5.2.7]{BD1}.

\section{More on \conjref{c:Nilp comp gen}} \label{s:more}

The material of this section will not be used in the rest of the paper. Here we 
record several more statements that are logically equivalent to \conjref{c:Nilp comp gen}. 

\ssec{The pro-left adjoint to the embedding} \label{ss:iota left}

\sssec{}

The embedding 
$$\iota:\Shv_{\on{Nilp}}(\Bun_G)\hookrightarrow \Shv(\Bun_G)$$
does \emph{not} admit a left adjoint. However, it admits a left adjoint with values in
the pro-category
$$\iota^L:\Shv(\Bun_G)\to \on{Pro}(\Shv_{\on{Nilp}}(\Bun_G)).$$

\sssec{} \label{sss:left adjoint via pro}

The functor $\iota^L$ is related to the functors $\sP_\CZ^{\on{enh}}$ as follows: 

\medskip

Let $\CZ$ be a prestack equipped with a map $f:\CZ\to \LocSys_\cG^{\on{restr}}(X)$. Assume that
$\CO_\CZ\in \QCoh(\CZ)$ is compact. Then it follows from \corref{c:right adj to projector}(c)
that the composition
$$\Shv(\Bun_G)\overset{\iota^L}\to \on{Pro}(\Shv_{\on{Nilp}}(\Bun_G)) \overset{\on{Pro}(\on{Id} \otimes f^*)}\longrightarrow
\on{Pro}\left(\Shv_{\on{Nilp}}(\Bun_G)\underset{\QCoh(\LocSys_\cG^{\on{restr}}(X))}\otimes \QCoh(\CZ)\right)$$
takes values in
$$\Shv_{\on{Nilp}}(\Bun_G)\underset{\QCoh(\LocSys_\cG^{\on{restr}}(X))}\otimes \QCoh(\CZ)
\subset \on{Pro}\left(\Shv_{\on{Nilp}}(\Bun_G)\underset{\QCoh(\LocSys_\cG^{\on{restr}}(X))}\otimes \QCoh(\CZ)\right)$$
and identifies with the functor $\sP_\CZ^{\on{enh}}$.

\sssec{}

Vice versa, we can express the functor $\iota^L$ in terms of the functors $\sP_\CZ^{\on{enh}}$:

\medskip

Let $Z_n\overset{f_n}\to \LocSys_\cG^{\on{restr}}(X)$ be as in \secref{sss:proof of Nilp comp gen abs}. 
We have:

\begin{lem}
For a \emph{compact} $\CF\in \Shv(\Bun_G)$, the object 
$$\iota^L(\CF)\in \on{Pro}(\Shv_{\on{Nilp}}(\Bun_G))$$
identifies canonically with 
$$\underset{n}{``\on{lim}"}\, (\on{Id}\otimes (f_n)_*)\circ \sP_{Z_n}^{\on{enh}}(\CF).$$
\end{lem} 

\begin{proof}

Follows from the fact that for $\CF'\in \Shv_{\on{Nilp}}(\Bun_G)$, we have a canonical isomorphism
$$\underset{n}{\on{colim}}\, (\on{Id}\otimes (f_n)_*)\circ (\on{Id}\otimes (f_n)^!)(\CF')\simeq \CF',$$
where $(\on{Id}\otimes (f_n)_*,\on{Id}\otimes (f_n)^!)$ are the adjoint functors
$$\Shv_{\on{Nilp}}(\Bun_G)\underset{\QCoh(\LocSys_\cG^{\on{restr}}(X))}\otimes \QCoh(Z_n)\rightleftarrows
\Shv_{\on{Nilp}}(\Bun_G).$$

\end{proof}

\ssec{The right adjoint} 

\sssec{} 

Let us now consider the right adjoint $\iota^R$ of $\iota$. Note, however, that since we do not
know \conjref{c:Nilp comp gen}, the functor $\iota^R$ is a priori discontinuous.

\medskip

We claim that there exists a natural transformation
\begin{equation} \label{e:from ra to dual}
\iota^R\to \sP
\end{equation}
as functors $\Shv(\Bun_G)\to \Shv_{\on{Nilp}}(\Bun_G)$. 

\medskip

Indeed, we start with the counit of the adjunction
$$\iota \circ \iota^R\to \on{Id}_{\Shv(\Bun_G)},$$
apply to both sides the functor $\sP$, and use the fact that
$\sP\circ \iota \simeq \on{Id}_{\Shv_{\on{Nilp}}(\Bun_G)}$. 

\sssec{}

We now claim:

\begin{prop} \label{p:projector as adjoint}
The following conditions are equivalent:

\smallskip

\noindent{\em(a)} \conjref{c:Nilp comp gen} holds;

\smallskip

\noindent{\em(b)} The functor $\sP$ provides a right adjoint to $\iota$;

\smallskip

\noindent{\em(c)} The natural transformation \eqref{e:from ra to dual} is an isomorphism.

\end{prop}

\begin{proof}

We have (c) $\Rightarrow$ (b) tautologically. Also, (b) implies that $\iota^R$ is continuous, which 
implies (a). Let us show that (a) implies (c).

\medskip

Let $\iota:\bC'\to \bC$ be a fully faithful functor of DG categories, and let $\sP$ be its left inverse.
In this case, as in \eqref{e:from ra to dual}, we construct a natural transformation
$$\iota^R\to \sP.$$

\medskip

We claim that this natural transformation is an isomorphism if and only if $\sP$ annihilates 
the subcategory $(\bC')^\perp\subset \bC$. Indeed, this follows by looking at the fiber sequence
$$\iota\circ \iota^R(\bc)\to \bc \to \bc'', \quad \bc\in \bC,$$
where $\bc''\in (\bC')^\perp$.

\medskip

Note also that above, ``annihilates" is equivalent to ``preserves" since $\sP(\bC)\subset \bC'$, while
$\bC'\cap (\bC')^\perp=0$. 

\medskip

Hence, we need to show that if (a) holds, then the functor $\sP$
preserves
$$(\Shv_{\on{Nilp}}(\Bun_G))^\perp \subset \Shv(\Bun_G).$$

We will now use the fact that $\iota\circ \sP$ is an integral Hecke functor,
see \secref{sss:proj int Hecke}. We claim that assumption (a) implies that \emph{any}
integral Hecke functor preserves $(\Shv_{\on{Nilp}}(\Bun_G))^\perp$. 

\medskip

Indeed, assumption (a) implies that $(\Shv_{\on{Nilp}}(\Bun_G))^\perp$ is closed under colimits. 
Therefore, it is sufficient to show that $(\Shv_{\on{Nilp}}(\Bun_G))^\perp$ is preserved by functors of the form
\begin{equation}  \label{e:int Hecke bis} 
\Shv(\Bun_G)\overset{\on{H}(V,-)}\longrightarrow \Shv(\Bun_G\times X^I) \overset{\text{\eqref{e:int Hecke}}}\longrightarrow \Shv(\Bun_G), 
\quad V\in (\Rep(\cG)^{\otimes I})^c, \quad \CS\in \Shv(X^I)^c.
\end{equation} 

The functor \eqref{e:int Hecke bis} admits a left adjoint, which is again a functor of the same form with $V$ replaced by its monoidal
dual and $\CS$ replaced by its Verdier dual. 

\medskip

Hence, we obtain that it is enough to show that functors of the form \eqref{e:int Hecke bis} preserve
the subcategory $\Shv_{\on{Nilp}}(\Bun_G)\subset \Shv(\Bun_G)$, which was already observed in \secref{sss:Ran and Nilp}. 

\end{proof} 

\ssec{\conjref{c:Nilp comp gen} and cohomological amplitudes}


\sssec{} \label{sss:base change fro coarse}

Recall again the stacks $Z_n$ from \secref{sss:proof of Nilp comp gen abs}. Note, however, that by 
Sects. \ref{sss:almost lift action a}-\ref{sss:almost lift action c}, we can assume that these substacks are 
obtained by base change of affine schemes $S_n$ equipped with regular embeddings $g_n$ into the 
coarse moduli spaces of connected components of 
$\LocSys_\cG^{\on{restr}}(X)$. 

\medskip

For the duration of this subsection we will assume that $Z_n$ has this form.

\sssec{}

First, we claim: 

\begin{lem} \label{l:left adj right exact}
For $(Z_n,f_n)$ as above, the composite functor
$$\Shv(\Bun_G)\overset{\sP_{Z_n}^{\on{enh}}}\longrightarrow 
\Shv_{\on{Nilp}}(\Bun_G)\underset{\QCoh(\LocSys_\cG^{\on{restr}}(X))}\otimes \QCoh(Z_n)
\overset{\on{Id} \otimes (f_n)_*}\longrightarrow \Shv_{\on{Nilp}}(\Bun_G)$$
is right t-exact.
\end{lem} 

\begin{proof}

By \secref{sss:left adjoint via pro}, it suffices to show that the endofunctor of $\Shv_{\on{Nilp}}(\Bun_G)$
given by
$$(\on{Id}\otimes (f_n)_*)\circ (\on{Id} \otimes (f_n)^*)$$
is right t-exact. 

\medskip

This endofunctor is given by the action on $\Shv_{\on{Nilp}}(\Bun_G)$ of the object
$$(f_n)_*(\CO_{Z_n})\in \QCoh(\LocSys_\cG^{\on{restr}}(X)).$$

However, by the construction of $\CO_{Z_n}$, it admits a left resolution with terms being
direct sums of copies of $\CO_{\CZ}$, where $\CZ$ is a connected component of 
$\LocSys_\cG^{\on{restr}}(X)$. So, it suffices to show that the endofunctor of
$\Shv_{\on{Nilp}}(\Bun_G)$ given by the action of $\CO_\CZ$ is right t-exact.

\medskip

However, the latter functor is t-exact, being a direct summand of the identity functor. 

\end{proof}

\sssec{}

By assumption, the maps $g_n$ from the affine schemes $S_n$ to coarse moduli spaces of connected components of 
$\LocSys_\cG^{\on{restr}}(X)$ are regular closed embeddings. For each $n$, let $m_n$ denote the length of the 
corresponding regular sequence. Note that $m_n$ only depends on the choice of a connected component of
$\LocSys_\cG^{\on{restr}}(X)$. 

\medskip

We claim:

\begin{prop} \label{p:bound generators}
The collection of integers $m_n$ is bounded.
\end{prop} 

\begin{proof}

As in \secref{sss:coarse change}, there exists a curve $X'$ over $\BC$, such that $\CZ^{\on{coarse}}$
is isomorphic to its counterpart in the (restricted) Betti context. Hence, we can assume that we are
in the Betti context. In this case, by \secref{ss:coarse Betti again}, our $\CZ^{\on{coarse}}$
can be realized as the completion of $\LocSys_\cG^{\on{Betti,coarse}}(X')$ at one of its closed points.

\medskip

This implies that we can take $m=\dim(\LocSys_\cG^{\on{Betti,coarse}}(X'))$: the schemes $S_n$
can be taken to be fat points around the corresponding closed point on $\LocSys_\cG^{\on{Betti,coarse}}(X')$, 
cut out by regular sequences on $\LocSys_\cG^{\on{Betti,coarse}}(X')$. 

\end{proof}

\begin{rem}

The version of \propref{p:bound generators} holds in the general context of \thmref{t:coarse restr},
i.e., for any target gentle Tannakian category $\bH$:

\medskip

Let $\CZ$ be a connected component of $\bMaps(\Rep(\sG),\bH))$, and let $\CZ^{\on{coarse}}$
be the corresponding coarse moduli space. Let $\sigma$ be the unique closed point of $\CZ$,
and let
$$m=\dim(\sG)\cdot \dim(H^0(T_\sigma(\CZ))).$$

We claim that we can find a map $\CZ^{\on{coarse}}\to \BA^m$ mapping the unique closed point to $0\in \BA^m$ 
such that the preimage of $0\in \BA^m$ is a scheme. 
 
\medskip

Indeed, let $n$ be such that the pro-algebraic group $\sH_{\phi\on{-isotyp}}$ is topologically
generated by $n$ elements. Then, according to the proof of \thmref{t:coarse groups pro}, there exists a map
$$\CZ^{\on{coarse}}\to \sG^n/\!/\on{Ad}(\sG),$$
such that the preimage of the point in $\sG^n/\!/\on{Ad}(\sG)$ equal to the image of $\brr(\sigma)$ 
is a scheme. Hence, it suffices to show that we can choose
$n\leq \dim(H^0(T_\sigma(\CZ)))$. 

\medskip

However, according the proof of \thmref{t:top fin gen}, we can choose $n$ be exactly $\dim(H^0(T_\sigma(\CZ)))$.

\end{rem} 

\medskip

Let $m$ denote the bound from \propref{p:bound generators}. 

\sssec{}

We claim:

\begin{prop} \label{p:projector estimate}
The functor $\sP$ has a cohomological amplitude 
bounded on the right by $M$.
\end{prop}

\begin{proof}

Recall that $\sP$ is given by the action of the object $\sR\in \Rep(\cG)_\Ran$. We write $\sR$ 
as a colimit \eqref{e:formula for R univ}. We will show that the terms in the resulting colimit
expression for $\sP$ have cohomological amplitudes
bounded on the right by $m$.

\medskip

However, each such term, up to a cohomological shift by $[-m_n]$ identifies with the functor
$(\on{Id}\otimes (f_n)_*)\circ \sP_{Z_n}^{\on{enh}}$.
Hence, it is right t-exact by \lemref{l:left adj right exact}.

\end{proof}

\sssec{}

We now claim:

\begin{prop} \label{p:Nilp comp gen equiv}
The following statements are equivalent:

\smallskip

\noindent{\em(a)} \conjref{c:Nilp comp gen} holds;

\smallskip

\noindent{\em(b)} The endofunctor \eqref{e:endofunctor point b}
has a cohomological amplitude bounded on the left;

\smallskip

\noindent{\em(b')} The endofunctor \eqref{e:endofunctor point b} has a cohomological 
amplitude bounded on the left by $M$;

\smallskip

\noindent{\em(c)} The functor $\sP$ has a cohomological amplitude 
bounded on the left;

\smallskip

\noindent{\em(c')} The functor $\sP$ is left t-exact;

\smallskip

\noindent{\em(d)} The functor $\iota^R$ has a cohomological amplitude bounded on the right; 

\smallskip

\noindent{\em(d')} The functor $\iota^R$ has a cohomological amplitude bounded on the right by $M$. 

\end{prop}

\begin{proof}

We have the tautological implications (b') $\Rightarrow$ (b), (c') $\Rightarrow$ (c), (d') $\Rightarrow$ (d) 
The implication (b) $\Rightarrow$ (a) was proved in \propref{p:cohom estim 1}. The implication
(a) $\Rightarrow$ (c') follows from \propref{p:projector as adjoint}. The implication (a) $\Rightarrow$ (d') follows
by combining Propositions \ref{p:projector as adjoint} and \ref{p:projector estimate}. 

\medskip

It remains to show the implications (c') $\Rightarrow$ (b'), (c) $\Rightarrow$ (b), (d) $\Rightarrow$ (a). 

\medskip

For (c') $\Rightarrow$ (b') and (c) $\Rightarrow$ (b), we note that the functor
$$\Shv(\Bun_G) \overset{\sP^{\on{enh}}_\CZ}\to 
\Shv_\Nilp(\Bun_G)\underset{\QCoh(\LocSys_\cG^{\on{restr}}(X))}\otimes \QCoh(Z_n)\overset{\on{Id} \otimes (f_n)_*}\longrightarrow 
\Shv_\Nilp(\Bun_G)$$
is the composite of $\sP$, and the endofunctor of $\Shv_\Nilp(\Bun_G)$, given by
the action of
$$(f_n)_*(\CO_{Z_n})\in \QCoh(\LocSys_\cG^{\on{restr}}(X)).$$

Now, $(f_n)_*(\CO_{Z_n})$ admits a left resolution by copies of $\CO_{\LocSys_\cG^{\on{restr}}(X)}$ of length
$m_n$, which is $\leq M$. 

\medskip

Finally, the implication  (d) $\Rightarrow$ (a) holds for any \emph{renormalization-adapted} pair of an algebraic stack 
$\CY$ and $\CN\subset T^*(\CY)$, see \secref{sss:non qc N constraccess}. Indeed, suppose that the right adjoint 
to the embedding
$$\iota:\Shv_\CN(\CY)\hookrightarrow \Shv(\CY)$$
has a cohomological amplitude bounded on the right by some $M'$. 

\medskip

We need to show that for a collection
of objects $\CF_\alpha\in \Shv(\CY)$, the map
$$\oplus\, \iota^R(\CF_\alpha)\to \iota^R(\oplus\, \CF_\alpha)$$
is an isomorphism. Since the t-structure on $\Shv_\CN(\CY)$ is separated, it suffices
to show that for every cohomological degree $n$, the map
$$\oplus\, H^n\left(\iota^R(\CF_\alpha)\right)\to H^n\left(\iota^R(\oplus\, \CF_\alpha)\right)$$
is an isomorphism. By the assumption on the cohomological amplitude of $\iota^R$, 
we can replace $\CF_\alpha$ by $\CF'_\alpha:=\tau^{\geq n-M'}(\CF_\alpha)$, so it is 
enough to check that the map
$$\oplus\, \iota^R(\CF'_\alpha)\to \iota^R(\oplus\, \CF'_\alpha)$$
is an isomorphism. However, the latter map takes place in the category
$$(\Shv_\CN(\CY))^{\geq n-M'}\subset (\Shv_\CN(\CY))^{>-\infty} \simeq (\Shv_\CN(\CY)^{\on{access}})^{>-\infty}.$$

Hence, we can test isomorphisms by mapping out of the compact generators of $\Shv_\CN(\CY)^{\on{access}}$.
Now, the assertion follows from the fact that the functor $\iota$, restricted to $\Shv_\CN(\CY)^{\on{access}}$
preserves compactness, by the assumption on $(\CY,\CN)$. 

\end{proof} 

\section{Spectral decomposition in the Betti context}  \label{s:spectral Betti}

In this section we will work over the ground field $k=\BC$, and we will consider the
sheaf-theoretic context of \emph{all sheaves} in the classical topology, denoted 
$\Shv^{\on{all}}(-)$. 

\medskip

We will establish analogs of the results proved in the preceding sections in this
context. 

\ssec{Sheaves locally constant for the Hecke action} 

\sssec{}

Consider the Hecke action in the Betti context, which is a compatible collection of functors
\begin{equation} \label{e:full Hecke Betti}
\on{H}(-,-):\Rep(\cG)^{\otimes I} \otimes \Shv^{\on{all}}(\Bun_G) \to \Shv^{\on{all}}(\Bun_G\times X^I).
\end{equation} 

Let 
$$\Shv^{\on{all}}(\Bun_G)^{\on{Hecke-loc.const.}}\subset  \Shv^{\on{all}}(\Bun_G)$$ be the full
subcategory, consisting of objects $\CF$ for which for all $V\in \Rep(\cG)$, we have
$$\on{H}(V,F)\in  \Shv^{\on{all}}(\Bun_G)\otimes \Shv_{\on{loc.const}}^{\on{all}}(X)\subset \Shv^{\on{all}}(\Bun_G\times X).$$

It is easy to see (see, e.g.,  \cite[Proposition C.2.5]{GKRV}) that the functors \eqref{e:full Hecke Betti} send 
\begin{equation} \label{e:full Hecke Betti loc const}
\Rep(\cG)^{\otimes I} \otimes \Shv^{\on{all}}(\Bun_G)^{\on{Hecke-loc.const.}}\to
\Shv^{\on{all}}(\Bun_G)^{\on{Hecke-loc.const.}} \otimes \Shv_{\on{loc.const}}^{\on{all}}(X)^{\otimes I}.
\end{equation}  

Hence, by \thmref{t:spectral Betti}, we obtain that the category $\Shv^{\on{all}}(\Bun_G)^{\on{Hecke-loc.const.}}$ 
carries a monoidal action of $\QCoh(\LocSys^{\on{Betti}}_\cG(X))$.

\sssec{}

Recall, following \cite{NY1} (see also \cite[Theorem B.5.2]{GKRV}), that the Hecke functor 
$$\on{H}(-,-):\Rep(\cG) \otimes \Shv^{\on{all}}(\Bun_G) \to \Shv^{\on{all}}(\Bun_G\times X)$$
sends
$$\Shv^{\on{all}}_{\on{Nilp}}(\Bun_G) \subset \Shv^{\on{all}}(\Bun_G)$$
to
$$\Shv^{\on{all}}_{\on{Nilp}\times \{0\}}(\Bun_G\times X) \subset  \Shv^{\on{all}}(\Bun_G\times X).$$

By \cite[Theorem A.3.8, case (a)]{GKRV}, the external tensor product functor
$$\Shv^{\on{all}}_{\on{Nilp}}(\Bun_G)  \otimes \on{Shv}_{\on{loc.const.}}(X)\to \Shv^{\on{all}}_{\on{Nilp}\times \{0\}}(\Bun_G\times X)$$
is an equivalence. 

\medskip

In particular, we obtain that
\begin{equation} \label{e:Nilp in loc const}
\Shv^{\on{all}}_{\on{Nilp}}(\Bun_G) \subset  \Shv^{\on{all}}(\Bun_G)^{\on{Hecke-loc.const.}}
\end{equation} 

\sssec{}

The functors \eqref{e:full Hecke Betti} give rise to a system of functors
\begin{equation} \label{e:full Hecke Betti Nilp}
\on{H}(-,-):\Rep(\cG)^{\otimes I} \otimes \Shv_\Nilp^{\on{all}}(\Bun_G) \to \Shv_\Nilp^{\on{all}}(\Bun_G)\otimes \Shv_{\on{loc.const}}^{\on{all}}(X)^{\otimes I}.
\end{equation}  
and the action of $\QCoh(\LocSys^{\on{Betti}}_\cG(X))$ on $\Shv^{\on{all}}(\Bun_G)^{\on{Hecke-loc.const.}}$ 
preserves the subcategory \eqref{e:Nilp in loc const}. In particular, we reproduce the following result of \cite{NY1}
(see also \cite[Corollary 5.4.5]{GKRV}):

\begin{thm} \label{t:NY LocSys}
The functors \eqref{e:full Hecke Betti Nilp} combine to a monoidal action of $\QCoh(\LocSys^{\on{Betti}}_\cG(X))$ on 
$\Shv^{\on{all}}_\Nilp(\Bun_G)$.
\end{thm} 

\sssec{}

The following assertion, which is the main result of this section, should be regarded as an analog
of \thmref{t:lisse prel} for $\Shv^{\on{all}}(-)$:

\begin{thm} \label{t:Nilp in loc const}
The inclusion \eqref{e:Nilp in loc const} is an equality.
\end{thm} 

\ssec{The left adjoint to the embedding} \label{ss:left adjoint Betti}

In this subsection we will assume \thmref{t:Nilp in loc const} and deduce some corollaries. 

\sssec{}

Consider the category $\Shv^{\on{all}}(\Bun_G)$ and its full subcategory 
\begin{equation} \label{e:embed Nilp Betti}
\Shv_{\on{Nilp}}^{\on{all}}(\Bun_G) \overset{\iota^{\on{all}}}\hookrightarrow \Shv^{\on{all}}(\Bun_G).
\end{equation} 

According to \corref{c:left adjoint all stack}, the embedding \eqref{e:embed Nilp Betti} admits a left adjoint.
We will now describe this left adjoint in terms of the Hecke action.

\sssec{}

Let $\wt\iota^{\on{all}}$ denote the embedding 
$$\Shv^{\on{all}}(\Bun_G)^{\on{Hecke-loc.const.}}\hookrightarrow  \Shv^{\on{all}}(\Bun_G).$$ 

\medskip

First, note that, parallel to Sects. \ref{ss:Hecke Ran}-\ref{ss:impl for eigen}, given a prestack $\CZ$ over $\sfe$ and a map
$f:\CZ\to \LocSys^{\on{Betti}}_\cG(X)$, we can consider the category
$$\on{Hecke}(\CZ,\Shv^{\on{all}}(\Bun_G)),$$
equipped with a pair of adjoint functors
$$\ind_{\on{Hecke},\CZ}:\Shv^{\on{all}}(\Bun_G)\otimes \QCoh(\CZ)\rightleftarrows \on{Hecke}(\CZ,\Shv^{\on{all}}(\Bun_G)):\oblv_{\on{Hecke},\CZ}.$$

\begin{rem}
We emphasize that unlike $\Rep(\cG)_\Ran$, the category $\Rep(\cG)^{\on{Betti}}_\Ran$ is \emph{not} rigid.
So, we are using \secref{sss:proj Betti} to establish the existence of the above adjoint pair.
\end{rem}

\sssec{}

Denote by $\sP_\CZ^{\on{enh}}$ the composition
$$\ind_{\on{Hecke},\CZ} \circ (-\otimes \CO_\CZ):\Shv^{\on{all}}(\Bun_G)\to \on{Hecke}(\CZ,\Shv^{\on{all}}(\Bun_G))$$
and by $\sP_\CZ$ the composition 
$$\oblv_{\on{Hecke},\CZ} \circ \sP_\CZ^{\on{enh}}:\Shv^{\on{all}}(\Bun_G)\to \Shv^{\on{all}}(\Bun_G)\otimes \QCoh(\CZ).$$

When $\CZ$ is such that $\CO_\CZ\in \QCoh(\CZ)$ is compact, the functor $\sP_\CZ^{\on{enh}}$ is the left adjoint of
$$\oblv_{\on{Hecke}}:=(\on{Id}\otimes \Gamma(\CZ,-))\circ \oblv_{\on{Hecke},\CZ}.$$

\sssec{}

Furthermore, the essential image of $\oblv_{\on{Hecke},\CZ}$ lies in 
$$\Shv^{\on{all}}(\Bun_G)^{\on{Hecke-loc.const.}}\otimes \QCoh(\CZ) \subset \Shv^{\on{all}}(\Bun_G)\otimes \QCoh(\CZ)$$
and we have a canonical identification
$$\on{Hecke}(\CZ,\Shv^{\on{all}}(\Bun_G)) \simeq
\Shv^{\on{all}}(\Bun_G)^{\on{Hecke-loc.const.}}\underset{\QCoh(\LocSys^{\on{Betti}}_\cG(X))}\otimes \QCoh(\CZ),$$
so that the functor $\oblv_{\on{Hecke},\CZ}$ identifies with the composition
\begin{multline*} 
\Shv^{\on{all}}(\Bun_G)^{\on{Hecke-loc.const.}}\underset{\QCoh(\LocSys^{\on{Betti}}_\cG(X))}\otimes \QCoh(\CZ) 
\overset{\on{Id}\otimes (\Delta_{\LocSys^{\on{Betti}}_\cG(X)})_*}\longrightarrow \\
\to \Shv^{\on{all}}(\Bun_G)^{\on{Hecke-loc.const.}}\otimes \QCoh(\CZ) \overset{\wt\iota^{\on{all}}}\hookrightarrow 
\Shv^{\on{all}}(\Bun_G)\otimes \QCoh(\CZ).
\end{multline*}

In particular, the functor $\sP_\CZ$ maps
$$\Shv^{\on{all}}(\Bun_G)\to \Shv^{\on{all}}(\Bun_G)^{\on{Hecke-loc.const.}}\otimes \QCoh(\CZ).$$

\medskip

If $\CO_\CZ\in \QCoh(\CZ)$ is compact (from which it formally follows that the functor $f_*$
is continuous), we can rewrite the functor $\oblv_{\on{Hecke}}$ as
\begin{multline*} 
\Shv^{\on{all}}(\Bun_G)^{\on{Hecke-loc.const.}}\underset{\QCoh(\LocSys^{\on{Betti}}_\cG(X))}\otimes \QCoh(\CZ) 
\overset{\on{Id}\otimes f_*}\longrightarrow \\
\to \Shv^{\on{all}}(\Bun_G)^{\on{Hecke-loc.const.}}\underset{\QCoh(\LocSys^{\on{Betti}}_\cG(X))}\otimes \QCoh(\LocSys^{\on{Betti}}_\cG(X)) \simeq \\
\simeq \Shv^{\on{all}}(\Bun_G)^{\on{Hecke-loc.const.}}\overset{\wt\iota^{\on{all}}}\hookrightarrow  \Shv^{\on{all}}(\Bun_G).
\end{multline*}

\sssec{}

Applying this to $\CZ=\LocSys^{\on{Betti}}_\cG(X)$ and $f$ being the identity 
map\footnote{Note that, unlike the restricted situation, $\LocSys^{\on{Betti}}_\cG(X)$ is a quasi-compact algebraic stack
with an affine diagonal, and hence $\CO_{\LocSys^{\on{Betti}}_\cG(X)}\in \QCoh(\LocSys^{\on{Betti}}_\cG(X))$ is compact.}, 
we obtain that the embedding $\wt\iota^{\on{all}}$ 
admits a left adjoint, given by
$$\sP_{\LocSys^{\on{Betti}}_\cG(X)}^{\on{enh}},$$
where 
$$\wt\iota^{\on{all}}\circ \sP_{\LocSys^{\on{Betti}}_\cG(X)}^{\on{enh}}\simeq
(\on{Id}\otimes \Gamma(\LocSys^{\on{Betti}}_\cG(X),-))\circ \sP_{\LocSys^{\on{Betti}}_\cG(X)}.$$

\begin{rem} \label{r:Hecke-loc.const dualiazab}
For future use we note that the functor $\wt\iota^{\on{all}}$ realizes $\Shv^{\on{all}}(\Bun_G)^{\on{Hecke-loc.const.}}$
as a retract of $\Shv^{\on{all}}(\Bun_G)$; in particular, $\Shv^{\on{all}}(\Bun_G)^{\on{Hecke-loc.const.}}$ is dualizable. 
\end{rem}

\sssec{}

Combining with \thmref{t:Nilp in loc const}, we obtain: 

\begin{cor} \label{c:projector Betti} \hfill 

\smallskip

\noindent{\em(a)}
The functor $\sP^{\on{enh}}_{\LocSys^{\on{Betti}}_\cG(X)}$ takes values in
$$\Shv^{\on{all}}_{\Nilp}(\Bun_G)\subset \Shv^{\on{all}}(\Bun_G)^{\on{Hecke-loc.const.}}$$
and identifies canonically with $(\iota^{\on{all}})^L$.

\smallskip

\noindent{\em(b)}
The monad $\iota^{\on{all}}\circ (\iota^{\on{all}})^L$ acting on $\Shv^{\on{all}}(\Bun_G)$ identifies canonically with
$$(\on{Id}\otimes \Gamma(\LocSys^{\on{Betti}}_\cG(X),-))\circ \sP_{\LocSys^{\on{Betti}}_\cG(X)}.$$
\end{cor}

\sssec{}

Note that \corref{c:projector Betti} contains the following assertion as a particular case: 

\medskip

Let $y$ be a $\BC$-point of $\Bun_G$, and consider the corresponding object 
$$\delta_y\in \Shv(\Bun_G)\subset \Shv^{\on{all}}(\Bun_G).$$

\begin{thm} \label{t:projector Betti pt}
The object 
$$\sP^{\on{enh}}_{\LocSys^{\on{Betti}}_\cG(X)}(\delta_y)\in \Shv^{\on{all}}(\Bun_G)^{\on{Hecke-loc.const.}}$$
belongs to $\Shv^{\on{all}}_{\Nilp}(\Bun_G)$.
\end{thm}

\begin{rem} 

In \secref{ss:proof of Hecke loc const}, we will show that the assertion of \thmref{t:projector Betti pt}
actually implies that of \thmref{t:Nilp in loc const}. This is how \thmref{t:Nilp in loc const} will be proved:
we will prove \thmref{t:projector Betti pt} directly in \secref{ss:proof projector Betti pt}, thereby establishing \thmref{t:Nilp in loc const}.

\end{rem} 

\sssec{}

Finally, as in \secref{ss:gen Nilp}, we obtain:

\begin{cor} \label{c:gen Nilp Betti}
Let $y_i\in \Bun_G(\BC)$ be points chosen as in \secref{sss:indicators for N}. 
Then the objects $\sP^{\on{enh}}_{\LocSys^{\on{Betti}}_\cG(X)}(\delta_{y_i})$ generate
$\Shv^{\on{all}}_{\Nilp}(\Bun_G)$.
\end{cor} 

\ssec{Comparing $\Shv_{\on{Nilp}}(\Bun_G)$ and $\Shv^{\on{all}}_{\on{Nilp}}(\Bun_G)$}

In this subsection we continue to assume the validity of \thmref{t:Nilp in loc const}, and we will
deduce some further consequences. 

\sssec{}

First, we claim:

\begin{prop} 
The functor 
\begin{equation} \label{e:constr to all}
\Shv_{\on{Nilp}}(\Bun_G)\to \Shv^{\on{all}}_{\on{Nilp}}(\Bun_G)
\end{equation}
preserves compactness and is fully faithful.
\end{prop}

\begin{proof}

By \conjref{c:Nilp comp gen}, which holds in the Betti context (\thmref{t:Nilp comp gen Betti}), the category $\Shv_{\on{Nilp}}(\Bun_G)$ 
is generated by objects that are compact in $\Shv(\Bun_G)$. Applying \propref{p:constr Betti stack}, we obtain that the functor
\eqref{e:constr to all} preserves compactness.

\medskip

Given this, in order to prove that \eqref{e:constr to all} is fully faithful, it suffices to show that it is fully faithful
when restricted to $\Shv_{\on{Nilp}}(\Bun_G)^c$. But this follows from the fact that (for any $\CY$)
the functor 
$$\Shv(\CY)^{\on{constr}}\to  \Shv^{\on{all}}(\CY)$$
is fully faithful.

\end{proof} 

\sssec{}

We will now explain how to single out \emph{ind-constructible} sheaves with nilpotent singular support among 
\emph{all} sheaves with nilpotent singular support in terms of the Hecke action. 

\medskip

Let 
\begin{equation} \label{e:shvs all fin mon}
(\Shv^{\on{all}}_{\on{Nilp}}(\Bun_G))^{\on{Hecke-fin.mon.}} \subset \Shv^{\on{all}}_{\on{Nilp}}(\Bun_G)
\end{equation}
be the full subcategory consisting of objects $\CF$ such that for all $V\in \Rep(\cG)$ we have
$$\on{H}(V,\CF)\in 
\Shv^{\on{all}}_{\on{Nilp}}(\Bun_G)  \otimes \qLisse(X)\subset
\Shv^{\on{all}}_{\on{Nilp}}(\Bun_G)  \otimes \on{Shv}_{\on{loc.const.}}(X),$$
cf. \secref{sss:fin mon}.

\medskip

As in \cite[Proposition C.2.5]{GKRV}, one easily shows that $(\Shv^{\on{all}}_{\on{Nilp}}(\Bun_G))^{\on{Hecke-fin.mon.}} $ is stable 
under the Hecke action.

\sssec{} \label{sss:shvs all fin mon} 

Note that by \propref{p:mon fin subcategory}, the subcategory \eqref{e:shvs all fin mon} 
equals 
\begin{multline*} 
\Shv^{\on{all}}_{\on{Nilp}}(\Bun_G)\underset{\QCoh(\LocSys^{\on{Betti}}_\cG(X))}\otimes \QCoh(\LocSys_\cG^{\on{restr}}(X))
\subset  \\
\subset 
\Shv^{\on{all}}_{\on{Nilp}}(\Bun_G)\underset{\QCoh(\LocSys^{\on{Betti}}_\cG(X))}\otimes \QCoh(\LocSys^{\on{Betti}}_\cG(X))=\Shv^{\on{all}}_{\on{Nilp}}(\Bun_G),
\end{multline*}
where we view
$$\QCoh(\LocSys_\cG^{\on{restr}}(X))\simeq \QCoh(\LocSys^{\on{Betti}}_\cG(X))_{\LocSys_\cG^{\on{restr}}(X)}$$
as a co-localization of $\QCoh(\LocSys^{\on{Betti}}_\cG(X))$. 

\sssec{}

Note that the essential image of the functor
$$\Shv_{\on{Nilp}}(\Bun_G)\to \Shv^{\on{all}}_{\on{Nilp}}(\Bun_G)$$
is contained in $(\Shv^{\on{all}}_{\on{Nilp}}(\Bun_G))^{\on{Hecke-fin.mon.}}$, 
see, e.g., \eqref{e:1 step Hecke bis}. 

\medskip

We claim: 

\begin{thm} \label{t:fin monod}
The inclusion 
$$\Shv_{\on{Nilp}}(\Bun_G) \hookrightarrow (\Shv^{\on{all}}_{\on{Nilp}}(\Bun_G))^{\on{Hecke-fin.mon.}} $$
is an equality.
\end{thm}

\begin{proof} 

By \secref{sss:shvs all fin mon}, we have to show that the essential image of
\begin{multline*}
\Shv^{\on{all}}(\Bun_G)^{\on{Hecke-loc.const.}}\underset{\QCoh(\LocSys^{\on{Betti}}_\cG(X))}\otimes
\QCoh(\LocSys^{\on{Betti}}_\cG(X))_{\LocSys_\cG^{\on{restr}}(X)} 
\hookrightarrow \\
\to \Shv^{\on{all}}(\Bun_G)^{\on{Hecke-loc.const.}} \underset{\QCoh(\LocSys^{\on{Betti}}_\cG(X))}\otimes\QCoh(\LocSys^{\on{Betti}}_\cG(X))  \simeq \\
\simeq \Shv^{\on{all}}(\Bun_G)^{\on{Hecke-loc.const.}}
\end{multline*}
equals
$$\Shv_{\on{Nilp}}(\Bun_G)\subset \Shv^{\on{all}}(\Bun_G)^{\on{Hecke-loc.const.}}.$$

\medskip

We will do so by exhibiting a set of (compact) generators of 
\begin{equation} \label{e:BC Betti restr}
\Shv^{\on{all}}(\Bun_G)^{\on{Hecke-loc.const.}}\underset{\QCoh(\LocSys^{\on{Betti}}_\cG(X))}\otimes
\QCoh(\LocSys^{\on{Betti}}_\cG(X))_{\LocSys_\cG^{\on{restr}}(X)}
\end{equation} 
and show that they belong to $\Shv_{\on{Nilp}}(\Bun_G)$.

\medskip

Namely, let $y_i$ be as in \secref{sss:indicators for N}. By \corref{c:gen Nilp Betti}, the objects 
$\sP^{\on{enh}}_{\LocSys^{\on{Betti}}_\cG(X)}(\delta_{y_i})$ generate
$\Shv^{\on{all}}_{\Nilp}(\Bun_G)$.

\medskip 

Let 
$$f_n:Z_n\to \LocSys^{\on{restr}}_\cG(X)$$
be as in \secref{sss:proof of Nilp comp gen abs}. Let $\wt{f}_n$ denote the composite map
$$Z_n\overset{f_n}\to \LocSys^{\on{restr}}_\cG(X)\to \LocSys^{\on{Betti}}_\cG(X).$$

We obtain that the objects 
$$\sP^{\on{enh}}_{\LocSys^{\on{Betti}}_\cG(X)}(\delta_{y_i})\otimes (\wt{f}_n)_*(\CO_{Z_n})$$
generate \eqref{e:BC Betti restr}. 

\medskip

However, diagram chase shows that these objects are isomorphic
to the objects \eqref{e:generators}, and so they indeed belong to $\Shv_{\on{Nilp}}(\Bun_G)$.

\end{proof}

\ssec{Proof of \thmref{t:Nilp in loc const}} \label{ss:proof of Hecke loc const} 

The rest of this section is devoted to the proof of \thmref{t:Nilp in loc const}. We will deduce it
from \thmref{t:projector Betti pt}.  In its turn, \thmref{t:projector Betti pt} will be proved
independently in \secref{ss:proof projector Betti pt}. 

\sssec{}

Recall that $\iota^{\on{all}}$ denotes the embedding
$$\Shv^{\on{all}}_{\on{Nilp}}(\Bun_G)\hookrightarrow \Shv^{\on{all}}(\Bun_G),$$
and $\wt\iota^{\on{all}}$ denotes the embedding
$$\Shv^{\on{all}}(\Bun_G)^{\on{Hecke-loc.const.}}\hookrightarrow \Shv^{\on{all}}(\Bun_G).$$

\medskip

Let $'\!\iota^{\on{all}}$ denote the embedding
$$\Shv^{\on{all}}_{\on{Nilp}}(\Bun_G)\hookrightarrow \Shv^{\on{all}}(\Bun_G)^{\on{Hecke-loc.const.}},$$
so that
$$\wt\iota^{\on{all}} \circ {}'\!\iota^{\on{all}}\simeq \iota^{\on{all}}.$$

Recall that $\iota^{\on{all}}$ admits a left adjoint (by \corref{c:left adjoint all stack}), and $\wt\iota^{\on{all}}$
admits a left adjoint, namely, $\sP_{\LocSys^{\on{Betti}}_\cG(X)}^{\on{enh}}$. Restricting the functor $(\iota^{\on{all}})^L$
to $\Shv^{\on{all}}(\Bun_G)^{\on{Hecke-loc.const.}}$, we obtain a left adjoint to $'\!\iota^{\on{all}}$, to be denoted 
$({}'\!\iota^{\on{all}})^L$, so that
$$(\iota^{\on{all}})^L \simeq ({}'\!\iota^{\on{all}})^L\circ (\wt\iota^{\on{all}})^L.$$

\medskip

We wish to show that $({}'\!\iota^{\on{all}})^L$ is conservative.

\sssec{} \label{sss:Hecke loc const self-dual}

Recall (see \corref{c:self-dual N all stk}) that the category $\Shv^{\on{all}}_{\on{Nilp}}(\Bun_G)$ is naturally self-dual,
so that with respect to the canonicaly self-duality of $\Shv^{\on{all}}(\Bun_G)$ (see \secref{sss:self-dual all stk}), 
we have
$$(\iota^{\on{all}})^L\simeq (\iota^{\on{all}})^\vee.$$


\medskip

We claim now that the category $\Shv^{\on{all}}(\Bun_G)^{\on{Hecke-loc.const.}}$ is also canonically self-dual, so that 
\begin{equation} \label{e:dual iota}
(\wt\iota^{\on{all}})^L\simeq (\iota^{\on{all}})^\vee.
\end{equation} 

\medskip

Assuming this for a moment, let us prove that the functor $({}'\!\iota^{\on{all}})^L$ is conservative.

\sssec{}

It follows formally from the above properties that with respect to the above self-dualities of 
$\Shv^{\on{all}}_{\on{Nilp}}(\Bun_G)$ and $\Shv^{\on{all}}(\Bun_G)^{\on{Hecke-loc.const.}}$, respectively,
we have
$$({}'\!\iota^{\on{all}})^L \simeq ({}'\!\iota^{\on{all}})^\vee.$$

So, it suffices to show that the functor $({}'\!\iota^{\on{all}})^\vee$ is conservative. Let $\CF$ be a non-zero
object of $\Shv^{\on{all}}(\Bun_G)^{\on{Hecke-loc.const.}}$. Let $y\in \Bun_G(\BC)$ be a point such that
$$\bi_y^*(\wt\iota^{\on{all}}(\CF))\neq 0,$$
where $\bi_y$ denotes the morphism $\on{pt}\to \Bun_G$ corresponding to $y$. 

\sssec{}

Note that 
$$\bi_y^*(\wt\iota^{\on{all}}(\CF)) \simeq
\on{counit}_{\Shv^{\on{all}}(\Bun_G)}\langle\wt\iota^{\on{all}}(\CF),\delta_y\rangle,$$
which we rewrite as
$$\on{counit}_{\Shv^{\on{all}}(\Bun_G)^{\on{Hecke-loc.const.}}}\langle\CF,(\wt\iota^{\on{all}})^\vee(\delta_y)\rangle,$$
and further as
$$\on{counit}_{\Shv^{\on{all}}(\Bun_G)^{\on{Hecke-loc.const.}}}\langle\CF,(\wt\iota^{\on{all}})^L(\delta_y)\rangle\simeq
\on{counit}_{\Shv^{\on{all}}(\Bun_G)^{\on{Hecke-loc.const.}}}\langle\CF,\sP_{\LocSys^{\on{Betti}}_\cG(X)}^{\on{enh}}(\delta_y)\rangle.$$

\sssec{}

Now, by \thmref{t:projector Betti pt},
$$\sP_{\LocSys^{\on{Betti}}_\cG(X)}^{\on{enh}}(\delta_y)\in \Shv^{\on{all}}_{\on{Nilp}}(\Bun_G),$$
hence, we have
\begin{multline*}
\on{counit}_{\Shv^{\on{all}}(\Bun_G)^{\on{Hecke-loc.const.}}}\langle\CF,\sP_{\LocSys^{\on{Betti}}_\cG(X)}^{\on{enh}}(\delta_y)\rangle\simeq \\
\simeq \on{counit}_{\Shv^{\on{all}}_\Nilp(\Bun_G)}\langle({}'\!\iota^{\on{all}})^\vee(\CF),\sP_{\LocSys^{\on{Betti}}_\cG(X)}^{\on{enh}}(\delta_y)\rangle.
\end{multline*}

Hence, we obtain 
$$({}'\!\iota^{\on{all}})^\vee(\CF) \neq 0,$$
as desired.

\qed[\thmref{t:Nilp in loc const}]

\ssec{Self-duality on $\Shv^{\on{all}}(\Bun_G)^{\on{Hecke-loc.const.}}$}

In this subsection we will construct a self-duality on $\Shv^{\on{all}}(\Bun_G)^{\on{Hecke-loc.const.}}$ 
with the properties specified in \secref{sss:Hecke loc const self-dual}. 

\sssec{}

We let the counit on $\Shv^{\on{all}}(\Bun_G)^{\on{Hecke-loc.const.}}$ be induced by the 
counit on $\Shv^{\on{all}}(\Bun_G)$, i.e., 
\begin{multline*}
\Shv^{\on{all}}(\Bun_G)^{\on{Hecke-loc.const.}} \otimes \Shv^{\on{all}}(\Bun_G)^{\on{Hecke-loc.const.}}
\overset{\wt\iota^{\on{all}}\otimes \wt\iota^{\on{all}}}\longrightarrow \\
\to \Shv^{\on{all}}(\Bun_G)\otimes \Shv^{\on{all}}(\Bun_G) \overset{\on{C}^\cdot_c(\Bun_G,-\overset{*}\otimes -)}\longrightarrow \Vect_\sfe.
\end{multline*}

Recall that the unit for the self-duality on $\Shv^{\on{all}}(\Bun_G)$ is given by
$$(\Delta_{\Bun_G})_!(\ul\sfe_{\Bun_G})\in \Shv^{\on{all}}(\Bun_G\times \Bun_G)\simeq
\Shv^{\on{all}}(\Bun_G)\otimes \Shv^{\on{all}}(\Bun_G),$$
see \secref{sss:self-dual all stk}. 

\medskip

\begin{prop}  \label{p:duality Hecke lo const}
The unit maps
\begin{multline*}
((\wt\iota^{\on{all}}\circ \sP_{\LocSys^{\on{Betti}}_\cG(X)}^{\on{enh}})\otimes \on{Id})((\Delta_{\Bun_G})_!(\ul\sfe_{\Bun_G}))\to \\
\to ((\wt\iota^{\on{all}}\circ \sP_{\LocSys^{\on{Betti}}_\cG(X)}^{\on{enh}})\otimes (\wt\iota^{\on{all}}\circ \sP_{\LocSys^{\on{Betti}}_\cG(X)}^{\on{enh}}))
((\Delta_{\Bun_G})_!(\ul\sfe_{\Bun_G}))
\end{multline*}
and
\begin{multline*} 
(\on{Id}\otimes (\wt\iota^{\on{all}}\circ \sP_{\LocSys^{\on{Betti}}_\cG(X)}^{\on{enh}}))((\Delta_{\Bun_G})_!(\ul\sfe_{\Bun_G})) \to \\
\to ((\wt\iota^{\on{all}}\circ \sP_{\LocSys^{\on{Betti}}_\cG(X)}^{\on{enh}})\otimes (\wt\iota^{\on{all}}\circ \sP_{\LocSys^{\on{Betti}}_\cG(X)}^{\on{enh}}))
((\Delta_{\Bun_G})_!(\ul\sfe_{\Bun_G}))
\end{multline*}
are isomorphisms.
\end{prop}

Assuming the proposition, we obtain that the object
$$(\sP_{\LocSys^{\on{Betti}}_\cG(X)}^{\on{enh}}\otimes \on{Id})((\Delta_{\Bun_G})_!(\ul\sfe_{\Bun_G}))\in 
\Shv^{\on{all}}(\Bun_G)^{\on{Hecke-loc.const.}} \otimes \Shv^{\on{all}}(\Bun_G)$$
in fact belongs to
\begin{multline}  \label{e:loc const ten 1}
\Shv^{\on{all}}(\Bun_G)^{\on{Hecke-loc.const.}} \otimes \Shv^{\on{all}}(\Bun_G)^{\on{Hecke-loc.const.}}  \subset \\
\subset \Shv^{\on{all}}(\Bun_G)^{\on{Hecke-loc.const.}} \otimes \Shv^{\on{all}}(\Bun_G),
\end{multline}
and the object 
$$(\on{Id}\otimes \sP_{\LocSys^{\on{Betti}}_\cG(X)}^{\on{enh}})((\Delta_{\Bun_G})_!(\ul\sfe_{\Bun_G}))\in
\Shv^{\on{all}}(\Bun_G)\otimes \Shv^{\on{all}}(\Bun_G)^{\on{Hecke-loc.const.}}$$
belongs to
\begin{multline} \label{e:loc const ten 2}
\Shv^{\on{all}}(\Bun_G)^{\on{Hecke-loc.const.}} \otimes \Shv^{\on{all}}(\Bun_G)^{\on{Hecke-loc.const.}}  \subset \\
\subset \Shv^{\on{all}}(\Bun_G)\otimes \Shv^{\on{all}}(\Bun_G)^{\on{Hecke-loc.const.}},
\end{multline}
and, moreover, the above two objects of $\Shv^{\on{all}}(\Bun_G)^{\on{Hecke-loc.const.}} \otimes \Shv^{\on{all}}(\Bun_G)^{\on{Hecke-loc.const.}}$
are isomorphic.

\medskip

(Note the functors in \eqref{e:loc const ten 1} and \eqref{e:loc const ten 2} are indeed inclusions of full subcategories, since 
$\Shv^{\on{all}}(\Bun_G)^{\on{Hecke-loc.const.}}$ is dualizable, see Remark \ref{r:Hecke-loc.const dualiazab}.) 

\medskip

This implies that the above object of $\Shv^{\on{all}}(\Bun_G)^{\on{Hecke-loc.const.}} \otimes \Shv^{\on{all}}(\Bun_G)^{\on{Hecke-loc.const.}}$
defines a unit for a self-duality of $\Shv^{\on{all}}(\Bun_G)^{\on{Hecke-loc.const.}}$ so that \eqref{e:dual iota} holds. 

\sssec{Proof of \propref{p:duality Hecke lo const}}

We will prove the first isomorphism; the second one will follow by symmetry. We need to show that the unit of the adjunction 
\begin{multline*}
((\wt\iota^{\on{all}}\circ \sP_{\LocSys^{\on{Betti}}_\cG(X)}^{\on{enh}})\otimes \on{Id})((\Delta_{\Bun_G})_!(\ul\sfe_{\Bun_G}))\to \\
\to (\on{Id}\otimes (\wt\iota^{\on{all}}\circ \sP_{\LocSys^{\on{Betti}}_\cG(X)}^{\on{enh}}))\circ ((\wt\iota^{\on{all}}\circ \sP_{\LocSys^{\on{Betti}}_\cG(X)}^{\on{enh}})\otimes \on{Id}) 
((\Delta_{\Bun_G})_!(\ul\sfe_{\Bun_G}))
\end{multline*}
is an isomorphism. 

\medskip

We will show that the object 
$$((\wt\iota^{\on{all}}\circ \sP_{\LocSys^{\on{Betti}}_\cG(X)}^{\on{enh}})\otimes \on{Id})((\Delta_{\Bun_G})_!(\ul\sfe_{\Bun_G}))\in
\Shv^{\on{all}}(\Bun_G)\otimes \Shv^{\on{all}}(\Bun_G)$$
already belongs to the essential image of 
$$\Shv^{\on{all}}(\Bun_G) \otimes \Shv^{\on{all}}(\Bun_G)^{\on{Hecke-loc.const.}},$$
and hence the unit for the $(\sP_{\LocSys^{\on{Betti}}_\cG(X)}^{\on{enh}},\wt\iota^{\on{all}})$-adjunction on it
is an isomorphism.

\medskip

In fact, we will show that
\begin{equation} \label{e:sym proj}
((\wt\iota^{\on{all}}\circ \sP_{\LocSys^{\on{Betti}}_\cG(X)}^{\on{enh}})\otimes \on{Id})((\Delta_{\Bun_G})_!(\ul\sfe_{\Bun_G}))\simeq
(\on{Id}\otimes (\wt\iota^{\on{all}}\circ \sP_{\LocSys^{\on{Betti}}_\cG(X)}^{\on{enh}}))((\Delta_{\Bun_G})_!(\ul\sfe_{\Bun_G}))
\end{equation} 
as objects of $\Shv^{\on{all}}(\Bun_G) \otimes \Shv^{\on{all}}(\Bun_G)$. 

\sssec{}

Let $\tau$ denote the Cartan involution on $\cG$. The Hecke functors \eqref{e:full Hecke Betti} have the basic property
that for $V\in \Rep(\cG)^{\otimes I}$,
$$(\on{H}(V,-)\otimes \on{Id})((\Delta_{\Bun_G})_!(\ul\sfe_{\Bun_G}))\simeq
(\on{Id}\otimes \on{H}(V^\tau,-))((\Delta_{\Bun_G})_!(\ul\sfe_{\Bun_G}))$$
as objects of $\Shv(\Bun_G\times \Bun_G)$, functorially in $V$ and $I\in \on{fSet}$. 

\medskip

This implies that for $\CV\in \Rep(\cG)^{\on{Betti}}_\Ran$, we have
$$((\CV\star -)\otimes \on{Id})((\Delta_{\Bun_G})_!(\ul\sfe_{\Bun_G}))\simeq
(\on{Id}\otimes (\CV^\tau\star -))((\Delta_{\Bun_G})_!(\ul\sfe_{\Bun_G})).$$

\sssec{}

Recall that the functor $\wt\iota^{\on{all}}\circ \sP_{\LocSys^{\on{Betti}}_\cG(X)}^{\on{enh}}$ identifies with
$$(\on{Id}\otimes \Gamma(\LocSys^{\on{Betti}}_\cG(X),-))\circ \sP_{\LocSys^{\on{Betti}}_\cG(X)},$$
while $\sP_{\LocSys^{\on{Betti}}_\cG(X)}$ is the functor
$$\Shv^{\on{all}}(\Bun_G) \to \Shv^{\on{all}}(\Bun_G) \otimes \QCoh(\LocSys^{\on{Betti}}_\cG(X)),$$
given by 
$$(\sR_{\LocSys^{\on{Betti}}_\cG(X)}\star -)\circ (\on{Id}\otimes \CO_{\LocSys^{\on{Betti}}_\cG(X)})$$
for the object
$$\sR_{\LocSys^{\on{Betti}}_\cG(X)} \in \Rep(\cG)^{\on{Betti}}_\Ran\otimes \QCoh(\LocSys^{\on{Betti}}_\cG(X)),$$
see Remark \ref{sss:R dr and Betti}. 

\medskip

Hence, in order to prove \eqref{e:sym proj}, it suffices to show that
\begin{equation} \label{e:sym R}
(\tau \otimes \on{Id})(\sR_{\LocSys^{\on{Betti}}_\cG(X)})\simeq (\on{Id}\otimes \tau^*)(\sR_{\LocSys^{\on{Betti}}_\cG(X)}),
\end{equation} 
where in the right-hand side, $\tau$ denotes the involution of $\LocSys^{\on{Betti}}_\cG(X)$, induced by $\tau$.

\medskip

However, \eqref{e:sym R} follows from the canonicity of the assignment 
$$\CC\mapsto \sR_{\CC}$$
with respect to $\CC$. 

\qed[\propref{p:duality Hecke lo const}]

\section{Preservation of nilpotence of singular support} \label{s:preserve sing}

In this section we will prove \thmref{t:preserve Nilp Sing Supp prel}. Let us indicate the main idea. 

\medskip

Let us ask the general question: how can we control the singular support of $f_*(\CF)$ for a morphism $f:\CY_1\to \CY_2$
and $\CF\in \Shv(\CY_1)$ in terms of the singular support of $\CF$?  One situation in which we can do it is when $f$ is proper. 
Namely, in this case, $\on{SingSupp}(f_*(\CF))$ is contained in the pull-push of $\on{SingSupp}(\CF)$ along the diagram
$$T^*(\CY_2) \leftarrow \CY_1\underset{\CY_2}\times T^*(\CY_2)\to T^*(\CY_1).$$

However, there is one more situation when this is possible: when $f$ is the open embedding of stacks of the form
$$\BP(E)\to E/\BG_m,$$
where $E$ is a vector bundle (over some base) and $\BP(E)$ is its projectivization. In fact, this situation can be essentially
reduced to one of a proper map, see \secref{ss:contractive}. We call an open embedding of this form \emph{contractive}. 

\medskip

The idea of the proof of \thmref{t:preserve Nilp Sing Supp prel}, borrowed from \cite{DrGa2}, is to find open substacks
$\CU_i$ so that we can can calculate the singular supports of *- (or !-) extensions by reducing to the contractive situation. 

\ssec{Statement of the result}

In this subsection we will give a more precise version of \thmref{t:preserve Nilp Sing Supp prel}, in which we will
specify what the open substacks $\CU_i$ are. 

\sssec{}

Denote $\Lambda^\BQ:=\Lambda\underset{\BZ}\otimes \BQ$. We denote by $\leq$ the partial order relation on $\Lambda^\BQ$ 
$$\lambda_1\leq \lambda_2 \, \Leftrightarrow\, \lambda_2-\lambda_1\in \{\text{Positive integral span of simple coroots}\}.$$

Let 
$$\Lambda^{\BQ,+}\subset \Lambda^\BQ$$
be the cone of dominant coweights. 

\medskip

Recall (see e.g., \cite[Theorem 7.4.3]{DrGa2}) that the stack $\Bun_G$ admits a decomposition into 
locally closed substacks (known as the Harder-Narasimhan stratification)
$$\Bun_G=\underset{\theta\in \Lambda^{\BQ,+}}\cup\, \Bun_G^{(\theta)},$$
where each $\Bun_G^{(\theta)}$ is quasi-compact.

\medskip

Moreover if a subset $S\subset \Lambda^+$ satisfies 
$$\theta\in S, \,\, \theta'\leq \theta\,\Rightarrow\, \theta'\in S,$$
then 
$$\underset{\theta\in S}\cup\, \Bun_G^{(\theta)}$$ is open in $\Bun_G$.

\sssec{}

For a fixed $\theta$, let 
$$\Bun_G^{(\leq \theta)}\overset{j^\theta}\hookrightarrow \Bun_G$$
denote the embedding of the open substack corresponding to $\underset{\theta'\leq \theta}\cup\, \Bun_G^{(\theta')}$.

\medskip

The goal of this section is to prove the following theorem:

\begin{thm}  \label{t:preserve Nilp Sing Supp}
There exists an integer $c$ \emph{(}depending on $G$, $\on{char}(k)$\emph{)}\footnote{For $\on{chark}=0$ one can take $c=0$.},
such that for 
$\theta$ satisfying 
\begin{equation} \label{e:estimate}
\langle \theta,\check\alpha_i\rangle\geq (2g-2)+c, \quad \forall\, i\in I,
\end{equation} 
the functor
$$j^\theta_*:\Shv(\Bun_G^{(\leq \theta)})^{\on{constr}}\to \Shv(\Bun_G)^{\on{constr}}$$
preserves the condition of having nilpotent singular support.
\end{thm} 

\sssec{Example}

Let $X$ have genus $1$ and $\on{char}(k)=0$, so $\theta=0$ satisfies \eqref{e:estimate}. Note that
$$\Bun_G^{\leq 0}:=\Bun_G^{\on{ss}}$$
is the semi-stable locus. 

\medskip

Objects of $\Shv_{\on{Nilp}}(\Bun_G^{\on{ss}})$ are known as character sheaves. So, in this case,
\thmref{t:preserve Nilp Sing Supp} says that the functor of $*$-extension from the semi-stable locus
sends character sheaves to sheaves with nilpotent singular support. 

\begin{rem}

Since the functors 
$$j^\theta_*:\Shv(\Bun_G^{(\leq \theta)})^{\on{constr}}\to \Shv(\Bun_G)^{\on{constr}}$$
and 
$$j^\theta_1:\Shv(\Bun_G^{(\leq \theta)})^{\on{constr}}\to \Shv(\Bun_G)^{\on{constr}}$$
are related by Verdier duality, we obtain that the assertion of \thmref{t:preserve Nilp Sing Supp} 
automatically applies to the functor
$j^\theta_!$ as well (Verdier duality preserves singular support). 

\medskip

In particular, it also applies
to the functor
$$j^\theta_{!*}: \Shv(\Bun_G^{(\leq \theta)})^\heartsuit\to \Shv(\Bun_G)^\heartsuit.$$

\end{rem}

\ssec{Set-up for the proof}

%

In this subsection we will explain how the calculation of extensions from the open substacks
specified in \thmref{t:preserve Nilp Sing Supp} can be reduced to a contractive situation. 

\sssec{}

Let $P$ be a parabolic in $G$ with Levi quotient $M$ and unipotent radical $N$. 
Let is call an open substack $\CU\subset \Bun_M$ \emph{good} if for $\CP_M\in U$, we have 
$$H^1(X,V^1_{\CP_M})=0 \text{ and } H^0(X,V^2_{\CP_M})=0$$
for any irreducible $M$-representation $V^1$ that appears as a subquotient of $\fg/\fp$ and an irreducible representation $V^2$
that appears as a subquotient of $\fn$ .

\medskip

Note that the above conditions guarantee that the map 
\begin{equation} \label{e:uniformize stratum}
\CU\underset{\Bun_M}\times \Bun_P \hookrightarrow \Bun_P \overset{\sfp}\to \Bun_G
\end{equation}
is smooth and 
$$\CU\underset{\Bun_M}\times \Bun_P \hookrightarrow \Bun_P \overset{\sfq}\to \Bun_M$$
is schematic, affine and smooth (see \cite[Proposition 11.1.4]{DrGa2}). In particular, the canonical map
$$\Bun_M\to \Bun_P$$
induces a closed embedding 
\begin{equation} \label{e:BunM U}
\CU \to \CU\underset{\Bun_M}\times \Bun_P.
\end{equation}

\sssec{}

We will use the following fact established in \cite[Proposition 9.2.2 and Sect. 9.3]{DrGa2}:

\begin{thm} \label{t:DrGa} There exists an integer $c$ such that for $\theta$ satisfying \eqref{e:estimate}, the closed
substack $\Bun_G-\Bun_G^{(\leq \theta)}$ can be decomposed into a (locally finite) union of locally closed substacks $\CY$ of the following form:

\medskip

\noindent There exists a parabolic $P$ with Levi quotient $M$ and a \emph{good} open substack $\CU\subset \Bun_M$ such that
the image $\CV$ of the map \eqref{e:uniformize stratum} contains $\CY$ as a closed substack, and the (closed) substack  
$$(\CU\underset{\Bun_M}\times \Bun_P) \underset{\Bun_G}\times \CY \subset \CU\underset{\Bun_M}\times \Bun_P$$
equals the (closed) substack 
$$\CU \subset \CU\underset{\Bun_M}\times \Bun_P$$
of \eqref{e:BunM U}. 
\end{thm} 

\sssec{}   \label{sss:sing supp ext}

Using a simple inductive argument, we obtain that in order to prove \thmref{t:preserve Nilp Sing Supp}, it suffices to 
prove the following:

\medskip

Let $\CY\subset \CV$ be as in \thmref{t:DrGa}; in particular $\CY$ is closed in $\CV$. Let $j$ denote the open embedding
$$\CV-\CY \overset{j}\hookrightarrow \CV.$$
Let $\CF$ be an object of $\Shv(\CV-\CY)^{\on{constr}}$ with nilpotent singular support. Then $j_*(\CF)\in \Shv(\CV)^{\on{constr}}$ also has 
nilpotent singular support. 

\ssec{What do we need to show?}

Let us put ourselves in the situation of \secref{sss:sing supp ext}. 

\medskip

In this subsection we will formulate
a general statement that estimates from the above the singular support of objects $j_*(\CF)$ in terms
of the singular support of $\CF$, see \secref{sss:estimate}. We will show how this estimate implies 
the preservation of nilpotence of singular support. 

\medskip

The statement from \secref{sss:estimate} will be proved in \secref{ss:contractive}. 

\sssec{}  \label{sss:estimate via blow up}

Let $\wt{j}$ denote the embedding 
$$\CU\underset{\Bun_M}\times \Bun_P-\CU \hookrightarrow \CU\underset{\Bun_M}\times \Bun_P,$$
and let $\wt\CF$ denote the pullback of $\CF$ along the (smooth) projection
$$\CU\underset{\Bun_M}\times \Bun_P-\CU\to \CV-\CY.$$

\medskip

Let $\CP_M$ be a point of $\CU\subset \Bun_M\subset \Bun_P$. Note that we have a canonical identification 
$$T^*_{\CP_M}(\Bun_P) \simeq \Gamma(X,\fp^*_{\CP_M}\otimes \omega) \simeq 
\Gamma(X,\fm^*_{\CP_M}\otimes \omega) \oplus \Gamma(X,\fn^*_{\CP_M}\otimes \omega).$$

For $\wt{A}\in T^*_{\CP_M}(\Bun_P)$, let $A^0$ and $A^-$ denote its components in $\Gamma(X,\fm^*_{\CP_M}\otimes \omega)$
and $\Gamma(X,\fn^*_{\CP_M}\otimes \omega)$, respectively.

\sssec{} \label{sss:A decomp}

Let
$\CP_M$ be a point of $\CU$, and let 
$$A\in \Gamma(X,\fg^*_{\CP_M}\otimes \omega) \simeq T^*_{\CP_M}(\Bun_G)$$
be an element contained in $\on{SingSupp}(j_*(\CF))$. We wish to show that $A$ is nilpotent. 

\medskip

Consider the map
\begin{equation} \label{e:proj n+}
\Gamma(X,\fg^*_{\CP_M}\otimes \omega)\to \Gamma(X,\fp^*_{\CP_M}\otimes \omega).
\end{equation}
Let $\wt{A}\in \Gamma(X,\fp^*_{\CP_M}\otimes \omega)$ denote the image of $A$ under the map
\eqref{e:proj n+}. Note that $A$ is nilpotent if and only if the component 
$$A^0\in \Gamma(X,\fm^*_{\CP_M}\otimes \omega)$$
of $\wt{A}$ is nilpotent.

\medskip

Indeed, identify $\fg$ with $\fg^*$ using an invariant bilinear form. Write 
$$
\Gamma(X,\fg^*_{\CP_M}\otimes \omega) \simeq \Gamma(X,\fg_{\CP_M}\otimes \omega)
\simeq \Gamma(X,\fn_{\CP_M}\otimes \omega)  \oplus 
\Gamma(X,\fm_{\CP_M}\otimes \omega) \oplus \Gamma(X,\fn^-_{\CP_M}\otimes \omega).$$

The projection \eqref{e:proj n+} corresponds to the projection on the last two factors. At the same
time, the assumption on $\CU$ implies that the first factor vanishes. So, we can think of $A$
as an element of
$$\Gamma(X,\fm_{\CP_M}\otimes \omega) \oplus \Gamma(X,\fn^-_{\CP_M}\otimes \omega)
\simeq  \Gamma(X,\fp^-_{\CP_M}\otimes \omega),$$
and it is nilpotent if and only if its Levi component is such. 

\sssec{} \label{sss:estimate}

We claim that it is enough to show the following: 

\medskip

Let $\wt\CF$ be \emph{an arbitrary} object of the category $\Shv(\CU\underset{\Bun_M}\times \Bun_P-\CU)^{\on{constr}}$,
and let $\wt{A}\in T^*_{\CP_M}(\Bun_P)$ belong to $\on{SingSupp}(\wt{j}_*(\wt\CF))$. Let $A^0$ be as in \secref{sss:estimate via blow up}. 
Then there exists a point 
$$\CP'_P\in \{\CP_M\}\underset{\Bun_M}\times \Bun_P-\{\CP_M\}$$
such that the image, denoted $\wt{A}'$, of $A^0$ along
$$\Gamma(X,\fm^*_{\CP_M}\otimes \omega) \simeq \Gamma(X,\fm^*_{\CP'_M}\otimes \omega)\hookrightarrow
\Gamma(X,\fp^*_{\CP'_M}\otimes \omega)\simeq T^*_{\CP'_P}(\Bun_P)$$
belongs to $\on{SingSupp}(\wt\CF)$. 

\sssec{}

Let us show how the claim in \secref{sss:estimate} implies the needed property in \secref{sss:sing supp ext}. 

\medskip

Let $\wt\CF$ be as in \secref{sss:estimate via blow up}, and let $A$ and $\wt{A}$ be as in \secref{sss:A decomp}.
Since the projection
\begin{equation} \label{e:U to V}
\CU \underset{\Bun_M}\times \Bun_P \to \CV
\end{equation}
is smooth, the element 
$$\wt{A}\in T^*_{\CP_M}(\Bun_P)$$
belongs to the singular support of $\on{SingSupp}(\wt{j}_*(\wt\CF))$.

\medskip

Let $\CP'_P$ be as in \secref{sss:estimate}. Using again the fact that \eqref{e:U to V} is smooth, 
we obtain that there exists
$$A'\in T^*_{\CP'_P}(\Bun_G)\simeq \Gamma(X,\fg^*_{\CP'_P}\otimes \omega)$$
that belongs to $\on{SingSupp}(\CF)$, and whose image along
$$\Gamma(X,\fg^*_{\CP'_P}\otimes \omega)\to \Gamma(X,\fp^*_{\CP'_P}\otimes \omega)$$
is contained in
$$\Gamma(X,\fm^*_{\CP'_P}\otimes \omega)\subset \Gamma(X,\fp^*_{\CP'_P}\otimes \omega)$$
and equals the image of $A^0$ under the identification
$$\Gamma(X,\fm^*_{\CP'_P}\otimes \omega) \simeq \Gamma(X,\fm^*_{\CP'_M}\otimes \omega)\simeq 
\Gamma(X,\fm^*_{\CP_M}\otimes \omega).$$

By assumption, $A'$ is nilpotent, and is contained in 
$$\Gamma(X,(\fg/\fn)^*_{\CP'_P}\otimes \omega)\subset \Gamma(X,\fg^*_{\CP'_P}\otimes \omega).$$

Hence, its projection along 
$$\Gamma(X,(\fg/\fn)^*_{\CP'_P}\otimes \omega)\to \Gamma(X,\fm^*_{\CP'_P}\otimes \omega)$$
is nilpotent as well, while the latter identifies with $A^0$.

\ssec{Singular support in a contractive situation} \label{ss:contractive}

In this subsection we will provide a general context for the proof of the claim in \secref{sss:estimate}.

\sssec{} \label{sss:contractive}

Let us be given a schematic affine map of stacks $\pi:\CW\to \CU$, equipped with a section $s:\CU\to \CW$.
Assume that $\CW$, viewed as a stack over $\CU$, is equipped with an action of the monoid $\BA^1$
(with respect to multiplication), such that the action of $0\in \BA^1$ on $\CW$ equals
$$\CW \overset{\pi}\to \CU \overset{s}\to \CW.$$

\medskip 

Denote by $j$ the open embedding $\CW-\CU\hookrightarrow \CW$. Let $\CF$ be an object of $\Shv(\CW-\CU)$.
Assume that $\CF$ is equivariant with respect to $\BG_m\subset \BA^1$, which acts on $\CW-\CU$. 

\medskip

Let $u$ be a point of $\CU$ and let $\xi$ be an element of 
$$T^*_u(\CU)\oplus T^*_u(\{u\}\underset{\CU}\times \CW)\simeq T^*_u(\CW).$$

Write $\xi^0$ and $\xi^-$ for its $T^*_u(\CU)$ and $T^*_u(\{u\}\underset{\CU}\times \CW)$
components, respectively. 

\medskip

We will prove: 

\begin{prop} \label{p:estim ss}
Suppose that $\xi$ belongs to $\on{SingSupp}(j_*(\CF))$. Then there exists a point $$w\in \{u\}\underset{\CU}\times \CW-\{u\}$$
and an element $\xi'\in T^*_w(\CW)$ that belongs to $\on{SingSupp}(\CF)$ and such that $\xi'$ equals the image of $\xi^0$
under the codifferential map 
$$T^*_u(\CU)\to T^*_w(\CW).$$
\end{prop} 

By \cite[Sect. 11.2]{DrGa2}, the set-up in \secref{sss:estimate} is a particular case of the situation
in \secref{sss:contractive}. Hence, the 
claim in \secref{sss:estimate} follows from \propref{p:estim ss}.

\medskip

The rest of this subsection is devoted to the proof of \propref{p:estim ss}.

\sssec{Reduction steps}

First, by performing a smooth base change along $\CU$, we can assume that $\CU$ is an affine scheme. 

\medskip

Second, choosing homogeneous generators (for the given $\BG_m$-action) of the ring of functions on $\CW$, we can assume that $\CW$ has the form
$\CU\times \on{Tot}(E)$, where $E$ is a vector space, on which $\BG_m$ acts via a collection of characters, which
we regard as a string of positive integers denoted $(d_1,...,d_n)=\ul{d}$. 

\sssec{}

We will consider a stacky blow-up, denoted $\wt{\on{Tot}}(E)_{\ul{d}}$ of $\on{Tot}(E)$ (it will be the usual blow up $E$ at the origin
for $(d_1,...,d_n)=(1,...,1)$).

\medskip

Namely, set:
$$\wt{\on{Tot}}(E)_{\ul{d}}:=(\BA^1 \times (\on{Tot}(E)-0))/\BG_m,$$
with respect to the anti-diagonal action (where the action on $\on{Tot}(E)-0$ is given by the specified set 
of characters). 

\medskip

We have a naturally defined map 
$$p:\wt{\on{Tot}}(E)_{\ul{d}}\to \on{Tot}(E),$$
given by the action of the monoid $\BA^1$ on $\on{Tot}(E)$.

\medskip

Denote also
$$\BP(E)_{\ul{d}}:=(\on{Tot}(E)-0)/\BG_m.$$

Note that $\BP(E)_{\ul{d}}$ identifies with a closed substack of $\wt{\on{Tot}}(E)_{\ul{d}}$ corresponding to
$0\in \BA^1$. Denote the embedding of the complement by $\wt{j}$. 

\medskip

Note that map $p$ induces an \emph{isomorphism} 
\begin{equation} \label{e:punctured blowup}
(\wt{\on{Tot}}(E)_{\ul{d}}-\BP(E)_{\ul{d}})\to (\on{Tot}(E)-0).
\end{equation}

Let $q$ denote the projection
$$\wt{\on{Tot}}(E)_{\ul{d}}\to \BP(E)_{\ul{d}}.$$

Note that $q$ realizes $\wt{\on{Tot}}(E)_{\ul{d}}$ as a line bundle over $\BP(E)_{\ul{d}}$, so that the embedding
$\BP(E)_{\ul{d}}\to \wt{\on{Tot}}(E)_{\ul{d}}$ is the zero section.

\medskip

By a slight abuse of notation, we will denote by the same characters $(p,\wt{j},q)$ the corresponding
morphisms after applying $\CU\times $. 

\medskip

Let $\pi$ denote the projection $\CU\times \on{Tot}(E)\to \CU$, and let $\wt\pi$ 
denote the projection $\CU\times \wt{\on{Tot}}(E)_{\ul{d}}\to \CU$, so that
$$\wt\pi=\pi\circ p.$$
Let $\ol\pi$ denote the projection $\CU\times \BP(E)_{\ul{d}}\to \CU$, so that
$$\wt\pi=\ol\pi\circ q.$$

\sssec{}

We claim: 

\begin{lem} \label{l:stacky proj}
The maps
$$\wt{\on{Tot}}(E)_{\ul{d}}\overset{p}\to \on{Tot}(E) \text{ and }
\BP(E)_{\ul{d}}\overset{\ol\pi}\to \on{pt}$$
are proper\footnote{Note that the notion of properness is applied here to maps of algebraic stacks
that are not necessarily schematic.}.
\end{lem}

The lemma will be proved in \secref{ss:DM}. 

\sssec{}

Let $\wt\CF$ denote the pullback of $\CF$ along the isomorphism \eqref{e:punctured blowup}. 
We have
$$j_*(\CF)\simeq p_*(\wt{j}_*(\wt\CF)).$$

We record the following lemma:

\begin{lem} \label{l:gen sing supp estim}
Let $f:\CY_1\to \CY_2$ be a proper map between algebraic stacks. 
Let $\CF_1\in \Shv(\CY_1)$ and denote $\CF_2:=f_*(\CF_1)$. 
Let $y_2\in \CY_2$ be a point and $\xi_2\in T^*_{y_2}(\CY_2)$ an element contained in $\on{SingSupp}(\CF_2)$.
Then there exists $y_1\in f^{-1}(y_2)$ such that
$$df^*(\xi_2)=:\xi_1\in T^*_{y_1}(Y_1)$$
belongs to $\on{SingSupp}(\CF_1)$.
\end{lem}

\begin{proof}

It follows from the definition of singular support in \cite{Be2} that the assertion of the lemma
holds for any separated morphism $f$, for which the canonical natural transformation $f_!\to f_*$
is an isomorphism. 

\medskip

The required property for proper maps follows from \cite[Theorem 1.1]{Ols}
(see, however Remark \ref{r:avoid Olsson} for an alternative argument in our specific case). 

\end{proof}

\sssec{}

We proceed with the proof of \propref{p:estim ss}. Let $\xi=(\xi^0,\xi^-)$
be as in the statement of the proposition. 

\medskip

We claim that we can assume that $\xi^-=0$. 

\medskip

Indeed, the action of $\BA^1$ (viewed as a monoid with respect to \emph{multiplication}) on 
$\CU\times \on{Tot}(E)$ induces an action of $\BA^1$ on 
$T^*_{(u,0)}(\CU\times \on{Tot}(E))$. Since 
$\CF$ is $\BG_m$-equivariant, the subset
$$\on{SingSupp}(j_*(\CF))\cap T^*_{(u,0)}(\CU\times \on{Tot}(E))$$
is $\BG_m$-invariant. Hence, it is invariant with respect to all of $\BA^1$, and in particular, with respect to 
the action of $0\in \BA^1$. However, the action of $0$ sends the pair 
$$(\xi^0,\xi^-)\in T^*_u(\CU)\oplus T^*_0(\on{Tot}(E))$$
to $(\xi^0,0)$. 

\medskip

Hence,
$$(\xi^0,0)\in \on{SingSupp}(j_*(\CF))\cap T^*_{(u,0)}(\CU\times \on{Tot}(E)).$$

\sssec{}

Thus, let $\xi^0\in T^*_u(\CU)$ be an element so that
$$d\pi^*(\xi^0)=(\xi^0,0)\in T^*_{(u,0)}(\CU\times \on{Tot}(E))$$
belongs to $\on{SingSupp}(j_*(\CF))$. 

\medskip

By Lemmas \ref{l:stacky proj} and \ref{l:gen sing supp estim}, we can find a point 
$\ol{e}\in \BP(E)_{\ul{d}}\subset \wt{\on{Tot}}(E)_{\ul{d}}$ such that the element
\begin{equation} \label{e:covector blowup}
dp^*\circ d\pi^*(\xi^0) \in T^*_{(u,\ol{e})}(\CU\times \wt{\on{Tot}}(E)_{\ul{d}})
\end{equation}
belongs to $\on{SingSupp}(\wt{j}_*(\wt\CF))$. 

\medskip

Let $e$ be a point of $\on{Tot}(E)-0$ that projects to $\ol{e}$. Set
$$\xi':=d\pi^*(\xi^0)\in T^*_{(u,e)}(\CU\times \on{Tot}(E)).$$

We will show that $\xi'\in \on{SingSupp}(\CF)$, which will prove \propref{p:estim ss}. 

\sssec{}

%

Since $\wt\CF$ is $\BG_m$-equivariant, it is of the form
$$q^*(\ol\CF)$$
for a canonically defined $\ol\CF\in \Shv(\CU\times \BP(E)_{\ul{d}})$. 

\medskip

As was mentioned above, $\CU\times \wt{\on{Tot}}(E)_{\ul{d}}$ is a total space of a line bundle over 
$\CU\times \BP(E)_{\ul{d}}$ by means of the projection $q$. We identify
$$T^*_{(u,\ol{e})}(\CU\times \wt{\on{Tot}}(E)_{\ul{d}})\simeq T^*_{(u,\ol{e})}(\CU\times \BP(E)_{\ul{d}}) \oplus k,$$
where the first summand is the image of $dq^*$.

\medskip

It is easy to see that 
$$\on{SingSupp}(\wt{j}_*(\wt\CF)) \cap T^*_{(u,\ol{e})}(\CU\times \BP(E)_{\ul{d}}) \subset T^*_{(u,\ol{e})}(\CU\times \wt{\on{Tot}}(E)_{\ul{d}})$$
equals the image of
$$\on{SingSupp}(\ol\CF) \cap T^*_{(u,\ol{e})}(\CU\times \BP(E)_{\ul{d}})$$
along $dq^*$.

\medskip

The condition in \eqref{e:covector blowup} reads as
$$d\wt\pi^*(\xi^0)\in  \on{SingSupp}(\wt{j}_*(\wt\CF)),$$
where we also note that
$$d\wt\pi^*=dq^*\circ d\ol\pi^*.$$

Hence, 
\begin{equation} \label{e:covector proj}
d\ol\pi^*(\xi^0)\in \on{SingSupp}(\ol\CF) \cap T^*_{(u,\ol{e})}(\CU\times \BP(E)_{\ul{d}}).
\end{equation}

Using the isomorphism \eqref{e:punctured blowup}, $q$ restricts to a map
$$\CU \times (\on{Tot}(E)-0)\to \CU\times \BP(E)_{\ul{d}}.$$

\medskip

In terms of this map,
$$\on{SingSupp}(\CF) \cap T^*_{(u,e)}(\CU\times \on{Tot}(E))$$
equals the image of 
$$\on{SingSupp}(\ol\CF) \cap T^*_{(u,\ol{e})}(\CU\times \BP(E)_{\ul{d}})$$
along $dq^*$, where we also note that
$$d\pi^*=dq^*\circ d\ol\pi^*.$$

\medskip

Hence, from \eqref{e:covector proj}, we obtain 
$$d\pi^*(\xi^0)\in \on{SingSupp}(\CF),$$
as desired.

\qed[\propref{p:estim ss}]

\ssec{The stacky weighted projective space} \label{ss:DM}

In this subsection we will prove \lemref{l:stacky proj}. 

\sssec{}

First, we observe that the morphism $p:\wt{\on{Tot}}(E)_{\ul{d}}\to \on{Tot}(E)$
factors as
$$\wt{\on{Tot}}(E)_{\ul{d}}\to \BP(E)_{\ul{d}}\times \on{Tot}(E) \to \on{Tot}(E),$$
where the first arrow is a finite morphism\footnote{Note, however, that unlike the usual blowup,
this map is not necessarily a closed embedding.}. 

\medskip

Hence, it is enough to prove the assertion of the lemma that concerns $\BP(E)_{\ul{d}}$.

\sssec{}

We first consider the case when the action of $\BG_m$ on $E$ is given by the $d$-th power
of the standard character, i.e., $\ul{d}=(d,...,d)$. We will denote the resulting stack 
$\BP(E)_{\ul{d}}$ simply by $\BP(E)_d$.

\medskip

%
%
%

In this case, the map
$$\BP(E)_{d}\to \BP(E)$$
(where $\BP(E)$ is the usual projectivization of $E$) is a Zariski locally trivial fibration with fiber
$\on{pt}/\mu_d$, where $\mu_d$ is the (finite) group-scheme of $d$-th roots of unity. 

\sssec{}

We now consider the case of a general $\ul{d}$. It suffices to find another vector space $E'$ and an
integer $d'$ so that we have a $\BG_m$-equivariant finite morphism
$$E\to E',$$
where $\BG_m$ acts on $E'$ by the $d'$-th power of the standard character.

\medskip

Write $E=(\BA^1)^n$, where $\BG_m$ acts on the $i$-th copy by the $d_i$-th power of the standard 
character. Set $d':=\on{lcm}(d_i)$ and $E':=(\BA^1)^n$. 

\medskip

The morphism $E\to E'$ is given by raising to the power $\frac{d'}{d_i}$ along the $i$-th coordinate.

\qed[\lemref{l:stacky proj}]

\begin{rem} \label{r:avoid Olsson}

The proof of \propref{p:estim ss} used \lemref{l:gen sing supp estim}, in whose proof one of the ingredients
was Olsson's theorem, which implies that for proper map between algebraic stacks, the natural transformation
$f_!\to f_*$ is an isomorphism.

\medskip

In our case, the morphism in question is
$$p:\CU'\times \wt{\on{Tot}}(E)_{\ul{d}}\to \CU'\times \on{Tot}(E),$$
(for some base $\CU'$), and we claim that this corresponding property can be established directly. 

\medskip

Indeed, tracing through the above proof of \lemref{ss:DM}, and using the fact that the direct image 
along a finite map is a conservative functor, we obtain that it is sufficient to establish the corresponding 
properties for the morphism
$$\CU' \times \BP(E)_{d} \to \CU'.$$

\medskip

However, this follows from the corresponding property for the morphisms
$$\CU' \times \BP(E)_{d}\to \CU' \times \BP(E),$$
(which is easy) and
$$\CU' \times \BP(E)\to \CU'$$
(which follows from the fact that $\BP(E)$ is a proper \emph{scheme}). 

\end{rem}

\section{Proof of \thmref{t:lisse}}  \label{s:proof lisse}

In this section we will prove \thmref{t:lisse}. 

\medskip

We first consider the case of $G=GL_2$, which explains the
main idea of the argument. We then implement this idea in a slightly more involved case of $G=GL_n$
(where it is sufficient consider the minuscule Hecke functors).  

\medskip

Finally, we treat the case of an arbitrary $G$;
the proof reduces to the analysis of the local Hitchin map and affine Springer fibers. 

\ssec{Estimating singular support from below}

In this subsection we will state a general result that allows us to guarantee that a certain cotangent vector
does belong to the singular support of a sheaf obtained as a direct image.

\sssec{}

Let $\CY$ be an algebraic stack. In this section we will be operating with the notion of singular support of objects
of $\Shv(\CY)$ that do not necessarily belong to $\Shv(\CY)^{\on{constr}}$.

\medskip

By definition, for $\CF\in \Shv(\CY)$, its singular support $\on{SingSupp}(\CF)$ is the 
subset of $T^*(\CY)$ equal to the set-theoretic union of singular supports of constructible subsheaves 
of each of its perverse cohomologies. 

\medskip

We refer the reader to \secref{ss:sing supp sing} for the explanation of what we mean by $T^*(\CY)$ when 
$\CY$ is a \emph{not necessarily} smooth scheme, and to Sects. \ref{sss:sing supp sing stack}-\ref{sss:sing supp stack qc}
for the generalization for stacks. The upshot is that in practice we can always assume that 
$\CY$ is a smooth scheme.

\medskip

We emphasize that with this definition, $\on{SingSupp}(\CF)$ is \emph{not necessarily closed} as a subset of $T^*(\CY)$. 

\medskip

That said, for a closed subset $\CN\subset T^*(\CY)$, we have
$$\on{SingSupp}(\CF) \subset \CN\, \Leftrightarrow\, \CF\in \Shv_\CN(\CY).$$

In particular, if $\CY$ is smooth, 
$$\on{SingSupp}(\CF) \subset \{0\}\, \Leftrightarrow\, \CF\in \qLisse(\CY).$$

\sssec{}

Let $f:\CY_1\to \CY_2$ be a schematic separated morphism between algebraic stacks with $\CY_2$ smooth. We denote by
$df^*$ the codifferential map
$$T^*(\CY_2)\underset{\CY_2}\times \CY_1\to T^*(\CY_1).$$

\begin{thm} \label{t:sing supp dir im}
Let
$\CF_1$ be an object of $\Shv(\CY_1)$ and let $\xi_2\neq 0$ be an element of $T^*_{y_2}(\CY_2)$ for some $y_2\in \CY_2$. 
Assume there exists a point $y_1\in f^{-1}(y_2)\subset \CY_1$ such that the following conditions hold:

\smallskip

\noindent{\em(i)}
The point 
$$(\xi_2,y_1)\in  T^*(\CY_2)\underset{\CY_2}\times \CY_1$$
satisfies
$$df^*(\xi_2)\in T^*_{y_1}(\CY_1)\cap \on{SingSupp}(\CF_1),$$ 
i.e., $(\xi_2,y_1)$ belongs to the intersection 
\begin{equation} \label{e:intersect sing supp bis}
(df^*)^{-1}(\on{SingSupp}(\CF_1))\cap (\{\xi_2\}\times f^{-1}(y_2)) \subset T^*(\CY_2)\underset{\CY_2}\times \CY_1. 
\end{equation} 

\smallskip

\noindent{\em(ii)} For every cohomological degree $m$, for every constructible sub-object $\CF'_1$
of $H^m(\CF_1)$ and for every irreducible component $\CN_1$ of $\on{SingSupp}(\CF'_1)$, if
$$(\xi_2,y_1)\in (df^*)^{-1}(\CN_1),$$ then the following 
conditions are satisfied:

\smallskip

{\em(iia)} The composite map 
$$(df^*)^{-1}(\CN_1) \hookrightarrow T^*(\CY_2)\underset{\CY_2}\times \CY_1\to T^*(\CY_2)$$
is quasi-finite on a neighborhood of the point $(\xi_2,y_1)$, or equivalently, the point $(\xi_2,y_1)$
is isolated in the intersection
$$(df^*)^{-1}(\CN_1)\cap (\xi_2\times f^{-1}(y_2));$$

\smallskip

{\em(iib)} The closed substack $(df^*)^{-1}(\CN_1)$ has dimension\footnote{This inequality is automatically an equality: 
since the assertion is local, we can assume that both $\CY_1$ and $\CY_2$ are smooth schemes; then by \cite{Be2}, every $\CN_1$ 
has dimension $\dim(\CY_1)$, and hence $(df^*)^{-1}(\CN_1)$, if non-empty, has dimension $\geq \dim(\CY_2)$.} 
$\leq \dim(\CY_2)$ at the point $(\xi_2,y_1)$.

\medskip

\noindent Finally, assume:

\begin{itemize}

\item Our sheaf-theoretic context is \'etale, Betti or ind-\emph{regular} holonomic. 

\end{itemize}

\medskip

\noindent Then $\xi_2$ belongs to $\on{SingSupp}(f_*(\CF_1))$.
\end{thm}

The proof will be given in \secref{s:proof sing supp}. Several remarks are in order:

\begin{rem}

The statement of \thmref{t:sing supp dir im} appeals to the notion of dimension of a 
Zariski-closed subset in 
\begin{equation} \label{e:rel cotan}
T^*(\CY_2)\underset{\CY_2}\times \CY_1.
\end{equation}

Note that, although for a non-smooth scheme/stack $\CY$, its cotangent bundle $T^*(\CY)$
is defined only up to a unipotent gerbe (see \secref{ss:sing supp sing}), so one cannot
a priori talk unambiguously about the dimension of its closed subsets, this difficulty is not
present for \eqref{e:rel cotan}, since $\CY_2$ is smooth, and so $T^*(\CY_2)$ is a well-defined
algebraic stack.

\end{rem}

\begin{rem} 
When $\on{char}(k)=0$ and we work either with holonomic D-modules, or when $k=\BC$
and we work with constructible sheaves in the classical topology, it is known that $\on{SingSupp}(\CF_1)$ 
is Lagrangian, and hence $f((df^*)^{-1}(\on{SingSupp}(\CF_1)))$ is isotropic.

\medskip

This implies that, given (iia), condition (iib) is automatic in this case. 

\end{rem}

\begin{rem} \label{r:lisse for Dmod gen}

One can ask whether a statement analogous to \thmref{t:sing supp dir im} with condition (iib) omitted
holds when instead of $\Shv(-)$ we work with entire category of D-modules (not necessarily holonomic ones).

\medskip

We believe that the answer is yes. In fact, when the object $\CF_1\in \Dmod(\CY_1)$ is coherent, the
proof was explained to us by P.~Schapira. 

\end{rem}

\begin{rem} \label{r:lisse for Dmod} 

Since the statement of \thmref{t:sing supp dir im} excludes the de Rham context (i.e., all D-modules
or even holonomic ones), we will not be able to apply it directly to prove \thmref{t:lisse prel} in this case.

\medskip

Instead, we will deduce \thmref{t:lisse prel} in the de Rham context as follows:

\medskip

The validity of \thmref{t:sing supp dir im} for $\Shv(-)$ in the Betti context implies, by
Lefschetz principle and Riemann-Hilbert, its validity in the context of \emph{regular holonomic} D-modules. 

\medskip

We will then formally deduce the assertion of \thmref{t:lisse prel} the entire category of D-modules from the regular holonomic case, 
see \secref{ss:lisse for D-mod}. The validity of \thmref{t:lisse} for holonomic D-modules would then follow from Remark \ref{r:lisse hol}. 

\medskip

(Recall, however, that we believe that (the stronger) \thmref{t:lisse} holds for the entire category 
$\Dmod(-)$, see Remark \ref{r:lisse other dR big}.) 

\end{rem}

\begin{rem}
We do not know\footnote{This was explained to us by P.~Schapira.} 
a viable analog of \thmref{t:sing supp dir im} for the category $\Shv^{\on{all}}(-)$.
This is why our method of proof of \thmref{t:Nilp in loc const} is indirect. 
\end{rem}


\ssec{The case of $G=GL_2$}  \label{ss:proof lisse GL2}

In this section we will assume that $\on{char}(k)>2$. 

\sssec{}

Take $G=GL_2$. To shorten the notation, we will write $\Bun_2$ instead of $\Bun_{GL_2}$. 
Let $\CF$ be an object in $\Shv(\Bun_2)^{\on{Hecke-lisse}}$. We will show that the singular
support of $\CF$ is contained in the nilpotent cone.

\medskip

Let 
$$\on{H}:\Shv(\Bun_2)\to \Shv(\Bun_2\times X)$$
be the basic Hecke functor, i.e., pull-push along the diagram
$$\Bun_2 \overset{\hl}\longleftarrow \CH_2 \overset{\hr\times s}\longrightarrow \Bun_2\times X,$$
where $\CH_2$ is the moduli space of triples
\begin{equation} \label{e:basic Hecke}
\CM\overset{\alpha}\hookrightarrow\CM',
\end{equation} 
where $\CM$ and $\CM'$ are vector bundles on $X$ and $\CM'/\CM$ is a torsion sheaf of length $1$ on $X$.
The maps $\hl$ and $\hr$ send the triple $\CM\overset{\alpha}\hookrightarrow \CM'$ to $\CM$ and $\CM'$, respectively, and $s$ 
sends it to the support of $\on{coker}(\alpha)$.  

\sssec{}

%

We will argue by contradiction, so assume that $\on{SingSupp}(\CF)$ is not contained in the nilpotent cone.

\medskip

Let 
$$\xi\in T^*_\CM(\Bun_2), \quad \CM\in \Bun_2$$
be an element contained in $\on{SingSupp}(\CF)$. Recall that the
cotangent space $T^*_\CM(\Bun_2)$ identifies with the space of
$$A\in \Hom(\CM,\CM\otimes \omega).$$

We wish to show that if $A$ corresponds to $\xi$, then $A$ is nilpotent.

\sssec{}

First, we claim that $\Tr(A)=0$ as an element of $\Gamma(X,\omega)$. Indeed, consider the action
\begin{equation} \label{e:Pic act}
\on{Pic}\times \Bun_2\to \Bun_2, \quad \CL,\CM\mapsto \CL\otimes \CM.
\end{equation}

As in \secref{sss:Hecke lisse Gm}, it is easy to see that the pullback of $\CF$ along \eqref{e:Pic act} belongs to
$$\Shv_{\{0\}\times T^*(\Bun_2)}(\on{Pic}\times \Bun_2)\subset \Shv(\on{Pic}\times \Bun_2).$$

\medskip

Hence, $A$ lies in the subspace of $T^*_\CM(\Bun_2)$ perpendicular to 
$$\on{Im}(T_\one(\on{Pic})\overset{\CL\mapsto \CL\otimes \CM}\longrightarrow T_\CM(\Bun_2))\subset T_\CM(\Bun_2),$$
and this subspace exactly consists of those $A$ that have trace $0$.

\sssec{}

Assume now that $A$ is non-nilpotent. This means that 
$\det(A)\neq 0$ as an element of $\Gamma(X,\omega^{\otimes 2})$. The conditions
\begin{equation} \label{e:gen ss}
\Tr(A)=0 \,\text{ and }\, \det(A)\neq 0
\end{equation} 
(plus the assumption that $\on{char}(k)>2$) imply that at the generic point of $X$, the operator $A$ is regular semi-simple.

\medskip

Let 
$$\wt{X}\subset T^*(X)$$
be the spectral curve corresponding to $A$. The fact that $A$ is generically regular semi-simple implies
that over the generic point of $X$, the projection
$$\wt{X}\to X$$
is generically \'etale.

\medskip

Let $x\in X$ be a point which has two distinct preimages in $\wt{X}$. Let $\wt{x}$ be one of them.
We can think of $\wt{x}$ as an element $T^*_x(X)$, which we will denote by $\xi_x$. 

\medskip

We will construct a point 
$\CM'\in \Bun_2$ and $A'\in T^*_{\CM'}(\Bun_2)$, such that the element
$$(A',\xi_x)\in T^*_{\CM',x}(\Bun_2\times X)$$
belongs to $\on{SingSupp}(\on{H}(\CF))$. 

\sssec{}

For a point \eqref{e:basic Hecke} of $\CH_2$, the 
intersection of
$$(d\hl^*)(T^*_\CM(\Bun_2)) \cap (d(\hr\times s)^*)(T^*_{\CM',x}(\Bun_2\times X))\subset T^*_{(\CM\overset{\alpha}\hookrightarrow \CM')}(\CH_2)$$
consists of commutative diagrams
\begin{equation} \label{e:Hecke spectral}
\CD
\CM'  @>{A'}>> \CM'\otimes \omega  \\
@A{\alpha}AA  @AA{\alpha\otimes \on{id}}A   \\
\CM  @>{A}>> \CM\otimes \omega,
\endCD
\end{equation} 
where the corresponding element of $T^*_x(X)$ is given by the induced map
$$\CM'/\CM \to (\CM'/\CM)\otimes \omega.$$

\sssec{}  \label{sss:modifications on spectral}

We can think of $\CM$ as a torsion-free sheaf $\CL$ on $\wt{X}$, which is generically a line bundle.
The possible diagrams \eqref{e:Hecke spectral} correspond to upper modifications of 
$$\CL \hookrightarrow \CL', \quad \on{supp}_X(\CL'/\CL)\subset \{x\}\underset{X}\times \wt{X}$$
as coherent sheaves on $\wt{X}$. 

\medskip

By the assumption on $x$, there are exactly two such modifications, corresponding to the two preimages of
$x$ in $\wt{X}$. We let $(\CM',A')$ be the modification corresponding to the chosen point $\wt{x}$, so 
$A'\in T^*_{\CM'}(\Bun_2)$. 

\sssec{}

We claim that $(A',\xi_x)\in T^*_{\CM',x}(\Bun_2\times X)$ indeed belongs to $\on{SingSupp}(\on{H}(X))$.
We will do so by applying \thmref{t:sing supp dir im} to 
$$\CY_1=\CH_2,\,\, \CY_2=\Bun_2\times X,\,\, f=(\hr\times s),\,\, \CF_1=\hl^*(\CF),$$
$$y_1=(x,\CM\overset{\alpha}\hookrightarrow \CM'),\,\, y_2=(\CM',x),\,\, \xi_2=(A',\xi_x).$$

\medskip

Note that since $\hl$ is smooth, 
$$\on{SingSupp}(\hl^*(\CF))\subset T^*(\CH_2)$$
equals the image of
$$\on{SingSupp}(\CF)\underset{\Bun_2,\hl}\times \CH_2$$ 
along the codifferential of $\hl$
$$\on{SingSupp}(\CF)\underset{\Bun_2,\hl}\times \CH_2\subset
T^*(\Bun_2)\underset{\Bun_2,\hl}\times \CH_2 \to T^*(\CH_2).$$

\sssec{} \label{sss:fibers finite i}

We first verify condition (i) of \thmref{t:sing supp dir im}. The fact that the point $((A',\xi_x),(\CM\overset{\alpha}\hookrightarrow\CM'))$ belongs to 
\begin{equation} \label{e:corr Hecke 2}
\on{SingSupp}(\hl^*(\CF)) \cap \left((A',\xi_x)\times (\hr\times s)^{-1}(\CM',x)\right) \subset 
T^*(\Bun_2\times X)\underset{\Bun_2\times X,(\hr\times s)}\times \CH_2
\end{equation}
follows from the construction.

\sssec{} \label{sss:fibers finite iia}

We now verify condition (iia). We have to show that the point 
$((A',\xi_x),(\CM\overset{\alpha}\hookrightarrow\CM'))$ is isolated in the intersection \eqref{e:corr Hecke 2}. For that
end, suffices to show that the intersection 
$$\left(T^*(\Bun_2)\underset{\Bun_2,\hl}\times \CH_2\right) \cap \left((A',\xi_x)\times (\hr\times s)^{-1}(\CM',x)\right)  \subset T^*(\CH_2)$$
is finite. 

\medskip

We will establish a slightly stronger assertion, namely that the intersection 
\begin{equation} \label{e:Higgs intersect 2}
\left(T^*(\Bun_2)\underset{\Bun_2,\hl_x}\times \CH_{2,x}\right) \cap \left(A'\times (\hr_x)^{-1}(\CM')\right)  \subset T^*(\CH_{2,x})
\end{equation}
is finite, where
$$\Bun_2 \overset{\hl_x}\longleftarrow \CH_{2,x}\overset{\hr_x}\longrightarrow \Bun_2$$
is the fiber of
$$\Bun_2 \overset{\hl}\longleftarrow \CH_{2}\overset{\hr}\longrightarrow \Bun_2$$
over $x\in X$. 

\medskip

The intersection \eqref{e:Higgs intersect 2} consists of diagrams \eqref{e:Hecke spectral} with fixed 
$(\CM',A',x)$. By \secref{sss:modifications on spectral},
such diagrams are in bijection with lower modifications of $\CL'$ as a coherent sheaf on
$\wt{X}$ supported at $x$, and there are exactly two of those. 

\begin{rem}
Note that most of the above argument would apply to $\Bun_n$ for $n\geq 2$, except for the last finiteness
assertion. The latter used the fact that $A$ is generically semi-simple, which in the case $n=2$ is guaranteed 
by the conditions \eqref{e:gen ss}.
\end{rem}  

\sssec{} \label{sss:fibers finite iib}

We now verify condition (iib) of \thmref{t:sing supp dir im}. Note that \cite{Be2}, for every cohomological degree $m$
and every constructible sub-object $\CF'$ of $H^m(\CF)$, all irreducible components of 
$\on{SingSupp}(\CF')$ have dimension equal to $\dim(\Bun_G)$. 

\medskip

Hence, it suffices to show that for every $\CF'$ as above, the fibers of the composite map 
$$\on{SingSupp}(\hl^*(\CF'))\underset{T^*(\CH_2)}\times \left(T^*(\Bun_2\times X)\underset{\Bun_2\times X,(\hr\times s)}\times \CH_2\right)\to
\on{SingSupp}(\hl^*(\CF')) \overset{\hl}\to \on{SingSupp}(\CF')$$
have dimension $\leq 1$ near $((A',\xi_x),(\CM\overset{\alpha}\hookrightarrow\CM'))$.

\medskip

We will show that the fibers of the map
$$\left(T^*(\Bun_2)\underset{\Bun_2,\hl}\times \CH_2\right)\underset{T^*(\CH_2)}\times 
\left(T^*(\Bun_2\times X)\underset{\Bun_2\times X,(\hr\times s)}\times \CH_2\right)\to $$
$$\to T^*(\Bun_2)\underset{\Bun_2,\hl}\times \CH_2\to T^*(\Bun_2)$$
have dimension $\leq 1$ near $((A',\xi_x),(\CM\overset{\alpha}\hookrightarrow\CM'))$.
It suffices to show that the map 
\begin{multline} \label{e:hl fiber 2}
\left(T^*(\Bun_2)\underset{\Bun_2,\hl_x}\times \CH_{2,x}\right)\underset{T^*(\CH_{2})}\times 
\left(T^*(\Bun_2\times X)\underset{\Bun_2\times X,\hr\times s}\times \CH_{2}\right)\simeq \\
\simeq \left(T^*(\Bun_2)\underset{\Bun_2,\hl_x}\times \CH_{2,x}\right)\underset{T^*(\CH_{2,x})}\times 
\left(T^*(\Bun_2)\underset{\Bun_2,\hr_x}\times \CH_{2,x}\right)\to \\
\to T^*(\Bun_2)\underset{\Bun_2,\hl_x}\times \CH_{2,x} \to T^*(\Bun_2)
\end{multline} 
is finite near $(A',(\CM\overset{\alpha}\hookrightarrow\CM'))$. 

\medskip

Since the map \eqref{e:hl fiber 2} is proper, it suffices to show that the point $(A',(\CM\overset{\alpha}\hookrightarrow\CM'))$ is isolated in its
fiber with respect to \eqref{e:hl fiber 2}. 

\medskip

However, this is a similar finiteness assertion to what we proved in \secref{sss:fibers finite iia}.

\ssec{The case of $G=GL_n$}  \label{ss:proof lisse GLn}

In this section we will assume that $\on{char}(k)>n$. 

\medskip

We will essentially follow the same argument as in the case of $n=2$, with the difference that we will have
to use all minuscule Hecke functors, and not just the basic one. 

\sssec{}

Let $G=GL_n$, and we will write $\Bun_n$ instead of $\Bun_{GL_n}$. Let $\CF$ be an object in $\Shv(\Bun_n)^{\on{Hecke-lisse}}$. 
We will show that the singular support of $\CF$ is contained in the nilpotent cone.

\medskip

For an integer $1\leq i\leq n$, let 
$$\on{H}^i:\Shv(\Bun_n)\to \Shv(\Bun_n\times X)$$
denote the $i$-th Hecke functor, i.e., pull-push along the diagram
$$\Bun_n \overset{\hl}\longleftarrow \CH^i_n \overset{\hr\times s}\longrightarrow \Bun_n\times X,$$
where $\CH_n^i$ is the moduli space of quadruples $(x,\CM\overset{\alpha}\hookrightarrow \CM')$, where:

\begin{itemize}

\item  $x$ is a point of $X$;

\item $\CM$ and $\CM'$ are rank $n$ bundles on $X$;

\item $\alpha$ is an injection of coherent sheaves 
\begin{equation} \label{e:Hecke i}
\CM\overset{\alpha}\hookrightarrow\CM',
\end{equation} 
such that $\on{coker}(\alpha)$ has length $i$ and is \emph{scheme-theoretically} supported at $x$.

\end{itemize} 

\medskip

For future use, let
$$\Bun_n \overset{\hl_x}\longleftarrow \CH^i_{n,x} \overset{\hr_x}\longrightarrow \Bun_n$$
denote the fiber of the above picture over a given $x\in X$. 

\sssec{}

We will argue by contradiction, so assume that $\on{SingSupp}(\CF)$ is not contained in the nilpotent cone.

\medskip

Let 
$$\xi_1\in T^*_\CM(\Bun_n), \quad \CM\in \Bun_n$$
be an element contained in $\on{SingSupp}(\CF)$. Thus $\xi_1$ corresponds to an element 
$$A\in \Hom(\CM,\CM\otimes \omega),$$
and assume that $A$ is non-nilpotent. Let $x\in X$ be a point such that
$$A_x\in \Hom(\CM_x,\CM_x\otimes T^*_x(X))$$
has a non-zero eigenvalue, to be denoted $\xi_x\in T^*_x(X)$. Let $i$ denote its multiplicity (as a generalized eigenvalue). We will construct a
point $\CM'\in \Bun_n$ and $\xi_2\in T^*_{\CM'}(\Bun_n)$, such that the element 
$$(\xi_2,i\cdot \xi_x)\in T^*_{\CM',x}(\Bun_n\times X)$$
belongs to $\on{SingSupp}(\on{H}^i(X))$ (it is here that we use the assumption that $\on{char}(k)>n$, namely that the integer
$i$ is non-zero in $k$).

\sssec{}

For a point $(x,\CM\overset{\alpha}\hookrightarrow \CM')$ of $\CH^i_n$, the intersection  
$$(d\hl^*)(T^*_\CM(\Bun_n)) \cap (d(\hr\times s)^*)(T^*_{\CM',x}(\Bun_n))\subset T^*_{(x,\CM\overset{\alpha}\hookrightarrow \CM')}(\CH^i_n)$$
consists of commutative diagrams
\begin{equation} \label{e:Hecke spectral i}
\CD
\CM'  @>{A'}>> \CM'\otimes \omega  \\
@A{\alpha}AA  @AA{\alpha\otimes \on{id}}A   \\
\CM  @>{A}>> \CM\otimes \omega,
\endCD
\end{equation} 
where the corresponding element of $T^*_x(X)$ is given by the \emph{trace} of the induced map
$$\CM'/\CM \to (\CM'/\CM)\otimes \omega.$$

\sssec{}

Let $\wt{X}\subset T^*(X)$ be the spectral curve corresponding to $A$. We can think of $\CM$ as a torsion-free sheaf $\CL$ on $\wt{X}$.
Its modifications 
$$\CM \overset{\alpha}\hookrightarrow \CM'$$
that fit into \eqref{e:Hecke spectral i} are in bijection with modifications
\begin{equation} \label{e:mod tors-free}
\CL \overset{\wt\alpha}\hookrightarrow \CL'
\end{equation} 
as torsion-free coherent sheaves on $\wt{X}$.

\sssec{}

Let $\cD_x$ be the formal disc around $x$, and set
$$\wt\cD_x:=\cD_x\underset{X}\times \wt{X}.$$

Modifications as in \eqref{e:mod tors-free} are in bijection with similar 
modifications of $\CL|_{\wt\cD_x}$. 

\medskip

The multi-disc $\wt\cD_x$ can be written as
$$\wt\cD_x:=\wt\cD^1_x\sqcup \wt\cD^2_x,$$
where $\wt\cD^1_x$ is the connected component containing the element $\xi_x\in T^*_x(X)\subset T^*(X)$.  
By assumption,
\begin{equation} \label{e:cover D1}
\wt\cD^1_x\to \cD_x
\end{equation}
is a finite flat ramified cover, such that the preimage of $x\in \cD_x$ is a ``fat point" of length $i$. Hence, the rank of
\eqref{e:cover D1} is $i$. 

\medskip

In particular, we obtain that $\CL|_{\wt\cD_x^1}$, viewed as a coherent sheaf on $\cD_x$ via the pushforward along \eqref{e:cover D1},
is a vector bundle of rank equal to $i$. (Note, however, that it is not in general true that $\CL|_{\wt\cD_x^1}$ itself is a line bundle on $\wt\cD_x^1$;
that only be the case if $\xi_x$ is a \emph{regular} eigenvalue, i.e., if the dimension of the actual eigenspace with eigenvalue $\xi_x$
is $1$.) 

\medskip

We let the sought-for modification of $\CL_{\cD_x}$
be given by
$$\CL'_{\cD_x}|_{\wt\cD^1_x}=\CL'_{\cD_x}(x)|_{\wt\cD^1_x} \text{ and } 
\CL'_{\cD_x}|_{\wt\cD^2_x}=\CL'_{\cD_x}|_{\wt\cD^2_x},$$
i.e., we leave $\CL$ intact on $\wt\cD^2_x$, and twist by the divisor equal to the preimage of $x$ on $\wt\cD^1_x$. 

\sssec{}

In order to show that the pair $(\xi_2,i\cdot \xi_x)$ indeed belongs to $\on{SingSupp}(\on{H}^i(X))$, we will apply 
\thmref{t:sing supp dir im} to 
$$\CY_1=\CH^i_n, \,\,\CY_2=\Bun_n\times X,\,\, f=(\hr\times s),\,\, \CF_1=\hl^*(\CF),$$
$$y_1=(x,\CM\overset{\alpha}\hookrightarrow \CM'),\,\, y_2=(\CM',x),\,\, \xi_2=(A',i\cdot \xi_x).$$

\medskip

Let us verify conditions (i) and (ii) of \thmref{t:sing supp dir im}. We start with condition (i).

\medskip

The point 
$$((A',i\cdot \xi_x),(x,\CM\overset{\alpha}\hookrightarrow \CM'))\in T^*_{(x,\CM\overset{\alpha}\hookrightarrow \CM')}(\CH^i_n)$$
belongs to $\on{SingSupp}(\hl^*(\CF))$ by assumption. 

\sssec{} \label{sss:isolated n}

Next we verify condition (iia). As in \secref{sss:fibers finite iia}, it suffices to show that the intersection 
\begin{equation} \label{e:Higgs intersect n}
\left(T^*(\Bun_2)\underset{\Bun_n,\hl_x}\times \CH^i_{n,x}\right) \cap \left(A'\times (\hr_x)^{-1}(\CM')\right)  \subset T^*(\CH^i_{n,x})
\end{equation}
is finite. 
%
%

\medskip

We interpret the pair $(\CM',A')$ as a torsion-free sheaf $\CL'$ on $\wt{X}$, and 
the intersection \eqref{e:Higgs intersect n} consists of its lower modifications \eqref{e:mod tors-free},
such that the quotient $\CL'/\CL$ , viewed as a coherent sheaf on $X$, is scheme-theoretically
supported at $x$ and has length $i$.

\medskip

Lower modifications of $\CL'$ on $\wt{X}$ over $x\in X$ are in bijection with lower modifications
of $\CL'|_{\wt\cD_x}$. Those split into connected components indexed by the length of the 
quotient $\CL'/\CL$ on \emph{each} connected component of $\wt\cD_x$.

\medskip

Take the connected component, where the length of the modification is $i$ over $\wt\cD^1_x$,
and $0$ on all other components. We claim that this connected component consists of a single
point, which corresponds to our $(x,\CM\overset{\alpha}\hookrightarrow \CM')$. 

\medskip

Indeed, the condition on the scheme-theoretic support of $\CL'/\CL$ implies that
$$\CL'(-x)\subset \CL,$$
while the requirement on the length implies that the above inclusion is an equality. 

\sssec{}

Let us verify condition (iib) in \thmref{t:sing supp dir im}. As in \secref{sss:fibers finite iib}, 
it suffices to show that the map
\begin{multline*}
\left(T^*(\Bun_n)\underset{\Bun_n,\hl}\times \CH^i_n\right) \underset{T^*(\CH^i_n)}\times
\left(T^*(\Bun_n\times X)\underset{\Bun_n\times X,\hr\times s}\times \CH^i_n\right) \to \\
\to \left(T^*(\Bun_n)\underset{\Bun_n,\hl}\times \CH^i_n\right)\to T^*(\Bun_n)
\end{multline*}
has fibers of dimension $\leq 1$ near $((A',i\cdot \xi_x),(x,\CM\overset{\alpha}\hookrightarrow \CM'))$.

\medskip

Furthermore, it suffices to show that the map 
\begin{multline}  \label{e:hl fiber n}
\left(T^*(\Bun_n)\underset{\Bun_n,\hl_x}\times \CH^i_{n,x}\right) \underset{T^*(\CH^i_{n})}\times
\left(T^*(\Bun_n\times X)\underset{\Bun_n\times X,\hr\times s}\times \CH^i_{n}\right) \simeq \\
\simeq \left(T^*(\Bun_n)\underset{\Bun_n,\hl_x}\times \CH^i_{n,x}\right)\underset{T^*(\CH^i_{n,x})}\times 
\left(T^*(\Bun_n)\underset{\Bun_n,\hr_x}\times \CH^i_{n,x}\right)
\to \\
\to \left(T^*(\Bun_n)\underset{\Bun_n,\hl_x}\times \CH^i_{n,x}\right)\to T^*(\Bun_n)
\end{multline}
is finite near $(A',(\CM\overset{\alpha}\hookrightarrow\CM'))$. 

\medskip

Since the map \eqref{e:hl fiber n} is proper, it suffices to show that the point $(A',(\CM\overset{\alpha}\hookrightarrow\CM')$
is isolated in its fiber with respect to \eqref{e:hl fiber n}. The latter is proved by the same consideration 
as in \secref{sss:isolated n}. 

\ssec{A digression: the notion of $(G,M)$-regularity}

Before we tackle \thmref{t:lisse} for a general reductive group $G$, we will need to make a digression
on the structure of Lie algebras.

\sssec{}

Fix a Cartan subgroup $T\subset G$, and a Levi subgroup $T\subset M\subset G$. We consider the affine schemes
$$\fa:=\fg/\!/\on{Ad}(G)\simeq \ft/\!/W \text{ and } \fa_M:=\fm/\!/\on{Ad}(M)\simeq \ft/\!/W_M,$$
and a natural map between them\footnote{Note that the above isomorphisms are part of our 
assumptions on $\on{char}(k)$.}.  

\medskip

Let $\oft_M\subset \ft$ be the open subset consisting of elements $\sft\in \ft$ for which $\check\alpha(\sft)\neq 0$
for all roots $\check\alpha$ that are \emph{not} roots of $M$. Since this subset is $W_M$-invariant,
it corresponds to an open subset
$$\ofa_M\subset \fa_M,$$
so that we have a Cartesian diagram
\begin{equation} \label{e:reg diag}
\CD
\oft_M @>>> \ft \\
@VVV  @VVV  \\
\ofa_M  @>>> \fa_M.
\endCD
\end{equation} 

We will refer to $\ofa_M$ as the $(G,M)$-\emph{regular} locus of $\fa_M$.

\sssec{}

We observe:

\begin{lem} \label{l:regularity}
For an element $A\in \fm$ the following conditions on are equivalent:  

\smallskip

\noindent{\em(i)} $Z_\fg(A)\subset \fm$;

\smallskip

\noindent{\em(i')} The adjoint action of $A$ on $\fg/\fm$ is invertible;

\smallskip

\noindent{\em(ii)} $Z_\fg(A^{\on{ss}})\subset \fm$, where $A^{\on{ss}}$ is the semi-simple part of $A$;

\smallskip

\noindent{\em(ii')} The adjoint action of $A^{\on{ss}}$ on $\fg/\fm$ is invertible;

\smallskip

\noindent{\em(iii)} The image of $A$ in $\fa_M$ belongs to $\ofa_M$.

\end{lem}

\begin{proof}

Clearly (i) $\Leftrightarrow$ (i') and (ii) $\Leftrightarrow$ (ii'). However, it is also clear that 
(i') $\Leftrightarrow$ (ii').  The equivalence (iii) $\Leftrightarrow$ (ii') is the fact that the diagram
\eqref{e:reg diag} is Cartesian. 

\end{proof} 

\sssec{} \label{sss:G M reg}

Let us say that an element $A\in \fm$ is $(G,M)$-\emph{regular} if it satisfies the equivalent 
conditions of \lemref{l:regularity}. 

\medskip

Elements of $\fm$ that are $(G,M)$-regular form a Zariski-open subset to be denoted
$\ofm$. We have a Cartesian diagram
\begin{equation} \label{e:am Cart}
\CD
\ofm @>>> \fm \\
@VVV  @VVV  \\
\ofa_M  @>>> \fa_M.
\endCD
\end{equation}

\sssec{}

We now claim:

\begin{lem}  \label{l:regularity bis} \hfill

\smallskip

\noindent{\em(a)} The open subset $\ofa_M\subset \fa_M$ is the locus of etaleness of the map 
$$\fa_M\to \fa.$$

\smallskip

\noindent{\em(b)} The open subset $\ofm\subset \fm$ is the locus of etaleness of the map 
$$\fm/\on{Ad}(M)\to \fg/\on{Ad}(G).$$

\smallskip

\noindent{\em(c)}
The diagram
$$
\CD
\ofm/\on{Ad}(M)  @>>> \fg/\on{Ad}(G) \\
@VVV  @VVV  \\
\ofa_M  @>>>  \fa
\endCD
$$
is Cartesian. 
\end{lem} 

\begin{proof}

Point (a) follows from the third assumption on $\on{char}(k)$ in \secref{sss:assump char}: 
an element $\sft\in \ft$ belongs to $\ofm$ if and only if 
its stabilizer in $W$ is contained in $W_M$. 

\medskip

Point (b) is a straightforward tangent space calculation. 

\medskip

For point (c), we note that by points (a) and (b), the map
$$\ofm/\on{Ad}(M) \to \fa_M \underset{\fa}\times  \fg/\on{Ad}(G)$$
is \'etale. So, it is sufficient to check that it is bijective at the level of
field-valued points, which follows from Jordan decomposition and \lemref{l:regularity}.  

\end{proof}

\ssec{The case of an arbitrary reductive group $G$} \label{ss:Hecke arb G}

The proof in the case of an arbitrary $G$ will follow the same idea as 
in the case of $GL_n$. What will be different is the local analysis:

\medskip

In the case of $GL_n$, to a cotangent vector to $\Bun_G$ (a.k.a. Higgs field), we attached
its spectral curve $\wt{X}$, and proved the theorem by analyzing the behavior of modifications
of sheaves on it.

\medskip

For an arbitrary $G$, there is no spectral curve.  Instead, our local analysis will amount
to studying the fibers of the affine (parabolic) Springer map. 

\sssec{}

Recall that the first assumption on $\on{char}(k)$ in \secref{sss:assump char}
says that there exists a non-degenerate $G$-equivariant pairing
\begin{equation} \label{e:almost Killing}
\fg\otimes \fg\to k,
\end{equation} 
whose restriction to the center of any Levi subalgebra remains non-degenerate. 

\medskip

We will use the pairing \eqref{e:almost Killing} to identify
$\fg^*$ with $\fg$ as $G$-modules, and also $\fm^*$ with $\fm$ for any Levi subgroup $M\subset G$. 

\sssec{}

Let $\CF\in \Shv(\Bun_G)$ be an object with non-nilpotent singular support. We will find an irreducible
representation $V^\lambda\in \Rep(\cG)$, such that the corresponding Hecke functor 
$$\on{H}(V^\lambda,-):\Shv(\Bun_G)\to \Shv(\Bun_G\times X),$$
sends $\CF$ to an object of $\Shv(\Bun_G\times X)$ whose singular support is \emph{not} contained
in 
$$T^*(\Bun_G)\times \{\text{zero-section}\}\subset T^*(\Bun_G)\times T^*(X)=T^*(\Bun_G\times X).$$

\sssec{}  \label{sss:specify M}

Using the pairing \eqref{e:almost Killing}, we can think of points of $T^*(\Bun_G)$ as pairs $(\CP_G,A)$,
where $\CP_G$ is a $G$-bundle on $X$ and $A$ is an element of
$$\Gamma(X,\fg_{\CP_G}\otimes \omega).$$

The Chevaley map attaches to $A$ above a global section $\on{ch}(A)$ of $\fa_\omega$, where the latter is the $\omega$-twist of 
$$\fa:=\fg/\!/G.$$ 

\medskip

By assumption, $\on{SingSupp}(\CF)$ contains a point $(\CP_G,A)$ for which $A$ is non-nilpotent, i.e.,
$\on{ch}(A)\neq 0$. Let $x\in X$ be a point for which the value 
$$\on{ch}(A)_x\in \fa_{\omega_x}\simeq \left(\fa \times (T^*_x(X)-0)\right)/\BG_m$$
of $\on{ch}(A)$ at $x$ is non-zero. 

\medskip

Choose a preimage $t_x\in \ft \otimes T^*_x(X)$ of $\on{ch}(A)_x$ along the projection 
$$\ft \otimes T^*_x(X)\to \fa_{\omega_x}.$$ 
Let $M$ be the Levi subgroup of $G$ equal to the centralizer of $t_x$. (Thus, if $\on{ch}(A)_x$ were zero, we would get
$M=G$, and if $A_x$ was regular semi-simple, we would get $M=T$, the Cartan subgroup.) 


\medskip

Let $\lambda$ be a coweight of $Z(M)$ that is $(G,M)$-regular (the latter means that the 
centralizer of $\lambda$ in $G$ is contained in $M$, see \secref{sss:G M reg}). By the non-degeneracy assumption on 
\eqref{e:almost Killing}, we can choose $\lambda$ so that the value of the pairing \eqref{e:almost Killing} on the 
pair $(t_x,\lambda)$ is non-zero. 

\medskip

Let $V^\lambda$ be an irreducible representation of $\cG$ corresponding to the conjugacy class of $\lambda$. 

\medskip

We claim that with this choice
of $\lambda$, the singular support of the object
$$\on{H}(V^\lambda,\CF)\in \Shv(\Bun_G\times X)$$
at the point $(\CP'_G,x)\in \Bun_G\times X$ will contain an element $(A',\xi_x)$, where $\CP'_G$ is 
the Hecke modification of $\CP_G$ at $x$ of type $\lambda$ specified in \secref{sss:specified Hecke} below, and
$$0\neq \xi_x\in T^*_x(X).$$

The element $A'$ will also be specified in \secref{sss:specified Hecke} below. 

\sssec{} \label{sss:specified Hecke}

By the choice of $M$, the fiber $(\CP_{G,x},A_x)$ of $(\CP_G,A)$ at $x$ admits a reduction $(\CP_{M,x},A_x)$ to $M$, so that
$$A_x\in \fm_{\CP_{M,x}}\otimes T^*_x(X)$$
is such that its semi-simple part lies in 
$$Z(\fm_{\CP_{M,x}}) \otimes T^*_x(X) \subset \fm_{\CP_{M,x}}\otimes T^*_x(X)$$
and is $(G,M)$-regular (see \secref{sss:G M reg}). 

\medskip

Note now that the map of the stack-theoretic quotients 
$$\fm/\on{Ad}(M)\to \fg/\on{Ad}(G)$$
is \'etale on the $(G,M)$-regular locus (see \lemref{l:regularity bis}).
This implies that the restriction $\CP^{\on{loc}}_G:=\CP_G|_{\cD_x}$ admits a \emph{unique} reduction to $M$,
to be denoted $\CP^{\on{loc}}_M$, such that:

\begin{itemize}

\item The value of $\CP^{\on{loc}}_M$ at $x$ is $\CP_{M,x}$;

\item $A^{\on{loc}}:=A|_{\cD_x}$ lies in $\Gamma(\cD_x,\fm_{\CP^{\on{loc}}_M}\otimes \omega)$;

\item $A^{\on{loc}}_x=A_x$ as elements of $\fm_{\CP_{M,x}}\otimes T^*_x(X)$. 

\end{itemize}

\medskip

Being a cocharacter of $Z(M)$, the element $\lambda$ defines a distinguished modification $\CP'{}_M^{\on{loc}}$ of $\CP^{\on{loc}}_M$.
We let $\CP'{}^{\on{loc}}_G$ be the induced modification of $\CP^{\on{loc}}_G:=\CP_G|_{\cD_x}$, and we let $\CP'_G$ denote the resulting
modification of $\CP_G$, i.e., 
$$
\begin{cases}
&\CP'_G|_{\cD_x}=\CP'{}^{\on{loc}}_G,\\
&\CP'_G|_{X-x}=\CP_G|_{X-x}.
\end{cases}
$$

\medskip

The centrality of $\lambda$ implies that we have a natural identification
$$\fm_{\CP'{}^{\on{loc}}_M}\simeq \fm_{\CP^{\on{loc}}_M},$$
and hence $A^{\on{loc}}$
gives rise to a section
$$A'{}^{\on{loc}}\in \Gamma(\cD_x,\fm_{\CP'{}^{\on{loc}}_M}\otimes \omega).$$

By a slight abuse of notation we will denote by the same symbol $A'{}^{\on{loc}}$ its image along 
$$\Gamma(\cD_x,\fm_{\CP'{}^{\on{loc}}_M}\otimes \omega)\to  
\Gamma(\cD_x,\fg_{\CP'{}^{\on{loc}}_M}\otimes \omega)=
\Gamma(\cD_x,\fg_{\CP'{}^{\on{loc}}_G}\otimes \omega)=
\Gamma(\cD_x,\fg_{\CP'_G}\otimes \omega).$$

Let 
$$A'\in \Gamma(X,\fg_{\CP'_G}\otimes \omega)$$
denote the element such that
$$
\begin{cases}
&A'|_{\cD_x}=A'{}^{\on{loc}},\\
&A'|_{X-x}=A|_{X-x}.
\end{cases}
$$

\sssec{} \label{sss:points a-d}

Consider the Hecke stack 
$$\Bun_G \overset{\hl}\longleftarrow \CH_G \overset{\hr\times s}\longrightarrow \Bun_G\times X.$$

For future use, denote by 
$$\Bun_G \overset{\hl_x}\longleftarrow \CH_{G,x} \overset{\hr_x}\longrightarrow \Bun_G$$
the fiber of this picture over a given $x\in X$.

\medskip

We will apply \thmref{t:sing supp dir im} to 
$$\CY_1=\CH_G, \,\, \CY_2=\Bun_G\times X,\,\, f=\hr\times s,\,\, \CF_1=\hl^*(\CF)\otimes \tau^*(\CV^\lambda),$$
where:

\begin{itemize}

\item $\tau: \CH_G\to \CH_G^{\on{loc}}$ is the projection on the local Hecke stack (see \cite[Sect.B.3.2]{GKRV});

\item $\CV^\lambda\in \Shv(\CH^{\on{loc}}_G)$ corresponds to $V^\lambda\in \Rep(\cG)$ by geometric Satake. 

\end{itemize}

\medskip

We take $y_2=(\CP'_G,x)$ and $y_1$ corresponding to the modification $\CP_G\overset{\alpha}\rightsquigarrow \CP'_G$. 

\begin{rem}

In what follows, we will appeal to the cotangent bundle of $\CH_G$ and related geometric objects, and
to the notion of singular support of sheaves on them. The apotropaic definitions that justify these manipulations
are spelled out in \secref{ss:sing supp sing}.

\end{rem}

\sssec{}

We will show the following: 

\medskip

\noindent(a) There exists \emph{some} $\xi_x\in T^*_x(X)$ such that 
$$((A',\xi_x),(x,\CP_G\rightsquigarrow \CP'_G))\in \left(T^*(\Bun_G\times X)\underset{\Bun_G\times X,\hr\times s}\times \ol\CH^\lambda_{G}\right) \cap 
\on{SingSupp}\Biggl(\hl^*(\CF)\otimes \tau^*(\CV^\lambda)\Biggr);$$

\medskip

\noindent(b) $\xi_x$ is the value of the pairing \eqref{e:almost Killing} on the pair $(A_x,\lambda)$, or equivalently, 
$(t_x,\lambda)$, and hence, is non-zero by the choice of $\lambda$; 

\medskip

\noindent(c) The point $((A',\xi_x),(x,\CP_G\rightsquigarrow \CP'_G))$ is isolated in the intersection
$$\left(T^*(\Bun_G\times \CH_G^{\on{loc}})\underset{\Bun_G\times \CH_G^{\on{loc}},\hl\times \tau}\times \ol\CH^\lambda_G\right)\cap 
\left((A',\xi_x) \times (\hr\times s)^{-1}(\CP'_G,x)\right) \subset T^*(\CH_G),$$
where $\ol\CH^\lambda_G$ is the closure of $\CH^\lambda_G\subset \CH_G$, the latter being the locus of modifications of type $\lambda$. 

\medskip

\noindent(d) The point $(A',(\CP_G\rightsquigarrow \CP'_G))$ is isolated in its fiber along the map
\begin{multline*}
\left(T^*(\Bun_G\times \CH_{G,x}^{\on{loc}})\underset{\Bun_G\times \CH_{G,x}^{\on{loc}},\hl_x\times \tau_x}\times \ol\CH^\lambda_{G,x}\right)
\underset{T^*(\ol\CH^\lambda_{G,x})}\times \left(T^*(\Bun_G)\underset{\Bun_G,\hr_x}\times \ol\CH^\lambda_{G,x}\right) \to \\
\to T^*(\Bun_G\times \CH_{G,x}^{\on{loc}})\underset{\Bun_G\times \CH_{G,x}^{\on{loc}},\hl_x\times \tau_x}\times \ol\CH^\lambda_{G,x} \to
T^*(\Bun_G),
\end{multline*}
where 
$$\tau_x:\CH_{G,x}\to \CH^{\on{loc}}_{G,x}, \quad \CV^\lambda_x\in \Shv(\CH^{\on{loc}}_{G,x})$$
are the counterparts of $(\tau,\CV^\lambda)$ at $x$.

\medskip

Arguing as in Sects. \ref{ss:proof lisse GL2} and \ref{ss:proof lisse GLn}, and using the fact that the map $\tau: \CH_G\to \CH_G^{\on{loc}}$
is pro-smooth, once we establish properties (a)-(d), the assertion of \thmref{t:lisse} will follow by applying \thmref{t:sing supp dir im}. 

\ssec{Proof of points (a) and (b)}

\sssec{} \label{sss:Sat loc const}

Since the map
$$\hl\times s: \CH_G\to \Bun_G\times X,$$
locally in the smooth topology, it can be isomorphed to the product situation 
$$\Bun_G\times X\times \Gr_G \to \Bun_G\times X,$$
so that $\hl^*(\CF)\otimes \tau^*(\CV)\in \Shv(\CH_G)$ identifies with
$$\CF \boxtimes \ul\sfe_X\boxtimes \CV'\in \Shv(\Bun_G\times X\times \Gr_G), \quad \CV'\in \Shv(\Gr_G),$$
in order to prove point (a) it is sufficient (in fact, equivalent) to show:

\medskip

\noindent(a') 
$$(A',(\CP_G\overset{\alpha}\rightsquigarrow \CP'_G))\in 
\left(T^*(\Bun_G)\underset{\Bun_G,\hr_x}\times \ol\CH^\lambda_{G,x}\right)\cap 
\on{SingSupp}(\hl_x^*(\CF)\otimes \tau_x^*(\CV^\lambda_x)).$$

\sssec{} \label{sss:cotangent modifications}

Recall (see, for example, \cite[Formula (B.23)]{GKRV}) that for a point
$$\CP^{\on{loc}}_G\overset{\alpha}\rightsquigarrow \CP'{}^{\on{loc}}_G$$
of $\CH^{\on{loc}}_{G,x}$, the cotangent space
$$T^*_{\CP^{\on{loc}}_G\overset{\alpha}\rightsquigarrow \CP'{}^{\on{loc}}_G}(\CH^{\on{loc}}_{G,x})$$
identifies with the set of pairs
\begin{equation} \label{e:identify cotan Hecke}
A^{\on{loc}}\in \Gamma(\cD_x,\fg_{\CP^{\on{loc}}_G}\otimes \omega),\,\,A'{}^{\on{loc}}\in \Gamma(\cD_x,\fg_{\CP'{}^{\on{loc}}_G}\otimes \omega),
\end{equation} 
such that 
$$\alpha(A^{\on{loc}})=A'{}^{\on{loc}}$$ as elements of 
$\Gamma(\ocD_x,\fg_{\CP'{}^{\on{loc}}_G}\otimes \omega_X)$. 

\medskip

Furthermore, given 
$$A\in T_{\CP_G}^*(\Bun_G)\simeq \Gamma(X,\fg_{\CP_G}\otimes \omega),\,\,
A'\in T_{\CP'_G}^*(\Bun_G)\simeq \Gamma(X,\fg_{\CP'_G}\otimes \omega)$$
their images in $T^*_{\CP_G\overset{\alpha}\rightsquigarrow \CP'_G}(\CH_{G,x})$ differ
by the image of an element in $T^*_{\CP_G|_{\cD_x}\overset{\alpha}\rightsquigarrow \CP'_G|_{\cD_x}}(\CH^{\on{loc}}_{G,x})$ 
if and only if 
$$\alpha(A|_{X-x})=A'|_{X-x},$$
and in this case the corresponding element of $T^*_{\CP_G|_{\cD_x}\overset{\alpha}\rightsquigarrow \CP'_G|_{\cD_x}}(\CH^{\on{loc}}_{G,x})$
is given in terms of \eqref{e:identify cotan Hecke} by
$$A^{\on{loc}}:=A|_{\cD_x},\,\, A'{}^{\on{loc}}=A'|_{\cD_x}.$$

\sssec{}

Hence, in order to prove (a'), it suffices to show that for a point
$$(\CP^{\on{loc}}_G\overset{\alpha}\rightsquigarrow \CP'{}^{\on{loc}}_G)\in \CH^{\on{loc}}_{G,x}$$
induced by a point
$$(\CP^{\on{loc}}_M\overset{\alpha}\rightsquigarrow \CP'{}^{\on{loc}}_M)\in \CH^{\on{loc}}_{M,x},$$
corresponding to $\lambda$ (see \secref{sss:specified Hecke}), \emph{any} pair
$$(A^{\on{loc}},A'{}^{\on{loc}})\in T^*_{\CP^{\on{loc}}_G\overset{\alpha}\rightsquigarrow \CP'{}^{\on{loc}}_G}(\CH^{\on{loc}}_{G,x})$$
belongs to $\on{SingSupp}(\CV^\lambda_x)$. 

\medskip

We identify 
$$\CH^{\on{loc}}_{G,x}=G\qqart\backslash G\ppart/G\qqart$$ 
so that the point $\CP^{\on{loc}}_G\overset{\alpha}\rightsquigarrow \CP'{}^{\on{loc}}_G$ corresponds to $t^\lambda$. 

\medskip

Recall that $\CV^\lambda_x$ is the IC-sheaf on the closure of the double coset of
$$t^\lambda \in G\qqart\backslash G\ppart/G\qqart.$$ 

\medskip

Hence, the fiber of $\on{SingSupp}(\CV^\lambda_x)$ at $t^\lambda$ is the conormal to this double coset,
and hence equals the entire cotangent space at this point. 

\sssec{}

To prove point (b), we mimic the argument of \cite[Sect. B.6.7]{GKRV}. We consider $\CH^{\on{loc}}_G$,
equipped with its natural crystal structure along $X$, and the corresponding splitting of the short exact sequence
$$0\to T^*_x(X) \to T^*_{\CP^{\on{loc}}_G\overset{\alpha}\rightsquigarrow \CP'{}^{\on{loc}}_G}(\CH^{\on{loc}}_G) \to
T^*_{\CP^{\on{loc}}_G\overset{\alpha}\rightsquigarrow \CP'{}^{\on{loc}}_G}(\CH^{\on{loc}}_{G,x})\to 0,$$
i.e.,
$$T^*_{\CP^{\on{loc}}_G\overset{\alpha}\rightsquigarrow \CP'{}^{\on{loc}}_G}(\CH^{\on{loc}}_G)\simeq 
T^*_x(X) \oplus T^*_{\CP^{\on{loc}}_G\overset{\alpha}\rightsquigarrow \CP'{}^{\on{loc}}_G}(\CH^{\on{loc}}_{G,x}).$$

It suffices to show that, in terms of this identification, for an element
$$(\xi_x,(A^{\on{loc}},A'{}^{\on{loc}}))\in T^*_{\CP^{\on{loc}}_G\overset{\alpha}\rightsquigarrow \CP'{}^{\on{loc}}_G}(\CH^{\on{loc}}_G)$$
that belongs to $\on{SingSupp}(\CV^\lambda)$, 
we have
\begin{equation} \label{e:xi x}
\xi_x:=\langle A^{\on{loc}}_x,\lambda \rangle,
\end{equation} 
where $A^{\on{loc}}_x$ is the value of $A^{\on{loc}}$ at $x$. 

\medskip

The assertion is local, so we can assume that $X$ is $\BA^1$, with coordinate $t$. This allows us to trivialize
the line $T^*_x(X)$. Further, we can assume that $\CP^{\on{loc}}_G$ is trivial. Then we can think of 
$$A^{\on{loc}}\in \Gamma(\cD_x,\fg_{\CP^{\on{loc}}_G}\otimes \omega)$$
as an element of $\fg\qqart$. 

\medskip

By \cite[Formula (B.33)]{GKRV}, the element $\xi_x$ equals 
$$\on{Res}_x(A^{\on{loc}}, \lambda\cdot \frac{dt}{t}),$$
whence \eqref{e:xi x}. 

\ssec{Proof of point (c) and affine Springer fibers}

To prove point (c), it suffices to show:

\medskip

\noindent(c') The point $(A',(\CP_G\rightsquigarrow \CP'_G))$ is isolated in the intersection
\begin{equation}  \label{e:Springer global lambda}
\left(T^*(\Bun_G\times \CH_{G,x}^{\on{loc}})\underset{\Bun_G\times \CH_{G,x}^{\on{loc}},\hl_x\times \tau_x}\times \ol\CH^\lambda_{G,x}\right)\cap 
\left(A' \times (\hr_x)^{-1}(\CP'_G)\right) \subset T^*(\CH_{G,x}).
\end{equation}

Point (d) in \secref{sss:points a-d} is proved similarly. 

\sssec{}

Consider first the larger intersection
\begin{equation}  \label{e:Springer global}
\left(T^*(\Bun_G\times \CH_{G,x}^{\on{loc}})\underset{\Bun_G\times \CH_{G,x}^{\on{loc}},\hl_x\times \tau_x}\times \CH_{G,x}\right)\cap 
\left(A' \times (\hr_x)^{-1}(\CP'_G)\right) \subset T^*(\CH_{G,x}).
\end{equation}

\medskip

By \secref{sss:cotangent modifications}, the scheme in \eqref{e:Springer global} is the space of modifications 
of $\CP'_G|_{\cD_x}\rightsquigarrow \CP_G^{\on{loc}}$, for which the element
$$A'|_{\cD_x}\in \Gamma(\cD_x,\fg_{\CP'_G}\otimes \omega)\subset 
\Gamma(\ocD_x,\fg_{\CP'_G}\otimes \omega)\simeq\Gamma(\ocD_x,\fg_{\CP^{\on{loc}}_G}\otimes \omega) $$
belongs to
$$\Gamma(\cD_x,\fg_{\CP^{\on{loc}}_G}\otimes \omega)\subset \Gamma(\ocD_x,\fg_{\CP^{\on{loc}}_G}\otimes \omega).$$

\medskip

Denote this space by $\on{Spr}_{G,A'}$: it is isomorphic to a parahoric affine Springer fiber over the element $A'$. 
Denote the intersection \eqref{e:Springer global lambda} by $\on{Spr}^{\leq\lambda}_{G,A'}$. 

\sssec{}

If we trivialize $\CP'{}^{\on{loc}}_G:=\CP'_G|_{\cD_x}$, we can think of $\on{Spr}_{G,A'}$ as a (closed) subscheme in $\Gr_G$, and we have
$$\on{Spr}^{\leq\lambda}_{G,A'}=\on{Spr}_{G,A'}\cap \ol\Gr^\lambda_G.$$

\medskip

We need to show that our particular point 
\begin{equation} \label{e:isol point}
\CP'_G|_{\cD_x}\rightsquigarrow \CP_G|_{\cD_x}
\end{equation} 
is isolated in $\on{Spr}^{\leq\lambda}_{G,A'}$. 

\sssec{}

By \secref{sss:specified Hecke}, the $G$-bundle $\CP'{}^{\on{loc}}_G$ on $\cD_x$ is equipped with a reduction to $M$, 
denoted $\CP'{}^{\on{loc}}_M$ and $A'|_{\cD_x}$
belongs to
$$\Gamma(\cD_x,\fm_{\CP'{}^{\on{loc}}_M}\otimes \omega).$$

So along with
$\on{Spr}_{G,A'}$, we can consider its variant for $M$, to be denoted $\on{Spr}_{M,A'}$. Since
$$\Gr_M\to \Gr_G$$
is a closed embedding, so is the embedding $\on{Spr}_{M,A'}\hookrightarrow \on{Spr}_{G,A'}$. 

\medskip

We claim:

\begin{prop}  \label{p:Springer}
The inclusion $\on{Spr}_{M,A'}\hookrightarrow \on{Spr}_{G,A'}$ is an equality.
\end{prop} 

\begin{rem}
For our purposes, which is proving that \eqref{e:isol point} is isolated in $\on{Spr}^{\leq\lambda}_{G,A'}$,
we only need the assertion \propref{p:Springer} at the level of sets of $k$-points. 
\end{rem}

\sssec{}

Let us show how \propref{p:Springer} implies that \eqref{e:isol point} is isolated in $\on{Spr}^{\leq\lambda}_{G,A'}$. 

\medskip

By \propref{p:Springer}, it suffices to show that the point $t^\lambda$ is isolated in
\begin{equation} \label{e:Gr M}
\ol\Gr_G^\lambda \cap \Gr_M.
\end{equation}

Since $t^\lambda$ belongs to 
\begin{equation} \label{e:Gr M bis}
\Gr_G^\lambda \cap \Gr_M,
\end{equation}
it suffices to show that it is isolated in \eqref{e:Gr M bis}.

\medskip

Note, however, that the intersection $\Gr_G^\lambda \cap \Gr_M$ is the union of $M\qqart$-orbits $\Gr_M^\mu$ 
over $M$-dominant coweights $\mu$ for which there exists $w\in W$ such that 
$$\mu=w(\lambda).$$

Note that the point $t^\lambda$ equals $\Gr_M^\lambda$, because $\lambda$ is a coweight of $Z(M)$. 
The assertion follows now from the regularity assumption on $\lambda$: the orbit $\Gr_M^\lambda$ 
belongs to a different connected component of $\Gr_M$ than the other $\Gr_M^\mu$ with $\mu=w(\lambda)$. 

\ssec{Proof of \propref{p:Springer}} 

\sssec{}

For a prestack $\CY$ denote by
$$\bMaps(\cD_x,\CY) \text{ and } \bMaps(\ocD_x,\CY)$$
the corresponding prestacks of arcs and loops into $\CY$, respectively:

$$\Maps(\Spec(R),\bMaps(\cD_x,\CY))=\Maps(\Spec(R\qqart),\CY)$$
and 
$$\Maps(\Spec(R),\bMaps(\ocD_x,\CY))=\Maps(\Spec(R\ppart),\CY).$$

%
%
%
%
%

\sssec{} 

Choose a trivialization of $\omega|_{\cD_x}$. Thus, we can think of the pair $(\CP'{}^{\on{loc}}_M,A')$ 
as a map 
\begin{equation}\label{e:Higgs'}
\cD_x\to \fm/\on{Ad}(M),
\end{equation}
which is $(G,M)$-regular; that is, it is in fact a map to $\ofm/\on{Ad}(M)$.
The map between the Springer fibers in \propref{p:Springer} can be written explicitly as

\medskip

\[\bMaps(\cD_x,\fm/\on{Ad}(M))\underset{\bMaps(\ocD_x,\fm/\on{Ad}(M))}
\times\{*\}\to
\bMaps(\cD_x,\fg/\on{Ad}(G))\underset{\bMaps(\ocD_x,\fg/\on{Ad}(G))}
\times\{*\};\]
we need to show that the map is an isomorphism. 
Here the map $\{*\}\to\bMaps(\ocD_x,\fm/\on{Ad}(M))$ is given by
the restriction of \eqref{e:Higgs'} to $\ocD_x$.

\sssec{} Notice that the map
\[\bMaps(\cD_x,\ofm/\on{Ad}(M))\underset{\bMaps(\ocD_x,\ofm/\on{Ad}(M))}
\times\{*\}\to\bMaps(\cD_x,\fm/\on{Ad}(M))\underset{\bMaps(\ocD_x,\fm/\on{Ad}(M))}
\times\{*\},\]
is an isomorphism. Indeed, using \eqref{e:am Cart}, it is obtained by base change from the map
\begin{equation}\label{e:Maps to a}
\bMaps(\cD_x,\ofa_M)\underset{\bMaps(\ocD_x,\ofa_M)}\times\{*\}\to
\bMaps(\cD_x,\fa_M)\underset{\bMaps(\ocD_x,\fa_M)}\times\{*\},
\end{equation}
which is an isomorphism. 

\sssec{} It remains to show that the map
\[\bMaps(\cD_x,\ofm/\on{Ad}(M))\underset{\bMaps(\ocD_x,\ofm/\on{Ad}(M))}
\times\{*\}\to
\bMaps(\cD_x,\fg/\on{Ad}(G))\underset{\bMaps(\ocD_x,\fg/\on{Ad}(G))}
\times\{*\}\]
is an isomorphism. By \lemref{l:regularity bis}(c), the map is obtained by base change from the map
\begin{equation}\label{e:Maps to a bis}
\bMaps(\cD_x,\ofa_M)\underset{\bMaps(\ocD_x,\ofa_M)}
\times\{*\}\to
\bMaps(\cD_x,\fa)\underset{\bMaps(\ocD_x,\fa)}
\times\{*\}.
\end{equation}
Therefore, it suffices to prove that \eqref{e:Maps to a bis} is an isomorphism. 
This is clear at the classical level: as both $\fa$ and $\ofa_M$ are separated, for a classical affine scheme $S$,
the sets of $S$-points in both the source and the target of \eqref{e:Maps to a bis} 
are singleton sets. (As was mentioned above, \propref{p:Springer} on the level
of sets of $k$-points actually suffices for our purposes.)

\medskip

To complete the proof, we notice that \eqref{e:Maps to a bis} is formally \'etale, because $\ofa_M\to\fa$ is \'etale.

\ssec{Proof of \thmref{t:lisse prel} for non-holonomic D-modules} \label{ss:lisse for D-mod} 

We will deduce the assertion of \thmref{t:lisse prel} for $\Dmod(-)$ from its validity for the subcategory
$\Shv(-)$ consisting of objects with regular holonomic cohomologies. 

\medskip

The proof is based on considering field extensions of the initial ground field $k$
(cf. the proof of Observation \ref{o:de Rham GLC}). 

\sssec{}

By \propref{p:lisse subcategory}, it suffices to show that the inclusion
\begin{multline*} 
\Shv_\Nilp(\Bun_G)\simeq 
\Shv_\Nilp(\Bun_G)\underset{\QCoh(\LocSys_\cG^{\on{restr}}(X))}\otimes \QCoh(\LocSys_\cG^{\on{restr}}(X)) \simeq \\
\simeq \Shv_\Nilp(\Bun_G)\underset{\QCoh(\LocSys_\cG(X))}\otimes \QCoh(\LocSys_\cG^{\on{restr}}(X)) \hookrightarrow \\
\hookrightarrow   \Dmod(\Bun_G)\underset{\QCoh(\LocSys_\cG(X))}\otimes \QCoh(\LocSys^{\on{restr}}_\cG(X)) 
\end{multline*}
is an equality.

\medskip

Let $S$ be an affine scheme mapping to $\LocSys_\cG^{\on{restr}}(X)$. It suffices to show
that the inclusion
\begin{equation} \label{e:holonom to all}
\Shv_\Nilp(\Bun_G)\underset{\QCoh(\LocSys_\cG(X))}\otimes \QCoh(S) \hookrightarrow 
\Dmod(\Bun_G)\underset{\QCoh(\LocSys_\cG(X))}\otimes \QCoh(S) 
\end{equation} 
is an equality\footnote{Indeed, $\QCoh(\LocSys_\cG(X))$ is rigid and $\Shv_\Nilp(\Bun_G)$ and 
$\Dmod(\Bun_G)$ are dualizable, so the operations $\Shv_\Nilp(\Bun_G)\underset{\QCoh(\LocSys_\cG(X))}\otimes -$
and $\Dmod(\Bun_G)\underset{\QCoh(\LocSys_\cG(X))}\otimes -$ commute with limits.}. 

\sssec{}

Let $k\subset k'$ be a field extension. Let $X'$ (resp., $S'$, $\Bun'_G$) be the base change of $X$ (resp., $S$, $\Bun_G$) 
from $k$ to $k'$. Note that for any prestack $\CY$ over $k$ and its base change $\CY'$ to $k'$, we have
\begin{equation} \label{e:base change D-mod}
\Dmod(\CY')\simeq \Dmod(\CY) \underset{\Vect_k}\otimes \Vect_{k'}
\end{equation} 
and 
\begin{equation} \label{e:base change LocSys}
\LocSys_\cG(X')\simeq \LocSys_\cG(X)\underset{\Spec(k)}\times \Spec(k').
\end{equation}

For a fixed $\CN\subset T^*(\CY)$, we have a fully faithful embedding
\begin{equation} \label{e:base change D-mod N}
\Dmod_\CN(\CY) \underset{\Vect_k}\otimes \Vect_{k'} \hookrightarrow 
\Dmod_{\CN'}(\CY'),
\end{equation} 
but which is no longer an equivalence (indeed, for example, for $\CN=\{0\}$, there are
more local systems over $k'$ than over $k$).  

\medskip

From \eqref{e:base change D-mod} and \eqref{e:base change LocSys} we obtain an equivalence
$$\left(\Dmod(\Bun_G)\underset{\QCoh(\LocSys_\cG(X))}\otimes \QCoh(S) \right)
\underset{\Vect_k}\otimes \Vect_{k'}\simeq
\Dmod(\Bun'_G)\underset{\QCoh(\LocSys_\cG(X'))}\otimes \QCoh(S') $$
and from \eqref{e:base change D-mod N} a fully faithful embedding 
\begin{multline} \label{e:Hecke base change}
\left(\Shv_\Nilp(\Bun_G)\underset{\QCoh(\LocSys_\cG(X))}\otimes \QCoh(S)\right)
\underset{\Vect_k}\otimes \Vect_{k'} \hookrightarrow \\
\to \Shv_\Nilp(\Bun'_G)\underset{\QCoh(\LocSys_\cG(X'))}\otimes \QCoh(S').
\end{multline}

However, we claim:

\begin{lem} \label{l:Hecke base change}
The fully faithful embedding \eqref{e:Hecke base change} is an equivalence. 
\end{lem}

\begin{proof}

We will show that the image of the functor \eqref{e:Hecke base change} contains the generators
of the target category. 

\medskip

Indeed, let $y_i\in \Bun_G$ be as \secref{sss:generators Nilp}. Let $y'_i$ be the corresponding $k'$-points
of $\Bun'_G$. Then the generators of 
$$\Shv_\Nilp(\Bun_G)\underset{\QCoh(\LocSys_\cG(X))}\otimes \QCoh(S)$$
are given by $\sP^{\on{enh}}_S(\delta_{y_i})$, and the generators of 
$$ \Shv_\Nilp(\Bun'_G)\underset{\QCoh(\LocSys_\cG(X'))}\otimes \QCoh(S')$$
are given by $\sP^{\on{enh}}_{S'}(\delta_{y'_i})$. Now these generators are sent to one another
by the functor \eqref{e:Hecke base change}. 

\end{proof} 

\sssec{}

We are now ready to prove that \eqref{e:holonom to all} is an equality. Let $\CF$ be an object in the right-hand
side, which is right-orthogonal to the left-hand side. By \lemref{l:Hecke base change} for any $k\subset k'$,
the pullback $\CF'$ to $\Bun'_G$ will have the same property.

\medskip

We now clam that for any $\CF$ as above, its image in $\Dmod(\Bun_G)$ is right-orthogonal to $\Shv(\Bun_G)$.
Indeed, for any $\CF_1\in \Shv(\Bun_G)$, we have
$$\CHom_{\Dmod(\Bun_G)}(\CF_1,\CF)
\simeq \CHom_{\Dmod(\Bun_G)\underset{\QCoh(\LocSys_\cG(X))}\otimes \QCoh(S) }(\sP^{\on{enh}}_S(\CF_1),\CF),$$
while 
$$\sP^{\on{enh}}_S(\CF_1)\in \Shv_\Nilp(\Bun_G)\underset{\QCoh(\LocSys_\cG(X))}\otimes \QCoh(S) .$$

\medskip

Hence, we obtain that for $\CF$ as above and any $k\subset k'$, the corresponding object $\CF'\in \Dmod(\Bun'_G)$
is right-orthogonal to $\Shv(\Bun'_G)$. 

\medskip

We wish to show that $\CF=0$. It suffices to show that the image of $\CF$ in $\Dmod(\Bun_G)$ is zero. This follows from the next assertion: 

\begin{lem}
Let $\CF\in \Dmod(\CY)$ be such that for any $k\subset k'$, the corresponding object $\CF'\in \Dmod(\CY')$
is right-orthogonal to $\Shv(\CY')$. Then $\CF=0$.
\end{lem} 

\begin{proof} 

Let $\CF\neq 0$. Consider the underlying object $\oblv_{\Dmod}(\CF)\in \QCoh(\CY)$. Then we can find a
geometric point
$$\Spec(k')\overset{\bi_y}\to \CY,$$
so that $\bi_y^*(\oblv_{\Dmod}(\CF))\neq 0$. 

\medskip

Let $\bi_{y'}$ denote the resulting geometric point of $\CY'$. Then $\bi_{y'}^*(\oblv_{\Dmod}(\CF'))\neq 0$.
However, the latter means that
$$\CHom_{\Dmod(\CY')}(\delta_{y'},\CF')\neq 0.$$

\end{proof}

\ssec{Proof of \thmref{t:projector Betti pt}} \label{ss:proof projector Betti pt}

We retain the notations of \secref{ss:left adjoint Betti}. 

\sssec{}

We need to show that $$\sP^{\on{enh}}_{\LocSys^{\on{Betti}}_\cG(X)}(\delta_y)\in \Shv^{\on{all}}(\Bun_G)^{\on{Hecke-loc.const.}}$$
belongs to $\Shv^{\on{all}}_{\Nilp}(\Bun_G)$. This is equivalent to showing that the object
$$\sP_{\LocSys^{\on{Betti}}_\cG(X)}(\delta_y)\in \Shv^{\on{all}}(\Bun_G) \otimes \QCoh(\LocSys^{\on{Betti}}_\cG(X))$$
belongs to
$$\Shv^{\on{all}}_{\Nilp}(\Bun_G) \otimes \QCoh(\LocSys^{\on{Betti}}_\cG(X)) \subset 
\Shv^{\on{all}}(\Bun_G) \otimes \QCoh(\LocSys^{\on{Betti}}_\cG(X)).$$

\medskip

Furthermore, the latter is equivalent to showing that the object 
$$\sP_{\LocSys^{\on{Betti},\on{rigid}_x}_\cG(X)}(\delta_y)\in 
\Shv^{\on{all}}(\Bun_G) \otimes \QCoh(\LocSys^{\on{Betti},\on{rigid}_x}_\cG(X))$$
belongs to
$$\Shv^{\on{all}}_{\Nilp}(\Bun_G) \otimes \QCoh(\LocSys^{\on{Betti},\on{rigid}_x}_\cG(X)) \subset 
\Shv^{\on{all}}(\Bun_G)\otimes \QCoh(\LocSys^{\on{Betti},\on{rigid}_x}_\cG(X)).$$

\sssec{}

Since $\LocSys^{\on{Betti},\on{rigid}_x}_\cG(X)$ is an eventually coconnective affine scheme, its $\QCoh$
is generated under colimits by objects of the form 
$$\wt\bi_*(\wt\sfe'),$$
where
$$\wt\bi:\Spec(\wt\sfe')\to \LocSys^{\on{Betti},\on{rigid}_x}_\cG(X)$$
and $\wt\sfe'\supset \sfe$ are the residue fields of scheme-theoretic points of $\LocSys^{\on{Betti},\on{rigid}_x}_\cG(X)$.

\medskip

In particular, $\CO_{\LocSys^{\on{Betti},\on{rigid}_x}_\cG(X)}$ can be expressed as a colimit of such objects. 

\medskip

Hence, it suffices to show that all 
$$\sP_{\LocSys^{\on{Betti},\on{rigid}_x}_\cG(X)}(\delta_y)\otimes \wt\bi_*(\wt\sfe')\in 
\Shv^{\on{all}}(\Bun_G) \otimes \QCoh(\LocSys^{\on{Betti},\on{rigid}_x}_\cG(X))$$
belong to 
$$\Shv^{\on{all}}_{\Nilp}(\Bun_G) \otimes \QCoh(\LocSys^{\on{Betti},\on{rigid}_x}_\cG(X)).$$

Note, however, that
$$\sP_{\LocSys^{\on{Betti},\on{rigid}_x}_\cG(X)}(\delta_y)\otimes \wt\bi_*(\wt\sfe')
\simeq (\on{Id}\otimes \wt\bi_*)(\sP_{\Spec(\wt\sfe')}(\delta_y)).$$

Hence, it suffices to show that the objects 
$$\sP_{\Spec(\wt\sfe')}(\delta_y)\in \Shv^{\on{all}}(\Bun_G) \otimes \Vect_{\wt\sfe'}$$
belong to
$$\Shv_\Nilp^{\on{all}}(\Bun_G)\otimes \Vect_{\wt\sfe'}
\subset \Shv^{\on{all}}(\Bun_G) \otimes \Vect_{\wt\sfe'}.$$

\sssec{}

Let $\sfe'\supset \wt\sfe'$ be the algebraic closure of $\wt\sfe'$. It is easy to see that it is sufficient to show that 
the objects 
$$\sP_{\Spec(\sfe')}(\delta_y)\in \Shv^{\on{all}}(\Bun_G)\otimes \Vect_{\sfe'}$$
belong to
$$\Shv_\Nilp^{\on{all}}(\Bun_G)\otimes \Vect_{\sfe'}
\subset \Shv^{\on{all}}(\Bun_G)\otimes \Vect_{\sfe'}.$$

\sssec{}

Note, however, that we have canonical identifications
\begin{equation} \label{e:BC all}
\Shv^{\sfe,\on{all}}(\Bun_G)\otimes \Vect_{\sfe'}\simeq \Shv^{\sfe',\on{all}}(\Bun_G), \quad
\Shv_\Nilp^{\sfe,\on{all}}(\Bun_G)\otimes \Vect_{\sfe'}\simeq \Shv_\Nilp^{\sfe',\on{all}}(\Bun_G)
\end{equation}
and
\begin{equation} \label{e:BC Betti}
\LocSys^{\on{Betti},\sfe}_\cG(X)\underset{\Spec(\sfe)}\times \Spec(\sfe')\simeq
\LocSys^{\on{Betti},\sfe'}_\cG(X),
\end{equation}
where the superscripts $\sfe$ and $\sfe'$ indicate the fields of coefficients of our sheaves.

\medskip

We can regard 
$$\bi:\Spec(\sfe')\to \LocSys^{\on{Betti},\sfe}_\cG(X)$$
as an $\sfe'$-point 
$$\bi':\Spec(\sfe')\to \LocSys^{\on{Betti},\sfe'}_\cG(X).$$

Under these identifications, the functor 
$$\sP_{\Spec(\sfe')}:\Shv^{\sfe,\on{all}}(\Bun_G)\to \Shv^{\on{all}}(\Bun_G)\otimes \Vect_{\sfe'}$$
identifies with
$$\Shv^{\sfe,\on{all}}(\Bun_G) \overset{-\otimes \sfe'}\longrightarrow \Shv^{\sfe',\on{all}}(\Bun_G)
\overset{\sP_{\Spec(\sfe')}}\longrightarrow \Shv^{\sfe',\on{all}}(\Bun_G),$$
where the second arrow is Beilinson's projector corresponding to the point $\bi'$.

\sssec{}

However, $\bi'$, being a rational point, factors as
$$\Spec(\sfe')\to \LocSys^{\on{restr},\on{rigid}_x,\sfe'}_\cG(X)\to \LocSys^{\on{Betti},\sfe'}_\cG(X).$$

Hence,
$$\sP_{\Spec(\sfe')}(\delta_y)\in \Shv_\Nilp^{\sfe'}(\Bun_G) \subset \Shv_\Nilp^{\sfe',\on{all}}(\Bun_G),$$
by \corref{c:P Z}. 

\begin{rem}
Note, however, that although we have an equivalence \eqref{e:BC all}, the functor 
$$\Shv_\Nilp^{\sfe}(\Bun_G)\otimes \Vect_{\sfe'}\subset \Shv_\Nilp^{\sfe'}(\Bun_G)$$
is a \emph{proper containment}.

\medskip

Similarly, although we have an isomorphism \eqref{e:BC Betti}, the map 
$$\LocSys^{\on{restr},\sfe}_\cG(X)\underset{\Spec(\sfe)}\times \Spec(\sfe')\to 
\LocSys^{\on{restr},\sfe'}_\cG(X),$$
is an embedding of a union of connected components, 
but \emph{not} an isomorphism.
\end{rem}

\newpage

\centerline{\bf Part IV: Langlands theory with nilpotent singular support} 

\bigskip

Let us make a brief overview of the contents of this Part.

\medskip

In \secref{c:GLC} we state the (categorical) Geometric Langlands conjecture with restricted variation:
$$\Shv_{\on{Nilp}}(\Bun_G) \simeq \IndCoh_{\on{Nilp}}(\on{LocSys}_\cG^{\on{restr}}(X)),$$
and compare it to the de Rham and Betti versions of the GLC. A priori the restricted version follows from these other two.
However, we show that the restricted version is actually equivalent to the full de Rham version (assuming
Hypothesis \ref{h:Lde Rham GLC}). 

\medskip

In \secref{s:Trace} we formulate one of the key points of this paper, the Trace Conjecture. We start by reviewing
the \emph{local term} map
$$\on{LT}:\Tr((\Frob_\CY)_*,\Shv(\CY))\to \on{Funct}_c(\CY(\BF_q)),$$
where $\CY$ is an algebraic stack defined over $\BF_q$, but considered over $\ol\BF_q$. 
The Trace Conjecture says that the composition
$$\Tr((\Frob_{\Bun_G})_*,\Shv_\Nilp(\Bun_G))\to \Tr((\Frob_{\Bun_G})_*,\Shv(\Bun_G))
\overset{\on{LT}}\to \on{Funct}_c(\Bun_G(\BF_q))=:\on{Autom}$$
is an isomorphism. We then discuss a generalization of the Trace Conjecture that recovers cohomology of
shtukas also as traces of functors acting on $\Shv_\Nilp(\Bun_G)$.

\medskip

In \secref{s:abelian varieties} we make a digression and proof a version of the Trace Conjecture for the
category of lisse sheaves on an abelian variety. 

\medskip

In \secref{s:spectral} we explain how the Trace Conjecture allows us to recover V.~Lafforgue's spectral decomposition
of $\on{Autom}$ with respect to (the coarse moduli space of) Langlands parameters. 

\medskip

We start by defining the (prestack) 
$\on{LocSys}_\cG^{\on{arithm}}(X)$ as Frobenius-fixed points on $\on{LocSys}_\cG^{\on{restr}}(X))$;
in \thmref{t:Frob-finite} we prove that $\on{LocSys}_\cG^{\on{arithm}}(X)$
is actually an algebraic stack. 

\medskip

We view $(\QCoh(\on{LocSys}_\cG^{\on{restr}}(X)),\Shv_\Nilp(\Bun_G))$ 
as a pair of a monoidal category with its module category, equipped with endofunctors (both given by Frobenius).
In this case we can consider
$$\on{cl}(\Shv_\Nilp(\Bun_G),(\Frob_{\Bun_G})_*)\in \on{HH}_\bullet(\QCoh(\on{LocSys}_\cG^{\on{restr}}(X)),\Frob^*)$$
attached to this data (see \cite[Sect. 3.8.1]{GKRV}). We identify 
$$\on{HH}_\bullet(\QCoh(\on{LocSys}_\cG^{\on{restr}}(X)),\Frob^*)\simeq \QCoh(\on{LocSys}_\cG^{\on{arithm}}(X))$$
(see \secref{sss:enhanced trace setup Y}), and denote the resulting object 
$$\Drinf\in \QCoh(\on{LocSys}_\cG^{\on{arithm}}(X)).$$
By design (see \thmref{t:enhanced Tr Y}), we have
$$\Gamma(\on{LocSys}_\cG^{\on{arithm}}(X),\Drinf)\simeq \Tr((\Frob_{\Bun_G})_*,\Shv_\Nilp(\Bun_G)),$$
where the right-hand side is naturally acted on by 
$$\Exc:=\Gamma(\on{LocSys}_\cG^{\on{arithm}}(X),\CO_{\on{LocSys}_\cG^{\on{arithm}}(X)}).$$

Combining with the Trace Conjecture, we obtain an action of $\Exc$ on $\on{Autom}$, i.e., a spectral
decomposition of $\on{Autom}$ with respect to the coarse moduli space of $\on{LocSys}_\cG^{\on{arithm}}(X)$. 

\medskip

Finally, assuming the Geometric Langlands Conjecture (plus a more elementary \conjref{c:ignore N}),
we deduce an equivalence
$$\Drinf\simeq \omega_{\on{LocSys}_\cG^{\on{arithm}}(X)},$$
as objects of $\QCoh(\on{LocSys}_\cG^{\on{arithm}}(X))$. Combining with the Trace Conjecture, we obtain 
a conjectural identification
$$\on{Autom}\simeq \Gamma(\on{LocSys}_\cG^{\on{arithm}}(X),\omega_{\on{LocSys}_\cG^{\on{arithm}}(X)}),$$
i.e., a description of the space of (unramified) automorphic functions purely in terms of the stack of
Langlands parameters. 

\medskip

In \secref{s:arithm} we prove \thmref{t:Frob-finite}. The key tools for the proof are the properties of the map
$\brr$ from \thmref{t:coarse restr}, combined with results from \cite{De} and \cite{LLaf}.

\bigskip

\section{Geometric Langlands Conjecture with nilpotent singular support} \label{c:GLC}

In this section we formulate a version of the Geometric Langlands Conjecture that involves
$\Shv_\Nilp(\Bun_G)$. Its main feature is that it makes sense for any sheaf theory from
our list. 

\medskip

We then explain the relationship between this version of the Geometric Langlands Conjecture
and the de Rham and Betti versions. We will show that both these versions imply the one
with nilpotent singular support. 

\medskip

Vice versa, we show (under a certain plausible assumption, see Hypothesis \ref{h:Lde Rham GLC}) that the
restricted version actually implies the full de Rham version. 

\ssec{Digression: coherent singular support} \label{ss:coh sing supp}

In this subsection we will show how to adapt the theory of singular support, developed in \cite{AG}
for quasi-smooth \emph{schemes}, to the case of quasi-smooth \emph{formal schemes}. We will assume
that the reader is familiar with the main tenets of the paper \cite{AG}. 

\sssec{} \label{sss:formal q-smooth}

Let $\CY$ be a formal affine scheme (see Remark \ref{r:formal affine}), locally almost
of finite type as a prestack. 

%
%
%
%
%
%
%
%

\medskip

We shall say that $\CY$ is \emph{quasi-smooth} if for every $\sfe$-point $y$ of $\CF$, the cotangent space
$T_y^*(\CY)$ is acyclic off degrees $0$ and $-1$. 

\medskip

Equivalently, $\CY$ is quasi-smooth if 
$$T^*(\CY)|_{^{\on{red}}\CY}\in \Coh({}^{\on{red}}\CY)$$
can be locally written as a 2-step complex of vector bundles
$\CE_{-1}\to \CE_0$. 

\sssec{}

We will denote by 
$$T(\CY)|_{^{\on{red}}\CY}\in  \Coh({}^{\on{red}}\CY)^{\leq 1}$$
the \emph{naive} dual of $T^*(\CY)|_{^{\on{red}}\CY}$, i.e., 
$$T(\CY)|_{^{\on{red}}\CY}=\ul\Hom(T^*(\CY)|_{^{\on{red}}\CY},\CO_{^{\on{red}}\CY}).$$

\medskip

We define the (reduced) scheme $\on{Sing}(\CY)$ to be the reduced scheme underlying
$$\Spec_{^{\on{red}}\CY}(\Sym_{\CO_{^{\on{red}}\CY}}(T(\CY)|_{^{\on{red}}\CY}[1])).$$

Explicitly, $\sfe$-points of $\on{Sing}(\CY)$ are pairs $(y,\xi)$, where $y\in \CY(k)$, 
and $\xi\in H^{-1}(T_y^*(\CY))$. 

\sssec{} \label{sss:q-smooth as ind}

Let $\CY$ be quasi-smooth. It follows from \thmref{t:formal} that we can write 
$\CY$ as a filtered colimit
\begin{equation} \label{e:q-smooth as colimit}
\CY=\underset{i}{\on{colim}}\, Y_n,
\end{equation} 
where:

\begin{itemize}

\item $Y_n$ are quasi-smooth affine schemes;

\item The maps $Y_{n_1}\to Y_{n_2}$ are closed embeddings that induce isomorphisms $^{\on{red}}Y_{n_1}\to {}^{\on{red}}Y_{n_2}$;

\item For every $n$, the map $Y_n\to \CY$ is a closed embedding such that the induced map
\begin{equation} \label{e:sing map}
\on{Sing}(\CY)\underset{\CY}\times Y_n\to \on{Sing}(Y_n)
\end{equation}
is a closed embedding. 

\end{itemize} 

Indeed, in the notations of \eqref{e:presentation A}, let 
$$\CY_n:=\CY\underset{\BA^m}\times \{0\},$$
where the map $\CY\to \BA^m$ is
$$\CY\to \Spec(R) \overset{\{f_1^n,...,f_m^n\}}\longrightarrow \BA^m.$$

Then, on the one hand, $\CY_n$ is a quasi-smooth formal affine scheme (since $\CY$ is such), and
$$\on{Sing}(\CY)\underset{\CY}\times \CY_n\to \on{Sing}(\CY_n)$$
is a closed embedding. 

\medskip

On the other hand, 
$\CY_n$ is actually an affine scheme isomorphic to $\Spec(R_n)$. So we can take $Y_n:=\CY_n$. 

\sssec{}

Recall that for a prestack $\CY$ locally almost of finite type it makes sense to talk about the category
$\IndCoh(\CY)$, which is a module category over the (symmetric) monoidal category $\QCoh(\CY)$.

\medskip

If $\CY$ is an ind-scheme, we have a well-defined (small) subcategory
$$\Coh(\CY) \subset \IndCoh(\CY)^c,$$
so that $\IndCoh(\CY)$ identifies with the ind-completion of $\Coh(\CY)$,
see \cite[Sect. 2.4.3 and Proposition 2.4.6]{GR3}. 

\medskip

For $\CY$ written as \eqref{e:q-smooth as colimit}, we have
\begin{equation} \label{e:IndCoh as lim}
\IndCoh(\CY) \simeq \underset{n}{\on{lim}}\, \IndCoh(Y_n),
\end{equation}
where the limit is formed using the functors $i_{n_1,n_2}^!$ for $Y_{n_1}\overset{i_{n_1,n_2}}\longrightarrow Y_{n_2}$,
and also
$$\IndCoh(\CY) \simeq \underset{n}{\on{colim}}\, \IndCoh(Y_n),$$
where the colimit is formed inside $\DGCat$ using the functors $(i_{n_1,n_2})^\IndCoh_*$,
see \cite[Sect. 2.4.2]{GR3}.

\medskip

In terms of this identification, we have 
\begin{equation} \label{e:q-smooth as colimit Coh}
\Coh(\CY) \simeq \underset{n}{\on{colim}}\, \Coh(Y_n),
\end{equation}
where the colimit is formed using the *-pushforward functors, but inside the
$\infty$-category of \emph{not-necessarily cocomplete} DG categories. 

\sssec{}

The theory of singular support for quasi-smooth \emph{schemes} developed in \cite{AG} applies
``as-is" in the case of formal affine schemes  that are quasi-smooth.

\medskip

In particular, to an object
$$\CM\in \Coh(\CY),$$
one can attach its singular support $\on{SingSupp}(\CM)$, which is a conical Zariski-closed subset in 
$\on{Sing}(\CY)$. 

\medskip

Explicitly, for a given $\CM\in \Coh(\CY)$, the fiber of $\on{SingSupp}(\CM)$ over a given $\Spec(\sfe)\overset{\bi_y}\longrightarrow \CY$
is the support of 
$$\underset{n}\oplus\, H^n(\bi_y^!(\CM)),$$ 
viewed as a module over the algebra
$$\Sym^n(H^1(T_y(\CY))),$$
where the action is defined as in \cite[Sect. 6.1.1]{AG}. 

\sssec{}

For a given conical Zariski-closed subset $\CN\subset \on{Sing}(\CY)$, we can talk about a full subcategory
$$\Coh_\CN(\CY)\subset \Coh(\CY),$$
consisting of objects whose singular support is contained in $\CN$. 
We denote by $\IndCoh_\CN(\CY)$ its ind-completion, which is a full subcategory in $\IndCoh(\CY)$. 

\medskip

One can describe the category $\IndCoh_\CN(\CY)$ in terms of \eqref{e:IndCoh as lim} as follows: 
Given $\CN\subset \on{Sing}(\CY)$, let 
$\CN_n\subset \on{Sing}(Y_n)$ be the image of 
$$\CN\underset{\CY}\times Y_n\subset \on{Sing}(\CY)\underset{\CY}\times Y_n$$
under the map \eqref{e:sing map}. Then for $Y_{n_1}\to Y_{n_2}$, the pullback functor
$$\IndCoh(Y_{n_2})\to \IndCoh(Y_{n_1})$$
sends
$$\IndCoh_{\CN_{n_2}}(Y_{n_2})\to \IndCoh_{\CN_{n_1}}(Y_{n_1})$$
(see \cite[Proposition 7.1.3.(a)]{AG}) and we have
\begin{equation} \label{e:IndCoh N as inv}
\IndCoh_\CN(\CY)\simeq \underset{n}{\on{lim}}\, \IndCoh_{\CN_n}(Y_n),
\end{equation} 
as subcategories in the two sides of \eqref{e:IndCoh as lim}.

\sssec{}

Finally, one checks, using \cite[Proposition 7.1.3]{AG} and base change, that for a pair of indices $n_1,n_2$ the composite functor
$$\IndCoh(Y_{n_1})\overset{\text{*-pushforward}}\longrightarrow \IndCoh(\CY)\overset{\text{!-pullback}}\longrightarrow \IndCoh(Y_{n_2})$$
sends
$$\IndCoh_{\CN_{n_1}}(Y_{n_1}) \to \IndCoh_{\CN_{n_2}}(Y_{n_2}).$$

This implies that the *-pushforward functors $\IndCoh(Y_n) \to \IndCoh(\CY)$
send
\begin{equation} \label{e:*-pushforward}
\IndCoh_{\CN_{n}}(Y_{n}) \to \IndCoh_{\CN}(\CY).
\end{equation} 

\medskip

This shows that the category $\IndCoh_\CN(\CY)$ is compactly generated, namely, by the essential
images of $\Coh_{\CN_{n}}(Y_{n})$ along the functors \eqref{e:*-pushforward}.

\begin{rem}
Note, however, that the individual functors
$$(i_{n_1,n_2})^\IndCoh_*:\IndCoh(Y_{n_1})\to \IndCoh(Y_{n_2})$$
do \emph{not} necessarily send $\IndCoh_{\CN_{n_1}}(Y_{n_1})$ to $\IndCoh_{\CN_{n_2}}(Y_{n_2})$.
\end{rem} 

\sssec{}

The action of $\QCoh(\CY)$ on $\IndCoh(\CY)$ preserves the subcategory $\IndCoh_\CN(\CY)$. 

\medskip

It follows formally from \eqref{e:IndCoh N as inv} that the action on the object $\omega_\CY\in \IndCoh(\CY)$
defines an equivalence from $\QCoh(\CY)$ onto the full subcategory
$$\IndCoh_{\{0\}}(\CY)\subset \IndCoh(\CY),$$
where $\{0\}\subset \Sing(\CY)$ is the zero-section. 

\sssec{} \label{sss:formal compl sing supp}

Suppose for a moment that we are given a quasi-smooth affine scheme $\CY'$ and a map
$\CY\to \CY'$, which is an ind-closed embedding and a formal isomorphism (see Remark \ref{r:union formal compl}).  

\medskip

Let $\CN'\subset \on{Sing}(\CY')$, and let $\CN:=\CN'|_\CY$. Then it is easy to see that the full subcategories
$$\IndCoh_\CN(\CY) \subset \IndCoh(\CY)$$
and
$$\IndCoh_{\CN'}(\CY')\cap  \IndCoh(\CY')_\CY\subset \IndCoh(\CY')_\CY$$
coincide under the identification
$$\IndCoh(\CY) \simeq \IndCoh(\CY')_\CY.$$

\ssec{Geometric Langlands Conjecture for $\Shv_{\on{Nilp}}(\Bun_G)$}

\sssec{} \label{sss:cotan LocSys}

Note that the identification of cotangent spaces of $\on{LocSys}_\cG^{\on{restr}}(X)$ given by \propref{p:def}(b) 
and \eqref{e:coinv as Cc} implies that for a $\sfe$-point of $\sigma$ of $\on{LocSys}_\cG^{\on{restr}}(X)$, we have
\begin{equation} \label{e:cotan LocSys}
T^*_\sigma(\on{LocSys}_\cG^{\on{restr}}(X))\simeq \on{C}^\cdot_c(X,\cg^\vee_\sigma\overset{*}\otimes \omega_X)[-1]
\simeq \on{C}^\cdot(X,\cg^\vee_\sigma)[1],
\end{equation} 
where the last isomorphism uses the fact that $X$ is a proper smooth curve. 

\medskip

In the above formula $\cg^\vee_\sigma$ is the local system on 
$X$ associated to $\sigma$ and $\cg^\vee\in \Rep(\cG)$.

\sssec{} \label{sss:LocSys q-sm}

In particular, we obtain that the cotangent spaces to $\on{LocSys}_\cG^{\on{restr}}(X)$ 
live in cohomological degrees $[-1,1]$. 

\medskip

This implies that the cotangent spaces to $\on{LocSys}_\cG^{\on{restr},\on{rigid}_x}(X)$ 
live in cohomological degrees $[-1,0]$. I.e., $\on{LocSys}_\cG^{\on{restr},\on{rigid}_x}(X)$ 
is a union of quasi-smooth formal affine schemes.

\medskip

This allows us to talk about 
$$\on{Sing}(\on{LocSys}_\cG^{\on{restr},\on{rigid}_x}(X)),$$
which is a (reduced) algebraic stack over $^{\on{red}}\!\on{LocSys}_\cG^{\on{restr},\on{rigid}_x}(X)$.

\medskip

Thus, to a given conical Zariski-closed $\CN'\subset \on{Sing}(\on{LocSys}_\cG^{\on{restr},\on{rigid}_x}(X))$ we can 
attach a full subcategory 
$$\IndCoh_{\CN'}(\on{LocSys}_\cG^{\on{restr},\on{rigid}_x}(X))\subset \IndCoh(\on{LocSys}_\cG^{\on{restr},\on{rigid}_x}(X)).$$

\sssec{}


The $\cG$-action on $\on{LocSys}_\cG^{\on{restr},\on{rigid}_x}(X)$ naturally extends to 
$\on{Sing}(\on{LocSys}_\cG^{\on{restr},\on{rigid}_x}(X))$.

\medskip

We define 
$$\on{Sing}(\on{LocSys}_\cG^{\on{restr}}(X)):=\on{Sing}(\on{LocSys}_\cG^{\on{restr},\on{rigid}_x}(X))/\cG,$$
which is a (reduced) algebraic stack over $^{\on{red}}\!\on{LocSys}_\cG^{\on{restr}}(X)$.

\medskip

Given a conical Zariski-closed $$\CN\subset \on{Sing}(\on{LocSys}_\cG^{\on{restr}}(X)),$$ let $\CN'$
be its preimage in $\on{Sing}(\on{LocSys}_\cG^{\on{restr},\on{rigid}_x}(X))$. Define 
$$\IndCoh_\CN(\on{LocSys}_\cG^{\on{restr}}(X))$$ 
as the full subcategory of $\IndCoh(\on{LocSys}_\cG^{\on{restr}}(X))$ consisting of objects whose
*- (or, equivalently, -!) pullback along the (smooth) projection
$$\on{LocSys}_\cG^{\on{restr}}(X) \simeq \on{LocSys}_\cG^{\on{restr},\on{rigid}_x}(X)/\cG$$
belongs to the full subcategory $\IndCoh_{\CN'}(\on{LocSys}_\cG^{\on{restr},\on{rigid}_x}(X))$
of $\IndCoh(\on{LocSys}_\cG^{\on{restr},\on{rigid}_x}(X))$. 

\sssec{} \label{sss:Arth restr}

Note that the identification of cotangent spaces of $\on{LocSys}_\cG^{\on{restr}}(X)$ given by \eqref{e:cotan LocSys}
allows us to identify $\sfe$-points of 
$$\on{Arth}_\cG(X):=\on{Sing}(\on{LocSys}_\cG^{\on{restr}}(X))$$ 
with pairs
$$(\sigma,A),$$
where:

\smallskip

\begin{itemize}

\item $\sigma$ is a $\cG$-point of $\on{LocSys}_\cG^{\on{restr}}(X)$;

\item $A$ is an element in $H^0(X,\cg^\vee_\sigma)$.

\end{itemize} 

\sssec{} \label{sss:Nilp spec}

Let 
$$\on{Nilp}\subset \on{Arth}_\cG(X)$$ be the closed subset whose $\sfe$-points consist of pairs $(\sigma,A)$ for which $A$ 
is nilpotent. Thus, we can consider the fullcategory
$$\IndCoh_{\on{Nilp}}(\on{LocSys}_\cG^{\on{restr}}(X)) \subset \IndCoh(\on{LocSys}_\cG^{\on{restr}}(X)).$$

\sssec{}

We propose the following ``restricted" version of the Geometric Langlands Conjecture:

\begin{mainconj}  \label{c:restr GLC}
There exists a canonical equivalence
$$\Shv_{\on{Nilp}}(\Bun_G) \simeq \IndCoh_{\on{Nilp}}(\on{LocSys}_\cG^{\on{restr}}(X)),$$
compatible with the action of $\QCoh(\LocSys_\cG^{\on{restr}}(X))$ on both sides. 
\end{mainconj} 

\sssec{Example}

Let $X$ have genus zero. Then the inclusion 
$$\Shv_{\on{Nilp}}(\Bun_G) \subset \Shv(\Bun_G),$$
is an equality, and $\on{LocSys}_\cG^{\on{restr}}(X)$ is actually an algebraic stack isomorphic to
$$(\on{pt}\underset{\cG}\times \on{pt})/\cG.$$

\medskip

In this case, the assertion of \conjref{c:restr GLC} is known: it follows from the (derived) geometric Satake.


\sssec{Example} \label{sss:ex Gm}

Let us see what \conjref{c:restr GLC} says for $G=\BG_m$. Note that in this case 
$$\on{Nilp}\subset \on{Arth}_\cG(X)$$
is the $0$-section, so
$$\IndCoh_{\on{Nilp}}(\on{LocSys}^{\on{restr}}_\cG(X))=\QCoh(\on{LocSys}^{\on{restr}}_{\BG_m}(X)).$$

Given a $\BG_m$-local system $\sigma$, let $\on{LocSys}^{\on{restr}}_{\BG_m}(X)_\sigma$
be the corresponding component of $\on{LocSys}^{\on{restr}}_{\BG_m}(X)$. Let
$$\qLisse(\on{Pic})_\sigma\subset \qLisse(\on{Pic})$$
be the corresponding direct factor of $\qLisse(\on{Pic})$, see \secref{sss:decomp Gm}.

\medskip

Thus, \conjref{c:restr GLC} says that we have an equivalence
\begin{equation} \label{e:restr CFT}
\qLisse(\on{Pic})_\sigma \simeq \QCoh(\on{LocSys}^{\on{restr}}_{\BG_m}(X)_\sigma).
\end{equation}

Let us show how to establish the isomorphism \eqref{e:restr CFT} directly. 

\medskip

Up to translation by $\sigma$ on $\on{LocSys}^{\on{restr}}_{\BG_m}(X)$ (resp.,
tensor product by $E_\sigma$ on $\on{Pic}$), we can assume that $\sigma$ is trivial, so
$E_\sigma$ is the constant sheaf $\ul{\sfe}_{\on{Pic}}$. The corresponding equivalence \eqref{e:restr CFT} is then the following
statement:

\medskip

Pick a point $x\in X$. Write
$$\on{LocSys}^{\on{restr}}_{\BG_m}(X)_\sigma\simeq \on{pt}/\BG_m \times \on{Tot}(H^1(X,\ul\sfe_X))^\wedge_{\{0\}} \times \on{Tot}(\sfe[-1]),$$
see \secref{sss:Gm}, and 
$$\on{Pic}\simeq \BZ\times \on{Jac}(X)\times \on{pt}/\BG_m.$$

So
$$\QCoh(\on{LocSys}^{\on{restr}}_{\BG_m}(X)_\sigma) \simeq \Vect_\sfe^\BZ \otimes 
\left(\Sym(H^1(X,\ul\sfe_X)[-1])\mod\right) \otimes \left(\sfe[\xi]\mod\right), \quad \deg(\xi)=-1,$$
and 
$$\qLisse(\on{Pic})_\sigma \simeq  \Vect_\sfe^\BZ \otimes \Shv(\on{Jac}(X))_0 \otimes \Shv(\on{pt}/\BG_m),$$
where $\on{Jac}(X)$ is the Jacobian \emph{variety} of $X$, and 
$\Shv(\on{Jac}(X))_0\subset  \Shv(\on{Jac}(X))$ is the full subcategory generated by the constant sheaf.

\medskip

Now the result follows from the canonical identifications
$$\Shv(\on{pt}/\BG_m)\simeq \on{C}_\cdot(\BG_m)\mod\simeq  \sfe[\xi]\mod, \quad \deg(\xi)=-1,$$
and
$$\Shv(\on{Jac}(X))_0\simeq \CEnd(\ul{\sfe}_{\on{Jac}(X)})\mod =
\Sym(H^1(X,\ul\sfe_X)[-1])\mod,$$
see \secref{sss:ident Ab} for the latter isomorphism. 

\ssec{Comparison to other forms of the Geometric Langlands Conjecture}

\sssec{}

Let us specialize to the de Rham context. Recall that in this case, we have a version of the geometric
Langlands conjecture from \cite[Conjecture 11.2.2]{AG}, which predicts the existence of a canonical equivalence
\begin{equation} \label{e:geom Langlands dR}
\Dmod(\Bun_G) \simeq \IndCoh_{\on{Nilp}}(\on{LocSys}^\dr_\cG(X)),
\end{equation}
compatible with the actions of $\QCoh(\on{LocSys}^\dr_\cG(X))$ on both sides.

\medskip

We note that \conjref{c:restr GLC} (in the de Rham context) is a formal corollary of this statement. Namely, tensoring both sides of
\eqref{e:geom Langlands dR} with $\QCoh(\LocSys_\cG^{\on{restr}}(X))$ over $\QCoh(\LocSys^\dr_\cG(X))$,
we obtain: 
$$
\xy
(0,0)*+{\QCoh(\LocSys_\cG^{\on{restr}}(X))\underset{\QCoh(\LocSys^\dr_\cG(X))}\otimes \Dmod(\Bun_G)}="A";
(75,0)*+{\Shv_{\on{Nilp}}(\Bun_G)}="B";
(0,-20)*+{\QCoh(\LocSys_\cG^{\on{restr}}(X))\underset{\QCoh(\LocSys^\dr_\cG(X))}\otimes \IndCoh_{\on{Nilp}}(\on{LocSys}^\dr_\cG(X)) }="C";
(75,-20)*+{\IndCoh_{\on{Nilp}}(\on{LocSys}_\cG^{\on{restr}}(X))}="D";
{\ar@{->} "A";"B"};
{\ar@{->} "C";"D"};
{\ar@{->}^{\sim} "A";"C"};
{\ar@{->}^{\sim} "B";"D"};
\endxy
$$
%
where the top horizontal arrow comes from \propref{p:nilp via spectral}, and the bottom horizontal arrow is an equivalence by 
\secref{sss:formal compl sing supp} and \eqref{e:qcoh restr via amb}.


\sssec{}

Let us now specialize to the Betti context. In this case, we have a version of the geometric
Langlands conjecture, proposed in \cite[Conjecture 1.5]{BN}, which says that there is an equivalence 
\begin{equation} \label{e:geom Langlands Betti}
\Shv^{\on{all}}_{\on{Nilp}}(\Bun_G) \simeq \IndCoh_{\on{Nilp}}(\on{LocSys}_\cG(X)),
\end{equation}
compatible with the actions of $\QCoh(\on{LocSys}_\cG(X))$ on both sides.

\medskip

We note that \conjref{c:restr GLC} (in the Betti context) is a formal corollary of this statement. 
Namely, tensoring both sides of
\eqref{e:geom Langlands Betti} with $\QCoh(\LocSys_\cG^{\on{restr}}(X))$ over $\QCoh(\LocSys^{\on{Betti}}_\cG(X))$,
we obtain: 
$$
\xy
(0,0)*+{\QCoh(\LocSys_\cG^{\on{restr}}(X))\underset{\QCoh(\LocSys^{\on{Betti}}_\cG(X))}\otimes \Shv^{\on{all}}_{\on{Nilp}}(\Bun_G)}="A";
(74,0)*+{\Shv_{\on{Nilp}}(\Bun_G)}="B";
(0,-20)*+{\QCoh(\LocSys_\cG^{\on{restr}}(X))\underset{\QCoh(\LocSys^{\on{Betti}}_\cG(X))}\otimes \IndCoh_{\on{Nilp}}(\on{LocSys}^{\on{Betti}}_\cG(X)) }="C";
(74,-20)*+{\IndCoh_{\on{Nilp}}(\on{LocSys}_\cG^{\on{restr}}(X))}="D";
{\ar@{->} "A";"B"};
{\ar@{->} "C";"D"};
{\ar@{->}^{\sim} "A";"C"};
{\ar@{->}^{\sim} "B";"D"};
\endxy
$$
%
where the top horizontal arrow comes from \thmref{t:fin monod} and \secref{sss:shvs all fin mon}, and the bottom horizontal arrow is an equivalence by 
\secref{sss:formal compl sing supp} and \eqref{e:qcoh restr via amb}. 

\ssec{A converse implication}

Above we have seen that the full de Rham version of the Geometric Langlands Conjecture implies the restricted version. 
Here we will show that the converse implication also takes place, under a plausible
hypothesis about the de Rham version.

\sssec{}

We place ourselves into the de Rham context of the Geometric Langlands Conjecture. Let us assume the following:

\begin{hypoth} \label{h:Lde Rham GLC}
There exists a functor
$$\BL:\IndCoh_{\on{Nilp}}(\on{LocSys}^\dr_\cG(X))\to \Dmod(\Bun_G)$$ 
that preserves compactness and is compatible with the actions of $\QCoh(\LocSys^\dr_\cG(X))$ on both sides.
\end{hypoth}

This hypothesis would be a theorem if one accepted Quasi-Theorems 6.7.2 and 9.5.3 from  \cite{Ga7}. 

\sssec{}

We now claim:

\begin{obs} \label{o:de Rham GLC}
Assume that the functor $\BL$ from Hypothesis \ref{h:Lde Rham GLC} induces an equivalence
\begin{multline} \label{e:de Rham GLC}
\IndCoh_{\on{Nilp}}(\on{LocSys}_\cG^{\on{restr}}(X)) \simeq \\
\simeq \QCoh(\LocSys_\cG^{\on{restr}}(X))\underset{\QCoh(\LocSys^\dr_\cG(X))}\otimes \IndCoh_{\on{Nilp}}(\on{LocSys}^\dr_\cG(X)) 
\overset{\on{Id}\otimes \BL}\to \\
\to \QCoh(\LocSys_\cG^{\on{restr}}(X))\underset{\QCoh(\LocSys^\dr_\cG(X))}\otimes \Dmod(\Bun_G)  \simeq 
\Shv_{\on{Nilp}}(\Bun_G)
\end{multline}
(for the corresponding objects for all field extensions of $k$). Then the functor $\BL$ itself is an equivalence.
\end{obs} 

The rest of this subsection is devoted to the proof this Observation.

\sssec{}

Since the functor $\BL$ preserves compactness, it admits a continuous right adjoint. Since 
$\QCoh(\LocSys^\dr_\cG(X))$ is rigid, this right adjoint is compatible with the action of $\QCoh(\LocSys^\dr_\cG(X))$.

\medskip

Consider the adjunction maps 
\begin{equation} \label{e:adj Langlands}
\on{Id}\to \BL^R\circ \BL \text{ and } \BL\circ \BL^R\to \on{Id}.
\end{equation}
We wish to show that they are isomorphisms. 

\medskip

We have the following general assertion:

\begin{lem} \label{l:check on alg stack}
Let $\bC$ be a category acted on by $\QCoh(\CY)$, where $\CY$ is a quasi-compact eventually coconnective
algebraic stack almost of finite type with affine diagonal. Then an object $\bc\in \bC$ is zero if and only if for every geometric
point 
$$\bi:\Spec(k')\to \CY,$$
the image of $\bc$ under
$$\bC \simeq \QCoh(\CY)\underset{\QCoh(\CY)}\otimes \bC \overset{\bi^*\otimes \on{Id}}\longrightarrow 
\Vect_{k'}\underset{\QCoh(\CY)}\otimes \bC$$
vanishes.
\end{lem}

The proof of the lemma is given below. Let us accept it temporarily. 

\sssec{}

Applying the lemma, we obtain that in order to prove that the maps \eqref{e:adj Langlands} are isomorphisms, it suffices to show that the
functor

\smallskip

\begin{multline} \label{e:Langlands tensor up}
\Vect_{k'}\underset{\QCoh(\LocSys^\dr_\cG(X))}\otimes \IndCoh_{\on{Nilp}}(\on{LocSys}^\dr_\cG(X))\overset{\on{Id}\otimes \BL}\longrightarrow \\
\to \Vect_{k'}\underset{\QCoh(\LocSys^\dr_\cG(X))}\otimes \Dmod(\Bun_G)
\end{multline} 
is an equivalence for all geometric points
$$\bi:\Spec(k')\to \LocSys^\dr_\cG(X).$$

\medskip

Let $X',G'$ denote the base change of $X,G$ along $k\rightsquigarrow k'$.  Let $\Bun'_G$ denote the
corresponding algebraic stack over $k'$. We have
$$\LocSys^\dr_{\cG'}(X')\simeq \Spec(k')\underset{\Spec(k)}\times \LocSys^\dr_{\cG}(X) \text{ and }
\Bun'_G\simeq \Spec(k')\underset{\Spec(k)}\times \Bun_G,$$
and hence
$$\QCoh(\LocSys^\dr_{\cG'}(X'))\simeq \Vect_{k'}\underset{\Vect_k}\otimes \QCoh(\LocSys^\dr_\cG(X)),$$
$$\IndCoh_{\on{Nilp}}(\on{LocSys}^\dr_{\cG'}(X')) \simeq \Vect_{k'}\underset{\Vect_k}\otimes \IndCoh_{\on{Nilp}}(\on{LocSys}^\dr_\cG(X))$$
and
$$\Dmod(\Bun'_G)\simeq \Vect_{k'}\underset{\Vect_k}\otimes \Dmod(\Bun_G).$$

\medskip

Hence, we can rewrite the map in \eqref{e:Langlands tensor up} as
$$\Vect_{k'}\underset{\QCoh(\LocSys^\dr_{\cG'}(X'))}\otimes \IndCoh_{\on{Nilp}}(\on{LocSys}^\dr_{\cG'}(X'))\to
\Vect_{k'}\underset{\QCoh(\LocSys_{\cG'}(X'))}\otimes \Dmod(\Bun'_G).$$

Thus, we have reduced the verification of the isomorphism \eqref{e:Langlands tensor up} to the case when $k'=k$.

\sssec{}

Note now that (for $k'=k$), the map
$\bi:\Spec(k)\to \LocSys^\dr_\cG(X)$ factors as
$$\Spec(k)\to \LocSys^{\on{restr}}_\cG(X)\to \LocSys^\dr_\cG(X).$$

Hence, the map \eqref{e:Langlands tensor up} is obtained by
$$\Vect_{k}\underset{\QCoh(\LocSys^{\on{restr}}_\cG(X))}\otimes-$$
from \eqref{e:de Rham GLC}, and hence is an equivalence.

\qed[Observation \ref{o:de Rham GLC}]

\begin{rem}
Note that whereas 
$$\LocSys^\dr_{\cG'}(X')\simeq \Spec(k')\underset{\Spec(k)}\times \LocSys^\dr_{\cG}(X),$$
the same \emph{no longer} holds for $\LocSys^{\on{restr}}_{\cG}(X)$. Similarly, while 
$$\Dmod(\Bun'_G)\simeq \Vect_{k'}\underset{\Vect_k}\otimes \Dmod(\Bun_G),$$
the same is no longer true for $\Shv_\Nilp(\Bun_G)$.
\end{rem}

\begin{rem}
A counterpart of Observation \ref{o:de Rham GLC} would apply also in the Betti version, using 
the identifications \eqref{e:BC all} and \eqref{e:BC Betti} (cf. \secref{ss:proof projector Betti pt}). 
We are just less confident of the status 
of the analog of Hypothesis \ref{h:Lde Rham GLC} in this case. 

\end{rem} 

\sssec{Proof of \lemref{l:check on alg stack}}

First, by \cite[Theorem 2.2.6]{Ga2}, $\CY$ is 1-affine\footnote{For our applications, we only need in the case
when $\CY$ is of the form $S/\sH$, where $S$ is an affine scheme and $\sH$ is an algebraic group, in which
case the assertion of \cite[Theorem 2.2.6]{Ga2} easily follows from the case of $\on{pt}/\sH$.},
hence we can replace $\CY$ by an affine scheme $S=\Spec(A)$.

\medskip

Since $\CY$ was assumed eventually coconnective, $S$ has the same property. Hence, $\QCoh(S)$
is generated under colimits by objects of the form $\wt\bi_*(\wt{k}')$, for 
$$\wt\bi:\Spec(\wt{k}')\to S,$$
$\wt{k}'$ are residue fields of scheme-theoretic points of $S$, see \lemref{l:descr ker}.
In particular $\CO_S$ can be expressed as a 
colimit of objects of this form.

\medskip

Hence, $\bc$ can be expressed as a colimit of objects of the form
$$\wt\bi_*(\wt{k}') \otimes \bc \simeq
(\wt\bi_*\otimes \on{Id})\circ (\wt\bi^*\otimes \on{Id})(\bc).$$

Hence, if all $(\wt\bi^*\otimes \on{Id})(\bc)$ vanish, then $\bc$ vanishes.

\medskip

Let $k'$ be the algebraic closure of $\wt{k}'$. It is easy to see that for any $\Vect_{\wt{k}'}$-linear category $\wt\bC$, 
the pullback functor
$$\wt\bC \to \Vect_{k'}\underset{\Vect_{\wt{k}'}}\otimes \wt\bC$$
is conservative. 

\medskip

Hence, for $\bi$ equal to
$$\Spec(k')\to \Spec(\wt{k}')\overset{\wt\bi}\to S,$$
if $(\bi^*\otimes \on{Id})(\bc)$ vanishes, then so does $(\wt\bi^*\otimes \on{Id})(\bc)$. 

\qed[\lemref{l:check on alg stack}]

\section{The trace conjecture} \label{s:Trace}

Throughout this section we will be working with schemes/algebraic stacks of finite type over $\ol\BF_q$, that
are defined over $\BF_q$, so that they carry the geometric Frobenius endomorphism. 

\medskip

Our sheaf-theoretic context will (by necessity) be that of $\ell$-adic sheaves, so $\sfe=\ol\BQ_\ell$. 

\medskip

This section contains what is the main point of this paper. We propose a conjecture that expresses the space
of automorphic functions as the categorical trace of Frobenius acting on the category $\Shv_{\on{Nilp}}(\Bun_G)$.

\ssec{The categorical trace of Frobenius}

\sssec{} \label{sss:Frob rad}

Let $\CY$ be an algebraic stack. While discussing general algebraic stacks in this and the next subsection,
we will assume that $\CY$ is locally a quotient (of a scheme by an algebraic group); this is an assumption under
which the results of \cite{GaVa} are established.

\medskip

We consider the Frobenius endomorphism $\Frob_\CY$ of $\CY$. One word of warning is that when $\CY$
is an Artin stack, $\Frob_\CY$ is not necessarily schematic. However, it is surjective and radicial, in the sense that it becomes
an isomorphism after we apply sheafificaton in the topology generated by surjective radicial maps. 

\medskip

Hence, the action of $\Frob_\CY$ on the category $\Shv(\CY)$ has properties of a surjective radicial map.
In particular, the functor 
$$(\on{Frob}_\CY)^*:\Shv(\CY)\to \Shv(\CY)$$
is an equivalence (and hence its right adjoint $(\on{Frob}_\CY)_*$ is also an equivalence). 

\medskip

Furthermore, the natural transformation
\begin{equation} \label{e:Frob rad}
(\on{Frob}_\CY)_!\to (\on{Frob}_\CY)_*
\end{equation}
(which is a priori well-defined due to the fact that $\on{Frob}_\CY$ is separated) is an isomorphism,
and both functors are equivalences. 

\medskip

From here it follows that left adjoint of $(\on{Frob}_\CY)_*$ is isomorophic to the right adjoint of $(\on{Frob}_\CY)_!$, i.e.,
$$(\on{Frob}_\CY)^* \simeq (\on{Frob}_\CY)^!.$$

\sssec{}

Assume first that $\CY$ is quasi-compact. 

\medskip

The category $\Shv(\CY)$ is compactly
generated, and hence dualizable. Hence, we can consider the categorical trace of $(\on{Frob}_\CY)_*$:
$$\Tr((\on{Frob}_\CY)_*,\Shv(\CY))\in \Vect_\sfe.$$

\medskip

To $\CF\in \Shv(\CY)^c$ equipped with a map 
\begin{equation} \label{e:lax Weil}
\CF\overset{\alpha}\to (\on{Frob}_\CY)_*(\CF),
\end{equation}
we can attach its class
$$\on{cl}(\CF,\alpha)\in \Tr((\on{Frob}_\CY)_*,\Shv(\CY)),$$
see \cite[Sect. 3.4.3]{GKRV}. 

\medskip

We will refer to the data of $\alpha$ as a \emph{lax Weil structure} on $\CF$, and to the pair
$(\CF,\alpha)$ as a \emph{lax Weil sheaf} on $\CY$.  

\sssec{}

We claim that there is a canonically defined map, called the \emph{Local Term},
\begin{equation} \label{e:Tr to funct}
\on{LT}:\Tr((\on{Frob}_\CY)_*,\Shv(\CY))\to \on{Funct}(\CY(\BF_q)),
\end{equation}
where $\on{Funct}(-)$ stands for the (classical) vector space of $\sfe$-valued functions on the set of isomorphism
classes of a given groupoid.

\medskip

In fact, there are two such maps, denoted $\on{LT}^{\on{naive}}$ and $\on{LT}^{\on{true}}$.

\sssec{} \label{sss:sheaf-funct}

The map $\on{LT}^{\on{naive}}$ is designed so that for a lax Weil sheaf $(\CF,\alpha)$, we will have 
$$\on{LT}^{\on{naive}}(\on{cl}(\CF,\alpha))=\on{funct}(\CF,\alpha),$$
where $\on{funct}(\CF,\alpha)$ is the usual function on $\CY(\BF_q)$ attached to $(\CF,\alpha)$
obtained by taking traces of the Frobenius on $\BF_q$-points:

\medskip

By adjunction, the datum of $\alpha$ is equivalent to the datum of a map
$$\alpha^L:(\on{Frob}_\CY)^*(\CF)\to \CF.$$

Now the value of $\on{funct}(\CF,\alpha)$ on a given $y\in \CY(\BF_q)$, 
$$\on{pt} \overset{\bi_y}\to \CY$$
equals the trace of 
$$\bi_y^*(\CF) \overset{y\text{ is Frob-invariant}}\longrightarrow (\on{Frob}_\CY\circ \bi_y)^*(\CF)\simeq
\bi_y^*\circ \on{Frob}_\CY^*(\CF) \overset{\alpha^L}\to \bi_y^*(\CF).$$

\sssec{}

The actual definition of $\on{LT}^{\on{naive}}$ proceeds as follows. Every $y$ as above defines a functor
$$\bi_y^*:\Shv(\CY)\to \Vect_\sfe,$$
which admits a continuous right adjoint (namely, $(\bi_y)_*$), and is
equipped with a morphism (in fact, an isomorphism) 
\begin{equation} \label{e:Frob nat trans ev}
\bi_y^*\circ (\on{Frob}_\CY)_*\to \bi_y^*.
\end{equation} 

\medskip

Hence, by \cite[Sect. 3.4.1]{GKRV}, it defines a map
$$\Tr((\on{Frob}_\CY)_*,\Shv(\CY))\to \Tr(\on{Id},\Vect_\sfe)\simeq \sfe.$$

This map is, by definition, the composition of $\on{LT}^{\on{naive}}$ with the evaluation map
$$\on{Funct}(\CY(\BF_q)) \overset{\on{ev}_y}\to \sfe.$$

\medskip

The map $\on{LT}^{\on{naive}}$ has the following features.

\sssec{} \label{sss:Frob ten prod}

For a lax Weil sheaf $(\CF_0,\alpha_0)$ on $\Shv(\CY)$, consider the functor
$$\Shv(\CY)\to \Shv(\CY), \quad \CF\mapsto \CF_0\overset{*}\otimes \CF.$$

This functor is endowed with a natural transformation
\begin{equation} \label{e:Frob nat trans ten}
\CF_0\overset{*}\otimes (\on{Frob}_\CY)_*(\CF)\to (\on{Frob}_\CY)_*(\CF_0\overset{*}\otimes \CF),
\end{equation} 
and it admits a (continuous) right adjoint, given by 
$$\CF\mapsto \BD(\CF_0)\sotimes \CF.$$

Hence, by \cite[Sect. 3.4.1]{GKRV}, it defines a map
\begin{equation} \label{e:Frob tensor product}
\Tr((\on{Frob}_\CY)_*,\Shv(\CY)) \to \Tr((\on{Frob}_\CY)_*,\Shv(\CY)).
\end{equation}

We claim that there is a commutative diagram
$$
\CD
\Tr((\on{Frob}_\CY)_*,\Shv(\CY)) @>{\on{LT}^{\on{naive}}}>> \on{Funct}(\CY(\BF_q))  \\
@V{\text{\eqref{e:Frob tensor product}}}VV @VV{\on{funct}(\CF_0,\alpha_0)\cdot -}V \\
\Tr((\on{Frob}_\CY)_*,\Shv(\CY)) @>{\on{LT}^{\on{naive}}}>> \on{Funct}(\CY(\BF_q)). 
\endCD
$$

Indeed, this follows from the fact that for a given $y\in \CY(\BF_q)$, we have a commutative diagram of functors
$$
\CD
\Shv(\CY) @>{\bi_y^*}>> \Vect_\sfe \\
@V{\CF_0\overset{*}\otimes}VV @VV{\bi_y^*(\CF_0)\otimes -}V \\
\Shv(\CY) @>{\bi_y^*}>> \Vect_\sfe
\endCD
$$
compatible with the natural transformations \eqref{e:Frob nat trans ev} and \eqref{e:Frob nat trans ten} via
the endomorphism on $\bi_y^*(\CF_0)$ given by $\alpha^L_0$. This implies that the resulting map
$$\sfe\simeq \Tr(\on{Id},\Vect_\sfe) \overset{(\bi_y^*(\CF_0),\alpha^L_0)}\longrightarrow  \Tr(\on{Id},\Vect_\sfe) \simeq \sfe$$
is given by multiplication by
$$\on{Tr}(\alpha^L_0,\bi_y^*(\CF_0))=\on{funct}(\CF_0,\alpha_0)(y),$$
as desired. 

\sssec{} \label{sss:Frob pull back}

Let $f:\CY_1\to \CY_2$ be a map. Consider the functor
$$f^*:\Shv(\CY_2)\to \Shv(\CY_1).$$

This functor is endowed with a natural transformation (in fact, an isomorphism)
\begin{equation} \label{e:Frob nat trans pullback}
f^*\circ (\on{Frob}_{\CY_2})_* \to (\on{Frob}_{\CY_1})_*\circ f^*.
\end{equation} 

The right adjoint of $f^*$ is the usual direct image functor $f_*$. However, for
morphisms between stacks, the functor $f_*$ is not necessarily continuous. 
Therefore, in what follows we will assume that $f$ is \emph{safe} in the sense
of \cite[Definition 10.2.2]{DrGa1}. Concretely, this condition means that for any
geometric point of any geometric fiber of $f$, the neutral connected component
of its group of automorphisms is unipotent. In particular, any schematic map 
between algebraic stacks is safe. One proves that if $f$ is safe, then the
functor $f_*$ is continuous, see \cite[Theorem 10.2.4]{DrGa1}.

\medskip

Assuming that $f_*$ is continuous, by \cite[Sect. 3.4.1]{GKRV}, the functor $f^*$ defines a map
\begin{equation} \label{e:Frob pull back}
\Tr((\on{Frob}_{\CY_2})_*,\Shv(\CY_2)) \to \Tr((\on{Frob}_{\CY_1})_*,\Shv(\CY_1)).
\end{equation}

We claim that there is a commutative diagram
$$
\CD
\Tr((\on{Frob}_{\CY_2})_*,\Shv(\CY_2)) @>{\on{LT}^{\on{naive}}}>> \on{Funct}(\CY_2(\BF_q))  \\
@V{\text{\eqref{e:Frob pull back}}}VV @VV{\text{pull back}}V \\
\Tr((\on{Frob}_{\CY_1})_*,\Shv(\CY_1)) @>{\on{LT}^{\on{naive}}}>> \on{Funct}(\CY_1(\BF_q)),
\endCD
$$
where the right vertical arrow is given by pullback of functions along the induced map
$$\CY_1(\BF_q)\to \CY_2(\BF_q).$$

This follows just from the fact that the *-pullback functor is compatible with compositions.  

\sssec{}

Finally, let $f:\CY_1\to \CY_2$ be as above. Consider the functor
$$f_!:\Shv(\CY_1)\to \Shv(\CY_2).$$

This functor is endowed with a natural transformation (in fact, an isomorphism)
\begin{equation} \label{e:Frob nat trans pushforward}
f_! \circ (\on{Frob}_{\CY_1})_* \to (\on{Frob}_{\CY_2})_*\circ f_!,
\end{equation} 
(coming from \eqref{e:Frob rad}), and it admits a (continuous) right adjoint, given by $f^!$. 

\medskip

Hence, by \cite[Sect. 3.4.1]{GKRV}, it defines a map
\begin{equation} \label{e:Frob pushforward}
\Tr((\on{Frob}_{\CY_1})_*,\Shv(\CY_1)) \to \Tr((\on{Frob}_{\CY_2})_*,\Shv(\CY_2)).
\end{equation} 

\begin{thm} \label{t:GLTF}
We have a commutative diagram
$$
\CD
\Tr((\on{Frob}_{\CY_1})_*,\Shv(\CY_1)) @>{\on{LT}^{\on{naive}}}>> \on{Funct}(\CY_1(\BF_q))  \\
@V{\text{\eqref{e:Frob pushforward}}}VV @VV{\text{\rm{push forward}}}V \\
\Tr((\on{Frob}_{\CY_2})_*,\Shv(\CY_2)) @>{\on{LT}^{\on{naive}}}>> \on{Funct}(\CY_2(\BF_q)),
\endCD
$$
where the right vertical arrow is given by (weighted)\footnote{We weigh each point by $\frac{1}{|\text{order of its group of automorphisms}|}$.} 
summation along the fiber of the induced map
$$\CY_1(\BF_q)\to \CY_2(\BF_q).$$
\end{thm}

This theorem is a version of the Grothendieck-Lefschetz trace formula. The proof is given in
\cite{GaVa}. 

\begin{rem}
It is easy to prove \thmref{t:GLTF} when $f$ is a locally closed embedding.
And this is the only case we will need in order to formulate \conjref{c:Trace conj}. 
\end{rem}

\sssec{}

Let now $\CY$ be an algebraic stack that is not necessarily quasi-compact. We write
\begin{equation} \label{e:Y as u U}
\CY:=\underset{\CU}\cup\, \CU,
\end{equation} 
where $\CU\overset{j}\hookrightarrow \CY$ is a filtered collection of 
quasi-compact open prestacks, so that
$$\Shv(\CY)\simeq \underset{\CU}{\on{lim}}\,\Shv(\CU),$$
with respect to the restriction maps and also
$$\Shv(\CY)\simeq \underset{\CU}{\on{colim}}\,\Shv(\CU),$$
with respect to !-pushforwards, see \cite[Proposition 1.7.5]{DrGa2}.

\medskip

We claim that the functors $j_!:\Shv(\CU)\to \Shv(\CY)$ induce an isomorphism
\begin{equation} \label{e:trace Frob as colim}
\underset{\CU}{\on{colim}}\, \Tr((\on{Frob}_\CU)_*,\Shv(\CU))\to \Tr((\on{Frob}_\CY)_*,\Shv(\CY)).
\end{equation} 

Indeed, we have 
$$\on{Id}_{\Shv(\CY)}\simeq \underset{\CU}{\on{colim}}\, j_!\circ j^*,$$
and hence
\begin{multline*} 
\Tr((\on{Frob}_\CY)_*,\Shv(\CY))\simeq \underset{\CU}{\on{colim}}\, \Tr((\on{Frob}_\CY)_*\circ j_!\circ j^*,\Shv(\CY))\simeq \\
\simeq \underset{\CU}{\on{colim}}\, \Tr(j_!\circ (\on{Frob}_\CU)_*\circ j^*,\Shv(\CY))\overset{\text{cyclicity of trace}}\simeq 
\underset{\CU}{\on{colim}}\, \Tr((\on{Frob}_\CU)_*\circ j^*\circ j_!,\Shv(\CU))\simeq \\
\simeq  \underset{\CU}{\on{colim}}\, \Tr((\on{Frob}_\CU)_*,\Shv(\CU)),
\end{multline*} 
as desired.

\medskip

From here, using \thmref{t:GLTF} for open embeddings, we obtain that the maps $\on{LT}^{\on{naive}}$ for $\CU$ 
give rise to a map
\begin{equation} \label{e:map to funct non-qc}
\on{LT}^{\on{naive}}: \Tr((\on{Frob}_\CY)_*,\Shv(\CY))\to \on{Funct}_c(\CY(\BF_q)),
\end{equation} 
where $\on{Funct}_c(-)$ stands for ``functions with finite support", so
$$\on{Funct}_c(\CY(\BF_q))\simeq \underset{\CU}{\on{colim}}\, \on{Funct}_c(\CU(\BF_q)).$$

\ssec{The true local term}

We now proceed to the definition of the map
\begin{equation} \label{e:Tr to funct true}
\on{LT}^{\on{true}}:\Tr((\on{Frob}_\CY)_*,\Shv(\CY))\to \on{Funct}(\CY(\BF_q)).
\end{equation}

\medskip

As in the previous subsection, on the first pass we will assume that $\CY$ is 
quasi-compact\footnote{This is essential for the construction because Verdier duality a priori
works only for quasi-compact stacks.}. 
We will also assume that $\CY$ is Verdier-compatible, see \secref{sss:duality adapted} for what this means. 
(According to \conjref{c:duality preserve compact}, all quasi-compact algebraic stacks with an affine
diagonal have this property; in \thmref{t:global quotient} it is shown that algebraic stacks that can 
locally be written as quotients are such.) 

\sssec{}

We recall that the algebraic stack $\CY^{\on{Frob}}$ is \emph{discrete}, i.e., has the form
$$\sqcup\, (\on{pt}/\Gamma), \quad \Gamma\in \text{Finite Groups},$$
so we can identify
$$\on{Funct}(\CY(\BF_q))\simeq \on{C}^\cdot(\CY^{\on{Frob}},\omega_{\CY^{\on{Frob}}}).$$

Let $\bi_\CY$ denote the forgetful map
$$\CY^{\on{Frob}}\to \CY.$$

Let us rewrite
$$\on{C}^\cdot(\CY^{\on{Frob}},\omega_{\CY^{\on{Frob}}})\simeq
\on{C}^\cdot(\CY,(\bi_\CY)_*(\omega_{\CY^{\on{Frob}}})).$$

Using base change along
$$
\CD
\CY^{\on{Frob}} @>{\bi_\CY}>> \CY \\
@V{\bi_\CY}VV @VV{(\on{Frob}_\CY,\on{id}_\CY)}V \\
\CY @>{\Delta_\CY}>> \CY\times \CY,
\endCD
$$
we can rewrite
$$(\bi_\CY)_*(\omega_{\CY^{\on{Frob}}})\simeq \Delta_\CY^! \circ (\on{Frob}_\CY\times \on{id}_\CY)_* \circ (\Delta_\CY)_*(\omega_\CY).$$

To summarize, we have
$$\on{Funct}(\CY(\BF_q)) \simeq 
\on{C}^\cdot(\CY,\Delta_\CY^! \circ (\on{Frob}_\CY\times \on{id}_\CY)_* \circ (\Delta_\CY)_*(\omega_\CY)).$$

\sssec{}

In order to compute $\Tr((\on{Frob}_\CY)_*,\Shv(\CY))$, we identify $\Shv(\CY)$ with its own dual, see 
\secref{sss:duality on stack}. We recall that the corresponding pairing 
$$\Shv(\CY)\otimes \Shv(\CY)\to \Vect_\sfe$$
is given by
$$\CF_1,\CF_2\mapsto \on{C}^\cdot_\blacktriangle(\CY,\CF_1\sotimes \CF_2)\simeq 
\on{C}^\cdot_\blacktriangle(\CY,\Delta_\CY^!(\CF_1\boxtimes \CF_2)),$$
where the notation $\on{C}^\cdot_\blacktriangle$ is as in \secref{sss:black triangle}.

\medskip 

Let $\on{u}_{\Shv(\CY)}\in \Shv(\CY)\otimes \Shv(\CY)$ be the unit of the self-duality on $\Shv(\CY)$. 

\medskip

We obtain that
$\Tr((\on{Frob}_\CY)_*,\Shv(\CY))$ is given by
$$\on{C}^\cdot_\blacktriangle(\CY,\Delta_\CY^! \circ \boxtimes \circ ((\on{Frob}_\CY)_*\otimes \on{Id}_{\Shv(\CY)})(\on{u}_{\Shv(\CY)})),$$
where $\boxtimes$ denotes the external tensor product functor
$$\Shv(\CY)\otimes \Shv(\CY)\to \Shv(\CY\times \CY).$$

\sssec{}
 
We note that
$$(\on{Frob}_\CY\times \on{id}_\CY)_* \circ \boxtimes \simeq \boxtimes \circ ((\on{Frob}_\CY)_*\otimes \on{Id}_{\Shv(\CY)}).$$

\medskip

Hence, in order to construct the map \eqref{e:Tr to funct true}, it suffices to construct a map
\begin{equation} \label{e:map from unit to diag}
\boxtimes(\on{u}_{\Shv(\CY)}) \to (\Delta_\CY)_*(\omega_\CY),
\end{equation} 
and a map
\begin{multline} \label{e:black triang Frob}
\on{C}^\cdot_\blacktriangle(\CY,\Delta_\CY^! \circ (\on{Frob}_\CY\times \on{id}_\CY)_* \circ (\Delta_\CY)_*(\omega_\CY))\to
\on{C}^\cdot(\CY,\Delta_\CY^! \circ (\on{Frob}_\CY\times \on{id}_\CY)_* \circ (\Delta_\CY)_*(\omega_\CY))\simeq \\
\simeq \on{Funct}(\CY(\BF_q)).
\end{multline}

The first arrow in \eqref{e:black triang Frob} is the map \eqref{e:from triangle}. In our case, it is in fact an isomorphism, 
because the object
$$ \Delta_\CY^! \circ (\on{Frob}_\CY\times \on{id}_\CY)_* \circ (\Delta_\CY)_*(\omega_\CY)
\simeq (\bi_\CY)_*(\omega_{\CY^{\on{Frob}}})\in \Shv(\CY)$$
is compact, since $\CY$ is Verdier-compatible (by \thmref{t:global quotient}). 

\medskip
 
We proceed to the construction of the map \eqref{e:map from unit to diag}, cf. Remark \ref{r:quasi-diag}. 

%
%
%
%
%
%

\sssec{} \label{sss:quasi-diag}

Let $\boxtimes^R$ denote the right adjoint of the functor
$$\boxtimes:\Shv(\CY)\otimes \Shv(\CY)\to \Shv(\CY\times \CY).$$

We claim that we have a canonical isomorphism
\begin{equation} \label{e:U Shv}
\on{u}_{\Shv(\CY)}\simeq \boxtimes^R((\Delta_\CY)_*(\omega_\CY)),
\end{equation}
which would then give rise to the desired map \eqref{e:map from unit to diag} by adjunction.

\medskip

To establish \eqref{e:U Shv} we note that for $\CF_1,\CF_2\in \Shv(\CY)^c$, we have
$$\CHom_{\Shv(\CY)\otimes \Shv(\CY)}(\CF_1\otimes \CF_2,\on{u}_{\Shv(\CY)})\simeq
\CHom_{\Shv(\CY)}(\CF_1,\BD(\CF_2)),$$
by the definition of the self-duality on $\Shv(\CY)$ (here $\BD$ is Verdier duality), while
\begin{multline*}
\CHom_{\Shv(\CY)\otimes \Shv(\CY)}(\CF_1\otimes \CF_2,\boxtimes^R((\Delta_\CY)_*(\omega_\CY)))\simeq
\CHom_{\Shv(\CY\times \CY)}(\CF_1\boxtimes \CF_2,(\Delta_\CY)_*(\omega_\CY))\simeq \\
\simeq \CHom_{\Shv(\CY)}(\CF_1,\BD(\CF_2))
\end{multline*} 
as well. 

\sssec{Example} \label{sss:true term Weil}

Let $(\CF,\alpha)$ be a lax Weil sheaf on $\CY$. Unwinding the construction, we obtain that the image of
$$\on{cl}(\CF,\alpha)\in \Tr((\on{Frob}_\CY)_*,\Shv(\CY))$$
along the map $\on{LT}^{\on{true}}$, thought of as an element of 
\begin{multline*}
\on{Funct}(\CY(\BF_q)) \simeq 
\on{C}^\cdot(\CY,\Delta_\CY^! \circ (\on{Frob}_\CY\times \on{id}_\CY)_* \circ (\Delta_\CY)_*(\omega_\CY))\simeq \\
\simeq \CHom_{\Shv(\CY)}(\ul\sfe_{\CY},\Delta_\CY^! \circ (\on{Frob}_\CY\times \on{id}_\CY)_* \circ (\Delta_\CY)_*(\omega_\CY)),
\end{multline*} 
equals
\begin{multline*}
\ul\sfe_{\CY}\to \CF\sotimes \BD(\CF) \simeq \Delta_\CY^!(\CF\boxtimes \BD(\CF)) \overset{\alpha\boxtimes \on{id}}\longrightarrow 
\Delta_\CY^!((\Frob_\CY)_*(\CF)\boxtimes \BD(\CF)) 
\simeq \Delta_\CY^! \circ (\Frob_\CY \times \on{id}_\CY)_*(\CF\boxtimes \BD(\CF)) \to \\
\to \Delta_\CY^! \circ (\Frob_\CY \times \on{id}_\CY)_*\circ (\Delta_\CY)_*\circ \Delta_\CY^*(\CF\boxtimes \BD(\CF)) \simeq
\Delta_\CY^! \circ (\Frob_\CY \times \on{id}_\CY)_*\circ (\Delta_\CY)_*(\CF\overset{*}\otimes \BD(\CF)) \to \\
\to \Delta_\CY^! \circ (\on{Frob}_\CY\times \on{id}_\CY)_* \circ (\Delta_\CY)_*(\omega_\CY). 
\end{multline*}

\begin{rem} \label{r:why naive trace}
The map \eqref{e:map from unit to diag} constructed above is, in general, \emph{not} an isomorphism. In fact it is an isomorphism
\emph{if and only if} the functor
$$\Shv(\CY)\otimes \Shv(\CY) \overset{\boxtimes}\to \Shv(\CY\times \CY)$$
is an equivalence, see \secref{sss:when tensor product}. 

\medskip

The fact that the map \eqref{e:map from unit to diag} is not in general an isomorphism prevents the map 
\eqref{e:Tr to funct true} from being an isomorphism.

\medskip

However (as was remarked in \secref{sss:when tensor product}), we obtain that \eqref{e:map from unit to diag} is an isomorphism for algebraic
stacks that have finitely many isomorphism classes of $\ol\BF_q$-points, e.g., for $N\backslash G/B$,
or a quasi-compact substack of $\Bun_G$ for a curve $X$ of genus $0$. Hence, \eqref{e:Tr to funct true} is an isomorphism in these cases
as well. 
\end{rem} 

\sssec{}

We claim:

\begin{thm} \label{t:local terms}
The maps 
$$\on{LT}^{\on{naive}}:\Tr((\on{Frob}_\CY)_*,\Shv(\CY))\to \on{Funct}(\CY(\BF_q))$$
and
$$\on{LT}^{\on{true}}:\Tr((\on{Frob}_\CY)_*,\Shv(\CY))\to \on{Funct}(\CY(\BF_q))$$
are canonically homotopic.
\end{thm}

The proof is given in \cite{GaVa}. From now on, we will just use the symbol 
$$\on{LT}:\Tr((\on{Frob}_\CY)_*,\Shv(\CY))\to \on{Funct}(\CY(\BF_q))$$
for the local term map. 

\begin{rem}
\thmref{t:local terms} implies that for a lax Weil sheaf $(\CF,\alpha)$, the images of
$\on{cl}(\CF,\alpha)\in \Tr((\on{Frob}_\CY)_*,\Shv(\CY))$ in $\on{Funct}(\CY(\BF_q))$
under the above two maps coincide.

\medskip

Interpreting the image of $\on{cl}(\CF,\alpha)$ along $\on{LT}^{\on{true}}$ as in \secref{sss:true term Weil}
and along $\on{LT}^{\on{naive}}$ as in \secref{sss:sheaf-funct}, 
the latter assertion becomes equivalent to one in \cite[Theorem 2.1.3]{Var2} (when $\CY$ is a scheme). 
The proof of \thmref{t:local terms} is an elaboration of the ideas from {\it loc. cit.}

\end{rem}

\begin{cor} \label{c:local terms} Let $f:\CY_1\to \CY_2$ be a map
between stacks. 

\smallskip

\noindent{\em(a)} 
If $f$ is safe, then the diagram 
$$
\CD
\Tr((\on{Frob}_{\CY_2})_*,\Shv(\CY_2)) @>{\on{LT}^{\on{true}}}>> \on{Funct}(\CY_2(\BF_q))  \\
@V{\text{\eqref{e:Frob pull back}}}VV @VV{\text{\rm{pull back}}}V \\
\Tr((\on{Frob}_{\CY_1})_*,\Shv(\CY_1)) @>{\on{LT}^{\on{true}}}>> \on{Funct}(\CY_1(\BF_q)),
\endCD
$$
commutes. 

\smallskip

\noindent{\em(b)} 
The diagram 
$$
\CD
\Tr((\on{Frob}_{\CY_1})_*,\Shv(\CY_1)) @>{\on{LT}^{\on{naive}}}>> \on{Funct}(\CY_1(\BF_q))  \\
@V{\text{\eqref{e:Frob pushforward}}}VV @VV{\text{\rm{push forward}}}V \\
\Tr((\on{Frob}_{\CY_2})_*,\Shv(\CY_2)) @>{\on{LT}^{\on{naive}}}>> \on{Funct}(\CY_2(\BF_q)),
\endCD
$$
commutes.

\end{cor}

\sssec{}

We now let $\CY$ be an arbitrary algebraic stack locally of finite type (i.e., not necessarily quasi-compact),
written as \eqref{e:Y as u U}

\medskip

Using \eqref{e:trace Frob as colim} and \corref{t:local terms}(a), we combine the maps $\on{LT}^{\on{true}}$ 
for the quasi-compact open substacks $\CU\subset \CY$ to a map
$$\on{LT}^{\on{true}}:\Tr((\on{Frob}_\CY)_*,\Shv(\CY))\to \on{Funct}_c(\CY(\BF_q)).$$

Furthermore, \thmref{t:local terms} implies that the above map equals the map
$$\on{LT}^{\on{naive}}:\Tr((\on{Frob}_\CY)_*,\Shv(\CY))\to \on{Funct}_c(\CY(\BF_q))$$
of \eqref{e:map to funct non-qc}. 

\ssec{From geometric to classical: the Trace Conjecture}

\sssec{}

We start with the following observation. Let $\CY$ be an algebraic stack over $\ol\BF_q$,
but defined over $\BF_q$. Consider the diagram
$$
\CD
\CY @>{\Frob^{\on{arithm}}_\CY}>> \CY @>{\Frob_\CY}>> \CY \\
@VVV @VVV \\
\Spec(\ol\BF_q) @>{\Frob^{\on{arithm}}}>> \Spec(\ol\BF_q),
\endCD
$$
where:

\begin{itemize}

\item The bottom horizontal arrow is the Frobenius automorphism of $\Spec(\ol\BF_q)$;

\item The square is Cartesian;

\item The composite top horizontal arrow is the absolute Frobenius on $\CY$.

\end{itemize}

\medskip

For a Zariski-closed conical $\CN\subset T^*(\CY)$, let $\CN'\subset T^*(\CY)$ denote the base-change of $\CN$ along 
$\Frob^{\on{arithm}}$. 

\medskip

We claim: 

\begin{lem} \label{l:Frob pres sing supp}
The functor $(\Frob_\CY)_*:\Shv(\CY)\to \Shv(\CY)$ 
sends $\Shv_\CN(\CY)\to \Shv_{\CN'}(\CY)$.
\end{lem}

\begin{proof}

The functor $(\Frob_\CY)_*$ is the inverse of the pullback functor $(\Frob_\CY)^*$, and the latter is the inverse 
of $(\Frob^{\on{arithm}}_\CY)^*$, since pullback by the absolute Frobenius acts as identity. 

\medskip

Hence, the functor $(\Frob_\CY)_*$ identifies with $(\Frob^{\on{arithm}}_\CY)^*$. Therefore,
it is sufficient to show that $(\Frob^{\on{arithm}}_\CY)^*$ sends $\Shv_\CN(\CY)\to \Shv_{\CN'}(\CY)$. 

\medskip

Now, for any map of fields $k\to k'$, the pullback functor along
$$\CY':=\Spec(k')\underset{\Spec(k)}\times \CY\to \CY$$
sends $\Shv_\CN(\CY)\subset \Shv(\CY)$ to $\Shv_{\CN'}(\CY')\subset \Shv(\CY')$.

\end{proof} 

\sssec{} We now come to one of the central ideas of this paper. 

\medskip

Recall (see \thmref{t:Nilp comp gen abs}) that the category $\Shv_{\on{Nilp}}(\Bun_G)$ is compactly generated, and in particular, dualizable.
By \lemref{l:Frob pres sing supp}, the action of $(\on{Frob}_{\Bun_G})_*$ on $\Shv(\Bun_G)$ preserves the subcategory
\begin{equation} \label{e:embed Nilp again}
\Shv_{\on{Nilp}}(\Bun_G)\overset{\iota}\hookrightarrow \Shv(\Bun_G);
\end{equation}
in particular, $(\on{Frob}_{\Bun_G})_*$ restricts to an endofunctor of $\Shv_{\on{Nilp}}(\Bun_G)$. 

\medskip

Hence, we can form the object
$$\Tr((\on{Frob}_{\Bun_G})_*,\Shv_{\on{Nilp}}(\Bun_G))\in \Vect_\sfe.$$

\sssec{}

Recall also the subcategory
\begin{equation} \label{e:embed Nilp acess again}
\Shv_{\on{Nilp}}(\Bun_G)^{\on{access}}\hookrightarrow \Shv_{\on{Nilp}}(\Bun_G),
\end{equation}
see \secref{sss:access BunG}\footnote{Note, however, that according to (the very plausible) \conjref{c:Nilp comp gen}, the inclusion 
\eqref{e:embed Nilp acess again} is an equivalence.}.

\medskip

By the combination of \thmref{t:preserve Nilp Sing Supp prel} and \corref{c:N preserved abs},
it is generated by objects compact in $\Shv(\Bun_G)$; in particular it is compactly generated. 
The endofunctor $(\on{Frob}_{\Bun_G})_*$ of $\Shv_{\on{Nilp}}(\Bun_G)$ preserves the subcategory
$\Shv_{\on{Nilp}}(\Bun_G)^{\on{access}}$ (indeed, its compact generators are those objects of
$\Shv_{\on{Nilp}}(\Bun_G)$ that are compact in $\Shv(\Bun_G)$).  

\medskip

Hence, we can form also the object
$$\Tr((\on{Frob}_{\Bun_G})_*,\Shv_{\on{Nilp}}(\Bun_G)^{\on{access}})\in \Vect_\sfe.$$

\sssec{}

The composite functor
\begin{equation} \label{e:embed Nilp acess into whole}
\Shv_{\on{Nilp}}(\Bun_G)^{\on{access}}\hookrightarrow 
\Shv_{\on{Nilp}}(\Bun_G)\overset{\iota}\hookrightarrow  \Shv(\Bun_G)
\end{equation}
preserves compactness. In particular, the functor \eqref{e:embed Nilp acess again} also
preserves compactness, and hence both functors \eqref{e:embed Nilp acess again} and
\eqref{e:embed Nilp acess into whole} admit continuous right adjoints. 

\medskip

Since the above functors commute with the action of $(\on{Frob}_{\Bun_G})_*$, by \cite[Sect. 3.4.1]{GKRV}, we obtain maps
$$\Tr((\on{Frob}_{\Bun_G})_*,\Shv_{\on{Nilp}}(\Bun_G)^{\on{access}})\to \Tr((\on{Frob}_{\Bun_G})_*,\Shv_\Nilp(\Bun_G))$$
and
$$\Tr((\on{Frob}_{\Bun_G})_*,\Shv_{\on{Nilp}}(\Bun_G)^{\on{access}})\to \Tr((\on{Frob}_{\Bun_G})_*,\Shv(\Bun_G)).$$

\sssec{}

We propose the following:

\begin{mainconj}  \label{c:Trace conj} \hfill

\smallskip

\noindent{\em(a)} There is a canonical isomorphism 
\begin{equation} \label{e:trace iso}
\Tr((\on{Frob}_{\Bun_G})_*,\Shv_{\on{Nilp}}(\Bun_G)) \simeq \on{Funct}_c(\Bun_G(\BF_q)).
\end{equation}

\smallskip

\noindent{\em(b)} The diagram
$$
\CD
\Tr((\on{Frob}_{\Bun_G})_*,\Shv_{\on{Nilp}}(\Bun_G)^{\on{access}}) @>>> \Tr((\on{Frob}_{\Bun_G})_*,\Shv_{\on{Nilp}}(\Bun_G)) \\
@VVV @VV{\text{\eqref{e:trace iso}}}V \\
\Tr((\on{Frob}_{\Bun_G})_*,\Shv_{\on{Nilp}}(\Bun_G)) @>{\on{LT}}>> \on{Funct}_c(\Bun_G(\BF_q))
\endCD
$$
commutes. 
\end{mainconj} 

In what follows we will use the notation
$$\on{Autom}:=\on{Funct}_c(\Bun_G(\BF_q))$$
and
$$\wt{\on{Autom}}:=\Tr((\on{Frob}_{\Bun_G})_*,\Shv_{\on{Nilp}}(\Bun_G)).$$

So, the statement of \conjref{c:Trace conj}(a) is that we have a canonical isomorphism
$$\wt{\on{Autom}}\simeq \on{Autom}.$$

\sssec{}

Let us explain a concrete meaning of point (b) of \conjref{c:Trace conj}. 

\medskip

Let $\CF$ be a compact object
of $\Shv_{\on{Nilp}}(\Bun_G)$, which is bounded below as an object of $\Shv(\Bun_G)$. In particular,
$\CF$ belongs to $\Shv_{\on{Nilp}}(\Bun_G)^{\on{access}}$ and is compact there, so it is compact
also as an object of $\Shv(\Bun_G)$, and, in particular, constructible. 

\medskip

Assume that $\CF$ is equipped with a lax Weil structure as in \eqref{e:lax Weil}. Then we can consider the elements 
$$\on{cl}(\CF,\alpha)\in  \Tr((\on{Frob}_{\Bun_G})_*,\Shv_{\on{Nilp}}(\Bun_G))$$
and
$$\on{funct}(\CF,\alpha)\in \on{Funct}_c(\Bun_G(\BF_q)).$$

Now point (b) of \conjref{c:Trace conj} says that these two elements match under the isomorphism 
\eqref{e:trace iso}. 

\begin{rem}

Note that if we knew \conjref{c:Nilp comp gen}, we could formulate \conjref{c:Trace conj} just as saying that the composite map
$$\Tr((\on{Frob}_{\Bun_G})_*,\Shv_{\on{Nilp}}(\Bun_G))\to \Tr((\on{Frob}_{\Bun_G})_*,\Shv(\Bun_G))
\overset{\on{LT}}\to \on{Funct}_c(\Bun_G(\BF_q))$$
is an isomorphism.

\end{rem}

\begin{rem}
Note that \conjref{c:Trace conj} defines a direct bridge from the geometric Langlands theory 
to the classical one, since it implies that the space of automorphic functions with compact support, can be expressed
as the categorical trace of the Frobenius endofunctor acting on the category of sheaves on $\Bun_G$ with nilpotent
singular support.

\medskip

Such a bridge allows us to transport structural assertions about $\Shv_{\on{Nilp}}(\Bun_G)$ as a category, to assertions 
about $\on{Autom}$ as a vector space. We will see some examples of this in \secref{s:spectral}.

\end{rem}

\sssec{}

There are several pieces of evidence towards the validity of \conjref{c:Trace conj}. 

\medskip

\noindent(I) It is true when $G$ is a torus. We will analyze this case in the next subsection. 

\medskip

\noindent(II) It is true when $X$ is of genus $0$. Indeed, in this case the inclusion
$$\Shv_{\on{Nilp}}(\Bun_G)\subset \Shv(\Bun_G)$$ 
is an equality, and the map $\on{LT}$ is an isomorphism because
\begin{equation} \label{e:ten on Bun_G}
\Shv(\Bun_G)\otimes \Shv(\Bun_G)\to \Shv(\Bun_G\times \Bun_G)
\end{equation}
is an equivalence, see Remark \ref{r:why naive trace}. 

\medskip

\noindent(III) We have seen in Remark \ref{r:why naive trace} that the failure of the map
$\on{LT}$ originates in the failure of the functor \eqref{e:ten on Bun_G} to be equivalence.
Now, this obstruction goes away for $\Shv_{\on{Nilp}}(\Bun_G)$ thanks to \thmref{t:tensor product}.

\medskip

\noindent(IV) The lisseness property of the cohomology of shtukas, recently established by C.~Xue
in \cite{Xue2}, see Remark \ref{r:Cong1}. 

\ssec{The case of $G=\BG_m$}

In this subsection we will verify by hand the assertion of \conjref{c:Trace conj}(a) for $G=\BG_m$. 
In \secref{s:abelian varieties} we will reprove it by a different method. 

\medskip

To simplify the notation, we will work not with the entire $\Bun_{\BG_m}\simeq \on{Pic}$, but with its
neutral connected component $\on{Pic}_0$.

\sssec{}

Recall that according to \secref{sss:ex Gm}, the category 
$$\Shv_{\{0\}}(\on{Pic}_0)=\qLisse(\on{Pic}_0)$$
is the direct sum over isomorphism classes of $\BG_m$-local systems $\sigma$ of copies of 
\begin{equation} \label{e:autom summand}
\left(\Sym(H^1(X,\ul\sfe_X)[-1])\mod\right) \otimes (\on{C}_\cdot(\BG_m)\mod),
\end{equation} 
where for every $\sigma$, we send the module
$$\Sym(H^1(X,\ul\sfe_X)[-1]) \otimes \sfe \in 
\left(\Sym(H^1(X,\ul\sfe_X)[-1])\mod\right) \otimes \left(\on{C}_\cdot(\BG_m)\mod\right)$$
(where $\sfe$ denotes the augmentation module over $\on{C}_\cdot(\BG_m)$)
to the irreducible Hecke eingensheaf $$E_\sigma\in \qLisse(\on{Pic}_0)$$ corresponding to $\sigma$.

\sssec{}

When we compute $\on{Tr}((\Frob_{\on{Pic}_0})_*,-)$ on this category, only the direct summands,
for which
\begin{equation} \label{e:Frob-inv local system}
\Frob_X^*(\sigma)\simeq \sigma
\end{equation}
can contribute. 

\medskip

For each such $\sigma$ choose an isomorphism in \eqref{e:Frob-inv local system}. This choice defines 
a Weil sheaf structure on the corresponding $E_\sigma$. Further, this choice identifies the action of
$\on{Tr}((\Frob_{\on{Pic}_0})_*,-)$ on the direct summand \eqref{e:autom summand} with the action induced
by the Frobenius automorphism of the algebra
$$A:=\Sym(H^1(X,\ul\sfe_X)[-1]) \otimes \on{C}_\cdot(\BG_m).$$

\medskip

We will show that
\begin{equation} \label{e:Tr Frob A mod}
\Tr(\Frob,A\mod)\simeq \sfe,
\end{equation} 
and that the induced map
\begin{multline} \label{e:Tr Frob eigensheaf}
\sfe \mapsto \Tr(\Frob,A\mod) \to  \on{Tr}((\Frob_{\on{Pic}_0})_*,\qLisse(\on{Pic}_0)) \to \\
\to \on{Tr}((\Frob_{\on{Pic}_0})_*,\Shv(\on{Pic}_0)) \overset{\on{LT}}\to \on{Funct}(\on{Pic}_0(\BF_q))
\end{multline} 
sends $1\in \sfe$ to $\on{funct}(E_\sigma)\cdot (1-q)$. 

\medskip

This will prove the required assertion since the functions $\on{funct}(E_\sigma)$ form a basis of 
$\on{Funct}(\on{Pic}_0(\BF_q))$, by Class Field Theory.

\sssec{} 

We have 
$$A=A_1\otimes A_2, \quad A_1=\Sym(H^1(X,\ul\sfe_X)[-1]) ,\,\, A_2=\on{C}_\cdot(\BG_m).$$

This corresponds to writing 
$$\on{Pic}_0\simeq \on{Jac}(X) \times \on{pt}/\BG_m.$$

We will perform the calculation for each factor separately.

\sssec{}

Note that if $A'$ is a polynomial algebra
$$A'\simeq \Sym(V),$$
where $V$ is a finite-dimensional vector space 
equipped with an endomorphism $F$ with no eigenvalue $1$, then the functor
$$\Vect_\sfe\to A'\mod, \quad \sfe\mapsto A'$$
defines an isomorphism
$$\sfe\simeq \Tr(\on{Id},\Vect_\sfe)\to \Tr(F,A'\mod).$$

Applying this to $A'=A_1$ and $A'=A_2$, we obtain the desired identifications
$$\Tr(\Frob,A_1\mod)\simeq \sfe \text{ and } \Tr(\Frob,A_2\mod)\simeq \sfe,$$
as required in \eqref{e:Tr Frob A mod}. 

\sssec{}

To prove \eqref{e:Tr Frob eigensheaf} for $A_1$, we consider the composite functor
$$\Vect_\sfe \overset{\sfe\mapsto A_1}\longrightarrow A_1\mod \to \qLisse(\on{Pic}_0)\to \Shv(\on{Pic}_0),$$
equipped with its datum of compatibility with the Frobenius. 

\medskip

It sends
$$\sfe \mapsto E_\sigma,$$
equipped with its Weil structure, to be denoted $\alpha$.

\medskip

Hence, the corresponding map
$$\sfe\simeq \Tr(\on{Id},\Vect_\sfe) \to \on{Tr}((\Frob_{\on{Pic}_0})_*,\Shv(\on{Pic}_0))$$
sends 
$$1\in \sfe\,\,  \mapsto \,\, \on{cl}(E_\sigma,\alpha)\in \on{Tr}((\Frob_{\on{Pic}_0})_*,\Shv(\on{Pic}_0)).$$

Hence, its image under $\on{LT}$ is $\on{funct}(E_\sigma)$.

\sssec{}

We now consider $A_2$. The composite functor
$$\Vect_\sfe \overset{\sfe\mapsto A_2}\longrightarrow A_2\mod \to \Shv(\on{pt}/\BG_m)\overset{\text{pullback}}\to \Shv(\on{pt})=\Vect_\sfe$$
sends 
$$\sfe\mapsto \on{C}_\cdot(\BG_m),$$
equipped with the natural datum of compatibility with the Frobenius. 

\medskip

Hence the resulting map
$$\sfe\simeq \Tr(\on{Id},\Vect_\sfe) \to \Tr((\on{Frob}_{\on{pt}/\BG_m})_*,\Shv(\on{pt}/\BG_m))\to \Tr((\Frob_{\on{pt}})_*,\Shv(\on{pt}))=
\Tr(\on{Id},\Vect_\sfe)\simeq \sfe$$
sends
$$1\in \sfe \,\, \mapsto \Tr(\Frob,\on{C}_\cdot(\BG_m))=1-q\in \sfe.$$

\ssec{A generalization: cohomologies of shtukas} \label{ss:shtukas}

In this subsection we will formulate a generalization of the Trace Conjecture, which gives a trace
interpretation to cohomologies of shtukas. 

\sssec{} \label{sss:shtukas}

Let us recall the construction of cohomologies of shtukas, following \cite{VLaf1} and \cite{Var1}. 

\medskip 

Let $I$ be a finite set and $V$ an object of $\Rep(\cG)^{\otimes I}$. To this data we attach an object
$$\on{Sht}_{I,V}\in \Shv(X^I)$$
as follows.

\medskip

We consider the $I$-legged Hecke stack
$$
\CD 
\Bun_G @<{\hl}<< \on{Hecke}_{X^I} @>{\hr}>> \Bun_G \\
& & @V{\pi}VV \\
& & X^I.
\endCD
$$

\medskip

The $I$-legged shtuka space is defined as the fiber product
$$
\CD
\on{Sht}_I @>>>  \on{Hecke}_{X^I}  \\
@VVV @VV{(\hl,\hr)}V \\
\Bun_G @>{(\on{Frob}_{\Bun_G},\on{Id})}>> \Bun_G\times \Bun_G. 
\endCD
$$

Let $\pi'$ denote the composite map
$$\on{Sht}_I \to \on{Hecke}_{X^I} \overset{\pi}\to X^I.$$

\medskip

Recall that (the classical\footnote{As opposed to derived.}) geometric Satake attaches to $V\in \Rep(\cG)^{\otimes I}$ an object
$$\CS_V\in \Shv(\on{Hecke}_{X^I}).$$

Let $\CS'_V\in \Shv(\on{Sht}_I)$ denote its *-restriction to $\Shv(\on{Sht}_I)$. Finally, we set
$$\on{Sht}_{I,V}:=\pi'_!(\CS'_V)\in \Shv(X^I).$$

\begin{rem} \label{r:Cong1}
A recent result of \cite{Xue2} says that the objects $\on{Sht}_{I,V}$ actually belong to
$$\qLisse(X^I)\subset \Shv(X^I).$$
\end{rem}

\sssec{Example} \label{sss:no leg shtuka}

Take $I=\emptyset$ and $V$ to be $\sfe$. Then 
$$\on{Sht}_\emptyset\simeq (\Bun_G)^{\on{Frob}}\simeq \Bun_G(\BF_q).$$

We obtain that
\begin{equation} \label{e:no leg shtuka}
\on{Sht}_{\emptyset,\sfe}=\on{C}^\cdot_c(\Bun_G(\BF_q),\ul\sfe_{\Bun_G(\BF_q)})\simeq 
\on{Funct}_c(\Bun_G(\BF_q))=\on{Autom}.
\end{equation}

\sssec{} 

We will now construct a different system of objects
$$\wt{\on{Sht}}_{I,V}\in \qLisse(X^I).$$

\sssec{}

Note that the categorical trace construction has the following variant. Let $\bC$ be a dualizable DG category
and let
$$F:\bC\to \bC\otimes \bD,$$
where $\bD$ is some other DG category. 

\medskip

Then we can consider an object 
$$\Tr(F,\bC)\in \bD.$$

Namely, $\Tr(F,\bC)$ is the composition
$$\Vect_\sfe\overset{\on{unit}}\longrightarrow \bC^\vee\otimes \bC
\overset{\on{Id}\otimes F}\longrightarrow \bC^\vee\otimes \bC\otimes \bD 
\overset{\on{counit}\otimes \on{Id}_\bD}\longrightarrow \bD.$$

(The usual trace construction is when $\bD=\Vect_\sfe$, so $\Tr(F,\bC)\in \Vect_\sfe$.)

\sssec{} \label{sss:sht tilde}

We apply this to $\bC:=\Shv_{\on{Nilp}}(\Bun_G)$, $\bD=\qLisse(X^I)$ and $F$ being the functor
$$\Shv_{\on{Nilp}}(\Bun_G) \overset{(\on{Frob}_{\Bun_G})_*}\longrightarrow 
\Shv_{\on{Nilp}}(\Bun_G)  \overset{\on{H}(V,-)}\longrightarrow \Shv_{\on{Nilp}}(\Bun_G)\otimes \qLisse(X^I).$$

\medskip

We set
$$\wt{\on{Sht}}_{I,V}:=\Tr(\on{H}(V,-) \circ (\on{Frob}_{\Bun_G})_*,\Shv_{\on{Nilp}}(\Bun_G))\in \qLisse(X^I).$$

\medskip

We propose:

\begin{mainconj} \label{c:Trace conj legs}
The objects $\on{Sht}_{I,V}$ and $\wt{\on{Sht}}_{I,V}$
are canonically isomorphic.
\end{mainconj}

\sssec{} \label{sss:vs comp supp}

Consider the case of $I=\emptyset$. As we have seen in \secref{sss:no leg shtuka}, 
$$\on{Sht}_{\emptyset,\sfe}=\on{C}^\cdot_c(\Bun_G(\BF_q),\ul\sfe_{\Bun_G(\BF_q)})\simeq 
\on{Funct}_c(\Bun_G(\BF_q))=\on{Autom}.$$

This is while, 
$$\wt{\on{Sht}}_{\emptyset,\sfe}=\Tr((\on{Frob}_{\Bun_G})_*,\Shv_{\on{Nilp}}(\Bun_G))=:\wt{\on{Autom}},$$
which, according to \conjref{c:Trace conj}(a), is isomorphic to $\on{Autom}$.

\medskip

So, \conjref{c:Trace conj}(a) is a special case of \conjref{c:Trace conj legs}. 
%
%
%
%
%
%

\begin{rem} \label{c:Cong}

A crucial piece of evidence for the validity of \conjref{c:Trace conj legs} is provided by the result of \cite{Xue2}
mentioned in Remark \ref{r:Cong1}.

\end{rem}

\sssec{Partial Frobeniuses} \label{sss:partial Frobs}

Recall (see \cite[Sect. 3]{VLaf1}) that the objects $\on{Sht}_{I,V}$ carry an additional structure, namely, equivariance
with respect to the \emph{partial Frobenius maps}.

\medskip

The construction from \cite[Sect. 5.3]{GKRV} endows the objects $\wt{\on{Sht}}_{I,V}$ with a similar structure. 
(See also Remark \ref{r:action by excursions} for a conceptual explanation of this structure.) 

\medskip

The statement of \conjref{c:Trace conj legs} should be strengthened as follows: the isomorphism
$$\on{Sht}_{I,V}\simeq \wt{\on{Sht}}_{I,V}$$
is compatible with the structure of equivariance with respect to the partial Frobenius maps. 

\section{The trace conjecture for abelian varieties} \label{s:abelian varieties} 

In this section we will prove a statement parallel to \conjref{c:Trace conj} for 
the category of lisse sheaves on an abelian variety. 

\medskip

The material of this section is not logically related to the contents of the rest of the paper. 

\ssec{Statement of the result}

\sssec{}

Let $A$ be an abelian variety. Consider the subcategory
\begin{equation} \label{e:iota A}
\qLisse(A)\overset{\iota_A}\hookrightarrow \Shv(A).
\end{equation} 

The category $\qLisse(A)$ is compactly generated by the character sheaves, to be denoted $E_\sigma$
(see \secref{sss:ident Ab} below). From here it follows that the pair $(A,\{0\})$ is constraccessible (see Definition \ref{d:constraccessible}). 
In particular, the embedding $\iota_A$ admits a continuous right adjoint. 

\sssec{}

We now take our ground field $k$ to be $\ol\BF_q$, but we assume that $A$ is defined over $\BF_q$, 
so it carries an action of the geometric Frobenius endomorphism $\Frob_A$. The map \eqref{e:iota A} induces a map
$$\Tr((\Frob_A)_*,\qLisse(A))\to \Tr((\Frob_A)_*,\Shv(A)).$$

\medskip

Consider the composition
\begin{equation} \label{e:Trace for Ab}
\Tr((\Frob_A)_*,\qLisse(A))\to \Tr((\Frob_A)_*,\Shv(A)) \overset{\on{LT}}\longrightarrow \on{Funct}(A(\BF_q)).
\end{equation}

The main result of this section reads:

\begin{thm} \label{t:Trace for Ab}
The map \eqref{e:Trace for Ab} is an isomorphism.
\end{thm} 

The rest of this section is devoted to the proof of this theorem.

\ssec{The right adjoint to the embedding of the lisse subcategory}

In this subsection we will study the right adjoint of the functor $\iota_A$ of \eqref{e:iota A}.
With future applications in mind, we will do this in greater generality than is actually needed for the proof of 
\thmref{t:Trace for Ab}. 

\sssec{} \label{sss:Lisse on X}

Deviating from the notations of the rest of the paper, for the duration of this section, 
we let $X$ denote a smooth scheme (not necessarily a curve). In the applications
to the proof of \thmref{t:Trace for Ab}, we will take $X$ to be the abelian variety $A$. 

\medskip

Consider the embedding
\begin{equation} \label{e:iota X}
\qLisse(X)\overset{\iota_X}\hookrightarrow \Shv(X).
\end{equation} 

We will assume that the pair $(X,\{0\})$ is constraccessible. I.e., $\qLisse(X)$ is generated by 
objects that are compact in $\Shv(X)$. In particular, the above functor $\iota_X$ admits a continuous
right adjoint\footnote{Even if $(X,\{0\})$ is not constraccessible, the discussion below applies if we 
replace the entire $\qLisse(X)$ by an arbitrary full subcategory of $\qLisse(X)$ which is generated by 
objects that are compact in $\Shv(X)$.}. 

\medskip 

We will now give an explicit formula for this right adjoint. 

\sssec{} 

For a pair of objects $E_1,E_2\in \Lisse(X)=\qLisse(X)^c$, consider the functor
$$\sP_{E_1,E_2}: \Shv(X)\to \qLisse(X), \quad \CF\mapsto E_1\otimes \CHom_{\Shv(X)}(E_2,\CF).$$

As $E_1,E_2$ vary we obtain a functor
\begin{equation} \label{e:to take coEnd}
\Lisse(X) \times \Lisse(X)^{\on{op}}  \to \on{Funct}(\Shv(X),\qLisse(X)).
\end{equation} 

Let $\sP_X\in  \on{Funct}(\Shv(X),\qLisse(X))$ be the coEnd of the functor \eqref{e:to take coEnd}, i.e.,
the colimit of \eqref{e:to take coEnd} over the index category $\on{TwArr}(\Lisse(X))$.

\medskip

The following is standard:

\begin{lem} \label{l:coEnd}
The functor $\sP_X$ identifies canonically with $\iota_X^R$.
\end{lem}

\sssec{} \label{sss:functors given by kernels}

Let $Y_1$ and $Y_2$ be a pair of quasi-compact schemes. For an object $\CQ\in \Shv(Y_1\times Y_2)$ we will
denote by $\sK_\CQ$ the functor
$$\Shv(Y_1)\to \Shv(Y_2), \quad \CF\mapsto (p_2)_*(p_1^!(\CF)\sotimes \CQ).$$

Let $Z$ be yet another scheme. We will denote by $\sK_\CQ \boxtimes \on{Id}_Z$ the functor
$$\Shv(Y_1\times Z)\to \Shv(Y_2\times Z)$$
equal to $\sK_{\CQ\boxtimes (\Delta_Z)_*(\omega_Z)}$, where 
$$\CQ\boxtimes (\Delta_Z)_*(\omega_Z)\in \Shv(Y_1\times Y_2\times Z\times Z)\simeq
\Shv((Y_1\times Z)\times (Y_2\times Z)).$$

\sssec{}

For $E_1,E_2\in \Lisse(X)$, set 
$$\CQ_{E_1,E_2}:=\BD(E_2)\boxtimes E_1\in \Shv(X\times X),$$
so that the functor $\sP_{E_1,E_2}$ above identifies with $\sK_{\CQ_{E_1,E_2}}$.

\medskip

The assignment
$$E_1,E_2 \mapsto \CQ_{E_1,E_2}$$ is a functor
$$\Lisse(X) \times \Lisse(X)^{\on{op}}  \to \Shv(X\times X),$$
and let $\CQ_{\qLisse(X)}$ be its coEnd. By construction, 
$$\sK_{\CQ_{\qLisse(X)}} \simeq \iota_X\circ \sP_X.$$

\sssec{}

Note that for a scheme $Z$, the essential image of the endofunctor 
$\sK_{\CQ_{\qLisse(X)}} \boxtimes \on{Id}_Z$ of $\Shv(X\times Z)$
is contained in the full subcategory
$$\qLisse(X)\otimes \Shv(Z)\hookrightarrow \Shv(X\times Z).$$

\medskip

We claim:

\begin{prop} \label{p:lisse projector product}
For a scheme $Z$, the endofunctor $\sK_{\CQ_{\qLisse(X)}} \boxtimes \on{Id}_Z$ of $\Shv(X\times Z)$
identifies with the (pre)composition of the fully faithful embedding
$$\qLisse(X)\otimes \Shv(Z)\hookrightarrow \Shv(X\times Z)$$
and its right adjoint.
\end{prop} 

\begin{proof}

We need to establish a functorial isomorphism
\begin{equation} \label{e:right adjoint prod}
\CHom_{\qLisse(X)\otimes \Shv(Z)}(E\boxtimes \CF_Z,(\sK_{\CQ_{\qLisse(X)}}\boxtimes \on{Id}_Z)(\CF))
\simeq \CHom_{\Shv(X\times Z)}(E\boxtimes \CF_Z,\CF)
\end{equation} 
for $E\in \Lisse(X),\,\, \CF_Z\in \Shv(Z)^c,\,\, \CF\in \Shv(X\times Z)$. 

\medskip

Set
$$\CF_X:=(p_X)_*(\CF\sotimes p_Z^!(\BD(\CF_Z))).$$

Then we can rewrite the left-hand side in \eqref{e:right adjoint prod} as
$$\CHom_{\qLisse(X)}(E,\sP_X(\CF_X)),$$
and the right-hand side as
$$\CHom_{\Shv(X)}(E,\CF_X).$$

So the assertion of the proposition follows from \lemref{l:coEnd}. 

\end{proof}

\ssec{Lisse sheaves an abelian variety}

Let us now specialize to the case when $X=A$ is an abelian variety. 

\sssec{} \label{sss:ident Ab prel}

Let $\sigma$ denote an isomorphism class of $1$-dimensional local systems on $A$, and let 
$$E_\sigma\in \qLisse(A)$$
be an object in the given isomorphism, canonically fixed by the requirement that its *-fiber at $1\in A$
is identified with $\sfe$.

\medskip

The following are standard facts about local systems on an abelian variety:

\begin{itemize} 

\item Each $E_\sigma$ has a (unique) structure of character sheaf, i.e., it is equipped
with an isomorphism
$$\on{mult}_A^*(E_\sigma)\simeq E_\sigma\boxtimes E_\sigma,$$
normalized so that its fiber at $1\times 1\in A\times A$ is the identity map
(it then automatically satisfies the associativity requirement);

\item Each irreducible object in $\Lisse(A)^\heartsuit$ is isomorphic to some $E_\sigma$;

\item $\CHom(E_{\sigma_1},E_{\sigma_2})=0$ if $\sigma_1\neq \sigma_2$;

\item $\CHom(E_\sigma,E_\sigma)\simeq \on{C}^\cdot(A,\sfe)\simeq \Sym(H^1(A,\sfe)[-1])$,
as associative algebras.

\end{itemize}

The second and the third points above imply that the category $\iLisse(A)$ splits as a direct sum
$$\underset{\sigma}\oplus\, \iLisse(A)_\sigma,$$
where each $\iLisse(A)_\sigma$ is compactly generated by $E_\sigma$.

\sssec{} \label{sss:ident Ab}

We now claim that the embedding
$$\iLisse(A) \hookrightarrow \qLisse(A)$$
is an equivalence, i.e., the objects $E_\sigma$ (compactly) generate all of $\qLisse(X)$.

\medskip

For that it suffices to show that each $\iLisse(A)_\sigma$ is left-complete in the t-structure 
(induced by the usual t-structure on $\qLisse(A)$.

\medskip

However, this t-structure, when viewed as a t-structure on
$$\iLisse(A)_\sigma \simeq \Sym(H^1(A,\sfe)[-1])\mod$$
translates via the Koszul duality
$$\Sym(H^1(A,\sfe)[-1])\mod\simeq \Sym(H^1(A,\sfe)^\vee)\mod_{\{0\}}$$
to the t-structure on $\Sym(H^1(A,\sfe)^\vee)\mod_{\{0\}}$ compatible
with the forgetful functor
$$\Sym(H^1(A,\sfe)^\vee)\mod_{\{0\}}\hookrightarrow \Sym(H^1(A,\sfe)^\vee)\mod\to \Vect_\sfe,$$
and the latter is manifestly left-complete. 

\sssec{}

The above description of $\qLisse(A)$ implies that 
we can rewrite the object $\sK_{\CQ_{\qLisse(A)}}$ more concisely. Namely, we have
$$\CQ_{\qLisse(A)} \simeq \underset{\sigma}\oplus\, \BD(E_\sigma)\underset{\on{C}^\cdot(A)}\boxtimes E_\sigma.$$

In particular, the functor $\iota_A^R$ acts as
\begin{equation} \label{e:projector for Abelian}
\CF \mapsto  \underset{\sigma}\oplus\, \left(E_\sigma \underset{\on{C}^\cdot(A)}\otimes \CHom_{\Shv(A)}(E_\sigma,\CF)\right).
\end{equation}

\sssec{}

The operation of convolution defines on $\Shv(A)$ a structure of symmetric monoidal category. We will denote
the corresponding binary operation by $\star$ (in order to distinguish it from the pointwise symmetric monoidal
structure, which is denoted by $\sotimes$).

\medskip

The full subcategory $\qLisse(A)$ is preserved by $\star$, and is in fact a monoidal ideal. Hence, the functor
$\iota^R_A$ acquires a right-lax symmetric monoidal structure. 

\medskip

We claim:

\begin{lem}
The right-lax symmetric monoidal structure on $\iota^R_A$ is strict.
\end{lem}

\begin{proof}
Follows from the fact that each $E_\sigma$ is a character sheaf.
\end{proof}

\sssec{} \label{sss:delta lisse}

Thus, we obtain that $\iota_A$ is a symmetric monoidal co-localization. Denote by
$$\delta_{1,\qLisse}\in \qLisse(A)$$
the object 
$$\iota_A^R(\delta_1).$$

It follows formally that the functor $\iota^R_A$ is given by convolution with $\delta_{1,\qLisse}$:
$$\iota^R_A(\CF) \simeq \CF\star \delta_{1,\qLisse}.$$

\medskip

Similarly, we obtain that the object
$$\CQ_{\qLisse(A)}\in \Shv(A\times A)$$
is obtained by convolving $(\Delta_A)_*(\omega_A)$ with $\delta_{1,\qLisse}$ along the 
second factor. 

\begin{rem}

By formula \eqref{e:projector for Abelian}, the object $\delta_{1,\qLisse}$ is a direct sum
$$\underset{\sigma}\oplus\,  \delta_{1,\sigma}, \quad 
\delta_{1,\sigma}\simeq E_\sigma\underset{\on{C}^\cdot(A)}\otimes\, \sfe,$$
where $\sfe$ is the augmentation module for $\on{C}^\cdot(A)$, corresponding to
the point $1\in A$.

\medskip

So, each $\delta_{1,\sigma}$ is an infinite Jordan block. It is naturally filtered and the
associated graded identifies with
$$E_\sigma \otimes \Sym(H^1(A,\sfe)).$$

\end{rem} 

\ssec{Calculating the trace of Frobenius}

We return to the general set-up of \secref{sss:Lisse on X}, where we now assume that the ground field is $\ol\BF_q$,
but the scheme $X$ is defined over $\BF_q$, so we can talk about
$$\Tr((\Frob_X)_*,\qLisse(X))\in \Vect_\sfe.$$

In this subsection we will produce an explicit formula for this trace. 

\sssec{}

Let $\CQ\in \Shv(Y\times Y)$ and $\sK_\CQ$ be as in \secref{sss:functors given by kernels}. We introduce the following notation
$$\Tr_{\on{geom}}(\CQ,Y):=\Gamma(Y,\Delta_Y^!(\CQ)).$$

Note that we have a canonically defined map
\begin{equation} \label{e:from geom to normal}
\Tr(\sK_\CQ,\Shv(Y))\to \Tr_{\on{geom}}(\CQ,Y) ,
\end{equation}
which is functorial in $\CQ$. 

\medskip

Indeed, this map comes from the map \eqref{e:map from unit to diag}.

\sssec{}

Note that for 
$$\CQ=\CQ_{\Frob_Y}:=(\on{Id} \times \Frob_Y)_*((\Delta_Y)_*(\omega_Y)) \simeq (\Frob_Y\times \on{Id})^!((\Delta_Y)_*(\omega_Y)),$$
so that
$$\sK_{\CQ_{\Frob_Y}}=(\Frob_Y)_*,$$
we have a canonical identification
$$\Tr_{\on{geom}}(\CQ_{\Frob_Y},Y) \simeq \on{Funct}(Y(\BF_q)).$$

By construction, the map $\on{LT}^{\on{true}}$ is the map \eqref{e:from geom to normal}
for $\CQ_{\Frob_Y}$.

\sssec{}

We take $Y=X$ and 
$$\CQ_1=\CQ_{\Frob_X}$$
and 
$$\CQ_2=\CQ_{\Frob_X,\qLisse}:=(\on{Id}_X \boxtimes \sK_{\qLisse(X)})(\CQ_1),$$
so that
$$\sK_{\CQ_2}=\sK_{\qLisse(X)}\circ (\Frob_X)_*.$$


The counit of the adjunction of \propref{p:lisse projector product} defines a map
$$\CQ_2\to \CQ_1.$$

In particular, we obtain a commutative diagram
\begin{equation} \label{e:com diag trace}
\CD
\Tr(\sK_{\qLisse(X)}\circ (\Frob_X)_*,\qLisse(X)) @>>> \Tr((\Frob_X)_*,\qLisse(X))  \\
@VVV @VVV \\
\Tr(\sK_{\qLisse(X)}\circ (\Frob_X)_*,\Shv(X)) @>>> \Tr((\Frob_X)_*,\Shv(X)) \\
@V{\text{\eqref{e:from geom to normal}}}VV @VV{\text{\eqref{e:from geom to normal}}}V \\
\Tr_{\on{geom}}(\CQ_{\Frob_X,\qLisse},X) @>>> \Tr_{\on{geom}}(\CQ_{\Frob_X},X) \\
& & @VV{\sim}V \\
& & \on{Funct}(X(\BF_q))
\endCD
\end{equation} 

Note that the right vertical composition in \eqref{e:com diag trace} is that map 
$$\Tr((\Frob_X)_*,\qLisse(X))\to \Tr((\Frob_X)_*,\Shv(X)) \overset{\on{LT}}\longrightarrow \on{Funct}(X(\BF_q))$$
of \eqref{e:Trace for Ab} when $X=A$. 

\sssec{}

The proof of \thmref{t:Trace for Ab} will amount to showing that:

\medskip

\begin{itemize}

\item The top horizontal arrow in \eqref{e:com diag trace} is an isomorphism (this is true for any $X$);

\item The upper left vertical arrow in \eqref{e:com diag trace} is an isomorphism (this is true for any $X$);

\item The lower left vertical arrow in \eqref{e:com diag trace} is an isomorphism (this is true for any $X$);

\item The bottom horizontal arrow in \eqref{e:com diag trace} is an isomorphism when $X=A$ is an 
abelian variety.

\end{itemize} 

Note that the combination of the first three isomorphisms implies that we have a canonical isomorphism
$$\Tr((\Frob_X)_*,\qLisse(X)) \simeq \Tr_{\on{geom}}(\CQ_{\Frob_X,\qLisse},X).$$

\sssec{}

The fact that the top horizontal arrow in \eqref{e:com diag trace} is an isomorphism is immediate:
by \propref{p:lisse projector product}, the functor $\sK_{\qLisse(X)}$ acts as identity on $\qLisse(X)$. 

\medskip

Similarly, the fact that the upper left vertical arrow in \eqref{e:com diag trace} is an isomorphism 
follows from the fact that the functor $\sK_{\qLisse(X)}$ identifies with the composition
$\iota_X\circ \iota_X^R$. 

\sssec{}

To prove that the lower left vertical arrow in \eqref{e:com diag trace} is an isomorphism, it suffices to show that
the map
\begin{equation} \label{e:sk map}
(\on{Id}_X \boxtimes \sK_{\qLisse(X)}) (\boxtimes (\on{u}_{\Shv(X)})) \to 
(\on{Id}_X \boxtimes \sK_{\qLisse(X)})  ((\Delta_X)_*(\omega_X)),
\end{equation} 
induced by  \eqref{e:map from unit to diag}, is an isomorphism. 

\medskip

Let us denote by $\Phi$ the $\boxtimes$ functor
$$\Shv(X)\otimes \Shv(X)\to \Shv(X\times X).$$ 

Then by \propref{p:lisse projector product}, the map \eqref{e:sk map} identifies with the value of the natural transformation 
\begin{equation} \label{e:sk map bis}
(\Phi \circ (\on{id}\otimes \iota_X)) \circ (\Phi \circ (\on{id}\otimes \iota_X))^R \circ \Phi \circ \Phi^R \to 
(\Phi \circ (\on{id}\otimes \iota_X)) \circ (\Phi \circ (\on{id}\otimes \iota_X))^R
\end{equation} 
on $(\Delta_X)_*(\omega_X)$. 

\medskip

We claim that the natural transformation \eqref{e:sk map bis} itself is an isomorphism. By definition, \eqref{e:sk map bis} is the map
$$\Phi \circ (\on{id}\otimes \iota_X) \circ  (\on{id}\otimes \iota_X)^R \circ \Phi^R \circ \Phi \circ \Phi^R \to
\Phi \circ (\on{id}\otimes \iota_X) \circ (\on{id}\otimes \iota_X)^R \circ \Phi^R,$$
given by the counit $\Phi \circ \Phi^R \to \on{Id}$. Hence, it suffices to show that 
$$\Phi^R \circ \Phi \circ \Phi^R \to \Phi^R$$
is an isomorphism, and the latter follows from the fact that $\Phi$ is fully faithful. 

\ssec{A calculation using Lang's isogeny}

In this subsection we will prove that the bottom horizontal arrow in \eqref{e:com diag trace} is an isomorphism.  

\sssec{}

Let $\CF$ be an arbitrary object of $\Shv(A)$. For any scheme $Z$, let $-\star \CF$
denote the endofunctor of $\Shv(Z\times A)$ obtained by convolving with $\CF$ along
the $A$ factor. 

\medskip

The following results from a diagram chase:

\begin{lem} \label{l:Lang}
There is a canonical isomorphism 
$$\Tr_{\on{geom}}(\CQ_{\Frob_A}\star \CF ,A)  \simeq \on{C}^\cdot(A,L^!(\CF)),$$ 
where $L:A\to A$ is the Lang isogeny.
\end{lem}

\sssec{}

Let us observe that the object $\CQ_{\Frob_X,\qLisse}$ identifies with $\CQ_{\Frob_X}\star \delta_{1,\qLisse}$
(see \secref{sss:delta lisse}). Under this identification, the map
$$\CQ_{\Frob_X,\qLisse} \to \CQ_{\Frob_X}$$
is obtained from counit of the adjunction.
$$\delta_{1,\qLisse}\to \delta_1.$$

Hence, by \lemref{l:Lang}, in order to show that the bottom horizontal arrow in \eqref{e:com diag trace} is an isomorphism,
we have to show that the map
\begin{equation} \label{e:Lang}
\on{C}^\cdot(A,L^!(\delta_{1,\qLisse}))\to \on{C}^\cdot(A,L^!(\delta_1))
\end{equation}
is an isomorphism.

\sssec{}

We rewrite the map \eqref{e:Lang} as
$$\CHom_{\Shv(A)}(\ul\sfe_A,L^!(\delta_{1,\qLisse}))\to \CHom_{\Shv(A)}(\ul\sfe_A,L^!(\delta_{1})),$$
and further by adjunction as 
$$\CHom_{\Shv(A)}(L_!(\ul\sfe_A),\delta_{1,\qLisse})\to \CHom_{\Shv(A)}(L_!(\ul\sfe_A),\delta_{1}).$$

Now, the latter map is indeed an isomorphism because
$$L_!(\ul\sfe_A)\in \qLisse(A),$$
since $L$ is a finite \'etale map. 

\qed[\thmref{t:Trace for Ab}]

\section{Localization of the space of automorphic functions} \label{s:spectral}

In this section we will introduce the space (in fact, a quasi-compact algebraic stack) of
\emph{arithmetic} Langlands parameters, denoted $\LocSys^{\on{arithm}} _\cG(X)$. 

\medskip

We will see how our Trace Conjecture leads to a \emph{localization} of the space of automorphic
forms onto $\LocSys^{\on{arithm}} _\cG(X)$.

\ssec{The arithmetic $\LocSys^{\on{restr}}_\cG$}


\sssec{} \label{sss:Frob on LocSys}

Consider the automorphism of the symmetric monoidal category $\qLisse(X)$, given by
pullback with respect to the Frobenius endomorphism of $X$:
$$\on{Frob}^*_X: \qLisse(X)\to \qLisse(X).$$

\medskip

By transport of structure, the prestack $\LocSys^{\on{restr}}_\cG(X)$ acquires an automorphism, which we will
denote simply by $\Frob$. 

\sssec{}

Consider the prestack
$$(\LocSys^{\on{restr}}_\cG(X))^{\on{Frob}}$$
of $\on{Frob}$-fixed points of $\LocSys^{\on{restr}}_\cG(X)$.

\medskip

Note that $\sfe$-points of $(\LocSys^{\on{restr}}_\cG(X))^{\on{Frob}}$ are Weil $\cG$-local systems on $X$, i.e., 
$\cG$-local systems on $X$ equipped with a Weil structure.

\sssec{}

In \secref{s:arithm} we will prove:

\begin{thm} \label{t:Frob-finite} 
The fixed-point locus
$(\LocSys^{\on{restr}}_\cG(X))^{\on{Frob}}$ is a quasi-compact, mock-affine\footnote{See \secref{sss:mock-affine} for what this means.} 
algebraic stack, locally almost of finite type.
\end{thm} 

\sssec{}

In the same \secref{s:arithm} we will also prove:

\begin{thm} \label{t:irred Weil} Assume that $\cG$ is semi-simple. Let 
$\sigma$ be an $\sfe$-point of $(\LocSys^{\on{restr}}_\cG(X))^{\on{Frob}}$, which is
irreducible as a Weil $\cG$-local system. Then the group of its 
automorphisms is finite, and the resulting map 
$$\on{pt}/\on{Aut}(\sigma)\to (\LocSys^{\on{restr}}_\cG(X))^{\on{Frob}}$$
is the embedding of a connected component.
\end{thm}

Combining with the quasi-compactness assertion from \thmref{t:Frob-finite}, we obtain:

\begin{cor}  \label{c:irred fin}
Let $\cG$ be semi-simple. Then there is only a finite number of irreducible Weil $\cG$-local systems on $X$.
\end{cor}

\sssec{}

We will think of $(\LocSys^{\on{restr}}_\cG(X))^{\on{Frob}}$ as the stack parameterizing $\cG$-local systems on $X$ equipped 
with a Weil structure, and henceforth denote it by
$$\LocSys^{\on{arithm}} _\cG(X).$$ 

\begin{rem}
We propose $\LocSys^{\on{arithm}} _\cG(X)$ as a candidate for the stack $\CS/\cG$, alluded to in \cite[Remark 8.5]{VLaf2}.

\medskip

Recently, P.~Scholze (unpublished) and X.~Zhu (in \cite{Zhu}) proposed two more definitions
of the stack of Weil $\cG$-local systems on $X$. Their definitions are different from each other,
and are of completely different flavor from ours. It is likely, however, that the resulting three versions
of $\LocSys^{\on{arithm}} _\cG(X)$ are actually equivalent. 

\end{rem} 

\begin{rem}

As we shall see in \secref{sss:arithm non-classical}, the stack 
$\LocSys^{\on{arithm}} _\cG(X)$ is \emph{non-classical}, i.e., its structure sheaf
has non-trivial negative cohomology. 

\end{rem}

\ssec{The excursion algebra} \label{ss:Exc}

\sssec{}

Denote
$$\Exc:=\Gamma(\LocSys^{\on{arithm}} _\cG(X),\CO_{\LocSys^{\on{arithm}} _\cG(X)}).$$

This is a commutative algebra object in $\Vect_\sfe$. Since $\LocSys^{\on{arithm}} _\cG(X)$ is mock-affine, 
the algebra $\Exc$ is connective. 

\medskip

Set
\begin{equation} \label{e:arithm coarse}
\LocSys^{\on{arithm,coarse}}_\cG(X):=\Spec(\Exc).
\end{equation} 

This is the coarse moduli space of arithmetic Langlands parameters. By construction, it is a \emph{derived}
affine scheme. 

\sssec{} \label{sss:naive Exc}

The algebra $\Exc$ is related to V.~Lafforgue's algebra of excursion operators as follows. 

\medskip

Let $\on{Weil}(X,x)^{\on{discr}}$ be the Weil group of $X$ (for some choice of a base point $x\in X$), considered as a discrete group. Set
$$\CX^{\on{discr}}:=B(\on{Weil}(X,x)^{\on{discr}})\in \Spc.$$

Consider the (mock-affine) algebraic stack
$$\LocSys^{\on{Betti}}_\cG(\CX^{\on{discr}})\simeq \LocSys^{\on{Betti,rigid}_x}_\cG(\CX^{\on{discr}})/\cG,$$
and set
$$\Exc^{\on{discr}}:=\Gamma(\LocSys^{\on{Betti}}_\cG(\CX^{\on{discr}}),\CO_{\LocSys_\cG(\CX^{\on{discr}})}),$$

\smallskip

$$\LocSys^{\on{arithm,coarse,discr}}_\cG(X):=\Spec(\Exc^{\on{discr}}).$$

\sssec{}

The algebra $\Exc^{\on{discr}}$ is the algebra of excursion operators attached to the group $\on{Weil}(X,x)^{\on{discr}}$, 
see \cite[Sect. 2.7]{GKRV}; it is denoted in {\it loc. cit.} by $\CEnd_{\CA^{\otimes Y}}(\one_{\CA^{\otimes Y}})$; in our case
$Y=\CX^{\on{discr}}$ and $\CA=\Rep(\cG)$. 

\medskip

The classical commutative algebra $H^0(\Exc^{\on{discr}})$ is the algebra of excursion operators in \cite{VLaf1}, so
$$^{\on{cl}}\!\LocSys^{\on{arithm,coarse,discr}}_\cG(X)$$
is the classical coarse moduli space of representations of $\on{Weil}(X,x)^{\on{discr}}$, considered as a discrete group. 

\sssec{} \label{sss:naive Exc bis}

We have a naturally defined closed embedding
\begin{equation} \label{e:into discrete}
\LocSys^{\on{arithm}} _\cG(X)\hookrightarrow \LocSys^{\on{Betti}}_\cG(\CX^{\on{discr}}),
\end{equation}
which induces a map
$$\Exc^{\on{discr}}\to \Exc,$$
surjective on $H^0$, and hence a closed embedding
$$\LocSys^{\on{arithm,coarse}}_\cG(X)\hookrightarrow \LocSys^{\on{arithm,coarse,discr}}_\cG(X).$$

\ssec{Enhanced trace and Drinfeld's object} \label{ss:enh Tr}

In this subsection we will explain how the procedure of \emph{2-categorical trace} produces from $\Shv_{\on{Nilp}}(\Bun_G)$,
equipped with the Frobenius endofunctor, an object of $\QCoh(\LocSys_\cG^{\on{arithm}}(X))$, to be denoted $\Drinf$. 

\sssec{} \label{sss:enhanced trace}

Recall the set-up of \secref{sss:enhanced trace setup}.

\medskip

Thus, we take $\bA$ to be the (symmetric) monoidal category $$\QCoh(\LocSys_\cG^{\on{restr}}(X)),$$ and we take $F_\bA$ to be
given by $\Frob^*$, where $\Frob$ is as in \secref{sss:Frob on LocSys}. 

\medskip

By \secref{sss:enhanced trace setup Y}, we have a canonical identification 
$$\on{HH}_\bullet(\Frob^*,\QCoh(\LocSys_\cG^{\on{restr}}(X)))
\simeq \QCoh((\LocSys_\cG^{\on{restr}}(X))^{\Frob})=:\QCoh(\LocSys_\cG^{\on{arithm}}(X))$$
as (symmetric) monoidal categories. 

\medskip
%
%

We take the module category $\bM$ to be $\Shv_{\on{Nilp}}(\Bun_G)$, and $F_\bM$ to be $(\Frob_{\Bun_G})_*$,
which is equipped with a natural structure of compatibility with $\Frob^*$. Since $\Shv_{\on{Nilp}}(\Bun_G)$ is
dualizable as a DG category (thanks to \thmref{t:Nilp comp gen abs}), we can consider the objects
$$\Tr((\Frob_{\Bun_G})_*,\Shv_{\on{Nilp}}(\Bun_G))=:\wt{\on{Autom}}\in \Vect_\sfe$$
and 
$$\Tr^{\on{enh}}_{\QCoh(\LocSys_\cG^{\on{restr}}(X))}((\Frob_{\Bun_G})_*,\Shv_{\on{Nilp}}(\Bun_G))=:\Drinf
\in \QCoh(\LocSys_\cG^{\on{arithm}}(X)).$$

In the next subsection we will explain that $\Drinf$ can be regarded as a ``universal shtuka", see
\propref{p:obsetvables}. 

\sssec{}

From \thmref{t:enhanced Tr Y} (combined with \corref{c:LocSys semi-rigid}), we obtain:

\begin{cor} \label{c:Autom tilde as sects}
There exists a canonical isomorphism 
\begin{equation} \label{e:Autom tilde as sects}
\wt{\on{Autom}}\simeq \Gamma(\LocSys_\cG^{\on{arithm}}(X),\Drinf).
\end{equation} 
\end{cor}

\begin{rem}
Note that according to \thmref{t:enhanced Tr Y}, a priori, in the right-hand side in \eqref{e:Autom tilde as sects},
we had to consider the functor $\Gamma_!(\LocSys_\cG^{\on{arithm}}(X),\Drinf)$. However, since 
$\LocSys_\cG^{\on{arithm}}(X)$ is actually an algebraic stack (thanks to \thmref{t:Frob-finite}), we have 
$$\Gamma_!(\LocSys_\cG^{\on{arithm}}(X),-)\simeq \Gamma(\LocSys_\cG^{\on{arithm}}(X),-).$$
\end{rem}

\sssec{}

In particular, from \corref{c:Autom tilde as sects}, 
we obtain an action of the algebra
$$\Exc:=\Gamma(\LocSys_\cG^{\on{arithm}}(X),\CO_{\LocSys_\cG^{\on{arithm}}(X)})$$
on $\wt{\on{Autom}}$. 

\sssec{}

Let us now combine \corref{c:Autom tilde as sects} with \conjref{c:Trace conj}. We obtain:

\begin{corconj} \label{c:loc of Autom}
There exists a canonical isomorphism of vector spaces
\begin{equation} \label{e:Autom as sects}
\on{Autom} \simeq \Gamma(\LocSys_\cG^{\on{arithm}}(X),\Drinf).
\end{equation} 
\end{corconj}

As a consequence, we obtain:

\begin{corconj}
There exists a canonically defined action of the algebra $\Exc$ on $\on{Autom}$.
\end{corconj}

Combining with Sects. \ref{sss:naive Exc}-\ref{sss:naive Exc bis}, we obtain an action of the algebra $\Exc^{\on{discr}}$ on $\on{Autom}$.
Thus, we obtain a spectral decomposition of $\on{Autom}$ over the affine scheme $\LocSys_\cG^{\on{arithm,coarse,discr}}(X)$.

\begin{rem}
As we will see in Remark \ref{r:action by excursions}, if we furthermore input \conjref{c:Trace conj legs} 
(along with its complement in \secref{sss:partial Frobs}), we will see that the resulting action of 
$\Exc^{\on{discr}}$ on $\on{Autom}$ equals the action defined in \cite{VLaf1} 
(with the extension by C.~Xue in \cite{Xue1}) by excursion operators.
\end{rem} 

\sssec{}

We can view the conclusion of \corref{c:loc of Autom} as ``localization" of the space $\on{Autom}$ of automorphic functions
onto the stack $\LocSys_\cG^{\on{arithm}}(X)$ of arithmetic Langlands parameters, in the sense that
we realize $\on{Autom}$ as the space of sections of a quasi-coherent sheaf on this stack.

\ssec{Relation to shtukas} \label{ss:shtukas again}

In this subsection we will explain that the Shtuka Conjecture  (i.e., \conjref{c:Trace conj legs}) implies that the object
$\Drinf$ constructed above, encodes the cohomology of shtukas. 

\sssec{}

Recall the objects 
$$\wt\Sht_{I,V}\in \qLisse(X),$$
see \secref{sss:sht tilde}.

\medskip

We will now show, following \cite[Sect. 5.2]{GKRV}, how the object
$$\Drinf\in \QCoh(\LocSys_\cG^{\on{arithm}}(X))$$
recovers these objects, and endows them with a structure of equivariance with respect to the partial
Frobenius maps. 

\sssec{}

For$I\in \on{fSet}$ and $V\in \Rep(\cG)^{\otimes I}$, let $\CE^I_V$ be the corresponding tautological object of
$$\QCoh(\LocSys_\cG^{\on{restr}}(X))\otimes \qLisse(X)^{\otimes I},$$
see \eqref{e:EV again}. 

\medskip

Namely, for $S\to \LocSys_\cG^{\on{restr}}(X)$, the value of $\CE^I_V$ on $S$, viewed as an object of 
$$\QCoh(S) \otimes \qLisse(X)^{\otimes I}$$ is the value on $V$ of the symmetric monoidal functor
$$\Rep(\cG)^{\otimes I}\to \QCoh(S)^{\otimes I} \otimes \qLisse(X)^{\otimes I} \to \QCoh(S)\otimes \qLisse(X)^{\otimes I},$$
where:

\begin{itemize}

\item The first arrow is the $I$ tensor power of the functor $\Rep(\cG)\to \QCoh(S) \otimes \qLisse(X)$
defining the map $S\to \LocSys_\cG^{\on{restr}}(X)$;

\item The second arrow uses the $I$-fold tensor product functor $\QCoh(S)^{\otimes I}\to \QCoh(S)$.

\end{itemize}

\medskip

In what follows, by a slight abuse of notation, we will denote by the same character $\CE^I_V$ the image
of $\CE^I_V$ under the fully faithful functor\footnote{This functor is fact an equivalence, by the combination
of \corref{c:curve Verdier compat} and \thmref{t:product thm 2 sch}.}
$$\QCoh(\LocSys_\cG^{\on{restr}}(X))\otimes \qLisse(X)^{\otimes I}\to 
\QCoh(\LocSys_\cG^{\on{restr}}(X))\otimes \qLisse(X^I).$$

\sssec{}

Let 
$$\CE^{I,\on{arithm}}_V\in  \QCoh(\LocSys_\cG^{\on{arithm}}(X))\otimes \qLisse(X^I)$$
denote the 
restriction of $\CE^I_V$ along 
$$\LocSys_\cG^{\on{arithm}}(X)\to \LocSys_\cG^{\on{restr}}(X).$$

\sssec{} 

We claim: 

\begin{prop} \label{p:obsetvables}
There exists a canonical isomorphism in $\qLisse(X^I)$
\begin{equation} \label{e:recover shtuka}
\wt\Sht_{I,V}\simeq
\left(\Gamma(\LocSys_\cG^{\on{arithm}}(X),-)\otimes \on{Id}\right)
(\Drinf\otimes \CE^{I,\on{arithm}}_V).
\end{equation}
\end{prop}

\begin{proof}

This is a variant of \cite[Theorem 4.4.4]{GKRV}, combined with \thmref{t:enhanced Tr Y}. 

\end{proof}

\begin{rem} \label{r:partial Frob}

For $I$ as above and $i\in I$, let $\Frob_{i,X^I}$ denote the Frobenius map along the $i$-th factor in $X^I$. By construction,
the object $\CE^{I,\on{arithm}}_V$ 
carries a natural structure of equivariance with respect to these endomorphisms:
$$((\on{Frob}_{X^I,i})^*\otimes \on{Id})(\CE^{I,\on{arithm}}_V) \simeq \CE^{I,\on{arithm}}_V.$$

\medskip

This structure endows the left-hand side in \eqref{e:recover shtuka} with a similar structure. Thus, we obtain a 
structure of equivariance with respect to the partial Frobenius maps on the objects $\wt\Sht_{I,V}$.

\medskip

It follows as in \cite[Proposition 5.3.3]{GKRV} that the resulting structure on $\wt\Sht_{I,V}$ can
be described by explicit excursion operators. 

\end{rem} 

\begin{rem} \label{r:action by excursions} One can use \propref{p:obsetvables} and \cite[Proposition 5.4.3]{GKRV} 
(combined with \thmref{t:enhanced Tr Y} along with its complement in
Remark \ref{r:enhanced Tr Y}), in order to describe the action of $\Exc^{\on{discr}}$ on $\wt{\on{Autom}}$ 
(see \corref{c:Autom tilde as sects} and Sects. \ref{sss:naive Exc}-\ref{sss:naive Exc bis}) by explicit excursion operators.

\medskip

Thus, if we assume \conjref{c:Trace conj legs} (along with its complement in \secref{sss:partial Frobs}), we obtain
that the above action matches under the isomorphism
$$\wt{\on{Autom}}\simeq \on{Autom}$$
with the action of $\Exc^{\on{discr}}$ on $\on{Autom}$, defined in V.~Lafforgue's work (with the extension by C.~Xue in \cite{Xue1}). 

\end{rem}

\begin{rem}
As has been remarked above, we propose our $\LocSys_\cG^{\on{arithm}}(X)$ as a candidate for the stack
sought-for in \cite[Remark 8.5]{VLaf2} and \cite[Sect. 6]{LafZh} (it was denoted $\CS/\cG$ in both these papers). 

\medskip

The space $\CS$ is supposed to be the affine scheme parameterizing homomorphisms from the Weil group 
of $X$ (based at $x$) to $\cG$. So, our proposal for $\CS$ itself is
$$\LocSys_\cG^{\on{arithm}}(X)\underset{\on{pt}/\cG}\times \on{pt}.$$

\medskip

Although in {\it loc.cit.} the space $\CS/\cG$ is only defined heuristically, it is designed so that it carries a collection
of quasi-coherent sheaves $\CE^{I,\on{arithm},\CS}_V$ for $(I,V)$ as above. 

\medskip

The goal of {\it loc.cit.} was to define an object
$$\Drinf^\CS\in \QCoh(\CS/\cG),$$
so that
$$\left(\Gamma(\CS/\cG,-)\otimes \on{Id}\right)
(\Drinf^\CS\otimes \CE^{I,\on{arithm},\CS}_V)\simeq \Sht_{I,V}$$

Thus, assuming \conjref{c:Trace conj legs}, our $\LocSys_\cG^{\on{arithm}}(X)$ with the object
$$\Drinf\in \QCoh(\LocSys_\cG^{\on{arithm}}(X))$$
achieves this goal.

\end{rem}

\ssec{Arithmetic Arthur parameters}

\sssec{}

Fix an $\sfe$-point $\sigma$ of $\LocSys^{\on{arithm}} _\sG(X)$. By \eqref{e:cotan LocSys} and Verdier duality, the tangent space 
$T_\sigma(\LocSys^{\on{restr}} _\sG(X))$ identifies with
$$\on{C}^\cdot(X,\cg_\sigma)[1].$$

Hence, the tangent space $T_\sigma(\LocSys^{\on{arithm}} _\sG(X))$ identifies with
$$\on{Fib}\left(T_\sigma(\LocSys^{\on{restr}}_\sG(X))\overset{\Frob-\on{id}}\longrightarrow T_\sigma(\LocSys^{\on{restr}}_\sG(X))\right)\simeq
\on{Fib}\left(\on{C}^\cdot(X,\cg_\sigma) \overset{\Frob-\on{id}}\longrightarrow \on{C}^\cdot(X,\cg_\sigma)\right)[1],$$
and thus is concentrated in the cohomological degrees $[-1,2]$. 

\medskip

This implies that $\LocSys^{\on{arithm}} _\sG(X)$ has a perfect cotangent complex and is \emph{quasi-quasi-smooth}.
The latter by definition means that it can be smoothly covered by an derived affine scheme, whose cotangent spaces 
are concentrated in the cohomological degrees $[-2,0]$. 

%
%

\sssec{}

We have
$$H^2(T_\sigma(\LocSys^{\on{arithm}} _\sG(X)))\simeq 
\on{coker}\left(H^2(X,\cg_\sigma) \overset{\Frob-\on{id}}\longrightarrow H^2(X,\cg_\sigma)\right).$$

Hence, by Verdier duality
$$H^{-2}(T^*_\sigma(\LocSys^{\on{arithm}} _\sG(X)))\simeq 
\left(H^0(X,\cg_\sigma(1))\right)^{\Frob},$$
where $(1)$ means Tate twist, and where we have identified $\cg$ with its dual using an invariant form.

\medskip

In other words, we can think of elements of $H^{-2}(T^*_\sigma(\LocSys^{\on{arithm}} _\sG(X)))$ as
elements
$$A\in H^0(X,\cg_\sigma)$$
such that
$$\Frob(A)=q\cdot A.$$

We note that such elements $A$ are necessarily nilpotent. 

\sssec{} \label{sss:arithm non-classical}

Note that $\LocSys^{\on{arithm}} _\sG(X)$ does contain points $\sigma$ which admit a non-zero element $A$ as above.
For example, take $\sigma$ to be geometrically trivial and fix an arbitrary non-zero nilpotent element 
$$A\in \cg\simeq H^0(X,\cg_\sigma).$$
Now let the Weil structure be given by the image of $q\in \BG_m$ under
$$\BG_m\to SL_2 \to \cG,$$
where $SL_2\to \cG$ is a Jacobson-Morozov map corresponding to $A$. 

\medskip

This implies that $\LocSys^{\on{arithm}} _\sG(X)$ is \emph{non-classical} and not even eventually coconnective: 
one can show that for a quasi-quasi smooth scheme that is \emph{not} quasi-smooth, its structure sheaf necessarily lives in infinitely 
many cohomological degrees. 

\sssec{}

Let $\CZ$ be a quasi-quasi-smooth algebraic stack. Following a suggestion of D.~Beraldo, one can mimic the construction of
\cite[Sect. 2.3.3]{AG} and produce a classical algebraic stack, denoted $\Sing_2(\CZ)$, whose $\sfe$-points are pairs
$$(z,\xi), \quad z\in \CZ,\,\xi\in H^{-2}(T^*_z(\CZ)).$$

\sssec{}

We will denote:
$$\on{Arth}^{\on{arithm}}(X):=\Sing_2(\LocSys^{\on{arithm}} _\sG(X)),$$
and refer to it as the stack of \emph{arithmetic Arthur parameters}. 

\medskip

Thus, the stack $\on{Arth}^{\on{arithm}}(X)$ projects to $\LocSys^{\on{arithm}} _\sG(X)$, and the fiber
over a given $\sigma$ is the vector space
$$A\in H^0(X,\cg_\sigma),\,\, \Frob(A)=q\cdot A.$$

\begin{rem}
The terminology ``Arthur parameters" is justified as follows: 

\medskip

If $\sigma$ is semi-simple (as a Weil local system), 
then using a Jacobson-Morozov argument, we can identify the \emph{set} 
$$\{A,\,\, \Frob(A)=q\cdot A\}/\on{Ad}(\on{Aut}(\sigma))$$
with the \emph{set} 
$$\{SL_2\to \on{Aut}(\sigma)\}/\on{Ad}(\on{Aut}(\sigma)).$$ 

\medskip

\noindent (Note, however, nilpotent elements have more automorphisms than $SL_2$-triples.)

\end{rem}

\ssec{A digression: categorical trace on $\IndCoh$ on stacks} \label{ss:Tr IndCoh}

In order to formulate our conjecture that expresses the space of automorphic forms 
explicitly in terms of the spectral side, we will need to make a digression and discuss
properties of the categorical trace construction applied to $\IndCoh(\CY)$, where
$\CY$ is a quasi-smooth stack.

\medskip

The material in this subsection was obtained as a result of communications with D.~Beraldo. 

\sssec{}  \label{sss:known cases QCoh}

Let $\CY$ be a quasi-compact algebraic stack equipped with an endomorphism $\phi$. Then according to \cite[Sect. 3.5.3]{GKRV},
we have
\begin{equation} \label{e:Tr QCoh}
\Tr(\phi^*,\QCoh(\CY))\simeq \Gamma(\CY^\phi,\CO_{\CY^\phi}).
\end{equation} 

\sssec{}

Assume now that $\CY$ locally almost of finite type. In this case, along with $\QCoh(\CY)$, we can consider the
category $\IndCoh(\CY)$, and the functor
$$\Upsilon_\CY:\QCoh(\CY)\to \IndCoh(\CY), \quad \CF\mapsto \CF\otimes \omega_\CY.$$

Recall that we have the canonical self-dualities
$$\QCoh(\CY)^\vee \simeq \QCoh(\CY) \text{ and } \IndCoh(\CY)^\vee\simeq \IndCoh(\CY),$$
with pairings given by
$$\CF_1,\CF_2\mapsto \Gamma(\CY,\CF_1\otimes \CF_2) \text{ and }
\CF_1,\CF_2\mapsto \Gamma^\IndCoh(\CY,\CF_1\sotimes \CF_2),$$
respectively. 

\medskip

With respect to these self-dualities, the functor $\Upsilon_\CY$ is the dual of the (tautological) functor
$$\on{un-ren}_\CY:\IndCoh(\CY)\to \QCoh(\CY).$$

Recall also that the functor $\on{un-ren}_\CY$ has the property that for a schematic map $f:\CY_1\to \CY_2$,
the diagram
$$
\CD
\IndCoh(\CY_1) @>{\on{un-ren}_{\CY_1}}>> \QCoh(\CY_1) \\
@V{f^\IndCoh_*}VV @VV{f_*}V \\
\IndCoh(\CY_2) @>{\on{un-ren}_{\CY_2}}>> \QCoh(\CY_2) 
\endCD
$$
commutes. 

\medskip

In particular, 
\begin{equation} \label{e:Gamma IndCoh}
\Gamma^\IndCoh(\CY,-)\simeq \Gamma(\CY,-)\circ \on{un-ren}_\CY.
\end{equation} 

\sssec{} \label{sss:known cases IndCoh}

A parallel computation to \eqref{e:Tr QCoh} shows that
\begin{equation} \label{e:Tr IndCoh}
\Tr(\phi^!,\IndCoh(\CY))\simeq \Gamma^\IndCoh(\CY^\phi,\omega_{\CY^\phi}).
\end{equation} 

\sssec{}

Furthermore, we can place ourselves in the paradigm of \secref{sss:enhanced trace}, and consider
$\QCoh(\CY)$ and $\IndCoh(\CY)$ as module categories over $\QCoh(\CY)$, equipped with compatible
endofunctors. 

\medskip

Thus, we can consider the objects 
$$\Tr^{\on{enh}}_{\QCoh(\CY)}(\phi^*,\QCoh(\CY)) \text{ and } \Tr^{\on{enh}}_{\QCoh(\CY)}(\phi^!,\IndCoh(\CY))$$
in $\QCoh(\CY^\phi)$. 

\medskip

Generalizing the computation \cite[Sect. 3.5.3]{GKRV}, one can show that
\begin{equation} \label{e:enh trace IndCoh}
\Tr^{\on{enh}}_{\QCoh(\CY)}(\phi^*,\QCoh(\CY)) \simeq \CO_{\CY^\phi} \text{ and }
\Tr^{\on{enh}}_{\QCoh(\CY)}(\phi^!,\IndCoh(\CY)) \simeq \on{un-ren}(\omega_{\CY^\phi}),
\end{equation}
as objects of $\QCoh(\CY^\phi)$ (note that the latter isomorphism is compatible with \eqref{e:Tr IndCoh}
via \eqref{e:Gamma IndCoh}). 

\sssec{} \label{sss:phi and N}

Assume now that $\CY$ is quasi-smooth. Let $\CN$ be a conical Zariski-closed subset 
in $\Sing(\CY)$. Assume that the codifferential map
$$\Sing(\phi):\CY\underset{\phi,\CY}\times  \Sing(\CY)\to \Sing(\CY),$$
sends $\CY \underset{\phi,\CY}\times \CN \subset \CY\underset{\phi,\CY}\times  \Sing(\CY)$ to 
$\CN\subset \Sing(\CY)$, 
so that the functor $\phi^!$ sends 
$$\IndCoh_\CN(\CY)\to \IndCoh_\CN(\CY),$$
see \cite[Proposition 7.1.3(a)]{AG}.

\medskip

Then it makes sense to consider
\begin{equation} \label{e:Tr on IndCoh N}
\Tr(\phi^!,\IndCoh_\CN(\CY))\in \Vect_\sfe.
\end{equation}

\medskip

Furthermore, we can regard $\IndCoh_\CN(\CY)$ as a module category over $\QCoh(\CY)$ and consider the object
$$\Tr^{\on{enh}}_{\QCoh(\CY)}(\phi^!,\IndCoh_\CN(\CY))\in \QCoh(\CY^\phi),$$
so that by \eqref{e:enhanced trace abs} we have
$$\Tr(\phi^!,\IndCoh_\CN(\CY))\simeq \Gamma\left(\CY^\phi,\Tr^{\on{enh}}_{\QCoh(\CY)}(\phi^!,\IndCoh_\CN(\CY))\right).$$

\begin{rem}

Unfortunately, we do not have an explicit answer for what the above object $\Tr^{\on{enh}}(\phi^!,\IndCoh_\CN(\CY))$
is in general. We expect, however, that one can give such an answer in terms of the subset 
$$\CN^\phi\subset \Sing_2(\CY^\phi),$$
defined in \eqref{e:N phi} below. 

\medskip

Yet, we know some particular cases: by \eqref{e:enh trace IndCoh}, we have
\begin{equation} \label{e:Tr enh on QCoh}
\Tr^{\on{enh}}_{\QCoh(\CY)}(\phi^!,\IndCoh_{\{0\}}(\CY))\simeq \CO_{\CY^\phi}
\end{equation}
and 
\begin{equation} \label{e:Tr enh on full IndCoh}
\Tr^{\on{enh}}_{\QCoh(\CY)}(\phi^!,\IndCoh(\CY))\simeq \on{un-ren}_{\CY^\phi}(\omega_{\CY^\phi}),
\end{equation}

\end{rem}

\sssec{}

Note now that for $\CY$ quasi-smooth, the stack $\CY^\phi$ is quasi-quasi-smooth and
\begin{equation} \label{e:sing fixed locus}
\Sing_2(\CY^\phi):=\{y\in \CY,\,\, \phi(y)\sim y,\,\, \xi\in H^{-1}(T^*_y(\CY)),\,\, \Sing(\phi)(\xi)=\xi\}.
\end{equation} 

\medskip

Let $\CN\subset \Sing(\CY)$ be as \secref{sss:phi and N}. Set
\begin{equation} \label{e:N phi}
\CN^\phi \subset \Sing_2(\CY^\phi)
\end{equation} 
be the subset that in terms of \eqref{e:sing fixed locus} corresponds to the condition that $\xi\in \CN\underset{\CY}\times \{y\}$. 

\sssec{}

We propose:

\begin{conj} \label{c:ignore N}
Suppose that for a pair of conical subsets $\CN_1\subset \CN_2$ as above, the inclusion
$$\CN_1^\phi \subset \CN_2^\phi$$
is an equality.
Then the inclusion functor
$$\IndCoh_{\CN_1}(\CY) \hookrightarrow \IndCoh_{\CN_2}(\CY)$$
defines an \emph{isomorphism}
$$\Tr^{\on{enh}}_{\QCoh(\CY)}(\phi^!,\IndCoh_{\CN_1}(\CY))\simeq \Tr^{\on{enh}}_{\QCoh(\CY)}(\phi^!,\IndCoh_{\CN_2}(\CY))$$
in $\QCoh(\CY^\phi)$.
\end{conj} 

This conjecture is not far-fetched, and might have been already established in the works of
D.~Beraldo. 

\sssec{}

As a particular case, and combining with \eqref{e:Tr enh on full IndCoh} we obtain: 

\begin{corconj} \label{cor:ignore N}
Suppose that for $\CN$ as above, the inclusion
$$\CN^\phi \subset \Sing_2(\CY^\phi)$$ is an equality.
Then the inclusion functor
$$\IndCoh_\CN(\CY) \hookrightarrow \IndCoh(\CY)$$
defines an \emph{isomorphism}
$$\Tr^{\on{enh}}_{\QCoh(\CY)}(\phi^!,\IndCoh_\CN(\CY))\simeq \Tr^{\on{enh}}_{\QCoh(\CY)}(\phi^!,\IndCoh(\CY))\simeq 
\on{un-ren}_{\CY^\phi}(\omega_{\CY^\phi})$$
in $\QCoh(\CY^\phi)$. In particular,
$$\Tr(\phi^!,\IndCoh_\CN(\CY))\simeq \Gamma^\IndCoh(\CY^\phi, \omega_{\CY^\phi}).$$
\end{corconj} 

\ssec{A digression: categorical trace on $\IndCoh$ on \emph{formal} stacks}  \label{ss:Tr IndCoh formal} 

We now generalize the discussion in \secref{ss:Tr IndCoh} to the case when instead of a quasi-compact algebraic stack $\CY$,
we have a \emph{formal algebraic stack} $\CY$ as in \secref{sss:formal union quot}. I.e., $\CY$ is a disjoint of
\'etale stacks, each of which is the quotient of a formal affine scheme by an action of an algebraic group. 

\sssec{}

First we observe that by Propositions \ref{p:self-duality for semi-pass} and \ref{p:counit self-duality for semi-pass}, 
we have a canonical identification
\begin{equation} \label{e:Tr QCoh formal}
\Tr(\phi^*,\QCoh(\CY))\simeq \Gamma_!(\CY^\phi,\CO_{\CY^\phi}).
\end{equation} 

(Note the difference with formula \eqref{e:Tr QCoh}: for a formal stack we have $\Gamma_!(\CY^\phi,-)$
instead of $\Gamma(\CY^\phi,-)$.)

\sssec{}

By contrast, the formula for $\Tr(\phi^!,\IndCoh(\CY))$ remains unchanged:
\begin{equation} \label{e:Tr IndCoh formal}
\Tr(\phi^!,\IndCoh(\CY))\simeq \Gamma^\IndCoh(\CY^\phi,\omega_{\CY^\phi}),
\end{equation} 
where $\Gamma^\IndCoh$ is (an equivariant version of) the functor in \cite[Chapter 3, Sect. 1.4]{GR2}.

\sssec{}

The functor 
$$\Upsilon_\CY:\QCoh(\CY)\to \IndCoh(\CY)$$
is defined as before:
$$\CF\mapsto \CF\otimes \omega_\CY$$
(in fact, this functor makes sense for \emph{any} prestack locally almost of finite type). 

\medskip

Let 
$$\on{un-ren}_\CY:\IndCoh(\CY)\to \QCoh(\CY)$$
be the functor dual to the functor $\Upsilon_\CY$. 

\medskip

This functor can be characterized by the property that the diagrams
$$
\CD
\IndCoh(S) @>{\on{un-ren}_{S}}>> \QCoh(S) \\
@V{f^\IndCoh_*}VV @VV{f_*}V \\
\IndCoh(\CY) @>{\on{un-ren}_{\CY}}>> \QCoh(\CY) 
\endCD
$$
are commutative for all $S\overset{f}\to \CY$, where $S$ is an affine scheme almost of finite type. 

\sssec{}

Then parallel to \eqref{e:enh trace IndCoh}, we have: 
\begin{equation} \label{e:enh trace IndCoh formal}
\Tr^{\on{enh}}_{\QCoh(\CY)}(\phi^*,\QCoh(\CY)) \simeq \CO_{\CY^\phi} \text{ and }
\Tr^{\on{enh}}_{\QCoh(\CY)}(\phi^!,\IndCoh(\CY)) \simeq \on{un-ren}_{\CY^\phi}(\omega_{\CY^\phi}),
\end{equation}

\sssec{}  \label{sss:Tr IndCoh formal 3}

Finally, we conjecture that a generalization of \conjref{c:ignore N}, stated ``as-is" 
holds in the case of formal stacks as well.

%
%
%
%
%
%
%
%
%
%

\ssec{Towards an explicit spectral description of the space of automorphic functions}

In this subsection we will assume two of our Main Conjectures, \ref{c:Trace conj} and \ref{c:restr GLC}
and (try to) deduce consequences for $\on{Autom}$.

\sssec{}

First, putting the above two conjectures together, we obtain:

\begin{mainconj} \label{c:autom via LocSys}
We have a canonical isomorphism
$$\on{Autom}\simeq \Tr(\Frob^*,\IndCoh_{\Nilp}(\LocSys^{\on{restr}}_\sG(X))).$$ 
\end{mainconj}

Since $\Frob$ is an automorphism of $\LocSys^{\on{restr}}_\sG(X)$, in the above conjecture 
we could replace the functor $\Frob^*$ by $\Frob^!$. 

\medskip

Thus, assuming the above conjecture, in order to describe $\on{Autom}$, we wish to 
have an explicit description of the object
$$\Tr(\Frob^!,\IndCoh_{\Nilp}(\LocSys^{\on{restr}}_\sG(X)))\in \Vect_\sfe.$$

\sssec{}

We apply the discussion in \secref{ss:Tr IndCoh formal} to $\CY=\LocSys^{\on{restr}}_\cG(X)$ with $\phi=\Frob$. We note that the inclusion
$$\Nilp^{\Frob}\hookrightarrow \Sing_2(\LocSys^{\on{arithm}} _\cG(X))=\on{Arth}^{\on{arithm}}(X)$$
is indeed an equality. 

\medskip

Hence, combining \conjref{c:restr GLC} with Corollary-of-Conjecture
\ref{cor:ignore N} (for formal stacks, see \secref{sss:Tr IndCoh formal 3}), we obtain:
\begin{mainconj}
We have a canonical isomorphism in $\QCoh(\LocSys^{\on{arithm}} _\cG(X))$:
$$\Drinf\simeq \on{un-ren}_{\LocSys^{\on{arithm}} _\cG(X)}(\omega_{\LocSys^{\on{arithm}} _\cG(X)}).$$ 
\end{mainconj}

\sssec{}

Taking global sections over $\LocSys^{\on{arithm}} _\cG(X)$, and taking into account Corollary-of-Conjecture
\ref{c:loc of Autom}, we obtain:

\begin{mainconj} \label{c:autom-omega}
We have a canonical isomorphism
$$\on{Autom}\simeq \Gamma^\IndCoh(\LocSys^{\on{arithm}}_\cG(X),\omega_{\LocSys^{\on{arithm}} _\cG(X)}).$$
\end{mainconj}

Note that \conjref{c:autom-omega} provides an explicit description of the space of automorphic functions
in terms of Galois representations. 

\begin{rem} 

Note also that the statement of \conjref{c:autom-omega} is close to the most naive guess for the expression
of $\on{Autom}$ in terms of the moduli space of Galois representations: the latter would say that we should 
take sections of the structure sheaf on $\LocSys^{\on{arithm}}_\cG(X)$, whereas \conjref{c:autom-omega}
says that we should rather take sections of the dualizing complex. 

\medskip

Note, however, that the objects
$$\CO_{\LocSys^{\on{arithm}}_\cG(X)} \text{ and } \on{un-ren}_{\LocSys^{\on{arithm}} _\cG(X)}(\omega_{\LocSys^{\on{arithm}} _\cG(X)})$$
of $\QCoh(\LocSys^{\on{arithm}} _\cG(X))$ are very far apart: 

\medskip

Since $\LocSys^{\on{arithm}} _\cG(X))$ is not quasi-smooth, 
the structure sheaf goes infinitely off in the connective direction, while the dualizing complex goes infinitely off in the coconnective direction. 

\end{rem}

\section{Proofs of Theorems \ref{t:Frob-finite} and \ref{t:irred Weil}} \label{s:arithm}

\ssec{Proof of \thmref{t:Frob-finite}}

In this subsection we will prove \thmref{t:Frob-finite}. We go back to the notations of Part I, and denote
the reductive group for which we are considering local systems by $\sG$. 

\medskip

The key ingredient will be provided by
\thmref{t:coarse restr}, which gives us a handle on ``how far is $\LocSys^{\on{restr}}_\sG(X)$
from being an algebraic stack", combined with some fundamental facts from algebraic geometry
pertaining to Weil sheaves on curves (specifically, Weil-II and L.~Lafforgue's theorem, which says that 
every irreducible Weil local system is pure). 

\sssec{} \label{sss:what to prove Frob finite}

First off, the assertion that $(\LocSys^{\on{restr}}_\sG(X))^{\on{Frob}}$ is mock-affine and locally almost of finite type
as a prestack follows from the corresponding property of $\LocSys^{\on{restr}}_\sG(X)$.

\medskip 

To prove the remaining assertions of the theorem, by \propref{p:repr of fibers}, it suffices to consider the
case of $\sG=GL_n$. 

\medskip

The assertion of the theorem can be broken into two parts:

\medskip

\noindent{(a)} There are only finitely many connected components of $\LocSys^{\on{restr}}_\sG(X)$
that are invariant under $\on{Frob}$.

\medskip

\noindent{(b)} For each connected component $\CZ$ of $\LocSys^{\on{restr}}_\sG(X)$, the fiber
product
$$(\LocSys^{\on{restr}}_\sG(X))^{\on{Frob}} \underset{\LocSys^{\on{restr}}_\sG(X)}\times \CZ$$
is an algebraic stack (as opposed to an ind-algebraic stack, see \secref{ss:ind alg}).

\sssec{}

We start by proving (a). Recall (see \propref{p:conn comps LocSys})
that connected components of $\LocSys^{\on{restr}}_\sG(X)$ are 
in bijection with isomorphism classes of semi-simple local systems. To such a local system 
we can associate a partition 
$$n=(n_1+...+n_1+n_2+...+n_2+....+n_k+...+n_k), \quad n_i\neq n_j$$
with $n_i$ appearing $m_i$ times, and a collection of irreducible local systems 
\begin{equation} \label{e:string loc syst}
(\sigma^1_{n_1},...,\sigma^{m_1}_{n_1},\sigma^1_{n_2},...,\sigma^{m_2}_{n_2},\sigma^1_{n_k},...,\sigma^{m_k}_{n_k}),
\end{equation} 
where each $\sigma^?_{n_i}$ has  rank $n_i$.

\medskip

We claim that there is only a finite number of possibilities for a string \eqref{e:string loc syst}, provided that
its isomorphism class is invariant under the Frobenius. 

\medskip

Indeed, the isomorphism class as above is invariant under the Frobenius if for every $j=1,...,k$ there
exists an element $g_j\in \Sigma_{m_j}$ (the symmetric group on $m_j$ letters) such that
$$(\on{Frob}_X(\sigma^1_{n_j}),...,\on{Frob}_X(\sigma^{m_j}_{n_j}))=
(\sigma^{g_j(1)}_{n_j},...,\sigma^{g_j(m_j)}_{n_j}).$$ 

\medskip

For every $j$ let $d_j:=\on{ord}(g_j)$. We obtain that all local systems $\sigma^?_{n_j}$ are invariant under $(\on{Frob}_X)^{d_j}$.
I.e., each such local system is an irreducible local system (over $\ol\BF_q$) that can be equipped with a Weil structure
(with respect to $\BF_{q^{d_j}}$). 

\medskip

We claim that the number of isomorphism classes of such local systems is finite. 

\medskip

To prove this, it suffices to show that the number of irreducible Weil local systems of a given rank $r$ and a fixed determinant
is finite. The latter follows from L.~Lafforgue's theorem (\cite{LLaf}), which says that such Weil local systems are in
bijection with unramified cuspidal automorphic representations of $GL_r$ with a fixed central character. Now, the 
number of such automorphic representations (for a given function field) is finite. 

\sssec{}

We now start tackling point (b) from \secref{sss:what to prove Frob finite}.

\medskip

Let $\CZ$ be a connected component of $\LocSys^{\on{restr}}_\sG(X)$ invariant under the Frobenius. Denote
$$\CZ^{\on{rigid}_x}:=\CZ\underset{\LocSys^{\on{restr}}_\sG(X)}\times \LocSys^{\on{restr},\on{rigid}_x}_\sG(X).$$

It is enough to show that
\begin{equation}  \label{e:Z Frob}
(\CZ^{\on{Frob}})^{\on{rigid}_x} :=\CZ^{\on{Frob}} \underset{\CZ}\times \CZ^{\on{rigid}_x}
\end{equation}
is an affine scheme; a priori we know that it is an ind-affine ind-scheme. 

\medskip

It follows from \corref{c:def}(a) that $(\CZ^{\on{Frob}})^{\on{rigid}_x}$ has a connective corepresentable deformation theory. 
Therefore, by \cite[Theorem 18.1.0.1]{Lu3}, it suffices to show that $^{\on{cl}}((\CZ^{\on{Frob}})^{\on{rigid}_x})$ is a classical affine scheme. Equivalently,
it suffices to show that the underlying classical prestack of $\CZ^{\on{Frob}}$ itself is a classical algebraic stack (as opposed to an ind-algebraic stack). 

\sssec{}

Set 
$$\CW^{\on{rigid}_x}:=\on{pt}\underset{\CZ^{\on{coarse}}}\times \CZ^{\on{rigid}_x},$$
where $\CZ^{\on{coarse}}$ is as in \thmref{t:coarse restr}, 
and set also
$$\CW:=\CW^{\on{rigid}_x}/\sG\simeq \on{pt}\underset{\CZ^{\on{coarse}}}\times \CZ.$$

We will prove:

\begin{prop} \label{p:non-deform coarse 1}
The map $\CW^{\on{Frob}}\to \CZ^{\on{Frob}}$ induces an isomorphism of the underlying classical prestacks. 
\end{prop} 

This proposition immediately implies that $^{\on{cl}}(\CZ^{\on{Frob}})$ is an algebraic stack, since we know that $\CW$ (and hence
$\CW^{\on{Frob}}$) is an algebraic stack, by \corref{c:coarse restr bis}. 

\sssec{}

Note that on the level of the underlying classical prestacks, the map
$$\on{pt}\to \CZ^{\on{coarse}}$$
is fully faithful (since $\CZ^{\on{coarse}}$ is a derived scheme). 

\medskip

Hence, the assertion of \propref{p:non-deform coarse 1} is equivalent to the following one:

\begin{prop} \label{p:non-deform coarse 2}
The composition
\begin{equation} \label{e:non-deform coarse 2}
\CZ^{\on{Frob}}\to \CZ \overset{\brr}\to \CZ^{\on{coarse}}
\end{equation}
factors as 
\begin{equation} \label{e:non-deform coarse 2 bis}
\CZ^{\on{Frob}}\to \on{pt}\to \CZ^{\on{coarse}}
\end{equation}
at the level of the underlying classical prestacks.
\end{prop}


\ssec{Proof of \propref{p:non-deform coarse 2}: the pure case} \label{ss:pure case}

\propref{p:non-deform coarse 2} says that all global functions on $\CZ$ become constant, when restricted to
$\CZ^{\Frob}$.

\medskip

In this subsection we will prove this assertion on a neighborhood of a point of $\CZ^{\Frob}$
that corresponds to a \emph{pure} local system. The proof will use Weil-II. 

\sssec{} \label{sss:factor Artinian}

Let 
\begin{equation} \label{e:deform Weil}
S\to \CZ^{\on{Frob}}
\end{equation} 
be a map, where $S=\Spec(A)$ with $A$ classical Artinian. 

\medskip

It suffices to show that for any such map, the composition
\begin{equation} \label{e:S factor 1}
S\to \CZ^{\on{Frob}}\to \CZ \overset{\brr}\to \CZ^{\on{coarse}}
\end{equation} 
factors as
\begin{equation} \label{e:S factor 2}
S\to \on{pt}\to \CZ^{\on{coarse}}.
\end{equation} 

\medskip

We can think of \eqref{e:deform Weil} as a local system $E_A$ on $X$, endowed with a Weil structure, 
and equipped with an action of $A$, whose fiber at $x\in X$ is a (locally) free $A$-module of rank $n$. 

\sssec{}

Let $E$ be the Weil local system corresponding to the composition
$$\on{pt}\to S\to \CZ^{\on{Frob}}.$$

Let us first consider the case when $E$ is \emph{pure of weight} $0$ (with respect to some
identification $\ol\BQ_\ell\simeq \BC$). 

\medskip 

Let $\ol{E}$ denote the underlying local system, when we forget the Weil structure. 
Let $\on{Aut}(\ol{E})$ denote the \emph{classical} algebraic group of automorphisms of $\ol{E}$.

\medskip

Varying the Weil structure on $\ol{E}$ defines a map
\begin{equation} \label{e:vary Weil}
\on{Aut}(\ol{E})/\on{Ad}_{\on{Frob}}(\on{Aut}(\ol{E}))\to \CZ^{\on{Frob}},
\end{equation}
where $\on{Ad}_{\on{Frob}}(\on{Aut}(\ol{E}))$ stands for the action of $\on{Aut}(\ol{E})$ on itself given by
$$g(g_1)=\on{Frob}(g)\cdot g_1\cdot g^{-1},$$
and where $\on{Frob}$ is the automorphism of $\on{Aut}(\ol{E})$ induced by 
$$\on{Aut}(\ol{E}) \overset{\on{functoriality}}\longrightarrow \on{Aut}(\on{Frob}(\ol{E})) 
\overset{\text{Weil structure}}\longrightarrow \on{Aut}(\ol{E}).$$

\sssec{}

We claim that the map \eqref{e:vary Weil} defines an isomorphism at the level
of formal completions at $E$. In order to prove this, since both sides admit deformation
theory, it suffices to show that the map \eqref{e:vary Weil} induces an isomorphism 
at the level of tangent spaces. 

\medskip

We have:
$$T_1(\on{Aut}(\ol{E}))\simeq H^0(X,\End(\ol{E})),$$
and 
$$T_1(\on{Aut}(\ol{E})/\on{Ad}_{\on{Frob}}(\on{Aut}(\ol{E}))) \simeq 
\on{coFib}\left(H^0(X,\End(\ol{E})) \overset{\on{Frob}-\on{Id}}\longrightarrow H^0(X,\End(\ol{E}))\right).$$

\medskip

We also have
$$T_E(\CZ^{\on{Frob}})\simeq \on{Fib}\left(T_{\ol{E}}(\CZ) \overset{\on{Frob}-\on{Id}}\longrightarrow  T_{\ol{E}}(\CZ)\right),$$
where
$$T_{\ol{E}}(\CZ)\simeq \on{C}^\cdot(X,\End(\ol{E}))[1].$$

The map that \eqref{e:vary Weil} induces at the level of tangent spaces corresponds to canonical map 
$$H^0(X,\End(\ol{E}))\to \on{C}^\cdot(X,\End(\ol{E})).$$

Hence, in order to show that
$$T_1(\on{Aut}(\ol{E})/\on{Ad}_{\on{Frob}}(\on{Aut}(\ol{E}))) \to T_E(\CZ^{\on{Frob}})$$
is an isomorphism, it suffices to show that $\on{Frob}-\on{Id}$ induces an isomorphism on $H^1(X,\End(\ol{E}))$ and $H^2(X,\End(\ol{E}))$.
In other words, we have to show that $\on{Frob}$ does not have eigenvalue $1$ on either $H^1(X,\End(\ol{E}))$ or $H^2(X,\End(\ol{E}))$. 

\sssec{}

We will now use the assumption that $E$ is pure of weight $0$. 

\medskip

This assumption implies that the induced Weil structure on
$\End(\ol{E})$ is also pure of weight $0$. Hence, the eigenvalues of $\on{Frob}$ on $H^1(X,\End(\ol{E}))$ (resp., $H^2(X,\End(\ol{E}))$)
are algebraic numbers that under any complex embedding have Archimedean absolute values $q^{\frac{1}{2}}$ (resp., $q$).

\medskip

In particular, these eigenvalues are different from $1$. 

\sssec{}

Thus, we have established that the map \eqref{e:vary Weil} is a formal isomorphism at $E$. Hence, by deformation
theory, the initial map
$$S\to \CZ^{\on{Frob}}$$
of \eqref{e:deform Weil} can be lifted to a map
$$S\to \on{Aut}(\ol{E})/\on{Ad}_{\on{Frob}}(\on{Aut}(\ol{E})).$$

\medskip

However, the composite map
$$\on{Aut}(\ol{E})/\on{Ad}_{\on{Frob}}(\on{Aut}(\ol{E}))\to \CZ^{\on{Frob}}\to \CZ$$
by definition factors as
$$\on{Aut}(\ol{E})/\on{Ad}_{\on{Frob}}(\on{Aut}(\ol{E}))\to \on{pt}/\on{Aut}(\ol{E}) \to \CZ,$$
while the composition
$$\on{pt}/\on{Aut}(\ol{E}) \to \CZ \to \CZ^{\on{coarse}}$$
factors as
$$\on{pt}/\on{Aut}(\ol{E}) \to \on{pt}\to \CZ^{\on{coarse}}.$$

\medskip

This proves the required factorization of  \eqref{e:S factor 1} as \eqref{e:S factor 2} (in the case when $E$ was pure of weight $0$). 

\ssec{Reduction to the pure case}

Above we have established the factorization of \eqref{e:S factor 1} as \eqref{e:S factor 2} when the Weil local system 
$E$ was pure of weight $0$. 

\medskip

In this subsection we will perform the reduction to this case. The source of pure local systems will be provided by
the theorem of L.~Lafforgue, which says that every irreducible Weil local system is pure. 

\sssec{} \label{sss:weight filtration}

Let us choose an isomorphism
\begin{equation} \label{e:Arch}
\ol\BQ_\ell\simeq \BC,
\end{equation}
so we can assign the Archimedean absolute value to elements of $\ol\BQ_\ell$. With this choice, we claim that 
every Weil local system $E'$ on $X$ acquires a \emph{canonical} (weight) filtration, indexed by real numbers
$$0\subset ...\subset E'_{r_1} \subset E'_{r_2} \subset ... \subset E',$$
such that each subquotient
$$\on{gr}_r(E')$$
is ``pure of weight $r$" in the sense that it is of the form
\begin{equation} \label{e:shape of pure}
E_0\otimes \ell_r,
\end{equation} 
where:

\begin{itemize}

\item $E_0$ is pure of weight $0$ (with respect to \eqref{e:Arch});

\item $\ell_r$ is a character of $\BZ=\on{Weil}(\ol\BF_q/\BF_q)$, on which the generator $1\in \BZ$ acts by a scalar
with Archimedean absolute value $q^{\frac{r}{2}}$. 

\end{itemize}

Moreover, this filtration is functorial in $E'$ and is compatible with tensor
products.

\medskip

The existence and properties of such a filtration follow from the combination of the following three facts:

\begin{itemize}

\item For two local systems $E^1_0\otimes \ell_{r^1}$ and $E^2_0\otimes \ell_{r^2}$ of the form \eqref{e:shape of pure},
$$r_1\neq r_2 \,\Rightarrow\, \Hom(E^1_0\otimes \ell_{r^1},E^2_0\otimes \ell_{r^2})=0.$$

\medskip

\item For two local systems $E^1_0\otimes \ell_{r^1}$ and $E^2_0\otimes \ell_{r^2}$ of the form \eqref{e:shape of pure},
$$r_1< r_2 \,\Rightarrow\, \on{Ext}^1(E^1_0\otimes \ell_{r^1},E^2_0\otimes \ell_{r^2})=0.$$
This follows from \cite{De}.

\medskip

\item Every irreducible Weil local system on $X$ is of the form \eqref{e:shape of pure}. This is a theorem
of L.~Lafforgue, \cite{LLaf}.

\end{itemize}

\sssec{}

Applying this construction to $E'=E_A$ (see \secref{sss:factor Artinian}), we obtain a filtration
$$0=E_{A,0}\subset E_{A,1} \subset ...\subset E_{A,k}=E_A$$
by Weil local systems, stable under the action of $A$.

\medskip

We claim that the fibers of $\on{gr}_i(E_A)$ at any $x\in X$ are (locally) free over $A$. For that end, it suffices to 
show that the induced filtration on $\on{ev}_x(E_A)$ canonically splits.

\medskip

Indeed, let $d$ be such that $x$ is defined over $\ol\BF_{q^d}$. Then $\on{Frob}^d$
acts on $\on{ev}_x(E_A)$, and its action on the different subquotients 
$$\on{gr}_i(E_A)$$ 
has distinct generalized eigenvalues. 

\sssec{}

Thus, we obtain that we can lift our initial $S$-point of $\CZ^{\on{Frob}}$ to a point of
$$(\LocSys_\sP^{\on{restr}}(X))^{\on{Frob}},$$
where $\sP$ is the parabolic corresponding to the ranks of $\on{gr}_i(E_A)$. 

\medskip

Let $\CZ_\sP$ denote the corresponding connected component of 
$\LocSys^{\on{restr}}_\sP(X)$, i.e., we now have a map
\begin{equation} \label{e:lift to P}
S\to (\CZ_\sP)^{\on{Frob}}.
\end{equation}

\medskip

It suffices to show that the composition
\begin{equation} \label{e:S factor 1 P}
S \overset{\text{\eqref{e:lift to P}}}\longrightarrow (\CZ_\sP)^{\on{Frob}}\to \CZ_\sP\to \CZ \to \CZ^{\on{coarse}}
\end{equation}
factors as 
\begin{equation} \label{e:S factor 2 P}
S\to \on{pt}\to \CZ^{\on{coarse}}.
\end{equation}

\sssec{}

In what follows we will want to consider the coarse moduli space corresponding to $\CZ_\sP$.
The slight inconvenience is that $\CZ_\sP$ is \emph{not} ind mock-affine (because $\sP$ is
not reductive). We will overcome this as follows.

\medskip

Set $$\CZ_\sP^{\on{unip-rigid}_x}:=\CZ_\sP\underset{\on{pt}/\sP}\times \on{pt}/\sM,$$
which we can also think of as
$$\CZ_\sP^{\on{rigid}_x}/\sM$$
for a choice of a Levi splitting $\sM\to \sP$.

\medskip

The map
$$\on{pt}/\sM\to \on{pt}/\sP$$
is smooth, so the map
$$S \overset{\text{\eqref{e:lift to P}}}\longrightarrow (\CZ_\sP)^{\on{Frob}}$$
can be lifted to a map
$$S\to (\CZ_\sP)^{\on{Frob}} \underset{\CZ_\sP}\times \CZ_\sP^{\on{unip-rigid}_x}.$$

It suffices to show that the composition
\begin{equation} \label{e:S factor 1 P rigid}
S\to (\CZ_\sP)^{\on{Frob}} \underset{\CZ_\sP}\times \CZ_\sP^{\on{unip-rigid}_x}\to  \CZ_\sP^{\on{unip-rigid}_x}\to  \CZ_\sP\to \CZ \to \CZ^{\on{coarse}}
\end{equation} 
factors as 
\begin{equation} \label{e:S factor 2 P rigid}
S\to \on{pt}\to \CZ^{\on{coarse}}.
\end{equation}

\sssec{}

Since 
$$\CZ^{\on{rigid}_x}_\sP:=\CZ_\sP\underset{\on{pt}/\sP}\times \on{pt}$$ 
is ind-affine ind-scheme, and $\sM$
is reductive, the ind-algebraic stack $\CZ_\sP^{\on{unip-rigid}_x}$ is ind mock-affine. Hence, we have the well-defined
ind-affine ind-scheme
$$\CZ_\sP^{\on{unip-rigid}_x,\on{coarse}},$$
and by construction, any map
$$\CZ_\sP^{\on{unip-rigid}_x}\to  \CU,$$
where $\CU$ is an ind-affine ind-scheme, factors as
$$\CZ_\sP^{\on{unip-rigid}_x}\to \CZ_\sP^{\on{unip-rigid}_x,\on{coarse}}\to \CU.$$

\sssec{}

In particular, the map
$$\CZ_\sP^{\on{unip-rigid}_x}\to  \CZ_\sP\to \CZ \to \CZ^{\on{coarse}}$$
that appears in \eqref{e:S factor 1 P rigid} factors as
$$\CZ_\sP^{\on{unip-rigid}_x}\to \CZ_\sP^{\on{unip-rigid}_x,\on{coarse}}\to \CZ^{\on{coarse}}.$$

\medskip

Hence, it suffices to show that the composition
\begin{equation} \label{e:S factor 1 P rigid bis}
S\to (\CZ_\sP)^{\on{Frob}} \underset{\CZ_\sP}\times \CZ_\sP^{\on{unip-rigid}_x}\to \CZ_\sP^{\on{unip-rigid}_x} \to
\CZ_\sP^{\on{unip-rigid}_x,\on{coarse}}
\end{equation}
factors as 
\begin{equation} \label{e:S factor 2 P rigid bis}
S\to \on{pt}\to \CZ_\sP^{\on{unip-rigid}_x,\on{coarse}},
\end{equation}
where 
$$\on{pt}\to \CZ_\sP^{\on{unip-rigid}_x,\on{coarse}}$$
is the unique $\sfe$-point of $\CZ_\sP^{\on{unip-rigid}_x,\on{coarse}}$, see isomorphism \eqref{e:P vs M} below. 

\sssec{}

Let $\CZ_\sM$ be the connected component of $\LocSys^{\on{restr}}_\sM(X)$, corresponding to $\CZ_\sP$. By the argument in 
\secref{sss:contr N}, the projection
$$\CZ_\sP^{\on{unip-rigid}_x}\to \CZ_\sM$$
induces an isomorphism
\begin{equation} \label{e:P vs M}
\CZ_\sP^{\on{unip-rigid}_x,\on{coarse}}\simeq \CZ_\sM^{\on{coarse}}.
\end{equation} 

Hence, it suffices to show that the composition
$$S \overset{\text{\eqref{e:lift to P}}}\longrightarrow (\CZ_\sP)^{\on{Frob}}\to (\CZ_\sM)^{\on{Frob}}\to \CZ_\sM \to \CZ_\sM^{\on{coarse}}$$
factors as 
$$S \to \on{pt}\to \CZ_\sM^{\on{coarse}}.$$

\sssec{}

Write 
$$\sM=\underset{i}\Pi\, GL_{n_i},$$
so it is enough to prove the corresponding factorization assertion for each of the $GL_{n_i}$ factors separately. 

\medskip

However, by the assumption on $\on{gr}_i(E_A)$, this reduces us to the pure of weight $0$ case considered in \secref{ss:pure case}. 
Indeed, the resulting local systems $\on{gr}_i(E_A)$ are pure of weight $0$ (up to a twist by a line). 

\qed[\thmref{t:Frob-finite}]

\ssec{Proof of \thmref{t:irred Weil}}

\sssec{}

First, we claim: 

\begin{lem} \label{l:irred arithm}
Let $\sigma\in (\LocSys^{\on{restr}}_\sG(X))^{\on{Frob}}(\sfe)$ be irreducible. Then the map
$$\on{pt}/\on{Aut}(\sigma)\to (\LocSys^{\on{restr}}_\sG(X))^{\on{Frob}}$$
is a closed embedding. 
\end{lem}

\begin{proof}

Using the closed embedding 
$$\LocSys^{\on{arithm}} _\cG(X)\hookrightarrow \LocSys^{\on{Betti}}_\cG(\CX^{\on{discr}})$$
of \eqref{e:into discrete}, suffices to show that the resulting composite map 
$$\on{pt}/\on{Aut}(\sigma)\to \LocSys^{\on{Betti}}_\cG(\CX^{\on{discr}})$$
is a closed embedding. 

\medskip

However, this follows from \propref{p:irred closed Betti}.

\end{proof}

\sssec{}

Given \lemref{l:irred arithm}, to prove \thmref{t:irred Weil}, it suffices to show that for an irreducible Weil local system $\sigma$, the tangent space
$$T_\sigma\left((\LocSys^{\on{restr}}_\sG(X))^{\on{Frob}}\right)=0.$$

We have:
$$T_\sigma\left((\LocSys^{\on{restr}}_\sG(X))^{\on{Frob}}\right)\simeq
\on{Fib}(T_\sigma\left(\LocSys^{\on{restr}}_\sG(X)\right)\overset{\on{Frob}-\on{Id}}\longrightarrow T_\sigma\left(\LocSys^{\on{restr}}_\sG(X)\right),$$
while
$$T_\sigma\left(\LocSys^{\on{restr}}_\sG(X)\right)\simeq \on{C}^\cdot(X,\sg_\sigma)[1].$$

So, it is enough to show that $\on{Frob}$ does not have fixed vectors when acting on $H^i(X,\sg_\sigma)$, $i=0,1,2$.

\sssec{}

We first consider the case of $i=0$.

\medskip

Note that 
$$(H^0(X,\sg_\sigma))^{\on{Frob}}$$
is the Lie algebra of the \emph{classical} group of automorphisms of $\sigma$ as a Weil local system.

\medskip

If this group has a non-trivial connected component, a standard argument implies that $\sigma$ can be reduced 
to a proper parabolic.

\sssec{}

To treat the cases $i=1,2$, it suffices to prove the following:

\begin{prop} \label{p:weights for G}
Let $\sG$ be a semi-simple group and let $\sigma$ be an irreducible Weil $\sG$-local system.
Then for any $V\in \Rep(\sG)^{c,\heartsuit}$, the associated Weil local system $V_\sigma$ is
pure of weight $0$.
\end{prop}

This proposition is likely well-known. We will provide a proof for completeness.

\sssec{} \label{p:tann filtr}

First, we recall the following general construction:

\medskip

Let 
$$\on{Fil}_\BR(\Vect_\sfe^{c,\heartsuit})$$ 
be the abelian symmetric monoidal category consisting of finite-dimensional vector spaces,
endowed with a filtration indexed by the poset of real numbers.

\medskip

Let us be given a $\sG$-torsor $\sigma$, thought of as a symmetric monoidal functor 
$$\sF_\sigma:\Rep(\sG)^{c,\heartsuit} \to \Vect_\sfe^{c,\heartsuit}.$$

\medskip

Note that the datum of a lift of $\sF$ 
to a symmetric monoidal functor
$$\sF_\sigma^{\on{Fil}_\BR}: \Rep(\sG)^{c,\heartsuit} \to \on{Fil}_\BR(\Vect_\sfe^{c,\heartsuit})$$ 
is equivalent to the datum of a reduction $\sigma_\sP$ 
of $\sigma$ to a parabolic $\sP$ (denote its Levi quotient by $\sM$)
and an element $\lambda\in \pi_{1,\on{alg}}(Z^0_\sM)\otimes \BR$, 
which is dominant and $(\sG,\sM)$-regular (see \secref{sss:specify M} for what this means).

\sssec{}

With respect to this bijection, for $V\in \Rep(\sG)^{c,\heartsuit}$, 
the filtration on $\sF_\sigma(V)$ is recovered as follows:

\medskip

The choice of $\sP$ defines a filtration on $V$ by $\sP$-subrepresentations, indexed by the poset
of characters of $Z^0_\sM$, 
$$V_\chi \subset V, \quad \chi\in \Hom(Z^0_\sM,\BG_m).$$
such that for a given character $\chi$, the action of $\sP$ on the subquotient $\on{gr}_\chi(V)$
factors through $\sM$ with $Z^0_\sM$ acting by $\chi$. Denote by
$$\sF_{\sigma_\sP}(V_\chi)\subset \sF_\sigma(V)$$
the induced filtration on $\sF_\sigma(V)$.

\medskip

Now for $r\in \BR$, the subspace 
$$(\sF_\sigma(V))_r\subset \sF_\sigma(V)$$ 
is the sum of the subspaces
$$\sF_{\sigma_\sP}(V_\chi), \quad \langle \chi,\lambda\rangle \leq r.$$

\medskip

In particular, if $\sG$ is semi-simple, then $\sP=\sG$ if and only the lift $\sF_\sigma^{\on{Fil}_\BR}$ is trivial, i.e., for every $V\in \Rep(\sG)^{c,\heartsuit}$
\begin{equation} \label{e:when pure}
\on{gr}_r(\sF_\sigma(V))=0 \text{ for } r\neq 0.
\end{equation}

\medskip

This construction is functorial. In particular, if a group acts on $\sF_\sigma$ in a way preserving its lift to
$\sF_\sigma^{\on{Fil}_\BR}$, then the action of this group on $\sigma$ is induced by its
action on $\sigma_\sP$. 

\begin{proof}[Proof of \propref{p:weights for G}]

Let $\on{Gal}(X,x)^W$ be the Tannakian pro-algebraic group corresponding to the 
(abelian) symmetric monoidal category of Weil local systems on $X$, equipped 
with the fiber functor given by $\on{ev}_x$.

\medskip

The datum of $\sigma$ can viewed as a $\sG$-torsor, 
acted on by $\on{Gal}(X,x)^W$.

\medskip

Recall the setting of \secref{sss:weight filtration}. We obtain that the canonical weight filtration on the Weil local systems $V_\sigma$ 
defines a reduction of $\sigma$ to a parabolic $\sP$. Since $\sigma$ was assumed irreducible, we obtain that $\sP=\sG$. By \eqref{e:when pure}, 
this implies
$$\on{gr}_r(V_\sigma)=0 \text{ for } r\neq 0.$$

I.e., all $V_\sigma$ are pure of weight $0$ as required.

\end{proof}

\qed[\thmref{t:irred Weil}]

\newpage

\appendix

\section{Formal affine schemes vs ind-schemes} \label{s:formal}

In this section we will supply the proof of \thmref{t:formal}. 

\ssec{Creating the ring}

In this subsection we will state (a particular case of) \cite[Theorem 18.2.3.2]{Lu3} and deduce
from it our \thmref{t:formal}.

\sssec{}

Ley $\CY$ be an ind-affine ind-scheme. Write it as
$$\CY\simeq \underset{\alpha\in A}{\on{colim}}\, Y_\alpha,$$
where:

\begin{itemize}

\item $A$ is a filtered index category;

\item $Y_\alpha=\Spec(R_\alpha)$'s are derived affine schemes almost of finite type;

\item The transition maps $Y_\alpha\to Y_\beta$ are closed embeddings, i.e., the corresponding maps
$R_\beta\to R_\alpha$ induce surjective maps $H^0(R_\beta)\twoheadrightarrow H^0(R_\alpha)$.

\end{itemize}

\medskip

We can form a commutative ring
$$R:=\underset{\alpha\in A}{\on{lim}}\, R_\alpha.$$

However, in general, we would not be able to say much about this $R$; in particular, we do not
know that it is connective. 

\sssec{}

Assume now that $\CY$ is as in \thmref{t:formal}, i.e.,

\begin{itemize}

\item $^{\on{red}}\CY$ is an affine scheme (to be denoted $Y_{\on{red}}=\Spec(R_{\on{red}})$)\footnote{Note that $Y_{\on{red}}={}^{\on{red}}Y_\alpha$
for $\alpha$ large, but it is \emph{not} $^{\on{red}}\!\Spec(R)$.};

\item $\CY$ admits a corepresentable deformation theory, i.e., 
for any $(S,y)\in \affSch_{/\CY}$, the cotangent space $T^*_y(\CY)\in \on{Pro}(\QCoh(S)^{\leq 0})$
actually belongs to $\QCoh(S)^{\leq 0}$.

\end{itemize} 

\medskip

In this case we claim:

\begin{thm} \label{t:Lurie}  \hfill

\smallskip

\noindent{\em(a)} The ring $R$ is connective. Furthermore, for every $n$, the natural map
$$\tau^{\geq -n}(R) \to \underset{\alpha\in A}{\on{lim}}\, \tau^{\geq -n}(R_\alpha)$$
is an isomorphism. 

\smallskip

\noindent{\em(b)} The ideal $I:=\on{ker}(H^0(R)\to {}R_{\on{red}})$ is finitely generated.

\smallskip

\noindent{\em(c)} The map
$$\CY\to \Spec(R)^\wedge_{\Spec({}R_{\on{red}})}$$
is an isomorphism.

\end{thm}

Let us remind the notation in the above formula: for a prestack $\CW$ and a classical reduced prestack $\CW^0\to {}^{\on{red}}\CW$, we denote by
$\CW^\wedge_{\CW^0}$
the completion of $\CW$ along $\CW_0$, i.e., 
$$\Maps(S,\CW^\wedge_{\CW^0})=\Maps(S,\CW) \underset{\Maps({}^{\on{red}}\!S,{}^{\on{red}}\CW)}\times \Maps({}^{\on{red}}\!S,\CW^0).$$

\begin{rem}
In the course of the proof, we will show that the ring $H^0(R)$ is Noetherian, and each 
$H^n(R)$ is finitely-generated as $H^0(R)$-module.
\end{rem}

\sssec{}

The assertion of \thmref{t:Lurie} implies that of \thmref{t:formal}. Indeed, the possibility to write $\CY$
as a colimit \eqref{e:presentation A} is the content of \cite[Proposition 6.7.4]{GR3}.

\medskip

The rest of this section is devoted to the proof of \thmref{t:Lurie}. 

\ssec{Analysis of the classical truncation} \label{ss:cl formal}

\sssec{}

With no restriction of generality, we can replace $A$ be a cofinal subcategory consisting of indices $\alpha$ for which
$$^{\on{red}}Y_\alpha\to \Spec(R_{\on{red}})$$
is an isomorphism.

\medskip

For an index $\alpha\in A$, let $I_\alpha$ denote the ideal
$$\on{ker}(H^0(R_\alpha)\to {}R_{\on{red}}).$$

For an integer $n$, we can consider its $n$-th power $I_\alpha^n\subset H^0(R_\alpha)$. 
We claim:

\begin{prop} \label{p:classical stabilization}
For every $n$, the $A$-family
$$\alpha\mapsto H^0(R_\alpha)/I_\alpha^n$$
stabilizes.
\end{prop}

\begin{proof}

It is clear that the assertion of the proposition for a given $n$ implies it for all $n'\leq n$. Hence, it
is enough to prove it for integers $n$ of the form $2^m$. The proof proceeds by induction on $m$.
We first consider the base case $m=1$, so $n=2$. 

\medskip

Thus, we wish to show that the family
$$\alpha \mapsto I_\alpha/I_\alpha^2$$
stabilizes. 

\medskip

For every $\alpha$, consider 
$$\on{Fib}(T^*(Y_\alpha)|_{Y_{\on{red}}}\to T^*(Y_{\on{red}}))\in \QCoh(Y_{\on{red}})^{\leq 0}.$$

By the assumption on $\CY$, the inverse system
$$\alpha\mapsto \on{Fib}(T^*(Y_\alpha)|_{Y_{\on{red}}}\to T^*(Y_{\on{red}}))$$
is equivalent to a constant object of $\QCoh(Y_{\on{red}})^{\leq 0}$. 

\medskip

Hence, the inverse system
$$\alpha\mapsto H^0\left(\on{Fib}(T^*(Y_\alpha)|_{Y_{\on{red}}}\to T^*(Y_{\on{red}}))\right)$$
is equivalent to a constant object of 
$\QCoh(Y_{\on{red}})^\heartsuit$.

\medskip

However
$$H^0\left(\on{Fib}(T^*(Y_\alpha)|_{Y_{\on{red}}}\to T^*(Y_{\on{red}}))\right)\simeq I_\alpha/I_\alpha^2$$
(see, e.g., \cite[Chapter 1, Lemma 5.4.3(a)]{GR2}), and the transition maps
$$I_\beta/I_\beta^2\to I_\alpha/I_\alpha^2$$
are surjective. 

\medskip

This implies the stabilization assertion for $n=2$. The induction step is carried 
out by the same argument: 

\medskip

Assume that the assertion holds for $n\leq 2^{m-1}$. Let $R_{n,\on{cl}}$ denote the resulting ring
(the eventual value of $H^0(R_\alpha)/I_\alpha^n$). Since $A$ is filtered, we can assume that $\Spec(R_{n,\on{cl}})$ 
maps to all the $Y_\alpha$. Now, in order to show that the assertion of the proposition holds for $n=2^m$,
we run the above argument with $Y_{\on{red}}$ replaced by 
$\Spec(R_{2^{m-1},\on{cl}})$. 

\end{proof} 

\sssec{}

Let $R_{n,\on{cl}}$ as above, i.e., the eventual value of $H^0(R_\alpha)/I_\alpha^n$.
Let $J_n:=\on{ker}(R_{n,\on{cl}}\to R_{\on{red}})$. 
By construction, for $m\leq n$, we have
$$R_{n,\on{cl}}/(J_n)^{m}\simeq R_{m,\on{cl}}.$$

\medskip

Set $$R_{\on{cl}}:=\underset{n}{\on{lim}}\, R_{n,\on{cl}}$$
Let $J$ denote the ideal $\on{ker}(R_{\on{cl}}\to R_{\on{red}})$. 

\medskip

By the almost of finite type assumption, the ideal
$$J_2\subset R_{2,\on{cl}}$$
is finitely generated; choose generators $\ol{f}_1,...,\ol{f}_m$. Let $f_1,...,f_m$ be their lifts to elements of $J$. 

\medskip

The following is a standard convergence argument:

\begin{lem} \label{l:adic} \hfill

\smallskip

\noindent{\em(a)} The elements $f_1,...,f_m$ generate $J$. 

\smallskip

\noindent{\em(b)} For any $n$, the ideal $J^n \subset R_{\on{cl}}$ is closed in the
$J$-adic topology on $R_{\on{cl}}$, and the inclusion
$$J^n \subset \on{ker}(R_{\on{cl}}\to {}R_{n,\on{cl}})$$
is an equality. 
\end{lem}

\begin{cor} \label{c:R is Noeth}
The ring $R_{\on{cl}}$ is Noetherian. 
\end{cor}

\begin{proof}

Let $\ol{g}_1,...,\ol{g}_n$ be generators of $R_{\on{red}}$ over $\sfe$. Let
$g_1,...,g_n$ be their lifts to elements of $R_{\on{cl}}$.

\medskip

Consider the algebra
$$\sfe[s_1,...,s_n][\![t_1,...,t_m]\!],$$
equipped with a map to $R_{\on{cl}}$ given by
$$s_i\mapsto g_i, \,\, t_j\mapsto f_j.$$

It is easy to see that this map is surjective. Hence, $R_{\on{cl}}$ is Noetherian, 
since $\sfe[s_1,...,s_n][\![t_1,...,t_m]\!]$ is such.

\end{proof}

\sssec{}

Let $^{\on{cl}}\CY$ denote the classical truncation of $\CY$. I.e., $^{\on{cl}}\CY$, viewed as a prestack
is the left Kan extension of the restriction of $\CY$ to classical affine schemes.
Explicitly,
$$^{\on{cl}}\CY=\underset{\alpha}{\on{colim}}\, {}^{\on{cl}}\!Y_\alpha.$$

\begin{rem}

Note that it is not a priori clear that $^{\on{cl}}\CY$ as defined above is an ind-scheme: 
we do not know that it is convergent (see \cite[Chapter 2, Sect. 1.4]{GR1} for what this means).
We do know, however, that its convergent
completion $^{\on{conv}}({}^{\on{cl}}\CY)$ is an ind-scheme
\cite[Chapter 2, Corollary 1.4.4]{GR2}.

\medskip

That said, \lemref{l:comp formal compl} below will imply that in our particular case, $^{\on{cl}}\CY$
is an ind-scheme. 

\end{rem}

\sssec{}

By construction, we obtain a map  
\begin{equation} \label{e:comp formal compl}
^{\on{cl}}\CY \to \Spec(R_{\on{cl}})^\wedge_{Y_{\on{red}}}. 
\end{equation}

\begin{lem} \label{l:comp formal compl}
The map \eqref{e:comp formal compl} is an isomorphism. 
\end{lem}

\begin{proof}

The fact that the map in question is an isomorphism when evaluated on classical affine schemes
follows from the construction. 

\medskip

Hence, it remains to show that the right-hand side is \emph{classical} as a prestack
(i.e., is isomorphic to the left Kan extension of its restriction to classical affine schemes).
This follows from \corref{c:R is Noeth} using \cite[Proposition 6.8.2]{GR3}.

\end{proof}

\ssec{Derived structure: a stabilization claim}

\sssec{}

For $k\geq 0$ consider the $k$-th coconnective truncation of $\CY$, denoted $^{\leq k}\CY$. Write
$$^{\leq k}\CY=\underset{\alpha\in A}{\on{colim}}\, Y_{\alpha,k},$$
where $Y_{\alpha,k}=\Spec(R_{\alpha,k})\in {}^{\leq k}\!\affSch_{/\sfe}$ and $R_{\alpha,k}=\tau^{\geq -k}(R_\alpha)$. 

\medskip

Set
$$R_k:=\underset{\alpha\in A}{\on{lim}}\, R_{\alpha,k}.$$

Using induction on $k$, we will prove the following statements:

\begin{itemize}

\item(i) The ring $R_k$ connective and for any $k'\leq k$, 
the map $\tau^{\geq -k'}(R_k)\to R_{k'}$ is an isomorphism.  In particular,
the map
$$H^0(R_k)\to R_{\on{cl}}$$
is an isomorphism. 

\smallskip

\item(ii) The map $^{\leq k}\CY\to \Spec(R_k)^\wedge_{Y_{\on{red}}}$
is an isomorphism.

\end{itemize}

\smallskip

Once we prove this, the assertion of \thmref{t:Lurie} will follow by taking the limit over $k$. 

\sssec{}

The base of the induction is the case $k=0$, which has been considered in \secref{ss:cl formal}. We proceed
to the induction step. Thus, we will assume that the statement is true for $k$ and prove it for $k+1$. 

\medskip

Thus, we write
$$^{\leq k+1}\CY=\underset{\alpha\in A}{\on{colim}}\, Y_{\alpha,k+1}.$$

\medskip

We have $Y_{\alpha,k}:={}^{\leq k}Y_{\alpha,k+1}$, and consider the corresponding fiber sequence
$$I_{\alpha,k+1}[k+1]\to R_{\alpha,k+1}\to R_{\alpha,k}, \quad I_{\alpha,k+1}\in \QCoh(Y_{\alpha,k})^\heartsuit.$$

\sssec{}

For each $\alpha$, let $J_{\alpha,k+1}$ denote the object in
$$\on{Pro}(R_{\alpha,k}\mod^\heartsuit)\simeq \on{Pro}(R_{\alpha,\on{cl}}\mod^\heartsuit)$$
given by 
$$\underset{\beta\geq \alpha}{``\on{lim}"}\, H^0\left(I_{\beta,k+1} \underset{R_{\beta,k}}\otimes R_{\alpha,k}\right).$$

We claim: 

\begin{prop} \label{p:J alpha stab}
The object $J_{\alpha,k+1}$ belongs to $R_{\alpha,k}\mod^\heartsuit\subset \on{Pro}(R_{\alpha,k}\mod^\heartsuit)$.
\end{prop} 

Note that \propref{p:J alpha stab} can be reformulated as saying that the canonical map in $\on{Pro}(R_{\alpha,k}\mod)$
$$\underset{\beta\geq \alpha}{\on{lim}}\, H^0\left(I_{\beta,k+1} \underset{R_{\beta,k}}\otimes R_{\alpha,k}\right)
\to \underset{\beta\geq \alpha}{``\on{lim}"}\, H^0\left(I_{\beta,k+1} \underset{R_{\beta,k}}\otimes R_{\alpha,k}\right)$$
is an isomorphism. 

\medskip

The next few subsections are devoted to the proof of \propref{p:J alpha stab}. 

\sssec{}

For a pair of indices $\beta\geq \alpha$, consider
\begin{equation} \label{e:fiber of normals}
\on{Fib}\left(T^*(Y_{\beta,k+1})|_{Y_{\alpha,k}}\to T^*(Y_{\beta,k})|_{Y_{\alpha,k}}\right).
\end{equation} 

By \cite[Chapter 1, Lemma 5.4.3(b)]{GR2}, the object \eqref{e:fiber of normals} lives in cohomological degrees 
$\leq -(k+1)$, and we have
$$H^{-(k+1)}\biggl(\on{Fib}(T^*(Y_{\beta,k+1})|_{Y_{\alpha,k}}\to T^*(Y_{\beta,k})|_{Y_{\alpha,k}})\biggr)\simeq 
H^0\left(I_{\beta,k+1} \underset{R_{\beta,k}}\otimes R_{\alpha,k}\right).$$

Hence, it suffices to show that the object
$$\underset{\beta\geq \alpha}{``\on{lim}"}\, H^{-(k+1)}\biggl(\on{Fib}(T^*(Y_{\beta,k+1})|_{Y_{\alpha,k}}\to T^*(Y_{\beta,k})|_{Y_{\alpha,k}})\biggr)\in
\on{Pro}(R_{\alpha,k}\mod^\heartsuit)$$
belongs to $R_{\alpha,k}\mod^\heartsuit\subset \on{Pro}(R_{\alpha,k}\mod^\heartsuit)$.

\sssec{}

Again by \cite[Chapter 1, Lemma 5.4.3(b)]{GR2}, the maps
$$\on{Fib}(T^*(Y_\beta)|_{Y_{\alpha,k}}\to T^*(Y_{\beta,k})|_{Y_{\alpha,k}})\to \on{Fib}(T^*(Y_{\beta,k+1})|_{Y_{\alpha,k}}\to T^*(Y_{\beta,k})|_{Y_{\alpha,k}})$$
induce isomorphisms on the cohomology in degree $-(k+1)$, 

\medskip

Hence, it suffices to show that the object
\begin{equation} \label{e:pro normals}
\underset{\beta\geq \alpha}{``\on{lim}"}\, \on{Fib}(T^*(Y_\beta)|_{Y_{\alpha,k}}\to T^*(Y_{\beta,k})|_{Y_{\alpha,k}})\in \on{Pro}(R_{\alpha,k}\mod)
\end{equation}
actually belongs to $\QCoh(Y_{\alpha,k})$.
 
\sssec{}

Note that the object \eqref{e:pro normals} identifies with
$$\on{Fib}(T^*(\CY)|_{Y_{\alpha,k}}\to T^*({}^{\leq k}\CY)|_{Y_{\alpha,k}}).$$

Now, $T^*(\CY)|_{Y_{\alpha,k}}$ belongs to $\QCoh(Y_{\alpha,k})$, by the assumption on $\CY$.

\medskip

The object 
$T^*({}^{\leq k}\CY)|_{Y_{\alpha,k}}$ also belongs to $\QCoh(Y_{\alpha,k})$, since $^{\leq k}\CY$ is
a formal completion of an affine scheme, by the inductive hypothesis.

\qed[\propref{p:J alpha stab}]

\sssec{}

Note that we can rewrite
$$H^0\left(I_{\beta,k+1} \underset{R_{\beta,k}}\otimes R_{\alpha,k}\right) \simeq
H^0\left(I_{\beta,k+1} \underset{R_{\beta,\on{cl}}}\otimes R_{\alpha,\on{cl}}\right),$$
and thus think of
$$J_{\alpha,k+1}\simeq \underset{\beta\geq \alpha}{``\on{lim}"}\, 
H^0\left(I_{\beta,k+1} \underset{R_{\beta,\on{cl}}}\otimes R_{\alpha,\on{cl}}\right)
\simeq  \underset{\beta\geq \alpha}{\on{lim}}\, 
H^0\left(I_{\beta,k+1} \underset{R_{\beta,\on{cl}}}\otimes R_{\alpha,\on{cl}}\right)$$
as an object of $R_{\alpha,\on{cl}}\mod^\heartsuit$. 

\begin{cor} \label{c:J alpha fg}
The $R_{\alpha,\on{cl}}$-module $J_{\alpha,k+1}$ is finitely generated.
\end{cor}

\begin{proof}

Follows formally from \propref{p:J alpha stab} using the fact that all
$H^0\left(I_{\beta,k+1} \underset{R_{\beta,\on{cl}}}\otimes R_{\alpha,\on{cl}}\right)$
are finitely generated (the latter follows from the laft assumption).

\end{proof}

\ssec{The induction step, assertion (i)}

\sssec{}

To prove assertion (i) in the induction step, we only have to show that
\begin{equation} \label{e:I alpha}
\underset{\alpha\in A}{\on{lim}}\, I_{\alpha,k+1}
\end{equation} 
lives in cohomological degree $0$. 

\sssec{}

First, we claim that the index category $A$ can be chosen to be the poset $\BN$ of natural numbers.
Indeed, this follows from \cite[Proposition 5.2.3]{GR3}.

\sssec{}

We have
$$\underset{\alpha\in A}{\on{lim}}\, I_{\alpha,k+1} \simeq \underset{\alpha\in A}{\on{lim}}\, \underset{\beta\geq \alpha}{\on{lim}}\, 
H^0\left(I_{\beta,k+1} \underset{R_{\beta,k}}\otimes R_{\alpha,k}\right)\simeq 
\underset{\alpha\in A}{\on{lim}}\, \underset{\beta\geq \alpha}{\on{lim}}\, 
H^0\left(I_{\beta,k+1} \underset{R_{\beta,\on{cl}}}\otimes R_{\alpha,\on{cl}}\right).$$

By \propref{p:J alpha stab}, we can rewrite this further as 
$$\underset{\alpha\in A}{\on{lim}}\, J_{\alpha,k+1},$$
and we claim that the latter object indeed lives in cohomological degree $0$. 

\sssec{}

Note that for 
$\alpha''\geq \alpha'$, we have
$$J_{\alpha',k+1}\simeq H^0(J_{\alpha'',k+1}\underset{R_{\alpha'',\on{cl}}}\otimes R_{\alpha',\on{cl}}).$$

Hence, for $\alpha''\geq \alpha'$, the transition map $J_{\alpha'',k+1}\to J_{\alpha',k+1}$ is surjective. Since 
the category of indices is $\BN$, this implies the desired assertion. 

\qed[Inductive assertion (i)]

\sssec{}

Let $I_{k+1}$ denote the $R_{\on{cl}}$-module 
$$\on{Fib}(R_{k+1}\to R_k)[-k-1] \simeq \underset{\alpha\in A}{\on{lim}}\, I_{\alpha,k+1}.$$

We claim:

\begin{lem} \label{l:stab mod} \hfill

\smallskip

\noindent{\em(a)} For every $\alpha$, the tautological map
$$H^0(I_{k+1}\underset{R_{\on{cl}}}\otimes R_{\alpha,\on{cl}}) \to J_{\alpha,k+1},$$
is an isomorphism, where $J_{\alpha,k+1}\in R_{\alpha,\on{cl}}\mod^\heartsuit$
is the object from \propref{p:J alpha stab}.

\smallskip

\noindent{\em(b)}
The $R_{\on{cl}}$-module $I_{k+1}$ is finitely generated.

\end{lem}

\begin{proof}

Since $R_{\on{cl}}$ is Noetherian, the ring $R_{\alpha,\on{cl}}$ is finitely presented as $R_{\on{cl}}$-module,
and hence the functor 
$$M\mapsto H^0(M\underset{R_{\on{cl}}}\otimes R_{\alpha,\on{cl}}), \quad R_{\on{cl}}\mod^\heartsuit\to \Vect_\sfe^\heartsuit$$
commutes with filtered limits. 

\medskip

Hence, since 
$$I_{k+1}\simeq \underset{\beta\geq \alpha}{\on{lim}}\, I_{\beta,k+1},$$
we can rewrite $H^0(I_{k+1}\underset{R_{\on{cl}}}\otimes R_{\alpha,\on{cl}})$ as
$$\underset{\beta\geq \alpha}{\on{lim}}\, H^0\left(I_{\beta,k+1} \underset{R_{\on{cl}}}\otimes R_{\alpha,\on{cl}}\right)
\simeq \underset{\beta\geq \alpha}{\on{lim}}\, H^0\left(I_{\beta,k+1} \underset{R_{\beta,\on{cl}}}\otimes R_{\alpha,\on{cl}}\right)$$
(where the limit taken in the abelian category $R_{\alpha,\on{cl}}\mod^\heartsuit$), 
while the latter is the same as $J_{\alpha,k+1}$, by \propref{p:J alpha stab}. This proves point (a). 

\medskip

For point (b), let us recall that the $R_{\alpha,\on{cl}}$-modules
$$H^0(I_{k+1}\underset{R_{\on{cl}}}\otimes R_{\alpha,\on{cl}}) \overset{\text{point (a)}}\simeq J_{\alpha,k+1}$$
are finitely generated, by \corref{c:J alpha fg}. 

\medskip

Choose some index $\alpha_0$, and choose a finite set of generators $\ol{m}_1,...,\ol{m}_l$
of $H^0(I_{k+1}\underset{R_{\on{cl}}}\otimes R_{\alpha_0,\on{cl}})$. Let $m_1,...,m_l$ be lifts
of these elements to $I_{k+1}$. We will show that the elements $m_1,...,m_l$ generate $I_{k+1}$
as a $R_{\on{cl}}$-module.

\medskip

First, since the ideals $\on{ker}(R_{\alpha,\on{cl}}\to R_{\alpha_0,\on{cl}})$ are nilpotent, we obtain that the images
of $m_1,...,m_l$ in every $H^0(I_{k+1}\underset{R_{\on{cl}}}\otimes R_{\alpha,\on{cl}})$ generate
it as $R_{\alpha,\on{cl}}$-module. 

\medskip

Now, we note that for every $\alpha$, the map
\begin{equation} \label{e:ideal as limit}
I_{k+1}\to \underset{\alpha}{\on{lim}}\, H^0(I_{k+1}\underset{R_{\on{cl}}}\otimes R_{\alpha,\on{cl}})
\end{equation}
is an isomorphism. Indeed, using point (a), we rewrite the right-hand side in \eqref{e:ideal as limit}
as
$$\underset{\alpha}{\on{lim}}\, J_{\alpha,k+1} \simeq
\underset{\alpha}{\on{lim}}\, \underset{\beta\geq \alpha}{\on{lim}}\, 
H^0(I_{\beta,k+1}\underset{R_{\beta,\on{cl}}}\otimes R_{\alpha,\on{cl}})\overset{\text{cofinality}}\simeq 
\underset{\alpha}{\on{lim}}\, 
H^0(I_{\alpha,k+1}\underset{R_{\alpha,\on{cl}}}\otimes R_{\alpha,\on{cl}}) \simeq 
\underset{\alpha}{\on{lim}}\, 
I_{\alpha,k+1},$$
which is the same as $I_{k+1}$.

\medskip

This implies the required generation property by \lemref{l:adic}(b) since the category of indices can be assumed to be $\BN$.

\end{proof}

\ssec{The induction step, assertion (ii)}

\sssec{}

We now proceed to the proof of assertion (ii) in the induction step, i.e., we wish to show that
\begin{equation} \label{e:trunk k+1}
^{\leq k+1}\CY\to \Spec(R_{k+1})^\wedge_{Y_{\on{red}}}
\end{equation}
is an isomorphism. The left-hand side is $(k+1)$-coconnective as a prestack, by definition.

\medskip

We claim that the right-hand side in \eqref{e:trunk k+1} is also $(k+1)$-coconnective 
as a prestack. Indeed, this follows by repeating the argument of \cite[Proposition 6.8.2]{GR3}
using the following:

\smallskip

\noindent{(i)} $R_{\on{cl}}$ is Noetherian. This is the content of \corref{c:R is Noeth}.

\smallskip

\noindent{(ii)} The $R_{\on{cl}}$-modules $I_{k'}:=\on{ker}(R_{k'}\to R_{k'-1})[-k']$
are finitely generated for $k'=1,...,k+1$. This is the content of \lemref{l:stab mod}(b). 

\medskip

Hence, it is enough to show that \eqref{e:trunk k+1} is an isomorphism 
when evaluated on $(k+1)$-coconnective affine schemes. By the inductive hypothesis,
we know that it is an isomorphism when evaluated on $k$-coconnective affine schemes.

\sssec{}

Consider the following situation. Let 
$$\CY_0\to \CY_1 \to \CY_2$$
be maps of prestacks that admit deformation theory. 

\medskip

Assume all three maps induce isomorphisms on $k$-coconnective affine schemes.
We claim:

\begin{lem} \label{l:k+1 coconn}
Suppose that for any $k$-coconnective affine scheme $S$ equipped with a map 
$S\to \CY_0$, the map
$$T^*(\CY_0/\CY_2)|_S \to T^*(\CY_0/\CY_1)|_S$$
in $\on{Pro}(\QCoh(S)^{<\infty})$
induces an isomorphism on the $\tau^{\geq -(k+2)}$ truncations.  Then 
the map $\CY_1\to \CY_2$ induces an isomorphism
on $(k+1)$-coconnective affine schemes.
\end{lem}

The proof of the lemma is given below. Let us show how the assertion of the lemma implies
the induction step.

\sssec{}

We will apply \lemref{l:k+1 coconn} to
$$\CY_0={}^{\leq k}\CY \simeq \Spec(R_k)^\wedge_{Y_{\on{red}}},\,\,
\CY_1:={}^{\on{conv}}({}^{\leq k+1}\CY),\,\, \CY_2=\Spec(R_{k+1})^\wedge_{Y_{\on{red}}}.$$

\medskip

Let us check that the condition of the lemma holds. By cofinality, we can assume that 
$S$ is one of the schemes $Y_{\alpha,k}$.

\medskip

We have
$$T^*({}^{\leq k}\CY/{}^{\on{conv}}({}^{\leq k+1}\CY))|_{Y_{\alpha,k}} \simeq
\underset{\beta\geq \alpha}{``\on{lim}"}\, 
T^*(Y_{\beta,k}/Y_{\beta,k+1})|_{Y_{\alpha,k}}.$$

For every $\beta$, we have
\begin{multline*}
\tau^{\geq -(k+2)}(T^*(Y_{\beta,k}/Y_{\beta,k+1}))=I_{\beta,k+1}[k+2]\in
R_{\beta,\on{cl}}\mod^\heartsuit[k+2] \simeq R_{\beta,k}\mod^\heartsuit[k+2] \simeq \\
\simeq \QCoh(Y_{\beta,k})^\heartsuit[k+2]  \subset \QCoh(Y_{\beta,k}).
\end{multline*}

Hence, 
$$\tau^{\geq -(k+2)}\left(T^*({}^{\leq k}\CY/{}^{\on{conv}}({}^{\leq k+1}\CY))|_{Y_{\alpha,k}}\right) \simeq J_{\alpha,k+1}[k+2].$$ 

\medskip

Similarly,
$$\tau^{\geq -(k+2)}\left(T^*(\Spec(R_k)^\wedge_{Y_{\on{red}}}/\Spec(R_{k+1})^\wedge_{Y_{\on{red}}})\right)
\simeq I_{k+1}[k+2],$$
and hence
$$\tau^{\geq -(k+2)}\left(T^*(\Spec(R_k)^\wedge_{Y_{\on{red}}}/\Spec(R_{k+1})^\wedge_{Y_{\on{red}}})|_{Y_{\alpha,k}}\right) 
\simeq H^0(I_{k+1}\underset{R_k}\otimes R_{\alpha,k})[k+2] \simeq H^0(I_{k+1}\underset{R_{\on{cl}}}\otimes R_{\alpha,\on{cl}})[k+2].$$

The desired isomorphism follows now from \lemref{l:stab mod}(a).

\qed[Inductive assertion (ii)]

\sssec{Proof of \lemref{l:k+1 coconn}} 

Thus, let us be given a pair $(S,S')$, where $S'$ is a $(k+1)$-coconnective affine scheme and $S$
is its $k$-coconnective truncation. The morphism $S\to S'$ has a canonical structure of square-zero
extension corresponding to an object $\CF\in \QCoh(S)^{\heartsuit}[k+1]$ and a map
$$T^*(S)\to \CF[1].$$

\medskip 

Let us be given a map $y:S\to \CY_0$. Denote by $y_1$ and $y_2$ the composite maps
to $\CY_1$ and $\CY_2$, respectively.

\medskip

The datum of extension of $y_1$ to a map $y'_1:S'\to \CY_1$ 
is equivalent to a null-homotopy of the composition
$$y_1^*(T^*(\CY_1))\to  y^*(T^*(\CY_0))\to T^*(S) \to  \CF[1],$$
and similarly for $y_2$.

\medskip

Hence, it is enough to show that restriction along
\begin{equation} \label{e:y pullback rel cotan}
y^*(T^*(\CY_0/\CY_2)) \to y^*(T^*(\CY_0/\CY_1))
\end{equation}
induces an isomorphism on spaces of maps to $\CF[1]$.

\medskip

Since $\CF[1]\in \QCoh(S)^\heartsuit[k+2]\subset \QCoh(S)^{\geq -(k+2)}$, it suffices to show that the map
\eqref{e:y pullback rel cotan} induces an isomorphism on the $\tau^{\geq -(k+2)}$ truncations. However,
the latter is the assumption in the lemma.

\qed[\lemref{l:k+1 coconn}]

\section{Colimits over $\on{TwArr}(\fSet)$} \label{s:twisted arr}

The goal of this section is to prove \lemref{l:fact homology} and related statements that involve
colimits over the category $\on{TwArr}(\fSet)$. 

\ssec{Operadic left Kan extensions}

\sssec{}
Let $O$ and $O'$ be symmetric monoidal categories.  We will denote by
$$ \on{Funct}^{\otimes}(O,O'), \ \on{Funct}^{\otimes\on{-rlax}}(O,O'), \ \on{Funct}^{\otimes\on{-llax}}(O,O') $$
respectively, for categories strict, right-lax, and left-lax symmetric monoidal functors.

\sssec{}
We now recall a particular aspect of the theory of operadic left Kan extensions \cite[Sect. 3.1]{Lu2}.

\medskip

In what follows, we will say that a symmetric monoidal category $\bA$ is a cocomplete
symmetric monoidal category if the underlying category $\bA$ is cocomplete and the tensor product commutes with colimits in each variable.

\medskip

The following is a special case of \cite[Corollary 3.1.3.5]{Lu2}\footnote{In the notation of \cite{Lu2}, the commutativity of the diagram \eqref{e:underlying of operadic LKE} follows from the fact that for a symmetric monoidal functor $O \to O'$, and any object $o'\in O'$, the functor $O/o' \to (O^{\otimes}_{\on{act}})/o'$ is cofinal.}.

\begin{thm}\label{t:operadic LKE}
Suppose that $F: O \to O'$ is a symmetric monoidal functor between (small) symmetric monoidal categories.  For any cocomplete symmetric monoidal category $\bA$, the restriction functor
$$\on{Res}_F: \on{Funct}^{\otimes\on{-rlax}}(O', \bA) \to \on{Funct}^{\otimes\on{-rlax}}(O, \bA)$$
admits a left adjoint $\on{LKE}_F^{\otimes}$ (the ``operadic left Kan extension'') and the
following diagram commutes
\begin{equation}\label{e:underlying of operadic LKE}
\xymatrix{
\on{Funct}^{\otimes\on{-rlax}}(O,\bA)\ar[r]^-{\on{LKE_F^{\otimes}}}\ar[d]_{\oblv} & \on{Funct}^{\otimes\on{-rlax}}(O', \bA) \ar[d]^{\oblv}\\
\on{Funct}(O, \bA) \ar[r]^-{\on{LKE}_F} & \on{Funct}(O', \bA)
}.
\end{equation}
In particular, the left Kan extension of a right-lax symmetric monoidal functor along a symmetric monoidal functor is canonically right-lax symmetric monoidal.
\end{thm}

\sssec{}
Unraveling the definitions, given a right-lax symmetric monoidal functor $\Phi: O \to \bA$, and $o_1, o_2\in O'$, the structure map
$$ (\on{LKE}_F \Phi)(o_1') \otimes (\on{LKE}_F \Phi)(o_2') \to (\on{LKE}_F \Phi)(o_1'\otimes o_2') $$
is the composite
$$ \underset{o_1 \in O/o_1'}{\on{colim}} \Phi(o_1) \otimes \underset{o_2 \in O/o_2'}{\on{colim}} \Phi(o_2) \overset{\sim}{\leftarrow} \underset{(o_1, o_2) \in O/o_1' \times O/o_2'}{\on{colim}} \Phi(o_1)\otimes \Phi(o_2) \to$$ 
$$\to \underset{(o_1, o_2) \in O/o_1' \times O/o_2'}{\on{colim}} \Phi(o_1 \otimes o_2) \to \underset{o \in O/(o_1' \otimes o_2')}{\on{colim}} \Phi(o),$$
where the first map is an isomorphism since the tensor product in $\bA$ commutes with colimits in each variable, the middle map is the right-lax structure of $\Phi$ and the last map is induced by the functor
$O/o_1' \times O/o_2' \to O/(o_1'\otimes o_2')$ given by tensor product (using the fact that $F$ is strictly symmetric monoidal).

\sssec{}

From the above discussion we obtain:

\begin{prop}\label{p:strict operadic LKE}
Let $F: O \to O'$ be a symmetric monoidal functor such that for any pair of objects $o_1', o_2' \in O'$
the functor
$$ O/o_1' \times O/o_2' \to O/(o_1'\otimes o_2') $$
given by tensor product is cofinal.  Then for any cocomplete symmetric monoidal category $\bA$ and any (strict) symmetric monoidal functor $\Phi: O \to \bA$, the operadic left Kan extension
$\on{LKE}^{\otimes}_F(\Phi): O'\to \bA $ is strictly symmetric monoidal.
\end{prop}

\sssec{}
We now specialize \thmref{t:operadic LKE} to the case that $O'=\{*\}$.
By definition, we have
$$ \on{ComAlg}(\bA) = \on{Funct}^{\otimes\on{-rlax}}(\{*\}, \bA) .$$
Given any symmetric monoidal category $O$, restriction along the terminal symmetric monoidal functor $O \to \{*\}$ gives a diagonal functor
\begin{equation}\label{e:diag alg}
 \on{ComAlg}(\bA) \simeq \on{Funct}^{\otimes\on{-rlax}}(\{*\}, \bA) \to \on{Funct}^{\otimes\on{-rlax}}(O, \bA)
\end{equation}

In this case, \thmref{t:operadic LKE} gives:

\begin{cor} \label{c:colim ten}
Let $O$ be a symmetric monoidal category.  For any cocomplete symmetric monoidal category $\bA$,
the diagonal functor \eqref{e:diag alg} admits a left adjoint
$$ \on{colim}_O^{\otimes}: \on{Funct}^{\otimes\on{-rlax}}(O, \bA) \to \on{ComAlg}(\bA) $$
which on underlying objects is given by colimit along $O$.
In particular,
the colimit of a right-lax symmetric monoidal functor is canonically a commutative algebra object.
\end{cor}

\ssec{Proof of \lemref{l:fact homology}} \label{ss:proof fact homology}

The proof we present was communicated to us by J.~Campbell.

\sssec{}\label{sss:adjoint setup}
Suppose we have a symmetric monoidal functor $\Phi: \bA' \to \bA$ between cocomplete symmetric monoidal categories.  The functor $\Phi$ induces a functor
\begin{equation}\label{e:induced comalg functor}
\on{ComAlg}(\bA') \to \on{ComAlg}(\bA) .
\end{equation}
Now suppose that the underlying functor fo $\Phi$ admits a left adjoint $\Phi^L$.  Our present goal is to
formulate and prove \propref{p:adjoint on algebras} which uses $\Phi^L$ to give a description of the left adjoint to \eqref{e:induced comalg functor}.  Note that in general, $\Phi^L$
itself is only a \emph{left-lax} symmetric monoidal functor and therefore does not induce a functor between commutative algebras.

\sssec{}
Recall \cite[Construction 2.2.4.1 and Proposition 2.2.4.9]{Lu2} that given a symmetric monoidal category $O$, there exist universal symmetric monoidal categories $\on{RLax}(O)$ and $\on{LLax}(O)$ equipped with, respectively, right-lax and left-lax
symmetric monoidal functors
$$ O \to \on{RLax}(O) \ \mbox{ and } \  O \to \on{LLax}(O) $$
which induce equivalences
$$ \on{Funct}^{\otimes}(\on{RLax}(O), \bA) \simeq \on{Funct}^{\otimes\on{-rlax}}(O, \bA) \mbox{ and }
 \on{Funct}^{\otimes}(\on{LLax}(O), \bA) \simeq \on{Funct}^{\otimes\on{-llax}}(O, \bA) $$
for any symmetric monoidal category $\bA$.  Evidently,
$$ \on{LLax}(O) \simeq \on{RLax}(O^{\on{op}})^{\on{op}} .$$

\sssec{}
We have $\on{RLax}(\{*\})=\fSet$, with the symmetric monoidal structure given by disjoint union.  In particular, this gives the equivalence
$$ \on{ComAlg}(\bA) \simeq \on{Funct}^{\otimes}(\fSet, \bA) $$
for any symmetric monoidal category $\bA$.

\sssec{}
Now, suppose we have a symmetric monoidal functor $\Phi: \bA' \to \bA$ which admits a left adjoint $\Phi^L$.  
Since $\Phi^L$ is canonically left-lax symmetric monoidal, we obtain a functor
\begin{multline}\label{e:llax left adjoint}
\on{ComAlg}(\bA) \simeq \on{Funct}^{\otimes\on{-rlax}}(\{*\}, \bA)\simeq
 \on{Funct}^{\otimes}(\fSet, \bA)\hookrightarrow \\
\hookrightarrow \on{Funct}^{\otimes\on{-llax}}(\fSet, \bA) \overset{\Phi^L \circ}{\to}
\on{Funct}^{\otimes\on{-llax}}(\fSet, \bA') .
\end{multline}

\sssec{}
Let $\on{TwArr}(\fSet)$ denote the twisted arrow category of $\fSet$ with tensor product given by disjoint union.  
Applying \cite[Construction 2.2.4.1]{Lu2} to (the opposite category of) $\fSet$, we obtain:

\begin{prop}\label{p:llax(fin)}
The left-lax symmetric monoidal functor
$$ \fSet \to \on{TwArr}(\fSet) $$
given by $I \mapsto (I \to \{*\})$ induces an equivalence
$$ \on{LLax}(\fSet) \simeq \on{TwArr}(\fSet) .$$
\end{prop}

\sssec{}
Composing the functor \eqref{e:llax left adjoint} with the equivalence of \propref{p:llax(fin)}, we obtain a functor
\begin{equation}\label{e:adjoint to twarr}
\on{ComAlg}(\bA) \to \on{Funct}^{\otimes}(\on{TwArr}(\fSet), \bA').
\end{equation}

Explicitly, given $R \in \on{ComAlg}(\bA)$, the corresponding functor
$$ \on{TwArr}(\fSet) \to \bA' $$
is given by
\begin{equation}
 (I \overset{\psi}{\to} J) \mapsto \underset{j\in J}{\otimes} \Phi^L(R^{\otimes \psi^{-1}(j)}).
\end{equation}

\sssec{}
The following is a more precise version of \lemref{l:fact homology}:

\begin{prop}\label{p:adjoint on algebras}
Let $\Phi: \bA' \to \bA$ be a symmetric monoidal functor between cocomplete symmetric monoidal categories which admits a left adjoint $\Phi^L$.
Then the induced functor $$\on{ComAlg}(\bA') \to \on{ComAlg}(\bA)$$ admits a left adjoint given by the composite
$$
\xymatrix{
\on{ComAlg}(\bA) \ar[r]^-{\eqref{e:adjoint to twarr}} & 
\on{Funct}^{\otimes}(\on{TwArr}(\fSet), \bA') \ar[r]^-{\on{colim}^{\otimes}} & \on{ComAlg}(\bA'),
}
$$
where $\on{colim}^{\otimes}$ is the composition
\begin{equation} \label{e:colim ten}
\on{Funct}^{\otimes}(\on{TwArr}(\fSet), \bA') \hookrightarrow \on{Funct}^{\otimes\on{-rlax}}(\on{TwArr}(\fSet), \bA') 
\overset{\on{colim}_{\on{TwArr}(\fSet)}}\to \on{ComAlg}(\bA'),
\end{equation}
where $\on{colim}_{\on{TwArr}(\fSet)}$ is as in \corref{c:colim ten}.
\end{prop}

\sssec{}
Before proving \propref{p:adjoint on algebras}, we establish the following:

\begin{prop}\label{p:twarr adjoint}
Let $\bA$ be a cocomplete symmetric monoidal category.  Then the inclusion functor
$$ \on{ComAlg}(\bA) \simeq \on{Funct}^{\otimes}(\fSet, \bA) 
\hookrightarrow \on{Funct}^{\otimes\on{-llax}}(\fSet, \bA) \simeq \on{Funct}^{\otimes}(\on{TwArr}(\fSet), \bA) $$
admits a left adjoint given by
$$ \on{colim}^{\otimes}: \on{Funct}^{\otimes}(\on{TwArr}(\fSet), \bA) \to \on{ComAlg}(\bA) ,$$
where $\on{colim}^{\otimes}$ is as in \eqref{e:colim ten}.
\end{prop}

\begin{proof}
The inclusion functor is given by restriction along the symmetric monoidal functor
$$s: \on{TwArr}(\fSet) \to \fSet $$
given by $(I \to J) \mapsto I$.  By \thmref{t:operadic LKE}, we have an adjunction
$$ \on{LKE}^{\otimes}_s:  \on{Funct}^{\otimes\on{-rlax}}(\on{TwArr}(\fSet), \bA) {}^{\longrightarrow}_{\longleftarrow}
 \on{Funct}^{\otimes\on{-rlax}}(\fSet, \bA) : \on{Res}_s. $$
 
 \medskip

Since $s$ is strictly symmetric monoidal, the restriction functor $\on{Res}_s$ preserves strictly symmetric monoidal functors.  
Furthermore, for every $I_1, I_2 \in \fSet$, the functor
$$ \on{TwArr}(\fSet)_{/I_1} \times \on{TwArr}(\fSet)_{/I_2} \to \on{TwArr}(\fSet)_{/I_1 \sqcup I_2} $$
is cofinal.  Therefore, by \propref{p:strict operadic LKE},
$\on{LKE}_s^{\otimes}$ also preserves strictly symmetric monoidal functors.  Thus,
the desired left adjoint is given by
$$\on{Funct}^{\otimes}(\on{TwArr}(\fSet), \bA) \overset{\on{LKE}^{\otimes}_s}\to
\on{Funct}^{\otimes}(\fSet, \bA) \simeq \on{ComAlg}(\bA) .$$

\medskip

It remains to show that this functor is canonically isomorphic to $\on{colim}^{\otimes}$.  The functor
$\on{colim}_{\on{TwArr}}^{\otimes}$ is given by the operadic left Kan extension along the composite
$$ \on{TwArr}(\fSet) \overset{s}{\to} \fSet \overset{p}{\to} \{*\} $$
Thus, it suffices to show that the
composite functor
$$ \on{ComAlg}(\bA) \simeq \on{Funct}^{\otimes}(\fSet, \bA) \hookrightarrow \on{Funct}^{\otimes\on{-rlax}}(\fSet, \bA) \overset{\on{LKE}^{\otimes}_p}{\to} \on{Funct}^{\otimes\on{-rlax}}(\{*\}, \bA) \simeq \on{ComAlg}(\bA) $$
is canonically isomorphic to the identity.  By \cite[Corollary 7.3.2.7]{Lu2}, since $p$ is symmetric monoidal, $\on{LKE}^{\otimes}_p$
is given by restriction along the right adjoint
$$ \{*\}\hookrightarrow \fSet,$$
which is the universal right-lax symmetric monoidal functor from $\{*\}$.  This gives the desired result.
\end{proof}

\sssec{Proof of \propref{p:adjoint on algebras}}
Since $\Phi$ is symmetric monoidal, by \cite[Corollary 7.3.2.7]{Lu2}, we obtain an adjunction
$$ \Phi^L \circ (-): \on{Funct}^{\otimes\on{-llax}}(\fSet, \bA) {}^{\longrightarrow}_{\longleftarrow} \on{Funct}^{\otimes\on{-llax}}(\fSet, \bA'): \Phi \circ (-). $$
We have a commutative diagram
$$ \xymatrix{
\on{Funct}^{\otimes}(\fSet, \bA') \ar[r]^-{\Phi \circ (-)}\ar@{^{(}->}[d] & \on{Funct}^{\otimes}(\fSet, \bA) \ar@{^{(}->}[d]\\
\on{Funct}^{\otimes\on{-llax}}(\fSet, \bA') \ar[r]^-{\Phi \circ (-)} & \on{Funct}^{\otimes\on{-llax}}(\fSet, \bA)
}$$
in which the vertical functors are fully faithful.  The desired result now follows from \propref{p:twarr adjoint}. \qed

\sssec{Proof of \lemref{l:coEnd abs}} \label{sss:proof coEnd}

We will show that the colimit expression \eqref{e:coEnd} with $C:=B^\vee$
possesses the stated universal property.

\medskip

Consider the symmetric monoidal functor
$$\Phi: (-) \otimes B: \bO \to B\mod(\bO).$$
It induces a functor
\begin{equation}\label{e:comalg tensor}
\on{ComAlg}(\bO) \to \on{ComAlg}(B\mod(\bO))
\end{equation}

Note that the value of $\Phi^L$ on an object of the form $M\otimes B$, $M\in \bO$ is $C\otimes M$. 
Hence, $\on{coEnd}(A,C)$ is the colimit expression in \propref{p:adjoint on algebras}
for the above functor $\Phi$ evaluated on $A\otimes B$. 

\medskip

Applying \propref{p:adjoint on algebras}, we obtain that the above colimit identifies with
the value of the left adjoint of \eqref{e:comalg tensor} evaluated on $A\otimes B\in \on{ComAlg}(B\mod(\bO))$. The 
latter has the required universal property by definition. 

\qed

\ssec{Proof of \lemref{l:maps from coend}} \label{ss:proof coEnd assoc}

\sssec{}

By definition, $\coHom(A,B)$ is value of the left adjoint of the functor
\begin{equation}\label{e:def of cohom}
(-) \otimes B: \on{ComAlg}(\bO) \to \on{ComAlg}(\bO)
\end{equation}
applied to $A \in \on{ComAlg}(\bO)$.

\medskip

Consider the right-lax symmetric monoidal functor
\begin{equation}
T_{B}: \bO \to \on{Funct}(\fSet, \bO)
\end{equation}
given by $X \mapsto (I \mapsto X \otimes B^{\otimes I})$. The functor \eqref{e:def of cohom} factors as
$$ \on{ComAlg}(\bO) \overset{\on{ComAlg}(T_{B})}{\longrightarrow} \on{ComAlg}(\on{Funct}(\fSet, \bO)) \simeq \on{Funct}(\fSet, \on{ComAlg}(\bO)) \overset{\on{ev}_{\{*\}}}{\to} \on{ComAlg}(\bO) ,$$
where the last functor is evaluation at $\{*\} \in \fSet$.  Thus we have that
$$ \coHom(A,B) \simeq \on{ComAlg}(T_{B})^L \circ \on{ev}_{\{*\}}^L(A) .$$
The functor $\on{ev}_{\{*\}}^L$ is given by left Kan extension along the inclusion $\{*\} \hookrightarrow \fSet$ and so
$$ \on{ev}_{\{*\}}^L(A) \simeq (I \mapsto A^{\otimes I}) \in \on{Funct}(\fSet, \on{ComAlg}(\bO)).$$

\sssec{}

Thus, we obtain that \lemref{l:maps from coend} follows from the following two assertions:

\begin{enumerate}
\item[(a)] The natural map
$$ \on{AssocAlg}(T_{B})^L(\oblv_{\on{Com}\to \on{Assoc}}) \to \oblv_{\on{Com}\to \on{Assoc}}(\on{ComAlg}(T_{B})^L) $$
of functors $\on{ComAlg}(\on{Funct}(\fSet, \bO)) \to \on{AssocAlg}(\bO)$ is an isomorphism.

\smallskip

\item[(b)] The natural map
$$ T_{B}^L (\oblv_{\on{Com}}) \to \oblv_{\on{Com}}(\on{ComAlg}(T_{B})^L) $$
of functors $\on{ComAlg}(\on{Funct}(\fSet, \bO)) \to \bO$ is an isomorphism.
\end{enumerate}

\medskip

To prove both assertions, it is enough to show that the left adjoint $T_{B}^L$, which is a priori left-lax symmetric monoidal, is strictly symmetric monoidal.  
Indeed, in this case, by \cite[Corollary 7.3.2.7]{Lu2}, we have
$$ \on{ComAlg}(T_{B})^L \simeq \on{ComAlg}(T_{B}^L)\  \mbox{ and } \on{AssocAlg}(T_{B})^L \simeq \on{AssocAlg}(T_{B}^L) .$$

\sssec{}

It remains to prove that the left-lax symmetric monoidal structure on $T_{B}^L$ is strict.  By definition, we have
$$ \on{Maps}_{\bO}(T_{B}^L(F), X) \simeq \on{Maps}_{\on{Funct}(\fSet, \bO)}(F, T_{B}(X)) .$$
However, (see e.g. \cite[Lemma 1.3.12]{GKRV}), the latter expression is canonically identified with
$$ \underset{(I\to J) \in \on{TwArr}(\fSet)^{\on{op}}}{\on{lim}} \on{Maps}_{\bO}(F(I), X \otimes B^{\otimes J}) \simeq \underset{(I\to J) \in \on{TwArr}(\fSet)^{\on{op}}}{\on{lim}} \on{Maps}_{\bO}(F(I) \otimes C^{\otimes J}, X),$$
where $C=B^\vee$. 

\medskip

Thus, we have
$$ T_{B}^L(F) \simeq \underset{(I \to J) \in \on{TwArr}(\fSet)}{\on{colim}} F(I) \otimes C^{\otimes J}. $$
Moreover, the left-lax structure map
\begin{multline}\label{e:left lax structure map}
 T_{B}^L(F_1 \otimes F_2) \simeq \underset{(I \to J) \in \on{TwArr}(\fSet)}{\on{colim}} F_1(I) \otimes F_2(I) \otimes C^{\otimes J} \to \\
 \to T_{B}^L(F_1)\otimes T_{B}^L(F_2) \simeq
 \underset{(I_1 \to J_1, I_2 \to J_2)\in \on{TwArr}(\fSet)^{\times 2}}{\on{colim}} F_1(I_1)\otimes F_2(I_2) \otimes C^{\otimes J_1}\otimes C^{\otimes J_2}
\end{multline}
is induced by the maps
$$\on{Id} \otimes \on{comult}^{\otimes J}:
F_1(I) \otimes F_2(I) \otimes C^{\otimes J} \to F_1(I) \otimes F_2(I) \otimes C^{\otimes J} \otimes C^{\otimes J}.$$

\sssec{}

Consider the category
$$ \on{TwArr}(\fSet) \underset{\fSet}{\times} \fSet^{\times 2},$$
where the functor $\fSet^{\times 2} \to \fSet$ is given by coproduct.  In other words, an object of this category
consists of three finite sets $I_1, I_2, J$ and a map $I_1 \sqcup I_2 \to J$.

\medskip

The functor $\on{TwArr}(\fSet)^{\times 2} \to \bO$ in the right-hand side of \eqref{e:left lax structure map} is given by restriction along the functor
\begin{equation}\label{e:twarr^2 reduction}
\on{TwArr}(\fSet)^{\times 2} \to \on{TwArr}(\fSet) \underset{\fSet}{\times} \fSet^{\times 2}
\end{equation}
given by $(I_1 \to J_1, I_2 \to J_2) \mapsto (I_1, I_2, I_1 \sqcup I_2 \to J_1 \sqcup J_2)$.  The functor \eqref{e:twarr^2 reduction}
admits a left adjoint given by $(I_1, I_2, I_1 \sqcup I_2 \to J) \mapsto (I_1 \to J, I_2 \to J)$ and is therefore cofinal.  

\medskip

Hence, we can rewrite 
the right-hand side of \eqref{e:left lax structure map} as
$$\underset{(I_1 \sqcup I_2 \to J)\in \on{TwArr}(\fSet)\underset{\fSet}{\times} \fSet^{\times 2}}{\on{colim}}\,\,
F_1(I_1)\otimes F_2(I_2) \otimes C^{\otimes J}.$$

\medskip

Furthermore, the map \eqref{e:left lax structure map} is induced by the functor 
$$ \on{TwArr}(\fSet) \to \on{TwArr}(\fSet) \underset{\fSet}{\times} \fSet^{\times 2} $$
given by $(I \to J) \mapsto (I \sqcup I \to J)$. This functor also 
admits a left adjoint and is therefore also cofinal.

\qed[\lemref{l:maps from coend}]

\section{Semi-rigid symmetric monoidal categories} \label{s:semi-rigid}

\ssec{Definition and examples}

\sssec{} \label{sss:semi-rigid defn}

Let $\bA$ be a (unital) symmetric monoidal DG category. We shall say that $\bA$ is \emph{semi-rigid} if the following 
two conditions are satisfied:

\begin{itemize}

\item(i) The functor $\on{mult}_\bA:\bA\otimes \bA\to \bA$ admits a continuous right adjoint
(to be denoted $\on{comult}_\bA$) and the structure on $\on{comult}_\bA$ of right-lax compatibility 
with the $\bA$-bimodule structure is strict.

\item(ii) $\bA$ is dualizable as a DG category.

\end{itemize}

\sssec{}

Of course, a rigid symmetric monoidal category is semi-rigid. 

\medskip

The one property that distinguishes rigid from semi-rigid is that in the former case we 
require that the unit object $\one_\bA\in \bA$ be compact.

\medskip

As we will see shortly, semi-rigid categories behave in a way very similar to rigid
ones with respect to 2-categorical properties, i.e., from the point of view of module
categories over them. 

\medskip

However, their internal structure is very different in that dualizable objects are not necessarily 
compact. However, in a compactly generated semi-rigid category,  compact objects are still 
dualizable, see \secref{ss:dual compact}. 

\sssec{}

A key example for this paper of a semi-rigid category is $\QCoh(\CY)$, where $\CY$ is a formal affine scheme.

\medskip

Condition (i) in \secref{sss:semi-rigid defn} holds because the diagonal morphism $\CY\to \CY\times \CY$
is affine. 

\medskip

Condition (ii) in \secref{sss:semi-rigid defn} holds because $\QCoh(\CY)$ is compactly generated.

\sssec{}

Let $\CY$ now be of the form $\CY'/\sG$, where $\CY'$ is a formal affine scheme. 
Then $\QCoh(\CY)$ is semi-rigid, see \lemref{l:QCoh s-rigid} and \corref{c:descent do quot form aff}(a). 

\sssec{}

Here is another way to deduce that $\QCoh(\CY)$ is semi-rigid when $\CY$ is a formal affine scheme,
realized as a formal completion of an affine scheme.

\medskip

Suppose that $\bA$ is a semi-rigid symmetric monoidal category, and let $\Phi:\bA\to \bA'$ be a
symmetric monoidal functor. Assume that $\Phi$ admits a left adjoint, to be denoted $\Phi^L$, such that:

\begin{itemize}

\item $\Phi^L$ is fully faithful;

\item The left-lax monoidal structure on $\Phi^L$ is strict (but not necessarily unital).

\end{itemize}

Then $\bA'$ is also semi-rigid.

\sssec{}

A completely different example of a semi-rigid category is $\Shv^{\on{all}}(Y)$ for a scheme $Y$ (or more generally, for a Hausdorff 
locally compact topological space).

\medskip

Indeed, in this case, the functor $\on{comult}_\bA$ identifies with
$$(\Delta_Y)_*\simeq (\Delta_Y)_!.$$

The fact that $\Shv^{\on{all}}(Y)$ is dualizable is also known, see \secref{sss:Shv Betti self-dual}. 

\ssec{Properties of semi-rigid categories}

In this subsection we let $\bA$ be a semi-rigid symmetric monoidal category.

\sssec{}

Set 
$$\sR_{\bA}:=\on{comult}_\bA(\one_\bA)\in \bA\otimes \bA.$$

This object has a canonical structure of commutative algebra in $\bA\otimes \bA$.

\sssec{}

Here is the first observation:

\begin{prop} \label{p:HH semi-rigid}
Let $\bA$ be semi-rigid\footnote{In fact we only need condition (i).}. Then for a $\bA$-bimodule category $\bP$, there is
a canonical equivalence
$$\on{Funct}_{(\bA\otimes \bA)\mod}(\bA,\bP)\simeq \bP\underset{\bA\otimes \bA}\otimes \bA.$$
\end{prop}

\begin{proof}

The adjunction
\begin{equation} \label{e:mult comult again}
\on{mult}_\bA:\bA\otimes \bA \rightleftarrows \bA:\on{comult}_\bA
\end{equation}
as $\bA$-bimodules gives rise to an adjunction
$$\bP \simeq \on{Funct}_{(\bA\otimes \bA)\mod}(\bA\otimes \bA,\bP) \rightleftarrows \on{Funct}_{(\bA\otimes \bA)\mod}(\bA,\bP),$$
where the right adjoint identifies with the forgetful functor
$$\on{Funct}_{(\bA\otimes \bA)\mod}(\bA,\bP)\to \on{Funct}_{\on{cont}}(\bA,\bP)\to \on{Funct}_{\on{cont}}(\Vect_\sfe,\bP)=\bP;$$
in particular, it is conservative, and hence monadic. Furthermore, the resulting monad on $\bP$ is given by the action of the object
$\sR_{\bA}\in \bA\otimes \bA$.

\medskip

Now, tensoring \eqref{e:mult comult again} we obtain an adjunction
$$\bP \simeq \bP\underset{\bA\otimes \bA}\otimes (\bA\otimes \bA)
\rightleftarrows  \bP\underset{\bA\otimes \bA}\otimes \bA,$$
where the essential image of the left adjoint generates the target category, and hence 
the right adjoint is conservative and hence monadic. The resulting monad on $\bP$ 
also identifies with one given the action of the object
$\sR_{\bA}\in \bA\otimes \bA$.

\medskip

Hence, we have identified both categories in the statement of the proposition with
$$\sR_{\bA}\mod(\bP).$$

\end{proof}

In the course of the proof, we have also shown:

\begin{cor} \label{c:dualizab persists}
Suppose that $\bP$ is dualizable as a plain DG category. Then so is 
$\bP\underset{\bA\otimes \bA}\otimes \bA$.
\end{cor}

\begin{proof}
The category $\bP\underset{\bA\otimes \bA}\otimes \bA$ is equivalent to that of modules 
over a monad in a dualizable category, and hence is dualizable (see \cite[Lemma 1.6.3]{GKRV}).
\end{proof}

\sssec{}

Let $\bM$ be an $\bA$-module category, dualizable as a plain DG category,
and consider $$\bM^\vee\simeq \on{Funct}_{\on{cont}}(\bM,\Vect_\sfe)$$
as an $\bA$-module via the action on the source. 

\begin{cor} \label{c:dualizable plain mod}
For $\bM$ as above and another $\bA$-module category $\bN$,
we have a canonical identification
$$\on{Funct}_{\bA\mod}(\bM,\bN)\simeq \bM^\vee\underset{\bA}\otimes \bN.$$
\end{cor}

\begin{proof}

Apply \propref{p:HH semi-rigid} to $\bP=\bM^\vee\otimes \bN$.

\end{proof}

\begin{cor} \label{c:two dualities}
If an $\bA$-module category $\bM$ is dualizable as a plain DG category, then 
it is dualizable as an $\bA$-module category; moreover we have a canonical
equivalence
\begin{equation} \label{e:two dualities}
\on{Funct}_{\bA\mod}(\bM,\bA)\simeq \bM^\vee.
\end{equation}
\end{cor}

\sssec{}

Finally, we claim:

\begin{lem} \label{l:dual module over s-rigid}
Let $\bM$ be an $\bA$-module category. Then it is dualizable as such if 
\emph{and only if} it is dualizable as a plain DG category.
\end{lem}

\begin{proof}

One direction has been proved in \corref{c:two dualities}. For the other direction,
this is a general property of algebras dualizable as objects in an ambient category. 

\end{proof}

\ssec{Self-duality}

\sssec{}

Taking $\bM=\bA$ in \eqref{e:two dualities} we obtain a canonical identification
\begin{equation} \label{e:semi-rigid self dual}
\bA\simeq \bA^\vee
\end{equation}
as $\bA$-modules.

\sssec{}

We claim:

\begin{lem} \label{l:semi-rigid self dual}
The unit of the self-duality \eqref{e:semi-rigid self dual} is given by the object
$\sR_{\bA}\in \bA\otimes \bA$.
\end{lem}

\begin{proof}

Let us denote by $\phi$ the identification \eqref{e:semi-rigid self dual}. By construction,
we have a commutative diagram
$$
\CD
\bA @>{\phi}>> \bA^\vee \\
@V{\sim}VV @VV{\sim}V \\
\on{Funct}_{\bA\mod}(\bA,\bA) & & \bA^\vee \underset{\bA}\otimes \bA \\
@VVV @VVV \\
\bA^\vee \otimes \bA @>{\on{Id}\otimes\on{Id}}>> \bA^\vee \otimes \bA,
\endCD
$$
where the lower right vertical arrow is the right adjoint to the tautological projection.  The image
of $\one_\bA$ along the composite left vertical arrow is the counit in $\bA^\vee \otimes \bA$. 

\medskip

Concatenating with the commutative diagram 
$$
\CD
\bA^\vee @>{\phi^{-1}}>> \bA \\
@VVV @VV{\on{comult}_\bA}V \\
\bA^\vee \otimes \bA @>{\phi^{-1}\otimes\on{Id}}>>  \bA \otimes \bA,
\endCD
$$
we obtain a diagram
$$
\CD
\bA @>{\on{Id}}>> \bA \\
@VVV @VV{\on{comult}_\bA}V \\
\bA^\vee \otimes \bA @>{\phi^{-1}\otimes\on{Id}}>>  \bA \otimes \bA.
\endCD
$$

The evaluating both routes on $\one_\bA\in \bA$, we obtain the desired result. 
\end{proof}

\sssec{}

We now claim:

\begin{lem} \label{l:dual of mult}
With respect to the above self-dualty of $\bA$, the functor $\on{comult}_\bA$ identifies 
canonically with the dual of $\on{mult}_\bA$.  
\end{lem}

\begin{proof}

By \lemref{l:semi-rigid self dual}, the unit of the self-duality on $\bA\otimes \bA$ is given by
$$\sigma_{2,3}(\sR_\bA\otimes \sR_\bA)\in \bA\otimes \bA\otimes \bA\otimes \bA,$$
where $\sigma_{2,3}$ denotes the transposition of the two inner factors. 

\medskip

Hence, we need to establish an isomorphism
$$(\on{mult}_\bA \otimes \on{Id} \otimes \on{Id})(\sigma_{2,3}(\sR_\bA\otimes \sR_\bA))
\simeq (\on{Id}\otimes \on{comult}_\bA)(\sR_\bA).$$

The latter is a diagram chase using the fact that $\on{comult}_\bA$ is compatible with the bimodule structure.

\end{proof}

\sssec{} \label{sss:Gamma ! abs}

Let $\Gamma_{!,\bA}:\bA\to \Vect_\sfe$ denote the functor dual to the unit functor
$$\Vect_\sfe \overset{\one_\bA}\to \bA$$
with respect to the above self-duality. 

\medskip

We claim

\begin{lem}  \label{l:semi-rigid self dual counit}
The counit of the self-duality on $\bA$ is given by
$$\bA\otimes \bA \overset{\on{mult}_\bA}\longrightarrow \bA \overset{\Gamma_{!,\bA}}\to \Vect_\sfe.$$
\end{lem}

\begin{proof}

Follows by duality from \lemref{l:dual of mult} above. 

\end{proof}

The above lemma says that $(\bA,\Gamma_{!,\bA})$ is 
a \emph{Frobenius algebra object} in the symmetric monoidal category $\DGCat$. 

\sssec{} \label{sss:Gamma ! is lax mon}

Consider the lax commutative diagram
\begin{equation} \label{e:to be dual of mult}
\xy
(0,0)*+{\bA\otimes \bA}="A";
(25,0)*+{\bA}="B";
(0,-25)*+{\Vect_\sfe}="C";
(25,-25)*+{\Vect_\sfe,}="D";
{\ar@{->}^{\on{comult}_\bA} "B";"A"};
{\ar@{->}^{\on{Id}} "D";"C"};
{\ar@{->}^{\one_\bA} "D";"B"};
{\ar@{->}^{\one_\bA\otimes \one_\bA} "C";"A"};
{\ar@{=>} "C";"B"};
\endxy
\end{equation}
obtained by passing to right adjoints along the horizontal arrows in the commutative diagram
$$
\CD
\bA\otimes \bA @>{\on{mult}_\bA}>> \bA \\
@A{\one_\bA\otimes \one_\bA}AA @AA{\one_\bA}A \\
\Vect_\sfe @>{\on{Id}}>> \Vect_\sfe.
\endCD
$$

\medskip

Passing to duals in \eqref{e:to be dual of mult}, and using \lemref{l:dual of mult}, we obtain a lax commutative diagram
\begin{equation} \label{e:dual of mult}
\xy
(0,0)*+{\bA\otimes \bA}="A";
(25,0)*+{\bA}="B";
(0,-25)*+{\Vect_\sfe}="C";
(25,-25)*+{\Vect_\sfe.}="D";
{\ar@{->}^{\on{mult}_\bA} "A";"B"};
{\ar@{->}^{\on{Id}} "C";"D"};
{\ar@{->}^{\Gamma_{!,\bA}} "B";"D"};
{\ar@{->}_{\Gamma_{!,\bA}\otimes \Gamma_{!,\bA}} "A";"C"};
{\ar@{=>} "C";"B"};
\endxy
\end{equation}

This diagram, and its analogs for higher order multiplication morphisms, 
endow the functor $\Gamma_{!,\bA}$ with a (non-unital) right-lax symmetric monoidal
structure. 

\ssec{Compactness and dualizability} \label{ss:dual compact}

The material in this subsection will not be needed in the sequel. 

\medskip

Let 
$\bA$ be a semi-rigid symmetric monoidal category.  For the duration 
of this subsection we will assume that $\bA$ is compactly generated
as a plain DG category. 

\sssec{}

Let $\ba\in \bA$ be a compact object. Let
$$\BD(\ba)\in \bA$$
be its abstract dual with respect to the canonical self-duality of $\bA$, i.e., 
\begin{equation} \label{e:abs dual defn}
\CHom_{\bA}(\ba,\bb)=\Gamma_{!,\bA}(\BD(\ba)\otimes \bb).
\end{equation}

Equivalently,
\begin{equation} \label{e:abs dual via R}
\BD(\ba) \simeq (\CHom_{\bA}(\ba,-)\otimes \on{Id}_\bA)(\sR_\bA).
\end{equation}

\sssec{}

We claim:

\begin{prop} \label{p:abs vs monoidal dual}
The object $\BD(\ba)$ identifies canonically with the monoidal dual of $\ba$.
\end{prop}

\begin{proof}

We need to establish a canonical isomorphism
\begin{equation} \label{e:check dual}
\CHom_{\bA}(\ba\otimes \bb,\bc)\simeq \CHom_{\bA}(\bb,\BD(\ba)\otimes \bc),\quad \bb,\bc\in \bA.
\end{equation}

With no restriction of generality, we can assume that $\bb$ is compact. We rewrite the left-hand side
as
$$\CHom_{\bA\otimes \bA}(\ba\boxtimes \bb,\on{comult}_{\bA}(\bc))\simeq 
\CHom_{\bA\otimes \bA}(\ba\boxtimes \bb,\sR_\bA\otimes (\one_\bA\boxtimes \bc)),$$
and further as 
\begin{multline*}
\CHom_{\bA}\left(\bb, (\CHom_{\bA}(\ba,-)\otimes \on{Id}_\bA)\circ (\on{Id}_\bA\otimes (-\otimes \bc))(\sR_\bA)\right)\simeq \\
\simeq \CHom_{\bA}\left(\bb,  (-\otimes \bc)\circ (\CHom_{\bA}(\ba,-)\otimes \on{Id}_\bA)(\sR_\bA)\right).
\end{multline*}

Using \eqref{e:abs dual via R}, we rewrite the latter expression as
$$\CHom_{\bA}(\bb, \BD(\ba)\otimes \bc),$$
as desired.

\end{proof}

\begin{cor} \label{c:comp dualzable abs}
In a compactly generated semi-rigid category, compact objects are dualizable with respect to
the monoidal structure.
\end{cor}

\sssec{}

Next, we claim: 

\begin{cor} \label{c:closed under ten abs}
The subcategory of compact objects in $\bA$ is closed under the monoidal operation.
\end{cor}

\begin{proof}

Follows from the fact that (in any monoidal category) the tensor product of
a compact object and a dualizable object is compact. 

\end{proof}

\sssec{}

Let $\ba$ be again a compact object. By \propref{p:abs vs monoidal dual}, we have 
$$\CHom_{\bA}(\BD(\ba),\one_\bA)\simeq \CHom_{\bA}(\one_\bA,\ba).$$

Combining with \eqref{e:abs dual defn}, we obtain:

\begin{cor} \label{c:Gamma ! as ren}
For a compact $\ba$, we have a canonical isomorphism
$$\Gamma_{!,\bA}(\ba)\simeq \CHom_{\bA}(\one_\bA,\ba).$$
\end{cor}

\begin{rem} \label{r:Gamma ! as ren}

The last corollary means that the functor $\Gamma_{!,\bA}$ can be thought of as a renormalized
version of the functor $\CHom_{\bA}(\one_\bA,-)$ in the following sense:

\medskip

The functor $\Gamma_{!,\bA}$ is the ind-extension of the restriction of $\CHom_{\bA}(\one_\bA,-)$ 
to the subcategory of compact objects.

\medskip

Furthermore, one can show that the right-lax symmetric monoidal structure on $\Gamma_{!,\bA}$ constructed in 
\secref{sss:Gamma ! is lax mon} agrees with one induced by the right-lax symmetric monoidal structure on
$\CHom_{\bA}(\one_\bA,-)$. 

\end{rem}

\ssec{Lax vs strict compatibility}

In this subsection we let $\bA$ be a semi-rigid symmetric monoidal DG category. 

\sssec{}

Let $\bM$ be an $\bA$-module category. Consider the action functor
$$\on{act}_\bM:\bA\otimes \bM\to \bM.$$

Let 
$$\on{coact}_\bM:\bM\to \bA\otimes \bM$$
denote the functor
$$\bM\overset{\sR_\bA\otimes \on{Id}_\bM}\longrightarrow \bA\otimes \bA\otimes \bM
\overset{\on{Id}_\bA \otimes \on{act}_\bM}\longrightarrow \bA\otimes \bM,$$
i.e., this is the $\bA$-dual map of $\on{act}_\bM$, with respect to the canonical self-duality of $\bA$.

\medskip

The functor $\on{act}_\bM$ is recovered from $\on{coact}_\bM$ as
$$\bA\otimes \bM \overset{\on{Id}_\bA \otimes \on{coact}_\bM}\longrightarrow 
\bA\otimes \bA\otimes \bM \overset{\on{counit}_\bA\otimes \on{Id}_\bM}\longrightarrow \bM.$$

\sssec{}

We claim:

\begin{lem} \label{l:dual to act}
The functor $\on{coact}_\bM$ is canonically isomorphic to the right adjoint of $\on{act}_\bM$.
\end{lem}

\begin{proof}

It suffices to establish the adjunction
$$\on{act}_\bA:\bA\otimes \bA\rightleftarrows \bA:\on{coact}_\bA$$
in a way compatible with the right action of $\bA$.  

\medskip

However, in this case $\on{act}_\bA=\on{mult}_\bA$, and it easy to see that 
$\on{coact}_\bA$ identifies with $\on{comult}_\bA$.

\end{proof}

\sssec{}

We now claim:

\begin{prop} \label{p:autom strict}
Let $T:\bM_1\to  \bM_2$ be a map of $\bA$-module categories. Suppose that $T$
admits a continuous right adjoint as a functor between plain DG categories.  
Then the right-lax structure of compatibility with $\bA$-actions on $T^R$ is strict.
\end{prop}

\begin{proof}

We need to show that the diagram
$$
\CD
\bA\otimes \bM_1 @>{\on{act}_{\bM_1}}>> \bM_1  \\
@A{\on{Id}_{\bA}\otimes T^R}AA @AA{T^R}A \\
\bA\otimes \bM_2 @>{\on{act}_{\bM_2}}>> \bM_2 
\endCD
$$
commutes.

\medskip

This is equivalent to the commutation of the $\bA$-dual diagram
$$
\CD
\bM_1 @>{\on{coact}_{\bM_1}}>> \bA\otimes \bM_1  \\
@A{T^R}AA @AA{\on{Id}_{\bA}\otimes T^R}A \\
\bM_2 @>{\on{coact}_{\bM_2}}>> \bA\otimes  \bM_2. 
\endCD
$$

However, by \lemref{l:dual to act}, the latter diagram can be obtained 
from the commutative diagram
$$
\CD
\bM_1 @<{\on{act}_{\bM_1}}<< \bA\otimes \bM_1  \\
@V{T}VV @VV{\on{Id}_{\bA}\otimes T}V \\
\bM_2 @<{\on{act}_{\bM_2}}<< \bA\otimes  \bM_2. 
\endCD
$$
by passing to right adjoints.

\end{proof}

\sssec{}

Finally, we claim:

\begin{prop} \label{p:autom strict left}
Let $T:\bM_1\to  \bM_2$ be a map of $\bA$-module categories. Suppose that $T$
admits a left adjoint as a functor between plain DG categories.  
Then the left-lax structure of compatibility with $\bA$-actions on $T^L$ is strict.
\end{prop}

\begin{proof}

We wish to show that the diagram 
$$
\CD
\bA\otimes \bM_1 @>{\on{act}_{\bM_1}}>> \bM_1  \\
@A{\on{Id}_{\bA}\otimes T^L}AA @AA{T^L}A \\
\bA\otimes \bM_2 @>{\on{act}_{\bM_2}}>> \bM_2 
\endCD
$$
commutes. 

\medskip

By passing to right adjoints, this is equivalent to the commutativity of the diagram
$$
\CD
\bA\otimes \bM_1 @<{\on{coact}_{\bM_1}}<< \bM_1  \\
@V{\on{Id}_{\bA}\otimes T}VV @V{T}VV \\
\bA\otimes \bM_2 @<{\on{coact}_{\bM_2}}<< \bM_2.
\endCD
$$

However, the latter follows from the fact that $T$ is compatible with $\bA$-actions.

\end{proof}

\ssec{Persistence of semi-rigidity}

\sssec{}

Let $\bA$ be a semi-rigid symmetric monoidal category. We claim:

\begin{prop} \label{p:HH semi-rigid adj}
The tautological functor
\begin{equation} \label{e:HH semi-rigid}
\bA\underset{\bA\otimes \bA}\otimes \bA\to \bA
\end{equation}
admits a continuous right adjoint, strictly compatible with the $\bA$-bimodule structures.
\end{prop}

\begin{proof}

First, we note that once we show that the right adjoint in question is continuous, 
the strict compatibility would follow by \propref{p:autom strict}. 

\medskip

Consider the projection 
$$\bA\otimes \bA\to \bA\underset{\bA\otimes \bA}\otimes \bA.$$

It admits a right adjoint, that is continuous and conservative
(say, by \propref{p:HH semi-rigid}). Hence in order to prove that the right adjoint to
\eqref{e:HH semi-rigid} is continuous, it suffices to show that 
the right adjoint of composite functor
$$\bA\otimes \bA\to \bA\underset{\bA\otimes \bA}\otimes \bA\to \bA$$
is continuous.

\medskip

However, the above composite functor is $\on{mult}_\bA$, so the assertion follows from
the definition of semi-rigidity.

\end{proof}

\sssec{}

We now claim:

\begin{prop} \label{p:semi-rigid persists}
Let $\bA_1 \leftarrow \bA\to \bA_2$ be a diagram of semi-rigid symmetric monoidal categories. 
Then the symmetric monoidal category $\bA_1 \underset{\bA}\otimes \bA_2$ is also semi-rigid.
\end{prop}

\begin{proof}

The fact that $\bA_1 \underset{\bA}\otimes \bA_2$ is dualizable follows from \corref{c:dualizab persists}. 

\medskip

We now show that $\on{mult}_{\bA_1 \underset{\bA}\otimes \bA_2}$ admits a continuous right adjoint,
strictly compatible with the $(\bA_1 \underset{\bA}\otimes \bA_2)$-bimodule structure. Note that for the latter, it suffices to show that it is strictly
compatible with the bimodule structure with respect to $\bA_1\otimes \bA_2$, and the latter would follow by \propref{p:autom strict},
once we establish the continuity.  

\medskip

We write $\on{mult}_{\bA_1 \underset{\bA}\otimes \bA_2}$ as 
$$((\bA_1\otimes \bA_1) \otimes (\bA_2\otimes \bA_2))\underset{(\bA\otimes \bA)\otimes (\bA\otimes \bA)}\otimes (\bA\otimes \bA)\to 
(\bA_1\otimes \bA_2) \underset{\bA\otimes \bA}\otimes \bA.$$

Denote
$$\bM_1=\bA_1\otimes \bA_1,\,\, \bM_2=\bA_2\otimes \bA_2,\,\, \bM'_1=\bA_1,\,\, \bM'_2=\bA_2,\,\,\wt{\bA}:=\bA\otimes \bA.$$

So the above functor is
$$(\bM_1\otimes \bM_2) \underset{\wt{\bA}\otimes \wt{\bA}}\otimes \wt{\bA} \to
(\bM'_1\otimes \bM'_2) \underset{\wt{\bA}\otimes \wt{\bA}}\otimes \wt{\bA} \simeq 
(\bM'_1\otimes \bM'_2) \underset{\bA\otimes \bA}\otimes (\bA\underset{\wt{\bA} }\otimes \bA)\to
(\bM'_1\otimes \bM'_2) \underset{\bA\otimes \bA}\otimes \bA.$$

We claim that both arrows in the above composition admit right adjoints with the required properties. 
Indeed, for the first arrow this follows from the fact that the corresponding property of the functors
$$\bM_1\to \bM'_1 \text{ and } \bM_2\to \bM'_2$$
(by the semi-rigidity of $\bA_1$ and $\bA_2$). 

\medskip

For the second arrow, this follows from \propref{p:HH semi-rigid adj}. 

\end{proof}

\ssec{Hochschild chains of semi-rigid categories}

\sssec{}

Let $\bA$ be a semi-rigid symmetric monoidal category, and let $F_\bA$ be a symmetric monoidal endofunctor 
of $\bA$. Consider the corresponding category of Hochschild chains
$$\on{HH}_\bullet(F_\bA,\bA):= \bA\underset{\on{mult},\bA\otimes \bA,\on{mult}\circ (\on{Id}\otimes F_\bA)}\otimes \bA.$$

Note that $\on{HH}_\bullet(F_\bA,\bA)$ is also semi-rigid, by \propref{p:semi-rigid persists}. 

\sssec{}

Let now $\bA_1$ and $\bA_2$ be a pair of semi-rigid symmetric monoidal categories, and let 
$$\Phi:\bA_1\to \bA_2$$
be a symmetric monoidal functor. 

\medskip

Let $F_{\bA_1}$ and $F_{\bA_2}$ be symmetric monoidal endofunctors of $\bA_1$ and $\bA_2$,
respectively, and let us be given an isomorphism
\begin{equation} \label{e:F Phi}
F_{\bA_2}\circ \Phi\simeq \Phi\circ F_{\bA_1}.
\end{equation}

Then we obtain a functor
\begin{equation} \label{e:HH Phi}
\on{HH}_\bullet(F_{\bA_1},\bA_1)\to\on{HH}_\bullet(F_{\bA_2},\bA_2),
\end{equation}
to be denoted $\on{HH}_\bullet(F,\Phi)$. 

\sssec{}

Let $\bM_2$ be an $\bA_2$-module category.
Let $F_\bM$ denote and endofunctor of $\bM_2$ compatible with $F_{\bA_2}$ (see \cite[Sect. 3.8.2]{GKRV}).

\medskip

Assume that $\bM_2$ is dualizable as a plain DG category. Then by
\corref{c:dualizable plain mod}, $\bM_2$ is dualizable also as an $\bA_2$-module, and by \cite[Sect. 3.8.2]{GKRV}
we can attach to it an object
$$\Tr_{\bA_2}^{\on{enh}}(F_\bM,\bM_2)\in \on{HH}_\bullet(F_{\bA_2},\bA_2).$$

\sssec{}

Let $\bM_1:=\Res_\Phi(\bM_2)\in \bA_1\mod$ 
be the $\bA_1$-module category obtained from $\bM_2$ be restriction along $\Phi$.

\medskip

The data of compatibility of $F_\bM$ and $F_{\bA_2}$, combined with \eqref{e:F Phi} defines 
a data of compatibility of $F_\bM$ and $F_{\bA_1}$. Hence, we can consider the object 
$$\Tr_{\bA_1}^{\on{enh}}(F_\bM,\bM_1)\in \on{HH}_\bullet(F_{\bA_1},\bA_1).$$

The main result of this section is the following:

\begin{thm} \label{t:Tr enh abstract}
There exists a canonical isomorphism
$$\Tr_{\bA_1}^{\on{enh}}(F_\bM,\bM_1) \simeq (\on{HH}_\bullet(F,\Phi))^\vee\left(\Tr_{\bA_2}^{\on{enh}}(F_\bM,\bM_2)\right),$$
where $(\on{HH}_\bullet(F,\Phi))^\vee$ is the functor dual to $\on{HH}_\bullet(F,\Phi)$ with respect to the canonical self-dualities
on $\on{HH}_\bullet(F_{\bA_i},\bA_i)$, $i=1,2$ of \eqref{e:semi-rigid self dual} for semi-rigid categories.
\end{thm}

\begin{rem}
This theorem is a generalization of \cite[Theorem 3.10.6]{GKRV}, where $\bA_1$ and $\bA_2$ were assumed rigid. 

\medskip

Note in {\it loc. cit.}, instead of the functor $(\on{HH}_\bullet(F,\Phi))^\vee$, one considered the functor 
$(\on{HH}_\bullet(F,\Phi))^R$. However, it follows from the proof of \thmref{t:Tr enh abstract} that when
$\bA_1$ and $\bA_2$ are rigid, we have a canonical equivalence
$$(\on{HH}_\bullet(F,\Phi))^R\simeq (\on{HH}_\bullet(F,\Phi))^\vee.$$

By contrast, in the semi-rigid case, the functor $(\on{HH}_\bullet(F,\Phi))^R$ may be discontinuous
(e.g., it corresponds to $\Gamma(\CY,-)$ on a formal affine scheme).

\end{rem}

\sssec{}

Consider the particular case when $\bA_1=\Vect_\sfe$ and $F_{\bA_1}=\on{Id}$. Denote $(\bA_2,F_{\bA_2})$ by
$(\bA,F_\bA)$ and $\bM_2$ by $\bM$. We obtain:

\begin{cor} \label{c:Tr enh abstract}
There exists a canonical isomorphism
$$\Tr(F_\bM,\bM) \simeq \Gamma_{!,\on{HH}_\bullet(F_\bA,\bA)}\left(\Tr_{\bA}^{\on{enh}}(F_\bM,\bM)\right).$$
\end{cor}

\sssec{} \label{sss:proof enhanced Tr Y}

Finally, we observe that \thmref{t:enhanced Tr Y} is a particular case of \corref{c:Tr enh abstract}. 

\ssec{Proof of \thmref{t:Tr enh abstract}} \label{ss:proof Tr enh abstract}

\sssec{}

First, we recall that if $\Psi:\bA'\to \bA''$ is a monoidal functor between monoidal categories, we have an adjoint pair of 
2-functors\footnote{We use the terminology ``2-functor" for 1-morphisms in the $(\infty,3)$-category of DG 2-categories.}
$$\on{Ind}_\Psi:\bA'\mod\rightleftarrows \bA''\mod:\Res_\Psi,$$
where 
$$\on{Ind}_\Psi(\bM)=\bA''\underset{\bA'}\otimes \bM.$$

In particular, the induction 2-functor $\on{Ind}_\Psi$ always admits a right adjoint.

\medskip

Suppose now that $\bA''$ is dualizable as a left $\bA'$-module category. Then the 2-functor $\Res_\Psi$
admits a right adjoint, denoted $\on{coInd}_\Psi$,
$$\on{coInd}_\Psi(\bM)=\on{Funct}_{\bA'}(\bA'',\bM).$$

\sssec{} \label{sss:smooth and proper}

Let $\bA$ be a symmetric monoidal category. Let as assume that:

\begin{itemize}

\item $\bA$ is dualizable as an $\bA\otimes \bA$-module;

\item $\bA$ is dualizable as a plain DG category.

\end{itemize}


\medskip

The first condition implies that the 2-functor
$$\Res_{\on{mult}_\bA}:\bA\mod \to (\bA\otimes \bA)\mod$$
admits a right adjoint, and the second condition implies that the 2-functor
$$\oblv_\bA:\bA\mod \to \DGCat$$
admits a right adjoint.

\medskip

In particular, the 2-functors
\begin{equation} \label{e:unit 2-categ}
\DGCat \overset{\bA}\to \bA\mod \overset{\Res_{\on{mult}_\bA}}\longrightarrow (\bA\otimes \bA)\mod
\end{equation}
and 
\begin{equation} \label{e:counit 2-categ}
(\bA\otimes \bA)\mod \overset{\on{Ind}_{\on{mult}_\bA}}\longrightarrow \bA\mod \overset{\oblv_\bA}\longrightarrow \DGCat
\end{equation}
admit right adjoints.
 
\sssec{} 

Let $F_\bA$ be a monoidal endofunctor of $\bA$. Then by \cite[Sects. 3.3.4 and 3.7.1-3.7.2]{GKRV}, the 
2-functors \eqref{e:unit 2-categ} and \eqref{e:counit 2-categ} give rise to functors
\begin{equation} \label{e:unit HH}
\Vect_\sfe\to \on{HH}_\bullet(F_\bA,\bA)\otimes \on{HH}_\bullet(F_\bA,\bA)
\end{equation}
and 
\begin{equation} \label{e:counit HH}
\on{HH}_\bullet(F_\bA,\bA)\otimes \on{HH}_\bullet(F_\bA,\bA) \to \Vect_\sfe.
\end{equation}

Furthermore, since the functors \eqref{e:unit 2-categ} and \eqref{e:counit 2-categ} define a unit 
and a counit of a self-duality on $\bA\mod$ (in the symmetric monoidal category of DG 2-categories,
see \cite[Sect. 3.6]{GKRV}\footnote{In \cite[Sect. 3.6]{GKRV}, 
this symmetric monoidal category is denoted $\on{Morita}(\DGCat)$.}), the functors 
\eqref{e:unit HH} and \eqref{e:counit HH} define a unit 
and a counit of a self-duality on $\on{HH}_\bullet(F_\bA,\bA)$.

\sssec{}

Let $\bA$ be a semi-rigid symmetric monoidal category. Note that it automatically satisfies the conditions
of \secref{sss:smooth and proper}.  We will prove:

\begin{prop} \label{p:self dual HH semi-rigid}
The data of self-duality on $\on{HH}_\bullet(F_\bA,\bA)$ defined by the functors 
\eqref{e:unit HH} and \eqref{e:counit HH} coincides with the data of self-duality 
on $\on{HH}_\bullet(F_\bA,\bA)$ as a semi-rigid symmetric monoidal category
of \eqref{e:semi-rigid self dual}.
\end{prop}

The proof of \propref{p:self dual HH semi-rigid} will be given in \secref{ss:HH semi-rigid}. 
Let us assume it for now, and use it in order to prove \thmref{t:Tr enh abstract}.

\sssec{}

Let $\Phi:\bA_1\to \bA_2$ be a monoidal functor between monoidal
categories. Let $F_{\bA_1}$ and $F_{\bA_2}$ be monoidal endofunctors of $\bA_1$
and $\bA_2$, respectively.  We will denote by the same symbol $F_{\bA_i}$ the resulting
2-endomorphism of $\bA_i\mod$, $i=1,2$. 

\medskip

Let us be given an isomorphism 
\begin{equation} \label{e:F Phi again}
F_{\bA_2}\circ \Phi\simeq \Phi\circ F_{\bA_1}
\end{equation}
as monoidal functors.

\medskip

The isomorphism \eqref{e:F Phi again} induces an isomorphism
$$F_{\bA_1} \circ \Res_\Phi \simeq \Res_\Phi \circ F_{\bA_2},$$
and by adjunction a morphism
$$\on{Ind}_\Phi \circ F_{\bA_1}\to F_{\bA_2} \circ \on{Ind}_\Phi.$$

\medskip

Then by \cite[Sects. 3.3.4 and 3.7.1-3.7.2]{GKRV}, the 2-functor
$$\on{Ind}_\Phi:\bA_1\mod \to \bA_2\mod$$
induces a functor
$$\Tr(F,\on{Ind}_\Phi):\on{HH}_\bullet(F_{\bA_1},\bA_1)\to \on{HH}_\bullet(F_{\bA_2},\bA_2).$$

\sssec{}

Assume now that $\bA_2$ is dualizable as a left $\bA_1$-module, so that the 
2-functor $\Res_\Phi$ also admits a right adjoint. Then again by \cite[Sects. 3.3.4 and 3.7.1-3.7.2]{GKRV},
we obtain a functor 
$$\on{HH}_\bullet(F_{\bA_2},\bA_2)\to \on{HH}_\bullet(F_{\bA_1},\bA_1),$$
which we will denote by $\Tr(F,\on{Res}_\Phi)$. 

\sssec{}

Suppose now that $\bA_1$ and $\bA_2$ satisfy the
assumptions of \secref{sss:smooth and proper}. (In particular, in this case,  
$\bA_2$ is automatically dualizable as an $\bA_1$-module.)

\medskip

We claim that we have a canonical identification
\begin{equation} \label{e:dual of Tr}
\Tr(F,\on{Res}_\Phi)\simeq (\Tr(F,\on{Ind}_\Phi))^\vee,
\end{equation}
with respect to the self-dualities of \eqref{e:unit HH} and \eqref{e:counit HH}.

\medskip

Indeed, this follows by taking traces from the identification of the 2-functors 
$$\on{Res}_\Phi \simeq (\on{Ind}_\Phi)^\vee$$
with respect to the self-dualities \eqref{e:unit 2-categ} and \eqref{e:counit 2-categ}.

\sssec{}

Assume now that $\bA_1,\bA_2$ are \emph{symmetric} monoidal and semi-rigid,
that the functors $F_{\bA_1},F_{\bA_2},\Phi$ are \emph{symmetric} monoidal, 
and that the data of compatibility \eqref{e:F Phi again} respects the 
symmetric monoidal structures. 

\medskip

In order to prove \thmref{t:Tr enh abstract}, it suffices to show that we have a canonical
identification
\begin{equation} \label{e:dual of HH}
\Tr(F,\Res_\Phi)\simeq (\on{HH}_\bullet(F,\Phi))^\vee,
\end{equation}
with respect to the canonical self-dualities
$\on{HH}_\bullet(F_{\bA_i},\bA_i)$, $i=1,2$ of \eqref{e:semi-rigid self dual} for semi-rigid categories.

\medskip

However, a straightforward calculation shows that we have a canonical identification
\begin{equation} \label{e:Tr Phi as HH}
\Tr(F,\on{Ind}_\Phi)\simeq \on{HH}_\bullet(F,\Phi).
\end{equation}

Hence, the identification \eqref{e:dual of HH} follows from \eqref{e:dual of Tr}, since by
\propref{p:self dual HH semi-rigid}, the above
self-dualities equal those given by \eqref{e:unit HH} and \eqref{e:counit HH}, 

\qed[\thmref{t:Tr enh abstract}]

\ssec{Proof of \propref{p:self dual HH semi-rigid}} \label{ss:HH semi-rigid}

\sssec{}

It suffices to show that the unit functor \eqref{e:unit HH} identifies canonically with 
$$\Vect_\sfe \to \on{HH}_\bullet(F_\bA,\bA) \overset{\on{comult}_{\on{HH}_\bullet(F_\bA,\bA)}}\longrightarrow 
\on{HH}_\bullet(F_\bA,\bA)\otimes \on{HH}_\bullet(F_\bA,\bA).$$

For this, it suffices to show that the latter functor is obtained from \eqref{e:unit 2-categ} by taking traces. 

\medskip

This is obvious for
the first arrow, i.e., 
$$\DGCat\to \bA\mod$$
(see \eqref{e:Tr Phi as HH}). 

\sssec{}

For the second arrow, i.e., 
$$\bA\mod \overset{\Res_{\on{mult}_\bA}}\longrightarrow (\bA\otimes \bA)\mod$$
we argue as follows:

\medskip

Note that the trace $\Tr(F, \on{Ind}_{\on{mult}_\bA})$ of the 2-functor
$$(\bA\otimes \bA)\mod\overset{\on{Ind}_{\on{mult}_\bA}}\longrightarrow \bA\mod$$
identifies with $\on{mult}_{\on{HH}_\bullet(F_\bA,\bA)}$, again by  \eqref{e:Tr Phi as HH}.

\medskip

Hence, it suffices to show that the functor
$$\Tr(F, \Res_{\on{mult}_\bA}):\on{HH}_\bullet(F_\bA,\bA) \to \on{HH}_\bullet(F_\bA,\bA)\otimes \on{HH}_\bullet(F_\bA,\bA)$$
identifies with the right adjoint of $\Tr(F, \on{Ind}_{\on{mult}_\bA})$. 

\sssec{} 

Recall the setting of \cite[Sect. 3.9]{GKRV}. Let $\fT_1$ and $\fT_2$ be a pair of DG 2-categories, each equipped
with an endofunctor $F_i$, $i=1,2$. Consider the corresponding categories
$$\Tr(F_1,\fT_1) \text{ and } \Tr(F_2,\fT_2).$$

Let $\bPhi:\fT_1\to \fT_2$ be a 2-functor that admits a right adjoint. Let us be given a natural transformation
\begin{equation} \label{e:alpha nat}
\alpha:\bPhi\circ F_1\to F_2\circ \bPhi.
\end{equation}

Then by \cite[Sect. 3.3.4]{GKRV}, we obtain a functor
$$\Tr(F,\bPhi): \Tr(F_1,\fT_1)\to \Tr(F_2,\fT_2).$$

\medskip

Suppose now that we are given two 2-functors $\bPhi',\bPhi'':\fT_1\to \fT_2$ as above, and let us be given 
a 2-morphism 
$$\beta:\bPhi'\to \bPhi'',$$
equipped with a natural 3-transformation $\gamma$ from 
$$\bPhi'\circ F_1 \overset{\alpha'}\to F_2\circ \bPhi' \overset{\beta}\to  F_2\circ \bPhi''$$
to
$$\bPhi'\circ F_1 \overset{\beta}\to \bPhi''\circ F_1  \overset{\alpha''}\to F_2\circ \bPhi''.$$

Finally assume that the 2-morphism $\beta$ admits a right adjoint. Then, by \cite[Sect. 3.9.4]{GKRV}, we obtain a natural 
transformation
\begin{equation} \label{e:Tr on 2-morph}
\Tr(F,\bPhi') \overset{\Tr(F,\beta)}\longrightarrow \Tr(F,\bPhi'').
\end{equation}

\sssec{} \label{sss:Tr adj}

Let $\bPhi:\fT_1\to \fT_2$ be as above, and suppose that the 2-functor $\bPhi$ is the left adjoint of a 
2-functor
$$\bPsi:\fT_2\to \fT_1,$$
equipped with an \emph{isomorphism}
\begin{equation} \label{e:alpha iso}
\bPsi\circ F_2\to F_1\circ \bPsi,
\end{equation}
so that \eqref{e:alpha nat} arises from \eqref{e:alpha iso} by adjunction.

\medskip

Note that in this case, the unit and the counit of the adjunction
\begin{equation} \label{e:Psi Phi unit}
\on{Id}_{\fT_1}\to \bPsi\circ \bPhi \text{ and } \bPhi\circ \bPsi\to \on{Id}_{\fT_2}
\end{equation}
automatically come equipped with the data of 3-morphisms $\gamma$ as above,
which are in fact \emph{isomorphisms}. 

\medskip

Assume that $\bPsi$ itself admits and a right adjoint (which is a 1-morphism) and that the
2-morphisms \eqref{e:Psi Phi unit} admit right adjoints (which are also 2-morphisms). We obtain that 
the natural transformation \eqref{e:Tr on 2-morph} applied to \eqref{e:Psi Phi unit} defines an adjunction
between $\Tr(F,\bPhi)$ and $\Tr(F,\bPsi)$. 

\sssec{}

We apply the material of \secref{sss:Tr adj} to the situation when
$$\fT_1:=\bA_1\mod,\,\, \fT_2:=\bA_2\mod, \,\, \bPhi:=\on{Ind}_\Phi,\,\, \bPsi:=\on{Res}_\Phi,$$
for a monoidal functor $\Phi:\bA_1\to \bA_2$. The datum of \eqref{e:alpha iso} is supplied by
\eqref{e:F Phi again}.

\medskip

Assume that $\bA_2$ is dualizable as an $\bA_1$-module. We obtain that the functors
$$\Tr(F,\on{Ind}_\Phi):\on{HH}_\bullet(F_{\bA_1},\bA_1)\rightleftarrows  \on{HH}_\bullet(F_{\bA_2},\bA_2):\Tr(F,\on{Res}_\Phi)$$
are an adjoint pair if the following conditions are satisfied:

\begin{itemize}

\item The functor $\Phi:\bA_1\to \bA_2$ admits a right adjoint as a map of $\bA_1$-bimodule categories;

\item The functor $\bA_2\underset{\bA_1}\otimes \bA_2 \to \bA_2$ 
admits a right adjoint as a map of $\bA_2$-bimodule categories.

\end{itemize}

\sssec{}

We apply this to $\bA_1=\bA\otimes \bA$, $\bA_2=\bA$ and $\Phi=\on{mult}_\bA$.

\medskip

Now, the existence of the right adjoint to $\on{mult}_\bA$ 
follows from the semi-rigidity condition. 

\medskip

The existence of the right adjoint to 
$$\bA\underset{\bA\otimes \bA}\otimes \bA\to \bA$$ 
follows from \propref{p:HH semi-rigid adj}.

\qed[\propref{p:self dual HH semi-rigid}]

\section{The dimension of the global nilpotent cone} \label{s:glob Nilp}

In this section we prove that under certain the restrictions on $\on{char}(k)$, the global
nilpotent cone
$$\Nilp\subset T^*(\Bun_G),$$
viewed as a \emph{classical} algebraic stack has dimension equal to $\dim(\Bun_G)=\dim(G)\cdot (g-1)$. 

\ssec{The Faltings-Ginzburg argument}

We first explain Faltings' proof that $\Nilp$ is isotropic (see \cite[Theorem III.2]{Fa}) which was 
conceptualized by V.~Ginzburg in \cite{Gi}. This argument is valid in characteristic $0$, and requires
a certain assumption in positive characteristic. This assumption will be satisfied if $\on{char}(k)$
is ``very good", see \secref{sss:JM}. 

\sssec{} \label{sss:nilp in par}

We first formulate the assumptionson $\fg$ and $\on{char}(k)$ that we need for the Faltings-Ginzburg
proof:

\begin{itemize}

\item The Lie algebra $\fg$ admits a non-degenerate $G$-invariant bilinear form;

\item For every (not necessarily algebraically) closed field extension
$\wt{k}\supset k$ and a nilpotent element $n\in \wt{k}\underset{k}\otimes \fg$, 
there exists a parabolic $\wt{P}\subset G$ defined over 
$\wt{k}$ such that $n$ belongs to the Lie algebra of its unipotent radical.

\end{itemize}

\medskip

When $\on{char}(k)=0$, the above conditions are automatically satisfied. Indeed, for the second condition,
we can take $P$ to be a Borel subgroup. Indeed, the element
$n$ generates a copy of $\BG_a\subset G$; take an arbitrary point on the flag variety, consider 
the resulting map $\BG_a\to G/B$, and complete it to a map $\BP^1\to G/B$. Then the image of
$\infty\in \BP^1$ is $\BG_a$-invariant, and hence the Lie algebra of the corresponding Borel subgroup
contains $n$. 

\medskip

When  $\on{char}(k)>0$, we were able to prove the validity of this assumption for
``very good" characteristics, using (a variant of) the Jacobson-Morozov theory, see \secref{sss:JM}. 

\sssec{}

Let $\CY$ be a smooth algebraic stack. We will regard $T^*(\CY)$ as a \emph{classical} algebraic stack,
see \secref{sss:cotan smooth stack} below. 

\medskip

For a Zariski locally closed subset $\CN \subset T^*(\CY)$, and a smooth map
$S\to \CY$, where $S$ is a scheme, denote by $\CN_S\subset T^*(S)$
the image of 
$$\CN\underset{\CY}\times S\subset T^*(\CY)\underset{\CY}\times S$$
under the codifferential map
$$T^*(\CY)\underset{\CY}\times S\to T^*(S).$$

We say that $\CN$ is half-dimensional if $\dim(\CN_S)\leq \dim(S)$ for every $S$ as above
(equivalently, for a collection of affine schemes $S$ that smoothly cover $\CY$). 

\medskip

We will say a Zariski locally closed subset of $Z\subset T^*(S)$ is isotropic if for every irreducible
component of $Z$, some non-empty smooth open subset $Z^\circ$ of this irreducible component 
is isotropic (i.e., the symplectic form vanishes on its tangent spaces). 

\medskip

We will say that a Zariski locally closed subset of $Z\subset T^*(S)$ is \emph{strongly} 
isotropic if for every point of $z\in Z$, the symplectic form vanishes on $H^0(T_z(Z))$. 

\medskip

We will say that a Zariski locally closed subset $\CN \subset T^*(\CY)$ is isotropic (resp., strongly isotropic) if $\CN_S$ is 
isotropic (resp., strongly isotropic) for every $S$ as above (equivalently, for a collection of schemes $S$ that smoothly cover $\CY$). 

\medskip

Clearly, if $\CN$ is isotropic then it is half-dimensional.

\begin{rem} \label{r:semi strong isotrop}
One can also consider a notion intermediate between isotropic and strongly isotropic: one can require that for
any smooth scheme $Z'$ mapping to $Z$, the pullback of the symplectic form to $Z'$ vanishes.
\end{rem}

\sssec{}

Let $f:\CY_1\to \CY_2$ be a schematic map between smooth stacks. Let $K_f$ denote the kernel of the codifferential, 
i.e., 
$$\on{Ker}\left(T^*(\CY_2)\underset{\CY_2}\times \CY_1\overset{df^*}\longrightarrow T^*(\CY_1)\right),$$
viewed as a Zariski-closed subset in $T^*(\CY_2)\underset{\CY_2}\times \CY_1$.

\medskip

Let $\CN\subset T^*(\CY_2)$ be a Zariski-closed subset. We have the following assertion:

\begin{prop} \label{p:isotrop}
Assume that the projection
$$T^*(\CY_2)\underset{\CY_2}\times \CY_1\to T^*(\CY_2)$$
maps $K_f$ to $\CN$, such that the following holds:

\medskip

\noindent There exists an open dense subset $\CN^\circ\subset \CN$ such that for every 
(not necessarily algebraically closed) field extension $k'\supset k$ and every $k'$-point
$n$ of $\CN^\circ$, there exists a finite separable field extension $k''\supset k'$ and a lift
of $n$ to a $k''$-point of $K_f$.

\medskip

Then $\CN$ is isotropic.

\end{prop}

\begin{rem}
If $\on{char}(k)=0$, the condition in \propref{p:isotrop} amounts to the requirement that the map
$K_f\to \CN$ be dominant.
\end{rem}

\begin{proof}

By base change, we can assume that $\CY_2=Y_2$ is a scheme. Then, since $f$ is schematic, $\CY_1=Y_1$ is a scheme
as well.  Consider the product
$$T^*(Y_1)\times T^*(Y_2)\simeq T^*(Y_1\times Y_2)$$
with its symplectic structure. 

\medskip

We have a natural embedding
$$T^*(Y_2)\underset{Y_2}\times Y_1 \hookrightarrow T^*(Y_1)\times T^*(Y_2),$$
whose image is a smooth Lagrangian. We can view $K_f$ as the intersection of this Lagrangian
with $\{0\}\times T^*(Y_2)$. 

\medskip

Hence, $K_f$ is \emph{strongly} isotropic as a subset of $T^*(Y_1)\times T^*(Y_2)$. We wish to show that
its image along the projection
$$T^*(Y_1)\times T^*(Y_2)\to T^*(Y_2)$$
is isotropic. 

\medskip

More precisely, let $Z$ be a smooth open subset of an irreducible component of $\CN^\circ$. We wish to show
that some non-empty open subset $Z^\circ$ of any such $Z$ is isotropic.

\medskip

Let $k'$ be the field of fractions of $Z$. By assumption, there exists a finite separable field extension $k''\supset k'$
and a $k''$-point of $K_f$ such that
$$\Spec(k'')\to K_f \to \CN$$
equals 
$$\Spec(k'')\to \Spec(k')\to Z \to \CN.$$

Hence, we can find a (non-empty) scheme $\wt{Z}$ equipped with an
\'etale map $\wt{Z}\to Z$ together with a lift of this map to a map
 $$\wt{Z} \to K_f.$$
 
Let $Z^\circ \subset Z$ denote the image of the map $\wt{Z}\to Z$. 
We claim that the symplectic form vanishes on $Z^\circ$.
 
 \medskip
 
Indeed, for every $z\in Z^\circ$, let $\wt{z}$ be its lift to a point of $\wt{Z}$.
Let $(w,z)$ denote the resulting point of $K_f\subset   T^*(Y_1)\times T^*(Y_2)$.

\medskip

By \'etaleness, the map $T_{\wt{z}}(\wt{Z})\to T_z(Z)$ is surjective. Hence, 
the map
$$H^0(T_{(w,z)}(K_f))\to T_z(Z)$$
is surjective. 

\medskip

Now, the restriction of the symplectic form on $T^*(Y_1)\times T^*(Y_2)$
along
$$H^0(T_{(w,z)}(K_f)) \to T_{(w,z)}(T^*(Y_1)\times T^*(Y_2))$$
equals the restriction of the symplectic form on $T^*(Y_2)$ along
$$H^0(T_{(w,z)}(K_f)) \to T_z(T^*(Y_2)).$$

Indeed, the two restrictions are already equal on $T_{(w,z)}(\{0\}\times T^*(Y_2))$.

\end{proof} 

\begin{rem}
The above proof shows that $\CN$ is actually semi-strongly isotropic, see Remark \ref{r:semi strong isotrop} for what this means.
\end{rem}
 
\sssec{}

We are now ready to prove that $\Nilp$ is half-dimensional. We are going to apply \propref{p:isotrop} in the following situation:

\medskip

We take $\CY_2=\Bun_G$ and $\CN=\Nilp$. We we take $\CY_1$ to be the union of $\Bun_P$
over the set of standard parabolics $P\subset G$. 

\medskip

Using an invariant form on $\fg$, we identify $T^*(\Bun_G)$ with the stack that classifies \
pairs $(\CP_G,A)$, where $\CP_G$ is a $G$-bundle on $X$ and $A$ is a section of $\fg_{\CP_G}\otimes \omega_X$.

\medskip

Then, for a given parabolic, the stack $K_f$ classifies pairs $(\CP_P,A)$, where $\CP_P$ is a $P$-bundle on
$X$, and  $A$ is a section of $\fn(P)_{\CP_G}\otimes \omega_X$.

\medskip

Clearly, the projection $K_f\to T^*(\Bun_G)$ has its image contained in $\Nilp$. Thus, in order to prove that
$\Nilp$ is isotropic, we have to show that the condition of \propref{p:isotrop} is satisfied. 

\medskip

Let $k'\supset k$ be a field extension. Let $(\CP'_G,A')$ be a $k'$-point of $T^*(\Bun_G)$ with $A'$ nilpotent. 
By \cite{DS}, there exists a separable field extension $k''\supset k'$ such that the pullback $\CP''_G$ of $\CP'_G$
to the curve 
$$X'':=X\underset{\Spec(k)}\times \Spec(k'')$$
can be trivialized at the generic $\eta$ point of $X''$. Let $\wt{k}$ denote the field of rational functions on $X''$. 
Denote $A'':=A'|_{X''}$ and $\wt{A}:=A''|_\eta$.  

\medskip

Up to trivializing $\omega_{X''}$ and $\CP''_G$ generically, we can think of $\wt{A}$ as a nilpotent element in 
$\wt{k}\underset{k}\otimes \fg$. Let $\wt{P}$ be a parabolic defined over $\wt{k}$ such that $\wt{A}\in \fn(\wt{P})$. 
It exists by the assumption in \secref{sss:nilp in par}.

\medskip

Let $P$ be the standard parabolic conjugate to
$\wt{P}$. Then we can think of $\wt{P}$ as a reduction $\CP''_P$ of $\CP''_G$ to $P$ at $\eta$, so that
$A''$ is a section of $\fn(P)_{\CP''_P}\otimes \omega_{X''}$. 

\medskip

By the valuative criterion, the reduction $\CP''_P$ of $\CP''_G$ (uniquely) extends to the entire $X''$,
and the section $A''$ belongs to $\fn(P)_{\CP''_P}\otimes \omega_{X'}$ (because it does so generically).

\medskip

The resulting pair $(\CP''_P,A'')$ is the sought-for lift of $(\CP''_G,A'')$ to a $k''$-point of $K_f$.

\qed[Isotropy of $\Nilp$]

\ssec{Adaptation of the Jacobson-Morozov theory} \label{ss:JL}

In this subsection, we will assume that the characteristic of $k$ is ``very good" for $\fg$
(this excludes very small primes for every isomorphism class of root data of $G$). 

%

\sssec{}

We first summarize the results that of \cite{Pre} that we will need.

\medskip

Let $n$ be a nilpotent element of $\fg$. Then there exists a homomorphism $\lambda:\BG_m\to G$ with the following properties: 

\medskip

Denote by
$$\fg=\underset{i}\oplus\, \fg_i$$
the weight decomposition of $\fg$ for the induced adjoint action of $\BG_m$. Denote by
$$\fg^i:=\underset{j\geq i}\oplus\,\fg_j$$
the corresponding filtration.

\medskip

We have:

\begin{itemize}

\item $n\in \fg_2$;

\item $\fg^0=:\fp$ is a Lie algebra of a parabolic subgroup (to be denoted $P$);

\item The map $\fp\overset{\on{ad}_n}\longrightarrow \fg^2$ is surjective.

\item $\fz_\fg(n)\subset \fp$ and $Z_G(n)\subset P$. 

\end{itemize}


\medskip

The above is (part of) the content of \cite[Theorem A]{Pre}.

\sssec{}

Note that the surjectivity of the map 
\begin{equation} \label{e:ad surj}
\fp\overset{\on{ad}_n}\longrightarrow \fg^2
\end{equation} 
implies that the $\on{Ad}_P$-orbit of $n$, denoted 
$\overset{\circ}\fg{}^2$, is a Zariski open subset in $\fg^2$. 

\sssec{}

Note also that the fact that $\fz_n(\fg)\subset \fp$ implies that the map
\begin{equation} \label{e:ad inj prel}
\on{ad}_n:\fg_{-1}\to \fg_1
\end{equation}
is injective. 

\medskip

Since $\dim(\fg_{-1})=\dim(\fg_1)$ (the bilinear form on $\fg$ restricts to a perfect pairing on
$\fg_{-1}\otimes \fg_1$), we obtain that \eqref{e:ad inj prel} is an isomorphism.

\sssec{}

Let $\bO\subset \fg$ denote the orbit of $n$ under the adjoint action. Denote by $Y$ the partial flag 
variety $Y=G/P$.

\medskip

Let $\wt\bO\to Y$ be the total space of a $G$-equivariant vector bundle over $Y$, whose fiber 
over $P\in Y$ is the space $\fg^2$, equipped with the natural action of $P$. 

\medskip

Let 
$$\wt\bO{}^\circ\subset \wt\bO$$
be the $G$-invariant open subscheme whose fiber over $P$ is
$$\overset{\circ}\fg{}^2\subset \fg^2.$$

The natural projection
$$\wt\bO\to\fg$$
restricts to a map
\begin{equation} \label{e:orb Spr}
\wt{\bO}^\circ\to \bO.
\end{equation}

\begin{lem} \label{l:orb Spr}
The map \eqref{e:orb Spr} is an isomorphism.
\end{lem}

\begin{proof}

Since $G$ acts transitively on both $\wt{\bO}^\circ$ and $\bO$, it suffices to show that
the stabilizers are equal. 

\medskip

However, this follows from the fact that $Z_G(n)\subset P$.

\end{proof}

\sssec{}   \label{sss:JM}

We are now ready to prove the property from \secref{sss:nilp in par}.

\medskip

The assumption that $\on{char}(k)$ is very good implies that $G$ acts on the nilpotent 
cone of $\fg$ with finitely many orbits. Hence, the nilpotent cone of $\fg$ is a finite union of 
its locally closed subsets $\bO$.

\medskip

Hence, given a field extension $\wt{k}\supset k$ and a nilpotent element $n\in \fg(\wt{k})$, 
there exists a nilpotent $G$-orbit $\bO$ defined over $k$ such that $n\in \bO(\wt{k})$.

\medskip

Now the assertion follows from \lemref{l:orb Spr} by taking $\wt{k}$-points.

\ssec{The Beilinson-Drinfeld argument}

In this subsection we will give another proof of the fact that $\Nilp$
is half-dimensional (under the same assumptions as above). 

\medskip

We will explicitly
write $\Nilp$ as a union of algebraic stacks of dimension $\leq \dim(G)\cdot (g-1)$. 

\medskip

This argument is wholly borrowed from \cite[Sect. 2.10.3]{BD2}. We include it here
for completeness. 

\sssec{} 

Let $\bO$ be a nilpotent conjugacy class. Let $\Nilp_\bO\subset \Nilp$ be the locally
closed substack consisting of pairs $(\CP_G,A)$, where $A$ generically belongs to $\bO$.

\medskip

By the same reasoning as in \secref{sss:JM} above, we have
$$\Nilp=\underset{\bO}\cup\, \Nilp_\bO.$$

\medskip

Therefore, it is enough show that each $\Nilp_\bO$ has dimension $\leq \dim(G)\cdot (g-1)$.

\sssec{}

Consider the algebraic stack $\bMaps(X,\wt\bO_{\omega_X}/G)$, where $\wt\bO_{\omega_X}$
denotes the twist of the constant bundle over $X$ with fiber $\wt\bO$ by $\omega_X$ viewed as
a $\BG_m$-torsor, using the $\BG_m$-action on $\wt\bO$ by fiber-wise dilations.

\medskip

Let 
$$\bMaps(X,\wt\bO_{\omega_X}/G)^\circ \subset \bMaps(X,\wt\bO_{\omega_X}/G)$$
be the open substack
consisting of maps that generically land in $\wt{\bO}^\circ$.

\medskip

The isomorphism of \lemref{l:orb Spr} and the valuative criterion imply that the 
map $$\bMaps(X,\wt\bO_{\omega_X}/G)^\circ\to \Nilp_\bO,$$
is bijective on geometric points. 

\medskip

Hence, it suffices to show that $\bMaps(X,\wt\bO_{\omega_X}/G)^\circ$ has dimension $\leq \dim(G)\cdot (g-1)$.

\medskip

We will show that $\bMaps(X,\wt\bO_{\omega_X}/G)^\circ$ is a smooth algebraic stack and that 
its (stacky) tangent spaces at $k$-points have Euler characteristics
$\leq \dim(G)\cdot (g-1)$. 

\sssec{}

Let $(\CP_P,A)$ be a $k$-point of $\bMaps(X,\wt\bO_{\omega_X}/G)^\circ$. Its stacky tangent space 
is given by
$$\Gamma(X,E_{-1}\overset{\on{ad}_A}\to E_0),$$
where $E_{-1}=\fg^0_{\CP_P}$ and $E_0=\fg^2_{\CP_P}\otimes \omega_X$.

\medskip

Note that since $A$ generically belongs to $\overset{\circ}\fg{}^2$ and 
the map \eqref{e:ad surj} is surjective, we obtain that the map
$$E_{-1}\overset{\on{ad}_A}\to E_0$$
is generically surjective, i.e., its cokernel is a torsion sheaf on $X$. 

\medskip

This implies that $T^*_{(\CP_P,A)}(\bMaps(X,\wt\bO_{\omega_X}/G)^\circ)$ is acyclic is cohomological 
degrees $>0$. This implies that $\bMaps(X,\wt\bO_{\omega_X}/G)^\circ$ is smooth. 

\sssec{}

We have
$$\chi\left(T^*_{(\CP_P,A)}(\bMaps(X,\wt\bO_{\omega_X}/G)^\circ)\right)=\chi(\Gamma(X,E_0))-\chi(\Gamma(X,E_{-1})).$$ 

Using the non-degenerate $G$-invariant form on $\fg$, we identify $\fg^2_{\CP_P}$ with the dual
vector bundle of $(\fg/\fg^{-1})_{\CP_G}$. Hence, by Serre duality
$$\chi(\Gamma(X,E_0))=-\chi(\Gamma(X,(\fg/\fg^{-1})_{\CP_G})).$$

Hence, we obtain 
$$\chi(T^*_{(\CP_P,A)}(\bMaps(X,\wt\bO_{\omega_X}/G)^\circ))=-\chi(\Gamma(X,(\fg/\fg^{-1})_{\CP_P}))-\chi(\Gamma(X,\fg^0_{\CP_P}))=$$
$$=-\chi(\Gamma(X,\fg_{\CP_P})) + \chi(\Gamma(X,(\fg^{-1}/\fg^0)_{\CP_P})=\dim(\Bun_G)+\chi(\Gamma(X,(\fg^{-1}/\fg^0)_{\CP_P}).$$

\medskip

It remains to show that $\chi(\Gamma(X,(\fg^{-1}/\fg^0)_{\CP_P})\leq 0$.

\sssec{}

Again, by Serre duality, we have
$$\chi(\Gamma(X,(\fg^{-1}/\fg^0)_{\CP_P})=-\chi(\Gamma(X,(\fg^{0}/\fg^2)_{\CP_P}\otimes\omega_X).$$

Hence,
$$2\chi(\Gamma(X,(\fg^{-1}/\fg^0)_{\CP_P})=\chi(\Gamma(X,(\fg^{-1}/\fg^0)_{\CP_P})-\chi(\Gamma(X,(\fg^{0}/\fg^2)_{\CP_P}\otimes\omega_X)=$$
$$=-\chi\left(\Gamma(X,(\fg^{-1}/\fg^0)_{\CP_P} \overset{\on{ad}_A}\to (\fg^{0}/\fg^2)_{\CP_P}\otimes\omega_X\right).$$

Hence, it is enough to show that
$$\chi\left(\Gamma(X,(\fg^{-1}/\fg^0)_{\CP_P} \overset{\on{ad}_A}\to (\fg^{0}/\fg^2)_{\CP_P}\otimes\omega_X\right) \geq 0.$$

Note, however, that the map
$$\on{ad}_n:\fg^{-1}/\fg^0\to \fg^1/\fg^2$$
is an isomorphism, since \eqref{e:ad inj prel} is an isomorphism.

\medskip

Therefore, the map
$$(\fg^{-1}/\fg^0)_{\CP_P} \overset{\on{ad}_A}\to (\fg^{0}/\fg^2)_{\CP_P}\otimes\omega_X$$
is generically an isomorphism. Hence, it is injective and its cokernel is torsion. Hence, the Euler characteristic
of its cone is non-negative. 

\qed

\section{Ind-constructible sheaves on schemes} \label{s:shvs on sch} 

Algebro-geometric objects 
in this section will be quasi-compact schemes over $k$, assumed almost\footnote{Since we are dealing with $\Shv(-)$,
we lose nothing by only considering classical schemes, i.e., derived algebraic geometry over $k$ will play no role.}
of finite type. Let $\Shv(-)^{\on{constr}}$ be one of the sheaf-theoretic contexts from \secref{sss:Shv}.

\ssec{The left completeness theorem} \label{ss:shvs on schemes}

\sssec{}

Recall that for a (quasi-compact) scheme $Y$ we define
$$\Shv(Y):=\on{Ind}(\Shv(Y)^{\on{constr}}).$$

The goal of this subsection is to prove \thmref{t:Shv left-complete}. The proof will be obtained
as a combination of the following two statements:

\begin{thm} \label{t:der perv}
The canonical functor
$$D^b(\on{Perv}(Y))\to \Shv(Y)^{\on{constr}}$$
is an equivalence.
\end{thm}

\begin{thm} \label{t:der left complete}
Let $\CA$ be a small abelian category of finite cohomological dimension.
Then the DG category $\on{Ind}(D^b(\CA))$ is left-complete in its t-structure.
\end{thm}

\thmref{t:der perv} is a theorem of A.~Beilinson, and it is proved in \cite{Be1}. The rest of this subsection
is devoted to the proof of \thmref{t:der left complete}.

\sssec{}

Recall that for a DG category $\bC$ equipped with a t-structure, we denote by $\bC^\wedge$ its
left completion, i.e.,
$$\bC^\wedge:=\underset{n}{\on{lim}}\, \bC^{\geq -n}.$$

We will think of objects of $\bC^\wedge$ as compatible collections
$$\{\bc^n\in \bC^{\geq -n}\} .$$

\medskip

We have the tautological functor
\begin{equation} \label{e:to left comp}
\bC\to \bC^\wedge, \quad \bc \mapsto \{\tau^{\geq -n}(\bc)\in \bC^{\geq -n}\}
\end{equation} 
and its right adjoint given by
\begin{equation} \label{e:from left comp}
\{\bc^n\in \bC^{\geq -n}\} \mapsto \underset{n}{\on{lim}}\, \bc^n,
\end{equation} 
where the limit is taken in $\bC$.

\sssec{}

We shall say that $\bC$ has \emph{convergent Postnikov towers} if \eqref{e:to left comp} is fully
faithful. Equivalently, if for $\bc\in \bC$, the natural map
$$\bc\to \underset{n}{\on{lim}}\, \tau^{\geq -n}(\bc)$$
is an isomorphism.

\sssec{}

We shall say that an object $\bc\in \bC$ has cohomological dimension $\leq n$ if
$$\Hom_\bC(\bc,\bc')=0 \text{ for all } \bc'\in \bC^{<-n}.$$

\medskip

We claim:

\begin{prop} \label{p:Postnikov}
Let $\bC$ be generated by compact objects of finite cohomological dimension.
Then $\bC$ has convergent Postnikov
towers. Furthermore, the right adjoint to \eqref{e:to left comp} is continuous.
\end{prop}

\begin{proof}

It is enough to show that for every $\bc_0\in \bC^c$, the functor
\begin{equation} \label{e:Hom into Postnikov}
\bC^\wedge \overset{\text{\eqref{e:from left comp}}}\longrightarrow \bC 
\overset{\CHom_\bC(\bc_0,-)}\longrightarrow \Vect_\sfe
\end{equation} 
is continuous, and that its precomposition with \eqref{e:to left comp} is isomorphic to
$\CHom_\bC(\bc_0,-)$.

\medskip

Now, the functor \eqref{e:Hom into Postnikov} sends 
$$\{\bc^n\in \bC^{\geq -n}\} \in \bC^\wedge$$ 
to 
$$\underset{n}{\on{lim}}\, \CHom_\bC(\bc_0,\bc^n),$$
while in the above limit, each individual cohomology group stabilizes due to the
assumption on $\bc_0.$

\medskip

This implies both claims.

\end{proof} 

\sssec{}

Applying \propref{p:Postnikov} to $\on{Ind}(D^b(\CA))$, we obtain that it has convergent Postnikov
towers. It remains to show that the functor \eqref{e:from left comp} is fully faithful. 

\medskip

We claim:

\begin{lem} \label{l:when left complete}
The functor \eqref{e:from left comp} is fully faithful if and only if the following condition
is satisfied: 

\medskip

\noindent For a family of objects of $\bC^{\leq 0}$ indexed by $\BN$
$$n\mapsto \bc_n$$
such that for every $N$ the family $\tau^{\geq -N}(\bc_n)$ stabilizes, we have
$$\underset{n}{\on{lim}}\, \bc_n\in \bC^{\leq 0}.$$

\end{lem}

\begin{proof} 

The ``only" if direction is obvious. We now prove the ``if" direction.

\medskip

Let $\{\bc^m\in \bC^{\geq -m}\}$ be an object of $\bC^\wedge$, and set
$$\bc:=\underset{m}{\on{lim}}\, \bc^m.$$

We need to show that for any $n$, the map
$$\tau^{\geq -n}(\bc)\to \bc^n$$
is an isomorphism. 

\medskip

We have a fiber sequence
$$\underset{m\geq n+1}{\on{lim}}\, \tau^{<-n}(\bc^m) \to \bc \to \underset{m\geq n+1}{\on{lim}}\, \tau^{\geq -n}(\bc^m),$$
where the left-most term belongs to $\bC^{<-n}$ by assumption, and the right-most term is $\bc^{n}$,
since the corresponding inverse family is constant with value $\bc^{n}$. 

\medskip

From here, we obtain that
$$\tau^{\geq -n}(\bc)\to \bc^n$$
is an isomorphism. 

\end{proof}

\sssec{}

Let us show that \lemref{l:when left complete} is applicable to $\on{Ind}(D^b(\CA))$. Let 
$$n\mapsto\bc_n$$
be a family of objects in $\on{Ind}(D^b(\CA))$ as in \lemref{l:when left complete}. 
Let $d$ be the cohomological dimension of $\CA$. Considering the fiber sequence
$$\tau^{\leq -d}(\bc_n) \to \bc_n\to \tau^{>-d}(\bc_n)$$
and taking into account that the family $n\mapsto  \tau^{> -d}(\bc_n)$ stabilizes, we 
obtain that we can 
assume that $\bc_n\in \on{Ind}(D^b(\CA))^{\leq -d}$.

\medskip

We claim that for any $a\in \CA$
$$\CHom_{\on{Ind}(D^b(\CA))}(a,\underset{n}{\on{lim}}\, \bc_n)\simeq 
\underset{n}{\on{lim}}\, \CHom_{\on{Ind}(D^b(\CA))}(a,\bc_n)\in \Vect_\sfe^{\leq 0}.$$

Indeed, in the family
$$n\mapsto \CHom_{\on{Ind}(D^b(\CA))}(a,\bc_n),$$
the terms belong to $\Vect^{\leq 0}_\sfe$, and for any $N$, the family 
$$n\mapsto \tau^{\geq -N}(\CHom_{\on{Ind}(D^b(\CA))}(a,\bc_n))$$
stabilizes. Hence the limit belongs to $\Vect_\sfe^{\leq 0}$ by the ``only if" direction in
\lemref{l:when left complete}, applied to $\Vect_\sfe$. 

\medskip

Now the fact that $\underset{n}{\on{lim}}\, \bc_n\in  \on{Ind}(D^b(\CA))^{\leq 0}$ 
follows from the next assertion: 

\begin{lem} \label{l:when below 0}
Let $\CA$ be a small abelian category of finite cohomological dimension. 
Let $\bc\in \on{Ind}(D^b(\CA))$ be an object such that
$$\CHom_{\on{Ind}(D^b(\CA))}(a,\bc)\in \Vect_\sfe^{\leq 0} \text{ for all } a\in \CA\subset \on{Ind}(D^b(\CA)).$$
Then $\bc\in \on{Ind}(D^b(\CA))^{\leq 0}$.
\end{lem}

%
%
%

\begin{proof}

Suppose that $\tau^{>0}(\bc)\neq 0$. Let $k>0$ be the smallest integer such that $H^k(\bc)\neq 0$. 
Then 
$$\underset{a}{\on{colim}}\, H^0\left(\CHom_{\on{Ind}(D^b(\CA))}(a[-k],\tau^{\geq k}(\bc))\right)\neq 0,$$
where the colimit goes over the (filtered) category, whose objects are objects of $\CA$, and whose morphisms
are surjections. 

\medskip

Note that for any fixed $n$, the map
$$\underset{a}{\on{colim}}\, H^0\left(\CHom_{\on{Ind}(D^b(\CA))}(a[-k],\tau^{\geq k-n}(\bc))\right)\to
\underset{a}{\on{colim}}\, H^0\left(\CHom_{\on{Ind}(D^b(\CA))}(a[-k],\tau^{\geq k}(\bc))\right)$$
is surjective (by the definition of the derived category). 

\medskip

Let $d$ be the cohomological dimension of $\CA$. Note that for any $a$, the map
$$H^0\left(\CHom_{\on{Ind}(D^b(\CA))}(a[-k],\bc)\right)\to 
H^0\left(\CHom_{\on{Ind}(D^b(\CA))}(a[-k],\tau^{\geq k-d}(\bc))\right)$$
is surjective. 

\medskip

Hence, we obtain that the map
$$\underset{a}{\on{colim}}\, H^0\left(\CHom_{\on{Ind}(D^b(\CA))}(a[-k],\bc)\right)\to 
\underset{a}{\on{colim}}\, H^0\left(\CHom_{\on{Ind}(D^b(\CA))}(a[-k],\tau^{\geq k}(\bc))\right)$$
is surjective. 

\medskip

In particular, we obtain that for some $a$,
$$H^0\left(\CHom_{\on{Ind}(D^b(\CA))}(a[-k],\bc)\right)\neq 0.$$

However, this contradicts the assumption on $\bc$. 

\end{proof} 

\ssec{Categorical $K(\pi,1)$'s}

\sssec{} \label{sss:Kpi1}

Recall the subcategories
$$\Lisse(Y)\subset \iLisse(Y)\subset \qLisse(Y),$$
see Sects. \ref{ss:lisse}-\ref{ss:iLisse}. 

\begin{defn}
We shall say that $Y$ is a \emph{categorical $K(\pi,1)$} if the naturally defined functor
$$D^b(\Lisse(Y)^\heartsuit)\to \Lisse(Y)$$
is an equivalence.
\end{defn}

Note that from \thmref{t:der left complete} and \secref{sss:qLisse left compl iLisse} we obtain:

\begin{cor} \label{c:Kpi1}
If $Y$ is a categorical $K(\pi,1)$, then the inclusion
\begin{equation} \label{e:access to all}
\iLisse(Y)\subset \qLisse(Y)
\end{equation} 
is an equality.
\end{cor}

\sssec{}

Let $\Lisse(Y)_0\subset \Lisse(Y)$ be the full subcategory, consisting of objects whose
cohomologies (with respect to the usual t-structure) are extensions of the constant sheaf
$\ul\sfe_Y$. Let
$$\iLisse(Y)_0\subset \iLisse(Y) \text{ and } \qLisse(Y)_0\subset \qLisse(Y)$$
denote the corresponding subcategories. 

\begin{defn}
We shall say that $Y$ is a \emph{unipotent categorical $K(\pi,1)$} if the naturally defined functor
$$D^b(\Lisse(Y)_0^\heartsuit)\to \Lisse(Y)_0$$
is an equivalence.
\end{defn}

\sssec{} \label{sss:analyze P1}

An easy example of $Y$, which is neither a categorical $K(\pi,1)$ nor a unipotent categorical $K(\pi,1)$
is $Y=\BP^1$. 

\medskip

First, note that in this case, the embedding
$$\Lisse(Y)_0 \hookrightarrow \Lisse(Y)$$
is an equivalence. So, the statements about categorical $K(\pi,1)$ vs. unipotent categorical $K(\pi,1)$
are equivalent. 

\medskip

We will show that the functor \eqref{e:access to all} is \emph{not} an equivalence.

\medskip

Indeed, the category $\iLisse(Y)$ is generated by one object, namely, $\ul\sfe_{\BP^1}$, whose
algebra of endomorphisms is
$$A:=\sfe[\eta]/\eta^2=0, \,\, \deg(\eta)=2.$$

Hence,
$$\iLisse(Y)\simeq A\mod.$$

By Koszul duality, we have
$$A\mod \simeq B\mod_{0},$$
where 
$$B=\sfe\langle \xi\rangle, \,\, \deg(\xi)=-1$$
is the free \emph{associative} algebra on one generator in degree $-1$, and 
\begin{equation} \label{e:P1 bad}
B\mod_{0}\subset B\mod
\end{equation}
is the full subcategory consisting of objects on which $\xi$ acts locally nilpotently.

\medskip

The t-structure on $\iLisse(Y)$ corresponds to the usual t-structure on 
$B\mod$, for which the forgetful functor to $\Vect_\sfe$ is t-exact.

\medskip

Now it is easy to see that the embedding \eqref{e:P1 bad} realizes $B\mod$ as
the left completion of $B\mod_{0}$.

\sssec{}  \label{sss:curves}

We now claim:

\begin{thm} \label{t:curves} \hfill

\smallskip

\noindent{\em(a)} All connected algebraic curves other than $\BP^1$ are categorical $K(\pi,1)$'s.

\smallskip

\noindent{\em(b)} All connected algebraic curves other than $\BP^1$ 
are unipotent categorical $K(\pi,1)$'s.
\end{thm} 

\sssec{}

We observe:

\begin{lem} \label{l:from heart to D}
Let $\bC_0$ be a small DG category equipped with a \emph{bounded} t-structure,
and consider the functor
\begin{equation} \label{e:from heart to D}
D^b(\bC_0^\heartsuit)\to \bC_0.
\end{equation}

\smallskip

\noindent{\em(a)} 
Suppose every object of $\bc_0\in \bC_0^\heartsuit$ admits a non-zero map
to an \emph{injective} object $\bc\in \on{Ind}(\bC^\heartsuit_0)$ that satisfies 
$$\Hom_{\on{Ind}(\bC^0)}(\bc'_0,\bc[k])=0, \quad \forall\, \bc'_0\in \bC_0^\heartsuit,\,\, \forall\, k>0.$$
Then \eqref{e:from heart to D} is an equivalence. 

\smallskip

\noindent{\em(b)} Suppose that 
$$\Hom_{\bC_0}(\bc'_0,\bc_0[k])=0 \text{ for } k> 2 \text{ for all } \bc_0,\bc_0'\in \bC_0^\heartsuit.$$
Then \eqref{e:from heart to D} is an equivalence if and only if for every $\bc_0,\bc'_0$ as above, the 
(a priori injective) map
$$\Ext^2_{\bC_0^\heartsuit}(\bc'_0,\bc_0)\to \Hom_{\bC_0}(\bc'_0,\bc_0[2])$$
is surjective.
\end{lem}

\begin{proof}

Point (a) is standard: the assumption allows us to compute $\Hom_{\bC^0}(\bc'_0,\bc_0)$
via (ind)-injective resolutions in $\bC_0$. Point (b) follows formally from point (a).

\end{proof} 

\begin{rem}

Note that \lemref{l:from heart to D}(b) implies the assertion of \thmref{t:curves} when $X$ is \emph{affine}, as in this case
$$\Hom_{\qLisse(Y)}(E,E'[2])=0$$
for any pair of local systems $E$ and $E'$. 

\end{rem}

\ssec{Proof of \thmref{t:curves} for complete curves}

Let $X$ be a complete algebraic curve of genus $>0$. Let $E_0\in \Lisse(X)^\heartsuit$ denote 
the trivial local system. 

\sssec{}

We will first show that point (b) of \thmref{t:curves} implies point (a). 

\medskip

By \lemref{l:from heart to D}(b), we have to show that for $E_1,E\in \Lisse(X)^\heartsuit$,
any element 
$$\alpha\in \Hom_{\Lisse(X)}(E_1,E[2])$$
can be written as a cup product of classes 
$$\beta\in \Ext^1(E_1,\wt{E}) \text{ and } \gamma\in \Ext^1(\wt{E},E)$$
for some $\wt{E}\in \Lisse(X)^\heartsuit$. 

\medskip

Dualizing $E_1$, we can assume that $E_1=E_0$, so we can think of $\alpha$ as an element of
$H^2(X,E)$. 

\medskip

Write
$$0\to E''\to E\to E'\to 0,$$
where $E'$ is an extension of copies of $E_0$, and $E''$ does not have trivial quotients. Note that the map
$$H^2(X,E)\to H^2(X,E')$$
is an isomorphism, since $H^2(X,E'')=0$. Let $\alpha'$ be the image of $\alpha$ in $H^2(X,E')$

\medskip

Assuming point (b), we can write $\alpha'$ as a cup product of classes
$$\beta\in H^1(X,\wt{E}) \text{ and } \gamma'\in \Ext^1(\wt{E},E')$$
for some $\wt{E}\in \Lisse(X)_0^\heartsuit$ (i.e., $\wt{E}$ is also an extension of
copies of $E_0$). 

\medskip

It suffices to show
that $\gamma'$ can be lifted to an element $\gamma\in \Ext^1(\wt{E},E)$. However, the obstruction to such a lift lies
in $\Ext^2(\wt{E},E'')$, which embeds into $\Hom_{\Lisse(X)}(\wt{E},E''[2])$, and the latter vanishes since 
$$\Hom_{\Lisse(X)}(E_0,E''[2])=H^2(X,E'')=0.$$

\sssec{}

We now prove point (b) of \thmref{t:curves}. By \lemref{l:from heart to D}, it suffices to construct an
object $E_0^{\on{cofree}_x}\in \iLisse(X)$ with the following properties:

\begin{itemize}

\item(i) $E_0^{\on{cofree}_x}\in \iLisse(X)^\heartsuit$;

\item(ii) $E_0^{\on{cofree}_x}\in \iLisse(X)_0$;

\item(iii) $\Hom(E_0,E_0^{\on{cofree}_x})\neq 0$;

\item(iv) $\Hom_{\iLisse(X)}(E_0,E_0^{\on{cofree}_x}[k])=0$ for $k>0$.

\end{itemize}

Indeed, note that point (iv) for $k=1$ implies that $E_0^{\on{cofree}_x}$ is injective as an object of 
$\iLisse^\heartsuit_0$, and combining with point (iii), we obtain that any object in $\Lisse^\heartsuit_0$
admits a non-zero map to $E_0^{\on{cofree}_x}$.

\sssec{}

Choose a point $x\in X$, and consider the corresponding augmentation map
$$\on{C}^\cdot(X)\to \sfe.$$

\medskip

Set 
$$E_0^{\on{cofree}_x}:=\sfe\underset{\on{C}^\cdot(X)}\otimes E_0,$$
where we regard $\on{C}^\cdot(X)$ as $\CHom_{\iLisse(X)}(E_0,E_0)$. 

\medskip

Item (ii) follows by construction. Items (iii) and (iv) follow from the fact that
$$\CHom_{\iLisse(X)}(E_0,E_0^{\on{cofree}_x}) \overset{\text{$E_0$ is compact}}\simeq
\sfe\underset{\on{C}^\cdot(X)}\otimes \CHom_{\iLisse(X)}(E_0,E_0) \simeq \sfe.$$

\medskip

It remains to establish item (i).

\sssec{}

Taking the fiber of $E_0^{\on{cofree}_x}$ at $x$, property (i) is equivalent to the fact that the object
\begin{equation} \label{e:Koszul}
\sfe\underset{\on{C}^\cdot(X)}\otimes \sfe \in \Vect_\sfe
\end{equation}
is acyclic off degree $0$. 

\medskip

This fact is probably well-known. We will supply a proof for completeness. 

\sssec{} \label{sss:reduce to Betti}

First, the manipulation in Sects. \ref{ss:spectral etale char 0}-\ref{ss:spectral etale char p} allows us to
reduce the assertion to the case when our sheaf-theoretic context is Betti (and $k=\BC$).

\begin{rem}

Note that in Sects. \ref{ss:spectral etale char 0}-\ref{ss:spectral etale char p} we appealed to \thmref{t:curves}
(which we are still in the process of proving) in order to compare the categories $\qLisse(X)$ in the Betti context to the the 
\'etale context in characteristic $0$ and further to the \'etale context in characteristic $p$.

\medskip

However, there is no circularity in the argument, because for the purposes of \secref{sss:reduce to Betti}, we 
only need to compare the algebras of cochains $\on{C}^\cdot(X)$ in the three contexts, and the fact that
the corresponding maps are isomorphisms is standard. 

\end{rem}

\sssec{}

In the Betti context, it is known that the algebra $\on{C}^\cdot(X)$ is formal (by \cite[Main Theorem, Sect. 6]{DGMS}). I.e., it is isomorphic
to the DG algebra $A$ with
$$A_0=\sfe, \,\, A_1=V, \,\, A_2=\sfe,\,\, A_n=0 \text{ for } n>2,$$
where $V$ is a symplectic vector space and the multiplication $V\otimes V\to \sfe$ is 
given by a symplectic form. 

\medskip

If $V\neq 0$ (it is here that we use the assumption that the genus of $X$ is $>0$), this
algebra is quadratic and hence Koszul\footnote{We are grateful for L.~Positselski for explaining this to us.},
see \cite[Chapter 5, Sect. 5(1)]{PP}. Hence, 
$$H^k(\sfe\underset{A}\otimes \sfe)=0 \text{ for } k\neq 0.$$

\qed[\thmref{t:curves}]

\ssec{The dual of $\qLisse(Y)$} \label{ss:Verdier compat}

In this subsection we let $Y$ be a \emph{smooth} scheme. 

\sssec{}

We give the following definition:

\begin{defn} 
We shall say that $\qLisse(Y)$ is duality-adapted if the functor
\begin{equation} \label{e:Verdier on lisse}
\qLisse(Y)\otimes \qLisse(Y)\to \Vect_\sfe, \quad E_1,E_2 \mapsto \on{C}^\cdot(Y,E_1\sotimes E_2).
\end{equation}
is the counit of a self-duality.
\end{defn}

\sssec{}

We claim:

\begin{prop} \label{p:duality-adapted}
Assume that $\iLisse(Y)\to \qLisse(Y)$
is an equivalence. Then $\qLisse(Y)$ is duality-adapted.
\end{prop}

\begin{proof}

Since $\iLisse(Y)\to \qLisse(Y)$ is an equivalence, the category 
$\qLisse(Y)$ is compactly generated by $\Lisse(Y)$. 
Now, \emph{naive} duality defines a contravariant equivalence
$$\Lisse(Y) \simeq \Lisse(Y)^{\on{op}}.$$

Since $Y$ is smooth, the above naive duality coincides with Verdier duality,
up to a shift. Hence, the latter defines an identification
$$\qLisse(Y)\simeq \qLisse(Y)^\vee.$$

Its counit is given by \eqref{e:Verdier on lisse} by definition.

\end{proof} 

\begin{rem} \label{r:Verdier compat}
Note that the above argument shows that for any smooth $Y$, the pairing
$$\CF_1,\CF_2 \mapsto \on{C}^\cdot(Y,\CF_1\sotimes \CF_2)$$
defines a self-duality on $\iLisse(Y)$.
\end{rem} 

\sssec{}

We claim:

\begin{cor}  \label{c:curve Verdier compat}
If $X$ is a smooth curve, then $\qLisse(X)$ is duality-adapted.
\end{cor}

\begin{proof}

The case of curves different from $\BP^1$ follows from \thmref{t:curves}(a) and \propref{p:duality-adapted}.

\medskip

The case of $\BP^1$ follows by direct inspection: in terms of the equivalence
$$\qLisse(\BP^1)\simeq B\mod$$
(see \secref{sss:analyze P1}), the pairing \eqref{e:Verdier on lisse} corresponds (up to a shift)
to the canonical pairing
$$B\mod\otimes B^{\on{op}}\mod\to \Vect_\sfe,$$
corresponding to the isomorphism
$$B\simeq B^{\on{op}}, \quad \xi\mapsto -\xi.$$

\end{proof} 

\ssec{Specifying singular support}

In this subsection we let $Y$ be a \emph{smooth} scheme. In \secref{ss:sing supp sing} below we will
explain how to extend the discussion to the case when $Y$ is not necessarily smooth. 

\sssec{} \label{sss:sing supp schemes perv}

Let $Y$ be a scheme and $\CN$ a conical Zariski-closed subset of $T^*(Y)$. In this case we have a well-defined
full abelian subcategory
$$\Perv_\CN(Y)\subset \Perv(Y),$$
see \cite[Sect. A.3.1]{GKRV}.

\medskip

A key property of $\Perv_\CN(Y)$ is that it is a Serre subcategory of $\Perv(Y)$, i.e., it is stable
under extensions and the operations of taking sub- and quotient objects. 

\begin{rem}
In the \'etale context, the notion of singular support for objects of $\Shv(Y)^{\on{constr}}$ was 
introduced in \cite{Be2} for \'etale sheaves with torsion coefficients. However, the theory applies ``as-is" for $\ell$-adic sheaves. Indeed,
the definition of singular support is based on the notion of \emph{universal local acyclicity}, which is equally applicable to $\ell$-adic sheaves.
\end{rem}

\sssec{}

Another basic feature of this subcategory is that the Verdier self-duality
$$\BD:\Perv(Y)\to \Perv(Y)$$
sends $\Perv_\CN(Y)$ to itself. 

\medskip

This follows from the geometric characterization of singular support in \cite{Be2}. 

\sssec{} \label{sss:sing supp schemes}

Consider the abelian category
$$\on{Ind}(\Perv_\CN(Y))\subset \on{Ind}(\Perv(Y))\simeq \Shv(Y)^\heartsuit.$$

\medskip

We let
$$\Shv_\CN(Y)\subset \Shv(Y)$$
be the full subcategory consisting of objects whose cohomologies belong to $\on{Ind}(\Perv_\CN(Y))$.

\medskip

Since the t-structure on $\Shv(Y)$ is compatible with filtered colimits, we obtain that $\Shv_\CN(Y)$
is closed under filtered colimits. 

\medskip 

By construction, $\Shv_\CN(Y)$ inherits a t-structure so that its embedding into $\Shv(Y)$
is t-exact. By \thmref{t:Shv left-complete}, the category $\Shv_\CN(Y)$ is left-complete in its t-structure. 

\sssec{} \label{sss:sing supp schemes constr}

Set
$$\Shv_\CN(Y)^{\on{constr}}:=\Shv_\CN(Y)\cap \Shv(Y)^{\on{constr}}\subset \Shv(Y).$$

This is the full subcategory of $\Shv(Y)^{\on{constr}}$ consisting of objects whose cohomologies
belong to $\Perv_\CN(Y)$. By construction, $\Shv_\CN(Y)^{\on{constr}}$ inherits a t-structure so that 
its embedding into $\Shv(Y)^{\on{constr}}$ is t-exact.

\medskip

Set 
$$\Shv_\CN(Y)^{\on{access}}:=\on{Ind}(\Shv_\CN(Y)^{\on{constr}}).$$

Ind-extension of the tautological embedding 
$$\Shv_\CN(Y)^{\on{constr}}\hookrightarrow \Shv_\CN(Y)$$
defines a functor
\begin{equation} \label{e:access to all N}
\Shv_\CN(Y)^{\on{access}}\to \Shv_\CN(Y).
\end{equation}

The functor \eqref{e:access to all N} is fully faithful: indeed, its composition with the embedding
$\Shv_\CN(Y)\hookrightarrow \Shv(Y)$ preserves compactness and is fully faithful on compacts. 

\medskip

The t-structure on $\Shv_\CN(Y)^{\on{constr}}$ extends to a unique t-structure on $\Shv_\CN(Y)^{\on{access}}$
compatible with filtered colimits. The functor \eqref{e:access to all N} is t-exact, since the t-structure on 
$\Shv(Y)$ (and hence $\Shv_\CN(Y)$) is also compatible with filtered colimits.

\medskip

It is easy to see that the functor \eqref{e:access to all N} 
induces an equivalence on the hearts. Hence, it induces an equivalence 
$$(\Shv_\CN(Y)^{\on{access}})^{\geq -n}\to (\Shv_\CN(Y))^{\geq -n}$$
for any $n$. From here, and the fact that $\Shv_\CN(Y)$ is left-complete, 
it follows that the functor \eqref{e:access to all N} identifies $\Shv_\CN(Y)$ with
the left completion of $\Shv_\CN(Y)^{\on{access}}$.

\sssec{Example}

Let $Y$ be smooth and take $\CN=\{0\}$. Then $\Shv_\CN(Y)$ is what we have previously denoted by $\qLisse(Y)$ and 
$\Shv_\CN(Y)^{\on{access}}=\iLisse(Y)$.

\begin{rem}  \label{r:bad P1}

Note that the process of left completion in \eqref{e:access to all N} is in general non-trivial,
i.e., the category $\Shv_\CN(Y)^{\on{access}}$ is not necessarily left-complete, see
\secref{sss:analyze P1}. 

\end{rem} 

\begin{rem}

Our conventions are different from those of \cite{GKRV}. In {\it loc.cit.} we denoted by
$\Shv_\CN(Y)$ what we denote here by $\Shv_\CN(Y)^{\on{access}}$.

\end{rem}  

\ssec{Singular support on schemes that are not necessarily smooth}  \label{ss:sing supp sing}

Let $Y$ be a scheme of finite type that is not necessarily smooth. In this subsection we explain what
we mean by a closed subset of $T^*(Y)$, and how to assign singular support to objects of $\Shv(Y)$.
(The same discussion applies when instead of $\Shv(Y)$ we consider $\Dmod(Y)$ and $\Shv^{\on{all}}(Y)$.)

\medskip

The discussion is Zariski-local, so we can assume that $Y$ is affine. 

\sssec{}

Let $Y$ be a scheme of finite type, and let $E$ be a (classical) coherent sheaf on $Y$. Write $E$ as a 
quotient of a map $E_{-1}\to E_0$, where $E_{-1}$ and $E_0$ are vector bundles. Let
$\on{Tot}(E)$ be the algebraic stack over $Y$ equal to 
$$\on{Tot}(E_0)/\on{Tot}(E_{-1}),$$
where we regard $\on{Tot}(E_{-1})$ as a group-scheme over $Y$ that acts on $\on{Tot}(E_0)$. 

\medskip

Of course, $\on{Tot}(E)$ as defined above depends on the presentation. However, it is well-defined
as an object of the category obtained by localizing algebraic stacks with respect to morphisms, whose
fibers are of the form $\on{pt}/V$, where $V$ is a vector group.

\medskip

In particular, the notion of a (locally) closed subset in $\on{Tot}(E)$ is well-defined. By definition, this is
the same as a (locally) closed subset in $\on{Tot}(E_0)$ that is $\on{Tot}(E_{-1})$-invariant for every 
$\on{Tot}(E_{-1})$ as above. 

\medskip

One can talk about (isomorphism classes of) $k$-points of $\on{Tot}(E)$. They are in bijection with pairs
$(y,\xi)$, where $y\in Y(k)$, and $\xi$ is an element in the (classical) fiber of $E$ at $y$. 

\sssec{}

We apply the above discussion to $E=\Omega(Y)$, the (classical) sheaf of K\"ahler differentials. Denote the
resulting stack $\on{Tot}(\Omega(Y))$ by $T^*(Y)$. 

\medskip

By definition, its $k$-points are pairs $(y,\xi)$, where $y\in Y(k)$ and $\xi\in T^*_y(Y)$. 

\medskip

Note that for every embedding $Y\overset{j}\hookrightarrow Y'$, where $Y'$ is smooth, we obtain a presentation 
of $T^*(Y)$ as a quotient of the vector bundle $T^*(Y')|_Y$. Denote by $dj^*$ the resulting projection
\begin{equation} \label{e:codiff sing}
T^*(Y')|_Y \to T^*(Y).
\end{equation}

\sssec{} \label{sss:half-dim sing}

We will say that a (locally) closed subset $\CN\subset T^*(Y)$ is half-dimensional/Largrangian
if, locally, its preimage along \eqref{e:codiff sing}, viewed as a subset of $T^*(Y')$, has the corresponding property. 

\medskip

It is easy to see that these notions are independent of the choice of an embedding $Y\overset{j}\hookrightarrow Y'$.

\sssec{}

We now claim that to every object $\CF\in \Shv(Y)^{\on{constr}}$ we can attach a closed subset
$$\on{SingSupp}(\CF)\subset T^*(Y)$$
with the following property: 

\medskip

For every closed embedding $Y\overset{j}\hookrightarrow Y'$ with $Y'$ is 
a smooth affine scheme, we have
$$(dj^*)^{-1}(\on{SingSupp}(\CF))=\on{SingSupp}(j_*(\CF)).$$

\begin{proof}

We only need to show that for every $y\in Y$, the intersection
$$T^*_y(Y')\cap \on{SingSupp}(j_*(\CF))$$
is invariant under translations by the elements of
$$\on{ker}(T^*_y(Y')\to T^*_y(Y)).$$

\medskip

Note that if 
$$Y\overset{j_1}\hookrightarrow Y'_1 \text{ and } Y\overset{j_2}\hookrightarrow Y'_2$$
are two embeddings as above, we can always complete them to a commutative diagram
$$
\CD
Y @>{j_1}>> Y_1 \\
@V{j_2}VV @VV{j'_2}V \\
Y_2 @>{j'_1}>> Y_{1,2},
\endCD
$$
where $j'_1$ and $j'_2$ are also closed embeddings and $Y_{1,2}$ is smooth. 

\medskip

From here it is easy to see that if for a given $(Y'_1,j_1)$ and a pair of cotangent
vectors $\xi'_1,\xi''_1\in T_y^*(Y_1)$ that project to the same vector in $T^*_y(Y)$
one can find a commutative diagram as above and a pair of cotangent
vectors $\xi'_{1,2},\xi''_{1,2}\in T_y^*(Y_{1,2})$ that project to $\xi'_1,\xi''_1$,
respectively under $T_y^*(Y_{1,2})\to T_y^*(Y_1)$ and that project to the same element,
to be denoted $\xi_2\in T_y^*(Y_2)$ under $T_y^*(Y_{1,2})\to T_y^*(Y_2)$. 

\medskip

Since our assertion holds for smooth $Y$, we have
$$\xi'_1\in \on{SingSupp}((j_1)_*(\CF))\, \Leftrightarrow\, 
\xi'_{1,2}\in \on{SingSupp}((j'_2)_*\circ (j_1)_*(\CF))\, \Leftrightarrow\, \xi'_{1,2}\in \on{SingSupp}((j'_1)_*\circ (j_2)_*(\CF))\, 
\Leftrightarrow\,$$
$$\Leftrightarrow\, \xi_{2}\in \on{SingSupp}((j_2)_*(\CF))\, \Leftrightarrow\,$$
$$\Leftrightarrow\,\xi''_{1,2}\in \on{SingSupp}((j'_1)_*\circ (j_2)_*(\CF))\,
\Leftrightarrow\, \xi''_{1,2}\in \on{SingSupp}((j'_2)_*\circ (j_1)_*(\CF))\, \Leftrightarrow\, 
\xi''_1\in \on{SingSupp}((j_1)_*(\CF)),$$
as required.

\end{proof}

\sssec{}

The above construction of $\on{SingSupp}(\CF)$ implies that it has the usual functoriality 
property of singular support. 

\medskip

In particular, if $f:Y_1\to Y_2$ is a closed embedding and $\CF_1\in \Shv(Y_1)$, we have 
$$\on{SingSupp}(f_*(\CF_1))=(df^*)^{-1}(\on{SingSupp}(\CF_1),$$
where $df^*$ denotes the map
$$T^*(Y_2)\supset Y_1\underset{Y_2}\times T^*(Y_2)\to T^*(Y_1).$$

\medskip

If $f:Y_1\to Y_2$ is a smooth morphism and  $\CF_2\in \Shv(Y_2)$, we have 
$$\on{SingSupp}(f^*(\CF_2)):=df^*\left(Y_1\underset{Y_2}\times \on{SingSupp}(\CF_2)\right).$$

\sssec{The case of ind-schemes}

In \secref{ss:Hecke arb G} we need to also consider the notion of singular support on ind-schemes
(of ind-finite type).

\medskip

If $\CY$ is an ind-scheme, and let $Z\subset \CY$ be a closed-subscheme. We will consider 
$T^*(\CY)|_Z$ as a pro-object in the above localization of the category of stacks. Namely,
$$T^*(\CY)|_Z=``\on{lim}"\, T^*(Y_i)|_Z,$$
where $Y_i$ is a (cofinal) family of closed subschemes of $\CY$ such that $Z\subset Y_i$ and
$$\CY= \on{colim}\, Y_i.$$

In particular, a $k$-point of $T^*(\CY)$ is a pair $(y,\xi)$, where $y\in Y(k)$ and $\xi$
is an element in the classical pro-cotangent space to $\CY$ at $y$, viewed as a pro-finite dimensional
vector space. 

\medskip

To an object $\CF\in \Shv(\CY)$ supported on $Z$, and \emph{constructible} when viewed as an object of $\Shv(Z)$,
we can attach its singular support, which 
is a subset in $T^*(\CY)|_Z$, which projects to a closed subset in each $T^*(Y_i)|_Z$ above.

\medskip

This justifies the manipulations with singular support of sheaves on ind-schemes/stacks in 
\secref{ss:Hecke arb G}. 

\ssec{The external tensor product functor}

For a pair of schemes $Y_1$ and $Y_2$, consider the external tensor product functor
\begin{equation} \label{e:box product Y}
\Shv(Y_1)\otimes \Shv(Y_2)\to \Shv(Y_1\times Y_2).
\end{equation} 

The functor \eqref{e:box product Y} sends compacts to compacts, and is fully faithful, but \emph{not an equivalence} (unless one of the schemes
is a disjoint union of set-theoretic points). 

\sssec{}

Recall that for a pair of DG categories, each equipped with a t-structure, their tensor product
acquires a t-structure, see \secref{sss:ten prod t}.

\medskip

We claim:

\begin{prop} \label{p:boxtimes t-exact}
The functor \eqref{e:box product Y} is t-exact. 
\end{prop} 

The proof is given in \secref{ss:proof boxtimes t-exact} below. 

\sssec{}

Given $\CN_i\subset T^*(Y_i)$ and $\CF_i\in \on{Perv}_{\CN_i}(Y_i)$, 
$$\CF_1\boxtimes \CF_2\in \Shv_{\CN_1\times \CN_2}(Y_1\times Y_2).$$

From here, it follows that the same is true for arbitrary 
$\CF_i\in \Shv_{\CN_i}(Y_i)$. 

\medskip

Hence, we obtain a functor
\begin{equation} \label{e:box product Y N}
\Shv_{\CN_1}(Y_1)\otimes \Shv_{\CN_2}(Y_2)\overset{\boxtimes}\to \Shv_{\CN_1\times \CN_2}(Y_1\times Y_2).
\end{equation}

If one of the categories $\Shv_{\CN_i}(Y_i)$, $i=1,2$, is dualizable, the functor \eqref{e:box product Y N} is fully faithful: 

\medskip

Indeed, if, say, $\Shv_{\CN_1}(Y_1)$ is dualizable, write the composition
$$\Shv_{\CN_1}(Y_1)\otimes \Shv_{\CN_2}(Y_2)\to \Shv_{\CN_1\times \CN_2}(Y_1\times Y_2)\hookrightarrow \Shv(Y_1\times Y_2)$$
as 
\begin{equation} \label{e:add N}
\Shv_{\CN_1}(Y_1)\otimes \Shv_{\CN_2}(Y_2)\to \Shv_{\CN_1}(Y_1)\otimes \Shv(Y_2)\to \Shv(Y_1)\otimes \Shv(Y_2)\to 
\Shv(Y_1\times Y_2).
\end{equation} 

\sssec{} \label{sss:comp gen t}

We will say that a t-structure on a DG category $\bC$ is \emph{compactly generated}, if it satisfies the following 
two conditions: 

\begin{itemize}

\item $\bC^{\leq 0}$ is generated under filtered colimits by objects in $\bC^{\leq 0}\cap \bC^c$;

\item All of $\bC$ is generated under filtered colimits by shifts of objects in $\bC^{\leq 0}\cap \bC^c$.

\end{itemize}

This is equivalent to:

\begin{itemize}

\item If $\CHom_\bC(c_0,c)\in \Vect_\sfe^{> 0}$ for all $c_0\in \bC^{\leq 0}\cap \bC^c$, then $c\in \bC^{>0}$;

\item If $\CHom_\bC(c_0,c)=0$ for all $c_0\in \bC^{\leq 0}\cap \bC^c$, then $c=0$.

\end{itemize} 

\medskip

For example, for any $Y$, the t-structure on $\Shv(Y)$ is compactly generated. 

\medskip

It is easy to see that if a t-structure on $\bC$ is compactly generated, then it is compatible
with filtered colimits and is right-complete. 

\medskip

We have:

\begin{lem} \label{l:comp gen t}
Let $\bC_i$ be DG categories, each equipped with a t-structure, and let $F:\bC_1\to \bC_2$ be a 
t-exact functor. Let $\bC$ be another DG category, equipped with a compactly generated t-structure.
Then the functor
$$F\otimes \on{Id}_\bC:\bC_1\otimes \bC\to \bC_2\otimes \bC$$
is t-exact.
\end{lem}

\sssec{}

Combining \propref{p:boxtimes t-exact} and \lemref{l:comp gen t} (applied to the maps in \eqref{e:add N}), we obtain:

\begin{cor} \label{c:comp gen t}
Suppose that for one of the categories $\Shv_{\CN_i}(Y_i)$, $i=1,2$ its t-structure is compactly generated.
Then the functor \eqref{e:box product Y N} is t-exact.
\end{cor} 

\sssec{Example}

Note that the conditions of the corollary are satisfied for $(Y,\CN)=(X,\{0\})$, where $X$ is a smooth curve.

\medskip

Indeed, when $X$ is a connected curve different from $\BP^1$, this follows from \thmref{t:curves}(a). 
When $X\simeq \BP^1$, this follows from the explicit description of the category $\qLisse(\BP^1)$ 
in \secref{sss:analyze P1}. 

\begin{rem}
In fact, one can show that the functor \eqref{e:box product Y N} is t-exact for any $(Y_i,\CN_i)$. 
In fact the assertion of \lemref{l:comp gen t} holds without the assumption that the t-structure on $\bC$
be compactly generated. This is a rather non-trivial assertion, proved in \cite[Proposition C.4.4.1]{Lu3}.
\end{rem}

\ssec{Proof of \propref{p:boxtimes t-exact}} \label{ss:proof boxtimes t-exact}

\sssec{}

The functor \eqref{e:box product Y} is right t-exact by construction. Hence, it remains to show that
it is left t-exact.

\sssec{Sheaves at the generic point}

Let $Z$ be an irreducible scheme of finite type, and let $\eta_Z$ be its generic point. 

\medskip 

Set $$\Shv(\eta_Z):=\underset{U}{\on{colim}}\, \Shv(U),$$
where $U$ runs the (filtered) category of non-empty open subschemes of $Z$, and the transition functors
$\Shv(U_1)\to \Shv(U_2)$ for $U_2\subset U_1$ are given by restriction. The category $\Shv(\eta_Z)$ carries
a naturally defined t-structure.

\medskip

Define $\iLisse(\eta_Z)$ by a similar procedure. 
We have a tautological functor 
\begin{equation} \label{e:Lisse to all gen}
\iLisse(\eta_Z)\to \Shv(\eta_Z),
\end{equation} 
and we claim that it is actually an equivalence. 

\medskip

Indeed, the functor \eqref{e:Lisse to all gen} is fully faithful because for every $U$, the category
$\iLisse(U)$ is compactly generated and the functor $\iLisse(U)\to \Shv(U)$ preserves compactness,
and the category of indices involved in the colimit is filtered. 

\medskip

The fact that \eqref{e:Lisse to all gen} is essentially surjective follows from the definition of
constructibility. 

\sssec{}

Let $Z'$ be another scheme. In a similar way, we define the category
$\Shv(\eta_Z\times Z')$, and if $Z'$ is also irreducible, the category $\Shv(\eta_Z\times \eta_{Z'})$. 

\sssec{}

Let $Y$ be a scheme of finite type. For an irreducible subvariety $Z\overset{\bi_Z}\hookrightarrow Y$ with generic point
$\eta_Z$, let $\bi_{\eta_Z}^!$ denote the functor 
$$\Shv(Y)\overset{\bi_Z^!}\to \Shv(Z)\to  \Shv(\eta_Z).$$

\medskip

It is easy to see that an object $\CF\in \Shv(Y)$ is coconnective if and only if $\bi_{\eta_Z}^!(\CF)\in \Shv(\eta_Z)$
is coconnective for every $Z$.  

\sssec{}

We now return to the proof of the fact that the functor \propref{p:boxtimes t-exact} is left t-exact.

\medskip

By the above, it suffices to show that for every irreducible $Z\overset{\bi_Z}\hookrightarrow Y_1\times Y_2$,
the composite functor 
\begin{equation} \label{e:boxtimes eta Z}
\Shv(Y_1)\otimes \Shv(Y_2) \overset{\boxtimes}\to \Shv(Y_1\times Y_2) \overset{\bi_{\eta_Z}^!}\longrightarrow \Shv(\eta_Z)
\end{equation}
is left t-exact.

\medskip

Let 
$$Z_1\overset{\bi_{Z_1}}\hookrightarrow Y_1 \text{ and } Z_2\overset{\bi_{Z_2}}\hookrightarrow Y_2$$
be the closures of the images of $Z$ in $Y_1$ and $Y_2$, respectively. 

\medskip

The functor \eqref{e:boxtimes eta Z} factors as 
$$\Shv(Y_1)\otimes \Shv(Y_2) \overset{\bi_{\eta_{Z_1}}^!\boxtimes \bi_{\eta_{Z_2}}^!}\longrightarrow 
\Shv(\eta_{Z_1})\otimes \Shv(\eta_{Z_2})\to  \Shv(\eta_Z),$$
where the first arrow is t-exact by \lemref{l:comp gen t}. Hence, it suffices to show that the functor 
\begin{equation} \label{e:boxtimes gen}
\Shv(\eta_{Z_1})\otimes \Shv(\eta_{Z_2})\to  \Shv(\eta_Z) 
\end{equation}
is left t-exact. 

\sssec{}

Using the equivalence \eqref{e:Lisse to all gen}, we rewrite the functor \eqref{e:boxtimes gen} as
$$\iLisse(\eta_{Z_1})\otimes \iLisse(\eta_{Z_2})\to  \iLisse(\eta_Z).$$

Hence, it is enough to show that for any open \emph{smooth} $U_1\subset Z_1$ and $U_2\subset Z_1$ and 
$$U\subset Z \cap (U_1\times U_2),$$
the functor 
\begin{equation} \label{e:boxtimes open}
\iLisse(U_1)\otimes \iLisse(U_2)\overset{\boxtimes}\to \iLisse(U_1\times U_2)\overset{\bi^!_Z}\to \iLisse(U)
\end{equation}
is left t-exact, where $\iLisse(-)$ are considered in the \emph{perverse} t-structure. (In the above formula, 
by a slight abuse of notation we denote by $\bi_Z$ the locally closed embedding $U\to U_1\times U_2$.)

\sssec{}

Let $\on{pt} \overset{\bi_z}\to U$ be the embedding corresponding to a closed point $z\in U$. By the definition
of the perverse t-structure on $\iLisse(U)$ and \secref{sss:fiber functor x}, it suffices to show that the
composition
\begin{equation} \label{e:boxtimes open pt}
\iLisse(U_1)\otimes \iLisse(U_2)\overset{\boxtimes}\to \iLisse(U_1\times U_2)\overset{\bi^!_Z}\to \iLisse(U)
\overset{\bi^!_z[\dim(Z)]}\to \Vect_\sfe
\end{equation}
is left t-exact. 

\medskip

Let $z_1$ and $z_2$ be the images of $z$ in $U_1$ and $U_2$, respectively. Let $\bi_{z_i}$, $i=1,2$
denote the corresponding embeddings $\on{pt}\to U_i$. The functor \eqref{e:boxtimes open pt} identifies
with
$$\iLisse(U_1)\otimes \iLisse(U_2)\overset{\bi_{z_1}^![\dim(Z_1)]\otimes \bi_{z_2}^![\dim(Z_2)]}\longrightarrow
\Vect_\sfe \overset{[\dim(Z)-\dim(Z_1)-\dim(Z_2)]} \longrightarrow \Vect_\sfe.$$

In the above composition, the first arrow is t-exact by \lemref{l:comp gen t}, and the second arrow is left t-exact
because $\dim(Z_1)+\dim(Z_2)\geq \dim(Z)$.

\qed[\propref{p:boxtimes t-exact}]

\ssec{The tensor product theorems}

In this subsection we will discuss several variants of the tensor product result \cite[Theorem A.3.8]{GKRV}.

\sssec{}

First, we have the following result, which is \cite[Theorem A.3.8]{GKRV}.

\begin{thm} \label{t:product thm 1 sch}
Assume that $X$ is smooth and proper. Let $\CN\subset T^*(Y)$ be half-dimensional. 
Then the resulting functor
\begin{equation} \label{e:product functor 1}
\iLisse(X)\otimes \Shv_{\CN}(Y)^{\on{access}}\to \Shv_{\{0\}\times \CN}(X\times Y)^{\on{access}}
\end{equation} 
is an equivalence.
\end{thm}

\begin{rem}
It is natural to ask whether the functor 
$$\qLisse(X)\otimes \Shv_{\CN}(Y)\to \Shv_{\{0\}\times \CN}(X\times Y)$$
is an equivalence. 

\medskip

Unfortunately, we do not have an answer to this, except in the cases covered by Theorems
\ref{t:product thm 1.5 sch} and \ref{t:product thm 2 sch} below. 
Namely, we did not find a way to determine when the tensor product
$$\qLisse(X)\otimes \Shv_{\CN}(Y)$$
is left-complete in its t-structure. 

\end{rem}

%
%

\sssec{}

Next, we claim: 

\begin{thm} \label{t:product thm 1.5 sch}
Let $X$ be smooth and proper. 
Then the functor 
\begin{equation} \label{e:product functor 1.5}
\qLisse(X)\otimes \Shv(Y)\to \Shv(X\times Y)
\end{equation} 
is an equivalence onto the full subcategory that consists of objects $\CF$ with the following property:

\smallskip

\noindent For every $m$
and every constructible sub-object $\CF'$ of $H^m(\CF)$, the singular support of $\CF'$
is contained in a subset of the form $\{0\}\times \CN$, where $\CN\subset T^*(Y)$ is half-dimensional.  
\end{thm}

The proof will use the following variant of \thmref{t:der left complete} (the proof is given in \secref{sss:left compl ten prod}):

\begin{thm} \label{t:der left complete prod}
Let $\CA$ be a small abelian category of finite
cohomological dimension. Let $\bC$ be a DG category equipped with a t-structure in 
which it is left-compete. Then 
$$\on{Ind}(D^b(\CA))\otimes \bC$$ 
is left-complete in its t-structure.
\end{thm}

\begin{proof}[Proof of \thmref{t:product thm 1.5 sch}]

First, we observe that the functor \eqref{e:product functor 1.5} is fully faithful, being a composition of
$$ \qLisse(X)\otimes \Shv(Y)\to \Shv(X)\otimes \Shv(Y)\to \Shv(X\times Y),$$
where the first arrow is fully faithful because $\Shv(Y)$ is dualizable.

\medskip

Thus, it remains to show that \eqref{e:product functor 1.5} is essentially surjective onto the
specified subcategory. 

\medskip

First, we claim that every \emph{bounded below} object in $\Shv(X\times Y)$ with the specified condition
belongs to the essential image of \eqref{e:product functor 1.5}. 

\medskip

Indeed, by devissage we can assume that the object
in question is also bounded above; then that it is in the heart of the t-structure, and then that it is contained in 
$\Perv(X\times Y)$, and has singular support of the form $\{0\}\times \CN$ with $\CN\subset T^*(Y)$ is half-dimensional.  
However, such an object is contained in the essential image of \eqref{e:product functor 1}, by \thmref{t:product thm 1 sch}.

\medskip

Now, the assertion of the theorem follows, as both sides are left-complete in their respective t-structures:
the right-hand side by \thmref{t:Shv left-complete}, and the left-hand side by \thmref{t:der left complete prod}.

\end{proof} 

\begin{cor} \label{c:1.75}
Suppose that $\on{char}(k)=0$. Then the functor \eqref{e:product functor 1.5} is an equivalence onto a subcategory 
consisting of objects whose singular support is contained in $\{0\}\times T^*(Y)\subset T^*(X\times Y)$.
\end{cor}

\begin{proof} 

It suffices to show for every constructible object $\CF\in \Shv(X\times Y)$, with 
$$\on{SingSupp}(\CF)\subset \{0\}\times T^*(Y)\subset T^*(X\times Y),$$
the singular support of $\CF$ is in fact contained in a subset of the form $\{0\} \times \CN$ for some half-dimensional $\CN\subset T^*(Y)$.

\medskip

However, since $\on{char}(k)=0$, the singular support of $\CF$ is a \emph{Lagrangian} subset of $T^*(X\times Y)$,

\medskip

We claim that any irreducible Lagrangian subset $\CL\subset T^*(X\times Y)$ contained in $\{0\}\times T^*(Y)$ is
of required form.

\medskip

Indeed, at its generic point, $\CL$ is the conormal of some $Z\subset X\times Y$. However, if 
$$N^*_{Z/X\times Y}\subset \{0\}\times T^*(Y),$$
then $Z$ is of the form $X\times Y'$ for $Y'\subset Y$. 

\end{proof}

\sssec{}

Finally, we claim:

\begin{thm} \label{t:product thm 2 sch}
Let $X$ be smooth and proper. Assume also that $\qLisse(X)$ is duality-adapted (see \secref{ss:Verdier compat}). 
Let $\CN\subset T^*(Y)$ be half-dimensional. 
Then the resulting functor
\begin{equation} \label{e:product functor 2}
\qLisse(X)\otimes \Shv_{\CN}(Y)\to \Shv_{\{0\}\times \CN}(X\times Y)
\end{equation} 
is an equivalence.
\end{thm}

\begin{proof} 

Since $\qLisse(X)$ is dualizable, the functor
\begin{equation} \label{e:lisse X Y}
\qLisse(X)\otimes \Shv_{\CN}(Y) \to \qLisse(X)\otimes \Shv(Y)
\end{equation} 
is fully faithful.

\medskip

Given \thmref{t:product thm 1.5 sch}, it suffices to show that any object
$$\CF\in \qLisse(X)\otimes \Shv(Y)$$
whose image $\CF'\in \Shv(X\times Y)$ has singular support in $\{0\}\times \CN$,
belongs to the essential image of \eqref{e:lisse X Y}. 

\medskip

Since $\qLisse(X)$ is duality-adapted, it suffices to show that for any $E\in \qLisse(X)$, the object
$$(\on{C}^\cdot(X,-)\otimes \on{Id})(E\sotimes \CF)\in \Shv(Y)$$
belongs to $\Shv_{\CN}(Y)$. 

\medskip

However, the latter object is the same as 
$$(p_Y)_*(p_X^!(E)\sotimes \CF'),$$
where $p_X$ and $p_Y$ are the two projections from $X\times Y$ to $X$ and $Y$, respectively.

\medskip

The latter object indeed belongs to $\Shv_{\CN}(Y)$, due to the assumption on the singular support of $\CF'$
and the fact that $X$ is proper (see \lemref{l:gen sing supp estim}).

\end{proof} 

\sssec{Proof of \thmref{t:der left complete prod}} \label{sss:left compl ten prod}

The proof repeats the argument of \thmref{t:der left complete}, using the following variants of
\propref{p:Postnikov} and \lemref{l:when below 0}, respectively:

\medskip

\begin{prop} \label{p:Postnikov prod}
Let $\bC$ be compactly generated by compact objects of finite cohomological dimension.
Then for any $\bC_1$ equipped with a t-structure in which it is left-complete, the functor
\begin{equation} \label{e:map to left compl, prod}
\bC\otimes \bC_1\to (\bC\otimes \bC_1)^\wedge
\end{equation} 
is fully faithful and its right adjoint is continuous.
\end{prop}

\begin{lem} \label{l:when below 0 prod}
Let $\CA$ be as in \thmref{t:der left complete prod}, and let $\bC_1$ be equipped with a t-structure.
Let $\bc\in \on{Ind}(D^b(\CA))\otimes \bC_1$ be an object satisfying
$$(\CHom_{\on{Ind}(D^b(\CA))}(a,-) \otimes \on{Id})(\bc)\in (\bC_1)^{\leq 0} \text{ for all } a\in \CA\subset \on{Ind}(D^b(\CA)).$$
Then $\bc\in (\on{Ind}(D^b(\CA))\otimes \bC_1)^{\leq 0}$.
\end{lem}

Both these statements are proved in a way mimicking the original arguments. 

\section{Constructible sheaves on an algebraic stack} \label{s:shvs on stacks} 

As in \secref{s:shvs on sch}, in this section we let $\Shv(-)^{\on{constr}}$ be a constructible sheaf theory. 
All algebro-geometric objects will be assumed (locally) of finite type over the ground field $k$. 

\ssec{Generalities}

\sssec{} 

Let $\CY$ be a prestack. Recall that we define
$$\Shv(\CY):=\underset{S}{\on{lim}}\, \Shv(S),$$
where the index category is that of affine schemes equipped with a map to $\CY$,
and the transition functors are given by !-pullback.

\medskip

Since we are in the constructible context, the !-pullback functor admits a left adjoint, given by !-pushforward.
Hence, using \cite[Chapter 1, Proposition 2.5.7]{GR1}, we can rewrite
\begin{equation} \label{e:on prestack as colim}
\Shv(\CY)\simeq \underset{S}{\on{colim}}\, \Shv(S),
\end{equation} 
where the transition functors are given by !-pushforward.

\medskip

In particular, we obtain that $\Shv(\CY)$ is compactly generated: the compact generators are of the form
\begin{equation} \label{e:ins f}
\on{ins}_{f_0}(\CF_{S_0}), \quad \CF_{S_0}\in \Shv(S_0)^c,
\end{equation} 
where for an affine scheme $S_0$ equipped with a map $f_0:S_0\to \CY$, we denote by
$\on{ins}_{f_0}$ the corresponding tautological functor 
$$\Shv(S_0)\to \underset{S}{\on{colim}}\, \Shv(S)\simeq \Shv(\CY).$$

\sssec{} \label{sss:comp gen stack}

Suppose for a moment that $\CY$ is an algebraic stack\footnote{As per our conventions, we will assume that $\CY$ 
has an affine diagonal.}. Then the above index category
can be replaced by its non-full subcategory, where we allow as objects affine schemes that
are smooth over $\CY$, and as morphisms smooth maps between those. 

\medskip

Furthermore, formula \eqref{e:ins f} has a more explicit meaning: for $S\overset{f}\to \CY$, where $S$
is an affine scheme, we have
$$\on{ins}_f\simeq f_!,$$
where the functor
$$f_!: \Shv(S)\to \Shv(\CY)$$
is defined because the morphism $f$ is schematic.

\medskip

Thus, $\Shv(\CY)$ is compactly generated by objects of the form 
$$f_!(\CF_S), \quad \CF_S\in \Shv(S)^c.$$

Moreover, as above, it is sufficient to consider only the pairs $(S,f)$ with $f$ smooth.

\sssec{}

Recall that for a quasi-compact scheme $Y$, Verdier duality defines a contravariant equivalence
$$(\Shv(Y)^{\on{constr}})^{\on{op}} \overset{\BD}\to \Shv(Y)^{\on{constr}}.$$

Since
$$\Shv(Y):=\on{Ind}(\Shv(Y)^{\on{constr}}),$$
we obtain that the category $\Shv(Y)$ is canonically self-dual with the counit 
$$\Shv(Y)\otimes \Shv(Y)\to \Vect_\sfe$$
given by
$$\CF_1,\CF_2\mapsto \on{C}^\cdot(Y,\CF_1\sotimes \CF_2).$$

\sssec{}

In particular, by \cite[Proposition 1.8.3]{DrGa2} and \eqref{e:on prestack as colim},
we obtain that for a prestack $\CY$, the category $\Shv(\CY)$ is dualizable, and
$$\Shv(\CY)^\vee\simeq \underset{S}{\on{colim}}\, \Shv(S),$$
where the transition functors are given by *-pushforward.

\begin{rem}
Note that there is no a priori reason for $\Shv(\CY)^\vee$ to be equivalent to the original
$\Shv(\CY)$. 

\medskip

We will see that there is a canonical such equivalence when $\CY$ is a quasi-compact algebraic stack
(at least when $\CY$ is locally a quotient). However, for more general $\CY$ (e.g., for non-quasi-compact 
algebraic stacks) such an equivalence would reflect a particular feature of $\CY$, for example its property 
of being \emph{miraculous}, see \cite[Sect. 6.7]{Ga3}. 
\end{rem}

\ssec{Constructible vs compact}

\sssec{}

Let $\CY$ be an algebraic stack. Let
$$\Shv(\CY)^{\on{constr}}\subset \Shv(\CY)$$
be the full subcategory consisting of objects that pullback to an object of
$$\Shv(S)^{\on{constr}}=\Shv(S)^c\subset \Shv(S)$$
for any affine scheme $S$ mapping to $\CY$. 

\medskip

It is easy to see that this condition is enough to test on 
smooth maps $S\to \CY$. In the latter case, we can use either !- or *- pullback,
as they differ by a cohomological shift. 
 
\sssec{} 

Using the definition of the constructible subcategory via *-pullbacks along smooth maps, 
we obtain that we have an inclusion
\begin{equation} \label{e:comp in constr}
\Shv(\CY)^c\subset \Shv(\CY)^{\on{constr}}.
\end{equation} 

Indeed, for $f:S\to \CY$, the functor $f^*$ sends compacts to compacts, since its right adjoint,
namely $f_*$, is continuous. 

\medskip

However, the inclusion \eqref{e:comp in constr}
is typically \emph{not} an equality. For example, the constant sheaf
$$\ul\sfe_\CY\in \Shv(\CY)^{\on{constr}}$$
is \emph{not} compact for $\CY=B(\BG_m)$. 

\sssec{}

That said, we have the following assertion:

\begin{prop} \label{p:almost compact}
Let $\CY$ be quasi-compact. Then an object $\CF\in \Shv(\CY)^{\on{constr}}\cap \Shv(\CY)^{\geq n}$ is compact as an object 
of $\Shv(\CY)^{\geq m}$ for any $m\leq n$.
\end{prop}

\begin{proof} 

It suffices to show that for any $k$, the functor
$$\CF'\mapsto \tau^{\leq k}\left(\CHom_{\Shv(\CY)}(\CF,\CF')\right)$$
commutes with filtered colimits as $\CF'$ ranges over $\Shv(\CY)^{\geq m}$ for some fixed $m$.

\medskip

Choose a smooth covering $f:S\to \CY$, where $S$ is an affine scheme, and let $S^\bullet$ be its \v{C}ech nerve;
let $f^n:S^n\to \CY$ denote the resulting maps.

\medskip

For $\CF'\in \Shv(\CY)$, we can calculate $\CHom_{\Shv(\CY)}(\CF,\CF')$ as the totalization of the 
cosimplicial complex whose $n$-simplices are
$$\CHom_{\Shv(S^n)}((f^n)^!(\CF),(f^n)^!(\CF')).$$

Assume that $\CF\in \Shv(\CY)^{\leq N}$. Then for all $n$, the objects $\CHom_{\Shv(S^n)}((f^n)^!(\CF),(f^n)^!(\CF'))$
belong to $\Vect_\sfe^{\geq -N+m}$. Hence,
$$\tau^{\leq k}\left(\underset{[n]}{\on{lim}}\, \CHom_{\Shv(S^n)}((f^n)^!(\CF),(f^n)^!(\CF'))\right)$$
maps isomorphically to the limit 
$$\tau^{\leq k}\left(\underset{[n],n\leq k+N-m+1}{\on{lim}}\, \CHom_{\Shv(S^n)}((f^n)^!(\CF),(f^n)^!(\CF'))\right)$$
(e.g., by \cite[Proposition 1.2.4.5(4)]{Lu2}), which is a \emph{finite} limit.  

\medskip

Since finite limits commute with filtered colimits, it suffices to show that for every $n$, the functor
$$\CF'\mapsto \CHom_{\Shv(S^n)}((f^n)^!(\CF),(f^n)^!(\CF'))$$
commutes with filtered colimits. However, this follows from the fact that $(f^n)^!(\CF)$ is compact.

\end{proof}

\sssec{}

Verdier duality defines a contravariant equivalence
$$(\Shv(\CY)^{\on{constr}})^{\on{op}} \overset{\BD}\to \Shv(\CY)^{\on{constr}}.$$

If $\CY$ is not quasi-compact, the functor $\BD$ will typically \emph{not} send
$\Shv(\CY)^c$ to $\Shv(\CY)^c$.

\sssec{} \label{sss:duality adapted}

Assume that $\CY$ is quasi-compact. 
We will say that  $\CY$ is \emph{Verdier-compatible} if the functor $\BD$ sends $\Shv(\CY)^c$ to $\Shv(\CY)^c$.

\medskip

Based in \cite[Corollary 8.4.2]{DrGa1}, we conjecture:

\begin{conj} \label{c:duality preserve compact}
Any quasi-compact algebraic stack with an affine diagonal 
is Verdier-compatible.
\end{conj}

We are going to prove:

\begin{thm} \label{t:global quotient}
Let $\CY$ be such that it can be covered by open subsets each of which has the form 
$Y/G$, where $Y$ is a quasi-compact scheme and $G$ is an algebraic group. 
Then $\CY$ is Verdier-compatible.
\end{thm} 

\ssec{Proof of \thmref{t:global quotient}}

\sssec{A reduction step} \label{sss:cover by quot}

Let us reduce the assertion to the case when $\CY$ is globally a quotient, i.e., is of the form $Y/G$.

\medskip

Indeed, suppose $\CY$ can be covered by open substacks $\CU_i\overset{j_i}\hookrightarrow \CY$, such that each $\CU_i$ is
Verdier-compatible. We will show that $\CY$ is Verdier-compatible.

\medskip

Since $\CY$ was assumed quasi-compact, we can assume that the above open cover is finite.
Now the assertion follows from the fact that an object $\CF\in \Shv(\CY)$ is compact if and only if
all $j_i^*(\CF)$ are compact. 

\medskip

Indeed, the implication
$$\CF\in \Shv(\CY)^c\,\, \Rightarrow \,\, j_i^*(\CF)\in \Shv(\CU_i)^c$$
follows from the fact that $j_i^*$ admits a continuous right adjoint, namely, $(j_i)_*$.

\medskip

The opposite implication follows from the fact that $\CHom_{\Shv(\CY)}(\CF,-)$
can be expressed as a finite limit in terms of $\CHom_{\Shv(\CU_i)}(j_i^*(\CF),-)$
and finite intersections of these opens. 

\sssec{Explicit generators for a global quotient} \label{sss:gen glob quot}

Thus, we can assume that $\CY$ has the form $Y/G$.

\medskip

It suffices to find a system of compact generators of $\Shv(\CY)$ that are sent to compact objects 
by the functor $\BD$.

\medskip

Let
$\pi_Y$ denote the map 
$$Y\to Y/G=\CY.$$

\medskip

Note that for any $\CF\in \Shv(Y/G)^{\on{constr}}$, the object 
$$(\pi_Y)_!\circ (\pi_Y)^*(\CF)$$ 
is compact. Hence, it suffices to show that:

\medskip

\noindent(I) Such objects generate $\Shv(Y/G)$;

\medskip

\noindent(II) They are sent to compact objects by Verdier duality.

\sssec{Digression: cochains on the group} \label{sss:cohomology of group}

Consider the map 
$$\pi_{\on{pt}}:\on{pt}\to \on{pt}/G$$
and the objects
$$(\pi_{\on{pt}})_*(\sfe),(\pi_{\on{pt}})_!(\sfe)\in \Shv(\on{pt}/G).$$
Note that that
$$(\pi_{\on{pt}})^*\circ (\pi_{\on{pt}})_*(\sfe)\simeq \on{C}^\cdot(G).$$

Note also that
\begin{equation} \label{e:group cohomology}
(\pi_{\on{pt}})_*(\sfe) \simeq (\pi_{\on{pt}})_!(\sfe)[d], 
\end{equation}
where for 
$$1\to G_{\on{unip}}\to G\to G_{\on{red}}\to 1,$$
we have
$$d=2\dim(G_{\on{unip}})+\dim(G_{\on{red}}).$$

The isomorphism \eqref{e:group cohomology} follows from the fact that for a reductive group $G$,
the DG algebra of cochains $\on{C}^\cdot(G)$ is a Frobenius algebra (in fact, a symmetric algebra on generators in odd degrees),
and so $\on{C}_c^\cdot(G)[2\dim(G)]\simeq (\on{C}^\cdot(G))^\vee$ is isomorphic to 
$\on{C}^\cdot(G)$ up to a shift by $[d]$. 

\sssec{Verification of Property II} \label{sss:prop glob quot new a}

Denote by $q$ the map $Y/G\to \on{pt}/G$. For any $\CF'\in \Shv(Y/G)$ we have:
$$(\pi_Y)_!\circ (\pi_Y)^*(\CF')\simeq \CF'\overset{*}\otimes q^*((\pi_{\on{pt}})_!(\sfe))$$
and
\begin{multline*}
(\pi_Y)_*\circ (\pi_Y)^!(\CF')\simeq \CF'\sotimes q^!((\pi_{\on{pt}})_*(\sfe))\simeq
\CF' \overset{*}\otimes q^*((\pi_{\on{pt}})_*(\sfe))[2\dim(G)]\simeq \\
\simeq \CF'\overset{*}\otimes q^*((\pi_{\on{pt}})_!(\sfe))[2\dim(G)+d]\simeq
(\pi_Y)_!\circ (\pi_Y)^*(\CF')[2\dim(G)+d]
\end{multline*}

Hence,
$$\BD((\pi_Y)_!\circ (\pi_Y)^*(\CF)) \simeq 
(\pi_Y)_*\circ (\pi_Y)^!(\BD(\CF)) \simeq (\pi_Y)_!\circ (\pi_Y)^*(\BD(\CF))[2\dim(G)+d].$$

This proves Property (II). 

\sssec{Verification of Property I} \label{sss:prop glob quot new b}

To prove Property (I), let $\CF'$ be a non-zero object of $\Shv(Y/G)$, and let us find $\CF\in \Shv(Y/G)^{\on{constr}}$
so that
$$\CHom_{\Shv(Y/G)}((\pi_Y)_!\circ (\pi_Y)^*(\CF),\CF')\neq 0.$$

\medskip

This can be done for $Y\overset{\pi_Y}\to Y/G$ replaced by any pair $Y\overset{f}\to \CZ$,
where $\CZ$ is an algebraic stack and $f$ is a smooth covering map. 

\medskip

Indeed, for any $\CF,\CF'\in \Shv(\CZ)$, we have
\begin{equation} \label{e:F F'}
\CHom_{\Shv(\CZ)}(f_!\circ f^*(\CF),\CF')\simeq \CHom_{\Shv(Y)}(f^*(\CF),f^!(\CF'))\simeq 
\CHom_{\Shv(\CZ)}(\CF,f_*\circ f^!(\CF')).
\end{equation} 

Since $\Shv(\CZ)$ is compactly generated, and  
$\Shv(\CZ)^c\subset \Shv(\CZ)^{\on{constr}}$, it suffices to show that if $\CF'\neq 0$, then $f_*\circ f^!(\CF')\neq 0$.

\medskip

Applying \eqref{e:F F'} to $\CF=\CF'$, it suffices to show that $\CHom_{\Shv(Y)}(f^*(\CF'),f^!(\CF'))\neq 0$.
However, $f^!$ is isomorphic to $f^*$ up to a shift, so the assertion follows from the fact that $f^*(\CF')\neq 0$. 

\qed[\thmref{t:global quotient}]

\ssec{Verdier duality on stacks}

In this subsection $\CY$ will be a Verdier-compatible quasi-compact algebraic stack. 

\sssec{} \label{sss:duality on stack}

The assumption that $\CY$ is Verdier-compatible implies that the Verdier duality functor defines a contravariant equivalence
$$(\Shv(\CY)^c)^{\on{op}} \to \Shv(\CY)^c.$$

Hence, we obtain a canonical identification
$$\Shv(\CY)^\vee\simeq \Shv(\CY).$$

By construction, the corresponding pairing
$$\Shv(\CY)^c\times \Shv(\CY)^c\to \Vect_\sfe$$
sends
$$\CF_1,\CF_2\to \on{C}^\cdot(\CY,\CF_1\sotimes \CF_2).$$

\sssec{} \label{sss:black triangle}

Let
$$\on{C}^\cdot_\blacktriangle(\CY,-):\Shv(\CY)\to \Vect$$
be the functor dual to the functor
$$\Vect_\sfe\to \Shv(\CY), \quad \sfe\mapsto \omega_{\CY},$$
see \cite[Sect. 9.1]{DrGa1}. This functor is the ind-extension of the restriction of the functor
$$\on{C}^\cdot(\CY,-):\Shv(\CY)\to \Vect$$
to $\Shv(\CY)^c\subset \Shv(\CY)$. 

\medskip

In particular, we have a natural transformation
\begin{equation} \label{e:from triangle}
\on{C}^\cdot_\blacktriangle(\CY,-)\to \on{C}^\cdot(\CY,-),
\end{equation} 
which is an equivalence when evaluated on compact objects.

\medskip

Furthermore, the duality pairing on all of $\Shv(\CY)\otimes \Shv(\CY)$ can be written as 
$$\CF_1,\CF_2 \mapsto \on{C}^\cdot_\blacktriangle(\CY,\CF_1\sotimes \CF_2).$$

Applying \eqref{e:from triangle}, we obtain a map
\begin{equation} \label{e:triang to star ten}
\on{C}^\cdot_\blacktriangle(\CY,\CF_1\sotimes \CF_2)\to \on{C}^\cdot(\CY,\CF_1\sotimes \CF_2). 
\end{equation}

\sssec{}

We observe:

\begin{lem} \label{l:ten prod by compact}
For $\CF\in \Shv(\CY)^c$ and $\CF'\in \Shv(\CY)^{\on{constr}}$, both
$$\CF\overset{*}\otimes \CF' \text{ and } \CF\sotimes \CF'$$ 
are compact.
\end{lem}

\begin{proof}
The assertion for  $\CF\overset{*}\otimes \CF'$ follows from the fact that
$$\CHom(\CF\overset{*}\otimes \CF',\CF'')\simeq \CHom(\CF,\BD(\CF')\sotimes \CF'').$$

The assertion for $\CF\sotimes \CF'$ follows by Verdier duality (and the assumption that $\CY$ is Verdier-compatible).
\end{proof} 

\begin{cor}
The map \eqref{e:triang to star ten} is an isomorphism when one of the objects $\CF_1$ or $\CF_2$ is compact.
\end{cor}

\begin{proof}

By \lemref{l:ten prod by compact}, 
the map \eqref{e:triang to star ten} is an isomorphism if both $\CF_1$ or $\CF_2$ are compact.

\medskip

Now, the assertion follows from the fact that if $\CF_1$ is compact, then both sides in \eqref{e:triang to star ten}
commute with colimits in $\CF_2$: this is the case for the left-hand side by definition, and for the right-hand side since
$$\on{C}^\cdot(\CY,\CF_1\sotimes \CF_2) \simeq \CHom(\BD(\CF_1),\CF_2),$$
and $\BD(\CF_1)$ by the assumption that $\CY$ is Verdier-compatible.

\end{proof}

%
%
%

\sssec{}

For future reference, we record the following properties of Verdier-compatible prestacks,
borrowed from \cite[Theorem 10.2.9]{DrGa1} (we will omit the proof as it repeats
verbatim the arguments from {\it loc. cit.}):

\begin{prop} \label{p:Drinf prop}
Assume that $\CY$ is Verdier-compatible. Then 
the following properties of an object $\CF\in \Shv(\CY)^{\on{constr}}$
are equivalent:

\smallskip

\noindent{\em(i)} $\CF$ is compact;

\smallskip

\noindent{\em(i')} $\BD(\CF)$ is compact;

\smallskip

\noindent{\em(ii)} $\CF$ belongs to the smallest (non cocomplete) DG subcategory of $\Shv(\CY)$
closed under taking direct summands that contains objects of the form $f_!(\CF_S)$,
where $f:S\to \CY$ with $S$ an affine scheme and $\CF_S\in \Shv(S)^{\on{constr}}$;

\smallskip

\noindent{\em(ii')} $\CF$ belongs to the smallest (non cocomplete) DG subcategory of $\Shv(\CY)$
closed under taking direct summands that contains objects of the form $f_*(\CF_S)$,
where $f:S\to \CY$ with $S$ an affine scheme and $\CF_S\in \Shv(S)^{\on{constr}}$;

\smallskip

\noindent{\em(iii)} The functor
$$\CF'\mapsto \on{C}^\cdot(\CY,\CF\sotimes \CF'), \quad \Shv(\CY)\to \Vect_\sfe$$
is continuous;

\smallskip

\noindent{\em(iv)} The functor
$$\CF'\mapsto \on{C}^\cdot(\CY,\CF\sotimes \CF'), \quad \Shv(\CY)\to \Vect_\sfe$$
is cohomologically bounded on the right;

\smallskip

\noindent{\em(v)} The functor
$$\CF'\mapsto \on{C}_c^\cdot(\CY,\CF\overset{*}\otimes \CF'), \quad \Shv(\CY)\to \Vect_\sfe$$
is cohomologically bounded on the left;

\smallskip

\noindent{\em(vi)} The functor
$$\CF'\mapsto \on{C}_\blacktriangle^\cdot(\CY,\CF\sotimes \CF'), \quad \Shv(\CY)\to \Vect_\sfe$$
is cohomologically bounded on the left;

\smallskip

\noindent{\em(vii)} The natural transformation
$$\on{C}_\blacktriangle^\cdot(\CY,\CF\sotimes \CF')\to  \on{C}^\cdot(\CY,\CF\sotimes \CF'), \quad \CF'\in \Shv(\CY)$$
is an isomorphism;

\smallskip

\noindent{\em(viii)} For any schematic quasi-compact morphism $g:\CY'\to \CY$ and $f:\CY'\to S$
where $S$ is a scheme, the object $f_*\circ g^!(\CF)$ is cohomologically bounded above;

\smallskip

\noindent{\em(viii')} For any schematic quasi-compact morphism $g:\CY'\to \CY$ and $f:\CY'\to S$
where $S$ is a scheme, the object $f_!\circ g^*(\CF)$ is cohomologically bounded below;

\smallskip

\noindent{\em(ix)} Same as \emph{(viii)} but $g$ is a finite \'etale map onto a locally closed substack of $\CY$;

\smallskip

\noindent{\em(ix')} Same as \emph{(viii')} but $g$ is a finite \'etale map onto a locally closed substack of $\CY$.

\end{prop} 

%
%
%
%
%
%
%
%

\ssec{The renormalized category of sheaves} \label{ss:ren sheaves}

\sssec{} \label{sss:ren}

Let $\CY$ be a quasi-compact algebraic stack.

\medskip

We define the renormalized version of the category of sheaves on $\CY$, denoted $\Shv(\CY)^{\on{ren}}$ to be
$$\on{Ind}(\Shv(\CY)^{\on{constr}}).$$

\medskip

The t-structure on $\Shv(\CY)^{\on{constr}}$ induces a unique t-structure on $\Shv(\CY)^{\on{ren}}$
compatible with filtered colimits. 

\medskip

Ind-extension of the tautological embedding $\Shv(\CY)^{\on{constr}}\hookrightarrow \Shv(\CY)$ defines a functor
$$\on{un-ren}_\CY:\Shv(\CY)^{\on{ren}}\to \Shv(\CY).$$
The functor $\on{un-ren}_\CY$ is t-exact by construction.

\sssec{} \label{sss:Shv ren on stack}

Note that for every fixed $n$, the functor 
\begin{equation} \label{e:un-ren n}
\on{un-ren}_\CY:(\Shv(\CY)^{\on{ren}})^{\geq -n}\to (\Shv(\CY))^{\geq -n}
\end{equation} 
preserves compactness (by \propref{p:almost compact}), and hence is fully faithful. 

\medskip
 
Moreover, since the functor $\on{un-ren}_\CY$ is essentially surjective on the hearts,
we obtain that \eqref{e:un-ren n} is actually an equivalence.

\medskip

From here, we obtain that the  functor $\on{un-ren}_\CY$ identifies $\Shv(\CY)$ with the left completion of $\Shv(\CY)^{\on{ren}}$
with respect to its t-structure. 

\sssec{}

This construction of the pair $(\Shv(\CY)^{\on{ren}},\on{un-ren}_\CY)$ 
mimics the construction of how one defines $\IndCoh(S)$ for an eventually coconnective
affine scheme (see \cite[Chapter 4, Sect. 1.2]{GR1}), and shares its formal properties:

\medskip

\begin{itemize}

\item Ind-extension of the tautological embedding $\Shv(\CY)^c\hookrightarrow \Shv(\CY)^{\on{constr}}$ defines a fully
faithful functor
$$\on{ren}_{\CY}:\Shv(\CY)\to \Shv(\CY)^{\on{ren}},$$
which is the left adjoint of $\on{un-ren}_\CY$.

\medskip

\item The functor $\on{un-ren}_{\CY}$ realizes $\Shv(\CY)$ as the co-localization of $\Shv(\CY)^{\on{ren}}$ with respect to
the subcategory consisting of objects all of whose cohomologies with respect to the above t-structure vanish.

\smallskip

\item The operation of *-tensor product makes $\Shv(\CY)^{\on{ren}}$ into a symmetric monoidal category,
and $\Shv(\CY)$ into a module category over it (see \lemref{l:ten prod by compact}). The same is true
for the !-tensor product provided that $\CY$ is Verdier-compatible. 

\end{itemize} 

\sssec{}

Note that Verdier duality
$$\BD:(\Shv(\CY)^{\on{constr}})^{\on{op}} \to \Shv(\CY)^{\on{constr}}$$
defines an identification
$$\Shv(\CY)^{\on{ren}}\simeq (\Shv(\CY)^{\on{ren}})^\vee.$$

Assume for a moment that $\CY$ is Verdier-compatible. In particular, we have also the identification
$$\Shv(\CY)\simeq \Shv(\CY)^\vee.$$

The functors $\on{ren}_\CY$ and $\on{un-ren}_\CY$ are mutually dual with respect to these identifications. 

\sssec{}

Let us consider the example of $\CY=\on{pt}/G$. In this case $\Shv(\on{pt}/G)^{\on{ren}}$ is compactly generated
by the object $\ul\sfe_{\on{pt}/G}$. Hence, we obtain a canonical equivalence
\begin{equation} \label{e:on BG}
\Shv(\on{pt}/G)^{\on{ren}}\simeq \on{C}^\cdot(\on{pt}/G)\mod.
\end{equation} 

\medskip

Under this equivalence, the symmetric monoidal structure on $\Shv(\on{pt}/G)^{\on{ren}}$ given by *-tensor
product corresponds to the usual symmetric monoidal structure on the category of modules over a
commutative algebra. 

\medskip

Recall that $\on{C}^\cdot(\on{pt}/G)$ is isomorphic to a polynomial algebra on generators in even degrees.
The canonical point
$$\pi_{\on{pt}}:\on{pt}\to \on{pt}/G$$
defines an augmentation module
$$\sfe\in  \on{C}^\cdot(\on{pt}/G)\mod.$$

\medskip

Note that under the equivalence \eqref{e:on BG}, we have
$$\sfe\in \on{C}^\cdot(\on{pt}/G)\mod\, \leftrightarrow\, \on{ren}_{\on{pt}/G}((\pi_{\on{pt}})_*(\sfe))\in \Shv(\on{pt}/G)^{\on{ren}}.$$

\medskip

Hence, under \eqref{e:on BG}, the (isomorphic) essential image of the functor $\on{ren}_{\on{pt}/G}$ corresponds to the full
subcategory 
$$\on{C}^\cdot(\on{pt}/G)\mod_0\subset \on{C}^\cdot(\on{pt}/G)\mod$$
be the full subcategory generated by the the augmentation module $\sfe$. 

\sssec{}

Let $\CY$ be of the form $Y/G$, where $Y$ is a quasi-compact scheme. The functor of *- (resp., !-)
pullback 
$$\Shv(\on{pt}/G)^{\on{ren}}\to \Shv(\CY)^{\on{ren}}$$
has a natural symmetric monoidal structure with respect to the *- (resp., !-) tensor product operation.

\medskip

We claim:

\begin{prop} 
The co-localization
$$\on{un-ren}_\CY:\Shv(\CY)^{\on{ren}}\leftrightarrows \Shv(\CY):\on{ren}_\CY$$
identifies with the co-localization 
$$\Shv(\CY)^{\on{ren}}\simeq \Shv(\CY)^{\on{ren}}\underset{\Shv(\on{pt}/G)^{\on{ren}}}\otimes \Shv(\on{pt}/G)^{\on{ren}}
\leftrightarrows \Shv(\CY)^{\on{ren}}\underset{\Shv(\on{pt}/G)^{\on{ren}}}\otimes \Shv(\on{pt}/G)$$
(for either *- or !- monoidal structures).
\end{prop}

\begin{proof}

The functor
$$\on{un-ren}_\CY: \Shv(\CY)^{\on{ren}}\to  \Shv(\CY)$$ 
clearly factors as
\begin{multline*}
\Shv(\CY)^{\on{ren}} 
\simeq \Shv(\CY)^{\on{ren}}\underset{\Shv(\on{pt}/G)^{\on{ren}}}\otimes \Shv(\on{pt}/G)^{\on{ren}}
\overset{\on{Id}\otimes \on{un-ren}_{\on{pt}/G}}\longrightarrow \\
\to \Shv(\CY)^{\on{ren}}\underset{\Shv(\on{pt}/G)^{\on{ren}}}\otimes \Shv(\on{pt}/G)\to \Shv(\CY).
\end{multline*}

Hence, to prove the proposition it suffices to show that the essential image of 
$$\Shv(\CY)^{\on{ren}}\underset{\Shv(\on{pt}/G)^{\on{ren}}}\otimes \Shv(\on{pt}/G)\overset{\on{Id}\otimes \on{ren}_{\on{pt}/G}}\longrightarrow 
\Shv(\CY)^{\on{ren}}\underset{\Shv(\on{pt}/G)^{\on{ren}}}\otimes \Shv(\on{pt}/G)^{\on{ren}}\simeq 
\Shv(\CY)^{\on{ren}}$$ 
is contained in that of
$$\Shv(\CY) \overset{\on{ren}_\CY}\hookrightarrow \Shv(\CY)^{\on{ren}}.$$

\medskip

For that end it suffices to show that for $\CF\in \Shv(\CY)^{\on{constr}}$, we have
$$\CF\overset{*}\otimes q^*((\pi_{\on{pt}})_*(\sfe))\in \Shv(\CY)^c,$$
where $q:\CY\to \on{pt}/G$. However, this follows from \eqref{e:group cohomology}. 

\end{proof} 

\sssec{}

Let now $\CY$ be a not necessarily quasi-compact algebraic stack. We let
\begin{equation} \label{e:ren on non qc}
\Shv(\CY)^{\on{ren}}:=\underset{\CU}{\on{lim}}\, \Shv(\CU)^{\on{ren}},
\end{equation}
where the limit is taken over the index category of quasi-compact open substacks $\CU\subset \CY$,
and the transition functors are given by restriction.

\medskip

The properties and structures listed in \secref{sss:ren} for the opens $\CU$ induce the corresponding 
properties and structures on $\CY$. In particular, we have an adjunction 
$$\on{un-ren}_\CY:\Shv(\CY)^{\on{ren}}\leftrightarrows \Shv(\CY):\on{ren}_\CY,$$
with $\on{ren}_\CY$ fully faithful, a t-structure on $\Shv(\CY)^{\on{ren}}$, etc.

\medskip

Note also that the transition functors in forming the limit \eqref{e:ren on non qc} admit left adjoints,
given by !-extension. Hence, we can rewrite $\Shv(\CY)^{\on{ren}}$ as
$$\underset{\CU}{\on{colim}}\, \Shv(\CU)^{\on{ren}},$$
where the transition functors are given by !-extension.

\medskip

In particular, we obtain that $\Shv(\CY)^{\on{ren}}$ is compactly generated by objects of the form 
$j_!(\CF)$, where 
$$j:\CU \hookrightarrow \CY$$
with $\CU$ quasi-compact and $\CF\in \Shv(\CU)^{\on{constr}}$. 

\ssec{Singular support on stacks} \label{ss:sing supp constr}

Let $\CY$ be an algebraic stack. 

\sssec{} \label{sss:cotan smooth stack}

Assume first that $\CY$ is smooth. Then one can talk about $T^*(\CY)$, which we regard
as a \emph{classical} algebraic stack. Namely, let $T(\CY)$ be the $\CO_\CY$-linear dual
of the cotangent complex of $\CY$, so that $T(\CY)$ can locally be written as 
$$\on{Cone}(E_{-1}\to E_0),$$
where $E_{-1},E_0$ are locally free sheaves on $\CY$, and 
$$T^*(\CY)=\Spec\left(H^0\left(\Sym_{\CO_\CY}(T(\CY))\right)\right).$$

\medskip

We could of course consider a derived enhancement of $T^*(\CY)$, namely, 
$$\Spec\left(\Sym_{\CO_\CY}(T(\CY))\right),$$
but we will not need it for the purposes of the present paper.

\sssec{} \label{sss:sing supp sing stack} 

Let now $\CY$ be arbitrary (i.e., not necessarily smooth). In this case,  
following \cite[Sect. A.3.6]{GKRV}, we will not even attempt to define 
$T^*(\CY)$ is an algebro-geometric object.
 
\medskip
 
Instead, we will talk about (isomorphism classes of) $k$-points of $T^*(\CY)$, which are, by definition, pairs
$(\xi,y)$, where $y\in \CY(k)$ and $\xi\in H^0(T^*_y(\CY))$.

\medskip

In addition, one can talk about \emph{Zariski-closed subsets} of $T^*(\CY)$.
By definition, such a subset $\CN$ is a compatible collection
of Zariski-closed subsets 
$$\CN_S\subset  T^*(S)$$
for schemes $S$ equipped with a smooth map $S\to \CY$ (see \secref{ss:sing supp sing}
for the notion of a Zariski-closed subset in a cotangent bundle of a non-smooth
scheme). 

\medskip

We will write symbolically that $\CN_S$ is the image of 
$$\CN\underset{\CY}\times S\subset T^*(\CY)\underset{\CY}\times S$$
under the codifferential 
$$T^*(\CY)\underset{\CY}\times S\to T^*(S).$$

\sssec{} \label{sss:sing supp stack qc}

Let $\CN$ be a conical Zariski-closed subset in $T^*(\CY)$, defined as above.

\medskip

We define the full subcategory
$$\Shv_{\CN}(\CY)\subset \Shv(\CY)$$
to consist of those $\CF\in \Shv(\CY)$ whose pullback under any smooth
map $S\to \CY$ (with $S$ a scheme) belongs to $\Shv_{\CN_S}(S)$ 
(see \secref{sss:sing supp schemes} in the case when $S$ is smooth and \secref{ss:sing supp sing} for a general $S$),
where 
$$\CN_S:=\CN\underset{\CY}\times S\subset T^*(\CY)\underset{\CY}\times S\subset T^*(S),$$
see \secref{sss:sing supp sing stack}.

\medskip

Vice versa, to $\CF\in \Shv(\CY)^{\on{constr}}$, we can associate its singular support,
which is a Zariski-closed subset of $T^*(\CY)$.

\sssec{}

Since  pullbacks with respect to smooth morphisms are t-exact up to a cohomological shift, we obtain that
the t-structures on $\Shv_{\CN_S}(S)$ induce a t-structure on $\Shv_{\CN}(\CY)$.
It follows automatically that $\Shv_{\CN}(\CY)$ is complete in its t-structure.

\medskip

It is easy to see that an object $\CF\in\Shv(\CY)$ belongs to $\Shv_{\CN}(\CY)$ if and only if for every $m$
and every constructible sub-object $\CF'$ of $H^m(\CF)$, the object $\CF'$ belongs to 
$$\on{Perv}_{\CN}(\CY):=\Shv_{\CN}(\CY)\cap \on{Perv}(\CY).$$

\medskip

Note also that $\on{Perv}_{\CN}(\CY)$ is a Serre subcategory of $\on{Perv}(\CY)$ and
$$\on{Ind}(\on{Perv}_{\CN}(\CY))\simeq \Shv_{\CN}(\CY)^\heartsuit.$$

\begin{rem}
The category $\Shv_{\CN}(\CY)$ defined above is different from the category denoted by the same symbol in 
\cite{GKRV}. In our current notations, the category in {\it loc.cit.} is
$$\underset{S}{\on{lim}}\, \Shv_{\CN_S}(S)^{\on{access}},$$
where the limit taken over the category of affine schemes $S$ smooth over $\CY$ and smooth maps between such,
and where $\Shv_{\CN_S}(S)^{\on{access}}$ is as in \secref{sss:sing supp schemes constr}. 

\medskip

Since for schemes, the functor $\Shv_{\CN_S}(S)^{\on{access}}\to \Shv_{\CN_S}(S)$ is fully faithful, the category
in \cite{GKRV} embeds fully faithfully into our $\Shv_{\CN}(\CY)$.
\end{rem} 

\sssec{}

We will now assume, for the duration of this subsection, that $\CY$ is quasi-compact.

\medskip

Set
$$\Shv_{\CN}(\CY)^{\on{constr}}:=\Shv_{\CN}(\CY)\cap \Shv(\CY)^{\on{constr}}$$
and define
$$\Shv_{\CN}(\CY)^{\on{ren}}:=\on{Ind}(\Shv_\CN(\CY)^{\on{constr}}).$$
Note that we have a tautologically defined functor
\begin{equation}  \label{e:unren N}
\Shv_{\CN}(\CY)^{\on{ren}}\to \Shv_{\CN}(\CY).
\end{equation} 

\medskip

The category $\Shv_{\CN}(\CY)^{\on{constr}}$ inherits a t-structure, and the latter uniquely extends to a t-structure
on $\Shv_{\CN}(\CY)^{\on{ren}}$, compatible with filtered colimits. The functor \eqref{e:unren N} is t-exact, by 
construction. 

\medskip

As in \secref{sss:Shv ren on stack} one shows that for every $n$, the functor \eqref{e:unren N} induces an equivalence
\begin{equation}  \label{e:unren N coconn}
(\Shv_{\CN}(\CY)^{\on{ren}})^{\geq -n} \simeq (\Shv_{\CN}(\CY))^{\geq -n}.
\end{equation} 

From here it follows that the functor \eqref{e:unren N} realizes $\Shv_{\CN}(\CY)$ 
as the left completion of $\Shv_{\CN}(\CY)^{\on{ren}}$ in its t-structure.

\sssec{}

Note also that we have a fully faithful t-exact functor
$$\Shv_{\CN}(\CY)^{\on{ren}}\to \Shv(\CY)^{\on{ren}}.$$

\begin{rem}

Note that when $\CY=Y$ is a quasi-compact scheme, the category that denoted above by 
$\Shv_{\CN}(\CY)^{\on{ren}}$ was denoted by $\Shv_{\CN}(Y)^{\on{access}}$
in \secref{sss:sing supp schemes constr}. 

\medskip

In the case of stacks, the notation $\Shv_{\CN}(\CY)^{\on{access}}$ will have a different
meaning, see \secref{sss:accessible} below.

\medskip

By contrast, when $\CN$ is all of $T^*(\CY)$, the category $\Shv_{\CN}(\CY)^{\on{ren}}$
is the category $\Shv(\CY)^{\on{ren}}$ introduced in \secref{ss:ren sheaves}. 

\medskip

Note also that, unlike the case of schemes, the functor \eqref{e:unren N} does not 
preserves compactness, and hence fails to be fully faithful. 

\end{rem}

\ssec{The accessible category on stacks}

We continue to assume that $\CY$ is a quasi-compact algebraic stack. 

\sssec{} \label{sss:accessible}

Denote by $\Shv_{\CN}(\CY)^{\on{access}}$ the full subcategory in $\Shv_{\CN}(\CY)$
generated under colimits by the essential image of the functor \eqref{e:unren N} above.

\medskip

Thus, we have the functors
\begin{equation} \label{e:access and ren}
\Shv_{\CN}(\CY)^{\on{ren}} \overset{\on{un-ren}_{\CY,\CN}}\longrightarrow \Shv_{\CN}(\CY)^{\on{access}} \hookrightarrow \Shv_{\CN}(\CY).
\end{equation} 

Since the functors \eqref{e:unren N coconn} are equivalences, we obtain that the subcategory
\begin{equation} \label{e:emb access}
\Shv_{\CN}(\CY)^{\on{access}}\hookrightarrow \Shv_{\CN}(\CY)
\end{equation} 
is preserved by the truncation functors acting on $\Shv_{\CN}(\CY)$. In particular, 
$\Shv_{\CN}(\CY)^{\on{access}}$ inherits a t-structure, and the embeddings 
$$(\Shv_{\CN}(\CY)^{\on{access}})^{\geq -n}\hookrightarrow (\Shv_{\CN}(\CY))^{\geq -n}$$
are equivalences. 

\medskip

In particular, the embedding \eqref{e:emb access} realizes $\Shv_{\CN}(\CY)$
as the left completion of $\Shv_{\CN}(\CY)^{\on{access}}$ in its t-structure.

\sssec{Examples}

When $\CY=Y$ is a scheme, the functor \eqref{e:unren N} is fully faithful, so the functor
$$\Shv_{\CN}(Y)^{\on{ren}} \overset{\on{un-ren}_{Y,\CN}}\longrightarrow \Shv_{\CN}(Y)^{\on{access}}$$
is an equivalence.

\medskip

When $\CY$ is a stack but $\CN=T^*(\CY)$, the embedding \eqref{e:emb access} is an equivalence. 

\sssec{} 

We give the following definitions:

\begin{defn}  \label{d:renormalization-adapted} 
We shall say that the pair $(\CY,\CN)$ is \emph{renormalization-adapted} if the category $\Shv_{\CN}(\CY)^{\on{access}}$
is generated by objects that are compact as objects of $\Shv(\CY)$.
\end{defn} 

\begin{defn} \label{d:constraccessible}
We shall say that the pair $(\CY,\CN)$ is \emph{constraccessible} if the inclusion \eqref{e:emb access}
is an equality. 
\end{defn}

\begin{rem} 
We emphasize again that a pair $(\CY,\CN)$ \emph{may not be} constraccessible even if $\CY=S$ is a scheme 
and $\CN=\{0\}$ (see Remark \ref{r:bad P1}). But it is tautologically renormalization-adapted. 
\end{rem}

\sssec{} \label{sss:constraccessible new new}

We make the following few observations: 

\medskip

\noindent(I) If a pair $(\CY,\CN)$ is renormalization-adapted, then $\Shv_{\CN}(\CY)^{\on{access}}$ 
equals the full subcategory of $\Shv_{\CN}(\CY)$ generated by $\Shv_{\CN}(\CY)\cap \Shv(\CY)^c$.

\medskip

\noindent(II) A pair $(\CY,\CN)$ is both renormalization-adapted and constraccessible if and
only if $\Shv_{\CN}(\CY)$ is generated by objects that are compact in $\Shv(\CY)$. 

\medskip

\noindent(III) If a pair $(\CY,\CN)$ is renormalization-adapted, the functor $\on{un-ren}_{\CY,\CN}$ admits a left adjoint,
to be denoted $\on{ren}_{\CY,\CN}$. In fact, this left adjoint is the ind-extension of the tautological embedding
$$\Shv_{\CN}(\CY)\cap \Shv(\CY)^c\hookrightarrow  \Shv_{\CN}(\CY)^{\on{constr}}.$$
In particular, $\on{ren}_{\CY,\CN}$ is fully faithful, so the adjunction
$$\on{ren}_{\CY,\CN}:\Shv_{\CN}(\CY)^{\on{access}} \rightleftarrows \Shv_{\CN}(\CY)^{\on{ren}}:\on{un-ren}_{\CY,\CN}$$
realizes $\Shv_{\CN}(\CY)^{\on{access}}$ as a co-localization of $\Shv_{\CN}(\CY)^{\on{ren}}$. 

\medskip

\noindent(IV) If $(\CY,\CN)$ is renormalization-adapted and $\CY$ is Verdier-compatible, then the category
$\Shv_{\CN}(\CY)^{\on{access}}$ is naturally self-dual, with the pairing
$$\Shv_{\CN}(\CY)^{\on{access}}\otimes \Shv_{\CN}(\CY)^{\on{access}}\to \Vect_\sfe$$
given by
$$\CF_1,\CF_2\mapsto \on{C}^\cdot_\blacktriangle(\CY,\CF_1\sotimes \CF_2),$$
and the corresponding contravariant equivalence on compact objects is given by the Verdier duality
functor. 

\sssec{}

We claim:

\begin{prop} \label{p:global quotient N}
Suppose that $\CY$ is a global quotient, i.e., $\CU=Y/G$, where $Y$ is a quasi-compact scheme and $G$ an algebraic group. 
Then $(\CY,\CN)$ is renormalization-adapted for any $\CN$.
\end{prop} 

\begin{proof}

Follows from the argument in \secref{sss:prop glob quot new b}. 

\end{proof} 

Based on the above proposition, we propose: 

\begin{conj} \label{c:ren adapt}
For any quasi-compact algebraic stack with an affine diagonal, 
and any $\CN\subset T^*(\CY)$, the pair $(\CY,\CN)$ is renormalization-adapted.
\end{conj} 

\ssec{Singular support condition for non quasi-compact stacks}

In this subsection we let $\CY$ be an algebraic stack, locally of finite type,
but not necessarily quasi-compact. 

\sssec{} \label{sss:non qc N}

We define 
$$\Shv_{\CN}(\CY):= \underset{\CU}{\on{lim}}\, \Shv_\CN(\CU),$$
where the index category is the poset of quasi-compact open substacks of $\CY$,
and the transition functors used in forming the limit are given by restriction. 

\medskip

The t-structures for the individual $\CU's$ induces a t-structure on $\Shv_{\CN}(\CY)$. 

\medskip

The fully faithful embeddings $\Shv_\CN(\CU)\to \Shv(\CU)$ give rise to a t-exact fully faithful
embedding
$$\Shv_\CN(\CY)\to \Shv(\CY).$$

\sssec{}  \label{sss:non qc N access}

We define the categories 
\begin{equation} \label{e:access non qc}
\Shv_{\CN}(\CY)^{\on{ren}} \text{ and } \Shv_{\CN}(\CY)^{\on{access}}
\end{equation} 
similarly:
$$\Shv_{\CN}(\CY)^{\on{ren}}:=\underset{\CU}{\on{lim}}\, \Shv_\CN(\CU)^{\on{ren}} \text{ and }
\Shv_{\CN}(\CY)^{\on{access}}:=\underset{\CU}{\on{lim}}\, \Shv_\CN(\CU)^{\on{access}}.$$

The t-structures for the individual $\CU's$ induces t-structures on the above categories.

\medskip

We have a t-exact functor
$$\on{un-ren}_{\CY,\CN}:\Shv_{\CN}(\CY)^{\on{ren}} \to \Shv_{\CN}(\CY)^{\on{access}}$$
and a t-exact fully faithful embedding
\begin{equation} \label{e:emb access non qc}
\Shv_{\CN}(\CY)^{\on{access}}\hookrightarrow \Shv_{\CN}(\CY).
\end{equation} 

Both these functors induce equivalences on the $n$-coconnective subcategories for any $n$. 
This implies both these functors identify $\Shv_{\CN}(\CY)$ with the left completion of
both $\Shv_{\CN}(\CY)^{\on{ren}}$ and $\Shv_{\CN}(\CY)^{\on{access}}$.

\sssec{}

The fully faithful embeddings $\Shv_\CN(\CU)^{\on{ren}} \hookrightarrow \Shv(\CU)^{\on{ren}}$
give rise to a t-exact fully faithful embedding 
$$\Shv_\CN(\CY)^{\on{ren}} \hookrightarrow \Shv(\CY)^{\on{ren}}.$$

\sssec{}  \label{sss:non qc N constraccess}

We shall say that $(\CY,\CN)$ is renormalization-adapted if $\Shv_{\CN}(\CY)^{\on{access}}$ is generated by objects
that are compact in $\Shv(\CY)$.

\medskip

We shall say that $(\CY,\CN)$ is constraccessible if the inclusion \eqref{e:emb access non qc} is an equality. 

\medskip

Observations (I), (II) and (III) from \secref{sss:constraccessible new new} remain valid for non quasi-compact
stacks as well. 

\begin{rem} 
Note that the notions of being renormalization-adapted and constraccessible, as defined above,
are of global nature, i.e., in general, they do not translate into statements that can be checked
on quasi-compact open substacks of $\CY$.

\medskip

We will now introduce a condition on $\CY$ that allows us to check these properties locally.
\end{rem} 

\sssec{} \label{sss:N preserved abs a}

We shall say that a pair $\CY$ is $\CN$-\emph{truncatable}\footnote{The terminology is borrowed from
\cite[Sect. 4]{DrGa2}.} if $\CY$ can be written as a filtered union of quasi-compact 
open substacks $\CU_i$, such that for every $i$ and the corresponding open embedding
$$\CU_i\overset{j_i}\hookrightarrow \CY,$$
the functor
$$(j_i)_!:\Shv(\CU_i)^{\on{constr}}\to \Shv(\CY)^{\on{constr}}$$
sends 
$$\Shv_\CN(\CU_i)^{\on{constr}}\to \Shv_\CN(\CY)^{\on{constr}}.$$

\sssec{} \label{sss:N preserved abs b}

Here are some equivalent ways to rewrite the $\CN$-truncatability condition: 

\medskip

\noindent(A) We can require that the functors $(j_{i_1,i_2})_!$ send 
$\Shv_\CN(\CU_{i_1})^{\on{constr}}\to \Shv_\CN(\CU_{i_2})^{\on{constr}}$ for $U_1\overset{j_{i_1,i_2}}\hookrightarrow U_2$.

\medskip

\noindent(B) We can require that the functors $(j_{i_1,i_2})_!$ send 
$\Shv_\CN(\CU_{i_1})\to \Shv_\CN(\CU_{i_2})$. Indeed, (B) $\Rightarrow$ (A) tautologically. For
the converse implication, using the fact that the functor $(j_{i_1,i_2})_!$ has a finite cohomological 
amplitude, it is enough to show that $(j_{i_1,i_2})_!$ sends $\on{Perv}_\CN(\CU_{i_1})$ to
$\Shv_\CN(\CU_{i_2})$, which follows from (A). 
 
\medskip

\noindent(C) We can require that the functors $(j_i)_!$ send $\Shv_\CN(\CU_i)\to \Shv_\CN(\CY)$
(indeed, this is easily seen to be equivalent to (B)); 

\medskip

\noindent(D) Same as all of the above with $(-)_!$ replaced by $(-)_*$. Indeed, for the constructible
categories, this follows by Verdier duality, and for the entire $\Shv_\CN(-)$, this follows as in 
the equivalence (A) $\Leftrightarrow$ (B) above. 

\sssec{} \label{sss:N preserved abs c}

Let $\CY$ and $\CU_i$ be as in \secref{sss:N preserved abs a}. 
Note that in this case, for $U_1\overset{j_{i_1,i_2}}\hookrightarrow U_2$ the functors 
$$(j_{i_1,i_2})^*:\Shv_\CN(\CU_{i_2})\to \Shv_\CN(\CU_{i_1}),$$ 
$$(j_{i_1,i_2})^*:\Shv_\CN(\CU_{i_2})^{\on{ren}}\to \Shv_\CN(\CU_{i_1})^{\on{ren}}$$
and 
$$(j_{i_1,i_2})^*:\Shv_\CN(\CU_{i_2})^{\on{access}}\to \Shv_\CN(\CU_{i_1})^{\on{access}}$$
all admit fully faithful left adjoints, to be denoted $(j_{i_1,i_2})_!$ in all three cases. 

\medskip

This implies that the corresponding categories 
$$\Shv_{\CN}(\CY), \,\, \Shv_{\CN}(\CY)^{\on{ren}} \text{ and } \Shv_{\CN}(\CY)^{\on{access}}$$
can also be written as \emph{colimits}
$$\Shv_{\CN}(\CY)\simeq \underset{i}{\on{colim}}\, \Shv_\CN(\CU_i),$$
$$\Shv_{\CN}(\CY)^{\on{ren}}\simeq \underset{i}{\on{colim}}\, \Shv_\CN(\CU_i)^{\on{ren}}$$
and 
$$\Shv_{\CN}(\CY)^{\on{access}}\simeq \underset{i}{\on{colim}}\, \Shv_\CN(\CU_i)^{\on{access}},$$
where the transition functors are given by !-extension. 

\sssec{}

The following assertion is proved by considering the $((j_{i_1,i_2})_!,(j_{i_1,i_2})^*)$
and $((j_{i_1,i_2})^*,(j_{i_1,i_2})_*)$ adjunctions:

\begin{lem}  \label{l:N preserved abs}
Let $\CY$ and $\CU_i$ be as in \secref{sss:N preserved abs a}. 

\smallskip

\noindent{\em(a)}
The category $\Shv_{\CN}(\CY)^{\on{ren}}$ is compactly generated. 

\smallskip

\noindent{\em(b)}
The category $\Shv_{\CN}(\CY)$ is compactly generated if and only if each 
$\Shv_{\CN}(\CU_i)$ is compactly generated. 

\smallskip

\noindent{\em(c)}
The pair $(\CY,\CN)$ is renormalization-adapted if and only if each $(\CU_i,\CN|_{\CU_i})$
is renormalization-adapted. 

\smallskip

\noindent{\em(d)}
The pair $(\CY,\CN)$ is constraccessible if and only if each $(\CU_i,\CN|_{\CU_i})$
is constraccessible. 
\end{lem} 

In particular, combining with \propref{p:global quotient N}, we obtain: 

\begin{cor} \label{c:N preserved abs}
Let $\CY$ and $\CU_i$ be as in \secref{sss:N preserved abs a}, and suppose that each $\CU_i$ is a global quotient. 
Then the pair $(\CY,\CN)$ is renormalization-adapted. 
\end{cor} 

\ssec{Product theorems for stacks}

In this subsection we will prove versions of Theorems \ref{t:product thm 1 sch} and \ref{t:product thm 2 sch}
for stacks.

\sssec{}

Let $\CY_1$ and $\CY_2$ be a pair of algebraic stacks. First, by \cite[Proposition A.2.10]{GKRV}, the external tensor product functor
\begin{equation} \label{e:prod stack ff}
\Shv(\CY_1)\otimes \Shv(\CY_2) \to \Shv(\CY_1\times \CY_2)
\end{equation} 
is fully faithful. It follows from \propref{p:boxtimes t-exact} and \lemref{l:comp gen t} that it is t-exact. 

\medskip

An argument parallel to {\it loc. cit.} shows that the functor
\begin{equation} \label{e:prod stack ff ren}
\Shv(\CY_1)^{\on{ren}}\otimes \Shv(\CY_2)^{\on{ren}} \to \Shv(\CY_1\times \CY_2)^{\on{ren}}
\end{equation} 
is also fully faithful. 

\sssec{}

Let $\CN\subset T^*(\CY)$ be a Zariski-closed conical subset. From the fact that \eqref{e:prod stack ff ren}
is fully faithful and using the fact that $\iLisse(X)$ is dualizable, we obtain that for a scheme $X$, the functor
\begin{equation} \label{e:prod stack ff ren N}
\iLisse(X)\otimes \Shv_\CN(\CY)^{\on{ren}} \to \Shv_{\{0\}\times \CN}(X\times \CY)^{\on{ren}}
\end{equation} 
is also fully faithful. 

\medskip

The functor \eqref{e:prod stack ff ren N} induces a functor
\begin{equation} \label{e:prod stack ff access N}
\iLisse(X)\otimes \Shv_\CN(\CY)^{\on{access}} \to \Shv_{\{0\}\times \CN}(X\times \CY)^{\on{access}}
\end{equation} 

\sssec{}

From now on let us assume that $\CN$ is half-dimensional. 

\medskip

We claim:

%
%
%

\begin{thm} \label{t:product thm stack 1b}
Let $X$ be smooth and proper. Then the functor
\begin{equation} \label{e:prod stack access N}
\iLisse(X)\otimes \Shv_\CN(\CY)^{\on{access}} \to \Shv_{\{0\}\times \CN}(X\times \CY)^{\on{access}}
\end{equation} 
is an equivalence.
\end{thm}

\begin{proof}

Since $\iLisse(X)$ is dualizable, by passing to the limit, we reduce to the case when $\CY$ is quasi-compact. 

\medskip

We have a commutative diagram
$$
\CD
\iLisse(X)\otimes \Shv_\CN(\CY)^{\on{access}}  @>>> \Shv_{\{0\}\times \CN}(X\times \CY)^{\on{access}} \\
@VVV @VVV \\
\iLisse(X)\otimes \Shv_\CN(\CY) & & \Shv_{\{0\}\times \CN}(X\times \CY) \\
@VVV @VVV  \\
\iLisse(X)\otimes \Shv(\CY)  @>>>\Shv(X\times \CY) 
\endCD
$$
with vertical arrows being fully faithful (for the left column this again uses the fact that $\iLisse(X)$ is dualizable).
The bottom horizontal arrow is also fully faithful (because $\Shv(\CY)$ is dualizable and the functor 
\eqref{e:prod stack ff} is fully faithful). This implies that the top horizontal arrow,
i.e., our functor \eqref{e:prod stack access N}, is fully faithful.

\medskip

Thus, it suffices to show that the right adjoint of \eqref{e:prod stack access N} is conservative.
Let us describe this right adjoint explicitly (it will be continuous by construction). 

\medskip

Recall that the category $\iLisse(X)$ is canonically self-dual (by Verdier), see Remark \ref{r:Verdier compat}. 
In terms of this self-duality, the right adjoint to \eqref{e:prod stack access N} corresponds to the functor
$$\iLisse(X)\otimes \Shv_{\{0\}\times \CN}(X\times \CY)^{\on{access}} \to \Shv_\CN(\CY)^{\on{access}}$$
given by
$$E,\CF\mapsto (p_Y)_*(p_X^!(E)\sotimes \CF).$$
where $p_X$ and $p_Y$ are the two projections from $X\times Y$ to $X$ and $Y$, respectively.

\medskip

This description of the right adjoint that it commutes with pullbacks for smooth maps $\CY'\to \CY$.
This allows us to replace $\CY$ by an affine scheme covering it. In the latter case, the assertion that
right adjoint to \eqref{e:prod stack access N} is conservative follows from \thmref{t:product thm 1 sch}.

\end{proof} 

\begin{rem}
A similar argument shows that if $X$ be smooth and proper, then the functor \eqref{e:prod stack ff ren N} is also 
an equivalence.
\end{rem}

\sssec{}

Finally, we claim:

\begin{thm} \label{t:product thm stack 2}
Let $X$ be smooth and proper. Assume also that $\qLisse(X)$ duality-adapted (see \secref{ss:Verdier compat}). 
Let $\CN\subset T^*(\CY)$ be half-dimensional. 
Then the resulting functor
$$\qLisse(X)\otimes \Shv_{\CN}(\CY)\to \Shv_{\{0\}\times \CN}(X\times \CY)$$
is an equivalence.
\end{thm}

\begin{proof}

Since $\qLisse(X)$ is assumed dualizable, the assertion reduces to the case of schemes.
In the latter case, it is given by \thmref{t:product thm 2 sch}. 

\end{proof}

\section{Sheaves in the Betti context} \label{s:sheaves Betti}

In this section we will work the ground field $k=\BC$. All our algebro-geometric
objects will be of finite type over $\BC$.

\medskip

The coefficients of our
sheaves will be an arbitrary algebraically closed field of characteristic $0$, denoted $\sfe$. 

\ssec{Betti sheaves on schemes}

\sssec{} \label{sss:all sheaves sch}

Let $Y$ be a scheme of finite type over $\BC$. We denote by $\Shv^{\on{all}}(Y)$ the category
of all sheaves of $\sfe$-vector spaces on the (topological space) underlying $S$.

\medskip

We will consider $\Shv^{\on{all}}(Y)$ as equipped with the usual t-structure, in which it is left-complete
(see \cite[Theorem 7.2.3.6 and Proposition 7.2.1.10]{Lu1}). 

\sssec{} \label{sss:prod Betti}

For a pair of affine schemes, we have a naturally defined functor
\begin{equation} \label{e:prod Betti}
\Shv^{\on{all}}(Y_1)\otimes \Shv^{\on{all}}(Y_2)\overset{\boxtimes}\to \Shv^{\on{all}}(Y_1\times Y_2).
\end{equation} 

This functor is known to be an equivalence (see, e.g., \cite[Theorem 7.3.3.9, Proposition 7.3.1.11]{Lu1} 
and \cite[Proposition 4.8.1.17]{Lu2}).

\sssec{} \label{sss:Shv Betti self-dual}

The category $\Shv^{\on{all}}(Y)$ is known to be dualizable, and in fact, self-dual, with the duality data given by
$$(\Delta_Y)_!(\ul\sfe_Y)\in \Shv^{\on{all}}(Y\times Y)\simeq \Shv^{\on{all}}(Y)\otimes \Shv^{\on{all}}(Y),$$
and
\begin{equation} \label{e:pairing Betti}
\Shv^{\on{all}}(Y)\otimes \Shv^{\on{all}}(Y) \overset{\boxtimes}\to \Shv^{\on{all}}(Y\times Y)
\overset{\Delta_Y^*}\to \Shv^{\on{all}}(Y)\overset{\Gamma_c(Y,-)}\to \Vect_\sfe.
\end{equation}

However, $\Shv^{\on{all}}(Y)$ is \emph{not} compactly generated (unless $Y$ is a disjoint union of points). 

\sssec{}

We have an embedding
\begin{equation} \label{e:constr sheaves to all}
\Shv(Y)^{\on{constr}}\hookrightarrow \Shv^{\on{all}}(Y),
\end{equation}
compatible with t-structures, whose ind-extension is a t-exact functor
\begin{equation} \label{e:sheaves to all}
\Shv(Y)\to  \Shv^{\on{all}}(Y).
\end{equation}

However, the functor \eqref{e:sheaves to all} is \emph{not} fully faithful, because \eqref{e:sheaves to all}
does not preserves compactness. 

\sssec{}

We have a fully faithful (continuous) embedding
\begin{equation} \label{e:loc const in all}
\Shv^{\on{all}}_{\on{loc.const.}}(Y)\hookrightarrow \Shv^{\on{all}}(Y).
\end{equation}

We claim:

\begin{lem}  \label{l:loc const in all}
The functor \eqref{e:loc const in all} admits a left adjoint.
\end{lem} 

\begin{proof}

Since the categories involved are presentable, the assertion of the proposition is
equivalent to the fact that the subcategory \eqref{e:loc const in all} is closed under limits.

\medskip

Thus, let $E_i$ be objects of $\Shv^{\on{all}}_{\on{loc.const.}}(Y)$, where $i$ runs over
some index category $I$, and we wish to show that
$$\CF:=\underset{i}{\on{lim}}\, E_i\in \Shv^{\on{all}}(Y)$$
also belongs to $\Shv^{\on{all}}_{\on{loc.const.}}(Y)$.

\medskip

For every $\BC$-point $y\in Y$, let $U_y\subset Y(\BC)$ be a contractible
\emph{analytic} neighborhood. We wish to show that $\CF|_{U_y}$ is 
constant. 

\medskip

The functor of restriction $\Shv^{\on{all}}(Y)\to \Shv^{\on{all}}(U_y)$ commutes
with limits, it is enough to show that the essential image of
$$\Vect_\sfe\simeq \Shv^{\on{all}}_{\on{loc.const.}}(U_y)\hookrightarrow \Shv^{\on{all}}(U_y)$$
is closed under limits. 

\medskip

However, the latter is evident: an inverse limit of constant sheaves is a constant sheaf.

\end{proof}

\begin{rem}
Let $U$ be a (real) ball of dimension $d$. Then the left adjoint to the embedding 
$$\Vect_\sfe\simeq \Shv^{\on{all}}_{\on{loc.const.}}(U)\hookrightarrow \Shv^{\on{all}}(U)$$
can be described explicitly as follows: it sends
$$\CF\mapsto \on{C}^\cdot_c(U,\CF)[d]\otimes \ell,$$
where $\ell$ is $\sfe$-line, induced from the $\pm 1$-torsor of orientations on $U$.
\end{rem}

\ssec{Singular support for Betti sheaves}

\sssec{} \label{sss:char sing supp Betti}

Let $Y$ be as above, a scheme of finite type over $\BC$. Let $\CN\subset T^*(Y)$ be a conical Lagrangian. 
Following  \cite[Sect. 8]{KS}, we introduce a full subcategory
\begin{equation} \label{e:emb N Betti}
\Shv_\CN^{\on{all}}(Y) \subset \Shv^{\on{all}}(Y).
\end{equation}

\sssec{} \label{sss:Whitney}

One can describe the subcategory \eqref{e:emb N Betti} more explicitly as follows. 

\medskip

thanks to \cite[Corollary 8.3.22]{KS}, we can choose a $\mu$-stratification
$$Y_i \overset{j_i}\hookrightarrow Y(\BC),$$
such that $\CN$ is contained in the union of the conormals to the strata (the latter
will be denote by $\CN_\mu$). 

\medskip

The property of being a $\mu$-stratification implies that for a pair of strata
$$Y_{i_1}\overset{j_{i_1}}\hookrightarrow Y \overset{j_{i_2}}\hookleftarrow Y_{i_2}$$
and $\CF_1\in \Shv^{\on{all}}_{\on{loc.const.}}(Y_{i_1})$, the object
$$(j_{i_2})^*\circ (j_{i_1})_*(\CF_1)\in \Shv^{\on{all}}(Y_{i_2})$$
belongs to $\Shv^{\on{all}}_{\on{loc.const.}}(Y_{i_2})$.

\medskip

This formally implies that for $\CF\in \Shv^{\on{all}}(Y)$ the following conditions are equivalent: 

\bigskip

$$(*)\hskip 1cm j_i^*(\CF)\in \Shv^{\on{all}}_{\on{loc.const.}}(Y_i),\, \forall\, i \,\,\,\Leftrightarrow\, \,\, 
j_i^!(\CF)\in \Shv^{\on{all}}_{\on{loc.const.}}(Y_i),\, \forall\, i;$$

\sssec{} \label{sss:Whitney bis}

Let 
\begin{equation} \label{e:emb mu}
\Shv^{\on{all}}_\mu(Y) \subset \Shv^{\on{all}}(Y)
\end{equation}
be the subcategory of sheaves satisfying (*). Then
\begin{equation} \label{e:emb N strat}
\Shv_\CN^{\on{all}}(Y)  \subset \Shv^{\on{all}}_\mu(Y).
\end{equation}

We will now explain, following \cite[Sect. A.1.3]{NY2}, 
how to single out the subcategory \eqref{e:emb N strat} by an explicit procedure. 

\sssec{} \label{sss:real van}

Choose a \emph{smooth} point $(\xi,y)\in \CN_\mu$. Let $B$ be a sufficiently small open
ball around $y$ in $Y(\BC)$, and let $g$ be a $C^\infty$ real-valued function on $B$ such 
that:

\begin{itemize}

\item $g(y)=0$;

\item $dg_y=\xi$;

\item $\on{Graph}(dg)\cap \CN_\mu=\{(\xi,y)\}$ and the intersection is transversal.

\end{itemize} 

\medskip

Consider the corresponding loci 
$$B_{\leq 0} \overset{q^-}\hookrightarrow B \overset{q^+}\hookleftarrow B_{\geq 0}.$$

Define the functors
$$\Phi^{+,!}_\xi,\Phi^{-,!}_\xi,\Phi^{+,*}_\xi,\Phi^{-,*}_\xi:\Shv(Y)\to \Vect_\sfe$$
by
$$\Phi^{+,!}_\xi(\CF):=\on{C}^\cdot_c(B_{\geq 0},(q^+)^*(\CF|_B)),\quad\Phi^{-,!}_\xi(\CF):=\on{C}^\cdot_c(B_{\leq 0},(q^-)^*(\CF|_B)),$$
$$\Phi^{+,*}_\xi(\CF)=\on{C}^\cdot(B_{\geq 0},(q^+)^!(\CF|_B)),\quad  \Phi^{-,*}_\xi(\CF):=\on{C}^\cdot(B_{\leq 0},(q^-)^!(\CF|_B)).$$

\medskip

For $\CF\in \Shv_\mu(Y)$, the above functors are independent of the choice of $g$
(i.e., they only depend on $\xi$) and satisfy
$$\Phi^{+,!}_\xi(\CF) \simeq \Phi^{-,*}_\xi(\CF) \text{ and }
\Phi^{-,!}_\xi(\CF) \simeq \Phi^{+,*}_\xi(\CF).$$

Note also that we have, tautologically,
$$\Phi^{+,!}_\xi(\CF)\simeq \Phi^{-,!}_{-\xi}(\CF) \text{ and } \Phi^{+,*}_\xi(\CF)\simeq \Phi^{-,*}_{-\xi}(\CF).$$

\sssec{}  \label{sss:N via van}

Now, an object $\CF\in \Shv_\mu(Y)$ belongs to $\Shv_\CN(Y)$ if and only if
$$\Phi^{+,*}_\xi(\CF)=0$$
for all smooth $(\xi,y)\in \CN_\mu-\CN$. 

\medskip

For future reference we note that since $\CN$ is a complex subvariety of $T^*(\CY)$
$$(\xi,y)\in \CN \,\, \Leftrightarrow\,\, (-\xi,y)\in \CN.$$

Hence, the condition that $\CF\in \Shv_\mu(Y)$ belongs to $\Shv_\CN(Y)$ is equivalent to all
$$\Phi^{+,!}_\xi(\CF),\,\, \Phi^{-,!}_\xi(\CF),\,\,\Phi^{+,*}_\xi(\CF),\,\, \Phi^{-,*}_\xi(\CF)$$
being $0$. 
 
\sssec{}

We now claim:

\begin{prop} \label{p:N all limits and colimits}
The subcategory \eqref{e:emb N Betti} is closed under limits and colimits.
\end{prop}

\begin{proof}

First, we claim that the subcategory $\Shv^{\on{all}}_\mu(Y)\subset \Shv(Y)$ 
is closed under limits and colimits. Indeed, this follows from \lemref{l:loc const in all},
combined with the fact that the functors $j_i^!$ commute with limits and the functors
$j_i^*$ commute with colimits.

\medskip

Now the assertion for $\Shv_\CN^{\on{all}}(Y)$ follows using its characterization in 
\secref{sss:N via van}, since the functors $\Phi^{\pm,*}_\xi$ commutes with limits
and the functors $\Phi^{\pm,!}_\xi$ commute with colimits (on all of $\Shv(Y)$). 

\end{proof} 

\sssec{}

Let $\iota^{\on{all}}$ denote the tautological embedding corresponding to \eqref{e:emb N Betti}.

\medskip

From \propref{p:N all limits and colimits} we see that the category $\Shv_\CN^{\on{all}}(Y)$
is cocomplete (and the functor $\iota^{\on{all}}$ is continuous). In addition, we obtain; 

\begin{cor} \label{c:left adjoint all sch}
The functor $\iota^{\on{all}}$ admits a left adjoint.
\end{cor} 

\ssec{Compact generators of the category $\Shv_\CN^{\on{all}}(Y)$}

\sssec{}

We now claim:

\begin{prop} \label{p:comp gen Betti}
The category $\Shv_\CN^{\on{all}}(Y)$ is compactly generated. Moreover, it admits a \emph{finite}
set of compact generators. 
\end{prop} 

The rest of this subsection is devoted to the proof of \propref{p:comp gen Betti}. In fact, we
will exhibit a particular finite set of compact generators. 

\sssec{}

For a point $y\in \CY$, let $\delta_y\in \Shv^{\on{all}}(\CY)$ be the !-$\delta$ function, i.e.,
$\delta_y:=(\bi_y)_!(\sfe)$, where
$$\bi_y:\on{pt}\to \CY$$
is the morphism corresponding to $y$.

\sssec{} \label{sss:mu strat}

Choose a $\mu$-stratification as in \secref{sss:Whitney}. It follows formally from its properties that the 
(\emph{discontinuous}) functors
$$j_i^!:\Shv^{\on{all}}(Y)\to \Shv^{\on{all}}(Y_i)$$
are \emph{continuous} when restricted to $\Shv_\CN^{\on{all}}(Y)$. 

\medskip

Choose a point $y_i$ on each connected component of each stratum. We will prove:

\begin{prop} \label{p:comp gen Betti bis}
The objects $(\iota^{\on{all}})^L(\delta_{y_i})$ compactly generate $\Shv_\CN^{\on{all}}(Y)$.
\end{prop}


\begin{proof}

We have to show that the functors
$$\CF\mapsto \CHom_{\Shv_\CN^{\on{all}}(Y)}((\iota^{\on{all}})^L(\delta_{y_i}),\CF), \quad \Shv_\CN^{\on{all}}(Y)\to \Vect_\sfe$$
are continuous, and of they all vanish on a given $\CF$, then $\CF=0$.

\medskip

By adjunction, we rewrite the above functors as
$$ \CF\mapsto \bi_{y_i}^!\circ \iota^{\on{all}}(\CF),$$
and further as 
$$\CF \mapsto \wt\bi_{y_i}^!\circ j_i^!\circ \iota^{\on{all}}(\CF),$$
where 
$$\on{pt}\overset{\wt\bi_{y_i}}\to Y_i.$$

Now, the assertion follows from the combination of the following facts:

\begin{itemize}

\item The functors $j_i^!:\Shv_\CN^{\on{all}}(Y)\to \Shv_{\on{loc.const.}}(Y_i)$ are continuous;

\item If $j_i^!(\CF)=0$ for some $\CF\in \Shv^{\on{all}}(Y)$ for all $i$, then $\CF=0$;

\item For every $i$, the functor $\wt\bi_{y_i}^!$ is continuous and conservative on $\Shv_{\on{loc.const.}}(Y_i)$.

\end{itemize} 

\end{proof}

\ssec{Constructible sheaves vs all sheaves}

\sssec{}

As was mentioned above, the functor \eqref{e:sheaves to all} does not preserve compactness. 
However, we claim:

\begin{prop} \label{p:constr Betti sch}
Let $\CN\subset T^*(Y)$ be a conical Lagrangian. Then the objects of $\Shv_\CN(Y)^{\on{constr}}$
are compact as objects of $\Shv_\CN^{\on{all}}(Y)$.
\end{prop}

As a corollary, we obtain:

\begin{cor}
The functor
$$\Shv_\CN^{\on{access}}(Y)\to \Shv_\CN^{\on{all}}(Y)$$
is fully faithful.
\end{cor}

\begin{rem}
We do not know whether the functor 
$$\Shv_\CN(Y)\to \Shv_\CN^{\on{all}}(Y)$$
preserves compactness or is fully faithful. (Note that in our usual counterexample
of \secref{sss:analyze P1}, this functor is actually an equivalence.) 
\end{rem} 

\sssec{Proof of \propref{p:constr Betti sch}}

Let $\CF$ be an object of $\Shv_\CN(Y)^{\on{constr}}$. Choose a stratification as in \secref{sss:mu strat}. 
Then for $\CF'\in \Shv^{\on{all}}_\CN(Y)$, we can express
$$\CHom_{\Shv_\CN^{\on{all}}(Y)}(\CF,\CF')$$
as a \emph{finite} colimit with terms 
\begin{equation} \label{e:Hom on stratum}
\CHom_{\Shv_\CN^{\on{all}}(Y_i)}(j_i^*(\CF),j_i^!(\CF')),
\end{equation}
functorially in $\CF'$. 

\medskip 

Hence, it suffices to show that the expression in \eqref{e:Hom on stratum}, viewed as a functor of $\CF'$,
is continuous. Since the functor $j_i^!$ is continuous on $\Shv_\CN^{\on{all}}(Y)$, and 
$$j_i^*(\CF)\in \Lisse(Y_i), \quad j_i^!(\CF')\in \Shv^{\on{all}}_{\on{loc.const.}}(Y_i),$$
it suffices to show that for a given $E\in \Lisse(Y_i)$, the functor
$$E'\mapsto \CHom_{\Shv^{\on{all}}_{\on{loc.const.}}(Y_i)}(E,E'), \quad E'\in \Shv^{\on{all}}_{\on{loc.const.}}(Y_i)$$
is continuous, which is obvious.

\qed[\propref{p:constr Betti sch}]

\ssec{Duality on the category of Betti sheaves with prescribed singular support} \label{ss:duality N all} 

\sssec{}

The goal of this subsection is to prove the following assertion: 

\begin{thm} \label{t:duality Betti N} \hfill

\smallskip

\noindent{\em(a)}
The pairing
$$\Shv_\CN^{\on{all}}(Y) \otimes \Shv_\CN^{\on{all}}(Y) \overset{\iota^{\on{all}}\otimes \iota^{\on{all}}}\longrightarrow
\Shv^{\on{all}}(Y) \otimes \Shv^{\on{all}}(Y)\overset{\on{C}^\cdot_c(Y,-\overset{*}\otimes -)}\longrightarrow \Vect_\sfe$$
is the counit of a self-duality. 

\smallskip

\noindent{\em(b)} With respect to the above identification $\Shv_\CN^{\on{all}}(Y)^\vee\simeq \Shv_\CN^{\on{all}}(Y)$
and the identification $\Shv^{\on{all}}(Y)^\vee\simeq \Shv^{\on{all}}(Y)$ of \secref{sss:Shv Betti self-dual}, the natural map
$$(\iota^{\on{all}})^\vee\to (\iota^{\on{all}})^L$$
is an isomorphism. 

\end{thm}

The assertion of \thmref{t:duality Betti N} is formally equivalent to the following: 

\begin{thm} \label{t:duality Betti N bis} 
For $\CF\in \Shv_\CN^{\on{all}}(Y)$ and $\CF'\in \Shv^{\on{all}}(Y)$, the unit of the $((\iota^{\on{all}})^L,\iota^{\on{all}})$-adjunction
defines an \emph{isomorphism} 
$$\on{C}^\cdot_c(Y,\iota^{\on{all}}(\CF)\overset{*}\otimes \CF')\to 
\on{C}^\cdot_c(Y,\iota^{\on{all}}(\CF)\overset{*}\otimes (\iota^{\on{all}}\circ (\iota^{\on{all}})^L(\CF'))).$$
\end{thm}

\sssec{}

Let 
$$\BD:(\Shv^{\on{all}}(Y))^{\on{op}}\to \Shv^{\on{all}}(Y),$$
defined by
$$\CHom_{\Shv^{\on{all}}(Y)}(\CF,\BD(\CF'))=\on{C}^\cdot_c(Y,\CF\overset{*}\otimes \CF')^\vee, \quad \CF,\CF'\in \Shv^{\on{all}}(Y).$$

\sssec{}

The assertion of \thmref{t:duality Betti N bis} is in turn equivalent to the following:

\begin{thm} \label{t:duality Betti N duality} 
The (contravariant) functor $\BD$ preserves the subcategory $\Shv_\CN^{\on{all}}(Y)$.
\end{thm}

The rest of this section is devoted to the proof of \thmref{t:duality Betti N duality}.

\sssec{Reduction to a $\mu$-stratification}

Choose a stratification as in \secref{sss:Whitney}. 
We claim that it suffices to show that the functor $\BD$ preserves the corresponding 
subcategory $\Shv_\mu^{\on{all}}(Y)$ (see \secref{sss:Whitney bis}). 

\medskip

Indeed, suppose that we know that for $\CF\in \Shv_\CN^{\on{all}}(Y)\subset \Shv_\mu^{\on{all}}(Y)$,
the object $\BD(\CF)$ belongs to $\Shv_\mu^{\on{all}}(Y)$. 

\medskip

By \secref{sss:N via van}, we have to show that for all smooth $\xi\in \CN_\mu-\CN$,
$$\Phi_\xi^{+,*}(\BD(\CF))=0.$$

However, for any $\CF\in \Shv^{\on{all}}(Y)$, we have
\begin{equation} \label{e:Phi duality}
\Phi_\xi^{+,*}(\BD(\CF)) \simeq (\Phi_\xi^{+,!}(\CF))^\vee.
\end{equation} 

Now, 
$$\Phi_\xi^{+,!}(\CF)\simeq \Phi_\xi^{-,*}(\CF) \simeq \Phi_{-\xi}^{+,*}(\CF)=0.$$

\sssec{Reduction to a stratum}

Let $Y_i\overset{j_i}\hookrightarrow Y$ be the embedding of a stratum. For any
$\CF\in \Shv^{\on{all}}(Y)$ we have
$$j_i^!(\BD_Y(\CF))\simeq \BD_{Y_i}(j_i^*(\CF)).$$

This reduces the assertion of \thmref{t:duality Betti N duality} to the case when $Y(\BC)$ is a manifold
and $$\Shv^{\on{all}}_\CN(Y)=\Shv^{\on{all}}_{\on{loc.const.}}(Y).$$

\sssec{Reduction to a ball}

Let $U\subset Y(\BC)$ be an open ball. For $\CF\in \Shv^{\on{all}}(U)$, we have 
$$\BD_Y(\CF)|_U\simeq \BD_U(\CF|_U).$$

Hence, it remains to show that the functor $\BD_U$ preserves the full subcategory

$$\Vect_\sfe\simeq \Shv^{\on{all}}_{\on{loc.const.}}(U)\subset \Shv^{\on{all}}(U).$$

However, for $V\in \Vect_\sfe$, 
$$\BD_U(\ul\sfe_U\otimes V) \simeq \ul\sfe_U\otimes V^\vee[2\dim(U)].$$

\qed[\thmref{t:duality Betti N duality}]

\ssec{A smoothness property of the category of Betti sheaves with restricted singular support}

The material in this subsection is not needed in the rest of the paper.

\sssec{}

Recall that a dualizable DG category $\bC$ is said to be \emph{smooth} if the unit
object
$$\on{u}_\bC\in \bC\otimes \bC^\vee$$
is compact. 

\medskip

We are going to prove:

\begin{prop} \label{p:sing supp smooth}
Let $Y$ be a scheme and let $\CN\subset T^*(Y)$ be a Zariski-closed conical Lagrangian 
subset. Then the category $\Shv_\CN^{\on{all}}(Y)$ is smooth. 
\end{prop} 

The rest of this subsection is devoted to the proof of this proposition\footnote{We have learned this assertion from D.~Nadler.}. 

\sssec{}

Let $\on{u}_{\Shv_\CN^{\on{all}}(Y)}$ denote the unit of the self-duality specified by
\thmref{t:duality Betti N}. We have have to check that the functor
\begin{equation} \label{e:Hom from unit}
(E\in \Shv_\CN^{\on{all}}(Y)\otimes \Shv_\CN^{\on{all}}(Y))\mapsto 
\CHom_{\Shv_\CN^{\on{all}}(Y)\otimes \Shv_\CN^{\on{all}}(Y)}(\on{u}_{\Shv_\CN^{\on{all}}(Y)},E)
\end{equation} 
is continuous. 

\medskip

By \thmref{t:duality Betti N}(b), object $\on{u}_{\Shv_\CN^{\on{all}}(Y)}$ can be explicitly described
as
$$((\iota^{\on{all}})^L\otimes (\iota^{\on{all}})^L)((\Delta_Y)_!(\ul\sfe_Y)).$$

Hence, the expression in \eqref{e:Hom from unit} can be rewritten as
\begin{equation} \label{e:Hom from unit again}
\CHom_{\Shv^{\on{all}}(Y)\otimes \Shv^{\on{all}}(Y)}((\Delta_Y)_!(\ul\sfe_Y),(\iota^{\on{all}}\otimes \iota^{\on{all}})(E)).
\end{equation} 

Denote by $\wt{E}$ the image of $E$ along
$$\Shv_\CN^{\on{all}}(Y)\otimes \Shv_\CN^{\on{all}}(Y) \overset{\boxtimes}\hookrightarrow \Shv^{\on{all}}_{\CN\times \CN}(Y\times Y).$$

Thus, we can rewrite \eqref{e:Hom from unit again} as
$$\CHom_{\Shv^{\on{all}}(Y\times Y)}((\Delta_Y)_!(\ul\sfe_Y),\iota_{Y\times Y}^{\on{all}}(\wt{E})) \simeq
\on{C}^\cdot(Y,\Delta_Y^!\circ \iota_{Y\times Y}^{\on{all}}(\wt{E})),$$
where $\iota_{Y\times Y}^{\on{all}}$ is the embedding
$$\Shv^{\on{all}}_{\CN\times \CN}(Y\times Y) \hookrightarrow \Shv^{\on{all}}(Y\times Y).$$

Thus, it is sufficient to show that the functor
\begin{equation} \label{e:Hom from unit again again}
\CF\mapsto \on{C}^\cdot(Y,\Delta_Y^!\circ \iota_{Y\times Y}^{\on{all}}(\CF)),
\quad \Shv^{\on{all}}_{\CN\times \CN}(Y\times Y)\to \Vect_\sfe
\end{equation} 
is continuous.  (Note, that the functor $\on{C}^\cdot(Y,\Delta_Y^!(-))$ itself is \emph{discontinuous}.) 

\sssec{}

Using a stratification as in \secref{sss:mu strat}, we can write the functor \eqref{e:Hom from unit again again} as a \emph{finite} limit
of functors of the form
$$\CF\mapsto \on{C}^\cdot(Y_i,\Delta_{Y_i}^! \circ (j_i\times j_i)^!(\CF)),$$
so it suffices to show that each of the latter functors is continuous on $\Shv^{\on{all}}_{\CN\times \CN}(Y\times Y)$. 

\medskip

As in \secref{sss:mu strat}, each of the functors $(j_i\times j_i)^!$ is continuous on $\Shv^{\on{all}}_{\CN\times \CN}(Y\times Y)$
and maps to
$$\Shv_{\on{loc.const}}(Y_i\times Y_i).$$

Next, the functor $\Delta_{Y_i}^!$, when restricted to $\Shv_{\on{loc.const}}(Y_i\times Y_i)$, is isomorphic to
$\Delta_{Y_i}^*$ up to a shift, and hence is also continuous, and maps to $\Shv_{\on{loc.const}}(Y_i)$. 

\medskip

Finally, the functor $\on{C}^\cdot(Y_i,-)$ is continuous, when restricted to $\Shv_{\on{loc.const}}(Y_i)$. 

\qed[\propref{p:sing supp smooth}]

\ssec{Betti sheaves on stacks}

For the purposes of this paper, we do \emph{not} need the category $\Shv^{\on{all}}(-)$ for 
arbitrary prestacks, but only for algebraic stacks. 

\sssec{}

Thus, let $\CY$ be an algebraic stack. We define 
\begin{equation} \label{e:all sheaves stack lim}
\Shv^{\on{all}}(\CY):=\underset{S}{\on{lim}}\, \Shv^{\on{all}}(S),
\end{equation} 
where the limit is taken over the category of affine schemes almost of finite type,
equipped with a smooth map to $\CY$, and where we allow only smooth maps
$f:S_1\to S_2$ over $\CY$. The transition functors are given by
$$\Shv^{\on{all}}(S_2) \overset{f^!}\to \Shv^{\on{all}}(S_1).$$
These functors are continuous because the morphisms $f$ were assumed smooth. 

\medskip

This limit can be rewritten as a colimit
\begin{equation} \label{e:all sheaves stack}
\Shv^{\on{all}}(\CY):=\underset{S}{\on{colim}}\, \Shv^{\on{all}}(S),
\end{equation} 
where the transition functors are given by !-puhsforward. 

\sssec{}

We define the functor
$$\on{C}^\cdot_c(\CY,-):\Shv^{\on{all}}(\CY)\to \Vect_\sfe$$
to correspond to the (compatible) system of functors
$$\on{C}^\cdot_c(S,-):\Shv^{\on{all}}(S)\to \Vect_\sfe$$
in the presentation \eqref{e:all sheaves stack}. 

\sssec{} \label{sss:Lagr stk}

Let $\CN\subset T^*(\CY)$ be a conical Lagrangian subset, which by definition means that for
any affine scheme equipped with a smooth map $S\to \CY$, the subset
$$\CN_S\subset T^*(\CY),$$
defined as in \secref{sss:sing supp stack qc}, is Lagrangian.

\medskip

Define the full subcategory
\begin{equation} \label{e:emb N Betti stack}
\Shv_\CN^{\on{all}}(\CY)\overset{\iota^{\on{all}}}\hookrightarrow \Shv^{\on{all}}(\CY)
\end{equation} 
as 
$$\underset{S}{\on{lim}}\, \Shv^{\on{all}}_{\CN_S}(S)\subset \underset{S}{\on{lim}}\, \Shv^{\on{all}}(S).$$

\sssec{}

Since the functors of pullback along smooth morphisms commute with limits and colimits,
from \propref{p:N all limits and colimits}:

\begin{cor} \label{c:N all limits and colimits stack}
The subcategory \eqref{e:emb N Betti stack} is closed under limits and colimits.
\end{cor}

\medskip

In particular, we obtain that the category $\Shv_\CN^{\on{all}}(\CY)$ is presentable, 
and the functor $\iota^{\on{all}}$ commutes with limits and colimits. In particular, we obtain:

\medskip

\begin{cor} \label{c:left adjoint all stack}
The functor $\iota^{\on{all}}$ of \eqref{e:emb N Betti stack} admits a left adjoint.
\end{cor} 

\sssec{}

We now claim:

\begin{prop} \label{p:comp gen Betti stack}
The category $\Shv_\CN^{\on{all}}(\CY)$ is compactly generated. If $\CY$ is quasi-compact, then 
is admits a finite set of compact generators.
\end{prop}

As in the case of \propref{p:comp gen Betti}, we will exhibit a particular set of compact generators
for $\Shv_\CN^{\on{all}}(\CY)$.

\begin{proof}[Proof of \propref{p:comp gen Betti stack}]

Let $f_\alpha:S_\alpha\to \CY$ be a covering collection of smooth maps, where $S_\alpha$ are affine schemes.
(Note that when $\CY$ is quasi-compact, we can take a single affine scheme $S$.)

\medskip

For each index $\alpha$ consider the category
$$\Shv^{\on{all}}_{\CN_{S_\alpha}}(S_\alpha),$$
and let $s_{\alpha,i}\in S_\alpha$ be a finite collection of points such that the objects
$$(\iota_{S_\alpha}^{\on{all}})^L(\delta_{s_{\alpha,i},S_\alpha})\in \Shv^{\on{all}}_{\CN_{S_\alpha}}(S_\alpha)$$
(compactly) generate $\Shv^{\on{all}}_{\CN_{S_\alpha}}(S_\alpha)$. 

\medskip

Let $y_{\alpha,i}$ denote the image of $s_{\alpha,i}$ in $\CY$. We claim that the objects 
$$(\iota^{\on{all}}_\CY)^L(\delta_{y_{\alpha,i},\CY})\in \Shv_\CN^{\on{all}}(\CY)$$
are compact and generate $\Shv_\CN^{\on{all}}(\CY)$.

\medskip

Indeed, for every $\alpha$ and $i$ and $\CF\in \Shv_\CN^{\on{all}}(\CY)$, we have
$$\CHom_{\Shv_\CN^{\on{all}}(\CY)}((\iota^{\on{all}}_\CY)^L(\delta_{y_{\alpha,i},\CY}),\CF)
\simeq \bi_{y_{\alpha,i}}^!\circ \iota_\CY^{\on{all}}(\CF) \simeq \bi_{s_{\alpha,i}}^!\circ f_\alpha^! \circ \iota_\CY^{\on{all}}(\CF)\simeq $$
$$\simeq \bi_{s_{\alpha,i}}^!\circ \iota_S^{\on{all}}\circ f_\alpha^!(\CF)
\simeq
\CHom_{\Shv^{\on{all}}_{\CN_{S_\alpha}}(S_\alpha)}((\iota_{S_\alpha}^{\on{all}})^L(\delta_{s_{\alpha,i},S_\alpha}),
f_\alpha^!(\CF)),$$
where in the second line, $f_\alpha^!$ is regarded as a functor 
\begin{equation} \label{e:f alpha}
\Shv_\CN^{\on{all}}(\CY)\to \Shv^{\on{all}}_{\CN_{S_\alpha}}(S_\alpha).
\end{equation} 

This implies the compactness statement, since the functors \eqref{e:f alpha} are continuous. This also
implies the generation statements as
$$\CHom_{\Shv_\CN^{\on{all}}(\CY)}((\iota^{\on{all}}_\CY)^L(\delta_{y_{\alpha,i},\CY}),\CF)=0,\,\,
\forall\, \alpha \text{ and } i$$
implies
$$f_\alpha^!(\CF)=0,\,\, \forall\, \alpha \,\, \Rightarrow \CF=0.$$

\end{proof}

\sssec{}

Finally, we have the following extension of \propref{p:constr Betti sch}:

\begin{prop} \label{p:constr Betti stack}
Suppose that $\CY$ can be covered by a filtered family of quasi-compact substacks,
each of which is a global quotient (i.e., of the form $Y/G$, where $Y$ is a quasi-compact
scheme). Then the functor 
$$\Shv_\CN(\CY)\to \Shv_\CN^{\on{all}}(\CY)$$
sends objects that are compact in $\Shv(\CY)$ to objects that are compact in $\Shv_\CN^{\on{all}}(\CY)$.
\end{prop} 

\begin{rem}
We expect that the assertion of \propref{p:constr Betti stack} holds for any algebraic
stack $\CY$, without the local global quotient condition. 
\end{rem}

\ssec{Proof of \propref{p:constr Betti stack}}

\sssec{}

Since compact objects in $\Shv(\CY)$ are !-extensions from quasi-compact substacks, we can 
assume that $\CY$ is quasi-compact.

\medskip

By assumption, we can assume that $\CY$ is of the form $Y/G$, where $Y$ is a quasi-compact 
scheme and $G$ is an algebraic group.
%
%

\sssec{}

For a pair of objects $\CF_1,\CF_2\in \Shv^{\on{all}}(Y/G)$, let $\CF'_1,\CF'_2$ denote their
*-pullbacks to $Y$. On the one hand, we can consider the object
$$\CHom_{\Shv^{\on{all}}(Y)}(\CF'_1,\CF'_2)\in \Vect_\sfe.$$

On the other hand, we can consider the object
$$\ul\CHom_{\Shv^{\on{all}}(Y)}(\CF'_1,\CF'_2)\in \Shv^{\on{all}}(\on{pt}/G)\simeq \Shv(\on{pt}/G)$$
defined by the requirement that
$$\CHom_{\Shv^{\on{all}}(\on{pt}/G)}(V,\ul\CHom_{\Shv^{\on{all}}(Y)}(\CF'_1,\CF'_2))\simeq
\CHom_{\Shv^{\on{all}}(Y/G)}(q^*(V)\overset{*}\otimes \CF_1,\CF_2), \quad V\in \Shv(\on{pt}/G),$$
where $q:Y/G\to \on{pt}/G$. 

\medskip

It is easy to see that 
$$\pi_{\on{pt}}^*(\ul\CHom_{\Shv^{\on{all}}(Y)}(\CF'_1,\CF'_2))\simeq 
\CHom_{\Shv^{\on{all}}(Y)}(\CF'_1,\CF'_2),$$
where $\pi_{\on{pt}}:\on{pt}\to \on{pt}/G$.

\medskip

Note that
\begin{equation} \label{e:abs vs inner Hom}
\CHom_{\Shv^{\on{all}}(Y/G)}(\CF_1,\CF_2)\simeq \on{C}^\cdot(\on{pt}/G,\ul\CHom_{\Shv^{\on{all}}(Y)}(\CF'_1,\CF'_2)).
\end{equation} 

\sssec{}

We wish to show that for $\CF_1\in \Shv_\CN(Y/G)\cap \Shv(Y/G)^c$, the functor
\begin{equation} \label{e:abs Hom}
\CF_2\mapsto \CHom_{\Shv^{\on{all}}(Y/G)}(\CF_1,\CF_2), \quad \Shv^{\on{all}}(Y/G)\to \Vect_\sfe
\end{equation} 
commutes with colimits when restricted to $\Shv_\CN^{\on{all}}(Y/G)$.

\medskip

First, we claim that for $\CF_1\in \Shv_\CN(Y/G)^{\on{constr}}$, the functor
\begin{equation} \label{e:inner Hom}
\CF_2\mapsto \ul\CHom_{\Shv^{\on{all}}(Y)}(\CF'_1,\CF'_2), \quad \Shv^{\on{all}}(Y/G)\to \Shv(\on{pt}/G)
\end{equation} 
commutes with colimits when restricted to $\Shv_\CN^{\on{all}}(Y/G)$.

\medskip

Indeed, since the functor $\pi^*_{\on{pt}}$ commutes with colimits and is conservative, it suffices to
show that the composition of \eqref{e:inner Hom} with $\pi^*_{\on{pt}}$ commutes with colimits
when restricted to $\Shv_\CN^{\on{all}}(Y/G)$.

\medskip

However, the latter composition is the functor
$$\CF_2\mapsto \CHom_{\Shv^{\on{all}}(Y)}(\CF'_1,\CF'_2),$$
and it commutes with colimits when restricted to $\Shv_\CN^{\on{all}}(Y/G)$ by 
\propref{p:constr Betti sch}.

\sssec{}

Recall now the natural transformation
$$\on{C}_\blacktriangle^\cdot(\on{pt}/G,-)\to \on{C}^\cdot(\on{pt}/G,-),$$
see \eqref{e:from triangle}. Given \eqref{e:abs vs inner Hom}, it suffices to show that
for $\CF_1\in \Shv(Y/G)^c$ and any $\CF_2$, the map
$$\on{C}_\blacktriangle^\cdot(\on{pt}/G,\ul\CHom_{\Shv^{\on{all}}(Y)}(\CF'_1,\CF'_2))\to
\on{C}^\cdot(\on{pt}/G,\ul\CHom_{\Shv^{\on{all}}(Y)}(\CF'_1,\CF'_2))$$
is an isomorphism. 

\sssec{}

By \secref{sss:gen glob quot}, we can assume that $\CF_1$ is of the form
$(\pi_Y)_!(\wt\CF)$ for $\wt\CF\in \Shv^{\on{all}}(Y)$. 
%
%
%
%
%
%

\medskip

In this case, it is easy to see that
$$\ul\CHom_{\Shv^{\on{all}}(Y)}(\CF'_1,\CF'_2)\simeq 
(\pi_{\on{pt}})_*(\CHom_{\Shv^{\on{all}}(Y)}(\wt\CF,(\pi_Y)^!(\CF_2))).$$

\sssec{}

Now, we claim that for any $V\in \Vect_\sfe$, the map
$$\on{C}_\blacktriangle^\cdot(\on{pt}/G,(\pi_{\on{pt}})_*(V))\to \on{C}^\cdot(\on{pt}/G,(\pi_{\on{pt}})_*(V))$$
is an isomorphism. 

\medskip

Indeed, it suffices to show that the right-hand side commutes with colimits in $V$. However, this
follows from the fact that
$$\on{C}^\cdot(\on{pt}/G,(\pi_{\on{pt}})_*(V))=
\CHom_{\Shv(\on{pt}/G)}(\ul\sfe_{\on{pt}/G},(\pi_{\on{pt}})_*(V))\simeq 
\CHom_{\Shv(\on{pt})}(\pi_{\on{pt}}^*(\ul\sfe_{\on{pt}/G}),V)=V.$$

\qed[\propref{p:constr Betti stack}]

\ssec{Duality for Betti sheaves on stacks}

\sssec{} \label{sss:self-dual all stk} 

Let $\CY_1,\CY_2$ be a pair of algebraic stacks. It follows from \eqref{e:all sheaves stack}
and the equivalence \eqref{e:prod Betti} for schemes that the functor
$$\Shv^{\on{all}}(\CY_1)\otimes \Shv^{\on{all}}(\CY_2) \overset{\boxtimes}\to 
\Shv^{\on{all}}(\CY_1 \times \CY_2)$$
is an equivalence. 

\medskip

From here it follows formally that for an algebraic stack $\CY$, the data of 
$$(\Delta_\CY)_!(\ul\sfe_\CY)\in \Shv^{\on{all}}(\CY\times \CY)\simeq \Shv^{\on{all}}(\CY)\otimes \Shv^{\on{all}}(\CY),$$
and
\begin{equation} \label{e:pairing Betti stk}
\Shv^{\on{all}}(\CY)\otimes \Shv^{\on{all}}(\CY) \overset{\boxtimes}\to \Shv^{\on{all}}(\CY\times \CY)
\overset{\Delta_\CY^*}\to \Shv^{\on{all}}(\CY)\overset{\Gamma_c(\CY,-)}\to \Vect_\sfe.
\end{equation}
define a self-duality on $\Shv^{\on{all}}(\CY)$.

\sssec{} 

Let $\CN\subset T^*(\CY)$ be as in \secref{sss:Lagr stk}. We are going to deduce from 
\thmref{t:duality Betti N duality} the following assertion:

\begin{cor} \label{c:self-dual N all stk}  \hfill

\smallskip

\noindent{\em(a)}
The pairing
$$\Shv_\CN^{\on{all}}(\CY) \otimes \Shv_\CN^{\on{all}}(\CY) \overset{\iota^{\on{all}}\otimes \iota^{\on{all}}}\longrightarrow
\Shv^{\on{all}}(\CY) \otimes \Shv^{\on{all}}(\CY)\overset{\on{C}^\cdot_c(\CY,-\overset{*}\otimes -)}\longrightarrow \Vect_\sfe$$
is the counit of a self-duality. 

\smallskip

\noindent{\em(b)} With respect to the above identification $\Shv_\CN^{\on{all}}(\CY)^\vee\simeq \Shv_\CN^{\on{all}}(\CY)$
and the identification $\Shv^{\on{all}}(\CY)^\vee\simeq \Shv^{\on{all}}(\CY)$ of \secref{sss:self-dual all stk}, the resulting map
$$(\iota^{\on{all}})^\vee\to (\iota^{\on{all}})^L$$
is an isomorphism. 

\end{cor}

\begin{proof} 

As in \secref{ss:duality N all}, the assertion of the corollary is formally equivalent to the fact that the functor
$$\BD:(\Shv^{\on{all}}(\CY))^{\on{op}}\to \Shv^{\on{all}}(\CY),$$
defined by
$$\CHom_{\Shv^{\on{all}}(\CY)}(\CF,\BD(\CF'))=\on{C}^\cdot_c(\CY,\CF\overset{*}\otimes \CF')^\vee, \quad \CF,\CF'\in \Shv^{\on{all}}(\CY)$$
preserves the subcategory $\Shv_\CN^{\on{all}}(\CY)$.

\medskip

Let $S$ be a scheme equipped with a smooth map $f:S\to \CF$. By the definition of $\Shv_\CN^{\on{all}}(\CY)$, we need to show
that for $\CF\in \Shv_\CN^{\on{all}}(\CY)$, we have
$$f^!(\BD_\CY(\CF))\in \Shv_{\CN_S}^{\on{all}}(S).$$

However, we have 
$$f^!(\BD_\CY(\CF)) \simeq \BD_S(f^*(\CF)),$$
for any $\CF\in \Shv(\CY)$. 

\medskip

Since $f$ is smooth, $f^*(\CF)\in \Shv_{\CN_S}^{\on{all}}(S)$. Now the fact that $\BD_S(f^*(\CF))$ belongs to
$\Shv_{\CN_S}^{\on{all}}(S)$ follows from \thmref{t:duality Betti N duality}.

\end{proof}

\section{Proof of \thmref{t:sing supp dir im}} \label{s:proof sing supp} 

In this section, we will prove \thmref{t:sing supp dir im}. We will work in a constructible \'etale, Betti
or regular holonomic sheaf-theoretic context. 

\medskip

Note, however, that in order to treat the Betti case, it is sufficient to treat the regular holonomic
one, by Riemann-Hilbert. So, from now on we will assume that our sheaf-theoretic context
is either \'etale or regular holonomic (this will be needed for a change of fields manipulation in 
\secref{ss:construct triple}.)

%
%
%
%
%
%
%
%

\ssec{Method of proof}

\sssec{Reduction to the case of a proper morphism}

First, with no restriction of generality, we can assume that $\CY_2$ is a smooth scheme
(and hence $\CY_1$ is a separated scheme as well, since $f$ was assumed schematic 
and separated). 

\medskip

The assumption on $\CF_1$ is local on $\CY_1$ around the point $y_1$. Hence, we 
can assume that $f$ is proper: indeed choose a relative compactification of $f$
$$\CY_1\overset{j}\hookrightarrow \ol\CY_1 \overset{\ol{f}}\to \CY_2,$$
and replace the initial $\CF_1$ with $j_*(\CF_1)$. 



\sssec{}

Recall that for a scheme $\CY$ and $\CF\in \Shv(\CF)$, we define 
$$\on{SingSupp}(\CF) \subset T^*(\CY)$$ 
as the union of irreducible closed subsets $\CN$ that appear as irreducible 
components of constructible sub-objects $\CF'$ of $H^m(\CF)$ for all $m$. 

\medskip

By \cite{Be2}, all such $\CN$ are half-dimensional\footnote{Note this makes sense whether
or not $\CY$ is smooth, see \secref{sss:half-dim sing}.}. (But in the
\'etale setting, when $\on{char}(k)\neq 0$, they are not necessarily Lagrangian.)

\medskip

Thus, an irreducible subvariety $Z$ of $T^*(\CY)$
belongs to $\on{SingSupp}(\CF)$ if it is contained in one of the irreducible components 
$\CN$ as above. 

\sssec{} \label{sss:choose N1}

By assumption (i) in \thmref{t:sing supp dir im}, we can find an irreducible 
half-dimensional subvariety $\CN_1\subset \on{SingSupp}(\CF_1)$ with $(df^*(\xi_2),y_1)\in \CN_1$. 
We will fix it from now on. 

\medskip

Consider the subscheme 
$$(df^*)^{-1}(\CN_1)\subset T^*(\CY_2)\underset{\CY_2}\times \CY_1.$$

Let $\wt\CN_1$ be its irreducible component that contains the point $(\xi_2,y_1)$.
By assumption, 
$$\dim(\wt\CN_1)=\dim(\CY_2),$$
and the map
$$\wt\CN_1\hookrightarrow T^*(\CY_2)\underset{\CY_2}\times \CY_1\to T^*(\CY_2)$$
is quasi-finite at $(\xi_2,y_1)$. 

\medskip

Let $\CN_2\subset T^*(\CY_2)$ be the closure of the image of the above map. This is
an irreducible subvariety of $T^*(\CY_2)$ of dimension equal to $\dim(\CY_2)$, and
$$(\xi_2,y_2)\in \CN_2.$$

\sssec{}

Let
$$\Shv_{\wt\CN_1\on{-q.f.}}(\CY_1)\subset \Shv(\CY_1)$$ 
be the full subcategory of objects $\CF'_1$ such that for every irreducible 
$\CN'_1\subset \on{SingSupp}(\CF'_1)$ such that 
$$\wt\CN_1 \subset (df^*)^{-1}(\CN'_1),$$
the dimension of $(df^*)^{-1}(\CN'_1)$
at the generic point of $\wt\CN_1$ equals $\dim(\CY_2)$. This is equivalent to requiring that
the map
$$(df^*)^{-1}(\CN'_1)\hookrightarrow T^*(\CY_2)\underset{\CY_2}\times \CY_1\to T^*(\CY_2)$$
be quasi-finite at the generic point of $\wt\CN_1$.  This is also equivalent to 
$\wt\CN_1$ being actually a whole irreducible component of $(df^*)^{-1}(\CN'_1)$ 
as a subset\footnote{I.e., we discard embedded components.}. 

\medskip

By assumption, our $\CF_1$ belongs to $\Shv_{\wt\CN_1\on{-q.f.}}(\CY_1)$.

\sssec{} \label{sss:Phi req}

We will construct a DG category $\bC$ and functors
$$\Phi_1:\Shv_{\wt\CN_1\on{-q.f.}}(\CY_1)\to \bC$$
and 
$$\Phi_2:\Shv(\CY_2)\to \bC$$
with the following properties:

\medskip

\begin{enumerate}

\item If $\CF_2\in \Shv(\CY_2)$ and $\CN_2 \not\subset \on{SingSupp}(\CF_2)$, then $\Phi_2(\CF_2)=0$;

\smallskip

\item If $\CF'_1\in \Shv_{\wt\CN_1\on{-q.f.}}(\CY_1)$ and $\CN_1\subset \on{SingSupp}(\CF'_1)$, then $\Phi_1(\CF'_1)\neq 0$;

\smallskip

\item $\Phi_1$ is canonically isomorphic to a direct summand of the composition
$$\Shv_{\wt\CN_1\on{-q.f.}}(\CY_1)\hookrightarrow \Shv(\CY_1) \overset{f_*}\longrightarrow \Shv(\CY_2)  \overset{\Phi_2}\longrightarrow \bC.$$

\end{enumerate} 

\smallskip

Clearly, the existence of a triple $(\bC,\Phi_1,\Phi_2)$ as above implies \thmref{t:sing supp dir im}.

\ssec{Isolated points and vanishing cycles, constructible case}

\sssec{}

Let $\CY$ be a scheme\footnote{The discussion in this subsection and the next is local, so to
check the validity of the results stated therein, we can assume that all the schemes involved are smooth.} 
and let $\CN\subset T^*(\CY)$ a Zariski-closed subset. 

\medskip

Let $g:\CY\to \BA^1$ be a function with a non-vanishing differential. For a $k$-point $y$ of $\CY$, we denote by
$0\neq dg_y\in T_y^*(\CY)$ the differential of $g$ at $y$. 

\medskip

We shall say that $g$ is $\CN$-\emph{characteristic} 
at $y\in \CY$ if the element $dg_y\in T_y^*(\CY)$ belongs to $\CN$. 

\medskip

We shall say that $g$ is \emph{non-characteristic with respect to} $\CN$ if it is not $\CN$-characteristic at all $y\in \CY$. 

\medskip

We shall say that a point $y\in \CY$ is an \emph{isolated} point for the pair $(\CN,g)$ if:

\begin{itemize}

\item $g(y)=0$ and $g$ is $\CN$-characteristic at $y$;

\smallskip

\item There exists a Zariski neighborhood $y\in U\subset \CY$ such that $g$ is 
non-characteristic with respect to $\CN$ on $U-\{y\}$.

\end{itemize}  

\sssec{}

We record the following two geometric observations: 

\begin{lem}  \label{l:codiff isol}
Let $f':\CY'_1\to \CY'_2$ be a map between two schemes.
Let $\CN'_1\subset T^*(\CY'_1)$ and $\CN'_2\subset T^*(\CY'_2)$ be
closed irreducible subvarieties and let $\wt\CN'_1$ be an irreducible component
of $(df'{}^*)^{-1}(\CN'_1)$ that projects to $\CN'_2$. 
Let $g$ be a function on $\CY'_2$, and let $y'_2$ be an isolated
point for $(\CN'_2,g)$. Let $y'_1\in \CY'_1$ be a point such that:

\smallskip 

\noindent{\em(i)} $(dg_{y'_2},y'_1)\in \wt\CN'_1$ and it does not 
belong to other irreducible components of $(df'{}^*)^{-1}(\CN'_1)$.

\smallskip 

\noindent{\em(ii)} The map $\wt\CN'_1\to \CN'_2$ is quasi-finite at $(dg_{y'_2},y'_1)$. 

\smallskip

Then $y'_1$ is an isolated point for $(g\circ f',\CN'_1)$.

\end{lem} 

The proof of the lemma is straightforward. 

\begin{prop} \label{p:isolated}
Assume that $\dim(\CN)\leq \dim(\CY)$. Then for any non-zero vector 
$\xi\in T^*_y(\CY)\cap \CN$ there exists a function $g$ defined on a Zariski
neighborhood of $y$, such that $dg_y=\xi$ and $y$ is an isolated point for the pair $(\CN,g)$.
\end{prop} 

The proof of the proposition is given in \secref{ss:proof isolated}.

\begin{rem} In fact, we will only need \propref{p:isolated} in
the case when $(\xi,y)$ is a smooth point on $\CN$ (see \secref{sss:gen isol}), 
in which case the assertion is easy (at least when $\on{char}(k)\neq 2$):

\medskip

Linearizing, we can assume being given a finite-dimensional vector space $V$
and a subspace $$W\subset V\oplus V^\vee,$$
of dimension $\leq \dim(V)$, and we have to find a quadratic
form $Q$ on $V$, such that for the associated bilinear form, viewed as a map
$V\to V^\vee$, its graph is transversal to $W$.

\end{rem}

\sssec{} \label{sss:use non-char}

Given a function $g$, consider the vanishing cycles functor
$$\Phi_g:\Shv(\CY)\to \Shv(\CY_0),$$
where $\CY_0:=\CY\underset{\BA^1}\times \{0\}$. 

\medskip

The following results from the definition of singular support:

\begin{lem} \label{l:non-char}
Let $\CN$ be conical and suppose that $g$ is non-characteristic with respect to $\CN$. Then the functor $\Phi_g$
vanishes on $\Shv_\CN(\CY)$.
\end{lem}

%

\sssec{} \label{sss:isolated}

Let now $y\in \CY$ be isolated for $(\CN,g)$. 

\medskip

Then by \lemref{l:non-char}, for $\CF\in \Shv_\CN(\CY)$, the object 
$$\Phi_g(\CF)\in \Shv(\CY_0),$$
canonically splits as a direct sum
\begin{equation} \label{e:Phi decomp N}
\Phi_g(\CF)\simeq \Phi_{g,y}(\CF) \oplus \Phi_{g,y\on{-disj}}(\CF),
\end{equation} 
where $\Phi_{g,y}(\CF)$ is supported at $\{y\}$ and $\Phi_{g,y\on{-disj}}(\CF)$ is supported on a closed subset of
$\CY_0$ disjoint from the point $y$. 

\sssec{}

We have the following fundamental result:

\begin{thm} \label{t:Sai}
In the situation of \secref{sss:isolated}, if for $\CF\in \on{Perv}(\CY)\cap \Shv_\CN(\CY)$, we have $(dg_y,y)\in \on{SingSupp}(\CF)$, then 
$$\Phi_{g,y}(\CF) \neq 0.$$
\end{thm} 

In the context of \'etale sheaves, this theorem follows from \cite[Theorem 5.9 (combined with Proposition 5.14)]{Sai}. 

\medskip

For constructible Betti sheaves, it follows easily from the definition of singular support in \cite{KS}, under 
the additional assumption that $(dg_y,y)\in T^*(\CY)$ is a \emph{smooth} point of $\CN$ (which will be 
the case in our situation, see \secref{sss:Phi 1 sat} below. 

\medskip

In the case of regular holonomic D-modules, 
it follows from the Betti case, by Riemann-Hilbert (under the same smoothness assumption). 

\begin{rem}
The reason we did not formulate \thmref{t:sing supp dir im} for holonomic (but not necessarily regular)
D-modules is that we are not certain about the status of an analog of \thmref{t:Sai} in this context.
\end{rem} 

\ssec{Isolated points and vanishing cycles--beyond the constructible}

\sssec{}

Given $y\in \CY$ and a function $g$ with $g(y)=0$, let 
\begin{equation} \label{e:isol}
\Shv_{(g,y)\on{-isol}}(\CY) \subset \Shv(\CY)
\end{equation}
denote the full subcategory consisting of objects $\CF$ such that, whenever $\CN$ is a conical closed irreducible subvariety satisfying 
$$\CN \subset \on{SingSupp}(\CF) \text{ and } (dg_y,y)\in \CN,$$
then $y$ is isolated for $(\CN,g)$. 

\medskip

The subcategory \eqref{e:isol} is compatible with the t-structure on the ambient category.
In particular, $\Shv_{(g,y)\on{-isol}}(\CY)$ acquires a t-structure.

\sssec{} \label{sss:isol nearby}

We claim that the restriction of the functor $\Phi_g$ to $\Shv_{(g,y)\on{-isol}}(\CY)$ also splits canonically
as a direct sum
\begin{equation} \label{e:Phi decomp infinite}
\Phi_g\simeq \Phi_{g,y}\oplus \Phi_{g,y\on{-disj}},
\end{equation} 
where $\Phi_{g,y}$ maps to
$$\Vect_\sfe\simeq \Shv(\{y\})\subset  \Shv(\CY_0),$$
and the range of the functor $\Phi_{g,y\on{-disj}}$ consists of objects $\CF'$ such that for every $m$, 
every constructible sub-object of $H^m(\CF')$ is supported on a closed subset of $\CY_0$ that does
not contain $y$.

\sssec{}

Indeed, as in \secref{sss:sing supp schemes constr}, the category $\Shv_{(g,y)\on{-isol}}(\CY)$ 
is the left-completion of its full subcategory
$$\Shv_{(g,y)\on{-isol}}(\CY)^{\on{access}}\subset  \Shv_{(g,y)\on{-isol}}(\CY)$$
equal to the ind-completion of 
\begin{equation} \label{e:access union}
\underset{\CN}\cup\, \Shv_{\CN}(\CY)^{\on{constr}}
\end{equation}  
for closed conical $\CN \subset T^*(\CY)$ such that
$$(dg_y,y)\in \CN\, \Rightarrow\, \text{ $y$ is isolated for $(\CN,g)$}.$$

\medskip

Since the functor $\Phi_g$ is t-exact, and the target category $\Shv(\CY_0)$ is left-complete 
in its t-structure, we obtain that in order to construct a decomposition \eqref{e:Phi decomp infinite}, it suffices
to construct it after restricting to $\Shv_{(g,y)\on{-isol}}(\CY)^{\on{access}}$. The latter amounts
to a similar decomposition after restriction to \eqref{e:access union}, which, in turn, follows
from \eqref{e:Phi decomp N}. 

\sssec{}

From \thmref{t:Sai} we deduce the following: 

\begin{cor}   \label{c:Sai}
Let $\CF$ be an object of $\Shv_{(g,y)\on{-isol}}(\CY)$ such that $(dg_y,y)\in \on{SingSupp}(\CF)$.
Then $$\Phi_{g,y}(\CF)\neq 0.$$
\end{cor}

\begin{proof}

Since the functor $\Phi_g$ is t-exact (on all of $\Shv(\CY)$), so is its restriction to $\Shv_{(g,y)\on{-isol}}(\CY)$,
and hence so is $\Phi_{g,y}$. Hence, it suffices to show that for some $m$ and for some constructible sub-object $\CF'$
of $H^m(\CF)$. we have $\Phi_{g,y}(\CF')\neq 0$.

\medskip

Let $m$ be such that $H^m(\CF)$ contains a constructible sub-object $\CF'$ with 
$(dg_y,y)\in \on{SingSupp}(\CF')$. Then $\Phi_{g,y}(\CF')\neq 0$ by \thmref{t:Sai}.

\end{proof} 

\ssec{Construction of $(\bC,\Phi_1,\Phi_2)$} \label{ss:construct triple}

\sssec{}

Let $\CN_2$ be as in \secref{sss:choose N1}. Let $k'$ be the algebraic closure of the field 
of rational functions on $\CN_2$.

\medskip

Let us denote by 
$$\CY'_1,\,\,\CY'_2,\,\,\CN'_2\subset T^*(\CY'_2),$$ 
etc., the base change of all of our geometric objects to $k'$. 

\medskip

Let $\on{BC}$ denote the corresponding base change functors
$$\Shv(\CY_1)\to \Shv(\CY'_1),\,\, \Shv(\CY_2)\to \Shv(\CY'_2),$$
etc. 

\medskip

The functor $\on{BC}$ preserves singular support, in the sense that if for $\CF_i\in \Shv(\CY_i)$
we have $\on{SingSupp}(\CF_i)=\CN_i$, then 
$$\on{SingSupp}(\on{BC}(\CF_i))=\CN'_i.$$

\sssec{} \label{sss:gen isol}

The generic point of $\CN_2$ gives rise to a $k'$-point of $\CN'_2$; denote it by
$(\xi'_2,y'_2)$. 

\medskip

After shrinking $\CY'_2$ around $y'_2$, choose a function $g$ such that $dg_{y'_2}=\xi'_2$ and $y'_2$ is an isolated point for
$(\CN'_2,g)$. Such $g$ exists by \propref{p:isolated}\footnote{Note that the point $(\xi'_2,y'_2)$ is a smooth point on $\CN'_2$,
so we are applying an easy case of \propref{p:isolated}.}. 

\medskip

Set
$$\CY'_{2,0}:=\CY'_2\underset{\BA^1_{k'}}\times \{0\}.$$

Let 
$$\bC:=\underset{y'_2\in U\subset \CY'_{2,0}}{\on{colim}}\, \Shv(U),$$
where the colimit is taken over the poset of Zariski neighborhoods of $y'_2$ in $\CY'_{2,0}$. 

\sssec{}

We let $\Phi_2$ be the composition
$$\Shv(\CY_2) \overset{\on{BC}}\to \Shv(\CY'_2) 
\overset{\Phi_g}\to \Shv(\CY'_{2,0}) \to \bC.$$

\medskip

We claim that $\Phi_2$ satisfies Property (1) in \secref{sss:Phi req}. 

\medskip

Indeed, it suffices to show that if $Z$ is a \emph{Zariski-closed} subset of $T^*(\CY_2)$
of dimension equal to $\dim(\CY_2)$ and
$$\CN_2\not\subset Z \text{ and } \CF_2\in \Shv_{Z}(\CY_2),$$ 
then 
$$\Phi_g(\on{BC}(\CF_2))$$
vanishes on some neighborhood of $y'_2$. 

\medskip

However, for $Z$ as above, 
$$(\xi'_2,y'_2)\notin Z',$$
hence $g$ is not $Z'$-characteristic at $y'_2$, and hence is non-characteristic with respect to $Z'$
on some neighborhood $U$ of $y'_2$. 

\medskip

Since $\on{SingSupp}(\on{BC}(\CF_2))\subset Z'$, we obtain that
$$\Phi_g(\on{BC}(\CF_2))|_U=0$$
by \lemref{l:non-char}. 

\sssec{}

Consider the function $g\circ f':\CY'_1\to \BA^1$. Set 
$$\CY'_{1,0}:=\CY'_1\underset{\BA^1_{k'}}\times \{0\},$$
and let $f'_0$ denote the induced map $\CY'_{1,0}\to \CY'_{2,0}$.

\medskip

Recall the subscheme $\wt\CN_1$, see \secref{sss:choose N1}, and consider its base change
$$\wt\CN'_1\subset T^*(\CY'_2)\underset{\CY'_2}\times \CY'_1.$$

Let $(\xi'_2,y'_1)$ be one (out of the finite and non-empty set) of $k'$-points of $\wt\CN'_1$
that projects to $(\xi'_2,y'_2)$ along 
$$\wt\CN'_1\to \CN'_2.$$

\medskip

By condition (iib) in \thmref{t:sing supp dir im}, the point $(\xi'_2,y'_1)$ does not belong to
other irreducible components of $(df'{}^*)^{-1}(\CN'_1)$. Applying \lemref{l:codiff isol}, we obtain 
that the functor $$\on{BC}:\Shv(\CY_1)\to \Shv(\CY'_1)$$ 
maps
$$\Shv_{\wt\CN_1\on{-q.f.}}(\CY_1)\to \Shv_{(g\circ f',y'_1)\on{-isol}}(\CY'_1).$$

We let $\Phi_1$ be the composition
$$\Shv_{\wt\CN_1\on{-q.f.}}(\CY_1) \overset{\on{BC}}\longrightarrow 
\Shv_{(g\circ f',y'_1)\on{-isol}}(\CY'_1) \overset{\Phi_{g\circ f',y'_1}}\longrightarrow \Shv(\{y'_1\}) \simeq 
\Vect_\sfe \simeq \Shv(\{y'_2\}) \hookrightarrow \Shv(Y'_{2,0})\to \bC.$$

\sssec{} \label{sss:Phi 1 sat}

We claim that $\Phi_1$ satisfies Property (2) in \secref{sss:Phi req}. 

\medskip

Indeed, if $\CF'_1\in \Shv_{\wt\CN_1\on{-q.f.}}(\CY_1)$ is such that 
$$\CN_1\subset \on{SingSupp}(\CF'_1),$$
then
$$\CN'_1\subset \on{SingSupp}(\on{BC}(\CF'_1)),$$
and hence $(df'{}^*(\xi'_2),y'_1)\in  \on{SingSupp}(\on{BC}(\CF'_1))$. 

\medskip

Hence, $\Phi_{g\circ f',y'_1}(\on{BC}(\CF'_1))\neq 0$ by \corref{c:Sai}. 

\sssec{}

We will now prove that $\Phi_1$ is canonically a direct summand of $\Phi_2\circ f_*$, which 
is Property (3) in \secref{sss:Phi req}. 

\medskip

We now use the fact that $f$ is proper. This implies that we have a canonical isomorphism
$$\Phi_g\circ f'_*\simeq (f'_0)_*\circ \Phi_{g\circ f'},$$
where $f'_0$ is the induced map $\CY'_{1,0}\to \CY'_{2,0}$.

\medskip

Hence, we can rewrite the composition 
$$\Shv_{\wt\CN_1\on{-q.f.}}(\CY_1)\hookrightarrow \Shv(\CY_1) \overset{f_*}\longrightarrow \Shv(\CY_1)  \overset{\Phi_2}\longrightarrow \bC$$
as
$$\Shv_{\wt\CN_1\on{-q.f.}}(\CY_1) \overset{\on{BC}}\longrightarrow 
\Shv_{(g\circ f',y'_1)\on{-isol}}(\CY'_1)   \overset{\Phi_{g\circ f'}}\longrightarrow \Shv(\CY'_{1,0}) \overset{(f'_0)_*}\longrightarrow
\Shv(\CY'_{2,0})\to \bC.$$

By \secref{sss:isol nearby}, the above functor contains a direct summand isomorphic to
$$\Shv_{\wt\CN_1\on{-q.f.}}(\CY_1) \overset{\on{BC}}\longrightarrow 
\Shv_{(g\circ f',y'_1)\on{-isol}}(\CY'_1)   \overset{\Phi_{g\circ f',y'_1}}\longrightarrow \Shv(\{y'_1\})\hookrightarrow 
\Shv(\CY'_{1,0}) \overset{(f'_0)_*}\longrightarrow \Shv(\CY'_{2,0})\to \bC.$$

However, the latter functor is the same as $\Phi_1$. Indeed, the composition
$$\Shv(\{y'_1\})\hookrightarrow \Shv(\CY'_{1,0}) \overset{(f'_0)_*}\longrightarrow \Shv(\CY'_{2,0})$$
is the same as
$$\Shv(\{y'_1\}) \simeq \Vect_\sfe \simeq \Shv(\{y'_2\}) \hookrightarrow \Shv(Y'_{2,0}).$$

\ssec{Proof of \propref{p:isolated}} \label{ss:proof isolated}

\sssec{}

The assertion is local, so by embedding $\CY$ into a smooth scheme, we can assume that
$\CY$ itself is smooth.

\medskip

The proof proceeds by induction on $\dim(\CY)$. The base of induction is when $\dim(\CY)=1$,
which is easy and is left to the reader (cf. proof of \lemref{l:exists h} below).

\medskip

Without loss of generality, we assume that 
assume that $\CY$ admits a smooth map
$$f:\CY\to \CY'$$
of relative dimension $1$, and that $\CY$ and $\CY'$ are affine. Denote $y'=f(y)$.

\medskip 

Let $h$ be a function on $\CY$, to be chosen later. 
We define an automorphism
$\on{shr}_h:T^*(\CY)\to T^*(\CY)$ by
\[\on{shr}_h(\eta)=\eta+dh_y,\qquad y\in \CY,\,\eta\in T^*_y(\CY).\]
Let $\tau_h$ denote the composite map
$$T^*(\CY')\underset{\CY'}\times \CY\hookrightarrow T^*(\CY) \overset{\on{shr}_h}\longrightarrow T^*(\CY).$$
Clearly, $\tau_h$ is a closed embedding. 

\medskip

Consider also the projection 
$$(\on{id}\times f):T^*(\CY')\underset{\CY'}\times \CY\to T^*(\CY').$$ 
It is a smooth map of relative dimension one. 

\medskip

We now specify the choice of $h$. Choose $0\neq \xi'\in T^*_{y'}(\CY')$. 

\begin{lem} \label{l:exists h}
There exists a function $h$ on $\CY$ with $h(y)=0$ such that:

\begin{itemize}

\item $\xi=dh_y+df^*(\xi')$, in particular, $(\xi,y)\in \on{Im}(\tau_h)$; 

\item The dimension of $\dim(\tau_h^{-1}(\CN))\leq \dim(\CY)-1$;

\item The restriction of the map $(\on{id}\times f)$ to $\tau_h^{-1}(\CN)$ 
is quasi-finite at the point $\tau_h^{-1}(\xi,y)$. 
\end{itemize} 

\end{lem}

\medskip

We prove \lemref{l:exists h} below. Let us assume it, and perform the induction step.

\sssec{}

Let $h$ be as in \lemref{l:exists h}. 

\medskip

Let $\CN'$ be the closure of the image $(\on{id}\times f)(\tau_h^{-1}(\CN))$.
By the choice of $h$, we have $(\xi',y')\in\CN'$, $\dim(\CN')\le\dim(\CY')$, and 
\[\tau_h^{-1}(\CN)\to\CN'\]
given by $(\on{id}\times f)$ is quasi-finite at $(\xi',y)$. 

\medskip

Applying the induction hypothesis to $\CN'$ and $\xi'$, we can find a function $g'$ on $\CY'$
such that the point $y'$ is isolated for $(\CN',g')$. 

\medskip

Take $g:=g'\circ f+h$. It is easy to see that the point $y$ is isolated for $(\CN,g)$.

\qed[\propref{p:isolated}]

\sssec{Proof of \lemref{l:exists h}} Let us impose the following restrictions on $h$:

\medskip

First require that 
$$dh_y=\xi-df^*(\xi'),$$
i.e., we impose the first condition of \lemref{l:exists h}.

\medskip

Consider the relative cotangent bundle $T^*(\CY/\CY')$. Let 
$$r:T^*(\CY)\to T^*(\CY/\CY')$$
be the natural projection; $r$ is smooth of relative dimension $n-1$. 

\medskip

For an irreducible component $\CN_i$ on $\CN$, consider two cases. 

\medskip

\noindent Case (a): $r(\CN_i)=\{r(\xi,y)\}$. In this case, $\dim(\tau^{-1}_h(\CN_i))\le n-1$
for any $h$.

\medskip

\noindent Case (b): There exists a point
$(\eta_i,y_i)\in r(\CN_i)$, $\eta_i\in T^*_{y_i}(\CY/\CY')$ such that $(\eta_i,y_i)\neq r(\xi,y)$. In this case,
we require that $\eta_i\neq r(dh_{y_i})$.
This restriction guarantees that $\dim(\tau^{-1}_h(\CN_i))\le n-1$.

\medskip

Imposing the above condition for for all irreducible components $\CN_i$ in case (b), we obtain that $h$ 
satisfies the second condition of \lemref{l:exists h}.

\medskip

Since $\dim(\CN)\leq n$, there exists a point $y_0$ that lies in the same connected component of $f^{-1}(y')$
as $y$ and an element $\xi_0\in T^*_{y_0}(\CY)$ such that 
$(\xi_0+df^*(\xi'),y_0)\not\in\CN$. Moreover, we can choose the pair $(y_0,\xi_0)$ so that $y_0\ne y$ and also 
$y_0\ne y_i$ for all of the points $y_i$ from the second step, case (b). We require that
$dh_{y_0}=\xi_0$. This implies that $\tau_h^{-1}(\CN)$ does not
entirely contain the connected component of $(\on{id}\times f)^{-1}(\xi',y')$ to which the point
$(\xi',y)$ belongs, which is equivalent to the third condition of \lemref{l:exists h} (since the
morphism $(\on{id}\times f)$ has relative dimension $1$). 

\medskip

It is easy to see that it is possible to find a function $h$ subject to the above restrictions. 

\qed[\lemref{l:exists h}]


%
%
%
%
%
%
%
%


\end{document}